\documentclass{report}
\usepackage{amssymb}
\usepackage{amsfonts}

\usepackage{amsmath}

\newtheorem{theorem}{Theorem}

\newtheorem{corollary}[theorem]{Corollary}

\newtheorem{lemma}[theorem]{Lemma}

\newtheorem{proposition}[theorem]{Proposition}
\newtheorem{remark}[theorem]{Remark}

\newenvironment{proof}[1][Proof]{\noindent\textbf{#1.} }{\ \rule{0.5em}{0.5em}}
\input{tcilatex}

\begin{document}

\title{\textbf{\ Advances in Inequalities of the Schwarz, Gr\"{u}ss and
Bessel Type in Inner Product Spaces}}
\author{\textit{Sever Silvestru Dragomir} \\
School of Computer Science \& Mathematics,\\
Victoria University, \\
Melbourne, Australia}
\date{December 2003}
\maketitle

\begin{abstract}
The main aim of this monograph is to survey some recent results obtained by
the author related to reverses of the Schwarz, triangle and Bessel
inequalities. Some Gr\"{u}ss' type inequalities for orthonormal families of
vectors in real or complex inner product spaces are presented as well.
Generalizations of \ the Boas-Bellman, Bombieri, Selberg, Heilbronn and Pe%
\v{c}ari\'{c} inequalities for finite sequences of vectors that are not
necessarily orthogonal are also provided. Two extensions of the celebrated
Ostrowski's inequalities for sequences or real numbers and the
generalization of Wagner's inequality in inner product spaces are pointed
out. Finally, some Gr\"{u}ss type inequalities for n-tuples of vectors in
inner product spaces and their natural applications for the approximation of
the discrete Fourier and Mellin transforms are given as well.
\end{abstract}

\tableofcontents

\chapter*{Preface}

The theory of Hilbert spaces plays a central role in contemporary
mathematics with numerous applications for Linear Operators, Partial
Differential Equations, in Nonlinear Analysis, Approximation Theory,
Optimization Theory, Numerical Analysis, Probability Theory, Statistics and
other fields.

The Schwarz, triangle, Bessel, Gram and most recently, Gr\"{u}ss type
inequalities have been frequently used as powerful tools in obtaining bounds
or estimating the errors for various approximation formulae occurring in the
domains mentioned above. Therefore, any new advancement related to these
fundamental facts will have a flow of important consequences in the
mathematical fields where these inequalities have been used before.

The main aim of this monograph is to survey some recent results obtained by
the author related to reverses of the Schwarz, triangle and Bessel
inequalities. Some Gr\"{u}ss type inequalities for orthonormal families of
vectors in real or complex inner product spaces are presented as well.
Generalizations of\ the Boas-Bellman, Bombieri, Selberg, Heilbronn and Pe%
\v{c}ari\'{c} inequalities for finite sequences of vectors that are not
necessarily orthogonal are also provided. Two extensions of the celebrated
Ostrowski inequalities for sequences of real numbers and the generalization
of Wagner's inequality in inner product spaces are pointed out. Finally,
some Gr\"{u}ss type inequalities for $n$-tuples of vectors in inner product
spaces and their natural applications for the approximation of the discrete
Fourier and Mellin transforms are given as well.

The monograph may be used by researchers in different branches of
Mathematical and Functional Analysis where the theory of Hilbert spaces is
of relevance. Since it is self-contained and all the results are completely
proved, the work may be also used by graduate students interested in Theory
of Inequalities and its Applications.

Every section is followed by the necessary references and in this way each
section may be read independently of the others.

\vspace{2in}

\textit{The Author},

December, 2003.

\chapter{Reverses for the Schwarz Inequality}

\section{Introduction}

\pagenumbering{arabic}Let $H$ be a linear space over the real or complex
number field $\mathbb{K}$. The functional $\left\langle \cdot ,\cdot
\right\rangle :H\times H\rightarrow \mathbb{K}$ is called an \textit{inner
product }on $H$ if it satisfies the conditions

\begin{enumerate}
\item[(i)] $\left\langle x,x\right\rangle \geq 0$ for any $x\in H$ and $%
\left\langle x,x\right\rangle =0$ iff $x=0;$

\item[(ii)] $\left\langle \alpha x+\beta y,z\right\rangle =\alpha
\left\langle x,z\right\rangle +\beta \left\langle y,z\right\rangle $ for any 
$\alpha ,\beta \in \mathbb{K}$ and $x,y,z\in H;$

\item[(iii)] $\left\langle y,x\right\rangle =\overline{\left\langle
x,y\right\rangle }$ for any $x,y\in H.$
\end{enumerate}

A first fundamental consequence of the properties (i)-(iii) above, is the 
\textit{Schwarz inequality:}%
\begin{equation}
\left\vert \left\langle x,y\right\rangle \right\vert ^{2}\leq \left\langle
x,x\right\rangle \left\langle y,y\right\rangle  \label{CBS}
\end{equation}%
for any $x,y\in H.$ The equality holds in (\ref{CBS}) if and only if the
vectors $x$ and $y$ are \textit{linearly dependent,} i.e., there exists a
nonzero constant $\alpha \in \mathbb{K}$ so that $x=\alpha y.$

If we denote $\left\Vert x\right\Vert :=\sqrt{\left\langle x,x\right\rangle }%
,x\in H,$ then one may state the following properties

\begin{enumerate}
\item[(n)] $\left\Vert x\right\Vert \geq 0$ for any $x\in H$ and $\left\Vert
x\right\Vert =0$ iff $x=0;$

\item[(nn)] $\left\Vert \alpha x\right\Vert =\left\vert \alpha \right\vert
\left\Vert x\right\Vert $ for any $\alpha \in \mathbb{K}$ and $x\in H;$

\item[(nnn)] $\left\Vert x+y\right\Vert \leq \left\Vert x\right\Vert
+\left\Vert y\right\Vert $ for any $x,y\in H$ (the triangle inequality);
\end{enumerate}

i.e., $\left\Vert \cdot \right\Vert $ is a \textit{norm} on $H.$

In this chapter we present some recent reverse inequalities for the Schwarz
and the triangle inequalities. More precisely, we point out upper bounds for
the nonnegative quantities%
\begin{equation*}
\left\Vert x\right\Vert \left\Vert y\right\Vert -\left\vert \left\langle
x,y\right\rangle \right\vert ,\text{ \ }\left\Vert x\right\Vert
^{2}\left\Vert y\right\Vert ^{2}-\left\vert \left\langle x,y\right\rangle
\right\vert ^{2}
\end{equation*}%
and 
\begin{equation*}
\left\Vert x\right\Vert +\left\Vert y\right\Vert -\left\Vert x+y\right\Vert
\end{equation*}%
under various assumptions for the vectors $x,y\in H.$

If the vectors $x,y\in H$ are not \textit{orthogonal}, i.e., $\left\langle
x,y\right\rangle \neq 0,$ then some upper bounds for the supra-unitary
quantities%
\begin{equation*}
\frac{\left\Vert x\right\Vert \left\Vert y\right\Vert }{\left\vert
\left\langle x,y\right\rangle \right\vert },\ \ \frac{\left\Vert
x\right\Vert ^{2}\left\Vert y\right\Vert ^{2}}{\left\vert \left\langle
x,y\right\rangle \right\vert ^{2}}
\end{equation*}%
are provided as well.

\section{An Additive Reverse of the Schwarz Inequality}

\subsection{Introduction}

Let $\overline{\mathbf{a}}=\left( a_{1},\dots ,a_{n}\right) $ and $\overline{%
\mathbf{b}}=\left( b_{1},\dots ,b_{n}\right) $ be two positive $n-$\textit{%
tuples} with%
\begin{equation}
0<m_{1}\leq a_{i}\leq M_{1}<\infty \text{ and }0<m_{2}\leq b_{i}\leq
M_{2}<\infty ;  \label{1.1a.1}
\end{equation}%
for each \ $i\in \left\{ 1,\dots ,n\right\} ,$ and some constants $%
m_{1},m_{2},M_{1},M_{2}.$

The following reverses of the Cauchy-Bunyakovsky-Schwarz inequality for
positive sequences of real numbers are well known:

\begin{enumerate}
\item \textit{P\'{o}lya-Szeg\"{o}'s inequality} \cite{6b.1}%
\begin{equation*}
\frac{\sum_{k=1}^{n}a_{k}^{2}\sum_{k=1}^{n}b_{k}^{2}}{\left(
\sum_{k=1}^{n}a_{k}b_{k}\right) ^{2}}\leq \frac{1}{4}\left( \sqrt{\frac{%
M_{1}M_{2}}{m_{1}m_{2}}}+\sqrt{\frac{m_{1}m_{2}}{M_{1}M_{2}}}\right) ^{2}.
\end{equation*}

\item \textit{Shisha-Mond's inequality} \cite{7b.1}%
\begin{equation*}
\frac{\sum_{k=1}^{n}a_{k}^{2}}{\sum_{k=1}^{n}a_{k}b_{k}}-\frac{%
\sum_{k=1}^{n}a_{k}b_{k}}{\sum_{k=1}^{n}b_{k}^{2}}\leq \left[ \left( \frac{%
M_{1}}{m_{2}}\right) ^{\frac{1}{2}}-\left( \frac{m_{1}}{M_{2}}\right) ^{%
\frac{1}{2}}\right] ^{2}.
\end{equation*}

\item \textit{Ozeki's inequality} \cite{5b.1}%
\begin{equation*}
\sum_{k=1}^{n}a_{k}^{2}\sum_{k=1}^{n}b_{k}^{2}-\left(
\sum_{k=1}^{n}a_{k}b_{k}\right) ^{2}\leq \frac{n^{2}}{4}\left(
M_{1}M_{2}-m_{1}m_{2}\right) ^{2}.
\end{equation*}

\item \textit{Diaz-Metcalf's inequality} \cite{1b.1}%
\begin{equation*}
\sum_{k=1}^{n}b_{k}^{2}+\frac{m_{2}M_{2}}{m_{1}M_{1}}\sum_{k=1}^{n}a_{k}^{2}%
\leq \left( \frac{M_{2}}{m_{1}}+\frac{m_{2}}{M_{1}}\right)
\sum_{k=1}^{n}a_{k}b_{k}.
\end{equation*}
\end{enumerate}

If $\overline{\mathbf{w}}=\left( w_{1},\dots ,w_{n}\right) $ is a positive
sequence, then the following weighted inequalities also hold:

\begin{enumerate}
\item \textit{Cassel's inequality} \cite{8b.1}. If the positive real
sequences $\overline{\mathbf{a}}=\left( a_{1},\dots ,a_{n}\right) $ and $%
\overline{\mathbf{b}}=\left( b_{1},\dots ,b_{n}\right) $ satisfy the
condition 
\begin{equation}
0<m\leq \frac{a_{k}}{b_{k}}\leq M<\infty \text{ for each }k\in \left\{
1,\dots ,n\right\} ,  \label{1.5a.1}
\end{equation}%
then%
\begin{equation*}
\frac{\left( \sum_{k=1}^{n}w_{k}a_{k}^{2}\right) \left(
\sum_{k=1}^{n}w_{k}b_{k}^{2}\right) }{\left(
\sum_{k=1}^{n}w_{k}a_{k}b_{k}\right) ^{2}}\leq \frac{\left( M+m\right) ^{2}}{%
4mM}.
\end{equation*}

\item \textit{Greub-Reinboldt's inequality} \cite{2b.1}. We have%
\begin{equation*}
\left( \sum_{k=1}^{n}w_{k}a_{k}^{2}\right) \left(
\sum_{k=1}^{n}w_{k}b_{k}^{2}\right) \leq \frac{\left(
M_{1}M_{2}+m_{1}m_{2}\right) ^{2}}{4m_{1}m_{2}M_{1}M_{2}}\left(
\sum_{k=1}^{n}w_{k}a_{k}b_{k}\right) ^{2},
\end{equation*}%
provided $\overline{\mathbf{a}}=\left( a_{1},\dots ,a_{n}\right) $ and $%
\overline{\mathbf{b}}=\left( b_{1},\dots ,b_{n}\right) $ satisfy the
condition $\left( \ref{1.1a.1}\right) .$

\item \textit{Generalised Diaz-Metcalf's inequality} \cite{1b.1}, see also 
\cite[p. 123]{4b.1}. If $u,v\in \left[ 0,1\right] $ and $v\leq u,$ $u+v=1$
and $\left( \ref{1.5a.1}\right) $ holds, then one has the inequality%
\begin{equation*}
u\sum_{k=1}^{n}w_{k}b_{k}^{2}+vMm\sum_{k=1}^{n}w_{k}a_{k}^{2}\leq \left(
vm+uM\right) \sum_{k=1}^{n}w_{k}a_{k}b_{k}.
\end{equation*}

\item \textit{Klamkin-McLenaghan's inequality} \cite{3b.1}. If $\overline{%
\mathbf{a}},\overline{\mathbf{b}}$ satisfy (\ref{1.5a.1}), then%
\begin{multline}
\left( \sum_{i=1}^{n}w_{i}a_{i}^{2}\right) \left(
\sum_{i=1}^{n}w_{i}b_{i}^{2}\right) -\left(
\sum_{i=1}^{n}w_{i}a_{i}b_{i}\right) ^{2}  \label{1.9.1} \\
\leq \left( M^{\frac{1}{2}}-m^{\frac{1}{2}}\right)
^{2}\sum_{i=1}^{n}w_{i}a_{i}b_{i}\sum_{i=1}^{n}w_{i}a_{i}^{2}.
\end{multline}
\end{enumerate}

For other recent results providing discrete reverse inequalities, see the
recent monograph online \cite{1bab.1}.

In this section, by following \cite{NSSD}, we point out a new reverse of
Schwarz's inequality in real or complex inner product spaces. Particular
cases for isotonic linear functionals, integrals and sequences are also
given.

\subsection{An Additive Reverse Inequality}

The following reverse of Schwarz's inequality in inner product spaces holds 
\cite{NSSD}.

\begin{theorem}
\label{t2.1.1}Let $A,a\in \mathbb{K}$ $\left( \mathbb{K}=\mathbb{C},\mathbb{R%
}\right) $ and $x,y\in H.$ If%
\begin{equation}
\func{Re}\left\langle Ay-x,x-ay\right\rangle \geq 0,  \label{2.1.1}
\end{equation}%
or equivalently,%
\begin{equation}
\left\Vert x-\frac{a+A}{2}\cdot y\right\Vert \leq \frac{1}{2}\left\vert
A-a\right\vert \left\Vert y\right\Vert ,  \label{2.1.a.1}
\end{equation}%
holds, then one has the inequality%
\begin{equation}
0\leq \left\Vert x\right\Vert ^{2}\left\Vert y\right\Vert ^{2}-\left\vert
\left\langle x,y\right\rangle \right\vert ^{2}\leq \frac{1}{4}\left\vert
A-a\right\vert ^{2}\left\Vert y\right\Vert ^{4}.  \label{2.2.1}
\end{equation}%
The constant $\frac{1}{4}$ is sharp in (\ref{2.2.1}).
\end{theorem}

\begin{proof}
The equivalence between $\left( \ref{2.1.1}\right) $ and $\left( \ref%
{2.1.a.1}\right) $ can be easily proved, see for example \cite{1ba.1}.

Let us define%
\begin{equation*}
I_{1}:=\func{Re}\left[ \left( A\left\Vert y\right\Vert ^{2}-\left\vert
\left\langle x,y\right\rangle \right\vert \right) \left( \overline{%
\left\langle x,y\right\rangle }-\overline{a}\left\Vert y\right\Vert
^{2}\right) \right]
\end{equation*}%
and 
\begin{equation*}
I_{2}:=\left\Vert y\right\Vert ^{2}\func{Re}\left\langle
Ay-x,x-ay\right\rangle .
\end{equation*}%
Then%
\begin{equation*}
I_{1}=\left\Vert y\right\Vert ^{2}\func{Re}\left[ A\overline{\left\langle
x,y\right\rangle }+\overline{a}\left\langle x,y\right\rangle \right]
-\left\vert \left\langle x,y\right\rangle \right\vert ^{2}-\left\Vert
y\right\Vert ^{4}\func{Re}\left( A\overline{a}\right)
\end{equation*}%
and%
\begin{equation*}
I_{2}=\left\Vert y\right\Vert ^{2}\func{Re}\left[ A\overline{\left\langle
x,y\right\rangle }+\overline{a}\left\langle x,y\right\rangle \right]
-\left\Vert x\right\Vert ^{2}\left\Vert y\right\Vert ^{2}-\left\Vert
y\right\Vert ^{4}\func{Re}\left( A\overline{a}\right) ,
\end{equation*}%
which gives%
\begin{equation*}
I_{1}-I_{2}=\left\Vert x\right\Vert ^{2}\left\Vert y\right\Vert
^{2}-\left\vert \left\langle x,y\right\rangle \right\vert ^{2}
\end{equation*}%
for any $x,y\in H$ and $a,A\in \mathbb{K}$. This is an interesting identity
in itself as well.

If (\ref{2.1.1}) holds, then $I_{2}\geq 0$ and thus%
\begin{equation}
\left\Vert x\right\Vert ^{2}\left\Vert y\right\Vert ^{2}-\left\vert
\left\langle x,y\right\rangle \right\vert ^{2}\leq \func{Re}\left[ \left(
A\left\Vert y\right\Vert ^{2}-\left\vert \left\langle x,y\right\rangle
\right\vert \right) \left( \overline{\left\langle x,y\right\rangle }-%
\overline{a}\left\Vert y\right\Vert ^{2}\right) \right] .  \label{2.4.1}
\end{equation}%
Further, if we use the elementary inequality for $u,v\in \mathbb{K}$ $\left( 
\mathbb{K}=\mathbb{C},\mathbb{R}\right) $%
\begin{equation*}
\func{Re}\left[ u\overline{v}\right] \leq \frac{1}{4}\left\vert
u+v\right\vert ^{2},
\end{equation*}%
then we have, for 
\begin{equation*}
u:=A\left\Vert y\right\Vert ^{2}-\left\langle x,y\right\rangle ,\ \
v:=\left\langle x,y\right\rangle -a\left\Vert y\right\Vert ^{2},
\end{equation*}%
that%
\begin{equation}
\func{Re}\left[ \left( A\left\Vert y\right\Vert ^{2}-\left\vert \left\langle
x,y\right\rangle \right\vert \right) \left( \overline{\left\langle
x,y\right\rangle }-\overline{a}\left\Vert y\right\Vert ^{2}\right) \right]
\leq \frac{1}{4}\left\vert A-a\right\vert ^{2}\left\Vert y\right\Vert ^{4}.
\label{2.6.1}
\end{equation}%
Making use of the inequalities $\left( \ref{2.4.1}\right) $ and $\left( \ref%
{2.6.1}\right) ,$ we deduce $\left( \ref{2.2.1}\right) $.

Now, assume that (\ref{2.2.1}) holds with a constant $C>0,$ i.e.,%
\begin{equation}
\left\Vert x\right\Vert ^{2}\left\Vert y\right\Vert ^{2}-\left\vert
\left\langle x,y\right\rangle \right\vert ^{2}\leq C\left\vert
A-a\right\vert ^{2}\left\Vert y\right\Vert ^{4},  \label{2.61.1}
\end{equation}%
where $x,y,a,A$ satisfy (\ref{2.1.1}).

Consider $y\in H,$ $\left\Vert y\right\Vert =1,$ $a\neq A$ and $m\in H,$ $%
\left\Vert m\right\Vert =1$ with $m\perp y.$ Define%
\begin{equation*}
x:=\frac{A+a}{2}y+\frac{A-a}{2}m.
\end{equation*}%
Then%
\begin{equation*}
\left\langle Ay-x,x-ay\right\rangle =\left\vert \frac{A-a}{2}\right\vert
^{2}\left\langle y-m,y+m\right\rangle =0,
\end{equation*}%
and thus the condition (\ref{2.1.1}) is fulfilled. From (\ref{2.61.1}) we
deduce%
\begin{equation}
\left\Vert \frac{A+a}{2}y+\frac{A-a}{2}m\right\Vert ^{2}-\left\vert
\left\langle \frac{A+a}{2}y+\frac{A-a}{2}m,y\right\rangle \right\vert
^{2}\leq C\left\vert A-a\right\vert ^{2},  \label{2.7.1}
\end{equation}%
and since%
\begin{equation*}
\left\Vert \frac{A+a}{2}y+\frac{A-a}{2}m\right\Vert ^{2}=\left\vert \frac{A+a%
}{2}\right\vert ^{2}+\left\vert \frac{A-a}{2}\right\vert ^{2}
\end{equation*}%
and%
\begin{equation*}
\left\vert \left\langle \frac{A+a}{2}y+\frac{A-a}{2}m,y\right\rangle
\right\vert ^{2}=\left\vert \frac{A+a}{2}\right\vert ^{2}
\end{equation*}%
then, by (\ref{2.7.1}), we obtain%
\begin{equation*}
\frac{\left\vert A-a\right\vert ^{2}}{4}\leq C\left\vert A-a\right\vert ^{2},
\end{equation*}%
which gives $C\geq \frac{1}{4},$ and the theorem is completely proved.
\end{proof}

\subsection{Applications for Isotonic Linear Functionals}

Let $F\left( T\right) $ be an algebra of real functions defined on $T$ and $%
L $ a subclass of $F\left( T\right) $ satisfying the conditions:

\begin{enumerate}
\item[(i)] $f,g\in L$ implies $f+g\in L;$

\item[(ii)] $f\in L,$ $\alpha \in \mathbb{R}$ implies $\alpha f\in L.$
\end{enumerate}

A functional $A$ defined on $L$ is an \textit{isotonic linear functional} on 
$L$ provided that

\begin{enumerate}
\item[(a)] $A\left( \alpha f+\beta g\right) =\alpha A\left( f\right) +\beta
A\left( g\right) $ for all $\alpha ,\beta \in \mathbb{R}$ and $f,g\in L;$

\item[(aa)] $f\geq g,$ that is, $f\left( t\right) \geq g\left( t\right) $
for all $t\in T,$ implies $A\left( f\right) \geq A\left( g\right) .$
\end{enumerate}

The functional $A$ is \textit{normalised} on $L,$ provided that $\mathbf{1}%
\in L,$ i.e., $\mathbf{1}\left( t\right) =1$ for all $t\in T,$ implies $%
A\left( \mathbf{1}\right) =1.$

Usual examples of isotonic linear functionals are integrals, sums, etc.

Now, suppose that $h\in F\left( T\right) ,$ $h\geq 0$ is given and satisfies
the properties that $fgh\in L,$ $fh\in L,$ $gh\in L$ for all $f,g\in L.$ For
a given isotonic linear functional $A:L\rightarrow \mathbb{R}$ with $A\left(
h\right) >0,$ define the mapping $\left( \cdot ,\cdot \right) _{A,h}:L\times
L\rightarrow \mathbb{R}$ by%
\begin{equation*}
\left( f,g\right) _{A,h}:=\frac{A\left( fgh\right) }{A\left( h\right) }.
\end{equation*}

This functional satisfies the following properties:

\begin{enumerate}
\item[(s)] $\left( f,f\right) _{A,h}\geq 0$ for all $f\in L;$

\item[(ss)] $\left( \alpha f+\beta g,k\right) _{A,h}=\alpha \left(
f,k\right) _{A,h}+\beta \left( g,k\right) _{A,h}$ for all $f,g,k\in L$ and $%
\alpha ,\beta \in \mathbb{R}$;

\item[(sss)] $\left( f,g\right) _{A,h}=\left( g,f\right) _{A,h}$ for all $%
f,g\in L.$
\end{enumerate}

The following reverse of Schwarz's inequality for positive linear
functionals holds \cite{NSSD}.

\begin{proposition}
\label{p3.1.1}Let $f,g,h\in F\left( T\right) $ be such that $fgh\in L,$ $%
f^{2}h\in L,$ $g^{2}h\in L.$ If $m,M$ are real numbers such that%
\begin{equation}
mg\leq f\leq Mg\text{ on }F\left( T\right) ,  \label{3.21.1}
\end{equation}%
then for any isotonic linear functional $A:L\rightarrow \mathbb{R}$ with $%
A\left( h\right) >0$ we have the inequality%
\begin{equation}
0\leq A\left( hf^{2}\right) A\left( hg^{2}\right) -\left[ A\left( hfg\right) %
\right] ^{2}\leq \frac{1}{4}\left( M-m\right) ^{2}A^{2}\left( hg^{2}\right) .
\label{3.3.1}
\end{equation}%
The constant $\frac{1}{4}$ in (\ref{3.3.1}) is sharp.
\end{proposition}

\begin{proof}
We observe that 
\begin{equation*}
\left( Mg-f,f-mg\right) _{A,h}=A\left[ h\left( Mg-f\right) \left(
f-mg\right) \right] \geq 0.
\end{equation*}%
Applying Theorem \ref{t2.1.1} for $\left( \cdot ,\cdot \right) _{A,h},$ we
get%
\begin{equation*}
0\leq \left( f,f\right) _{A,h}\left( g,g\right) _{A,h}-\left( f,g\right)
_{A,h}^{2}\leq \frac{1}{4}\left( M-m\right) ^{2}\left( g,g\right) _{A,h}^{2},
\end{equation*}%
which is clearly equivalent to (\ref{3.3.1}).
\end{proof}

The following corollary holds.

\begin{corollary}
\label{c3.2.1}Let $f,g\in F\left( T\right) $ such that $fg,$ $f^{2},g^{2}\in
F\left( T\right) .$ If $m,M$ are real numbers such that (\ref{3.21.1})
holds, then%
\begin{equation}
0\leq A\left( f^{2}\right) A\left( g^{2}\right) -A^{2}\left( fg\right) \leq 
\frac{1}{4}\left( M-m\right) ^{2}A^{2}\left( g^{2}\right) .  \label{3.4.1}
\end{equation}%
The constant $\frac{1}{4}$ is sharp in (\ref{3.4.1}).
\end{corollary}

\begin{remark}
\label{r3.3.1}The condition (\ref{3.21.1}) may be replaced with the weaker
assumption%
\begin{equation*}
\left( Mg-f,f-mg\right) _{A,h}\geq 0.
\end{equation*}
\end{remark}

\subsection{Applications for Integrals}

Let $\left( \Omega ,\Sigma ,\mu \right) $ be a measure space consisting of a
set $\Omega ,$ $\Sigma $ a $\sigma -$algebra of subsets of $\Omega $ and $%
\mu $ a countably additive and positive measure on $\Sigma $ with values in $%
\mathbb{R}\cup \left\{ \infty \right\} .$

Denote by $L_{\rho }^{2}\left( \Omega ,\mathbb{K}\right) $ the Hilbert space
of all $\mathbb{K}$-valued functions $f$ defined on $\Omega $ that are $%
2-\rho -$integrable on $\Omega ,$ i.e., $\int_{\Omega }\rho \left( t\right)
\left\vert f\left( s\right) \right\vert ^{2}d\mu \left( s\right) <\infty ,$
where $\rho :\Omega \rightarrow \lbrack 0,\infty )$ is a measurable function
on $\Omega .$

The following proposition contains a reverse of the weighted
Cauchy-Bunyakovsky-Schwarz's integral inequality \cite{NSSD}.

\begin{proposition}
\label{p4.1.1}Let $A,a\in \mathbb{K}$ $\left( \mathbb{K}=\mathbb{C},\mathbb{R%
}\right) $ and $f,g\in L_{\rho }^{2}\left( \Omega ,\mathbb{K}\right) .$ If 
\begin{equation}
\int_{\Omega }\func{Re}\left[ \left( Ag\left( s\right) -f\left( s\right)
\right) \left( \overline{f\left( s\right) }-\overline{a}\text{ }\overline{g}%
\left( s\right) \right) \right] \rho \left( s\right) d\mu \left( s\right)
\geq 0  \label{4.1.1}
\end{equation}%
or equivalently,%
\begin{equation*}
\int_{\Omega }\rho \left( s\right) \left\vert f\left( s\right) -\frac{a+A}{2}%
g\left( s\right) \right\vert ^{2}d\mu \left( s\right) \leq \frac{1}{4}%
\left\vert A-a\right\vert ^{2}\int_{\Omega }\rho \left( s\right) \left\vert
g\left( s\right) \right\vert ^{2}d\mu \left( s\right) ,
\end{equation*}%
holds, then one has the inequality%
\begin{align*}
0& \leq \int_{\Omega }\rho \left( s\right) \left\vert f\left( s\right)
\right\vert ^{2}d\mu \left( s\right) \int_{\Omega }\rho \left( s\right)
\left\vert g\left( s\right) \right\vert ^{2}d\mu \left( s\right) -\left\vert
\int_{\Omega }\rho \left( s\right) f\left( s\right) \overline{g\left(
s\right) }d\mu \left( s\right) \right\vert ^{2} \\
& \leq \frac{1}{4}\left\vert A-a\right\vert ^{2}\left( \int_{\Omega }\rho
\left( s\right) \left\vert g\left( s\right) \right\vert ^{2}d\mu \left(
s\right) \right) ^{2}.
\end{align*}%
The constant $\frac{1}{4}$ is best possible.
\end{proposition}

\begin{proof}
Follows by Theorem \ref{t2.1.1} applied for the inner product $\left\langle
\cdot ,\cdot \right\rangle _{\rho }:=L_{\rho }^{2}\left( \Omega ,\mathbb{K}%
\right) \times L_{\rho }^{2}\left( \Omega ,\mathbb{K}\right) \rightarrow 
\mathbb{K}$,%
\begin{equation*}
\left\langle f,g\right\rangle _{\rho }:=\int_{\Omega }\rho \left( s\right)
f\left( s\right) \overline{g\left( s\right) }d\mu \left( s\right) .
\end{equation*}
\end{proof}

\begin{remark}
\label{r4.2.1}A sufficient condition for (\ref{4.1.1}) to hold is 
\begin{equation*}
\func{Re}\left[ \left( Ag\left( s\right) -f\left( s\right) \right) \left( 
\overline{f\left( s\right) }-\overline{a}\text{ }\overline{g}\left( s\right)
\right) \right] \geq 0,\text{ \ for \ }\mu -\text{a.e. }s\in \Omega .
\end{equation*}
\end{remark}

In the particular case $\rho =1,$ we have the following reverse of the
Cauchy-Bunyakovsky-Schwarz inequality.

\begin{corollary}
\label{c4.3.1}Let $a,A\in \mathbb{K}$ $\left( \mathbb{K}=\mathbb{C},\mathbb{R%
}\right) $ and $f,g\in L^{2}\left( \Omega ,\mathbb{K}\right) .$ If%
\begin{equation}
\int_{\Omega }\func{Re}\left[ \left( Ag\left( s\right) -f\left( s\right)
\right) \left( \overline{f\left( s\right) }-\overline{a}\text{ }\overline{g}%
\left( s\right) \right) \right] d\mu \left( s\right) \geq 0,  \label{4.4.1}
\end{equation}%
or equivalently,%
\begin{equation*}
\int_{\Omega }\left\vert f\left( s\right) -\frac{a+A}{2}g\left( s\right)
\right\vert ^{2}d\mu \left( s\right) \leq \frac{1}{4}\left\vert
A-a\right\vert ^{2}\int_{\Omega }\left\vert g\left( s\right) \right\vert
^{2}d\mu \left( s\right) ,
\end{equation*}%
holds, then one has the inequality%
\begin{align*}
0& \leq \int_{\Omega }\left\vert f\left( s\right) \right\vert ^{2}d\mu
\left( s\right) \int_{\Omega }\left\vert g\left( s\right) \right\vert
^{2}d\mu \left( s\right) -\left\vert \int_{\Omega }f\left( s\right) 
\overline{g\left( s\right) }d\mu \left( s\right) \right\vert ^{2} \\
& \leq \frac{1}{4}\left\vert A-a\right\vert ^{2}\left( \int_{\Omega
}\left\vert g\left( s\right) \right\vert ^{2}d\mu \left( s\right) \right)
^{2}.
\end{align*}%
The constant $\frac{1}{4}$ is best possible
\end{corollary}

\begin{remark}
\label{r4.4.1}If $\mathbb{K}=\mathbb{R},$ then a sufficient condition for
either (\ref{4.1.1}) or (\ref{4.4.1}) \ to hold is 
\begin{equation*}
ag\left( s\right) \leq f\left( s\right) \leq Ag\left( s\right) ,\text{ \ for
\ }\mu -\text{a.e. }s\in \Omega ,
\end{equation*}%
where, in this case, $a,A\in \mathbb{R}$ with $A>a>0.$
\end{remark}

\subsection{Applications for Sequences}

For a given sequence $\left( w_{i}\right) _{i\in \mathbb{N}}$ of nonnegative
real numbers, consider the Hilbert space $\ell _{w}^{2}\left( \mathbb{K}%
\right) ,$ $\left( \mathbb{K}=\mathbb{C},\mathbb{R}\right) ,$ where 
\begin{equation*}
\ell _{w}^{2}\left( \mathbb{K}\right) :=\left\{ \overline{\mathbf{x}}=\left(
x_{i}\right) _{i\in \mathbb{N}}\subset \mathbb{K}\left\vert
\sum_{i=0}^{\infty }w_{i}\left\vert x_{i}\right\vert ^{2}<\infty \right.
\right\} .
\end{equation*}

The following proposition that provides a reverse of the weighted
Cauchy-Bunyakovsky-Schwarz inequality for complex numbers holds.

\begin{proposition}
\label{p5.1.1}Let $a,A\in \mathbb{K}$ and $\overline{\mathbf{x}},\overline{%
\mathbf{y}}\in \ell _{w}^{2}\left( \mathbb{K}\right) .$ If%
\begin{equation}
\sum_{i=0}^{\infty }w_{i}\func{Re}\left[ \left( Ay_{i}-x_{i}\right) \left( 
\overline{x_{i}}-\overline{a}\text{ }\overline{y_{i}}\right) \right] \geq 0,
\label{5.2.1}
\end{equation}%
then one has the inequality%
\begin{equation*}
0\leq \sum_{i=0}^{\infty }w_{i}\left\vert x_{i}\right\vert
^{2}\sum_{i=0}^{\infty }w_{i}\left\vert y_{i}\right\vert ^{2}-\left\vert
\sum_{i=0}^{\infty }w_{i}x_{i}\overline{y_{i}}\right\vert ^{2}\leq \frac{1}{4%
}\left\vert A-a\right\vert ^{2}\left( \sum_{i=0}^{\infty }w_{i}\left\vert
y_{i}\right\vert ^{2}\right) ^{2}.
\end{equation*}%
The constant $\frac{1}{4}$ is sharp.
\end{proposition}

\begin{proof}
Follows by Theorem \ref{t2.1.1} applied for the inner product $\left\langle
\cdot ,\cdot \right\rangle _{w}:\ell _{w}^{2}\left( \mathbb{K}\right) \times
\ell _{w}^{2}\left( \mathbb{K}\right) \rightarrow \mathbb{K}$,%
\begin{equation*}
\left\langle \overline{\mathbf{x}},\overline{\mathbf{y}}\right\rangle
_{w}:=\sum_{i=0}^{\infty }w_{i}x_{i}\overline{y_{i}}.
\end{equation*}
\end{proof}

\begin{remark}
\label{r5.2.1}A sufficient condition for (\ref{5.2.1}) to hold is%
\begin{equation*}
\func{Re}\left[ \left( Ay_{i}-x_{i}\right) \left( \overline{x_{i}}-\overline{%
a}\overline{y_{i}}\right) \right] \geq 0,\text{ \ for all }i\in \mathbb{N}.
\end{equation*}
\end{remark}

In the particular case $w_{i}=1,$ $i\in \mathbb{N}$, we have the following
reverse of the Cauchy-Bunyakovsky-Schwarz inequality.

\begin{corollary}
\label{c5.3.1}Let $a,A\in \mathbb{K}$ $\left( \mathbb{K}=\mathbb{C},\mathbb{R%
}\right) $ and $\overline{\mathbf{x}},\overline{\mathbf{y}}\in \ell
^{2}\left( \mathbb{K}\right) .$ If%
\begin{equation}
\sum_{i=0}^{\infty }\func{Re}\left[ \left( Ay_{i}-x_{i}\right) \left( 
\overline{x_{i}}-\overline{a}\overline{y_{i}}\right) \right] \geq 0,
\label{5.5.1}
\end{equation}%
then one has the inequality%
\begin{equation*}
0\leq \sum_{i=0}^{\infty }\left\vert x_{i}\right\vert ^{2}\sum_{i=0}^{\infty
}\left\vert y_{i}\right\vert ^{2}-\left\vert \sum_{i=0}^{\infty }x_{i}%
\overline{y_{i}}\right\vert ^{2}\leq \frac{1}{4}\left\vert A-a\right\vert
^{2}\left( \sum_{i=0}^{\infty }\left\vert y_{i}\right\vert ^{2}\right) ^{2}.
\end{equation*}
\end{corollary}

\begin{remark}
\label{r5.4.1}If $\mathbb{K}=\mathbb{R}$, then a sufficient condition for
either (\ref{5.2.1}) or (\ref{5.5.1}) to hold is%
\begin{equation*}
ay_{i}\leq x_{i}\leq Ay_{i}\text{ \ for each \ }i\in \mathbb{N},
\end{equation*}%
with $A>a>0.$
\end{remark}

\newpage

\section{A Generalisation of the Cassels and Greub-Reinboldt Inequalities}

\subsection{Introduction}

The following result was proved by J.W.S. Cassels in 1951 (see Appendix 1 of 
\cite{8b.2}).

\begin{theorem}
\label{t1.1.2}Let $\overline{\mathbf{a}}=\left( a_{1},\dots ,a_{n}\right) ,$ 
$\overline{\mathbf{b}}=\left( b_{1},\dots ,b_{n}\right) $ be sequences of
positive real numbers and $\overline{\mathbf{w}}=\left( w_{1},\dots
,w_{n}\right) $ a sequence of nonnegative real numbers. Suppose that%
\begin{equation}
m=\min_{i=\overline{1,n}}\left\{ \frac{a_{i}}{b_{i}}\right\} \text{ \ and \ }%
M=\max_{i=\overline{1,n}}\left\{ \frac{a_{i}}{b_{i}}\right\} .  \label{1.1.2}
\end{equation}%
Then one has the inequality%
\begin{equation}
\frac{\sum_{i=1}^{n}w_{i}a_{i}^{2}\sum_{i=1}^{n}w_{i}b_{i}^{2}}{\left(
\sum_{i=1}^{n}w_{i}a_{i}b_{i}\right) ^{2}}\leq \frac{\left( m+M\right) ^{2}}{%
4mM}.  \label{1.2.2}
\end{equation}%
The equality holds in (\ref{1.2.2}) when $w_{1}=\frac{1}{a_{1}b_{1}},$ $%
w_{n}=\frac{1}{a_{n}b_{n}},$ $w_{2}=\cdots =w_{n-1}=0,$ $m=\frac{a_{n}}{b_{1}%
}$ and $M=\frac{a_{1}}{b_{n}}.$
\end{theorem}

If one assumes that $0<a\leq a_{i}\leq A<\infty $ and $0<b\leq b_{i}\leq
B<\infty $ for each $i\in \left\{ 1,\dots ,n\right\} ,$ then by (\ref{1.2.2}%
) we may obtain \textit{Greub-Reinboldt's inequality} \cite{2b.2}%
\begin{equation*}
\frac{\sum_{i=1}^{n}w_{i}a_{i}^{2}\sum_{i=1}^{n}w_{i}b_{i}^{2}}{\left(
\sum_{i=1}^{n}w_{i}a_{i}b_{i}\right) ^{2}}\leq \frac{\left( ab+AB\right) ^{2}%
}{4abAB}.
\end{equation*}%
The following \textquotedblleft unweighted\textquotedblright\ Cassels'
inequality also holds%
\begin{equation*}
\frac{\sum_{i=1}^{n}a_{i}^{2}\sum_{i=1}^{n}b_{i}^{2}}{\left(
\sum_{i=1}^{n}a_{i}b_{i}\right) ^{2}}\leq \frac{\left( m+M\right) ^{2}}{4mM},
\end{equation*}%
provided $\overline{\mathbf{a}}$ and $\overline{\mathbf{b}}$ satisfy (\ref%
{1.1.2}). This inequality will produce the well known \textit{P\'{o}lya-Szeg%
\"{o} inequality} \cite[pp. 57, 213-114]{5b.2}, \cite[pp. 71-72, 253-255]%
{3b.2}:%
\begin{equation*}
\frac{\sum_{i=1}^{n}a_{i}^{2}\sum_{i=1}^{n}b_{i}^{2}}{\left(
\sum_{i=1}^{n}a_{i}b_{i}\right) ^{2}}\leq \frac{\left( ab+AB\right) ^{2}}{%
4abAB},
\end{equation*}%
provided $0<a\leq a_{i}\leq A<\infty $ and $0<b\leq b_{i}\leq B<\infty $ for
each $i\in \left\{ 1,\dots ,n\right\} .$

In \cite{4ba.2}, C.P. Niculescu proved, amongst others, the following
generalisation of Cassels' inequality:

\begin{theorem}
\label{t1.2.2}Let $E$ be a vector space endowed with a Hermitian product $%
\left\langle \cdot ,\cdot \right\rangle .$ Then%
\begin{equation}
\frac{\func{Re}\left\langle x,y\right\rangle }{\left\langle x,x\right\rangle
^{\frac{1}{2}}\left\langle y,y\right\rangle ^{\frac{1}{2}}}\geq \frac{2}{%
\sqrt{\frac{\omega }{\Omega }}+\sqrt{\frac{\Omega }{\omega }}}  \label{1.6.2}
\end{equation}%
for every $x,y\in E$ and every $\omega ,\Omega >0$ for which $\func{Re}%
\left\langle x-\omega y,x-\Omega y\right\rangle \leq 0.$
\end{theorem}

For other reverses of the Cauchy-Bunyakovsky-Schwarz inequality, see the
references \cite{1b.2} -- \cite{8b.2}.

In this section, by following \cite{1NSSD}, we establish a generalisation of
(\ref{1.6.2}) for complex numbers $\omega $ and $\Omega $ for which $\func{Re%
}\left( \overline{\omega }\Omega \right) >0.$ Applications for isotonic
linear functionals, integrals and sequences are also given.

\subsection{An Inequality in Real or Complex Inner Product Spaces}

The following reverse of Schwarz's inequality in inner product spaces holds 
\cite{1NSSD}.

\begin{theorem}
\label{t2.1.2}Let $a,A\in \mathbb{K}$ $\left( \mathbb{K}=\mathbb{C},\mathbb{R%
}\right) $ so that $\func{Re}\left( \overline{a}A\right) >0.$ If $x,y\in H$
are such that%
\begin{equation}
\func{Re}\left\langle Ay-x,x-ay\right\rangle \geq 0,  \label{2.1.2}
\end{equation}%
then one has the inequality%
\begin{equation}
\left\Vert x\right\Vert \left\Vert y\right\Vert \leq \frac{1}{2}\cdot \frac{%
\func{Re}\left[ A\overline{\left\langle x,y\right\rangle }+\overline{a}%
\left\langle x,y\right\rangle \right] }{\left[ \func{Re}\left( \overline{a}%
A\right) \right] ^{\frac{1}{2}}}\leq \frac{1}{2}\cdot \frac{\left\vert
A\right\vert +\left\vert a\right\vert }{\left[ \func{Re}\left( \overline{a}%
A\right) \right] ^{\frac{1}{2}}}\left\vert \left\langle x,y\right\rangle
\right\vert .  \label{2.2.2}
\end{equation}%
The constant $\frac{1}{2}$ is sharp in both inequalities.
\end{theorem}

\begin{proof}
We have, obviously, that%
\begin{align*}
I& :=\func{Re}\left\langle Ay-x,x-ay\right\rangle \\
& =\func{Re}\left[ A\overline{\left\langle x,y\right\rangle }+\overline{a}%
\left\langle x,y\right\rangle \right] -\left\Vert x\right\Vert ^{2}-\left[ 
\func{Re}\left( \overline{a}A\right) \right] \left\Vert y\right\Vert ^{2}
\end{align*}%
and, thus, by (\ref{2.1.2}), one has%
\begin{equation*}
\left\Vert x\right\Vert ^{2}+\left[ \func{Re}\left( \overline{a}A\right) %
\right] \cdot \left\Vert y\right\Vert ^{2}\leq \func{Re}\left[ A\overline{%
\left\langle x,y\right\rangle }+\overline{a}\left\langle x,y\right\rangle %
\right] ,
\end{equation*}%
which gives%
\begin{equation}
\frac{1}{\left[ \func{Re}\left( \overline{a}A\right) \right] ^{\frac{1}{2}}}%
\left\Vert x\right\Vert ^{2}+\left[ \func{Re}\left( \overline{a}A\right) %
\right] ^{\frac{1}{2}}\left\Vert y\right\Vert ^{2}\leq \frac{\func{Re}\left[
A\overline{\left\langle x,y\right\rangle }+\overline{a}\left\langle
x,y\right\rangle \right] }{\left[ \func{Re}\left( \overline{a}A\right) %
\right] ^{\frac{1}{2}}}.  \label{2.3.2}
\end{equation}%
On the other hand, by the elementary inequality%
\begin{equation*}
\alpha p^{2}+\frac{1}{\alpha }q^{2}\geq 2pq,
\end{equation*}%
valid for $p,q\geq 0$ and $\alpha >0,$ we deduce%
\begin{equation}
2\left\Vert x\right\Vert \left\Vert y\right\Vert \leq \frac{1}{\left[ \func{%
Re}\left( \overline{a}A\right) \right] ^{\frac{1}{2}}}\left\Vert
x\right\Vert ^{2}+\left[ \func{Re}\left( \overline{a}A\right) \right] ^{%
\frac{1}{2}}\left\Vert y\right\Vert ^{2}.  \label{2.4.2}
\end{equation}%
Utilizing (\ref{2.3.2}) and (\ref{2.4.2}) we deduce the first part of (\ref%
{2.2.2}).

The second part is obvious by the fact that for $z\in \mathbb{C},$ $%
\left\vert \func{Re}\left( z\right) \right\vert \leq \left\vert z\right\vert
.$

Now, assume that the first inequality in (\ref{2.2.2}) holds with a constant 
$c>0,$ i.e., 
\begin{equation}
\left\Vert x\right\Vert \left\Vert y\right\Vert \leq c\frac{\func{Re}\left[ A%
\overline{\left\langle x,y\right\rangle }+\overline{a}\left\langle
x,y\right\rangle \right] }{\left[ \func{Re}\left( \overline{a}A\right) %
\right] ^{\frac{1}{2}}},  \label{2.5.2}
\end{equation}%
where $a,A,x$ and $y$ satisfy (\ref{2.2.2}).

If we choose $a=A=1,$ $y=x\neq 0,$ then obviously (\ref{2.1.2}) holds and
from (\ref{2.5.2}) we obtain%
\begin{equation*}
\left\Vert x\right\Vert ^{2}\leq 2c\left\Vert x\right\Vert ^{2},
\end{equation*}%
giving $c\geq \frac{1}{2}.$

The theorem is completely proved.
\end{proof}

The following corollary is a natural consequence of the above theorem \cite%
{1NSSD}.

\begin{corollary}
\label{c2.2.2}Let $m,M>0.$ If $x,y\in H$ are such that%
\begin{equation*}
\func{Re}\left\langle My-x,x-my\right\rangle \geq 0,
\end{equation*}%
then one has the inequality%
\begin{equation}
\left\Vert x\right\Vert \left\Vert y\right\Vert \leq \frac{1}{2}\cdot \frac{%
M+m}{\sqrt{mM}}\func{Re}\left\langle x,y\right\rangle \leq \frac{1}{2}\cdot 
\frac{M+m}{\sqrt{mM}}\left\vert \left\langle x,y\right\rangle \right\vert .
\label{2.7.2}
\end{equation}%
The constant $\frac{1}{2}$ is sharp in (\ref{2.7.2}).
\end{corollary}

\begin{remark}
\label{r2.2.1.2}The inequality (\ref{2.7.2}) is equivalent to Niculescu's
inequality (\ref{1.6.2}).
\end{remark}

The following corollary is also obvious \cite{1NSSD}.

\begin{corollary}
\label{c2.3.2}With the assumptions of Corollary \ref{c2.2.2}, we have%
\begin{align}
0& \leq \left\Vert x\right\Vert \left\Vert y\right\Vert -\left\vert
\left\langle x,y\right\rangle \right\vert \leq \left\Vert x\right\Vert
\left\Vert y\right\Vert -\func{Re}\left\langle x,y\right\rangle
\label{2.8.2} \\
& \leq \frac{\left( \sqrt{M}-\sqrt{m}\right) ^{2}}{2\sqrt{mM}}\func{Re}%
\left\langle x,y\right\rangle \leq \frac{\left( \sqrt{M}-\sqrt{m}\right) ^{2}%
}{2\sqrt{mM}}\left\vert \left\langle x,y\right\rangle \right\vert  \notag
\end{align}%
and%
\begin{align}
0& \leq \left\Vert x\right\Vert ^{2}\left\Vert y\right\Vert ^{2}-\left\vert
\left\langle x,y\right\rangle \right\vert ^{2}\leq \left\Vert x\right\Vert
^{2}\left\Vert y\right\Vert ^{2}-\left[ \func{Re}\left\langle
x,y\right\rangle \right] ^{2}  \label{2.9.2} \\
& \leq \frac{\left( M-m\right) ^{2}}{4mM}\left[ \func{Re}\left\langle
x,y\right\rangle \right] ^{2}\leq \frac{\left( M-m\right) ^{2}}{4mM}%
\left\vert \left\langle x,y\right\rangle \right\vert ^{2}.  \notag
\end{align}%
The constants $\frac{1}{2}$ and $\frac{1}{4}$ are sharp.
\end{corollary}

\begin{proof}
If we subtract $\func{Re}\left\langle x,y\right\rangle \geq 0$ from the
first inequality in (\ref{2.7.2}), we get%
\begin{align*}
\left\Vert x\right\Vert \left\Vert y\right\Vert -\func{Re}\left\langle
x,y\right\rangle & \leq \left( \frac{1}{2}\cdot \frac{M+m}{\sqrt{mM}}%
-1\right) \func{Re}\left\langle x,y\right\rangle \\
& =\frac{\left( \sqrt{M}-\sqrt{m}\right) ^{2}}{2\sqrt{mM}}\func{Re}%
\left\langle x,y\right\rangle
\end{align*}%
which proves the third inequality in (\ref{2.8.2}). The other ones are
obvious.

Now, if we square the first inequality in (\ref{2.7.2}) and then subtract $%
\left[ \func{Re}\left\langle x,y\right\rangle \right] ^{2},$ we get%
\begin{align*}
\left\Vert x\right\Vert ^{2}\left\Vert y\right\Vert ^{2}-\left[ \func{Re}%
\left\langle x,y\right\rangle \right] ^{2}& \leq \left[ \frac{\left(
M+m\right) ^{2}}{4mM}-1\right] \left[ \func{Re}\left\langle x,y\right\rangle %
\right] ^{2} \\
& =\frac{\left( M-m\right) ^{2}}{4mM}\left[ \func{Re}\left\langle
x,y\right\rangle \right] ^{2}
\end{align*}%
which proves the third inequality in (\ref{2.9.2}). The other ones are
obvious.
\end{proof}

\subsection{Applications for Isotonic Linear Functionals}

The following proposition holds \cite{1NSSD}.

\begin{proposition}
\label{p3.1.2}Let $f,g,h\in F\left( T\right) $ be such that $fgh\in L,$ $%
f^{2}h\in L,$ $g^{2}h\in L.$ If $m,M>0$ are such that%
\begin{equation}
mg\leq f\leq Mg\text{ on }F\left( T\right) ,  \label{3.2.2}
\end{equation}%
then for any isotonic linear functional $A:L\rightarrow \mathbb{R}$ with $%
A\left( h\right) >0,$ we have the inequality%
\begin{equation}
1\leq \frac{A\left( f^{2}h\right) A\left( g^{2}h\right) }{A^{2}\left(
fgh\right) }\leq \frac{\left( M+m\right) ^{2}}{4mM}.  \label{3.3.2}
\end{equation}%
The constant $\frac{1}{4}$ in (\ref{3.3.2}) is sharp.
\end{proposition}

\begin{proof}
We observe that%
\begin{equation*}
\left( Mg-f,f-mg\right) _{A,h}=A\left[ h\left( Mg-f\right) \left(
f-mg\right) \right] \geq 0.
\end{equation*}%
Applying Corollary \ref{c2.2.2} for $\left( \cdot ,\cdot \right) _{A,h}$ we
get%
\begin{equation*}
1\leq \frac{\left( f,f\right) _{A,h}\left( g,g\right) _{A,h}}{\left(
f,g\right) _{A,h}^{2}}\leq \frac{\left( M+m\right) ^{2}}{4mM},
\end{equation*}%
which is clearly equivalent to (\ref{3.3.2}).
\end{proof}

The following additive versions of (\ref{3.3.2}) also hold \cite{1NSSD}.

\begin{corollary}
\label{c3.2.2}With the assumption in Proposition \ref{p3.1.2}, one has%
\begin{align*}
0& \leq \left[ A\left( f^{2}h\right) A\left( g^{2}h\right) \right] ^{\frac{1%
}{2}}-A\left( hfg\right) \\
& \leq \frac{\left( \sqrt{M}-\sqrt{m}\right) ^{2}}{2\sqrt{mM}}A\left(
hfg\right)
\end{align*}%
and 
\begin{align*}
0& \leq A\left( f^{2}h\right) A\left( g^{2}h\right) -A^{2}\left( fgh\right)
\\
& \leq \frac{\left( M-m\right) ^{2}}{4mM}A^{2}\left( fgh\right) .
\end{align*}%
The constants $\frac{1}{2}$ and $\frac{1}{4}$ are sharp.
\end{corollary}

\begin{remark}
\label{r3.3.2}The condition (\ref{3.2.2}) may be replaced with the weaker
assumption%
\begin{equation}
\left( Mg-f,f-mg\right) _{A,h}\geq 0.  \label{3.6.2}
\end{equation}
\end{remark}

\begin{remark}
\label{r3.4.2}With the assumption (\ref{3.2.2}) or (\ref{3.6.2}) and if $%
f,g\in F\left( T\right) $ with $fg,f^{2},g^{2}\in L,$ then one has the
inequalities%
\begin{equation*}
1\leq \frac{A\left( f^{2}\right) A\left( g^{2}\right) }{A^{2}\left(
fg\right) }\leq \frac{\left( M+m\right) ^{2}}{4mM},
\end{equation*}%
\begin{align*}
0& \leq \left[ A\left( f^{2}\right) A\left( g^{2}\right) \right] ^{\frac{1}{2%
}}-A\left( fg\right) \\
& \leq \frac{\left( \sqrt{M}-\sqrt{m}\right) ^{2}}{2\sqrt{mM}}A\left(
fg\right)
\end{align*}%
and%
\begin{equation*}
0\leq A\left( f^{2}\right) A\left( g^{2}\right) -A^{2}\left( fg\right) \leq 
\frac{\left( M-m\right) ^{2}}{4mM}A^{2}\left( fg\right) .
\end{equation*}
\end{remark}

\subsection{Applications for Integrals}

The following proposition contains a reverse of the weighted
Cauchy-Bunyakovsky-Schwarz integral inequality.

\begin{proposition}
\label{p4.1.2}Let $A,a\in \mathbb{K}$ $\left( \mathbb{K}=\mathbb{C},\mathbb{R%
}\right) $ with $\func{Re}\left( \overline{a}A\right) >0$ and $f,g\in
L_{\rho }^{2}\left( \Omega ,\mathbb{K}\right) .$ If 
\begin{equation}
\int_{\Omega }\func{Re}\left[ \left( Ag\left( s\right) -f\left( s\right)
\right) \left( \overline{f\left( s\right) }-\overline{a}\text{ }\overline{g}%
\left( s\right) \right) \right] \rho \left( s\right) d\mu \left( s\right)
\geq 0,  \label{4.1.2}
\end{equation}%
then one has the inequality%
\begin{align}
& \left[ \int_{\Omega }\rho \left( s\right) \left\vert f\left( s\right)
\right\vert ^{2}d\mu \left( s\right) \int_{\Omega }\rho \left( s\right)
\left\vert g\left( s\right) \right\vert ^{2}d\mu \left( s\right) \right] ^{%
\frac{1}{2}}  \label{4.2.2} \\
& \leq \frac{1}{2}\cdot \frac{\int_{\Omega }\rho \left( s\right) \func{Re}%
\left[ A\overline{f\left( s\right) }g\left( s\right) +\overline{a}f\left(
s\right) \overline{g\left( s\right) }\right] d\mu \left( s\right) }{\left[ 
\func{Re}\left( \overline{a}A\right) \right] ^{\frac{1}{2}}}  \notag \\
& \leq \frac{1}{2}\cdot \frac{\left\vert A\right\vert +\left\vert
a\right\vert }{\left[ \func{Re}\left( \overline{a}A\right) \right] ^{\frac{1%
}{2}}}\left\vert \int_{\Omega }\rho \left( s\right) f\left( s\right) 
\overline{g\left( s\right) }d\mu \left( s\right) \right\vert .  \notag
\end{align}%
The constant $\frac{1}{2}$ is sharp in (\ref{4.2.2}).
\end{proposition}

\begin{proof}
Follows by Theorem \ref{t2.1.2} applied for the inner product $\left\langle
\cdot ,\cdot \right\rangle _{\rho }:=L_{\rho }^{2}\left( \Omega ,\mathbb{K}%
\right) \times L_{\rho }^{2}\left( \Omega ,\mathbb{K}\right) \rightarrow 
\mathbb{K}$,%
\begin{equation*}
\left\langle f,g\right\rangle :=\int_{\Omega }\rho \left( s\right) f\left(
s\right) \overline{g\left( s\right) }d\mu \left( s\right) .
\end{equation*}
\end{proof}

\begin{remark}
\label{r4.2.2}A sufficient condition for (\ref{4.1.2}) to hold is 
\begin{equation*}
\func{Re}\left[ \left( Ag\left( s\right) -f\left( s\right) \right) \left( 
\overline{f\left( s\right) }-\overline{a}\overline{g\left( s\right) }\right) %
\right] \geq 0,\text{ \ for }\mu \text{-a.e. \ }s\in \Omega .
\end{equation*}
\end{remark}

In the particular case $\rho =1,$ we have the following result.

\begin{corollary}
\label{c4.3.2}Let $a,A\in \mathbb{K}$ $\left( \mathbb{K}=\mathbb{C},\mathbb{R%
}\right) $ with $\func{Re}\left( \overline{a}A\right) >0$ and $f,g\in
L^{2}\left( \Omega ,\mathbb{K}\right) .$ If%
\begin{equation}
\int_{\Omega }\func{Re}\left[ \left( Ag\left( s\right) -f\left( s\right)
\right) \left( \overline{f\left( s\right) }-\overline{a}\overline{g\left(
s\right) }\right) \right] d\mu \left( s\right) \geq 0,  \label{4.4.2}
\end{equation}%
then one has the inequality%
\begin{align*}
& \left[ \int_{\Omega }\left\vert f\left( s\right) \right\vert ^{2}d\mu
\left( s\right) \int_{\Omega }\left\vert g\left( s\right) \right\vert
^{2}d\mu \left( s\right) \right] ^{\frac{1}{2}} \\
& \leq \frac{1}{2}\cdot \frac{\int_{\Omega }\func{Re}\left[ A\overline{%
f\left( s\right) }g\left( s\right) +\overline{a}f\left( s\right) \overline{%
g\left( s\right) }\right] d\mu \left( s\right) }{\left[ \func{Re}\left( 
\overline{a}A\right) \right] ^{\frac{1}{2}}} \\
& \leq \frac{1}{2}\cdot \frac{\left\vert A\right\vert +\left\vert
a\right\vert }{\left[ \func{Re}\left( \overline{a}A\right) \right] ^{\frac{1%
}{2}}}\left\vert \int_{\Omega }f\left( s\right) \overline{g\left( s\right) }%
d\mu \left( s\right) \right\vert .
\end{align*}
\end{corollary}

\begin{remark}
\label{r4.4.2}If $\mathbb{K}=\mathbb{R},$ then a sufficient condition for
either (\ref{4.1.2}) or (\ref{4.4.2}) to hold is 
\begin{equation*}
ag\left( s\right) \leq f\left( s\right) \leq Ag\left( s\right) ,\text{ \ for
\ }\mu \text{-a.e. }s\in \Omega ,
\end{equation*}%
where, in this case, $a,A\in \mathbb{R}$ with $A>a>0.$
\end{remark}

If $a,A$ are real positive constants, then the following proposition holds.

\begin{proposition}
\label{p4.5.2}Let $m,M>0.$ If $f,g\in L_{\rho }^{2}\left( \Omega ,\mathbb{K}%
\right) $ such that%
\begin{equation*}
\int_{\Omega }\rho \left( s\right) \func{Re}\left[ \left( Mg\left( s\right)
-f\left( s\right) \right) \left( \overline{f\left( s\right) }-m\overline{g}%
\left( s\right) \right) \right] d\mu \left( s\right) \geq 0,
\end{equation*}%
then one has the inequality%
\begin{multline*}
\left[ \int_{\Omega }\rho \left( s\right) \left\vert f\left( s\right)
\right\vert ^{2}d\mu \left( s\right) \int_{\Omega }\rho \left( s\right)
\left\vert g\left( s\right) \right\vert ^{2}d\mu \left( s\right) \right] ^{%
\frac{1}{2}} \\
\leq \frac{1}{2}\cdot \frac{M+m}{\sqrt{mM}}\int_{\Omega }\rho \left(
s\right) \func{Re}\left[ f\left( s\right) \overline{g\left( s\right) }\right]
d\mu \left( s\right) .
\end{multline*}
\end{proposition}

The proof follows by Corollary \ref{c2.2.2} applied for the inner product 
\begin{equation*}
\left\langle f,g\right\rangle _{\rho }:=\int_{\Omega }\rho \left( s\right)
f\left( s\right) \overline{g\left( s\right) }d\mu \left( s\right) .
\end{equation*}

The following additive versions also hold \cite{1NSSD}.

\begin{corollary}
\label{c4.6.2}With the assumptions in Proposition \ref{p4.5.2}, one has the
inequalities%
\begin{align*}
0& \leq \left[ \int_{\Omega }\rho \left( s\right) \left\vert f\left(
s\right) \right\vert ^{2}d\mu \left( s\right) \int_{\Omega }\rho \left(
s\right) \left\vert g\left( s\right) \right\vert ^{2}d\mu \left( s\right) %
\right] ^{\frac{1}{2}} \\
& \ \ \ \ \ \ \ \ \ \ \ \ \ \ \ \ \ \ \ \ \ \ \ \ \ -\int_{\Omega }\rho
\left( s\right) \func{Re}\left[ f\left( s\right) \overline{g\left( s\right) }%
\right] d\mu \left( s\right) \\
& \leq \frac{\left( \sqrt{M}-\sqrt{m}\right) ^{2}}{2\sqrt{mM}}\int_{\Omega
}\rho \left( s\right) \func{Re}\left[ f\left( s\right) \overline{g\left(
s\right) }\right] d\mu \left( s\right)
\end{align*}%
and 
\begin{align*}
0& \leq \int_{\Omega }\rho \left( s\right) \left\vert f\left( s\right)
\right\vert ^{2}d\mu \left( s\right) \int_{\Omega }\rho \left( s\right)
\left\vert g\left( s\right) \right\vert ^{2}d\mu \left( s\right) \\
& \ \ \ \ \ \ \ \ \ \ \ \ \ \ \ \ \ \ \ \ \ \ \ \ -\left( \int_{\Omega }\rho
\left( s\right) \func{Re}\left[ f\left( s\right) \overline{g\left( s\right) }%
\right] d\mu \left( s\right) \right) ^{2} \\
& \leq \frac{\left( M-m\right) ^{2}}{4mM}\left( \int_{\Omega }\rho \left(
s\right) \func{Re}\left[ f\left( s\right) \overline{g\left( s\right) }\right]
d\mu \left( s\right) \right) ^{2}.
\end{align*}
\end{corollary}

\begin{remark}
\label{r4.7.2}If $\mathbb{K}=\mathbb{R}$, a sufficient condition for (\ref%
{4.1.2}) to hold is 
\begin{equation*}
mg\left( s\right) \leq f\left( s\right) \leq Mg\left( s\right) ,\text{ \ for
\ }\mu \text{-a.e. }s\in \Omega ,
\end{equation*}%
where $M>m>0.$
\end{remark}

\subsection{Applications for Sequences}

For a given sequence $\left( w_{i}\right) _{i\in \mathbb{N}}$ of nonnegative
real numbers, consider the Hilbert space $\ell _{w}^{2}\left( \mathbb{K}%
\right) ,$ $\left( \mathbb{K}=\mathbb{C},\mathbb{R}\right) ,$ where 
\begin{equation*}
\ell _{w}^{2}\left( \mathbb{K}\right) :=\left\{ \overline{\mathbf{x}}=\left(
x_{i}\right) _{i\in \mathbb{N}}\subset \mathbb{K}\left\vert
\sum_{i=0}^{\infty }w_{i}\left\vert x_{i}\right\vert ^{2}<\infty \right.
\right\} .
\end{equation*}

The following proposition that provides a reverse of the weighted
Cauchy-Bunyakovsky-Schwarz inequality for complex numbers holds \cite{1NSSD}.

\begin{proposition}
\label{p5.1.2}Let $a,A\in \mathbb{K}$ with $\func{Re}\left( \overline{a}%
A\right) >0$ and $\overline{\mathbf{x}},\overline{\mathbf{y}}\in \ell
_{w}^{2}\left( \mathbb{K}\right) .$ If%
\begin{equation}
\sum_{i=0}^{\infty }w_{i}\func{Re}\left[ \left( Ay_{i}-x_{i}\right) \left( 
\overline{x_{i}}-\overline{a}\overline{y_{i}}\right) \right] \geq 0,
\label{5.2.2}
\end{equation}%
then one has the inequality%
\begin{align}
\left[ \sum_{i=0}^{\infty }w_{i}\left\vert x_{i}\right\vert
^{2}\sum_{i=0}^{\infty }w_{i}\left\vert y_{i}\right\vert ^{2}\right] ^{\frac{%
1}{2}}& \leq \frac{1}{2}\cdot \frac{\sum_{i=0}^{\infty }w_{i}\func{Re}\left[
A\overline{x_{i}}y_{i}+\overline{a}x_{i}\overline{y_{i}}\right] }{\left[ 
\func{Re}\left( \overline{a}A\right) \right] ^{\frac{1}{2}}}  \label{5.3.2}
\\
& \leq \frac{1}{2}\cdot \frac{\left\vert A\right\vert +\left\vert
a\right\vert }{\left[ \func{Re}\left( \overline{a}A\right) \right] ^{\frac{1%
}{2}}}\left\vert \sum_{i=0}^{\infty }w_{i}x_{i}\overline{y_{i}}\right\vert .
\notag
\end{align}%
The constant $\frac{1}{2}$ is sharp in (\ref{5.3.2}).
\end{proposition}

\begin{proof}
Follows by Theorem \ref{t2.1.2} applied for the inner product $\left\langle
\cdot ,\cdot \right\rangle :\ell _{w}^{2}\left( \mathbb{K}\right) \times
\ell _{w}^{2}\left( \mathbb{K}\right) \rightarrow \mathbb{K}$,%
\begin{equation*}
\left\langle \overline{\mathbf{x}},\overline{\mathbf{y}}\right\rangle
_{w}:=\sum_{i=0}^{\infty }w_{i}x_{i}\overline{y_{i}}.
\end{equation*}
\end{proof}

\begin{remark}
\label{r5.2.2}A sufficient condition for (\ref{5.2.2}) to hold is%
\begin{equation}
\func{Re}\left[ \left( Ay_{i}-x_{i}\right) \left( \overline{x_{i}}-\overline{%
a}\overline{y_{i}}\right) \right] \geq 0\text{ \ for all }i\in \mathbb{N}.
\label{5.4.2}
\end{equation}
\end{remark}

In the particular case $\rho =1,$ we have the following result.

\begin{corollary}
\label{c5.3.2}Let $a,A\in \mathbb{K}$ with $\func{Re}\left( \overline{a}%
A\right) >0$ and $\overline{\mathbf{x}},\overline{\mathbf{y}}\in \ell
^{2}\left( \mathbb{K}\right) .$ If%
\begin{equation*}
\sum_{i=0}^{\infty }\func{Re}\left[ \left( Ay_{i}-x_{i}\right) \left( 
\overline{x_{i}}-\overline{a}\overline{y_{i}}\right) \right] \geq 0,
\end{equation*}%
then one has the inequality%
\begin{align*}
\left[ \sum_{i=0}^{\infty }\left\vert x_{i}\right\vert
^{2}\sum_{i=0}^{\infty }\left\vert y_{i}\right\vert ^{2}\right] ^{\frac{1}{2}%
}& \leq \frac{1}{2}\cdot \frac{\sum_{i=0}^{\infty }\func{Re}\left[ A%
\overline{x_{i}}y_{i}+\overline{a}x_{i}\overline{y_{i}}\right] }{\left[ 
\func{Re}\left( \overline{a}A\right) \right] ^{\frac{1}{2}}} \\
& \leq \frac{1}{2}\cdot \frac{\left\vert A\right\vert +\left\vert
a\right\vert }{\left[ \func{Re}\left( \overline{a}A\right) \right] ^{\frac{1%
}{2}}}\left\vert \sum_{i=0}^{\infty }x_{i}\overline{y_{i}}\right\vert .
\end{align*}
\end{corollary}

\begin{remark}
\label{r5.4.2}If $\mathbb{K}=\mathbb{R}$, then a sufficient condition for
either (\ref{5.2.2}) or (\ref{5.4.2}) to hold is%
\begin{equation*}
ay_{i}\leq x_{i}\leq Ay_{i}\text{ \ for each \ }i\in \left\{ 1,\dots
,n\right\} ,
\end{equation*}%
where, in this case, $a,A\in \mathbb{R}$ with $A>a>0.$
\end{remark}

For $a=m,$ $A=M$, then the following proposition also holds.

\begin{proposition}
\label{p5.5.2}Let $m,M>0.$ If $\overline{\mathbf{x}},\overline{\mathbf{y}}%
\in \ell _{w}^{2}\left( \mathbb{K}\right) $ such that 
\begin{equation}
\sum_{i=0}^{\infty }w_{i}\func{Re}\left[ \left( My_{i}-x_{i}\right) \left( 
\overline{x_{i}}-m\overline{y_{i}}\right) \right] \geq 0,  \label{5.7.2}
\end{equation}%
then one has the inequality%
\begin{equation*}
\left[ \sum_{i=0}^{\infty }w_{i}\left\vert x_{i}\right\vert
^{2}\sum_{i=0}^{\infty }w_{i}\left\vert y_{i}\right\vert ^{2}\right] ^{\frac{%
1}{2}}\leq \frac{1}{2}\cdot \frac{M+m}{\sqrt{mM}}\sum_{i=0}^{\infty }w_{i}%
\func{Re}\left( x_{i}\overline{y_{i}}\right) .
\end{equation*}
\end{proposition}

The proof follows by Corollary \ref{c2.2.2} applied for the inner product 
\begin{equation*}
\left\langle \overline{\mathbf{x}},\overline{\mathbf{y}}\right\rangle
_{w}:=\sum_{i=0}^{\infty }w_{i}x_{i},\overline{y_{i}}.
\end{equation*}

The following additive version also holds \cite{1NSSD}.

\begin{corollary}
\label{c5.6.2}With the assumptions in Proposition \ref{p5.5.2}, one has the
inequalities%
\begin{align*}
0& \leq \left[ \sum_{i=0}^{\infty }w_{i}\left\vert x_{i}\right\vert
^{2}\sum_{i=0}^{\infty }w_{i}\left\vert y_{i}\right\vert ^{2}\right] ^{\frac{%
1}{2}}-\sum_{i=0}^{\infty }w_{i}\func{Re}\left( x_{i}\overline{y_{i}}\right)
\\
& \leq \frac{\left( \sqrt{M}-\sqrt{m}\right) ^{2}}{2\sqrt{mM}}%
\sum_{i=0}^{\infty }w_{i}\func{Re}\left( x_{i}\overline{y_{i}}\right)
\end{align*}%
and%
\begin{align*}
0& \leq \sum_{i=0}^{\infty }w_{i}\left\vert x_{i}\right\vert
^{2}\sum_{i=0}^{\infty }w_{i}\left\vert y_{i}\right\vert ^{2}-\left[
\sum_{i=0}^{\infty }w_{i}\func{Re}\left( x_{i}\overline{y_{i}}\right) \right]
^{2} \\
& \leq \frac{\left( M-m\right) ^{2}}{4mM}\left[ \sum_{i=0}^{\infty }w_{i}%
\func{Re}\left( x_{i}\overline{y_{i}}\right) \right] ^{2}.
\end{align*}
\end{corollary}

\begin{remark}
\label{r5.7.2}If $\mathbb{K}=\mathbb{R}$, a sufficient condition for (\ref%
{5.7.2}) to hold is%
\begin{equation*}
my_{i}\leq x_{i}\leq My_{i}\text{ \ for each \ }i\in \mathbb{N},
\end{equation*}%
where $M>m>0.$
\end{remark}

\newpage

\section{Quadratic Reverses of Schwarz's Inequality}

\subsection{Two Better Reverse Inequalities}

It has been proven in \cite{2NSSD}, that%
\begin{equation}
0\leq \left\Vert x\right\Vert ^{2}-\left\vert \left\langle x,e\right\rangle
\right\vert ^{2}\leq \frac{1}{4}\left\vert \phi -\varphi \right\vert
^{2}-\left\vert \frac{\phi +\varphi }{2}-\left\langle x,e\right\rangle
\right\vert ^{2},  \label{1.7.3}
\end{equation}%
provided, either 
\begin{equation}
\func{Re}\left\langle \phi e-x,x-\varphi e\right\rangle \geq 0,
\label{1.8.3}
\end{equation}%
or equivalently,%
\begin{equation}
\left\Vert x-\frac{\phi +\varphi }{2}e\right\Vert \leq \frac{1}{2}\left\vert
\phi -\varphi \right\vert ,  \label{1.9.3}
\end{equation}%
holds, where $e=H,$ $\left\Vert e\right\Vert =1.$ The constant $\frac{1}{4}$
in (\ref{1.7.3}) is best possible.

If we choose $e=\frac{y}{\left\Vert y\right\Vert },$ $\phi =\Gamma
\left\Vert y\right\Vert ,$ $\varphi =\gamma \left\Vert y\right\Vert $ $%
\left( y\neq 0\right) ,$ $\Gamma ,\gamma \in \mathbb{K}$, then by (\ref%
{1.8.3}) and (\ref{1.9.3}) we have,%
\begin{equation}
\func{Re}\left\langle \Gamma y-x,x-\gamma y\right\rangle \geq 0,
\label{1.10.3}
\end{equation}%
or equivalently,%
\begin{equation}
\left\Vert x-\frac{\Gamma +\gamma }{2}y\right\Vert \leq \frac{1}{2}%
\left\vert \Gamma -\gamma \right\vert \left\Vert y\right\Vert ,
\label{1.11.3}
\end{equation}%
implying the following reverse of Schwarz's inequality:%
\begin{equation}
0\leq \left\Vert x\right\Vert ^{2}\left\Vert y\right\Vert ^{2}-\left\vert
\left\langle x,y\right\rangle \right\vert ^{2}\leq \frac{1}{4}\left\vert
\Gamma -\gamma \right\vert ^{2}\left\Vert y\right\Vert ^{4}-\left\vert \frac{%
\Gamma +\gamma }{2}\left\Vert y\right\Vert ^{2}-\left\langle
x,y\right\rangle \right\vert ^{2}.  \label{1.12.3}
\end{equation}%
The constant $\frac{1}{4}$ in (\ref{1.12.3}) is sharp.

Note that, this inequality is an improvement of (\ref{2.2.1}), but it may
not be very convenient for applications.

In \cite{5NSSD}, it has also been proven that 
\begin{equation}
0\leq \left\Vert x\right\Vert ^{2}-\left\vert \left\langle x,e\right\rangle
\right\vert ^{2}\leq \frac{1}{4}\left\vert \phi -\varphi \right\vert ^{2}-%
\func{Re}\left\langle \Phi e-x,x-\varphi e\right\rangle  \label{1.12.3a}
\end{equation}%
provided either (\ref{1.8.3}) or (\ref{1.9.3}) holds true.

If we make the same choice for $e,\Phi $ and $\varphi $ as above, then we
deduce the inequality%
\begin{equation}
0\leq \left\Vert x\right\Vert ^{2}\left\Vert y\right\Vert ^{2}-\left\vert
\left\langle x,y\right\rangle \right\vert ^{2}\leq \frac{1}{4}\left\vert
\Gamma -\gamma \right\vert ^{2}\left\Vert y\right\Vert ^{4}-\left\Vert
y\right\Vert ^{2}\func{Re}\left\langle \Gamma y-x,x-\gamma y\right\rangle
\label{1.12.3b}
\end{equation}%
provided either (\ref{1.10.3}) or (\ref{1.11.3}) holds true.

The constant $\frac{1}{4}$ is best possible in (\ref{1.12.3b}).

\subsection{A Reverse of Schwarz's Inequality Under More General Assumptions}

The following result holds \cite{3NSSD}.

\begin{theorem}
\label{t2.1.3}Let $\left( H;\left\langle \cdot ,\cdot \right\rangle \right) $
be an inner product space over the real or complex number field $\mathbb{K}$ 
$\left( \mathbb{K}=\mathbb{R},\ \mathbb{K}=\mathbb{C}\right) $ and $x,a\in
H, $ $r>0$ are such that%
\begin{equation*}
x\in \overline{B}\left( a,r\right) :=\left\{ z\in H|\left\Vert
z-a\right\Vert \leq r\right\} .
\end{equation*}

\begin{enumerate}
\item[(i)] If $\left\Vert a\right\Vert >r,$ then we have the inequalities%
\begin{equation}
0\leq \left\Vert x\right\Vert ^{2}\left\Vert a\right\Vert ^{2}-\left\vert
\left\langle x,a\right\rangle \right\vert ^{2}\leq \left\Vert x\right\Vert
^{2}\left\Vert a\right\Vert ^{2}-\left[ \func{Re}\left\langle
x,a\right\rangle \right] ^{2}\leq r^{2}\left\Vert x\right\Vert ^{2}.
\label{2.2.3}
\end{equation}%
The constant $C=1$ in front of $r^{2}$ is best possible in the sense that it
cannot be replaced by a smaller one.

\item[(ii)] If $\left\Vert a\right\Vert =r,$ then%
\begin{equation}
\left\Vert x\right\Vert ^{2}\leq 2\func{Re}\left\langle x,a\right\rangle
\leq 2\left\vert \left\langle x,a\right\rangle \right\vert .  \label{2.3.3}
\end{equation}%
The constant $2$ is best possible in both inequalities.

\item[(iii)] If $\left\Vert a\right\Vert <r,$ then%
\begin{equation}
\left\Vert x\right\Vert ^{2}\leq r^{2}-\left\Vert a\right\Vert ^{2}+2\func{Re%
}\left\langle x,a\right\rangle \leq r^{2}-\left\Vert a\right\Vert
^{2}+2\left\vert \left\langle x,a\right\rangle \right\vert .  \label{2.4.3}
\end{equation}%
Here the constant $2$ is also best possible.
\end{enumerate}
\end{theorem}

\begin{proof}
Since $x\in \overline{B}\left( a,r\right) ,$ then obviously $\left\Vert
x-a\right\Vert ^{2}\leq r^{2},$ which is equivalent to 
\begin{equation}
\left\Vert x\right\Vert ^{2}+\left\Vert a\right\Vert ^{2}-r^{2}\leq 2\func{Re%
}\left\langle x,a\right\rangle .  \label{2.5.3}
\end{equation}

\begin{enumerate}
\item[(i)] If $\left\Vert a\right\Vert >r,$ then we may divide (\ref{2.5.3})
by $\sqrt{\left\Vert a\right\Vert ^{2}-r^{2}}>0$ getting 
\begin{equation}
\frac{\left\Vert x\right\Vert ^{2}}{\sqrt{\left\Vert a\right\Vert ^{2}-r^{2}}%
}+\sqrt{\left\Vert a\right\Vert ^{2}-r^{2}}\leq \frac{2\func{Re}\left\langle
x,a\right\rangle }{\sqrt{\left\Vert a\right\Vert ^{2}-r^{2}}}.  \label{2.6.3}
\end{equation}%
Using the elementary inequality%
\begin{equation*}
\alpha p+\frac{1}{\alpha }q\geq 2\sqrt{pq},\ \ \ \alpha >0,\ \ p,q\geq 0,
\end{equation*}%
we may state that%
\begin{equation}
2\left\Vert x\right\Vert \leq \frac{\left\Vert x\right\Vert ^{2}}{\sqrt{%
\left\Vert a\right\Vert ^{2}-r^{2}}}+\sqrt{\left\Vert a\right\Vert ^{2}-r^{2}%
}.  \label{2.7.3}
\end{equation}%
Making use of (\ref{2.6.3}) and (\ref{2.7.3}), we deduce%
\begin{equation}
\left\Vert x\right\Vert \sqrt{\left\Vert a\right\Vert ^{2}-r^{2}}\leq \func{%
Re}\left\langle x,a\right\rangle ,  \label{2.8.3}
\end{equation}%
which is an interesting inequality in itself as well.

Taking the square in (\ref{2.8.3}) and re-arranging the terms, we deduce the
third inequality in (\ref{2.2.3}). The others are obvious.

To prove the sharpness of the constant, assume, under the hypothesis of the
theorem, that, there exists a constant $c>0$ such that%
\begin{equation}
\left\Vert x\right\Vert ^{2}\left\Vert a\right\Vert ^{2}-\left[ \func{Re}%
\left\langle x,a\right\rangle \right] ^{2}\leq cr^{2}\left\Vert x\right\Vert
^{2},  \label{2.9.3}
\end{equation}%
provided $x\in \overline{B}\left( a,r\right) $ and $\left\Vert a\right\Vert
>r.$

Let $r=\sqrt{\varepsilon }>0,$ $\varepsilon \in \left( 0,1\right) ,$ $a,e\in
H$ with $\left\Vert a\right\Vert =\left\Vert e\right\Vert =1$ and $a\perp e.$
Put $x=a+\sqrt{\varepsilon }e.$ Then obviously $x\in \overline{B}\left(
a,r\right) ,$ $\left\Vert a\right\Vert >r$ and $\left\Vert x\right\Vert
^{2}=\left\Vert a\right\Vert ^{2}+\varepsilon \left\Vert e\right\Vert
^{2}=1+\varepsilon $, $\func{Re}\left\langle x,a\right\rangle =\left\Vert
a\right\Vert ^{2}=1,$ and thus $\left\Vert x\right\Vert ^{2}\left\Vert
a\right\Vert ^{2}-\left[ \func{Re}\left\langle x,a\right\rangle \right]
^{2}=\varepsilon .$ Using (\ref{2.9.3}), we may write that%
\begin{equation*}
\varepsilon \leq c\varepsilon \left( 1+\varepsilon \right) ,\ \ \varepsilon
>0
\end{equation*}%
giving 
\begin{equation}
c+c\varepsilon \geq 1\text{ \ for any }\varepsilon >0.  \label{2.10.3}
\end{equation}%
Letting $\varepsilon \rightarrow 0+,$ we get from (\ref{2.10.3}) that $c\geq
1,$ and the sharpness of the constant is proved.

\item[(ii)] The inequality (\ref{2.3.3}) is obvious by (\ref{2.5.3}) since $%
\left\Vert a\right\Vert =r.$ The best constant follows in a similar way to
the above.

\item[(iii)] The inequality (\ref{2.4.3}) is obvious. The best constant may
be proved in a similar way to the above. We omit the details.
\end{enumerate}
\end{proof}

The following reverse of Schwarz's inequality holds \cite{3NSSD}.

\begin{theorem}
\label{t2.2.3}Let $\left( H;\left\langle \cdot ,\cdot \right\rangle \right) $
be an inner product space over $\mathbb{K}$ and $x,y\in H,$ $\gamma ,\Gamma
\in \mathbb{K}$ such that either%
\begin{equation}
\func{Re}\left\langle \Gamma y-x,x-\gamma y\right\rangle \geq 0,
\label{2.11.3}
\end{equation}%
or equivalently,%
\begin{equation}
\left\Vert x-\frac{\Gamma +\gamma }{2}y\right\Vert \leq \frac{1}{2}%
\left\vert \Gamma -\gamma \right\vert \left\Vert y\right\Vert ,
\label{2.12.3}
\end{equation}%
holds.

\begin{enumerate}
\item[(i)] If $\func{Re}\left( \Gamma \overline{\gamma }\right) >0,$ then we
have the inequalities%
\begin{align}
\left\Vert x\right\Vert ^{2}\left\Vert y\right\Vert ^{2}& \leq \frac{1}{4}%
\cdot \frac{\left\{ \func{Re}\left[ \left( \overline{\Gamma }+\overline{%
\gamma }\right) \left\langle x,y\right\rangle \right] \right\} ^{2}}{\func{Re%
}\left( \Gamma \overline{\gamma }\right) }  \label{2.13.3} \\
& \leq \frac{1}{4}\cdot \frac{\left\vert \Gamma +\gamma \right\vert ^{2}}{%
\func{Re}\left( \Gamma \overline{\gamma }\right) }\left\vert \left\langle
x,y\right\rangle \right\vert ^{2}.  \notag
\end{align}%
The constant $\frac{1}{4}$ is best possible in both inequalities.

\item[(ii)] If $\func{Re}\left( \Gamma \overline{\gamma }\right) =0,$ then 
\begin{equation*}
\left\Vert x\right\Vert ^{2}\leq \func{Re}\left[ \left( \overline{\Gamma }+%
\overline{\gamma }\right) \left\langle x,y\right\rangle \right] \leq
\left\vert \Gamma +\gamma \right\vert \left\vert \left\langle
x,y\right\rangle \right\vert .
\end{equation*}

\item[(iii)] If $\func{Re}\left( \Gamma \overline{\gamma }\right) <0,$ then 
\begin{align*}
\left\Vert x\right\Vert ^{2}& \leq -\func{Re}\left( \Gamma \overline{\gamma }%
\right) \left\Vert y\right\Vert ^{2}+\func{Re}\left[ \left( \overline{\Gamma 
}+\overline{\gamma }\right) \left\langle x,y\right\rangle \right] \\
& \leq -\func{Re}\left( \Gamma \overline{\gamma }\right) \left\Vert
y\right\Vert ^{2}+\left\vert \Gamma +\gamma \right\vert \left\vert
\left\langle x,y\right\rangle \right\vert .
\end{align*}
\end{enumerate}
\end{theorem}

\begin{proof}
The proof of the equivalence between the inequalities (\ref{2.11.3}) and (%
\ref{2.12.3}) follows by the fact that in an inner product space $\func{Re}%
\left\langle Z-x,x-z\right\rangle \geq 0$ for $x,z,Z\in H$ is equivalent
with $\left\Vert x-\frac{z+Z}{2}\right\Vert \leq \frac{1}{2}\left\Vert
Z-z\right\Vert $ (see for example \cite{SSD3.3}).

Consider, for $y\neq 0,$ $a=\frac{\gamma +\Gamma }{2}y$ and $r=\frac{1}{2}%
\left\vert \Gamma -\gamma \right\vert \left\Vert y\right\Vert .$ Then%
\begin{equation*}
\left\Vert a\right\Vert ^{2}-r^{2}=\frac{\left\vert \Gamma +\gamma
\right\vert ^{2}-\left\vert \Gamma -\gamma \right\vert ^{2}}{4}\left\Vert
y\right\Vert ^{2}=\func{Re}\left( \Gamma \overline{\gamma }\right)
\left\Vert y\right\Vert ^{2}.
\end{equation*}

\begin{enumerate}
\item[(i)] If $\func{Re}\left( \Gamma \overline{\gamma }\right) >0,$ then
the hypothesis of (i) in Theorem \ref{t2.1.3} is satisfied, and by the
second inequality in (\ref{2.2.3}) we have%
\begin{equation*}
\left\Vert x\right\Vert ^{2}\frac{\left\vert \Gamma +\gamma \right\vert ^{2}%
}{4}\left\Vert y\right\Vert ^{2}-\frac{1}{4}\left\{ \func{Re}\left[ \left( 
\overline{\Gamma }+\overline{\gamma }\right) \left\langle x,y\right\rangle %
\right] \right\} ^{2}\leq \frac{1}{4}\left\vert \Gamma -\gamma \right\vert
^{2}\left\Vert x\right\Vert ^{2}\left\Vert y\right\Vert ^{2}
\end{equation*}%
from where we derive%
\begin{equation*}
\frac{\left\vert \Gamma +\gamma \right\vert ^{2}-\left\vert \Gamma -\gamma
\right\vert ^{2}}{4}\left\Vert x\right\Vert ^{2}\left\Vert y\right\Vert
^{2}\leq \frac{1}{4}\left\{ \func{Re}\left[ \left( \overline{\Gamma }+%
\overline{\gamma }\right) \left\langle x,y\right\rangle \right] \right\}
^{2},
\end{equation*}%
giving the first inequality in (\ref{2.13.3}).

The second inequality is obvious.

To prove the sharpness of the constant $\frac{1}{4},$ assume that the first
inequality in (\ref{2.13.3}) holds with a constant $c>0,$ i.e., 
\begin{equation}
\left\Vert x\right\Vert ^{2}\left\Vert y\right\Vert ^{2}\leq c\cdot \frac{%
\left\{ \func{Re}\left[ \left( \overline{\Gamma }+\overline{\gamma }\right)
\left\langle x,y\right\rangle \right] \right\} ^{2}}{\func{Re}\left( \Gamma 
\overline{\gamma }\right) },  \label{2.16.3}
\end{equation}%
provided $\func{Re}\left( \Gamma \overline{\gamma }\right) >0$ and either (%
\ref{2.11.3}) or (\ref{2.12.3}) holds.

Assume that $\Gamma ,\gamma >0,$ and let $x=\gamma y.$ Then (\ref{2.11.3})
holds and by (\ref{2.16.3}) we deduce%
\begin{equation*}
\gamma ^{2}\left\Vert y\right\Vert ^{4}\leq c\cdot \frac{\left( \Gamma
+\gamma \right) ^{2}\gamma ^{2}\left\Vert y\right\Vert ^{4}}{\Gamma \gamma }
\end{equation*}%
giving%
\begin{equation}
\Gamma \gamma \leq c\left( \Gamma +\gamma \right) ^{2}\text{ \ for any \ }%
\Gamma ,\gamma >0.  \label{2.17.3}
\end{equation}%
Let $\varepsilon \in \left( 0,1\right) $ and choose in (\ref{2.17.3}), $%
\Gamma =1+\varepsilon ,$ $\gamma =1-\varepsilon >0$ to get $1-\varepsilon
^{2}\leq 4c$ for any $\varepsilon \in \left( 0,1\right) .$ Letting $%
\varepsilon \rightarrow 0+,$ we deduce $c\geq \frac{1}{4},$ and the
sharpness of the constant is proved.

(ii) and (iii) are obvious and we omit the details.
\end{enumerate}
\end{proof}

\begin{remark}
We observe that the second bound in (\ref{2.13.3}) for $\left\Vert
x\right\Vert ^{2}\left\Vert y\right\Vert ^{2}$ is better than the second
bound provided by (\ref{2.2.2}).
\end{remark}

The following corollary provides a reverse inequality for the additive
version of Schwarz's inequality \cite{3NSSD}.

\begin{corollary}
\label{c2.3.3}With the assumptions of Theorem \ref{t2.2.3} and if $\func{Re}%
\left( \Gamma \overline{\gamma }\right) >0,$ then we have the inequality:%
\begin{equation}
0\leq \left\Vert x\right\Vert ^{2}\left\Vert y\right\Vert ^{2}-\left\vert
\left\langle x,y\right\rangle \right\vert ^{2}\leq \frac{1}{4}\cdot \frac{%
\left\vert \Gamma -\gamma \right\vert ^{2}}{\func{Re}\left( \Gamma \overline{%
\gamma }\right) }\left\vert \left\langle x,y\right\rangle \right\vert ^{2}.
\label{2.18.3}
\end{equation}%
The constant $\frac{1}{4}$ is best possible in (\ref{2.18.3}).
\end{corollary}

The proof is obvious from (\ref{2.13.3}) on subtracting in both sides the
same quantity $\left\vert \left\langle x,y\right\rangle \right\vert ^{2}.$
The sharpness of the constant may be proven in a similar manner to the one
incorporated in the proof of (i), Theorem \ref{t2.2.3}. We omit the details.

For other recent results in connection to Schwarz's inequality, see \cite%
{ADR.3}, \cite{DM.3} and \cite{GH.3}.

\subsection{Reverses of the Triangle Inequality}

The following reverse of the triangle inequality holds \cite{3NSSD}.

\begin{proposition}
\label{p2.4.3}Let $\left( H;\left\langle \cdot ,\cdot \right\rangle \right) $
be an inner product space over the real or complex number field $\mathbb{K}$ 
$\left( \mathbb{K}=\mathbb{R},\mathbb{C}\right) $ and $x,a\in H,$ $r>0$ are
such that 
\begin{equation*}
\left\Vert x-a\right\Vert \leq r<\left\Vert a\right\Vert .
\end{equation*}%
Then we have the inequality%
\begin{align}
0& \leq \left\Vert x\right\Vert +\left\Vert a\right\Vert -\left\Vert
x+a\right\Vert  \label{2.20.3} \\
& \leq \sqrt{2}r\cdot \sqrt{\frac{\func{Re}\left\langle x,a\right\rangle }{%
\sqrt{\left\Vert a\right\Vert ^{2}-r^{2}}\left( \sqrt{\left\Vert
a\right\Vert ^{2}-r^{2}}+\left\Vert a\right\Vert \right) }}.  \notag
\end{align}
\end{proposition}

\begin{proof}
Using the inequality (\ref{2.8.3}), we may write that%
\begin{equation*}
\left\Vert x\right\Vert \left\Vert a\right\Vert \leq \frac{\left\Vert
a\right\Vert \func{Re}\left\langle x,a\right\rangle }{\sqrt{\left\Vert
a\right\Vert ^{2}-r^{2}}},
\end{equation*}%
which gives%
\begin{align}
0& \leq \left\Vert x\right\Vert \left\Vert a\right\Vert -\func{Re}%
\left\langle x,a\right\rangle  \label{2.21.3} \\
& \leq \frac{\left\Vert a\right\Vert -\sqrt{\left\Vert a\right\Vert
^{2}-r^{2}}}{\sqrt{\left\Vert a\right\Vert ^{2}-r^{2}}}\func{Re}\left\langle
x,a\right\rangle  \notag \\
& =\frac{r^{2}\func{Re}\left\langle x,a\right\rangle }{\sqrt{\left\Vert
a\right\Vert ^{2}-r^{2}}\left( \sqrt{\left\Vert a\right\Vert ^{2}-r^{2}}%
+\left\Vert a\right\Vert \right) }.  \notag
\end{align}%
Since%
\begin{equation*}
\left( \left\Vert x\right\Vert +\left\Vert a\right\Vert \right)
^{2}-\left\Vert x+a\right\Vert ^{2}=2\left( \left\Vert x\right\Vert
\left\Vert a\right\Vert -\func{Re}\left\langle x,a\right\rangle \right) ,
\end{equation*}%
then by (\ref{2.21.3}), we have%
\begin{align*}
\left\Vert x\right\Vert +\left\Vert a\right\Vert & \leq \sqrt{\left\Vert
x+a\right\Vert ^{2}+\frac{2r^{2}\func{Re}\left\langle x,a\right\rangle }{%
\sqrt{\left\Vert a\right\Vert ^{2}-r^{2}}\left( \sqrt{\left\Vert
a\right\Vert ^{2}-r^{2}}+\left\Vert a\right\Vert \right) }} \\
& \leq \left\Vert x+a\right\Vert +\sqrt{2}r\cdot \sqrt{\frac{\func{Re}%
\left\langle x,a\right\rangle }{\sqrt{\left\Vert a\right\Vert ^{2}-r^{2}}%
\left( \sqrt{\left\Vert a\right\Vert ^{2}-r^{2}}+\left\Vert a\right\Vert
\right) }},
\end{align*}%
giving the desired inequality (\ref{2.20.3}).
\end{proof}

The following proposition providing a simpler reverse for the triangle
inequality also holds \cite{3NSSD}.

\begin{proposition}
\label{p2.5.3}Let $\left( H;\left\langle \cdot ,\cdot \right\rangle \right) $
be an inner product space over $\mathbb{K}$ and $x,y\in H,$ $M>m>0$ such
that either%
\begin{equation*}
\func{Re}\left\langle My-x,x-my\right\rangle \geq 0,
\end{equation*}%
or equivalently,%
\begin{equation*}
\left\Vert x-\frac{M+m}{2}\cdot y\right\Vert \leq \frac{1}{2}\left(
M-m\right) \left\Vert y\right\Vert ,
\end{equation*}%
holds. Then we have the inequality%
\begin{equation}
0\leq \left\Vert x\right\Vert +\left\Vert y\right\Vert -\left\Vert
x+y\right\Vert \leq \frac{\sqrt{M}-\sqrt{m}}{\sqrt[4]{mM}}\sqrt{\func{Re}%
\left\langle x,y\right\rangle }.  \label{2.24.3}
\end{equation}
\end{proposition}

\begin{proof}
Choosing in (\ref{2.8.3}), $a=\frac{M+m}{2}y,$ $r=\frac{1}{2}\left(
M-m\right) \left\Vert y\right\Vert $ we get%
\begin{equation*}
\left\Vert x\right\Vert \left\Vert y\right\Vert \sqrt{Mm}\leq \frac{M+m}{2}%
\func{Re}\left\langle x,y\right\rangle ,
\end{equation*}%
giving 
\begin{equation*}
0\leq \left\Vert x\right\Vert \left\Vert y\right\Vert -\func{Re}\left\langle
x,y\right\rangle \leq \frac{\left( \sqrt{M}-\sqrt{m}\right) ^{2}}{2\sqrt{mM}}%
\func{Re}\left\langle x,y\right\rangle .
\end{equation*}%
Following the same argument as in the proof of Proposition \ref{p2.4.3}, we
deduce the desired inequality (\ref{2.24.3}).
\end{proof}

For some results related to triangle inequality in inner product spaces, see 
\cite{JBDFTM.3}, \cite{SMK.3}, \cite{PMM.3} and \cite{DKR.3}.

\subsection{Integral Inequalities}

Let $\left( \Omega ,\Sigma ,\mu \right) $ be a measurable space consisting
of a set $\Omega ,$ $\Sigma $ a $\sigma -$algebra of parts and $\mu $ a
countably additive and positive measure on $\Sigma $ with values in $\mathbb{%
R\cup }\left\{ \infty \right\} .$ Let $\rho \geq 0$ be a $g-$measurable
function on $\Omega $ with $\int_{\Omega }\rho \left( s\right) d\mu \left(
s\right) =1.$ Denote by $L_{\rho }^{2}\left( \Omega ,\mathbb{K}\right) $ the
Hilbert space of all real or complex valued functions defined on $\Omega $
and $2-\rho -$integrable on $\Omega ,$ i.e.,%
\begin{equation*}
\int_{\Omega }\rho \left( s\right) \left\vert f\left( s\right) \right\vert
^{2}d\mu \left( s\right) <\infty .
\end{equation*}%
It is obvious that the following inner product%
\begin{equation*}
\left\langle f,g\right\rangle _{\rho }:=\int_{\Omega }\rho \left( s\right)
f\left( s\right) \overline{g\left( s\right) }d\mu \left( s\right) ,
\end{equation*}%
generates the norm $\left\Vert f\right\Vert _{\rho }:=\left( \int_{\Omega
}\rho \left( s\right) \left\vert f\left( s\right) \right\vert ^{2}d\mu
\left( s\right) \right) ^{\frac{1}{2}}$ of $L_{\rho }^{2}\left( \Omega ,%
\mathbb{K}\right) ,$ and all the above results may be stated for integrals.

It is important to observe that, if 
\begin{equation*}
\func{Re}\left[ f\left( s\right) \overline{g\left( s\right) }\right] \geq 0,%
\text{ \ for }\mu -\text{a.e. }s\in \Omega ,
\end{equation*}%
then, obviously,%
\begin{align}
\func{Re}\left\langle f,g\right\rangle _{\rho }& =\func{Re}\left[
\int_{\Omega }\rho \left( s\right) f\left( s\right) \overline{g\left(
s\right) }d\mu \left( s\right) \right]  \label{7.4.3} \\
& =\int_{\Omega }\rho \left( s\right) \func{Re}\left[ f\left( s\right) 
\overline{g\left( s\right) }\right] d\mu \left( s\right) \geq 0.  \notag
\end{align}%
The reverse is evidently not true in general.

Moreover, if the space is real, i.e., $\mathbb{K=R}$, then a sufficient
condition for (\ref{7.4.3}) to hold is:%
\begin{equation*}
f\left( s\right) \geq 0,\ \ g\left( s\right) \geq 0,\text{ \ for }\mu -\text{%
a.e. }s\in \Omega .
\end{equation*}

We now provide, by the use of certain results obtained above, some integral
inequalities that may be used in practical applications.

\begin{proposition}
\label{p7.1.3}Let $f,g\in L_{\rho }^{2}\left( \Omega ,\mathbb{K}\right) $
and $r>0$ with the properties that%
\begin{equation}
\left\vert f\left( s\right) -g\left( s\right) \right\vert \leq r\leq
\left\vert g\left( s\right) \right\vert ,\ \text{\ for }\mu -\text{a.e. }%
s\in \Omega .  \label{7.6.3}
\end{equation}%
Then we have the inequalities%
\begin{align}
0& \leq \int_{\Omega }\rho \left( s\right) \left\vert f\left( s\right)
\right\vert ^{2}d\mu \left( s\right) \int_{\Omega }\rho \left( s\right)
\left\vert g\left( s\right) \right\vert ^{2}d\mu \left( s\right)
\label{7.7.3} \\
& \ \ \ \ \ \ \ \ \ \ \ \ \ \ \ \ \ \ \ \ \ \ \ \ \ \ \ -\left\vert
\int_{\Omega }\rho \left( s\right) f\left( s\right) \overline{g\left(
s\right) }d\mu \left( s\right) \right\vert ^{2}  \notag \\
& \leq \int_{\Omega }\rho \left( s\right) \left\vert f\left( s\right)
\right\vert ^{2}d\mu \left( s\right) \int_{\Omega }\rho \left( s\right)
\left\vert g\left( s\right) \right\vert ^{2}d\mu \left( s\right)  \notag \\
& \ \ \ \ \ \ \ \ \ \ \ \ \ \ \ \ \ \ \ \ \ \ \ \ \ \ \ -\left[ \int_{\Omega
}\rho \left( s\right) \func{Re}\left( f\left( s\right) \overline{g\left(
s\right) }\right) d\mu \left( s\right) \right] ^{2}  \notag \\
& \leq r^{2}\int_{\Omega }\rho \left( s\right) \left\vert g\left( s\right)
\right\vert ^{2}d\mu \left( s\right) .  \notag
\end{align}%
The constant $c=1$ in front of $r^{2}$ is best possible.
\end{proposition}

The proof follows by Theorem \ref{t2.1.3} and we omit the details \cite%
{3NSSD}.

\begin{proposition}
\label{p7.2.3}Let $f,g\in L_{\rho }^{2}\left( \Omega ,\mathbb{K}\right) $
and $\gamma ,\Gamma \in \mathbb{K}$ such that $\func{Re}\left( \Gamma 
\overline{\gamma }\right) >0$ and%
\begin{equation*}
\func{Re}\left[ \left( \Gamma g\left( s\right) -f\left( s\right) \right)
\left( \overline{f\left( s\right) }-\overline{\gamma }\overline{g\left(
s\right) }\right) \right] \geq 0,\text{ \ for }\mu -\text{a.e. }s\in \Omega .
\end{equation*}%
Then we have the inequalities%
\begin{align}
& \int_{\Omega }\rho \left( s\right) \left\vert f\left( s\right) \right\vert
^{2}d\mu \left( s\right) \int_{\Omega }\rho \left( s\right) \left\vert
g\left( s\right) \right\vert ^{2}d\mu \left( s\right)  \label{7.9.3} \\
& \leq \frac{1}{4}\cdot \frac{\left\{ \func{Re}\left[ \left( \overline{%
\Gamma }+\overline{\gamma }\right) \int_{\Omega }\rho \left( s\right)
f\left( s\right) \overline{g\left( s\right) }d\mu \left( s\right) \right]
\right\} ^{2}}{\func{Re}\left( \Gamma \overline{\gamma }\right) }  \notag \\
& \leq \frac{1}{4}\cdot \frac{\left\vert \Gamma +\gamma \right\vert ^{2}}{%
\func{Re}\left( \Gamma \overline{\gamma }\right) }\left\vert \int_{\Omega
}\rho \left( s\right) f\left( s\right) \overline{g\left( s\right) }d\mu
\left( s\right) \right\vert ^{2}.  \notag
\end{align}%
The constant $\frac{1}{4}$ is best possible in both inequalities.
\end{proposition}

The proof follows by Theorem \ref{t2.2.3} and we omit the details.

\begin{corollary}
\label{c7.3.3}With the assumptions of Proposition \ref{p7.2.3}, we have the
inequality%
\begin{align}
0& \leq \int_{\Omega }\rho \left( s\right) \left\vert f\left( s\right)
\right\vert ^{2}d\mu \left( s\right) \int_{\Omega }\rho \left( s\right)
\left\vert g\left( s\right) \right\vert ^{2}d\mu \left( s\right)
\label{7.10.3} \\
& \ \ \ \ \ \ \ \ \ \ \ \ \ \ \ \ \ \ \ \ \ \ \ \ \ -\left\vert \int_{\Omega
}\rho \left( s\right) f\left( s\right) \overline{g\left( s\right) }d\mu
\left( s\right) \right\vert ^{2}  \notag \\
& \leq \frac{1}{4}\cdot \frac{\left\vert \Gamma -\gamma \right\vert ^{2}}{%
\func{Re}\left( \Gamma \overline{\gamma }\right) }\left\vert \int_{\Omega
}\rho \left( s\right) f\left( s\right) \overline{g\left( s\right) }d\mu
\left( s\right) \right\vert ^{2}.  \notag
\end{align}%
The constant $\frac{1}{4}$ is best possible.
\end{corollary}

\begin{remark}
If the space is real and we assume, for $M>m>0,$ that%
\begin{equation*}
mg\left( s\right) \leq f\left( s\right) \leq Mg\left( s\right) ,\text{ \ for 
}\mu -\text{a.e. }s\in \Omega ,
\end{equation*}%
then, by (\ref{7.9.3}) and (\ref{7.10.3}), we deduce the inequalities%
\begin{multline}
\int_{\Omega }\rho \left( s\right) \left[ f\left( s\right) \right] ^{2}d\mu
\left( s\right) \int_{\Omega }\rho \left( s\right) \left[ g\left( s\right) %
\right] ^{2}d\mu \left( s\right)  \label{7.12.3} \\
\leq \frac{1}{4}\cdot \frac{\left( M+m\right) ^{2}}{mM}\left[ \int_{\Omega
}\rho \left( s\right) f\left( s\right) g\left( s\right) d\mu \left( s\right) %
\right] ^{2}
\end{multline}%
and 
\begin{align}
0& \leq \int_{\Omega }\rho \left( s\right) \left[ f\left( s\right) \right]
^{2}d\mu \left( s\right) \int_{\Omega }\rho \left( s\right) \left[ g\left(
s\right) \right] ^{2}d\mu \left( s\right)  \label{7.13.3} \\
& \ \ \ \ \ \ \ \ \ \ \ \ \ \ \ \ \ \ \ -\left[ \int_{\Omega }\rho \left(
s\right) f\left( s\right) g\left( s\right) d\mu \left( s\right) \right] ^{2}
\notag \\
& \leq \frac{1}{4}\cdot \frac{\left( M-m\right) ^{2}}{mM}\left[ \int_{\Omega
}\rho \left( s\right) f\left( s\right) g\left( s\right) d\mu \left( s\right) %
\right] ^{2}.  \notag
\end{align}%
The inequality (\ref{7.12.3}) is known in the literature as Cassel's
inequality.
\end{remark}

\newpage

\section{More Reverses of Schwarz's Inequality}

\subsection{General Results}

The following result holds \cite{4NSSD}.

\begin{theorem}
\label{t2.1.4}Let $\left( H;\left\langle \cdot ,\cdot \right\rangle \right) $
be an inner product space over the real or complex number field $\mathbb{K}$%
, $x,a\in H$ and $r>0.$ If%
\begin{equation}
x\in \bar{B}\left( a,r\right) :=\left\{ z\in H|\left\Vert z-a\right\Vert
\leq r\right\} ,  \label{2.1.4}
\end{equation}%
then we have the inequalities:%
\begin{align}
0& \leq \left\Vert x\right\Vert \left\Vert a\right\Vert -\left\vert
\left\langle x,a\right\rangle \right\vert \leq \left\Vert x\right\Vert
\left\Vert a\right\Vert -\left\vert \func{Re}\left\langle x,a\right\rangle
\right\vert  \label{2.2.4} \\
& \leq \left\Vert x\right\Vert \left\Vert a\right\Vert -\func{Re}%
\left\langle x,a\right\rangle \leq \frac{1}{2}r^{2}.  \notag
\end{align}%
The constant $\frac{1}{2}$ is best possible in (\ref{2.2.4}) in the sense
that it cannot be replaced by a smaller constant.
\end{theorem}

\begin{proof}
The condition (\ref{2.1.4}) is clearly equivalent to%
\begin{equation}
\left\Vert x\right\Vert ^{2}+\left\Vert a\right\Vert ^{2}\leq 2\func{Re}%
\left\langle x,a\right\rangle +r^{2}.  \label{2.3.4}
\end{equation}%
Using the elementary inequality%
\begin{equation*}
2\left\Vert x\right\Vert \left\Vert a\right\Vert \leq \left\Vert
x\right\Vert ^{2}+\left\Vert a\right\Vert ^{2},\ \ \ \ a,x\in H
\end{equation*}%
and (\ref{2.3.4}), we deduce%
\begin{equation*}
2\left\Vert x\right\Vert \left\Vert a\right\Vert \leq 2\func{Re}\left\langle
x,a\right\rangle +r^{2},
\end{equation*}%
giving the last inequality in (\ref{2.2.4}). The other inequalities are
obvious.

To prove the sharpness of the constant $\frac{1}{2},$ assume that%
\begin{equation}
0\leq \left\Vert x\right\Vert \left\Vert a\right\Vert -\func{Re}\left\langle
x,a\right\rangle \leq cr^{2}  \label{2.3a.4}
\end{equation}%
for any $x,a\in H$ and $r>0$ satisfying (\ref{2.1.4}).

Assume that $a,e\in H,$ $\left\Vert a\right\Vert =\left\Vert e\right\Vert =1$
and $e\perp a.$ If $r=\sqrt{\varepsilon },$ $\varepsilon >0$ and if we
define $x=a+\sqrt{\varepsilon }e,$ then $\left\Vert x-a\right\Vert =\sqrt{%
\varepsilon }=r$ showing that the condition (\ref{2.1.4}) is fulfilled.

On the other hand,%
\begin{align*}
\left\Vert x\right\Vert \left\Vert a\right\Vert -\func{Re}\left\langle
x,a\right\rangle & =\sqrt{\left\Vert a+\sqrt{\varepsilon }e\right\Vert ^{2}}-%
\func{Re}\left\langle a+\sqrt{\varepsilon }e,a\right\rangle \\
& =\sqrt{\left\Vert a\right\Vert ^{2}+\varepsilon \left\Vert e\right\Vert
^{2}}-\left\Vert a\right\Vert ^{2} \\
& =\sqrt{1+\varepsilon }-1.
\end{align*}%
Utilising (\ref{2.3a.4}), we conclude that%
\begin{equation}
\sqrt{1+\varepsilon }-1\leq c\varepsilon \text{ \ for any }\varepsilon >0.
\label{2.4b.4}
\end{equation}%
Multiplying (\ref{2.4b.4}) by $\sqrt{1+\varepsilon }+1>0$ and then dividing
by $\varepsilon >0,$ we get%
\begin{equation}
\left( \sqrt{1+\varepsilon }+1\right) c\geq 1\text{ \ for any \ }\varepsilon
>0.  \label{2.5.4}
\end{equation}%
Letting $\varepsilon \rightarrow 0+$ in (\ref{2.5.4}), we deduce $c\geq 
\frac{1}{2},$ and the theorem is proved.
\end{proof}

The following result also holds \cite{4NSSD}.

\begin{theorem}
\label{t2.2.4}Let $\left( H;\left\langle \cdot ,\cdot \right\rangle \right) $
be an inner product space over $\mathbb{K}$ and $x,y\in H,$ $\gamma ,\Gamma
\in \mathbb{K}$ \ $\left( \Gamma \neq -\gamma \right) $ so that either%
\begin{equation}
\func{Re}\left\langle \Gamma y-x,x-\gamma y\right\rangle \geq 0,
\label{2.6.4}
\end{equation}%
or equivalently,%
\begin{equation}
\left\Vert x-\frac{\gamma +\Gamma }{2}y\right\Vert \leq \frac{1}{2}%
\left\vert \Gamma -\gamma \right\vert \left\Vert y\right\Vert ,
\label{2.7.4}
\end{equation}%
holds. Then we have the inequalities%
\begin{align}
0& \leq \left\Vert x\right\Vert \left\Vert y\right\Vert -\left\vert
\left\langle x,y\right\rangle \right\vert  \label{2.8.4} \\
& \leq \left\Vert x\right\Vert \left\Vert y\right\Vert -\left\vert \func{Re}%
\left[ \frac{\bar{\Gamma}+\bar{\gamma}}{\left\vert \Gamma +\gamma
\right\vert }\left\langle x,y\right\rangle \right] \right\vert  \notag \\
& \leq \left\Vert x\right\Vert \left\Vert y\right\Vert -\func{Re}\left[ 
\frac{\bar{\Gamma}+\bar{\gamma}}{\left\vert \Gamma +\gamma \right\vert }%
\left\langle x,y\right\rangle \right]  \notag \\
& \leq \frac{1}{4}\cdot \frac{\left\vert \Gamma -\gamma \right\vert ^{2}}{%
\left\vert \Gamma +\gamma \right\vert }\left\Vert y\right\Vert ^{2}.  \notag
\end{align}%
The constant $\frac{1}{4}$ in the last inequality is best possible.
\end{theorem}

\begin{proof}
Consider for $a,y\neq 0,$ $a=\frac{\Gamma +\gamma }{2}\cdot y$ and $r=\frac{1%
}{2}\left\vert \Gamma -\gamma \right\vert \left\Vert y\right\Vert .$ Thus
from (\ref{2.2.4}), we get%
\begin{align*}
0& \leq \left\Vert x\right\Vert \left\vert \frac{\Gamma +\gamma }{2}%
\right\vert \left\Vert y\right\Vert -\left\vert \frac{\Gamma +\gamma }{2}%
\right\vert \left\vert \left\langle x,y\right\rangle \right\vert \\
& \leq \left\Vert x\right\Vert \left\vert \frac{\Gamma +\gamma }{2}%
\right\vert \left\Vert y\right\Vert -\left\vert \func{Re}\left[ \frac{\bar{%
\Gamma}+\bar{\gamma}}{\left\vert \Gamma +\gamma \right\vert }\left\langle
x,y\right\rangle \right] \right\vert \\
& \leq \left\Vert x\right\Vert \left\vert \frac{\Gamma +\gamma }{2}%
\right\vert \left\Vert y\right\Vert -\func{Re}\left[ \frac{\bar{\Gamma}+\bar{%
\gamma}}{\left\vert \Gamma +\gamma \right\vert }\left\langle
x,y\right\rangle \right] \\
& \leq \frac{1}{8}\cdot \left\vert \Gamma -\gamma \right\vert ^{2}\left\Vert
y\right\Vert ^{2}.
\end{align*}%
Dividing by $\frac{1}{2}\left\vert \Gamma +\gamma \right\vert \geq 0,$ we
deduce the desired inequality (\ref{2.8.4}).

To prove the sharpness of the constant $\frac{1}{4},$ assume that there
exists a $c>0$ such that:%
\begin{equation}
\left\Vert x\right\Vert \left\Vert y\right\Vert -\func{Re}\left[ \frac{\bar{%
\Gamma}+\bar{\gamma}}{\left\vert \Gamma +\gamma \right\vert }\left\langle
x,y\right\rangle \right] \leq c\cdot \frac{\left\vert \Gamma -\gamma
\right\vert ^{2}}{\left\vert \Gamma +\gamma \right\vert }\left\Vert
y\right\Vert ^{2},  \label{2.9.4}
\end{equation}%
provided either (\ref{2.6.4}) or (\ref{2.7.4}) holds.

Consider the real inner product space $\left( \mathbb{R}^{2},\left\langle
\cdot ,\cdot \right\rangle \right) $ with $\left\langle \mathbf{\bar{x}},%
\mathbf{\bar{y}}\right\rangle =x_{1}y_{1}+x_{2}y_{2},$ $\mathbf{\bar{x}}%
=\left( x_{1},x_{2}\right) ,$ $\mathbf{\bar{y}}=\left( y_{1},y_{2}\right)
\in \mathbb{R}^{2}.$ Let $\mathbf{\bar{y}}=\left( 1,1\right) $ and $\Gamma
,\gamma >0$ with $\Gamma >\gamma .$ Then, by (\ref{2.9.4}), we deduce%
\begin{equation}
\sqrt{2}\sqrt{x_{1}^{2}+x_{2}^{2}}-\left( x_{1}+x_{2}\right) \leq 2c\cdot 
\frac{\left( \Gamma -\gamma \right) ^{2}}{\Gamma +\gamma }.  \label{2.10.4}
\end{equation}%
If $x_{1}=\Gamma ,$ $x_{2}=\gamma ,$ then 
\begin{equation*}
\left\langle \Gamma \mathbf{\bar{y}}-\mathbf{\bar{x}},\mathbf{\bar{x}}%
-\gamma \mathbf{\bar{y}}\right\rangle =\left( \Gamma -x_{1}\right) \left(
x_{1}-\gamma \right) +\left( \Gamma -x_{2}\right) \left( x_{2}-\gamma
\right) =0,
\end{equation*}%
showing that the condition (\ref{2.6.4}) is valid. Replacing \thinspace $%
x_{1}$ and $x_{2}$ in (\ref{2.10.4}), we deduce%
\begin{equation}
\sqrt{2}\sqrt{\Gamma ^{2}+\gamma ^{2}}-\left( \Gamma +\gamma \right) \leq 2c%
\frac{\left( \Gamma -\gamma \right) ^{2}}{\Gamma +\gamma }.  \label{2.11.4}
\end{equation}%
If in (\ref{2.11.4}) we choose $\Gamma =1+\varepsilon ,$ $\gamma
=1-\varepsilon $ with $\varepsilon \in \left( 0,1\right) ,$ then we have 
\begin{equation*}
2\sqrt{1+\varepsilon ^{2}}-2\leq 2c\frac{4\varepsilon ^{2}}{2},
\end{equation*}%
giving%
\begin{equation}
\sqrt{1+\varepsilon ^{2}}-1\leq 2c\varepsilon ^{2}.  \label{2.12.4}
\end{equation}%
Finally, multiplying (\ref{2.12.4}) with $\sqrt{1+\varepsilon ^{2}}+1>0$ and
thus dividing by $\varepsilon ^{2},$ we deduce%
\begin{equation}
1\leq 2c\left( \sqrt{1+\varepsilon ^{2}}+1\right) \text{ \ for any \ }%
\varepsilon \in \left( 0,1\right) .  \label{2.13.4}
\end{equation}%
Letting $\varepsilon \rightarrow 0+$ in (\ref{2.13.4}) we get $c\geq \frac{1%
}{4},$ and the sharpness of the constant is proved.
\end{proof}

\subsection{Reverses of the Triangle Inequality}

The following reverse of the triangle inequality in inner product spaces
holds \cite{4NSSD}.

\begin{proposition}
\label{p2.3.4}Let $\left( H;\left\langle \cdot ,\cdot \right\rangle \right) $
be an inner product space over the real or complex number field $\mathbb{K}$%
, $x,a\in H$ and $r>0.$ If $\left\Vert x-a\right\Vert \leq r,$ then we have
the inequality%
\begin{equation}
0\leq \left\Vert x\right\Vert +\left\Vert a\right\Vert -\left\Vert
x+a\right\Vert \leq r.  \label{2.14.4}
\end{equation}
\end{proposition}

\begin{proof}
Since%
\begin{equation*}
\left( \left\Vert x\right\Vert +\left\Vert a\right\Vert \right)
^{2}-\left\Vert x+a\right\Vert ^{2}\leq 2\left( \left\Vert x\right\Vert
\left\Vert a\right\Vert -\func{Re}\left\langle x,a\right\rangle \right) ,
\end{equation*}%
then by Theorem \ref{t2.1.4} we deduce%
\begin{equation*}
\left( \left\Vert x\right\Vert +\left\Vert a\right\Vert \right)
^{2}-\left\Vert x+a\right\Vert ^{2}\leq r^{2},
\end{equation*}%
from where we obtain%
\begin{equation*}
\left\Vert x\right\Vert +\left\Vert a\right\Vert \leq \sqrt{r^{2}+\left\Vert
x+a\right\Vert ^{2}}\leq r+\left\Vert x+a\right\Vert ,
\end{equation*}%
giving the desired result (\ref{2.14.4}).
\end{proof}

We may state the following result \cite{4NSSD}.

\begin{proposition}
\label{p2.4.4}Let $\left( H;\left\langle \cdot ,\cdot \right\rangle \right) $
be an inner product space over $\mathbb{K}$ and $x,y\in H,$ $M>m>0$ such
that either%
\begin{equation*}
\func{Re}\left\langle My-x,x-my\right\rangle \geq 0,
\end{equation*}%
or equivalently,%
\begin{equation*}
\left\Vert x-\frac{M+m}{2}y\right\Vert \leq \frac{1}{2}\left( M-m\right)
\left\Vert y\right\Vert ,
\end{equation*}%
holds. Then we have the inequality%
\begin{equation}
0\leq \left\Vert x\right\Vert +\left\Vert y\right\Vert -\left\Vert
x+y\right\Vert \leq \frac{\sqrt{2}}{2}\cdot \frac{\left( M-m\right) }{\sqrt{%
M+m}}\left\Vert y\right\Vert .  \label{2.19.4}
\end{equation}
\end{proposition}

\begin{proof}
By Theorem \ref{t2.2.4} for $\Gamma =M,$ $\gamma =m,$ we have the inequality%
\begin{equation*}
\left\Vert x\right\Vert \left\Vert y\right\Vert -\func{Re}\left\langle
x,y\right\rangle \leq \frac{1}{4}\cdot \frac{\left( M-m\right) ^{2}}{\left(
M+m\right) }\left\Vert y\right\Vert ^{2}.
\end{equation*}%
Then we may state that%
\begin{align*}
\left( \left\Vert x\right\Vert +\left\Vert y\right\Vert \right)
^{2}-\left\Vert x+y\right\Vert ^{2}& =2\left( \left\Vert x\right\Vert
\left\Vert y\right\Vert -\func{Re}\left\langle x,y\right\rangle \right) \\
& \leq \frac{1}{2}\cdot \frac{\left( M-m\right) ^{2}}{M+m}\left\Vert
y\right\Vert ^{2},
\end{align*}%
from where we get%
\begin{align*}
\left\Vert x\right\Vert +\left\Vert y\right\Vert & \leq \sqrt{\frac{1}{2}%
\cdot \frac{\left( M-m\right) ^{2}}{M+m}\left\Vert y\right\Vert
^{2}+\left\Vert x+y\right\Vert ^{2}} \\
& \leq \left\Vert x+y\right\Vert +\frac{\left( M-m\right) }{\sqrt{2\left(
M+m\right) }}\left\Vert y\right\Vert ,
\end{align*}%
giving the desired inequality (\ref{2.19.4}).
\end{proof}

\subsection{Integral Inequalities}

We provide now, by the use of certain results obtained above, some integral
inequalities that may be used in practical applications.

\begin{proposition}
\label{p7.1.4}Let $f,g\in L_{\rho }^{2}\left( \Omega ,\mathbb{K}\right) $
and $r>0$ with the property that%
\begin{equation*}
\left\vert f\left( s\right) -g\left( s\right) \right\vert \leq r\text{ \ \
for, \ }\mu -\text{a.e. \ }s\in \Omega .
\end{equation*}%
Then we have the inequalities%
\begin{align}
0& \leq \left[ \int_{\Omega }\rho \left( s\right) \left\vert f\left(
s\right) \right\vert ^{2}d\mu \left( s\right) \int_{\Omega }\rho \left(
s\right) \left\vert g\left( s\right) \right\vert ^{2}d\mu \left( s\right) %
\right] ^{\frac{1}{2}}  \label{7.7.4} \\
& \qquad \qquad -\left\vert \int_{\Omega }\rho \left( s\right) f\left(
s\right) \overline{g\left( s\right) }d\mu \left( s\right) \right\vert  \notag
\\
& \leq \left[ \int_{\Omega }\rho \left( s\right) \left\vert f\left( s\right)
\right\vert ^{2}d\mu \left( s\right) \int_{\Omega }\rho \left( s\right)
\left\vert g\left( s\right) \right\vert ^{2}d\mu \left( s\right) \right] ^{%
\frac{1}{2}}  \notag \\
& \qquad \qquad -\left\vert \int_{\Omega }\rho \left( s\right) \func{Re}%
\left[ f\left( s\right) \overline{g\left( s\right) }\right] d\mu \left(
s\right) \right\vert  \notag \\
& \leq \left[ \int_{\Omega }\rho \left( s\right) \left\vert f\left( s\right)
\right\vert ^{2}d\mu \left( s\right) \int_{\Omega }\rho \left( s\right)
\left\vert g\left( s\right) \right\vert ^{2}d\mu \left( s\right) \right] ^{%
\frac{1}{2}}  \notag \\
& \qquad \qquad -\int_{\Omega }\rho \left( s\right) \func{Re}\left[ f\left(
s\right) \overline{g\left( s\right) }\right] d\mu \left( s\right)  \notag \\
& \leq \frac{1}{2}r^{2}.  \notag
\end{align}%
The constant $\frac{1}{2}$ is best possible in (\ref{7.7.4}).
\end{proposition}

The proof follows by Theorem \ref{t2.1.4}, and we omit the details.

\begin{proposition}
\label{p7.2.4}Let $f,g\in L_{\rho }^{2}\left( \Omega ,\mathbb{K}\right) $
and $\gamma ,\Gamma \in \mathbb{K}$ so that $\Gamma \neq -\gamma ,$ and%
\begin{equation*}
\func{Re}\left[ \left( \Gamma g\left( s\right) -f\left( s\right) \right)
\left( \overline{f\left( s\right) }-\overline{\gamma }\overline{g\left(
s\right) }\right) \right] \geq 0,\text{ \ \ for \ }\mu -\text{a.e. \ }s\in
\Omega .
\end{equation*}%
Then we have the inequalities%
\begin{align}
0& \leq \left[ \int_{\Omega }\rho \left( s\right) \left\vert f\left(
s\right) \right\vert ^{2}d\mu \left( s\right) \int_{\Omega }\rho \left(
s\right) \left\vert g\left( s\right) \right\vert ^{2}d\mu \left( s\right) %
\right] ^{\frac{1}{2}}  \label{7.9.4} \\
& \qquad \qquad -\left\vert \int_{\Omega }\rho \left( s\right) f\left(
s\right) \overline{g\left( s\right) }d\mu \left( s\right) \right\vert  \notag
\\
& \leq \left[ \int_{\Omega }\rho \left( s\right) \left\vert f\left( s\right)
\right\vert ^{2}d\mu \left( s\right) \int_{\Omega }\rho \left( s\right)
\left\vert g\left( s\right) \right\vert ^{2}d\mu \left( s\right) \right] ^{%
\frac{1}{2}}  \notag \\
& \qquad \qquad -\left\vert \func{Re}\left[ \frac{\bar{\Gamma}+\bar{\gamma}}{%
\left\vert \Gamma +\gamma \right\vert }\int_{\Omega }\rho \left( s\right)
f\left( s\right) \overline{g\left( s\right) }d\mu \left( s\right) \right]
\right\vert  \notag
\end{align}%
\begin{align*}
& \leq \left[ \int_{\Omega }\rho \left( s\right) \left\vert f\left( s\right)
\right\vert ^{2}d\mu \left( s\right) \int_{\Omega }\rho \left( s\right)
\left\vert g\left( s\right) \right\vert ^{2}d\mu \left( s\right) \right] ^{%
\frac{1}{2}} \\
& \ \ \ \ \ \ \ \ \ \ \ \ \ \ \ \ \ \ \ \ \ \ \ \ \ \ \ -\func{Re}\left[ 
\frac{\bar{\Gamma}+\bar{\gamma}}{\left\vert \Gamma +\gamma \right\vert }%
\int_{\Omega }\rho \left( s\right) f\left( s\right) \overline{g\left(
s\right) }d\mu \left( s\right) \right] \\
& \leq \frac{1}{4}\cdot \frac{\left\vert \Gamma -\gamma \right\vert ^{2}}{%
\left\vert \Gamma +\gamma \right\vert }\int_{\Omega }\rho \left( s\right)
\left\vert g\left( s\right) \right\vert ^{2}d\mu \left( s\right) .
\end{align*}%
The constant $\frac{1}{4}$ is best possible.
\end{proposition}

\begin{remark}
If the space is real and we assume, for $M>m>0,$ that%
\begin{equation}
mg\left( s\right) \leq f\left( s\right) \leq Mg\left( s\right) ,\text{ \ \
for \ }\mu -\text{a.e. \ }s\in \Omega ,  \label{7.10.4}
\end{equation}%
then, by (\ref{7.9.4}), we deduce the inequality:%
\begin{eqnarray*}
0 &\leq &\left[ \int_{\Omega }\rho \left( s\right) \left\vert f\left(
s\right) \right\vert ^{2}d\mu \left( s\right) \int_{\Omega }\rho \left(
s\right) \left\vert g\left( s\right) \right\vert ^{2}d\mu \left( s\right) %
\right] ^{\frac{1}{2}} \\
&&\qquad \qquad -\left\vert \int_{\Omega }\rho \left( s\right) f\left(
s\right) \overline{g\left( s\right) }d\mu \left( s\right) \right\vert \\
&\leq &\frac{1}{4}\cdot \frac{\left( M-m\right) ^{2}}{M+m}\int_{\Omega }\rho
\left( s\right) \left\vert g\left( s\right) \right\vert ^{2}d\mu \left(
s\right) .
\end{eqnarray*}%
The constant $\frac{1}{4}$ is best possible.
\end{remark}

The following reverse of the triangle inequality for integrals holds.

\begin{proposition}
\label{p.7.3.4}Assume that the functions $f,g\in L_{\rho }^{2}\left( \Omega ,%
\mathbb{K}\right) $ satisfy (\ref{7.10.4}). Then we have the inequality%
\begin{align*}
0& \leq \left( \int_{\Omega }\rho \left( s\right) \left\vert f\left(
s\right) \right\vert ^{2}d\mu \left( s\right) \right) ^{\frac{1}{2}}+\left(
\int_{\Omega }\rho \left( s\right) \left\vert g\left( s\right) \right\vert
^{2}d\mu \left( s\right) \right) ^{\frac{1}{2}} \\
& \ \ \ \ \ \ \ \ \ \ \ \ \ \ \ \ \ \ \ \ \ \ \ \ \ \ \ -\left( \int_{\Omega
}\rho \left( s\right) \left\vert f\left( s\right) +g\left( s\right)
\right\vert ^{2}d\mu \left( s\right) \right) ^{\frac{1}{2}} \\
& \leq \frac{\sqrt{2}}{2}\cdot \frac{M-m}{\sqrt{M+m}}\left( \int_{\Omega
}\rho \left( s\right) \left\vert g\left( s\right) \right\vert ^{2}d\mu
\left( s\right) \right) ^{\frac{1}{2}}.
\end{align*}
\end{proposition}

The proof follows by Proposition \ref{p2.4.4}.

\chapter{Inequalities of the Gr\"{u}ss Type}

\section{Introduction}

Over the last five years, the development of Gr\"{u}ss type inequalities has
experienced a surge, having been stimulated by their applications in
different branches of Applied Mathematics including: in perturbed quadrature
rules (see for example \cite{CL5.0}, \cite{CDb5.0}) and in the approximation
of integral transforms (see \cite{DK5.0}, \cite{HDR5.0}) and the references
therein.

For two Lebesgue integrable functions $f,g:\left[ a,b\right] \rightarrow 
\mathbb{R}$, consider the \v{C}eby\v{s}ev functional: 
\begin{equation*}
T\left( f,g\right) :=\frac{1}{b-a}\int_{a}^{b}f\left( t\right) g\left(
t\right) dt-\frac{1}{b-a}\int_{a}^{b}f\left( t\right) dt\cdot \frac{1}{b-a}%
\int_{a}^{b}g\left( t\right) dt.
\end{equation*}

In 1934, G. Gr\"{u}ss \cite{4b5.0} showed that 
\begin{equation}
\left\vert T\left( f,g\right) \right\vert \leq \frac{1}{4}\left( M-m\right)
\left( N-n\right) ,  \label{1.25.0}
\end{equation}%
provided $m,M,n,N$ are real numbers with the property 
\begin{equation}
-\infty <m\leq f\leq M<\infty ,\;-\infty <n\leq g\leq N<\infty \text{ 
\hspace{0.05in}a.e. on }\left[ a,b\right] .  \label{1.35.0}
\end{equation}

The constant $\frac{1}{4}$ is best possible in (\ref{1.25.0}) in the sense
that it cannot be replaced by a smaller one. Another less well known
inequality for $T\left( f,g\right) $ was derived in 1882 by \v{C}eby\v{s}ev 
\cite{3b5.0} under the assumption that $f^{\prime }$, $g^{\prime }$ exist
and are continuous in $\left[ a,b\right] $ and is given by%
\begin{equation}
\left\vert T\left( f,g\right) \right\vert \leq \frac{1}{12}\left\Vert
f^{\prime }\right\Vert _{\infty }\left\Vert g^{\prime }\right\Vert _{\infty
}\left( b-a\right) ^{2},  \label{1.45.0}
\end{equation}%
where $\left\Vert f^{\prime }\right\Vert _{\infty }:=\sup\limits_{t\in \left[
a,b\right] }\left\vert f^{\prime }\left( t\right) \right\vert .$

The constant $\frac{1}{12}$ cannot be improved in the general case.

\v{C}eby\v{s}ev's inequality (\ref{1.45.0}) also holds if $f,g:\left[ a,b%
\right] \rightarrow \mathbb{R}$ are assumed to be absolutely continuous and $%
f^{\prime },g^{\prime }\in L_{\infty }\left[ a,b\right] .$

In 1970, A.M. Ostrowski \cite{5b5.0} proved, amongst others, the following
result that is in a sense a combination of the \v{C}eby\v{s}ev and Gr\"{u}ss
results 
\begin{equation*}
\left\vert T\left( f,g\right) \right\vert \leq \frac{1}{8}\left( b-a\right)
\left( M-m\right) \left\Vert g^{\prime }\right\Vert _{\infty },
\end{equation*}%
provided $f$ is Lebesgue integrable on $\left[ a,b\right] $ and satisfying (%
\ref{1.35.0}) with $g:\left[ a,b\right] \rightarrow \mathbb{R}$ being
absolutely continuous and $g^{\prime }\in L_{\infty }\left[ a,b\right] .$
Here the constant $\frac{1}{8}$ is also sharp.

Finally, let us recall a result by Lupa\c{s} (see for example \cite[p. 210]%
{1b5.0}), which states that: 
\begin{equation*}
\left\vert T\left( f,g\right) \right\vert \leq \frac{1}{\pi ^{2}}\left\Vert
f^{\prime }\right\Vert _{2}\left\Vert g^{\prime }\right\Vert _{2}\left(
b-a\right) ,
\end{equation*}%
provided $f,g$ are absolutely continuous and $f^{\prime },g^{\prime }\in
L_{2}\left[ a,b\right] $. The constant $\frac{1}{\pi ^{2}}$ is the best
possible here also.

For other Gr\"{u}ss type integral inequalities, see the books \cite{4bb5.0}, 
\cite{1b5.0}, and the papers \cite{4b15.0} -- \cite{1bbb5.0}, where further
references are given.

In \cite{1bb5.0}, P. Cerone has obtained the following identity that
involves a Stieltjes integral (Lemma 2.1, p. 3):

\begin{lemma}
\label{l2.1b5.0}Let $f,g:\left[ a,b\right] \rightarrow \mathbb{R}$, where $f$
is of bounded variation and $g$ is continuous on $\left[ a,b\right] ,$ then
the $T\left( f,g\right) $ satisfies the identity, 
\begin{equation}
T\left( f,g\right) =\frac{1}{\left( b-a\right) ^{2}}\int_{a}^{b}\Psi \left(
t\right) df\left( t\right) ,  \label{a.1b5.0}
\end{equation}%
where 
\begin{equation*}
\Psi \left( t\right) :=\left( t-a\right) A\left( t,b\right) -\left(
b-t\right) A\left( a,t\right) ,
\end{equation*}%
with 
\begin{equation*}
A\left( c,d\right) :=\int_{c}^{d}g\left( x\right) dx.
\end{equation*}
\end{lemma}

Using this representation and the properties of Stieltjes integrals he
obtained the following result in bounding the functional $T\left( \cdot
,\cdot \right) $ (Theorem 2.5, p. 4):

\begin{theorem}
\label{t2.2b5.0}With the assumptions in Lemma \ref{l2.1b5.0}, we have: 
\begin{multline*}
\left\vert T\left( f,g\right) \right\vert \\
\leq \frac{1}{\left( b-a\right) ^{2}}\cdot \left\{ 
\begin{array}{ll}
\sup\limits_{t\in \left[ a,b\right] }\left\vert \Psi \left( t\right)
\right\vert \bigvee_{a}^{b}\left( f\right) , &  \\ 
&  \\ 
L\int_{a}^{b}\left\vert \Psi \left( t\right) \right\vert dt, & \text{for }L-%
\text{Lipschitzian;} \\ 
&  \\ 
\int_{a}^{b}\left\vert \Psi \left( t\right) \right\vert df\left( t\right) ,
& \text{for }f\text{ monotonic nondecreasing,}%
\end{array}%
\right.
\end{multline*}%
where $\bigvee_{a}^{b}\left( f\right) $ denotes the total variation of $f$
on $\left[ a,b\right] .$
\end{theorem}

Now, if we use the function $\varphi :\left( a,b\right) \rightarrow \mathbb{R%
}$, 
\begin{equation}
\varphi \left( t\right) :=D\left( g;a,t,b\right) =\frac{\int_{t}^{b}g\left(
x\right) dx}{b-t}-\frac{\int_{a}^{t}g\left( x\right) dx}{t-a},
\label{2.0b5.0}
\end{equation}%
then by (\ref{a.1b5.0}) we may obtain the identity: 
\begin{equation*}
T\left( f,g\right) =\frac{1}{\left( b-a\right) ^{2}}\int_{a}^{b}\left(
t-a\right) \left( b-t\right) \varphi \left( t\right) df\left( t\right) .
\end{equation*}

In \cite{CDb5.0} various upper bounds for $\left\vert T\left( f,g\right)
\right\vert $ have been given, from which we would like to mention only the
following ones

\begin{theorem}
\label{t2.5b5.0}Let $f:\left[ a,b\right] \rightarrow \mathbb{R}$ be a
function of bounded variation and $g:\left[ a,b\right] \rightarrow \mathbb{R}
$ an absolutely continuous function so that $\varphi $ is bounded on $\left(
a,b\right) .$ Then one has the inequality: 
\begin{equation*}
\left\vert T\left( f,g\right) \right\vert \leq \frac{1}{4}\left\Vert \varphi
\right\Vert _{\infty }\bigvee_{a}^{b}\left( f\right) ,
\end{equation*}%
where $\varphi $ is as given by (\ref{2.0b5.0}) and 
\begin{equation*}
\left\Vert \varphi \right\Vert _{\infty }:=\sup\limits_{t\in \left(
a,b\right) }\left\vert \varphi \left( t\right) \right\vert .
\end{equation*}
\end{theorem}

The case of Lipschitzian functions $f:\left[ a,b\right] \rightarrow \mathbb{R%
}$ is embodied in the following theorem as well \cite{CDb5.0}.

\begin{theorem}
\label{t2.6b5.0}Let $f:\left[ a,b\right] \rightarrow \mathbb{R}$ be an $L-$%
Lipschitzian function on $\left[ a,b\right] $ and $g:\left[ a,b\right]
\rightarrow \mathbb{R}$ an absolutely continuous function on $\left[ a,b%
\right] .$ Then 
\begin{eqnarray*}
&&\left\vert T\left( f,g\right) \right\vert \\
&\leq &\left\{ 
\begin{array}{ll}
L\frac{\left( b-a\right) ^{3}}{6}\left\Vert \varphi \right\Vert _{\infty } & 
\text{if \hspace{0.05in}}\varphi \in L_{\infty }\left[ a,b\right] ; \\ 
&  \\ 
L\left( b-a\right) ^{\frac{1}{q}}\left[ B\left( q+1,q+1\right) \right] ^{%
\frac{1}{q}}\left\Vert \varphi \right\Vert _{p}, & p>1,\;\frac{1}{p}+\frac{1%
}{q}=1 \\ 
& \text{if \hspace{0.05in}}\varphi \in L_{p}\left[ a,b\right] ; \\ 
\frac{L}{4}\left\Vert \varphi \right\Vert _{1}, & \text{if \hspace{0.05in}}%
\varphi \in L_{1}\left[ a,b\right] \text{,}%
\end{array}%
\right.
\end{eqnarray*}%
where $\left\Vert \cdot \right\Vert _{p}$ are the usual Lebesgue $p-$norms
on $\left[ a,b\right] $ and $B\left( \cdot ,\cdot \right) $ is Euler's Beta
function.
\end{theorem}

Finally, the following result containing Stieltjes integral holds \cite%
{CDb5.0}:

\begin{theorem}
\label{t2.7b5.0}Let $f:\left[ a,b\right] \rightarrow \mathbb{R}$ be a
monotonic nondecreasing function on $\left[ a,b\right] .$ If $g$ is \hspace{%
0.05in}continuous, then one has the inequality: 
\begin{align*}
& \left\vert T\left( f,g\right) \right\vert \\
& \leq \left\{ 
\begin{array}{l}
\dfrac{1}{4}\dint_{a}^{b}\left\vert \varphi \left( t\right) \right\vert
df\left( t\right) \\ 
\\ 
\dfrac{1}{\left( b-a\right) ^{2}}\left( \dint_{a}^{b}\left[ \left(
b-t\right) \left( t-a\right) \right] ^{q}df\left( t\right) \right) ^{\frac{1%
}{q}}\left( \dint_{a}^{b}\left\vert \varphi \left( t\right) \right\vert
^{p}df\left( t\right) \right) ^{\frac{1}{p}},\; \\ 
\hfill p>1,\;\frac{1}{p}+\frac{1}{q}=1; \\ 
\\ 
\dfrac{1}{\left( b-a\right) ^{2}}\sup\limits_{t\in \left[ a,b\right]
}\left\vert \varphi \left( t\right) \right\vert \dint_{a}^{b}\left(
t-a\right) \left( b-t\right) df\left( t\right) .%
\end{array}%
\right.
\end{align*}
\end{theorem}

In \cite{1abb5.0}, the authors have considered the following functional 
\begin{equation*}
D\left( f;u\right) :=\int_{a}^{b}f\left( x\right) du\left( x\right) -\left[
u\left( b\right) -u\left( a\right) \right] \cdot \frac{1}{b-a}%
\int_{a}^{b}f\left( t\right) dt,
\end{equation*}%
provided that the involved integrals exist.

In the same paper, the following result in estimating the above functional
has been obtained.

\begin{theorem}
\label{ta.1b5.0}Let $f,u:\left[ a,b\right] \rightarrow \mathbb{R}$ be such
that $u$ is Lipschitzian on $\left[ a,b\right] ,$ i.e., 
\begin{equation*}
\left\vert u\left( x\right) -u\left( y\right) \right\vert \leq L\left\vert
x-y\right\vert \text{ \hspace{0.05in}for any }x,y\in \left[ a,b\right]
\;\;\left( L>0\right)
\end{equation*}%
and $f$ is Riemann integrable on $\left[ a,b\right] .$ If $m,M\in \mathbb{R}$
are such that 
\begin{equation*}
m\leq f\left( x\right) \leq M\text{ \hspace{0.05in}for any }x,y\in \left[ a,b%
\right] ,
\end{equation*}%
then we have the inequality 
\begin{equation*}
\left\vert D\left( f;u\right) \right\vert \leq \frac{1}{2}L\left( M-m\right)
\left( b-a\right) .
\end{equation*}%
The constant $\frac{1}{2}$ is sharp in the sense that it cannot be replaced
by a smaller constant.
\end{theorem}

In \cite{1bbb5.0}, the following result complementing the above one was
obtained.

\begin{theorem}
\label{ta.2b5.0}Let $f,u:\left[ a,b\right] \rightarrow \mathbb{R}$ be such
that $u:\left[ a,b\right] \rightarrow \mathbb{R}$ is of bounded variation in 
$\left[ a,b\right] $ and $f:\left[ a,b\right] \rightarrow \mathbb{R}$ is $K-$%
Lipschitzian $\left( K>0\right) .$ Then we have the inequality 
\begin{equation*}
\left\vert D\left( f;u\right) \right\vert \leq \frac{1}{2}K\left( b-a\right)
\bigvee_{a}^{b}\left( u\right) .
\end{equation*}%
The constant $\frac{1}{2}$ is sharp in the above sense.
\end{theorem}

The main aim of this section is to survey some recent inequalities of the Gr%
\"{u}ss type holding in the general setting of inner product spaces. Natural
applications for Lebesgue integrals in measure spaces are presented as well.

\newpage

\section{Gr\"{u}ss' Inequality in Inner Product Spaces}

\subsection{Introduction}

In \cite{SSD.5}, the author has proved the following Gr\"{u}ss' type
inequality in real or complex inner product spaces.

\begin{theorem}
\label{1.1.5}Let $\left( H,\left\langle \cdot ,\cdot \right\rangle \right) $
be an inner product space over $\mathbb{K}\left( \mathbb{K=R}\text{, }%
\mathbb{C}\right) $ and $e\in H,\left\Vert e\right\Vert =1.$ If $\varphi
,\gamma ,\Phi ,\Gamma $ are real or complex numbers and $x,y$ are vectors in 
$H$ such that the conditions 
\begin{equation}
\func{Re}\left\langle \Phi e-x,x-\varphi e\right\rangle \geq 0\text{ and }%
\func{Re}\left\langle \Gamma e-y,y-\gamma e\right\rangle \geq 0
\label{i.1.5}
\end{equation}%
hold, then we have the inequality 
\begin{equation}
\left\vert \left\langle x,y\right\rangle -\left\langle x,e\right\rangle
\left\langle e,y\right\rangle \right\vert \leq \frac{1}{4}\left\vert \Phi
-\varphi \right\vert \cdot \left\vert \Gamma -\gamma \right\vert .
\label{i.2.5}
\end{equation}%
The constant $\frac{1}{4}$ is best possible in the sense that it cannot be
replaced by a smaller constant.
\end{theorem}

Some particular cases of interest for integrable functions with real or
complex values and the corresponding discrete versions are listed below.

\begin{corollary}
\label{1.2.5}Let $f,g:\left[ a,b\right] \rightarrow \mathbb{K}\left( \mathbb{%
K=R}\text{, }\mathbb{C}\right) $ be Lebesgue integrable and such that 
\begin{equation*}
\func{Re}\left[ \left( \Phi -f\left( x\right) \right) \left( \overline{%
f\left( x\right) }-\overline{\varphi }\right) \right] \geq 0,\;\;\;\func{Re}%
\left[ \left( \Gamma -g\left( x\right) \right) \left( \overline{g\left(
x\right) }-\overline{\gamma }\right) \right] \geq 0
\end{equation*}%
for a.e. $x\in \left[ a,b\right] ,$ where $\varphi ,\gamma ,\Phi ,\Gamma $
are real or complex numbers and $\bar{z}$ denotes the complex conjugate of $%
z.$ Then we have the inequality 
\begin{multline*}
\left\vert \frac{1}{b-a}\int_{a}^{b}f\left( x\right) \overline{g\left(
x\right) }dx-\frac{1}{b-a}\int_{a}^{b}f\left( x\right) dx\cdot \frac{1}{b-a}%
\int_{a}^{b}\overline{g\left( x\right) }dx\right\vert \\
\leq \frac{1}{4}\left\vert \Phi -\varphi \right\vert \cdot \left\vert \Gamma
-\gamma \right\vert .
\end{multline*}%
The constant $\frac{1}{4}$ is best possible.
\end{corollary}

The discrete case is embodied in

\begin{corollary}
\label{1.3.5}Let $\mathbf{x,y\in }\mathbb{K}^{n}$ and $\varphi ,\gamma ,\Phi
,\Gamma $ are real or complex numbers such that 
\begin{equation*}
\func{Re}\left[ \left( \Phi -x_{i}\right) \left( \overline{x_{i}}-\overline{%
\varphi }\right) \right] \geq 0,\text{ \ }\func{Re}\left[ \left( \Gamma
-y_{i}\right) \left( \overline{y_{i}}-\overline{\gamma }\right) \right] \geq
0
\end{equation*}%
for each $i\in \left\{ 1,\dots ,n\right\} .$ Then we have the inequality 
\begin{equation*}
\left\vert \frac{1}{n}\sum_{i=1}^{n}x_{i}\overline{y_{i}}-\frac{1}{n}%
\sum_{i=1}^{n}x_{i}\cdot \frac{1}{n}\sum_{i=1}^{n}\overline{y_{i}}%
\right\vert \leq \frac{1}{4}\left\vert \Phi -\varphi \right\vert \cdot
\left\vert \Gamma -\gamma \right\vert .
\end{equation*}%
The constant $\frac{1}{4}$ is best possible.
\end{corollary}

For other applications of Theorem \ref{1.1.5}, see the recent paper \cite%
{SSDIG.5}.

In the present section, by following \cite{05NSSD}, we show that the
condition $\left( \ref{i.1.5}\right) $ may be replaced by an equivalent but
simpler assumption and a new proof of Theorem \ref{1.1.5} is produced. A
refinement of the Gr\"{u}ss type inequality $\left( \ref{i.2.5}\right) ,$
some companions and applications for integrals are pointed out as well.

\subsection{An Equivalent Assumption}

The following lemma holds \cite{05NSSD}.

\begin{lemma}
\label{l2.1.5} Let $a,x,A$ be vectors in the inner product space $\left(
H,\left\langle \cdot ,\cdot \right\rangle \right) $ over $\mathbb{K\ }\left( 
\mathbb{K=R}\text{,}\mathbb{C}\right) $ with $a\neq A.$ Then 
\begin{equation*}
\func{Re}\left\langle A-x,x-a\right\rangle \geq 0
\end{equation*}%
if and only if 
\begin{equation*}
\left\Vert x-\frac{a+A}{2}\right\Vert \leq \frac{1}{2}\left\Vert
A-a\right\Vert .
\end{equation*}
\end{lemma}

\begin{proof}
Define 
\begin{equation*}
I_{1}:=\func{Re}\left\langle A-x,x-a\right\rangle ,\ \ \ \ \ I_{2}:=\frac{1}{%
4}\left\Vert A-a\right\Vert ^{2}-\left\Vert x-\frac{a+A}{2}\right\Vert ^{2}.
\end{equation*}%
A simple calculation shows that 
\begin{equation*}
I_{1}=I_{2}=\func{Re}\left[ \left\langle x,a\right\rangle +\left\langle
A,x\right\rangle \right] -\func{Re}\left\langle A,a\right\rangle -\left\Vert
x\right\Vert ^{2}
\end{equation*}%
and thus, obviously, $I_{1}\geq 0$ iff $I_{2}\geq 0$ showing the required
equivalence.
\end{proof}

The following corollary is obvious

\begin{corollary}
\label{c2.2.5} Let $x,e\in H$ with $\left\Vert e\right\Vert =1$ and $\delta
,\Delta \in \mathbb{K}$ with $\delta \neq \Delta .$ Then 
\begin{equation*}
\func{Re}\left\langle \Delta e-x,x-\delta e\right\rangle \geq 0
\end{equation*}%
iff 
\begin{equation*}
\left\Vert x-\frac{\delta +\Delta }{2}\cdot e\right\Vert \leq \frac{1}{2}%
\left\vert \Delta -\delta \right\vert .
\end{equation*}
\end{corollary}

\begin{remark}
\label{r2.3.5} If $H=\mathbb{C}$, then 
\begin{equation*}
\func{Re}\left[ \left( A-x\right) \left( \bar{x}-\bar{a}\right) \right] \geq
0
\end{equation*}%
if and only if 
\begin{equation*}
\left\vert x-\frac{a+A}{2}\right\vert \leq \frac{1}{2}\left\vert
A-a\right\vert ,
\end{equation*}%
where $a,x,A\in \mathbb{C}$. If $H=\mathbb{R}$, and $A>a$ then $a\leq x\leq
A $ if and only if $\left\vert x-\frac{a+A}{2}\right\vert \leq \frac{1}{2}%
\left\vert A-a\right\vert .$
\end{remark}

The following lemma also holds \cite{05NSSD}.

\begin{lemma}
\label{l2.2.5} Let $x,e\in H$ with $\left\Vert e\right\Vert =1.$ Then one
has the following representation 
\begin{equation}
0\leq \left\Vert x\right\Vert ^{2}-\left\vert \left\langle x,e\right\rangle
\right\vert ^{2}=\inf_{\lambda \in \mathbb{K}}\left\Vert x-\lambda
e\right\Vert ^{2}.  \label{a.1.5}
\end{equation}
\end{lemma}

\begin{proof}
Observe, for any $\lambda \in \mathbb{K}$, that 
\begin{align*}
\left\langle x-\lambda e,x-\left\langle x,e\right\rangle e\right\rangle &
=\left\Vert x\right\Vert ^{2}-\left\vert \left\langle x,e\right\rangle
\right\vert ^{2}-\lambda \left[ \left\langle e,x\right\rangle -\left\langle
e,x\right\rangle \left\Vert e\right\Vert ^{2}\right] \\
& =\left\Vert x\right\Vert ^{2}-\left\vert \left\langle x,e\right\rangle
\right\vert ^{2}.
\end{align*}%
Using Schwarz's inequality, we have 
\begin{align*}
\left[ \left\Vert x\right\Vert ^{2}-\left\vert \left\langle x,e\right\rangle
\right\vert ^{2}\right] ^{2}& =\left\vert \left\langle x-\lambda
e,x-\left\langle x,e\right\rangle e\right\rangle \right\vert ^{2} \\
& \leq \left\Vert x-\lambda e\right\Vert ^{2}\left\Vert x-\left\langle
x,e\right\rangle e\right\Vert ^{2} \\
& =\left\Vert x-\lambda e\right\Vert ^{2}\left[ \left\Vert x\right\Vert
^{2}-\left\vert \left\langle x,e\right\rangle \right\vert ^{2}\right] ,
\end{align*}%
giving the bound 
\begin{equation}
\left\Vert x\right\Vert ^{2}-\left\vert \left\langle x,e\right\rangle
\right\vert ^{2}\leq \left\Vert x-\lambda e\right\Vert ^{2},\ \ \ \ \
\lambda \in \mathbb{K}\text{.}  \label{a.2.5}
\end{equation}%
Taking the infimum in $\left( \ref{a.2.5}\right) $ over $\lambda \in \mathbb{%
K}$, we deduce 
\begin{equation*}
\left\Vert x\right\Vert ^{2}-\left\vert \left\langle x,e\right\rangle
\right\vert ^{2}\leq \inf_{\lambda \in \mathbb{K}}\left\Vert x-\lambda
e\right\Vert ^{2}.
\end{equation*}%
Since, for $\lambda _{0}=\left\langle x,e\right\rangle ,$ we get $\left\Vert
x-\lambda _{0}e\right\Vert ^{2}=\left\Vert x\right\Vert ^{2}-\left\vert
\left\langle x,e\right\rangle \right\vert ^{2},$ then the representation $%
\left( \ref{a.1.5}\right) $ is proved.
\end{proof}

We are now able to provide a different proof for the Gr\"{u}ss type
inequality in inner product spaces (mentioned in the Introduction), than the
one from the paper \cite{SSD.5}.

\begin{theorem}
\label{t2.1.5}Let $\left( H,\left\langle \cdot ,\cdot \right\rangle \right) $
be an inner product space over $\mathbb{K}\left( \mathbb{K=R}\text{,}\mathbb{%
C}\right) $ and $e\in H,\left\Vert e\right\Vert =1.$ If $\varphi ,\gamma
,\Phi ,\Gamma $ are real or complex numbers and $x,y$ are vectors in $H$
such that the conditions $\left( \ref{i.1.5}\right) $ hold, or equivalently,
the following assumptions 
\begin{equation}
\left\Vert x-\frac{\varphi +\Phi }{2}\cdot e\right\Vert \leq \frac{1}{2}%
\left\vert \Phi -\varphi \right\vert ,\text{ \ }\left\Vert y-\frac{\gamma
+\Gamma }{2}\cdot e\right\Vert \leq \frac{1}{2}\left\vert \Gamma -\gamma
\right\vert  \label{a.3.5}
\end{equation}%
are valid, then one has the inequality 
\begin{equation}
\left\vert \left\langle x,y\right\rangle -\left\langle x,e\right\rangle
\left\langle e,y\right\rangle \right\vert \leq \frac{1}{4}\left\vert \Phi
-\varphi \right\vert \cdot \left\vert \Gamma -\gamma \right\vert .
\label{a.4.5}
\end{equation}%
The constant $\frac{1}{4}$ is best possible.
\end{theorem}

\begin{proof}
It can be easily shown (see for example the proof of Theorem 1 from \cite%
{SSD.5}) that 
\begin{equation}
\left\vert \left\langle x,y\right\rangle -\left\langle x,e\right\rangle
\left\langle e,y\right\rangle \right\vert \leq \left[ \left\Vert
x\right\Vert ^{2}-\left\vert \left\langle x,e\right\rangle \right\vert ^{2}%
\right] ^{\frac{1}{2}}\left[ \left\Vert y\right\Vert ^{2}-\left\vert
\left\langle y,e\right\rangle \right\vert ^{2}\right] ^{\frac{1}{2}},
\label{a.5.5}
\end{equation}%
for any $x,y\in H$ and $e\in H,\left\Vert e\right\Vert =1.$ Using Lemma \ref%
{l2.2.5} and the conditions $\left( \ref{a.3.5}\right) $ we obviously have
that 
\begin{equation*}
\left[ \left\Vert x\right\Vert ^{2}-\left\vert \left\langle x,e\right\rangle
\right\vert ^{2}\right] ^{\frac{1}{2}}=\inf_{\lambda \in \mathbb{K}%
}\left\Vert x-\lambda e\right\Vert \leq \left\Vert x-\frac{\varphi +\Phi }{2}%
\cdot e\right\Vert \leq \frac{1}{2}\left\vert \Phi -\varphi \right\vert
\end{equation*}%
and 
\begin{equation*}
\left[ \left\Vert y\right\Vert ^{2}-\left\vert \left\langle y,e\right\rangle
\right\vert ^{2}\right] ^{\frac{1}{2}}=\inf_{\lambda \in \mathbb{K}%
}\left\Vert y-\lambda e\right\Vert \leq \left\Vert y-\frac{\gamma +\Gamma }{2%
}\cdot e\right\Vert \leq \frac{1}{2}\left\vert \Gamma -\gamma \right\vert
\end{equation*}%
and by $\left( \ref{a.5.5}\right) $ the desired inequality $\left( \ref%
{a.4.5}\right) $ is obtained.

The fact that $\frac{1}{4}$ is the best possible constant, has been shown in 
\cite{SSD.5} and we omit the details.
\end{proof}

\subsection{A Refinement of the Gr\"{u}ss Inequality}

The following result improving $\left( \ref{i.1.5}\right) $ holds \cite%
{05NSSD}.

\begin{theorem}
\label{t3.1.5}Let $\left( H,\left\langle \cdot ,\cdot \right\rangle \right) $
be an inner product space over $\mathbb{K}\left( \mathbb{K=R}\text{,}\mathbb{%
C}\right) $ and $e\in H,\left\Vert e\right\Vert =1.$ If $\varphi ,\gamma
,\Phi ,\Gamma $ are real or complex numbers and $x,y$ are vectors in $H$
such that the conditions $\left( \ref{i.1.5}\right) $, or equivalently, $%
\left( \ref{a.3.5}\right) $ hold, then we have the inequality 
\begin{align}
& \left\vert \left\langle x,y\right\rangle -\left\langle x,e\right\rangle
\left\langle e,y\right\rangle \right\vert  \label{3.1.5} \\
& \leq \frac{1}{4}\left\vert \Phi -\varphi \right\vert \cdot \left\vert
\Gamma -\gamma \right\vert -\left[ \func{Re}\left\langle \Phi e-x,x-\varphi
e\right\rangle \right] ^{\frac{1}{2}}\left[ \func{Re}\left\langle \Gamma
e-y,y-\gamma e\right\rangle \right] ^{\frac{1}{2}}  \notag \\
& \leq \left( \frac{1}{4}\left\vert \Phi -\varphi \right\vert \cdot
\left\vert \Gamma -\gamma \right\vert \right) .  \notag
\end{align}%
The constant $\frac{1}{4}$ is best possible.
\end{theorem}

\begin{proof}
As in \cite{SSD.5}, we have 
\begin{equation}
\left\vert \left\langle x,y\right\rangle -\left\langle x,e\right\rangle
\left\langle e,y\right\rangle \right\vert ^{2}\leq \left[ \left\Vert
x\right\Vert ^{2}-\left\vert \left\langle x,e\right\rangle \right\vert ^{2}%
\right] \left[ \left\Vert y\right\Vert ^{2}-\left\vert \left\langle
y,e\right\rangle \right\vert ^{2}\right] ,  \label{3.2.5}
\end{equation}%
\begin{multline}
\left\Vert x\right\Vert ^{2}-\left\vert \left\langle x,e\right\rangle
\right\vert ^{2}  \label{3.3.5} \\
=\func{Re}\left[ \left( \Phi -\left\langle x,e\right\rangle \right) \left( 
\overline{\left\langle x,e\right\rangle }-\overline{\varphi }\right) \right]
-\func{Re}\left\langle \Phi e-x,x-\varphi e\right\rangle
\end{multline}%
and 
\begin{multline}
\left\Vert y\right\Vert ^{2}-\left\vert \left\langle y,e\right\rangle
\right\vert ^{2}  \label{3.4.5} \\
=\func{Re}\left[ \left( \Gamma -\left\langle y,e\right\rangle \right) \left( 
\overline{\left\langle y,e\right\rangle }-\overline{\gamma }\right) \right] -%
\func{Re}\left\langle \Gamma e-x,x-\gamma e\right\rangle .
\end{multline}%
Using the elementary inequality 
\begin{equation*}
4\func{Re}\left( a\overline{b}\right) \leq \left\vert a+b\right\vert ^{2};\
\ \ a,b\in \mathbb{K\ }\left( \mathbb{K=R}\text{,}\mathbb{C}\right) ,
\end{equation*}%
we may state that 
\begin{equation}
\func{Re}\left[ \left( \Phi -\left\langle x,e\right\rangle \right) \left( 
\overline{\left\langle x,e\right\rangle }-\overline{\varphi }\right) \right]
\leq \frac{1}{4}\left\vert \Phi -\varphi \right\vert ^{2}  \label{3.5.5}
\end{equation}%
and 
\begin{equation}
\func{Re}\left[ \left( \Gamma -\left\langle y,e\right\rangle \right) \left( 
\overline{\left\langle y,e\right\rangle }-\overline{\gamma }\right) \right]
\leq \frac{1}{4}\left\vert \Gamma -\gamma \right\vert ^{2}.  \label{3.6.5}
\end{equation}%
Consequently, by $\left( \ref{3.2.5}\right) -\left( \ref{3.6.5}\right) $ we
may state that 
\begin{multline}
\left\vert \left\langle x,y\right\rangle -\left\langle x,e\right\rangle
\left\langle e,y\right\rangle \right\vert ^{2}  \label{3.7.5} \\
\leq \left[ \frac{1}{4}\left\vert \Phi -\varphi \right\vert ^{2}-\left( %
\left[ \func{Re}\left\langle \Phi e-x,x-\varphi e\right\rangle \right] ^{%
\frac{1}{2}}\right) ^{2}\right] \\
\times \left[ \frac{1}{4}\left\vert \Gamma -\gamma \right\vert ^{2}-\left( %
\left[ \func{Re}\left\langle \Gamma e-y,y-\gamma e\right\rangle \right] ^{%
\frac{1}{2}}\right) ^{2}\right] .
\end{multline}%
Finally, using the elementary inequality for positive real numbers 
\begin{equation*}
\left( m^{2}-n^{2}\right) \left( p^{2}-q^{2}\right) \leq \left( mp-nq\right)
^{2},
\end{equation*}%
we have 
\begin{multline*}
\left[ \frac{1}{4}\left\vert \Phi -\varphi \right\vert ^{2}-\left( \left[ 
\func{Re}\left\langle \Phi e-x,x-\varphi e\right\rangle \right] ^{\frac{1}{2}%
}\right) ^{2}\right] \\
\times \left[ \frac{1}{4}\left\vert \Gamma -\gamma \right\vert ^{2}-\left( %
\left[ \func{Re}\left\langle \Gamma e-y,y-\gamma e\right\rangle \right] ^{%
\frac{1}{2}}\right) ^{2}\right] \\
\leq \left( \frac{1}{4}\left\vert \Phi -\varphi \right\vert \cdot \left\vert
\Gamma -\gamma \right\vert -\left[ \func{Re}\left\langle \Phi e-x,x-\varphi
e\right\rangle \right] ^{\frac{1}{2}}\left[ \func{Re}\left\langle \Gamma
e-y,y-\gamma e\right\rangle \right] ^{\frac{1}{2}}\right) ^{2},
\end{multline*}%
giving the desired inequality $\left( \ref{3.1.5}\right) .$
\end{proof}

\subsection{Some Companion Inequalities}

The following companion of the Gr\"{u}ss inequality in inner product spaces
holds \cite{05NSSD}.

\begin{theorem}
\label{t4.1.5} Let $\left( H,\left\langle \cdot ,\cdot \right\rangle \right) 
$ be an inner product space over $\mathbb{K\ }\left( \mathbb{K=R}\text{,}%
\mathbb{C}\right) $ and $e\in H,\left\Vert e\right\Vert =1.$ If $\gamma
,\Gamma \in \mathbb{K}$ and $x,y\in H$ are such that 
\begin{equation}
\func{Re}\left\langle \Gamma e-\frac{x+y}{2},\frac{x+y}{2}-\gamma
e\right\rangle \geq 0  \label{4.1.5}
\end{equation}%
or equivalently, 
\begin{equation*}
\left\Vert \frac{x+y}{2}-\frac{\gamma +\Gamma }{2}\cdot e\right\Vert \leq 
\frac{1}{2}\left\vert \Gamma -\gamma \right\vert ,
\end{equation*}%
then we have the inequality 
\begin{equation}
\func{Re}\left[ \left\langle x,y\right\rangle -\left\langle x,e\right\rangle
\left\langle e,y\right\rangle \right] \leq \frac{1}{4}\left\vert \Gamma
-\gamma \right\vert ^{2}.  \label{4.2.a.5}
\end{equation}%
The constant $\frac{1}{4}$ is best possible in the sense that it cannot be
replaced by a smaller constant.
\end{theorem}

\begin{proof}
Start with the obvious inequality 
\begin{equation}
\func{Re}\left\langle z,u\right\rangle \leq \frac{1}{4}\left\Vert
z+u\right\Vert ^{2};\ \ \ z,u\in H.  \label{4.3.5}
\end{equation}%
Since 
\begin{equation*}
\left\langle x,y\right\rangle -\left\langle x,e\right\rangle \left\langle
e,y\right\rangle =\left\langle x-\left\langle x,e\right\rangle
e,y-\left\langle y,e\right\rangle e\right\rangle ,
\end{equation*}%
then using $\left( \ref{4.3.5}\right) $ we may write 
\begin{align}
\func{Re}\left[ \left\langle x,y\right\rangle -\left\langle x,e\right\rangle
\left\langle e,y\right\rangle \right] & =\func{Re}\left[ \left\langle
x-\left\langle x,e\right\rangle e,y-\left\langle y,e\right\rangle
e\right\rangle \right]  \label{4.4.5} \\
& \leq \frac{1}{4}\left\Vert x-\left\langle x,e\right\rangle
e+y-\left\langle y,e\right\rangle e\right\Vert ^{2}  \notag \\
& =\left\Vert \frac{x+y}{2}-\left\langle \frac{x+y}{2},e\right\rangle \cdot
e\right\Vert ^{2}  \notag \\
& =\left\Vert \frac{x+y}{2}\right\Vert ^{2}-\left\vert \left\langle \frac{x+y%
}{2},e\right\rangle \right\vert ^{2}.  \notag
\end{align}%
If we apply Gr\"{u}ss' inequality in inner product spaces for, say, $a=b=%
\frac{x+y}{2},$ we get 
\begin{equation}
\left\Vert \frac{x+y}{2}\right\Vert ^{2}-\left\vert \left\langle \frac{x+y}{2%
},e\right\rangle \right\vert ^{2}\leq \frac{1}{4}\left\vert \Gamma -\gamma
\right\vert ^{2}.  \label{4.5.5}
\end{equation}%
Making use of $\left( \ref{4.4.5}\right) $ and $\left( \ref{4.5.5}\right) $
we deduce $\left( \ref{4.2.a.5}\right) .$

The fact that $\frac{1}{4}$ is the best possible constant in $\left( \ref%
{4.2.a.5}\right) $ follows by the fact that if in $\left( \ref{4.1.5}\right) 
$ we choose $x=y,$ then it becomes $\func{Re}\left\langle \Gamma
e-x,x-\gamma e\right\rangle \geq 0,$ implying that $0\leq \left\Vert
x\right\Vert ^{2}-\left\vert \left\langle x,e\right\rangle \right\vert
^{2}\leq \frac{1}{4}\left\vert \Gamma -\gamma \right\vert ^{2},$ for which,
by Gr\"{u}ss' inequality in inner product spaces, we know that the constant $%
\frac{1}{4}$ is best possible.
\end{proof}

The following corollary might be of interest if one wanted to evaluate the
absolute value of 
\begin{equation*}
\func{Re}\left[ \left\langle x,y\right\rangle -\left\langle x,e\right\rangle
\left\langle e,y\right\rangle \right] .
\end{equation*}

\begin{corollary}
\label{c4.2.5}Let $\left( H,\left\langle \cdot ,\cdot \right\rangle \right) $
be an inner product space over $\mathbb{K\ }\left( \mathbb{K=R}\text{,}%
\mathbb{C}\right) $ and $e\in H,\left\Vert e\right\Vert =1.$ If $\gamma
,\Gamma \in \mathbb{K}$ and $x,y\in H$ are such that 
\begin{equation*}
\func{Re}\left\langle \Gamma e-\frac{x\pm y}{2},\frac{x\pm y}{2}-\gamma
e\right\rangle \geq 0,
\end{equation*}%
or equivalently, 
\begin{equation*}
\left\Vert \frac{x\pm y}{2}-\frac{\gamma +\Gamma }{2}\cdot e\right\Vert \leq 
\frac{1}{2}\left\vert \Gamma -\gamma \right\vert ,
\end{equation*}%
holds, then we have the inequality 
\begin{equation}
\left\vert \func{Re}\left[ \left\langle x,y\right\rangle -\left\langle
x,e\right\rangle \left\langle e,y\right\rangle \right] \right\vert \leq 
\frac{1}{4}\left\vert \Gamma -\gamma \right\vert ^{2}.  \label{4.8.5}
\end{equation}%
If the inner product space $H$ is real, then $\left( \text{for }m,M\in 
\mathbb{R}\text{, }M>m\right) $ 
\begin{equation*}
\left\langle Me-\frac{x\pm y}{2},\frac{x\pm y}{2}-me\right\rangle \geq 0,
\end{equation*}%
or equivalently, 
\begin{equation*}
\left\Vert \frac{x\pm y}{2}-\frac{m+M}{2}\cdot e\right\Vert \leq \frac{1}{2}%
\left( M-m\right) ,
\end{equation*}%
implies 
\begin{equation}
\left\vert \left\langle x,y\right\rangle -\left\langle x,e\right\rangle
\left\langle e,y\right\rangle \right\vert \leq \frac{1}{4}\left( M-m\right)
^{2}.  \label{4.11.5}
\end{equation}%
In both inequalities $\left( \ref{4.8.5}\right) $ and $\left( \ref{4.11.5}%
\right) ,$ the constant $\frac{1}{4}$ is best possible.
\end{corollary}

\begin{proof}
We only remark that, if 
\begin{equation*}
\func{Re}\left\langle \Gamma e-\frac{x-y}{2},\frac{x-y}{2}-\gamma
e\right\rangle \geq 0
\end{equation*}%
holds, then by Theorem \ref{t4.1.5}, we get 
\begin{equation*}
\func{Re}\left[ -\left\langle x,y\right\rangle +\left\langle
x,e\right\rangle \left\langle e,y\right\rangle \right] \leq \frac{1}{4}%
\left\vert \Gamma -\gamma \right\vert ^{2},
\end{equation*}%
showing that 
\begin{equation}
\func{Re}\left[ \left\langle x,y\right\rangle -\left\langle x,e\right\rangle
\left\langle e,y\right\rangle \right] \geq -\frac{1}{4}\left\vert \Gamma
-\gamma \right\vert ^{2}.  \label{4.12.5}
\end{equation}%
Making use of $\left( \ref{4.2.a.5}\right) $ and $\left( \ref{4.12.5}\right) 
$ we deduce the desired result $\left( \ref{4.8.5}\right) .$
\end{proof}

Finally, we may state and prove the following dual result as well \cite%
{05NSSD}.

\begin{proposition}
\label{p4.1.5}Let $\left( H,\left\langle \cdot ,\cdot \right\rangle \right) $
be an inner product space over $\mathbb{K\ }\left( \mathbb{K=R}\text{,}%
\mathbb{C}\right) $ and $e\in H,\left\Vert e\right\Vert =1.$ If $\varphi
,\Phi \in \mathbb{K}$ and $x,y\in H$ are such that 
\begin{equation}
\func{Re}\left[ \left( \Phi -\left\langle x,e\right\rangle \right) \left( 
\overline{\left\langle x,e\right\rangle }-\overline{\varphi }\right) \right]
\leq 0,  \label{4.13.5}
\end{equation}%
then we have the inequalities 
\begin{align}
\left\Vert x-\left\langle x,e\right\rangle e\right\Vert & \leq \left[ \func{%
Re}\left\langle x-\Phi e,x-\varphi e\right\rangle \right] ^{\frac{1}{2}}
\label{4.14.5} \\
& \leq \frac{\sqrt{2}}{2}\left[ \left\Vert x-\Phi e\right\Vert
^{2}+\left\Vert x-\varphi e\right\Vert ^{2}\right] ^{\frac{1}{2}}.  \notag
\end{align}
\end{proposition}

\begin{proof}
We know that the following identity holds true $($see $\left( \ref{3.3.5}%
\right) )$%
\begin{multline}
\left\Vert x\right\Vert ^{2}-\left\vert \left\langle x,e\right\rangle
\right\vert ^{2}  \label{4.15.5} \\
=\func{Re}\left[ \left( \Phi -\left\langle x,e\right\rangle \right) \left( 
\overline{\left\langle x,e\right\rangle }-\overline{\varphi }\right) \right]
+\func{Re}\left\langle x-\Phi e,x-\varphi e\right\rangle .
\end{multline}%
Using the assumption $\left( \ref{4.13.5}\right) $ and the fact that 
\begin{equation*}
\left\Vert x\right\Vert ^{2}-\left\vert \left\langle x,e\right\rangle
\right\vert ^{2}=\left\Vert x-\left\langle x,e\right\rangle e\right\Vert
^{2},
\end{equation*}%
by $\left( \ref{4.15.5}\right) ,$ we deduce the first inequality in $\left( %
\ref{4.14.5}\right) .$

The second inequality in $\left( \ref{4.14.5}\right) $ follows by the fact
that for any $v,w\in H$ one has 
\begin{equation*}
\func{Re}\left\langle w,v\right\rangle \leq \frac{1}{2}\left( \left\Vert
w\right\Vert ^{2}+\left\Vert v\right\Vert ^{2}\right) .
\end{equation*}%
The proposition is thus proved.
\end{proof}

\subsection{Integral Inequalities}

The following proposition holds \cite{05NSSD}.

\begin{proposition}
\label{p5.1.5} If $f,g,h\in L^{2}\left( \Omega ,\mathbb{K}\right) $ and $%
\varphi ,\Phi ,\gamma ,\Gamma \in \mathbb{K}$, are such that $\int_{\Omega
}\left\vert h\left( s\right) \right\vert ^{2}d\mu \left( s\right) =1$ and 
\begin{align}
\int_{\Omega }\func{Re}\left[ \left( \Phi h\left( s\right) -f\left( s\right)
\right) \left( \overline{f\left( s\right) }-\overline{\varphi }\overline{%
h\left( s\right) }\right) \right] d\mu \left( s\right) & \geq 0,
\label{5.1.5} \\
\int_{\Omega }\func{Re}\left[ \left( \Gamma h\left( s\right) -g\left(
s\right) \right) \left( \overline{g\left( s\right) }-\overline{\gamma }%
\overline{h\left( s\right) }\right) \right] d\mu \left( s\right) & \geq 0 
\notag
\end{align}%
or, equivalently 
\begin{align*}
\left( \int_{\Omega }\left\vert f\left( s\right) -\frac{\Phi +\varphi }{2}%
h\left( s\right) \right\vert ^{2}d\mu \left( s\right) \right) ^{\frac{1}{2}%
}& \leq \frac{1}{2}\left\vert \Phi -\varphi \right\vert , \\
\left( \int_{\Omega }\left\vert g\left( s\right) -\frac{\Gamma +\gamma }{2}%
h\left( s\right) \right\vert ^{2}d\mu \left( s\right) \right) ^{\frac{1}{2}%
}& \leq \frac{1}{2}\left\vert \Gamma -\gamma \right\vert ,
\end{align*}%
hold, then we have the following refinement of the Gr\"{u}ss integral
inequality 
\begin{multline*}
\left\vert \int_{\Omega }f\left( s\right) \overline{g\left( s\right) }d\mu
\left( s\right) -\int_{\Omega }f\left( s\right) \overline{h\left( s\right) }%
d\mu \left( s\right) \int_{\Omega }h\left( s\right) \overline{g\left(
s\right) }d\mu \left( s\right) \right\vert \\
\leq \frac{1}{4}\left\vert \Phi -\varphi \right\vert \cdot \left\vert \Gamma
-\gamma \right\vert -\left[ \int_{\Omega }\func{Re}\left[ \left( \Phi
h\left( s\right) -f\left( s\right) \right) \left( \overline{f\left( s\right) 
}-\overline{\varphi }\overline{h\left( s\right) }\right) \right] d\mu \left(
s\right) \right. \\
\times \left. \int_{\Omega }\func{Re}\left[ \left( \Gamma h\left( s\right)
-g\left( s\right) \right) \left( \overline{g\left( s\right) }-\overline{%
\gamma }\overline{h\left( s\right) }\right) \right] d\mu \left( s\right) %
\right] ^{\frac{1}{2}}.
\end{multline*}%
The constant $\frac{1}{4}$ is best possible.
\end{proposition}

The proof follows by Theorem \ref{t3.1.5} on choosing $H=L^{2}\left( \Omega ,%
\mathbb{K}\right) $ with the inner product 
\begin{equation*}
\left\langle f,g\right\rangle :=\int_{\Omega }f\left( s\right) \overline{%
g\left( s\right) }d\mu \left( s\right) .
\end{equation*}%
We omit the details.

\begin{remark}
\label{r.1.5} It is obvious that a sufficient condition for $\left( \ref%
{5.1.5}\right) $ to hold is 
\begin{equation*}
\func{Re}\left[ \left( \Phi h\left( s\right) -f\left( s\right) \right)
\left( \overline{f\left( s\right) }-\overline{\varphi }\overline{h\left(
s\right) }\right) \right] \geq 0,
\end{equation*}%
and 
\begin{equation*}
\func{Re}\left[ \left( \Gamma h\left( s\right) -g\left( s\right) \right)
\left( \overline{g\left( s\right) }-\overline{\gamma }\overline{h\left(
s\right) }\right) \right] \geq 0,
\end{equation*}%
for $\mu -$a.e.$\;s\in \Omega ,$ or equivalently, 
\begin{align*}
\left\vert f\left( s\right) -\frac{\Phi +\varphi }{2}h\left( s\right)
\right\vert & \leq \frac{1}{2}\left\vert \Phi -\varphi \right\vert
\left\vert h\left( s\right) \right\vert \text{ \ \ and} \\
\left\vert g\left( s\right) -\frac{\Gamma +\gamma }{2}h\left( s\right)
\right\vert & \leq \frac{1}{2}\left\vert \Gamma -\gamma \right\vert
\left\vert h\left( s\right) \right\vert ,
\end{align*}%
for $\mu -$a.e.$\;s\in \Omega .$
\end{remark}

The following result may be stated as well \cite{05NSSD}.

\begin{corollary}
\label{c5.1.5}If $z,Z,t,T\in \mathbb{K}$, $\mu \left( \Omega \right) <\infty 
$ and $f,g\in L^{2}\left( \Omega ,\mathbb{K}\right) $ are such that: 
\begin{align*}
\func{Re}\left[ \left( Z-f\left( s\right) \right) \left( \overline{f\left(
s\right) }-\bar{z}\right) \right] & \geq 0, \\
\func{Re}\left[ \left( T-g\left( s\right) \right) \left( \overline{g\left(
s\right) }-\bar{t}\right) \right] & \geq 0,\text{ \hspace{0.05in}for a.e. }%
s\in \Omega ,
\end{align*}%
or, equivalently 
\begin{align*}
\left\vert f\left( s\right) -\frac{z+Z}{2}\right\vert & \leq \frac{1}{2}%
\left\vert Z-z\right\vert , \\
\left\vert g\left( s\right) -\frac{t+T}{2}\right\vert & \leq \frac{1}{2}%
\left\vert T-t\right\vert ,\text{ \hspace{0.05in}for a.e. }s\in \Omega ,
\end{align*}%
then we have the inequality 
\begin{multline}
\left\vert \frac{1}{\mu \left( \Omega \right) }\int_{\Omega }f\left(
s\right) \overline{g\left( s\right) }d\mu \left( s\right) \right. -\left. 
\frac{1}{\mu \left( \Omega \right) }\int_{\Omega }f\left( s\right) d\mu
\left( s\right) \cdot \frac{1}{\mu \left( \Omega \right) }\int_{\Omega }%
\overline{g\left( s\right) }d\mu \left( s\right) \right\vert  \label{5.6.5}
\\
\leq \frac{1}{4}\left\vert Z-z\right\vert \left\vert T-t\right\vert -\frac{1%
}{\mu \left( \Omega \right) }\left[ \int_{\Omega }\func{Re}\left[ \left(
Z-f\left( s\right) \right) \left( \overline{f\left( s\right) }-\bar{z}%
\right) \right] d\mu \left( s\right) \right. \\
\times \left. \int_{\Omega }\func{Re}\left[ \left( T-g\left( s\right)
\right) \left( \overline{g\left( s\right) }-\bar{t}\right) \right] d\mu
\left( s\right) \right] ^{\frac{1}{2}}.
\end{multline}
\end{corollary}

Using Theorem \ref{t4.1.5} we may state the following result as well \cite%
{05NSSD}.

\begin{proposition}
\label{p5.3.5}If $f,g,h\in L^{2}\left( \Omega ,\mathbb{K}\right) $ and $%
\gamma ,\Gamma \in \mathbb{K}$ are such that $\int_{\Omega }\left\vert
h\left( s\right) \right\vert ^{2}d\mu \left( s\right) =1$ and 
\begin{multline}
\int_{\Omega }\func{Re}\left\{ \left[ \Gamma h\left( s\right) -\frac{f\left(
s\right) +g\left( s\right) }{2}\right] \right.  \label{5.7.5} \\
\times \left. \left[ \frac{\overline{f\left( s\right) }+\overline{g\left(
s\right) }}{2}-\bar{\gamma}\bar{h}\left( s\right) \right] \right\} d\mu
\left( s\right) \geq 0
\end{multline}%
or equivalently, 
\begin{equation}
\left( \int_{\Omega }\left\vert \frac{f\left( s\right) +g\left( s\right) }{2}%
-\frac{\gamma +\Gamma }{2}h\left( s\right) \right\vert ^{2}d\mu \left(
s\right) \right) ^{\frac{1}{2}}\leq \frac{1}{2}\left\vert \Gamma -\gamma
\right\vert ,  \label{5.8.5}
\end{equation}%
holds, then we have the inequality 
\begin{align*}
I& :=\int_{\Omega }\func{Re}\left[ f\left( s\right) \overline{g\left(
s\right) }\right] d\mu \left( s\right) \\
& \qquad \qquad \qquad -\func{Re}\left[ \int_{\Omega }f\left( s\right) 
\overline{h\left( s\right) }d\mu \left( s\right) \cdot \int_{\Omega }h\left(
s\right) \overline{g\left( s\right) }d\mu \left( s\right) \right] \\
& \leq \frac{1}{4}\left\vert \Gamma -\gamma \right\vert ^{2}.
\end{align*}%
If (\ref{5.7.5}) and (\ref{5.8.5}) hold with \textquotedblleft\ $\pm $
\textquotedblright\ instead of \textquotedblleft\ $+$ \textquotedblright ,
then 
\begin{equation*}
\left\vert I\right\vert \leq \frac{1}{4}\left\vert \Gamma -\gamma
\right\vert ^{2}.
\end{equation*}
\end{proposition}

\begin{remark}
It is obvious that a sufficient condition for (\ref{5.7.5}) to hold is 
\begin{equation*}
\func{Re}\left\{ \left[ \Gamma h\left( s\right) -\frac{f\left( s\right)
+g\left( s\right) }{2}\right] \cdot \left[ \frac{\overline{f\left( s\right) }%
+\overline{g\left( s\right) }}{2}-\bar{\gamma}\bar{h}\left( s\right) \right]
\right\} \geq 0,
\end{equation*}%
for a.e. $s\in \Omega ,$ or equivalently, 
\begin{equation*}
\left\vert \frac{f\left( s\right) +g\left( s\right) }{2}-\frac{\gamma
+\Gamma }{2}h\left( s\right) \right\vert \leq \frac{1}{2}\left\vert \Gamma
-\gamma \right\vert \left\vert h\left( s\right) \right\vert ,\text{ \hspace{%
0.05in}for a.e. }s\in \Omega .
\end{equation*}
\end{remark}

Finally, the following corollary holds.

\begin{corollary}
\label{c5.4.5}If $Z,z\in \mathbb{K}$, $\mu \left( \Omega \right) <\infty $
and $f,g\in L^{2}\left( \Omega ,\mathbb{K}\right) $ are such that 
\begin{equation}
\func{Re}\left[ \left( Z-\frac{f\left( s\right) +g\left( s\right) }{2}%
\right) \left( \frac{\overline{f\left( s\right) }+\overline{g\left( s\right) 
}}{2}-\overline{z}\right) \right] \geq 0,\text{ \hspace{0.05in}for a.e. }%
s\in \Omega  \label{5.13.5}
\end{equation}%
or, equivalently 
\begin{equation}
\left\vert \frac{f\left( s\right) +g\left( s\right) }{2}-\frac{z+Z}{2}%
\right\vert \leq \frac{1}{2}\left\vert Z-z\right\vert ,\text{ \hspace{0.05in}%
for a.e. }s\in \Omega ,  \label{5.14.5}
\end{equation}%
then we have the inequality 
\begin{align*}
J& :=\frac{1}{\mu \left( \Omega \right) }\int_{\Omega }\func{Re}\left[
f\left( s\right) \overline{g\left( s\right) }\right] d\mu \left( s\right) \\
& \qquad \qquad \qquad -\func{Re}\left[ \frac{1}{\mu \left( \Omega \right) }%
\int_{\Omega }f\left( s\right) d\mu \left( s\right) \cdot \frac{1}{\mu
\left( \Omega \right) }\int_{\Omega }\overline{g\left( s\right) }d\mu \left(
s\right) \right] \\
& \leq \frac{1}{4}\left\vert Z-z\right\vert ^{2}.
\end{align*}%
If (\ref{5.13.5}) and (\ref{5.14.5}) hold with \textquotedblleft\ $\pm $
\textquotedblright\ instead of \textquotedblleft\ $+$ \textquotedblright $,$
then 
\begin{equation*}
\left\vert J\right\vert \leq \frac{1}{4}\left\vert Z-z\right\vert ^{2}.
\end{equation*}
\end{corollary}

\begin{remark}
It is obvious that if one chooses the discrete measure above, then all the
inequalities in this section may be written for sequences of real or complex
numbers. We omit the details.
\end{remark}

\newpage

\section{Companions of Gr\"{u}ss' Inequality}

\subsection{A General Result}

The following Gr\"{u}ss type inequality in inner product spaces holds \cite%
{6NSSD}.

\begin{theorem}
\label{t2.1.6}Let $x,y,e\in H$ with $\left\Vert e\right\Vert =1,$ and the
scalars $a,A,b,B\in \mathbb{K}$ $\left( \mathbb{K}=\mathbb{C},\mathbb{R}%
\right) $ such that $\func{Re}\left( \bar{a}A\right) >0$ and $\func{Re}%
\left( \bar{b}B\right) >0.$ If 
\begin{equation*}
\func{Re}\left\langle Ae-x,x-ae\right\rangle \geq 0\text{ \hspace{0.05in}and 
\hspace{0.05in}}\func{Re}\left\langle Be-y,y-be\right\rangle \geq 0
\end{equation*}%
or equivalently, 
\begin{equation*}
\left\Vert x-\frac{a+A}{2}e\right\Vert \leq \frac{1}{2}\left\vert
A-a\right\vert \text{ \hspace{0.05in}and \hspace{0.05in}}\left\Vert y-\frac{%
b+B}{2}e\right\Vert \leq \frac{1}{2}\left\vert B-b\right\vert ,
\end{equation*}%
holds, then we have the inequality 
\begin{equation}
\left\vert \left\langle x,y\right\rangle -\left\langle x,e\right\rangle
\left\langle e,y\right\rangle \right\vert \leq \frac{1}{4}M\left( a,A\right)
M\left( b,B\right) \left\vert \left\langle x,e\right\rangle \left\langle
e,y\right\rangle \right\vert ,  \label{2.3.6}
\end{equation}%
where $M\left( \cdot ,\cdot \right) $ is defined by 
\begin{equation*}
M\left( a,A\right) :=\left[ \frac{\left( \left\vert A\right\vert -\left\vert
a\right\vert \right) ^{2}+4\left[ \left\vert A\bar{a}\right\vert -\func{Re}%
\left( A\bar{a}\right) \right] }{\func{Re}\left( \bar{a}A\right) }\right] ^{%
\frac{1}{2}}.
\end{equation*}%
The constant $\frac{1}{4}$ is best possible in the sense that it cannot be
replaced by a smaller constant.
\end{theorem}

\begin{proof}
Start with the inequality 
\begin{equation}
\left\vert \left\langle x,y\right\rangle -\left\langle x,e\right\rangle
\left\langle e,y\right\rangle \right\vert ^{2}\leq \left( \left\Vert
x\right\Vert ^{2}-\left\vert \left\langle x,e\right\rangle \right\vert
^{2}\right) \left( \left\Vert y\right\Vert ^{2}-\left\vert \left\langle
y,e\right\rangle \right\vert ^{2}\right) .  \label{2.5.6}
\end{equation}%
Now, assume that $u,v\in H,$ and $c,C\in \mathbb{K}$ with $\func{Re}\left( 
\bar{c}C\right) >0$ and $\func{Re}\left\langle Cv-u,u-cv\right\rangle \geq
0. $ This last inequality is equivalent to 
\begin{equation*}
\left\Vert u\right\Vert ^{2}+\func{Re}\left( \bar{c}C\right) \left\Vert
v\right\Vert ^{2}\leq \func{Re}\left[ C\overline{\left\langle
u,v\right\rangle }+\bar{c}\left\langle u,v\right\rangle \right] .
\end{equation*}%
Dividing this inequality by $\left[ \func{Re}\left( C\bar{c}\right) \right]
^{\frac{1}{2}}>0,$ we deduce 
\begin{equation}
\frac{1}{\left[ \func{Re}\left( \bar{c}C\right) \right] ^{\frac{1}{2}}}%
\left\Vert u\right\Vert ^{2}+\left[ \func{Re}\left( \bar{c}C\right) \right]
^{\frac{1}{2}}\left\Vert v\right\Vert ^{2}\leq \frac{\func{Re}\left[ C%
\overline{\left\langle u,v\right\rangle }+\bar{c}\left\langle
u,v\right\rangle \right] }{\left[ \func{Re}\left( \bar{c}C\right) \right] ^{%
\frac{1}{2}}}.  \label{2.7.6}
\end{equation}%
On the other hand, by the elementary inequality 
\begin{equation*}
\alpha p^{2}+\frac{1}{\alpha }q^{2}\geq 2pq,\,\;\;\;\alpha >0,\;p,q\geq 0,
\end{equation*}%
we deduce 
\begin{equation}
2\left\Vert u\right\Vert \left\Vert v\right\Vert \leq \frac{1}{\left[ \func{%
Re}\left( \bar{c}C\right) \right] ^{\frac{1}{2}}}\left\Vert u\right\Vert
^{2}+\left[ \func{Re}\left( \bar{c}C\right) \right] ^{\frac{1}{2}}\left\Vert
v\right\Vert ^{2}.  \label{2.8.6}
\end{equation}%
Making use of (\ref{2.7.6}) and (\ref{2.8.6}) and the fact that for any $%
z\in \mathbb{C}$, $\func{Re}\left( z\right) \leq \left\vert z\right\vert ,$
we get 
\begin{equation*}
\left\Vert u\right\Vert \left\Vert v\right\Vert \leq \frac{\func{Re}\left[ C%
\overline{\left\langle u,v\right\rangle }+\bar{c}\left\langle
u,v\right\rangle \right] }{2\left[ \func{Re}\left( \bar{c}C\right) \right] ^{%
\frac{1}{2}}}\leq \frac{\left\vert c\right\vert +\left\vert C\right\vert }{2%
\left[ \func{Re}\left( \bar{c}C\right) \right] ^{\frac{1}{2}}}\left\vert
\left\langle u,v\right\rangle \right\vert .
\end{equation*}%
Consequently 
\begin{align}
\left\Vert u\right\Vert ^{2}\left\Vert v\right\Vert ^{2}-\left\vert
\left\langle u,v\right\rangle \right\vert ^{2}& \leq \left[ \frac{\left(
\left\vert c\right\vert +\left\vert C\right\vert \right) ^{2}}{4\left[ \func{%
Re}\left( \bar{c}C\right) \right] }-1\right] \left\vert \left\langle
u,v\right\rangle \right\vert ^{2}  \label{2.9.6} \\
& =\frac{1}{4}\frac{\left( \left\vert c\right\vert -\left\vert C\right\vert
\right) ^{2}+4\left[ \left\vert \bar{c}C\right\vert -\func{Re}\left( \bar{c}%
C\right) \right] }{\func{Re}\left( \bar{c}C\right) }\left\vert \left\langle
u,v\right\rangle \right\vert ^{2}  \notag \\
& =\frac{1}{4}M^{2}\left( c,C\right) \left\vert \left\langle
u,v\right\rangle \right\vert ^{2}.  \notag
\end{align}%
Now, if we write (\ref{2.9.6}) for the choices $u=x,$ $v=e$ and $u=y,$ $v=e$
respectively and use (\ref{2.5.6}), we deduce the desired result (\ref{2.3.6}%
). The sharpness of the constant will be proved in the case where $H$ is a
real inner product space.
\end{proof}

The following corollary which provides a simpler Gr\"{u}ss type inequality
for real constants (and in particular, for real inner product spaces) holds 
\cite{6NSSD}.

\begin{corollary}
\label{c2.2.6}With the assumptions of Theorem \ref{t2.1.6} and if $%
a,b,A,B\in \mathbb{R}$ are such that $A>a>0,$ $B>b>0$ and 
\begin{equation}
\left\Vert x-\frac{a+A}{2}e\right\Vert \leq \frac{1}{2}\left( A-a\right) 
\text{ \hspace{0.05in}and \hspace{0.05in}}\left\Vert y-\frac{b+B}{2}%
e\right\Vert \leq \frac{1}{2}\left( B-b\right) ,  \label{2.10.6}
\end{equation}%
then we have the inequality 
\begin{equation}
\left\vert \left\langle x,y\right\rangle -\left\langle x,e\right\rangle
\left\langle e,y\right\rangle \right\vert \leq \frac{1}{4}\cdot \frac{\left(
A-a\right) \left( B-b\right) }{\sqrt{abAB}}\left\vert \left\langle
x,e\right\rangle \left\langle e,y\right\rangle \right\vert .  \label{2.11.6}
\end{equation}%
The constant $\frac{1}{4}$ is best possible.
\end{corollary}

\begin{proof}
The prove the sharpness of the constant $\frac{1}{4}$ assume that the
inequality (\ref{2.11.6}) holds in real inner product spaces with $x=y$ and
for a constant $k>0,$ i.e., 
\begin{equation}
\left\Vert x\right\Vert ^{2}-\left\vert \left\langle x,e\right\rangle
\right\vert ^{2}\leq k\cdot \frac{\left( A-a\right) ^{2}}{aA}\left\vert
\left\langle x,e\right\rangle \right\vert ^{2}\;\;\;\;\left( A>a>0\right) ,
\label{2.11.a.6}
\end{equation}%
provided that $\left\Vert x-\frac{a+A}{2}e\right\Vert \leq \frac{1}{2}\left(
A-a\right) ,$ or equivalently, $\left\langle Ae-x,x-ae\right\rangle \geq 0.$

We choose $H=\mathbb{R}^{2},$ $x=\left( x_{1},x_{2}\right) \in \mathbb{R}%
^{2},$ $e=\left( \frac{1}{\sqrt{2}},\frac{1}{\sqrt{2}}\right) .$ Then we
have 
\begin{align*}
\left\Vert x\right\Vert ^{2}-\left\vert \left\langle x,e\right\rangle
\right\vert ^{2}& =x_{1}^{2}+x_{2}^{2}-\frac{\left( x_{1}+x_{2}\right) ^{2}}{%
2}=\frac{\left( x_{1}-x_{2}\right) ^{2}}{2}, \\
\left\vert \left\langle x,e\right\rangle \right\vert ^{2}& =\frac{\left(
x_{1}+x_{2}\right) ^{2}}{2},
\end{align*}%
and by (\ref{2.11.a.6}) we get 
\begin{equation}
\frac{\left( x_{1}-x_{2}\right) ^{2}}{2}\leq k\cdot \frac{\left( A-a\right)
^{2}}{aA}\cdot \frac{\left( x_{1}+x_{2}\right) ^{2}}{2}.  \label{2.11.b.6}
\end{equation}%
Now, if we let $x_{1}=\frac{a}{\sqrt{2}},$ $x_{2}=\frac{A}{\sqrt{2}}$ $%
\left( A>a>0\right) ,$ then obviously 
\begin{equation*}
\left\langle Ae-x,x-ae\right\rangle =\sum_{i=1}^{2}\left( \frac{A}{\sqrt{2}}%
-x_{i}\right) \left( x_{i}-\frac{a}{\sqrt{2}}\right) =0,
\end{equation*}%
which shows that the condition (\ref{2.10.6}) is fulfilled, and by (\ref%
{2.11.b.6}) we get 
\begin{equation*}
\frac{\left( A-a\right) ^{2}}{4}\leq k\cdot \frac{\left( A-a\right) ^{2}}{aA}%
\cdot \frac{\left( a+A\right) ^{2}}{4}
\end{equation*}%
for any $A>a>0.$ This implies 
\begin{equation}
\left( A+a\right) ^{2}k\geq aA  \label{2.11.c.6}
\end{equation}%
for any $A>a>0.$

Finally, let $a=1-q$, $A=1+q,$ $q\in \left( 0,1\right) .$ Then from (\ref%
{2.11.c.6}) we get $4k\geq 1-q^{2}$ for any $q\in \left( 0,1\right) $ which
produces $k\geq \frac{1}{4}.$
\end{proof}

\begin{remark}
\label{r2.3.6}If $\left\langle x,e\right\rangle ,\left\langle
y,e\right\rangle $ are assumed not to be zero, then the inequality (\ref%
{2.3.6}) is equivalent to 
\begin{equation*}
\left\vert \frac{\left\langle x,y\right\rangle }{\left\langle
x,e\right\rangle \left\langle e,y\right\rangle }-1\right\vert \leq \frac{1}{4%
}M\left( a,A\right) M\left( b,B\right) ,
\end{equation*}%
while the inequality (\ref{2.11.6}) is equivalent to 
\begin{equation*}
\left\vert \frac{\left\langle x,y\right\rangle }{\left\langle
x,e\right\rangle \left\langle e,y\right\rangle }-1\right\vert \leq \frac{1}{4%
}\cdot \frac{\left( A-a\right) \left( B-b\right) }{\sqrt{abAB}}.
\end{equation*}%
The constant $\frac{1}{4}$ is best possible in both inequalities.
\end{remark}

\subsection{Some Related Results}

The following result holds \cite{6NSSD}.

\begin{theorem}
\label{t3.1.6}Let $\left( H;\left\langle \cdot ,\cdot \right\rangle \right) $
be an inner product space over $\mathbb{K}$ $\left( \mathbb{K}=\mathbb{C},%
\mathbb{R}\right) .$ If $\gamma ,\Gamma \in \mathbb{K}$, $e,x,y\in H$ with $%
\left\Vert e\right\Vert =1$ and $\lambda \in \left( 0,1\right) $ are such
that 
\begin{equation}
\func{Re}\left\langle \Gamma e-\left( \lambda x+\left( 1-\lambda \right)
y\right) ,\left( \lambda x+\left( 1-\lambda \right) y\right) -\gamma
e\right\rangle \geq 0,  \label{3.1.6}
\end{equation}%
or equivalently, 
\begin{equation*}
\left\Vert \lambda x+\left( 1-\lambda \right) y-\frac{\gamma +\Gamma }{2}%
e\right\Vert \leq \frac{1}{2}\left\vert \Gamma -\gamma \right\vert ,
\end{equation*}%
then we have the inequality 
\begin{equation}
\func{Re}\left[ \left\langle x,y\right\rangle -\left\langle x,e\right\rangle
\left\langle e,y\right\rangle \right] \leq \frac{1}{16}\cdot \frac{1}{%
\lambda \left( 1-\lambda \right) }\left\vert \Gamma -\gamma \right\vert ^{2}.
\label{3.3.6}
\end{equation}%
The constant $\frac{1}{16}$ is the best possible constant in (\ref{3.3.6})
in the sense that it cannot be replaced by a smaller one.
\end{theorem}

\begin{proof}
We know that for any $z,u\in H$ one has 
\begin{equation*}
\func{Re}\left\langle z,u\right\rangle \leq \frac{1}{4}\left\Vert
z+u\right\Vert ^{2}.
\end{equation*}%
Then for any $a,b\in H$ and $\lambda \in \left( 0,1\right) $ one has 
\begin{equation}
\func{Re}\left\langle a,b\right\rangle \leq \frac{1}{4\lambda \left(
1-\lambda \right) }\left\Vert \lambda a+\left( 1-\lambda \right)
b\right\Vert ^{2}.  \label{3.4.6}
\end{equation}%
Since 
\begin{equation*}
\left\langle x,y\right\rangle -\left\langle x,e\right\rangle \left\langle
e,y\right\rangle =\left\langle x-\left\langle x,e\right\rangle
e,y-\left\langle y,e\right\rangle e\right\rangle \;\;\;\;\;\text{(as }%
\left\Vert e\right\Vert =1\text{),}
\end{equation*}%
using (\ref{3.4.6}), we have 
\begin{align}
& \func{Re}\left[ \left\langle x,y\right\rangle -\left\langle
x,e\right\rangle \left\langle e,y\right\rangle \right]  \label{3.5.6} \\
& =\func{Re}\left[ \left\langle x-\left\langle x,e\right\rangle
e,y-\left\langle y,e\right\rangle e\right\rangle \right]  \notag \\
& \leq \frac{1}{4\lambda \left( 1-\lambda \right) }\left\Vert \lambda \left(
x-\left\langle x,e\right\rangle e\right) +\left( 1-\lambda \right) \left(
y-\left\langle y,e\right\rangle e\right) \right\Vert ^{2}  \notag \\
& =\frac{1}{4\lambda \left( 1-\lambda \right) }\left\Vert \lambda x+\left(
1-\lambda \right) y-\left\langle \lambda x+\left( 1-\lambda \right)
y,e\right\rangle e\right\Vert ^{2}.  \notag
\end{align}%
Since, for $m,e\in H$ with $\left\Vert e\right\Vert =1,$ one has the
equality 
\begin{equation*}
\left\Vert m-\left\langle m,e\right\rangle e\right\Vert ^{2}=\left\Vert
m\right\Vert ^{2}-\left\vert \left\langle m,e\right\rangle \right\vert ^{2},
\end{equation*}%
then by (\ref{3.5.6}) we deduce the inequality 
\begin{multline}
\func{Re}\left[ \left\langle x,y\right\rangle -\left\langle x,e\right\rangle
\left\langle e,y\right\rangle \right]  \label{3.7.6} \\
\leq \frac{1}{4\lambda \left( 1-\lambda \right) }\left[ \left\Vert \lambda
x+\left( 1-\lambda \right) y\right\Vert ^{2}-\left\vert \left\langle \lambda
x+\left( 1-\lambda \right) y,e\right\rangle \right\vert ^{2}\right] .
\end{multline}%
Now, if we apply Gr\"{u}ss' inequality 
\begin{equation*}
0\leq \left\Vert a\right\Vert ^{2}-\left\vert \left\langle a,e\right\rangle
\right\vert ^{2}\leq \frac{1}{4}\left\vert \Gamma -\gamma \right\vert ^{2},
\end{equation*}%
provided 
\begin{equation*}
\func{Re}\left\langle \Gamma e-a,a-\gamma e\right\rangle \geq 0,
\end{equation*}%
for $a=\lambda x+\left( 1-\lambda \right) y,$ we have 
\begin{equation}
\left\Vert \lambda x+\left( 1-\lambda \right) y\right\Vert ^{2}-\left\vert
\left\langle \lambda x+\left( 1-\lambda \right) y,e\right\rangle \right\vert
^{2}\leq \frac{1}{4}\left\vert \Gamma -\gamma \right\vert ^{2}.
\label{3.8.6}
\end{equation}%
Utilising (\ref{3.7.6}) and (\ref{3.8.6}) we deduce the desired inequality (%
\ref{3.3.6}).

To prove the sharpness of the constant $\frac{1}{16}$, assume that (\ref%
{3.3.6}) holds with a constant $C>0,$ provided that (\ref{3.1.6}) is valid,
i.e., 
\begin{equation}
\func{Re}\left[ \left\langle x,y\right\rangle -\left\langle x,e\right\rangle
\left\langle e,y\right\rangle \right] \leq C\cdot \frac{1}{\lambda \left(
1-\lambda \right) }\left\vert \Gamma -\gamma \right\vert ^{2}.  \label{3.9.6}
\end{equation}%
If in (\ref{3.9.6}) we choose $x=y,$ given that (\ref{3.1.6}) holds with $%
x=y $ and $\lambda \in \left( 0,1\right) ,$ then 
\begin{equation}
\left\Vert x\right\Vert ^{2}-\left\vert \left\langle x,e\right\rangle
\right\vert ^{2}\leq C\cdot \frac{1}{\lambda \left( 1-\lambda \right) }%
\left\vert \Gamma -\gamma \right\vert ^{2},  \label{3.10.6}
\end{equation}%
provided 
\begin{equation*}
\func{Re}\left\langle \Gamma e-x,x-\gamma e\right\rangle \geq 0.
\end{equation*}%
Since we know, in Gr\"{u}ss' inequality, that the constant $\frac{1}{4}$ is
best possible, then by (\ref{3.10.6}), one has 
\begin{equation*}
\frac{1}{4}\leq \frac{C}{\lambda \left( 1-\lambda \right) }\text{ \hspace{%
0.05in}for \hspace{0.05in}}\lambda \in \left( 0,1\right) ,
\end{equation*}%
giving, for $\lambda =\frac{1}{2},$ $C\geq \frac{1}{16}.$

The theorem is completely proved.
\end{proof}

The following corollary is a natural consequence of the above result \cite%
{6NSSD}.

\begin{corollary}
\label{c3.2.6}Assume that $\gamma ,\Gamma ,e,x,y$ and $\lambda $ are as in
Theorem \ref{t3.1.6}. If 
\begin{equation*}
\func{Re}\left\langle \Gamma e-\left( \lambda x\pm \left( 1-\lambda \right)
y\right) ,\left( \lambda x\pm \left( 1-\lambda \right) y\right) -\gamma
e\right\rangle \geq 0,
\end{equation*}%
or equivalently, 
\begin{equation*}
\left\Vert \lambda x\pm \left( 1-\lambda \right) y-\frac{\gamma +\Gamma }{2}%
e\right\Vert \leq \frac{1}{2}\left\vert \Gamma -\gamma \right\vert ^{2},
\end{equation*}%
then we have the inequality 
\begin{equation}
\left\vert \func{Re}\left[ \left\langle x,y\right\rangle -\left\langle
x,e\right\rangle \left\langle e,y\right\rangle \right] \right\vert \leq 
\frac{1}{16}\cdot \frac{1}{\lambda \left( 1-\lambda \right) }\left\vert
\Gamma -\gamma \right\vert ^{2}.  \label{3.13.6}
\end{equation}%
The constant $\frac{1}{16}$ is best possible in (\ref{3.13.6}).
\end{corollary}

\begin{proof}
Using Theorem \ref{t3.1.6} for $\left( -y\right) $ instead of $y$, we have
that 
\begin{equation*}
\func{Re}\left\langle \Gamma e-\left( \lambda x-\left( 1-\lambda \right)
y\right) ,\left( \lambda x-\left( 1-\lambda \right) y\right) -\gamma
e\right\rangle \geq 0,
\end{equation*}%
which implies that 
\begin{equation*}
\func{Re}\left[ -\left\langle x,y\right\rangle +\left\langle
x,e\right\rangle \left\langle e,y\right\rangle \right] \leq \frac{1}{16}%
\cdot \frac{1}{\lambda \left( 1-\lambda \right) }\left\vert \Gamma -\gamma
\right\vert ^{2}
\end{equation*}%
giving 
\begin{equation}
\func{Re}\left[ \left\langle x,y\right\rangle -\left\langle x,e\right\rangle
\left\langle e,y\right\rangle \right] \geq -\frac{1}{16}\cdot \frac{1}{%
\lambda \left( 1-\lambda \right) }\left\vert \Gamma -\gamma \right\vert ^{2}.
\label{3.14.6}
\end{equation}%
Consequently, by (\ref{3.3.6}) and (\ref{3.14.6}) we deduce the desired
inequality (\ref{3.13.6}).
\end{proof}

\begin{remark}
\label{r1.6}If $M,m\in \mathbb{R}$ with $M>m$ and, for $\lambda \in \left(
0,1\right) ,$%
\begin{equation}
\left\Vert \lambda x+\left( 1-\lambda \right) y-\frac{M+m}{2}e\right\Vert
\leq \frac{1}{2}\left( M-m\right) ,  \label{3.15.6}
\end{equation}%
then 
\begin{equation*}
\left\langle x,y\right\rangle -\left\langle x,e\right\rangle \left\langle
e,y\right\rangle \leq \frac{1}{16}\cdot \frac{1}{\lambda \left( 1-\lambda
\right) }\left( M-m\right) ^{2}.
\end{equation*}%
If (\ref{3.15.6}) holds with \textquotedblleft $\pm $\textquotedblright\
instead of \hspace{0.05in}\textquotedblleft $+$\textquotedblright \hspace{%
0.05in}, then 
\begin{equation*}
\left\vert \left\langle x,y\right\rangle -\left\langle x,e\right\rangle
\left\langle e,y\right\rangle \right\vert \leq \frac{1}{16}\cdot \frac{1}{%
\lambda \left( 1-\lambda \right) }\left( M-m\right) ^{2}.
\end{equation*}
\end{remark}

\begin{remark}
\label{r2.6}If $\lambda =\frac{1}{2}$ in (\ref{3.1.6}), then we obtain the
result from \cite{SSD1.6}, i.e., 
\begin{equation*}
\func{Re}\left\langle \Gamma e-\frac{x+y}{2},\frac{x+y}{2}-\gamma
e\right\rangle \geq 0,
\end{equation*}%
or equivalently, 
\begin{equation*}
\left\Vert \frac{x+y}{2}-\frac{\gamma +\Gamma }{2}e\right\Vert \leq \frac{1}{%
2}\left\vert \Gamma -\gamma \right\vert ,
\end{equation*}%
implies 
\begin{equation}
\func{Re}\left[ \left\langle x,y\right\rangle -\left\langle x,e\right\rangle
\left\langle e,y\right\rangle \right] \leq \frac{1}{4}\left\vert \Gamma
-\gamma \right\vert ^{2}.  \label{3.19.6}
\end{equation}%
The constant $\frac{1}{4}$ is best possible in (\ref{3.19.6}).
\end{remark}

For $\lambda =\frac{1}{2},$ Corollary \ref{c3.2.6} and Remark \ref{r1.6}
will produce the corresponding results obtained in \cite{SSD1.6}. We omit
the details.

\subsection{Integral Inequalities}

The following proposition holds \cite{6NSSD}.

\begin{proposition}
\label{p5.1.6} If $f,g,h\in L^{2}\left( \Omega ,\mathbb{K}\right) $ and $%
\varphi ,\Phi ,\gamma ,\Gamma \in \mathbb{K}$, are such that $\func{Re}%
\left( \Phi \overline{\varphi }\right) >0,\func{Re}\left( \Gamma \overline{%
\gamma }\right) >0,$ $\int_{\Omega }\left\vert h\left( s\right) \right\vert
^{2}d\mu \left( s\right) =1$ and 
\begin{align}
\int_{\Omega }\func{Re}\left[ \left( \Phi h\left( s\right) -f\left( s\right)
\right) \left( \overline{f\left( s\right) }-\overline{\varphi }\overline{%
h\left( s\right) }\right) \right] d\mu \left( s\right) & \geq 0,
\label{5.1.6} \\
\int_{\Omega }\func{Re}\left[ \left( \Gamma h\left( s\right) -g\left(
s\right) \right) \left( \overline{g\left( s\right) }-\overline{\gamma }%
\overline{h\left( s\right) }\right) \right] d\mu \left( s\right) & \geq 0, 
\notag
\end{align}%
or equivalently, 
\begin{align*}
\left( \int_{\Omega }\left\vert f\left( s\right) -\frac{\Phi +\varphi }{2}%
h\left( s\right) \right\vert ^{2}d\mu \left( s\right) \right) ^{\frac{1}{2}%
}& \leq \frac{1}{2}\left\vert \Phi -\varphi \right\vert , \\
\left( \int_{\Omega }\left\vert g\left( s\right) -\frac{\Gamma +\gamma }{2}%
h\left( s\right) \right\vert ^{2}d\mu \left( s\right) \right) ^{\frac{1}{2}%
}& \leq \frac{1}{2}\left\vert \Gamma -\gamma \right\vert ,
\end{align*}%
then we have the following Gr\"{u}ss type integral inequality 
\begin{multline}
\left\vert \int_{\Omega }f\left( s\right) \overline{g\left( s\right) }d\mu
\left( s\right) -\int_{\Omega }f\left( s\right) \overline{h\left( s\right) }%
d\mu \left( s\right) \int_{\Omega }h\left( s\right) \overline{g\left(
s\right) }d\mu \left( s\right) \right\vert  \label{5.3.6} \\
\leq \frac{1}{4}M\left( \varphi ,\Phi \right) M\left( \gamma ,\Gamma \right)
\left\vert \int_{\Omega }f\left( s\right) \overline{h\left( s\right) }d\mu
\left( s\right) \int_{\Omega }h\left( s\right) \overline{g\left( s\right) }%
d\mu \left( s\right) \right\vert ,
\end{multline}%
where%
\begin{equation*}
M\left( \varphi ,\Phi \right) :=\left[ \frac{\left( \left\vert \Phi
\right\vert -\left\vert \varphi \right\vert \right) ^{2}+4\left[ \left\vert
\Phi \overline{\varphi }\right\vert -\func{Re}\left( \Phi \overline{\varphi }%
\right) \right] }{\func{Re}\left( \Phi \overline{\varphi }\right) }\right] ^{%
\frac{1}{2}}.
\end{equation*}%
The constant $\frac{1}{4}$ is best possible.
\end{proposition}

The proof follows by Theorem \ref{t3.1.6} on choosing $H=L^{2}\left( \Omega ,%
\mathbb{K}\right) $ with the inner product 
\begin{equation*}
\left\langle f,g\right\rangle :=\int_{\Omega }f\left( s\right) \overline{%
g\left( s\right) }d\mu \left( s\right) .
\end{equation*}%
We omit the details.

\begin{remark}
\label{r.1.6} It is obvious that a sufficient condition for $\left( \ref%
{5.1.6}\right) $ to hold is 
\begin{equation*}
\func{Re}\left[ \left( \Phi h\left( s\right) -f\left( s\right) \right)
\left( \overline{f\left( s\right) }-\overline{\varphi }\overline{h\left(
s\right) }\right) \right] \geq 0,
\end{equation*}%
and 
\begin{equation*}
\func{Re}\left[ \left( \Gamma h\left( s\right) -g\left( s\right) \right)
\left( \overline{g\left( s\right) }-\overline{\gamma }\overline{h\left(
s\right) }\right) \right] \geq 0,
\end{equation*}%
for $\mu -$a.e.$\;s\in \Omega ,$ or equivalently, 
\begin{equation*}
\left\vert f\left( s\right) -\frac{\Phi +\varphi }{2}h\left( s\right)
\right\vert \leq \frac{1}{2}\left\vert \Phi -\varphi \right\vert \left\vert
h\left( s\right) \right\vert
\end{equation*}%
\ and%
\begin{equation*}
\left\vert g\left( s\right) -\frac{\Gamma +\gamma }{2}h\left( s\right)
\right\vert \leq \frac{1}{2}\left\vert \Gamma -\gamma \right\vert \left\vert
h\left( s\right) \right\vert ,
\end{equation*}%
for $\mu -$a.e.$\;s\in \Omega .$
\end{remark}

The following result may be stated as well.

\begin{corollary}
\label{c5.1.6}If $z,Z,t,T\in \mathbb{K}$, $\mu \left( \Omega \right) <\infty 
$ and $f,g\in L^{2}\left( \Omega ,\mathbb{K}\right) $ are such that: 
\begin{align*}
\func{Re}\left[ \left( Z-f\left( s\right) \right) \left( \overline{f\left(
s\right) }-\bar{z}\right) \right] & \geq 0, \\
\func{Re}\left[ \left( T-g\left( s\right) \right) \left( \overline{g\left(
s\right) }-\bar{t}\right) \right] & \geq 0,\text{ \hspace{0.05in}for a.e. }%
s\in \Omega
\end{align*}%
or equivalently, 
\begin{align*}
\left\vert f\left( s\right) -\frac{z+Z}{2}\right\vert & \leq \frac{1}{2}%
\left\vert Z-z\right\vert , \\
\left\vert g\left( s\right) -\frac{t+T}{2}\right\vert & \leq \frac{1}{2}%
\left\vert T-t\right\vert ,\text{ \hspace{0.05in}for a.e. }s\in \Omega ;
\end{align*}%
then we have the inequality 
\begin{multline}
\left\vert \frac{1}{\mu \left( \Omega \right) }\int_{\Omega }f\left(
s\right) \overline{g\left( s\right) }d\mu \left( s\right) \right. -\left. 
\frac{1}{\mu \left( \Omega \right) }\int_{\Omega }f\left( s\right) d\mu
\left( s\right) \cdot \frac{1}{\mu \left( \Omega \right) }\int_{\Omega }%
\overline{g\left( s\right) }d\mu \left( s\right) \right\vert  \label{5.6.6}
\\
\leq \frac{1}{4}M\left( z,Z\right) M\left( t,T\right) \left\vert \frac{1}{%
\mu \left( \Omega \right) }\int_{\Omega }f\left( s\right) d\mu \left(
s\right) \cdot \frac{1}{\mu \left( \Omega \right) }\int_{\Omega }\overline{%
g\left( s\right) }d\mu \left( s\right) \right\vert .
\end{multline}
\end{corollary}

\begin{remark}
\label{r.3.6} The case of real functions incorporates the following
interesting inequality 
\begin{equation*}
\left\vert \frac{\mu \left( \Omega \right) \int_{\Omega }f\left( s\right)
g\left( s\right) d\mu \left( s\right) }{\int_{\Omega }f\left( s\right) d\mu
\left( s\right) \int_{\Omega }g\left( s\right) d\mu \left( s\right) }%
-1\right\vert \leq \frac{1}{4}\cdot \frac{\left( Z-z\right) \left(
T-t\right) }{\sqrt{ztZT}},
\end{equation*}%
provided $\mu \left( \Omega \right) <\infty ,$%
\begin{equation*}
z\leq f\left( s\right) \leq Z,\ \ t\leq g\left( s\right) \leq T
\end{equation*}%
for $\mu -$a.e. $s\in \Omega ,$ where $z,t,Z,T$ are real numbers and the
integrals at the denominator are not zero. Here the constant $\frac{1}{4}$
is best possible in the sense mentioned above.
\end{remark}

Using Theorem \ref{t3.1.6} we may state the following result as well \cite%
{6NSSD}.

\begin{proposition}
\label{p5.3.6}If $f,g,h\in L^{2}\left( \Omega ,\mathbb{K}\right) $ and $%
\gamma ,\Gamma \in \mathbb{K}$ are such that $\int_{\Omega }\left\vert
h\left( s\right) \right\vert ^{2}d\mu \left( s\right) =1$ and%
\begin{multline}
\int_{\Omega }\left\{ \func{Re}\left[ \Gamma h\left( s\right) -\left(
\lambda f\left( s\right) +\left( 1-\lambda \right) g\left( s\right) \right) %
\right] \right.  \label{5.7.6} \\
\times \left. \left[ \lambda \overline{f\left( s\right) }+\left( 1-\lambda
\right) \overline{g\left( s\right) }-\bar{\gamma}\bar{h}\left( s\right) %
\right] \right\} d\mu \left( s\right) \geq 0,
\end{multline}%
or equivalently, 
\begin{equation}
\left( \int_{\Omega }\left\vert \lambda f\left( s\right) +\left( 1-\lambda
\right) g\left( s\right) -\frac{\gamma +\Gamma }{2}h\left( s\right)
\right\vert ^{2}d\mu \left( s\right) \right) ^{\frac{1}{2}}\leq \frac{1}{2}%
\left\vert \Gamma -\gamma \right\vert ,  \label{5.8.6}
\end{equation}%
then we have the inequality 
\begin{align*}
I& :=\int_{\Omega }\func{Re}\left[ f\left( s\right) \overline{g\left(
s\right) }\right] d\mu \left( s\right) \\
& \qquad \qquad \qquad -\func{Re}\left[ \int_{\Omega }f\left( s\right) 
\overline{h\left( s\right) }d\mu \left( s\right) \cdot \int_{\Omega }h\left(
s\right) \overline{g\left( s\right) }d\mu \left( s\right) \right] \\
& \leq \frac{1}{16}\cdot \frac{1}{\lambda \left( 1-\lambda \right) }%
\left\vert \Gamma -\gamma \right\vert ^{2}.
\end{align*}%
The constant $\frac{1}{16}$ is best possible.

If (\ref{5.7.6}) and (\ref{5.8.6}) hold with \textquotedblleft\ $\pm $
\textquotedblright\ instead of \textquotedblleft\ $+$ \textquotedblright\
(see Corollary \ref{c3.2.6}), then 
\begin{equation*}
\left\vert I\right\vert \leq \frac{1}{16}\cdot \frac{1}{\lambda \left(
1-\lambda \right) }\left\vert \Gamma -\gamma \right\vert ^{2}.
\end{equation*}
\end{proposition}

\begin{remark}
It is obvious that a sufficient condition for (\ref{5.7.6}) to hold is 
\begin{equation*}
\func{Re}\left\{ \left[ \Gamma h\left( s\right) -\left( \lambda f\left(
s\right) +\left( 1-\lambda \right) g\left( s\right) \right) \right] \cdot %
\left[ \lambda \overline{f\left( s\right) }+\left( 1-\lambda \right) 
\overline{g\left( s\right) }-\bar{\gamma}\bar{h}\left( s\right) \right]
\right\} \geq 0
\end{equation*}%
for a.e. $s\in \Omega ,$ or equivalently, 
\begin{equation*}
\left\vert \lambda f\left( s\right) +\left( 1-\lambda \right) g\left(
s\right) -\frac{\gamma +\Gamma }{2}h\left( s\right) \right\vert \leq \frac{1%
}{2}\left\vert \Gamma -\gamma \right\vert \left\vert h\left( s\right)
\right\vert ,\text{ }
\end{equation*}%
for a.e. $s\in \Omega .$
\end{remark}

Finally, the following corollary holds.

\begin{corollary}
\label{c5.4.6}If $Z,z\in \mathbb{K}$, $\mu \left( \Omega \right) <\infty $
and $f,g\in L^{2}\left( \Omega ,\mathbb{K}\right) $ are such that 
\begin{equation}
\func{Re}\left[ \left( Z-\left( \lambda f\left( s\right) +\left( 1-\lambda
\right) g\left( s\right) \right) \right) \left( \lambda \overline{f\left(
s\right) }+\left( 1-\lambda \right) \overline{g\left( s\right) }-\overline{z}%
\right) \right] \geq 0  \label{5.13.6}
\end{equation}%
for a.e. $s\in \Omega ,$ or, equivalently 
\begin{equation}
\left\vert \lambda f\left( s\right) +\left( 1-\lambda \right) g\left(
s\right) -\frac{z+Z}{2}\right\vert \leq \frac{1}{2}\left\vert Z-z\right\vert 
\text{,\hspace{0.05in}}  \label{5.14.6}
\end{equation}%
for a.e. $s\in \Omega ,$ then we have the inequality 
\begin{align*}
J& :=\frac{1}{\mu \left( \Omega \right) }\int_{\Omega }\func{Re}\left[
f\left( s\right) \overline{g\left( s\right) }\right] d\mu \left( s\right) \\
& \qquad \qquad \qquad -\func{Re}\left[ \frac{1}{\mu \left( \Omega \right) }%
\int_{\Omega }f\left( s\right) d\mu \left( s\right) \cdot \frac{1}{\mu
\left( \Omega \right) }\int_{\Omega }\overline{g\left( s\right) }d\mu \left(
s\right) \right] \\
& \leq \frac{1}{16}\cdot \frac{1}{\lambda \left( 1-\lambda \right) }%
\left\vert Z-z\right\vert ^{2}.
\end{align*}%
If (\ref{5.13.6}) and (\ref{5.14.6}) hold with \textquotedblleft\ $\pm $
\textquotedblright\ instead of \textquotedblleft\ $+$ \textquotedblright $,$
then 
\begin{equation*}
\left\vert J\right\vert \leq \frac{1}{16}\cdot \frac{1}{\lambda \left(
1-\lambda \right) }\left\vert Z-z\right\vert ^{2}.
\end{equation*}
\end{corollary}

\begin{remark}
\label{r.4.6} It is obvious that if one chooses the discrete measure above,
then all the inequalities in this section may be written for sequences of
real or complex numbers. We omit the details.
\end{remark}

\newpage

\section{Other Gr\"{u}ss Type Inequalities}

\subsection{General Results}

We may state the following result \cite{7NSSD}.

\begin{theorem}
\label{t4.1.7a}Let $\left( H;\left\langle \cdot ,\cdot \right\rangle \right) 
$ be an inner product space over the real or complex number field $\mathbb{K}
$ $\left( \mathbb{K}=\mathbb{R},\mathbb{K}=\mathbb{C}\right) $ and $x,y,e\in
H$ with $\left\Vert e\right\Vert =1.$ If $r_{1},r_{2}\in \left( 0,1\right) $
and 
\begin{equation*}
\left\Vert x-e\right\Vert \leq r_{1},\ \ \ \ \left\Vert y-e\right\Vert \leq
r_{2},
\end{equation*}%
then we have the inequality%
\begin{equation}
\left\vert \left\langle x,y\right\rangle -\left\langle x,e\right\rangle
\left\langle e,y\right\rangle \right\vert \leq r_{1}r_{2}\left\Vert
x\right\Vert \left\Vert y\right\Vert .  \label{4.2.7a}
\end{equation}%
The inequality (\ref{4.2.7a}) is sharp in the sense that the constant $c=1$
in front of $r_{1}r_{2}$ cannot be replaced by a smaller constant.
\end{theorem}

\begin{proof}
Start with the inequality%
\begin{equation}
\left\vert \left\langle x,y\right\rangle -\left\langle x,e\right\rangle
\left\langle e,y\right\rangle \right\vert ^{2}\leq \left( \left\Vert
x\right\Vert ^{2}-\left\vert \left\langle x,e\right\rangle \right\vert
^{2}\right) \left( \left\Vert y\right\Vert ^{2}-\left\vert \left\langle
y,e\right\rangle \right\vert ^{2}\right) .  \label{4.3.7a}
\end{equation}%
Using Theorem \ref{t2.1.3} for $a=e,$ we may state that%
\begin{equation}
\left\Vert x\right\Vert ^{2}-\left\vert \left\langle x,e\right\rangle
\right\vert ^{2}\leq r_{1}^{2}\left\Vert x\right\Vert ^{2},\ \ \ \ \ \
\left\Vert y\right\Vert ^{2}-\left\vert \left\langle y,e\right\rangle
\right\vert ^{2}\leq r_{2}^{2}\left\Vert y\right\Vert ^{2}.  \label{4.4.7a}
\end{equation}%
Utilizing (\ref{4.3.7a}) and (\ref{4.4.7a}), we deduce%
\begin{equation*}
\left\vert \left\langle x,y\right\rangle -\left\langle x,e\right\rangle
\left\langle e,y\right\rangle \right\vert ^{2}\leq
r_{1}^{2}r_{2}^{2}\left\Vert x\right\Vert ^{2}\left\Vert y\right\Vert ^{2},
\end{equation*}%
which is clearly equivalent to the desired inequality (\ref{4.2.7a}).

The sharpness of the constant follows by the fact that for $x=y,$ $%
r_{1}=r_{2}=r,$ we get from (\ref{4.2.7a}) that%
\begin{equation}
\left\Vert x\right\Vert ^{2}-\left\vert \left\langle x,e\right\rangle
\right\vert ^{2}\leq r^{2}\left\Vert x\right\Vert ^{2},  \label{4.6.7a}
\end{equation}%
provided $\left\Vert e\right\Vert =1$ and $\left\Vert x-e\right\Vert \leq
r<1.$ The inequality (\ref{4.6.7a}) is sharp, as shown in Theorem \ref%
{t2.1.3}, and the proof is completed.
\end{proof}

Another companion of the Gr\"{u}ss inequality may be stated as well \cite%
{7NSSD}.

\begin{theorem}
\label{t4.2.7a}Let $\left( H;\left\langle \cdot ,\cdot \right\rangle \right) 
$ be an inner product space over $\mathbb{K}$ and $x,y,e\in H$ with $%
\left\Vert e\right\Vert =1.$ Suppose also that $a,A,b,B\in \mathbb{K}$ $%
\left( \mathbb{K}=\mathbb{R},\mathbb{C}\right) $ such that $\func{Re}\left( A%
\overline{a}\right) ,$ $\func{Re}\left( B\overline{b}\right) >0.$ If either%
\begin{equation*}
\func{Re}\left\langle Ae-x,x-ae\right\rangle \geq 0,\ \ \func{Re}%
\left\langle Be-y,y-be\right\rangle \geq 0,\ 
\end{equation*}%
or equivalently,%
\begin{equation*}
\left\Vert x-\frac{a+A}{2}e\right\Vert \leq \frac{1}{2}\left\vert
A-a\right\vert ,\ \ \left\Vert y-\frac{b+B}{2}e\right\Vert \leq \frac{1}{2}%
\left\vert B-b\right\vert ,
\end{equation*}%
holds, then we have the inequality%
\begin{equation}
\left\vert \left\langle x,y\right\rangle -\left\langle x,e\right\rangle
\left\langle e,y\right\rangle \right\vert \leq \frac{1}{4}\cdot \frac{%
\left\vert A-a\right\vert \left\vert B-b\right\vert }{\sqrt{\func{Re}\left( A%
\overline{a}\right) \func{Re}\left( B\overline{b}\right) }}\left\vert
\left\langle x,e\right\rangle \left\langle e,y\right\rangle \right\vert .
\label{4.9.7a}
\end{equation}%
The constant $\frac{1}{4}$ is best possible.
\end{theorem}

\begin{proof}
We know, that%
\begin{equation}
\left\vert \left\langle x,y\right\rangle -\left\langle x,e\right\rangle
\left\langle e,y\right\rangle \right\vert ^{2}\leq \left( \left\Vert
x\right\Vert ^{2}-\left\vert \left\langle x,e\right\rangle \right\vert
^{2}\right) \left( \left\Vert y\right\Vert ^{2}-\left\vert \left\langle
y,e\right\rangle \right\vert ^{2}\right) .  \label{4.10.7a}
\end{equation}%
If we use Corollary \ref{c2.3.3}, then we may state that%
\begin{equation}
\left\Vert x\right\Vert ^{2}-\left\vert \left\langle x,e\right\rangle
\right\vert ^{2}\leq \frac{1}{4}\cdot \frac{\left\vert A-a\right\vert ^{2}}{%
\func{Re}\left( A\overline{a}\right) }\left\vert \left\langle
x,e\right\rangle \right\vert ^{2}  \label{4.11.7a}
\end{equation}%
and%
\begin{equation}
\left\Vert y\right\Vert ^{2}-\left\vert \left\langle y,e\right\rangle
\right\vert ^{2}\leq \frac{1}{4}\cdot \frac{\left\vert B-b\right\vert ^{2}}{%
\func{Re}\left( B\overline{b}\right) }\left\vert \left\langle
y,e\right\rangle \right\vert ^{2}.  \label{4.12.7a}
\end{equation}%
Utilizing (\ref{4.10.7a}) -- (\ref{4.12.7a}), we deduce%
\begin{equation*}
\left\vert \left\langle x,y\right\rangle -\left\langle x,e\right\rangle
\left\langle e,y\right\rangle \right\vert ^{2}\leq \frac{1}{16}\cdot \frac{%
\left\vert A-a\right\vert ^{2}\left\vert B-b\right\vert ^{2}}{\func{Re}%
\left( A\overline{a}\right) \func{Re}\left( B\overline{b}\right) }\left\vert
\left\langle x,e\right\rangle \left\langle e,y\right\rangle \right\vert ^{2},
\end{equation*}%
which is clearly equivalent to the desired inequality (\ref{4.9.7a}).

The sharpness of the constant follows from Corollary \ref{c2.3.3}, and we
omit the details.
\end{proof}

\begin{remark}
With the assumptions of Theorem \ref{t4.2.7a} and if $\left\langle
x,e\right\rangle ,\left\langle y,e\right\rangle \neq 0$ (that is actually
the interesting case), then one has the inequality%
\begin{equation*}
\left\vert \frac{\left\langle x,y\right\rangle }{\left\langle
x,e\right\rangle \left\langle e,y\right\rangle }-1\right\vert \leq \frac{1}{4%
}\cdot \frac{\left\vert A-a\right\vert \left\vert B-b\right\vert }{\sqrt{%
\func{Re}\left( A\overline{a}\right) \func{Re}\left( B\overline{b}\right) }}.
\end{equation*}%
The constant $\frac{1}{4}$ is best possible.
\end{remark}

We may state the following result \cite{8NSSD}.

\begin{theorem}
\label{t4.1.7b}Let $\left( H;\left\langle \cdot ,\cdot \right\rangle \right) 
$ be an inner product space over the real or complex number field $\mathbb{K}
$ and $x,y,e\in H$ with $\left\Vert e\right\Vert =1.$ If $r_{1},r_{2}>0$ and 
\begin{equation*}
\left\Vert x-e\right\Vert \leq r_{1},\ \ \ \ \left\Vert y-e\right\Vert \leq
r_{2},
\end{equation*}%
then we have the inequalities%
\begin{align}
\left\vert \left\langle x,y\right\rangle -\left\langle x,e\right\rangle
\left\langle e,y\right\rangle \right\vert & \leq \frac{1}{2}r_{1}r_{2}\sqrt{%
\left\Vert x\right\Vert +\left\vert \left\langle x,e\right\rangle
\right\vert }\cdot \sqrt{\left\Vert y\right\Vert +\left\vert \left\langle
y,e\right\rangle \right\vert }  \label{4.2.7b} \\
& \leq r_{1}r_{2}\left\Vert x\right\Vert \left\Vert y\right\Vert .  \notag
\end{align}%
The constant $\frac{1}{2}$ is best possible in the sense that it cannot be
replaced by a smaller constant.
\end{theorem}

\begin{proof}
Start with the inequality%
\begin{equation}
\left\vert \left\langle x,y\right\rangle -\left\langle x,e\right\rangle
\left\langle e,y\right\rangle \right\vert ^{2}\leq \left( \left\Vert
x\right\Vert ^{2}-\left\vert \left\langle x,e\right\rangle \right\vert
^{2}\right) \left( \left\Vert y\right\Vert ^{2}-\left\vert \left\langle
y,e\right\rangle \right\vert ^{2}\right) .  \label{4.3.7b}
\end{equation}%
Using Theorem \ref{t2.1.4} for $a=e,$ we have%
\begin{align}
0& \leq \left\Vert x\right\Vert ^{2}-\left\vert \left\langle
x,e\right\rangle \right\vert ^{2}  \label{4.4.7b} \\
& =\left( \left\Vert x\right\Vert -\left\vert \left\langle x,e\right\rangle
\right\vert \right) \left( \left\Vert x\right\Vert +\left\vert \left\langle
x,e\right\rangle \right\vert \right)  \notag \\
& \leq \frac{1}{2}r_{1}^{2}\left( \left\Vert x\right\Vert +\left\vert
\left\langle x,e\right\rangle \right\vert \right) \leq r_{1}^{2}\left\Vert
x\right\Vert ,  \notag
\end{align}%
and, in a similar way%
\begin{align}
0& \leq \left\Vert y\right\Vert ^{2}-\left\vert \left\langle
y,e\right\rangle \right\vert ^{2}  \label{4.5.7b} \\
& \leq \frac{1}{2}r_{2}^{2}\left( \left\Vert y\right\Vert +\left\vert
\left\langle y,e\right\rangle \right\vert \right) \leq r_{2}^{2}\left\Vert
y\right\Vert .  \notag
\end{align}%
Utilising (\ref{4.3.7b}) -- (\ref{4.5.7b}), we may state that%
\begin{align}
\left\vert \left\langle x,y\right\rangle -\left\langle x,e\right\rangle
\left\langle e,y\right\rangle \right\vert ^{2}& \leq \frac{1}{4}%
r_{1}^{2}r_{2}^{2}\left( \left\Vert x\right\Vert +\left\vert \left\langle
x,e\right\rangle \right\vert \right) \left( \left\Vert y\right\Vert
+\left\vert \left\langle y,e\right\rangle \right\vert \right)  \label{4.6.7b}
\\
& \leq r_{1}^{2}r_{2}^{2}\left\Vert x\right\Vert \left\Vert y\right\Vert , 
\notag
\end{align}%
giving the desired inequality (\ref{4.2.7b}).

To prove the sharpness of the constant $\frac{1}{2}$, let us assume that $%
x=y $ in (\ref{4.2.7b}), to get%
\begin{equation}
\left\Vert x\right\Vert ^{2}-\left\vert \left\langle x,e\right\rangle
\right\vert ^{2}\leq \frac{1}{2}r_{1}^{2}\left( \left\Vert x\right\Vert
+\left\vert \left\langle x,e\right\rangle \right\vert \right) ,
\label{4.7.7b}
\end{equation}%
provided $\left\Vert x-e\right\Vert \leq r_{1}.$ If $x\neq 0,$ then dividing
(\ref{4.7.7b}) with $\left\Vert x\right\Vert +\left\vert \left\langle
x,e\right\rangle \right\vert >0$ we get%
\begin{equation}
\left\Vert x\right\Vert -\left\vert \left\langle x,e\right\rangle
\right\vert \leq \frac{1}{2}r_{1}^{2}  \label{4.8.7b}
\end{equation}%
provided $\left\Vert x-e\right\Vert \leq r_{1},$ $\left\Vert e\right\Vert
=1. $ However, (\ref{4.8.7b}) is in fact (\ref{2.2.4}) for $a=e,$ for which
we have shown that $\frac{1}{2}$ is the best possible constant.
\end{proof}

The following result also holds \cite{8NSSD}.

\begin{theorem}
\label{t4.2.7b}With the assumptions of Theorem \ref{t4.1.7b}, we have the
inequality%
\begin{equation}
\left\vert \left\langle x,y\right\rangle -\left\langle x,e\right\rangle
\left\langle e,y\right\rangle \right\vert \leq r_{1}r_{2}\sqrt{\frac{1}{4}%
r_{1}^{2}+\left\vert \left\langle x,e\right\rangle \right\vert }\cdot \sqrt{%
\frac{1}{4}r_{2}^{2}+\left\vert \left\langle y,e\right\rangle \right\vert }.
\label{4.9.7b}
\end{equation}
\end{theorem}

\begin{proof}
Note that, from Theorem \ref{t2.2.4}, we have%
\begin{equation}
\left\Vert x\right\Vert \left\Vert a\right\Vert \leq \left\vert \left\langle
x,a\right\rangle \right\vert +\frac{1}{2}r^{2},  \label{4.10.7b}
\end{equation}%
provided $\left\Vert x-a\right\Vert \leq r.$

Taking the square of (\ref{4.10.7b}) and re-arranging the terms, we obtain:%
\begin{equation*}
0\leq \left\Vert x\right\Vert ^{2}\left\Vert a\right\Vert ^{2}-\left\vert
\left\langle x,a\right\rangle \right\vert ^{2}\leq r^{2}\left( \frac{1}{4}%
r^{2}+\left\vert \left\langle x,a\right\rangle \right\vert \right) ,
\end{equation*}%
provided $\left\Vert x-a\right\Vert \leq r.$

Using the assumption of the theorem, we then have%
\begin{equation}
0\leq \left\Vert x\right\Vert ^{2}-\left\vert \left\langle x,e\right\rangle
\right\vert ^{2}\leq r_{1}^{2}\left( \frac{1}{4}r_{1}^{2}+\left\vert
\left\langle x,e\right\rangle \right\vert \right) ,  \label{4.12.7b}
\end{equation}%
and 
\begin{equation}
0\leq \left\Vert y\right\Vert ^{2}-\left\vert \left\langle y,e\right\rangle
\right\vert ^{2}\leq r_{2}^{2}\left( \frac{1}{4}r_{2}^{2}+\left\vert
\left\langle y,e\right\rangle \right\vert \right) .  \label{4.13.7b}
\end{equation}%
Utilising (\ref{4.3.7b}), (\ref{4.12.7b}) and (\ref{4.13.7b}), we deduce the
desired inequality (\ref{4.9.7b}).
\end{proof}

The following result may be stated as well \cite{8NSSD}.

\begin{theorem}
\label{t4.3.7b}Let $\left( H;\left\langle \cdot ,\cdot \right\rangle \right) 
$ be an inner product space over $\mathbb{K}$ and $x,y,e\in H$ with $%
\left\Vert e\right\Vert =1.$ Suppose also that $a,A,b,B\in \mathbb{K}$ $%
\left( \mathbb{K}=\mathbb{C},\mathbb{R}\right) $ such that $A\neq -a,B\neq
-b.$ If either%
\begin{equation*}
\func{Re}\left\langle Ae-x,x-ae\right\rangle \geq 0,\ \ \func{Re}%
\left\langle Be-y,y-be\right\rangle \geq 0,
\end{equation*}%
or equivalently,%
\begin{equation*}
\left\Vert x-\frac{a+A}{2}e\right\Vert \leq \frac{1}{2}\left\vert
A-a\right\vert ,\ \ \ \left\Vert y-\frac{b+B}{2}e\right\Vert \leq \frac{1}{2}%
\left\vert B-b\right\vert ,
\end{equation*}%
holds, then we have the inequality%
\begin{align}
& \left\vert \left\langle x,y\right\rangle -\left\langle x,e\right\rangle
\left\langle e,y\right\rangle \right\vert  \label{4.16.7b} \\
& \leq \frac{1}{4}\cdot \frac{\left\vert A-a\right\vert \left\vert
B-b\right\vert }{\sqrt{\left\vert A+a\right\vert \left\vert B+b\right\vert }}%
\sqrt{\left\Vert x\right\Vert +\left\vert \left\langle x,e\right\rangle
\right\vert }\cdot \sqrt{\left\Vert y\right\Vert +\left\vert \left\langle
y,e\right\rangle \right\vert }  \notag \\
& \leq \frac{1}{2}\cdot \frac{\left\vert A-a\right\vert \left\vert
B-b\right\vert }{\sqrt{\left\vert A+a\right\vert \left\vert B+b\right\vert }}%
\sqrt{\left\Vert x\right\Vert \left\Vert y\right\Vert }.  \notag
\end{align}%
The constant $\frac{1}{4}$ is best possible in (\ref{4.16.7b}).
\end{theorem}

\begin{proof}
From Theorem \ref{t2.2.4}, we may state that%
\begin{align}
0& \leq \left\Vert x\right\Vert ^{2}-\left\vert \left\langle
x,e\right\rangle \right\vert ^{2}  \label{4.17.7b} \\
& =\left( \left\Vert x\right\Vert -\left\vert \left\langle x,e\right\rangle
\right\vert \right) \left( \left\Vert x\right\Vert +\left\vert \left\langle
x,e\right\rangle \right\vert \right)  \notag \\
& \leq \frac{1}{4}\cdot \frac{\left\vert A-a\right\vert ^{2}}{\left\vert
A+a\right\vert }\left( \left\Vert x\right\Vert +\left\vert \left\langle
x,e\right\rangle \right\vert \right) ,  \notag
\end{align}%
and%
\begin{equation}
0\leq \left\Vert y\right\Vert ^{2}-\left\vert \left\langle y,e\right\rangle
\right\vert ^{2}\leq \frac{1}{4}\cdot \frac{\left\vert B-b\right\vert ^{2}}{%
\left\vert B+b\right\vert }\left( \left\Vert y\right\Vert +\left\vert
\left\langle y,e\right\rangle \right\vert \right) .  \label{4.18.7b}
\end{equation}%
Making use of (\ref{4.3.7b}) and (\ref{4.17.7b}), (\ref{4.18.7b}), we deduce
the first inequality in (\ref{4.16.7b}).

The best constant follows by the use of Theorem \ref{t2.2.4}, and we omit
the details.
\end{proof}

Finally, we may state the following theorem as well \cite{8NSSD}.

\begin{theorem}
\label{t4.4.7b}With the assumptions of Theorem \ref{t4.3.7b}, we have the
inequality%
\begin{multline}
\left\vert \left\langle x,y\right\rangle -\left\langle x,e\right\rangle
\left\langle e,y\right\rangle \right\vert  \label{4.19.7b} \\
\leq \frac{1}{2}\cdot \frac{\left\vert A-a\right\vert \left\vert
B-b\right\vert }{\sqrt{\left\vert A+a\right\vert \left\vert B+b\right\vert }}%
\sqrt{\frac{1}{8}\cdot \frac{\left\vert A-a\right\vert ^{2}}{\left\vert
A+a\right\vert }+\left\vert \left\langle x,e\right\rangle \right\vert } \\
\times \sqrt{\frac{1}{8}\cdot \frac{\left\vert B-b\right\vert ^{2}}{%
\left\vert B+b\right\vert }+\left\vert \left\langle y,e\right\rangle
\right\vert }.
\end{multline}
\end{theorem}

\begin{proof}
Using Theorem \ref{t2.2.4}, we may state that%
\begin{equation*}
0\leq \left\Vert x\right\Vert -\left\vert \left\langle x,e\right\rangle
\right\vert \leq \frac{1}{4}\cdot \frac{\left\vert A-a\right\vert ^{2}}{%
\left\vert A+a\right\vert }.
\end{equation*}%
This inequality implies that%
\begin{equation*}
\left\Vert x\right\Vert ^{2}\leq \left\vert \left\langle x,e\right\rangle
\right\vert ^{2}+\frac{1}{2}\left\vert \left\langle x,e\right\rangle
\right\vert \cdot \frac{\left\vert A-a\right\vert ^{2}}{\left\vert
A+a\right\vert }+\frac{1}{16}\cdot \frac{\left\vert A-a\right\vert ^{4}}{%
\left\vert A+a\right\vert ^{2}},
\end{equation*}%
giving%
\begin{equation}
0\leq \left\Vert x\right\Vert ^{2}-\left\vert \left\langle x,e\right\rangle
\right\vert ^{2}\leq \frac{1}{2}\cdot \frac{\left\vert A-a\right\vert ^{2}}{%
\left\vert A+a\right\vert }\left[ \left\vert \left\langle x,e\right\rangle
\right\vert +\frac{1}{8}\cdot \frac{\left\vert A-a\right\vert ^{2}}{%
\left\vert A+a\right\vert }\right] .  \label{4.20.7b}
\end{equation}%
Similarly, we have%
\begin{equation}
0\leq \left\Vert y\right\Vert ^{2}-\left\vert \left\langle y,e\right\rangle
\right\vert ^{2}\leq \frac{1}{2}\cdot \frac{\left\vert B-b\right\vert ^{2}}{%
\left\vert B+b\right\vert }\left[ \left\vert \left\langle y,e\right\rangle
\right\vert +\frac{1}{8}\cdot \frac{\left\vert B-b\right\vert ^{2}}{%
\left\vert B+b\right\vert }\right] .  \label{4.21.7b}
\end{equation}%
By making use of (\ref{4.3.7b}) and (\ref{4.20.7b}), (\ref{4.21.7b}), we
deduce the desired inequality (\ref{4.19.7b}).
\end{proof}

For some recent results on Gr\"{u}ss type inequalities in inner product
spaces, see \cite{SSD0.7}, \cite{SSD00.7} and \cite{PFR.7}.

\subsection{Integral Inequalities}

The following Gr\"{u}ss type integral inequality for real or complex-valued
functions also holds \cite{7NSSD}.

\begin{proposition}
\label{p.7.3.7a}Let $f,g,h\in L_{\rho }^{2}\left( \Omega ,\mathbb{K}\right) $
with $\int_{\Omega }\rho \left( s\right) \left\vert h\left( s\right)
\right\vert ^{2}d\mu \left( s\right) =1$ and $a,A,b,B\in \mathbb{K}$ such
that $\func{Re}\left( A\overline{a}\right) ,\func{Re}\left( B\overline{b}%
\right) >0$ and%
\begin{eqnarray*}
\func{Re}\left[ \left( Ah\left( s\right) -f\left( s\right) \right) \left( 
\overline{f\left( s\right) }-\overline{a}\overline{h\left( s\right) }\right) %
\right] &\geq &0, \\
\func{Re}\left[ \left( Bh\left( s\right) -g\left( s\right) \right) \left( 
\overline{g\left( s\right) }-\overline{b}\overline{h\left( s\right) }\right) %
\right] &\geq &0\text{,}
\end{eqnarray*}%
for $\mu -$a.e. $s\in \Omega .$ Then we have the inequalities%
\begin{multline*}
\left\vert \int_{\Omega }\rho \left( s\right) f\left( s\right) \overline{%
g\left( s\right) }d\mu \left( s\right) -\int_{\Omega }\rho \left( s\right)
f\left( s\right) \overline{h\left( s\right) }d\mu \left( s\right)
\int_{\Omega }\rho \left( s\right) h\left( s\right) \overline{g\left(
s\right) }d\mu \left( s\right) \right\vert \\
\leq \frac{1}{4}\cdot \frac{\left\vert A-a\right\vert \left\vert
B-b\right\vert }{\sqrt{\func{Re}\left( A\overline{a}\right) \func{Re}\left( B%
\overline{b}\right) }}\left\vert \int_{\Omega }\rho \left( s\right) f\left(
s\right) \overline{h\left( s\right) }d\mu \left( s\right) \int_{\Omega }\rho
\left( s\right) h\left( s\right) \overline{g\left( s\right) }d\mu \left(
s\right) \right\vert .
\end{multline*}%
The constant $\frac{1}{4}$ is best possible.
\end{proposition}

The proof follows by Theorem \ref{t4.2.7a}.

By making use of Theorem \ref{t4.3.7b}, we may also state

\begin{proposition}
\label{p7.3.7b}Let $f,g,h\in L_{\rho }^{2}\left( \Omega ,\mathbb{K}\right) $
be such that $\int_{\Omega }\rho \left( s\right) \left\vert h\left( s\right)
\right\vert ^{2}d\mu \left( s\right) =1.$ Suppose also that $a,A,b,B\in 
\mathbb{K}$ with $A\neq -a,B\neq -b$ and%
\begin{eqnarray*}
\func{Re}\left[ \left( Ah\left( s\right) -f\left( s\right) \right) \left( 
\overline{f\left( s\right) }-\overline{a}\overline{h\left( s\right) }\right) %
\right] &\geq &0\text{,} \\
\func{Re}\left[ \left( Bh\left( s\right) -g\left( s\right) \right) \left( 
\overline{g\left( s\right) }-\overline{b}\overline{h\left( s\right) }\right) %
\right] &\geq &0,\text{ \ \ }
\end{eqnarray*}%
for\ $\mu -$a.e.\ $s\in \Omega .$ Then we have the inequality%
\begin{multline*}
\left\vert \int_{\Omega }\rho \left( s\right) f\left( s\right) \overline{%
g\left( s\right) }d\mu \left( s\right) -\int_{\Omega }\rho \left( s\right)
f\left( s\right) \overline{h\left( s\right) }d\mu \left( s\right)
\int_{\Omega }\rho \left( s\right) h\left( s\right) \overline{g\left(
s\right) }d\mu \left( s\right) \right\vert \\
\leq \frac{1}{4}\cdot \frac{\left\vert A-a\right\vert \left\vert
B-b\right\vert }{\sqrt{\left\vert A+a\right\vert \left\vert B+b\right\vert }}%
\qquad \qquad \qquad \qquad \qquad \qquad \qquad \\
\times \sqrt{\left( \int_{\Omega }\rho \left( s\right) \left\vert f\left(
s\right) \right\vert ^{2}d\mu \left( s\right) \right) ^{\frac{1}{2}%
}+\left\vert \int_{\Omega }\rho \left( s\right) f\left( s\right) \overline{%
h\left( s\right) }d\mu \left( s\right) \right\vert } \\
\times \sqrt{\left( \int_{\Omega }\rho \left( s\right) \left\vert g\left(
s\right) \right\vert ^{2}d\mu \left( s\right) \right) ^{\frac{1}{2}%
}+\left\vert \int_{\Omega }\rho \left( s\right) g\left( s\right) \overline{%
h\left( s\right) }d\mu \left( s\right) \right\vert }.
\end{multline*}%
The constant $\frac{1}{4}$ is best possible.
\end{proposition}

\newpage

\chapter{Reverses of Bessel's Inequality}

\section{Introduction}

Let $\left( H,\left\langle \cdot ,\cdot \right\rangle \right) $ be an inner
product space over $\mathbb{K}$ $\left( \mathbb{K}=\mathbb{R},\mathbb{C}%
\right) $ and $\left\{ e_{i}\right\} _{i\in I}$ a finite or infinite family
of \textit{orthonormal vectors} in $H$, i.e.,%
\begin{equation*}
\left\langle e_{i},e_{j}\right\rangle =\left\{ 
\begin{array}{ll}
0 & \text{if }i\neq j \\ 
&  \\ 
1 & \text{if }i=j%
\end{array}%
\right. ;\ \ i,j\in I.
\end{equation*}%
It is well known that, the following inequality due to Bessel, holds 
\begin{equation*}
\sum_{i\in I}\left\vert \left\langle x,e_{i}\right\rangle \right\vert
^{2}\leq \left\Vert x\right\Vert ^{2}
\end{equation*}%
\ for any $x\in H,$ where the meaning of the sum is:%
\begin{equation*}
\sum_{i\in I}\left\vert \left\langle x,e_{i}\right\rangle \right\vert
^{2}:=\sup_{F\subset I}\left\{ \sum_{i\in F}\left\vert \left\langle
x,e_{i}\right\rangle \right\vert ^{2},\ \ \ F\text{ is a finite part of }%
I\right\} .
\end{equation*}%
If $\left( H,\left\langle \cdot ,\cdot \right\rangle \right) $ is an
infinite dimensional Hilbert space and $\left\{ e_{i}\right\} _{i\in \mathbb{%
N}}$ an orthonormal family in $H,$ then we also have%
\begin{equation*}
\sum_{i=0}^{\infty }\left\vert \left\langle x,e_{i}\right\rangle \right\vert
^{2}\leq \left\Vert x\right\Vert ^{2}
\end{equation*}%
for any $x\in H.$ Here the meaning of the series is the usual one.

In this chapter we establish reverses of the Bessel inequality and some Gr%
\"{u}ss type inequalities for orthonormal families, namely, upper bounds for
the expressions%
\begin{equation*}
\left\Vert x\right\Vert ^{2}-\sum_{i\in I}\left\vert \left\langle
x,e_{i}\right\rangle \right\vert ^{2},\text{ \ }\left\Vert x\right\Vert
-\left( \sum_{i\in I}\left\vert \left\langle x,e_{i}\right\rangle
\right\vert ^{2}\right) ^{\frac{1}{2}},\ \ \ x\in X
\end{equation*}%
and%
\begin{equation*}
\left\vert \left\langle x,y\right\rangle -\sum_{i\in I}\left\langle
x,e_{i}\right\rangle \left\langle e_{i},y\right\rangle \right\vert ,\ \ \
x,y\in H,
\end{equation*}%
under various assumptions for the vectors $x,y$ and the orthonormal family $%
\left\{ e_{i}\right\} _{i\in I}.$

\section{Reverses of Bessel's Inequality}

\subsection{Introduction}

In \cite{SSD1.8}, the author has proved the following Gr\"{u}ss type
inequality in real or complex inner product spaces.

\begin{theorem}
\label{t1.8}Let $\left( H,\left\langle \cdot ,\cdot \right\rangle \right) $
be an inner product space over $\mathbb{K}$ $\left( \mathbb{K}=\mathbb{R},%
\mathbb{C}\right) $ and $e\in H,$ $\left\Vert e\right\Vert =1.$ If $\phi
,\Phi ,\gamma ,\Gamma $ are real or complex numbers and $x,y$ are vectors in 
$H$ such that the conditions%
\begin{equation}
\func{Re}\left\langle \Phi e-x,x-\phi e\right\rangle \geq 0\text{ \ and \ }%
\func{Re}\left\langle \Gamma e-y,y-\gamma e\right\rangle \geq 0
\label{1.1.8}
\end{equation}%
hold, then we have the inequality%
\begin{equation}
\left\vert \left\langle x,y\right\rangle -\left\langle x,e\right\rangle
\left\langle e,y\right\rangle \right\vert \leq \frac{1}{4}\left\vert \Phi
-\phi \right\vert \left\vert \Gamma -\gamma \right\vert .  \label{1.2.8}
\end{equation}%
The constant $\frac{1}{4}$ is best possible in the sense that it cannot be
replaced by a smaller constant.
\end{theorem}

In \cite{SSD2.8}, the following refinement of (\ref{1.2.8}) has been pointed
out.

\begin{theorem}
\label{t2.8}Let $H,$ $\mathbb{K}$ and $e$ be as in Theorem \ref{t1.8}. If $%
\phi ,\Phi ,\gamma ,\Gamma ,x,y$ satisfy (\ref{1.1.8}) or equivalently 
\begin{equation*}
\left\Vert x-\frac{\phi +\Phi }{2}e\right\Vert \leq \frac{1}{2}\left\vert
\Phi -\phi \right\vert ,\ \ \ \left\Vert y-\frac{\gamma +\Gamma }{2}%
e\right\Vert \leq \frac{1}{2}\left\vert \Gamma -\gamma \right\vert ,
\end{equation*}%
then%
\begin{multline}
\left\vert \left\langle x,y\right\rangle -\left\langle x,e\right\rangle
\left\langle e,y\right\rangle \right\vert  \label{1.4.8} \\
\leq \frac{1}{4}\left\vert \Phi -\phi \right\vert \left\vert \Gamma -\gamma
\right\vert -\left[ \func{Re}\left\langle \Phi e-x,x-\phi e\right\rangle %
\right] ^{\frac{1}{2}}\left[ \func{Re}\left\langle \Gamma e-y,y-\gamma
e\right\rangle \right] ^{\frac{1}{2}}.
\end{multline}
\end{theorem}

In \cite{NU.8}, N. Ujevi\'{c} has generalised Theorem \ref{t1.8} for the
case of real inner product spaces as follows.

\begin{theorem}
\label{t3.8}Let $\left( H,\left\langle \cdot ,\cdot \right\rangle \right) $
be an inner product space over the real number field $\mathbb{R}$, and $%
\left\{ e_{i}\right\} _{i\in \left\{ 1,\dots ,n\right\} }$ an orthornormal
family in $H.$ If $\phi _{i},\gamma _{i},\Phi _{i},\Gamma _{i}\in \mathbb{R}$%
, $i\in \left\{ 1,\dots ,n\right\} $ satisfy the condition%
\begin{equation*}
\left\langle \sum_{i=1}^{n}\Phi _{i}e_{i}-x,x-\sum_{i=1}^{n}\phi
_{i}e_{i}\right\rangle \geq 0,\ \ \ \ \ \ \left\langle \sum_{i=1}^{n}\Gamma
_{i}e_{i}-y,y-\sum_{i=1}^{n}\gamma _{i}e_{i}\right\rangle \geq 0,
\end{equation*}%
then one has the inequality:%
\begin{equation}
\left\vert \left\langle x,y\right\rangle -\sum_{i=1}^{n}\left\langle
x,e_{i}\right\rangle \left\langle e_{i},y\right\rangle \right\vert \leq 
\frac{1}{4}\left[ \sum_{i=1}^{n}\left( \Phi _{i}-\phi _{i}\right) ^{2}\cdot
\sum_{i=1}^{n}\left( \Gamma _{i}-\gamma _{i}\right) ^{2}\right] ^{\frac{1}{2}%
}.  \label{1.6.8}
\end{equation}%
The constant $\frac{1}{4}$ is best possible in the sense that it cannot be
replaced by a smaller constant.
\end{theorem}

We note that the key point in his proof is the following identity:%
\begin{multline*}
\sum_{i=1}^{n}\left( \left\langle x,e_{i}\right\rangle -\phi _{i}\right)
\left( \Phi _{i}-\left\langle x,e_{i}\right\rangle \right) -\left\langle
x-\sum_{i=1}^{n}\phi _{i}e_{i},\sum_{i=1}^{n}\Phi _{i}e_{i}-x\right\rangle \\
=\left\Vert x\right\Vert ^{2}-\sum_{i=1}^{n}\left\langle
x,e_{i}\right\rangle ^{2},
\end{multline*}%
holding for $x\in H,$ $\phi _{i},\Phi _{i}\in \mathbb{R}$, $i\in \left\{
1,\dots ,n\right\} $ and $\left\{ e_{i}\right\} _{i\in \left\{ 1,\dots
,n\right\} }$ an orthornormal family of vectors in the real inner product
space $H.$

In this section, by following \cite{9NSSD}, we point out a reverse of
Bessel's inequality in both real and complex inner product spaces. This
result will then be employed to provide a refinement of the Gr\"{u}ss type
inequality $\left( \ref{1.6.8}\right) $ for real or complex inner products.
Related results as well as integral inequalities for general measure spaces
are also given.

\subsection{A General Result}

The following lemma holds \cite{9NSSD}.

\begin{lemma}
\label{l2.1.8}Let $\left\{ e_{i}\right\} _{i\in I}$ be a family of
orthornormal vectors in $H,$ $F$ a finite part of $I$ and $\phi _{i},\Phi
_{i}$ $\left( i\in F\right) ,$ real or complex numbers. The following
statements are equivalent for $x\in H:$

\begin{enumerate}
\item[(i)] $\func{Re}\left\langle \sum_{i\in F}\Phi _{i}e_{i}-x,x-\sum_{i\in
F}\phi _{i}e_{i}\right\rangle \geq 0,$

\item[(ii)] $\left\Vert x-\sum_{i\in F}\frac{\phi _{i}+\Phi _{i}}{2}%
e_{i}\right\Vert \leq \frac{1}{2}\left( \sum_{i\in F}\left\vert \Phi
_{i}-\phi _{i}\right\vert ^{2}\right) ^{\frac{1}{2}}.$
\end{enumerate}
\end{lemma}

\begin{proof}
It is easy to see that for $y,a,A\in H,$ the following are equivalent (see 
\cite[Lemma 1]{SSD2.8})

\begin{enumerate}
\item[(b)] $\func{Re}\left\langle A-y,y-a\right\rangle \geq 0$ and

\item[(bb)] $\left\Vert y-\frac{a+A}{2}\right\Vert \leq \frac{1}{2}%
\left\Vert A-a\right\Vert .$
\end{enumerate}

Now, for $a=\sum_{i\in F}\phi _{i}e_{i},$ $A=\sum_{i\in F}\Phi _{i}e_{i},$
we have%
\begin{align*}
\left\Vert A-a\right\Vert & =\left\Vert \sum_{i\in F}\left( \Phi _{i}-\phi
_{i}\right) e_{i}\right\Vert =\left( \left\Vert \sum_{i\in F}\left( \Phi
_{i}-\phi _{i}\right) e_{i}\right\Vert ^{2}\right) ^{\frac{1}{2}} \\
& =\left( \sum_{i\in F}\left\vert \Phi _{i}-\phi _{i}\right\vert
^{2}\left\Vert e_{i}\right\Vert ^{2}\right) ^{\frac{1}{2}}=\left( \sum_{i\in
F}\left\vert \Phi _{i}-\phi _{i}\right\vert ^{2}\right) ^{\frac{1}{2}},
\end{align*}%
giving, for $y=x,$ the desired equivalence.
\end{proof}

The following reverse of Bessel's inequality holds \cite{9NSSD}.

\begin{theorem}
\label{t2.2.8}Let $\left\{ e_{i}\right\} _{i\in I},$ $F,$ $\phi _{i},\Phi
_{i},$ $i\in F$ and $x\in H$ such that either (i) or (ii) of Lemma \ref%
{l2.1.8} holds. Then we have the inequality:%
\begin{align}
0& \leq \left\Vert x\right\Vert ^{2}-\sum_{i\in F}\left\vert \left\langle
x,e_{i}\right\rangle \right\vert ^{2}  \label{2.4.8} \\
& \leq \frac{1}{4}\sum_{i\in F}\left\vert \Phi _{i}-\phi _{i}\right\vert
^{2}-\func{Re}\left\langle \sum_{i\in F}\Phi _{i}e_{i}-x,x-\sum_{i\in F}\phi
_{i}e_{i}\right\rangle  \notag \\
& \leq \frac{1}{4}\sum_{i\in F}\left\vert \Phi _{i}-\phi _{i}\right\vert
^{2}.  \notag
\end{align}%
The constant $\frac{1}{4}$ is best in both inequalities.
\end{theorem}

\begin{proof}
Define%
\begin{equation*}
I_{1}:=\sum_{i\in H}\func{Re}\left[ \left( \Phi _{i}-\left\langle
x,e_{i}\right\rangle \right) \left( \overline{\left\langle
x,e_{i}\right\rangle }-\overline{\phi _{i}}\right) \right]
\end{equation*}%
and%
\begin{equation*}
I_{2}:=\func{Re}\left[ \left\langle \sum_{i\in H}\Phi
_{i}e_{i}-x,x-\sum_{i\in H}\phi _{i}e_{i}\right\rangle \right] .
\end{equation*}%
Observe that%
\begin{equation*}
I_{1}=\sum_{i\in H}\func{Re}\left[ \Phi _{i}\overline{\left\langle
x,e_{i}\right\rangle }\right] +\sum_{i\in H}\func{Re}\left[ \overline{\phi
_{i}}\left\langle x,e_{i}\right\rangle \right] -\sum_{i\in H}\func{Re}\left[
\Phi _{i}\overline{\phi _{i}}\right] -\sum_{i\in H}\left\vert \left\langle
x,e_{i}\right\rangle \right\vert ^{2}
\end{equation*}%
and%
\begin{align*}
I_{2}& =\func{Re}\left[ \sum_{i\in H}\Phi _{i}\overline{\left\langle
x,e_{i}\right\rangle }+\sum_{i\in H}\overline{\phi _{i}}\left\langle
x,e_{i}\right\rangle -\left\Vert x\right\Vert ^{2}-\sum_{i\in H}\sum_{j\in
H}\Phi _{i}\overline{\phi _{i}}\left\langle e_{i},e_{j}\right\rangle \right]
\\
& =\sum_{i\in H}\func{Re}\left[ \Phi _{i}\overline{\left\langle
x,e_{i}\right\rangle }\right] +\sum_{i\in H}\func{Re}\left[ \overline{\phi
_{i}}\left\langle x,e_{i}\right\rangle \right] -\left\Vert x\right\Vert
^{2}-\sum_{i\in H}\func{Re}\left[ \Phi _{i}\overline{\phi _{i}}\right] .
\end{align*}%
Consequently, subtracting $I_{2}$ from $I_{1},$ we deduce the following
equality that is interesting in its turn%
\begin{multline}
\left\Vert x\right\Vert ^{2}-\sum_{i\in F}\left\vert \left\langle
x,e_{i}\right\rangle \right\vert ^{2}=\sum_{i\in H}\func{Re}\left[ \left(
\Phi _{i}-\left\langle x,e_{i}\right\rangle \right) \left( \overline{%
\left\langle x,e_{i}\right\rangle }-\overline{\phi _{i}}\right) \right]
\label{2.5.8} \\
-\func{Re}\left[ \left\langle \sum_{i\in H}\Phi _{i}e_{i}-x,x-\sum_{i\in
H}\phi _{i}e_{i}\right\rangle \right] .
\end{multline}%
Using the following elementary inequality for complex numbers%
\begin{equation*}
\func{Re}\left[ a\overline{b}\right] \leq \frac{1}{4}\left\vert
a+b\right\vert ^{2},\ \ \ \ a,b\in \mathbb{K},
\end{equation*}%
for the choices $a=\Phi _{i}-\left\langle x,e_{i}\right\rangle ,$ $%
b=\left\langle x,e_{i}\right\rangle -\phi _{i}$ \ $\left( i\in F\right) ,$
we deduce%
\begin{equation}
\sum_{i\in H}\func{Re}\left[ \left( \Phi _{i}-\left\langle
x,e_{i}\right\rangle \right) \left( \overline{\left\langle
x,e_{i}\right\rangle }-\overline{\phi _{i}}\right) \right] \leq \frac{1}{4}%
\sum_{i\in H}\left\vert \Phi _{i}-\phi _{i}\right\vert ^{2}.  \label{2.6.8}
\end{equation}%
Making use of (\ref{2.5.8}), (\ref{2.6.8}) and the assumption (i), we deduce
(\ref{2.4.8}).

The sharpness of the constant $\frac{1}{4}$ was proved for a single element $%
e,$ $\left\Vert e\right\Vert =1$ in \cite{SSD1.8}, or for the real case in 
\cite{NU.8}.

We can give here a simple proof as follows.

Assume that there is a $c>0$ such that%
\begin{align}
0& \leq \left\Vert x\right\Vert ^{2}-\sum_{i\in F}\left\vert \left\langle
x,e_{i}\right\rangle \right\vert ^{2}  \label{2.7.8} \\
& \leq c\sum_{i\in F}\left\vert \Phi _{i}-\phi _{i}\right\vert ^{2}-\func{Re}%
\left\langle \sum_{i\in F}\Phi _{i}e_{i}-x,x-\sum_{i\in F}\phi
_{i}e_{i}\right\rangle ,  \notag
\end{align}%
provided $\phi _{i},\Phi _{i},$ $x$ and $F$ satisfy (i) or (ii).

We choose $F=\left\{ 1\right\} ,$ $e_{1}=e_{2}=\left( \frac{1}{\sqrt{2}},%
\frac{1}{\sqrt{2}}\right) \in \mathbb{R}^{2},$ $x=\left( x_{1},x_{2}\right)
\in \mathbb{R}^{2},\ \Phi _{1}=\Phi =m>0,$ $\phi _{1}=\phi =-m,$ $H=\mathbb{R%
}^{2}$ to get from (\ref{2.7.8}) that%
\begin{align}
0& \leq x_{1}^{2}+x_{2}^{2}-\frac{\left( x_{1}+x_{2}\right) ^{2}}{2}
\label{2.8.9a} \\
& \leq 4cm^{2}-\left( \frac{m}{\sqrt{2}}-x_{1}\right) \left( x_{1}+\frac{m}{%
\sqrt{2}}\right) -\left( \frac{m}{\sqrt{2}}-x_{2}\right) \left( x_{2}+\frac{m%
}{\sqrt{2}}\right) ,  \notag
\end{align}%
provided%
\begin{align}
0& \leq \left\langle me-x,x+me\right\rangle  \label{2.9.9} \\
& =\left( \frac{m}{\sqrt{2}}-x_{1}\right) \left( x_{1}+\frac{m}{\sqrt{2}}%
\right) +\left( \frac{m}{\sqrt{2}}-x_{2}\right) \left( x_{2}+\frac{m}{\sqrt{2%
}}\right) .  \notag
\end{align}%
If we choose $x_{1}=\frac{m}{\sqrt{2}},$ $x_{2}=-\frac{m}{\sqrt{2}},$ then (%
\ref{2.9.9}) is fulfilled and by (\ref{2.8.9a}) we get $m^{2}\leq 4cm^{2},$
giving $c\geq \frac{1}{4}.$
\end{proof}

\subsection{A Refinement of the Gr\"{u}ss Inequality for Orthonormal Families%
}

The following result holds \cite{9NSSD}.

\begin{theorem}
\label{t3.1.8}Let $\left\{ e_{i}\right\} _{i\in I}$ be a family of
orthornormal vectors in $H,$ $F$ a finite part of $I$ and $\phi _{i},\Phi
_{i},$ $\gamma _{i},\Gamma _{i}\in \mathbb{K},\ i\in F$ and $x,y\in H.$ If
either%
\begin{align*}
\func{Re}\left\langle \sum_{i\in F}\Phi _{i}e_{i}-x,x-\sum_{i\in F}\phi
_{i}e_{i}\right\rangle & \geq 0, \\
\func{Re}\left\langle \sum_{i\in F}\Gamma _{i}e_{i}-y,y-\sum_{i\in F}\gamma
_{i}e_{i}\right\rangle & \geq 0,
\end{align*}%
or equivalently,%
\begin{align*}
\left\Vert x-\sum_{i\in F}\frac{\Phi _{i}+\phi _{i}}{2}e_{i}\right\Vert &
\leq \frac{1}{2}\left( \sum_{i\in F}\left\vert \Phi _{i}-\phi
_{i}\right\vert ^{2}\right) ^{\frac{1}{2}}, \\
\left\Vert y-\sum_{i\in F}\frac{\Gamma _{i}+\gamma _{i}}{2}e_{i}\right\Vert
& \leq \frac{1}{2}\left( \sum_{i\in F}\left\vert \Gamma _{i}-\gamma
_{i}\right\vert ^{2}\right) ^{\frac{1}{2}},
\end{align*}%
hold, then we have the inequalities%
\begin{align}
& \left\vert \left\langle x,y\right\rangle -\sum_{i\in F}\left\langle
x,e_{i}\right\rangle \left\langle e_{i},y\right\rangle \right\vert
\label{3.3.8} \\
& \leq \frac{1}{4}\left( \sum_{i\in F}\left\vert \Phi _{i}-\phi
_{i}\right\vert ^{2}\right) ^{\frac{1}{2}}\cdot \left( \sum_{i\in
F}\left\vert \Gamma _{i}-\gamma _{i}\right\vert ^{2}\right) ^{\frac{1}{2}} 
\notag \\
& \ \ \ \ \ \ -\left[ \func{Re}\left\langle \sum_{i\in F}\Phi
_{i}e_{i}-x,x-\sum_{i\in F}\phi _{i}e_{i}\right\rangle \right] ^{\frac{1}{2}}%
\left[ \func{Re}\left\langle \sum_{i\in F}\Gamma _{i}e_{i}-y,y-\sum_{i\in
F}\gamma _{i}e_{i}\right\rangle \right] ^{\frac{1}{2}}  \notag \\
& \leq \frac{1}{4}\left( \sum_{i\in F}\left\vert \Phi _{i}-\phi
_{i}\right\vert ^{2}\right) ^{\frac{1}{2}}\cdot \left( \sum_{i\in
F}\left\vert \Gamma _{i}-\gamma _{i}\right\vert ^{2}\right) ^{\frac{1}{2}}. 
\notag
\end{align}%
The constant $\frac{1}{4}$ is best possible.
\end{theorem}

\begin{proof}
Using Schwartz's inequality in the inner product space $\left(
H,\left\langle \cdot ,\cdot \right\rangle \right) $ one has%
\begin{multline}
\left\vert \left\langle x-\sum_{i\in F}\left\langle x,e_{i}\right\rangle
e_{i},y-\sum_{i\in F}\left\langle y,e_{i}\right\rangle e_{i}\right\rangle
\right\vert ^{2}  \label{3.4.8} \\
\leq \left\Vert x-\sum_{i\in F}\left\langle x,e_{i}\right\rangle
e_{i}\right\Vert ^{2}\left\Vert y-\sum_{i\in F}\left\langle
y,e_{i}\right\rangle e_{i}\right\Vert ^{2}
\end{multline}%
and since a simple calculation shows that 
\begin{equation*}
\left\langle x-\sum_{i\in F}\left\langle x,e_{i}\right\rangle
e_{i},y-\sum_{i\in F}\left\langle y,e_{i}\right\rangle e_{i}\right\rangle
=\left\langle x,y\right\rangle -\sum_{i\in F}\left\langle
x,e_{i}\right\rangle \left\langle e_{i},y\right\rangle
\end{equation*}%
and 
\begin{equation*}
\left\Vert x-\sum_{i\in F}\left\langle x,e_{i}\right\rangle e_{i}\right\Vert
^{2}=\left\Vert x\right\Vert ^{2}-\sum_{i\in F}\left\vert \left\langle
x,e_{i}\right\rangle \right\vert ^{2}
\end{equation*}%
for any $x,y\in H,$ then by (\ref{3.4.8}) and by the reverse of Bessel's
inequality in Theorem \ref{t2.2.8}, we have%
\begin{align}
& \left\vert \left\langle x,y\right\rangle -\sum_{i\in F}\left\langle
x,e_{i}\right\rangle \left\langle e_{i},y\right\rangle \right\vert ^{2}
\label{3.5.8} \\
& \leq \left( \left\Vert x\right\Vert ^{2}-\sum_{i\in F}\left\vert
\left\langle x,e_{i}\right\rangle \right\vert ^{2}\right) \left( \left\Vert
y\right\Vert ^{2}-\sum_{i\in F}\left\vert \left\langle y,e_{i}\right\rangle
\right\vert ^{2}\right)  \notag \\
& \leq \left[ \frac{1}{4}\sum_{i\in F}\left\vert \Phi _{i}-\phi
_{i}\right\vert ^{2}-\func{Re}\left\langle \sum_{i\in F}\Phi
_{i}e_{i}-x,x-\sum_{i\in F}\phi _{i}e_{i}\right\rangle \right]  \notag \\
& \quad \times \left[ \frac{1}{4}\sum_{i\in F}\left\vert \Gamma _{i}-\gamma
_{i}\right\vert ^{2}-\func{Re}\left\langle \sum_{i\in F}\Gamma
_{i}e_{i}-y,y-\sum_{i\in F}\gamma _{i}e_{i}\right\rangle \right]  \notag \\
& \leq \left[ \frac{1}{4}\left( \sum_{i\in F}\left\vert \Phi _{i}-\phi
_{i}\right\vert ^{2}\right) ^{\frac{1}{2}}\cdot \left( \sum_{i\in
F}\left\vert \Gamma _{i}-\gamma _{i}\right\vert ^{2}\right) ^{\frac{1}{2}%
}\right.  \notag \\
& \quad -\left. \left[ \func{Re}\left\langle \sum_{i\in F}\Phi
_{i}e_{i}-x,x-\sum_{i\in F}\phi _{i}e_{i}\right\rangle \right] ^{\frac{1}{2}}%
\left[ \func{Re}\left\langle \sum_{i\in F}\Gamma _{i}e_{i}-y,y-\sum_{i\in
F}\gamma _{i}e_{i}\right\rangle \right] ^{\frac{1}{2}}\right]  \notag
\end{align}%
where, for the last inequality, we have made use of the inequality%
\begin{equation*}
\left( m^{2}-n^{2}\right) \left( p^{2}-q^{2}\right) \leq \left( mp-nq\right)
^{2},
\end{equation*}%
holding for any $m,n,p,q>0.$

Taking the square root in (\ref{3.5.8}) and observing that the quantity in
the last square bracket is nonnegative (see for example (\ref{2.4.8})), we
deduce the desired result (\ref{3.3.8}).

The best constant has been proved in \cite{SSD1.8} for one element and we
omit the details.
\end{proof}

\subsection{Some Companion Inequalities}

The following companion of the Gr\"{u}ss inequality also holds \cite{9NSSD}.

\begin{theorem}
\label{t4.1.8}Let $\left\{ e_{i}\right\} _{i\in I}$ be a family of
orthornormal vectors in $H,$ $F$ a finite part of $I$ and $\phi _{i},\Phi
_{i}\in \mathbb{K},\ i\in F$ and $x,y\in H$ such that%
\begin{equation}
\func{Re}\left\langle \sum_{i\in F}\Phi _{i}e_{i}-\frac{x+y}{2},\frac{x+y}{2}%
-\sum_{i\in F}\phi _{i}e_{i}\right\rangle \geq 0  \label{4.1.8}
\end{equation}%
or equivalently,%
\begin{equation*}
\left\Vert \frac{x+y}{2}-\sum_{i\in F}\frac{\Phi _{i}+\phi _{i}}{2}\cdot
e_{i}\right\Vert \leq \frac{1}{2}\left( \sum_{i\in F}\left\vert \Phi
_{i}-\phi _{i}\right\vert ^{2}\right) ^{\frac{1}{2}},
\end{equation*}%
holds, then we have the inequality%
\begin{equation}
\func{Re}\left[ \left\langle x,y\right\rangle -\sum_{i\in F}\left\langle
x,e_{i}\right\rangle \left\langle e_{i},y\right\rangle \right] \leq \frac{1}{%
4}\sum_{i\in F}\left\vert \Phi _{i}-\phi _{i}\right\vert ^{2}.  \label{4.3.8}
\end{equation}%
The constant $\frac{1}{4}$ is best possible.
\end{theorem}

\begin{proof}
Start with the well known inequality 
\begin{equation}
\func{Re}\left\langle z,u\right\rangle \leq \frac{1}{4}\left\Vert
z+u\right\Vert ^{2},\ \ \ z,u\in H.  \label{4.4.8}
\end{equation}%
Since%
\begin{equation*}
\left\langle x,y\right\rangle -\sum_{i\in F}\left\langle
x,e_{i}\right\rangle \left\langle e_{i},y\right\rangle =\left\langle
x-\sum_{i\in F}\left\langle x,e_{i}\right\rangle e_{i},y-\sum_{i\in
F}\left\langle y,e_{i}\right\rangle e_{i}\right\rangle ,
\end{equation*}%
for any \thinspace $x,y\in H,$ then, by (\ref{4.4.8}), we get%
\begin{align}
& \func{Re}\left[ \left\langle x,y\right\rangle -\sum_{i\in F}\left\langle
x,e_{i}\right\rangle \left\langle e_{i},y\right\rangle \right]  \label{4.5.8}
\\
& =\func{Re}\left[ \left\langle x-\sum_{i\in F}\left\langle
x,e_{i}\right\rangle e_{i},y-\sum_{i\in F}\left\langle y,e_{i}\right\rangle
e_{i}\right\rangle \right]  \notag \\
& \leq \frac{1}{4}\left\Vert x-\sum_{i\in F}\left\langle
x,e_{i}\right\rangle e_{i}+y-\sum_{i\in F}\left\langle y,e_{i}\right\rangle
e_{i}\right\Vert ^{2}  \notag \\
& =\left\Vert \frac{x+y}{2}-\sum_{i\in F}\left\langle \frac{x+y}{2}%
,e_{i}\right\rangle e_{i}\right\Vert ^{2}  \notag \\
& =\left\Vert \frac{x+y}{2}\right\Vert ^{2}-\sum_{i\in F}\left\vert
\left\langle \frac{x+y}{2},e_{i}\right\rangle \right\vert ^{2}.  \notag
\end{align}%
If we apply the reverse of Bessel's inequality in Theorem \ref{t2.2.8} for $%
\frac{x+y}{2},$ we may state that%
\begin{equation}
\left\Vert \frac{x+y}{2}\right\Vert ^{2}-\sum_{i\in F}\left\vert
\left\langle \frac{x+y}{2},e_{i}\right\rangle \right\vert ^{2}\leq \frac{1}{4%
}\left( \sum_{i\in F}\left\vert \Phi _{i}-\phi _{i}\right\vert ^{2}\right) ^{%
\frac{1}{2}}.  \label{4.6.8}
\end{equation}%
Now, by making use of (\ref{4.5.8}) and (\ref{4.6.8}), we deduce (\ref{4.3.8}%
).

The fact that $\frac{1}{4}$ is the best constant in (\ref{4.3.8}) follows by
the fact that if in (\ref{4.1.8}) we choose $x=y,$ then it becomes (i) of
Lemma \ref{l2.1.8}, implying (\ref{2.4.8}), for which, we have shown that $%
\frac{1}{4}$ was the best constant.
\end{proof}

The following corollary may be of interest if we wish to evaluate the
absolute value of%
\begin{equation*}
\func{Re}\left[ \left\langle x,y\right\rangle -\sum_{i\in F}\left\langle
x,e_{i}\right\rangle \left\langle e_{i},y\right\rangle \right] .
\end{equation*}

\begin{corollary}
\label{c4.2.8}With the assumptions of Theorem \ref{t4.1.8} and if 
\begin{equation*}
\func{Re}\left\langle \sum_{i\in F}\Phi _{i}e_{i}-\frac{x\pm y}{2},\frac{%
x\pm y}{2}-\sum_{i\in F}\phi _{i}e_{i}\right\rangle \geq 0
\end{equation*}%
or equivalently,%
\begin{equation*}
\left\Vert \frac{x\pm y}{2}-\sum_{i\in F}\frac{\Phi _{i}+\phi _{i}}{2}\cdot
e_{i}\right\Vert \leq \frac{1}{2}\left( \sum_{i\in F}\left\vert \Phi
_{i}-\phi _{i}\right\vert ^{2}\right) ^{\frac{1}{2}},
\end{equation*}%
holds, then we have the inequality%
\begin{equation}
\left\vert \func{Re}\left[ \left\langle x,y\right\rangle -\sum_{i\in
F}\left\langle x,e_{i}\right\rangle \left\langle e_{i},y\right\rangle \right]
\right\vert \leq \frac{1}{4}\sum_{i\in F}\left\vert \Phi _{i}-\phi
_{i}\right\vert ^{2}.  \label{4.9.8}
\end{equation}
\end{corollary}

\begin{proof}
We only remark that, if%
\begin{equation*}
\func{Re}\left\langle \sum_{i\in F}\Phi _{i}e_{i}-\frac{x-y}{2},\frac{x-y}{2}%
-\sum_{i\in F}\phi _{i}e_{i}\right\rangle \geq 0
\end{equation*}%
holds, then by Theorem \ref{t4.1.8} for $\left( -y\right) $ instead of $y,$
we have%
\begin{equation*}
\func{Re}\left[ -\left\langle x,y\right\rangle +\sum_{i\in F}\left\langle
x,e_{i}\right\rangle \left\langle e_{i},y\right\rangle \right] \leq \frac{1}{%
4}\sum_{i\in F}\left\vert \Phi _{i}-\phi _{i}\right\vert ^{2},
\end{equation*}%
showing that%
\begin{equation}
\func{Re}\left[ \left\langle x,y\right\rangle -\sum_{i\in F}\left\langle
x,e_{i}\right\rangle \left\langle e_{i},y\right\rangle \right] \geq -\frac{1%
}{4}\sum_{i\in F}\left\vert \Phi _{i}-\phi _{i}\right\vert ^{2}.
\label{4.10.8}
\end{equation}%
Making use of (\ref{4.3.8}) and (\ref{4.10.8}), we deduce the desired
inequality (\ref{4.9.8}).
\end{proof}

\begin{remark}
\label{r4.3.8}If $H$ is a real inner product space and $m_{i},M_{i}\in 
\mathbb{R}$ with the property that%
\begin{equation*}
\left\langle \sum_{i\in F}M_{i}e_{i}-\frac{x\pm y}{2},\frac{x\pm y}{2}%
-\sum_{i\in F}m_{i}e_{i}\right\rangle \geq 0
\end{equation*}%
or equivalently,%
\begin{equation*}
\left\Vert \frac{x\pm y}{2}-\sum_{i\in F}\frac{M_{i}+m_{i}}{2}\cdot
e_{i}\right\Vert \leq \frac{1}{2}\left( \sum_{i\in F}\left(
M_{i}-m_{i}\right) ^{2}\right) ^{\frac{1}{2}},
\end{equation*}%
then we have the Gr\"{u}ss type inequality%
\begin{equation*}
\left\vert \left\langle x,y\right\rangle -\sum_{i\in F}\left\langle
x,e_{i}\right\rangle \left\langle e_{i},y\right\rangle \right\vert \leq 
\frac{1}{4}\sum_{i\in F}\left( M_{i}-m_{i}\right) ^{2}.
\end{equation*}
\end{remark}

\subsection{Integral Inequalities}

Let $\left( \Omega ,\Sigma ,\mu \right) $ be a measure space consisting of a
set $\Omega ,$ $\Sigma $ a $\sigma -$algebra of parts and $\mu $ a countably
additive and positive measure on $\Sigma $ with values in $\mathbb{R}\cup
\left\{ \infty \right\} .$ Let $\rho \geq 0$ be a $\mu -$measurable function
on $\Omega .$ Denote by $L_{\rho }^{2}\left( \Omega ,\mathbb{K}\right) $ the
Hilbert space of all real or complex valued functions defined on $\Omega $
and $2-\rho -$integrable on $\Omega ,$ i.e.,%
\begin{equation*}
\int_{\Omega }\rho \left( s\right) \left\vert f\left( s\right) \right\vert
^{2}d\mu \left( s\right) <\infty .
\end{equation*}

Consider the family $\left\{ f_{i}\right\} _{i\in I}$ of functions in $%
L_{\rho }^{2}\left( \Omega ,\mathbb{K}\right) $ with the properties that%
\begin{equation*}
\int_{\Omega }\rho \left( s\right) f_{i}\left( s\right) \overline{f_{j}}%
\left( s\right) d\mu \left( s\right) =\delta _{ij},\ \ \ i,j\in I,
\end{equation*}%
where $\delta _{ij}$ is $0$ if $i\neq j$ and $\delta _{ij}=1$ if $i=j.$ $%
\left\{ f_{i}\right\} _{i\in I}$ is an orthornormal family in $L_{\rho
}^{2}\left( \Omega ,\mathbb{K}\right) .$

The following proposition holds \cite{9NSSD}.

\begin{proposition}
\label{p5.1.8}Let $\left\{ f_{i}\right\} _{i\in I}$ be an orthornormal
family of functions in $L_{\rho }^{2}\left( \Omega ,\mathbb{K}\right) ,$ $F$
a finite subset of $I,$ $\phi _{i},\Phi _{i}\in \mathbb{K}$ $\left( i\in
F\right) $ and $f\in L_{\rho }^{2}\left( \Omega ,\mathbb{K}\right) ,$ such
that either%
\begin{multline}
\int_{\Omega }\rho \left( s\right) \func{Re}\left[ \left( \sum_{i\in F}\Phi
_{i}f_{i}\left( s\right) -f\left( s\right) \right) \right.  \label{5.3.8} \\
\times \left. \left( \overline{f}\left( s\right) -\sum_{i\in F}\overline{%
\phi _{i}}\text{ }\overline{f_{i}}\left( s\right) \right) \right] d\mu
\left( s\right) \geq 0
\end{multline}%
or equivalently,%
\begin{equation*}
\int_{\Omega }\rho \left( s\right) \left\vert f\left( s\right) -\sum_{i\in F}%
\frac{\Phi _{i}+\phi _{i}}{2}f_{i}\left( s\right) \right\vert ^{2}d\mu
\left( s\right) \leq \frac{1}{4}\sum_{i\in F}\left\vert \Phi _{i}-\phi
_{i}\right\vert ^{2}.
\end{equation*}%
Then we have the inequality%
\begin{align}
0& \leq \int_{\Omega }\rho \left( s\right) \left\vert f\left( s\right)
\right\vert ^{2}d\mu \left( s\right) -\sum_{i\in F}\left\vert \int_{\Omega
}\rho \left( s\right) f\left( s\right) \overline{f_{i}}\left( s\right) d\mu
\left( s\right) \right\vert ^{2}  \label{5.5.8} \\
& \leq \frac{1}{4}\sum_{i\in F}\left\vert \Phi _{i}-\phi _{i}\right\vert ^{2}
\notag \\
& \ \ \ \ \ -\int_{\Omega }\rho \left( s\right) \func{Re}\left[ \left(
\sum_{i\in F}\Phi _{i}f_{i}\left( s\right) -f\left( s\right) \right) \left( 
\overline{f}\left( s\right) -\sum_{i\in F}\overline{\phi _{i}}\text{ }%
\overline{f_{i}}\left( s\right) \right) \right] d\mu \left( s\right)  \notag
\\
& \leq \frac{1}{4}\sum_{i\in F}\left\vert \Phi _{i}-\phi _{i}\right\vert
^{2}.  \notag
\end{align}%
The constant $\frac{1}{4}$ is best possible in both inequalities.
\end{proposition}

The proof follows by Theorem \ref{t2.2.8} applied for the Hilbert space $%
L_{\rho }^{2}\left( \Omega ,\mathbb{K}\right) $ and the orthornormal family $%
\left\{ f_{i}\right\} _{i\in I}.$

The following Gr\"{u}ss type inequality also holds \cite{9NSSD}.

\begin{proposition}
\label{p5.2.8}Let $\left\{ f_{i}\right\} _{i\in I}$ and $F$ be as in
Proposition \ref{p5.1.8}. If $\phi _{i},\Phi _{i},\gamma _{i},\Gamma _{i}\in 
\mathbb{K}$ $\left( i\in F\right) $ and $f,g\in L_{\rho }^{2}\left( \Omega ,%
\mathbb{K}\right) $ so that either%
\begin{align*}
\int_{\Omega }\rho \left( s\right) \func{Re}\left[ \left( \sum_{i\in F}\Phi
_{i}f_{i}\left( s\right) -f\left( s\right) \right) \left( \overline{f}\left(
s\right) -\sum_{i\in F}\overline{\phi _{i}}\text{ }\overline{f_{i}}\left(
s\right) \right) \right] d\mu \left( s\right) & \geq 0, \\
\int_{\Omega }\rho \left( s\right) \func{Re}\left[ \left( \sum_{i\in
F}\Gamma _{i}f_{i}\left( s\right) -g\left( s\right) \right) \left( \overline{%
g}\left( s\right) -\sum_{i\in F}\overline{\gamma _{i}}\text{ }\overline{f_{i}%
}\left( s\right) \right) \right] d\mu \left( s\right) & \geq 0,
\end{align*}%
or equivalently,%
\begin{align*}
\int_{\Omega }\rho \left( s\right) \left\vert f\left( s\right) -\sum_{i\in F}%
\frac{\Phi _{i}+\phi _{i}}{2}f_{i}\left( s\right) \right\vert ^{2}d\mu
\left( s\right) & \leq \frac{1}{4}\sum_{i\in F}\left\vert \Phi _{i}-\phi
_{i}\right\vert ^{2}, \\
\int_{\Omega }\rho \left( s\right) \left\vert g\left( s\right) -\sum_{i\in F}%
\frac{\Gamma _{i}+\gamma _{i}}{2}f_{i}\left( s\right) \right\vert ^{2}d\mu
\left( s\right) & \leq \frac{1}{4}\sum_{i\in F}\left\vert \Gamma _{i}-\gamma
_{i}\right\vert ^{2},
\end{align*}%
hold, then we have the inequalities%
\begin{multline}
\left\vert \int_{\Omega }\rho \left( s\right) f\left( s\right) g\left(
s\right) d\mu \left( s\right) \right.  \label{5.8.8} \\
-\left. \sum_{i\in F}\int_{\Omega }\rho \left( s\right) f\left( s\right) 
\overline{f_{i}}\left( s\right) d\mu \left( s\right) \int_{\Omega }\rho
\left( s\right) f_{i}\left( s\right) \overline{g\left( s\right) }d\mu \left(
s\right) \right\vert
\end{multline}%
\vspace{-0.2in}%
\begin{align*}
& \leq \frac{1}{4}\left( \sum_{i\in F}\left\vert \Phi _{i}-\phi
_{i}\right\vert ^{2}\right) ^{\frac{1}{2}}\left( \sum_{i\in F}\left\vert
\Gamma _{i}-\gamma _{i}\right\vert ^{2}\right) ^{\frac{1}{2}} \\
& \quad -\left[ \int_{\Omega }\rho \left( s\right) \func{Re}\left[ \left(
\sum_{i\in F}\Phi _{i}f_{i}\left( s\right) -f\left( s\right) \right) \left( 
\overline{f}\left( s\right) -\sum_{i\in F}\phi \overline{_{i}}\overline{f_{i}%
}\left( s\right) \right) \right] d\mu \left( s\right) \right] ^{\frac{1}{2}}
\\
& \quad \times \left[ \int_{\Omega }\rho \left( s\right) \func{Re}\left[
\left( \sum_{i\in F}\Gamma _{i}f_{i}\left( s\right) -g\left( s\right)
\right) \left( \overline{g}\left( s\right) -\sum_{i\in F}\overline{\gamma
_{i}}\overline{f_{i}}\left( s\right) \right) \right] d\mu \left( s\right) %
\right] ^{\frac{1}{2}} \\
& \leq \frac{1}{4}\left( \sum_{i\in F}\left\vert \Phi _{i}-\phi
_{i}\right\vert ^{2}\right) ^{\frac{1}{2}}\left( \sum_{i\in F}\left\vert
\Gamma _{i}-\gamma _{i}\right\vert ^{2}\right) ^{\frac{1}{2}}.
\end{align*}%
The constant $\frac{1}{4}$ is the best possible.
\end{proposition}

The proof follows by Theorem \ref{t3.1.8} and we omit the details.

\begin{remark}
\label{r5.3.8}Similar results may be stated if we apply the other
inequalities obtained above. We omit the details.
\end{remark}

In the case of real spaces, the following corollaries provide much simpler
sufficient conditions for the reverse of Bessel's inequality (\ref{5.5.8})
or for the Gr\"{u}ss type inequality (\ref{5.8.8}) to hold.

\begin{corollary}
\label{c5.4.8}Let $\left\{ f_{i}\right\} _{i\in I}$ be an orthornormal
family of functions in the real Hilbert space $L_{\rho }^{2}\left( \Omega
\right) ,$ $F$ a finite part of $I,$ $M_{i},m_{i}\in \mathbb{R}$ \ $\left(
i\in F\right) $ and $f\in L_{\rho }^{2}\left( \Omega \right) $ such that%
\begin{equation*}
\sum_{i\in F}m_{i}f_{i}\left( s\right) \leq f\left( s\right) \leq \sum_{i\in
F}M_{i}f_{i}\left( s\right) \text{ \ for \ }\mu -\text{a.e. }s\in \Omega .
\end{equation*}%
Then we have the inequalities%
\begin{align*}
0& \leq \int_{\Omega }\rho \left( s\right) f^{2}\left( s\right) d\mu \left(
s\right) -\sum_{i\in F}\left[ \int_{\Omega }\rho \left( s\right) f\left(
s\right) f_{i}\left( s\right) d\mu \left( s\right) \right] ^{2} \\
& \leq \frac{1}{4}\sum_{i\in F}\left( M_{i}-m_{i}\right) ^{2} \\
& \ \ \ \ \ \ -\int_{\Omega }\rho \left( s\right) \left( \sum_{i\in
F}M_{i}f_{i}\left( s\right) -f\left( s\right) \right) \left( f\left(
s\right) -\sum_{i\in F}m_{i}f_{i}\left( s\right) \right) d\mu \left( s\right)
\\
& \leq \frac{1}{4}\sum_{i\in F}\left( M_{i}-m_{i}\right) ^{2}.
\end{align*}%
The constant $\frac{1}{4}$ is best possible.
\end{corollary}

\begin{corollary}
\label{c5.5.8}Let $\left\{ f_{i}\right\} _{i\in I}$ and $F$ be as in
Corollary \ref{c5.4.8}. If $M_{i},m_{i},N_{i},n_{i}\in \mathbb{R}$ $\left(
i\in F\right) $ and $f,g\in L_{\rho }^{2}\left( \Omega \right) $ are such
that%
\begin{equation*}
\sum_{i\in F}m_{i}f_{i}\left( s\right) \leq f\left( s\right) \leq \sum_{i\in
F}M_{i}f_{i}\left( s\right)
\end{equation*}%
and%
\begin{equation*}
\sum_{i\in F}n_{i}f_{i}\left( s\right) \leq g\left( s\right) \leq \sum_{i\in
F}N_{i}f_{i}\left( s\right) ,\text{ \ for \ }\mu -\text{a.e. }s\in \Omega ,
\end{equation*}%
then we have the inequalities%
\begin{multline*}
\left\vert \int_{\Omega }\rho \left( s\right) f\left( s\right) g\left(
s\right) d\mu \left( s\right) \right. \\
-\left. \sum_{i\in F}\int_{\Omega }\rho \left( s\right) f\left( s\right)
f_{i}\left( s\right) d\mu \left( s\right) \int_{\Omega }\rho \left( s\right)
g\left( s\right) f_{i}\left( s\right) d\mu \left( s\right) \right\vert
\end{multline*}%
\vspace{-0.2in}%
\begin{align*}
& \leq \frac{1}{4}\left( \sum_{i\in F}\left( M_{i}-m_{i}\right) ^{2}\right)
^{\frac{1}{2}}\left( \sum_{i\in F}\left( N_{i}-n_{i}\right) ^{2}\right) ^{%
\frac{1}{2}} \\
& \quad -\left[ \int_{\Omega }\rho \left( s\right) \left( \sum_{i\in
F}M_{i}f_{i}\left( s\right) -f\left( s\right) \right) \left( f\left(
s\right) -\sum_{i\in F}m_{i}f_{i}\left( s\right) \right) d\mu \left(
s\right) \right] ^{\frac{1}{2}} \\
& \quad \times \left[ \int_{\Omega }\rho \left( s\right) \left( \sum_{i\in
F}N_{i}f_{i}\left( s\right) -g\left( s\right) \right) \left( g\left(
s\right) -\sum_{i\in F}n_{i}f_{i}\left( s\right) \right) d\mu \left(
s\right) \right] ^{\frac{1}{2}} \\
& \leq \frac{1}{4}\left( \sum_{i\in F}\left( M_{i}-m_{i}\right) ^{2}\right)
^{\frac{1}{2}}\left( \sum_{i\in F}\left( N_{i}-n_{i}\right) ^{2}\right) ^{%
\frac{1}{2}}.
\end{align*}
\end{corollary}

\newpage

\section{Another Reverse for Bessel's Inequality}

\subsection{A General Result}

The following lemma holds \cite{10NSSD}.

\begin{lemma}
\label{l2.1.9}Let $\left\{ e_{i}\right\} _{i\in I}$ be a family of
orthornormal vectors in $H,$ $F$ a finite part of $I,$ $\lambda _{i}\in 
\mathbb{K},$ $i\in F$, $r>0$ and $x\in H.$ If%
\begin{equation*}
\left\Vert x-\sum_{i\in F}\lambda _{i}e_{i}\right\Vert \leq r,
\end{equation*}%
then we have the inequality%
\begin{equation}
0\leq \left\Vert x\right\Vert ^{2}-\sum_{i\in F}\left\vert \left\langle
x,e_{i}\right\rangle \right\vert ^{2}\leq r^{2}-\sum_{i\in F}\left\vert
\lambda _{i}-\left\langle x,e_{i}\right\rangle \right\vert ^{2}.
\label{2.2.9}
\end{equation}
\end{lemma}

\begin{proof}
Consider%
\begin{align*}
I_{1}& :=\left\Vert x-\sum_{i\in F}\lambda _{i}e_{i}\right\Vert
^{2}=\left\langle x-\sum_{i\in F}\lambda _{i}e_{i},x-\sum_{j\in F}\lambda
_{j}e_{j}\right\rangle \\
& =\left\Vert x\right\Vert ^{2}-\sum_{i\in F}\lambda _{i}\overline{%
\left\langle x,e_{i}\right\rangle }-\sum_{i\in F}\overline{\lambda _{i}}%
\left\langle x,e_{i}\right\rangle +\sum_{i\in F}\sum_{j\in F}\lambda _{i}%
\overline{\lambda _{j}}\left\langle e_{i},e_{j}\right\rangle \\
& =\left\Vert x\right\Vert ^{2}-\sum_{i\in F}\lambda _{i}\overline{%
\left\langle x,e_{i}\right\rangle }-\sum_{i\in F}\overline{\lambda _{i}}%
\left\langle x,e_{i}\right\rangle +\sum_{i\in F}\left\vert \lambda
_{i}\right\vert ^{2}
\end{align*}%
and%
\begin{align*}
I_{2}& :=\sum_{i\in F}\left\vert \lambda _{i}-\left\langle
x,e_{i}\right\rangle \right\vert ^{2}=\sum_{i\in F}\left( \lambda
_{i}-\left\langle x,e_{i}\right\rangle \right) \left( \overline{\lambda _{i}}%
-\overline{\left\langle x,e_{i}\right\rangle }\right) \\
& =\sum_{i\in F}\left[ \left\vert \lambda _{i}\right\vert ^{2}+\left\vert
\left\langle x,e_{i}\right\rangle \right\vert ^{2}-\overline{\lambda _{i}}%
\left\langle x,e_{i}\right\rangle -\lambda _{i}\overline{\left\langle
x,e_{i}\right\rangle }\right] \\
& =\sum_{i\in F}\left\vert \lambda _{i}\right\vert ^{2}+\sum_{i\in
F}\left\vert \left\langle x,e_{i}\right\rangle \right\vert ^{2}-\sum_{i\in F}%
\overline{\lambda _{i}}\left\langle x,e_{i}\right\rangle -\sum_{i\in
F}\lambda _{i}\overline{\left\langle x,e_{i}\right\rangle }.
\end{align*}%
If we subtract $I_{2}$ from $I_{1}$ we deduce the following identity that is
interesting in its own right%
\begin{equation*}
\left\Vert x-\sum_{i\in F}\lambda _{i}e_{i}\right\Vert ^{2}-\sum_{i\in
F}\left\vert \lambda _{i}-\left\langle x,e_{i}\right\rangle \right\vert
^{2}=\left\Vert x\right\Vert ^{2}-\sum_{i\in F}\left\vert \left\langle
x,e_{i}\right\rangle \right\vert ^{2},
\end{equation*}%
from which we easily deduce (\ref{2.2.9}).
\end{proof}

The following reverse of Bessel's inequality holds \cite{10NSSD}.

\begin{theorem}
\label{t2.2.9}Let $\left\{ e_{i}\right\} _{i\in I}$ be a family of
orthornormal vectors in $H,$ $F$ a finite part of $I,$ $\phi _{i},$ $\Phi
_{i},$ $i\in I$ real or complex numbers. For $x\in H,$ if either

\begin{enumerate}
\item[(i)] $\func{Re}\left\langle \sum_{i\in F}\Phi _{i}e_{i}-x,x-\sum_{i\in
F}\phi _{i}e_{i}\right\rangle \geq 0;$\newline

or equivalently,

\item[(ii)] $\left\Vert x-\sum_{i\in F}\frac{\phi _{i}+\Phi _{i}}{2}%
e_{i}\right\Vert \leq \frac{1}{2}\left( \sum_{i\in F}\left\vert \Phi
_{i}-\phi _{i}\right\vert ^{2}\right) ^{\frac{1}{2}};$
\end{enumerate}

holds, then the following reverse of Bessel's inequality%
\begin{align}
0& \leq \left\Vert x\right\Vert ^{2}-\sum_{i\in F}\left\vert \left\langle
x,e_{i}\right\rangle \right\vert ^{2}  \label{2.4.9} \\
& \leq \frac{1}{4}\sum_{i\in F}\left\vert \Phi _{i}-\phi _{i}\right\vert
^{2}-\sum_{i\in F}\left\vert \frac{\phi _{i}+\Phi _{i}}{2}-\left\langle
x,e_{i}\right\rangle \right\vert ^{2}  \notag \\
& \leq \frac{1}{4}\sum_{i\in F}\left\vert \Phi _{i}-\phi _{i}\right\vert
^{2},  \notag
\end{align}%
is valid.

The constant $\frac{1}{4}$ is best possible in both inequalities.
\end{theorem}

\begin{proof}
If we apply Lemma \ref{l2.1.9} for $\lambda _{i}=\frac{\phi _{i}+\Phi _{i}}{2%
}$ and 
\begin{equation*}
r:=\frac{1}{2}\left( \sum_{i\in F}\left\vert \Phi _{i}-\phi _{i}\right\vert
^{2}\right) ^{\frac{1}{2}},
\end{equation*}%
we deduce the first inequality in (\ref{2.4.9}).

Let us prove that $\frac{1}{4}$ is best possible in the second inequality in
(\ref{2.4.9}).

Assume that there is a $c>0$ such that%
\begin{equation}
0\leq \left\Vert x\right\Vert ^{2}-\sum_{i\in F}\left\vert \left\langle
x,e_{i}\right\rangle \right\vert ^{2}\leq c\sum_{i\in F}\left\vert \Phi
_{i}-\phi _{i}\right\vert ^{2}-\sum_{i\in F}\left\vert \frac{\phi _{i}+\Phi
_{i}}{2}-\left\langle x,e_{i}\right\rangle \right\vert ^{2},  \label{2.7.9}
\end{equation}%
provided that $\phi _{i},$ $\Phi _{i},x$ and $F$ satisfy (i) and (ii).

We choose $F=\left\{ 1\right\} ,$ $e_{1}=e=\left( \frac{1}{\sqrt{2}},\frac{1%
}{\sqrt{2}}\right) \in \mathbb{R}^{2},$ $x=\left( x_{1},x_{2}\right) \in 
\mathbb{R}^{2},$ $\Phi _{1}=\Phi =m>0,$ $\phi _{1}=\phi =-m,$ $H=\mathbb{R}%
^{2}$ to get from (\ref{2.7.9}) that%
\begin{align}
0& \leq x_{1}^{2}+x_{2}^{2}-\frac{\left( x_{1}+x_{2}\right) ^{2}}{2}
\label{2.8.9} \\
& \leq 4cm^{2}-\frac{\left( x_{1}+x_{2}\right) ^{2}}{2},  \notag
\end{align}%
provided%
\begin{align}
0& \leq \left\langle me-x,x+me\right\rangle  \label{2.9.9a} \\
& =\left( \frac{m}{\sqrt{2}}-x_{1}\right) \left( x_{1}+\frac{m}{\sqrt{2}}%
\right) +\left( \frac{m}{\sqrt{2}}-x_{2}\right) \left( x_{2}+\frac{m}{\sqrt{2%
}}\right) .  \notag
\end{align}%
From (\ref{2.8.9}) we get%
\begin{equation}
x_{1}^{2}+x_{2}^{2}\leq 4cm^{2}  \label{2.10.9}
\end{equation}%
provided (\ref{2.9.9a}) holds.

If we choose $x_{1}=\frac{m}{\sqrt{2}},$ $x_{2}=-\frac{m}{\sqrt{2}},$ then (%
\ref{2.9.9a}) is fulfilled and by (\ref{2.10.9}) we get $m^{2}\leq 4cm^{2},$
giving $c\geq \frac{1}{4}.$
\end{proof}

\begin{remark}
\label{r2.3.9}If $F=\left\{ 1\right\} ,$ $e_{1}=1,$ $\left\Vert e\right\Vert
=1$ and for $\phi ,\Phi \in \mathbb{K}$ and $x\in H$ one has either%
\begin{equation*}
\func{Re}\left\langle \Phi e-x,x-\phi e\right\rangle \geq 0
\end{equation*}%
or equivalently,%
\begin{equation*}
\left\Vert x-\frac{\phi +\Phi }{2}e\right\Vert \leq \frac{1}{2}\left\vert
\Phi -\phi \right\vert ,
\end{equation*}%
then%
\begin{align*}
0& \leq \left\Vert x\right\Vert ^{2}-\left\vert \left\langle
x,e\right\rangle \right\vert ^{2} \\
& \leq \frac{1}{4}\left\vert \Phi -\phi \right\vert ^{2}-\left\vert \frac{%
\phi +\Phi }{2}-\left\langle x,e\right\rangle \right\vert ^{2}\leq \frac{1}{4%
}\left\vert \Phi -\phi \right\vert ^{2}.
\end{align*}%
The constant $\frac{1}{4}$ is best possible in both inequalities.
\end{remark}

\begin{remark}
\label{r2.4.9}It is important to compare the bounds provided by Theorem \ref%
{t2.2.8} and Theorem \ref{t2.2.9}.

For this purpose, consider%
\begin{equation*}
B_{1}\left( x,e,\phi ,\Phi \right) :=\frac{1}{4}\left( \Phi -\phi \right)
^{2}-\left\langle \Phi e-x,x-\phi e\right\rangle
\end{equation*}%
and%
\begin{equation*}
B_{2}\left( x,e,\phi ,\Phi \right) :=\frac{1}{4}\left( \Phi -\phi \right)
^{2}-\left( \frac{\phi +\Phi }{2}-\left\langle x,e\right\rangle \right) ^{2},
\end{equation*}%
where $H$ is a real inner product, $e\in H,$ $\left\Vert e\right\Vert =1,$ $%
x\in H,$ $\phi ,\Phi \in \mathbb{R}$ with 
\begin{equation*}
\left\langle \Phi e-x,x-\phi e\right\rangle \geq 0,
\end{equation*}%
or equivalently, 
\begin{equation*}
\left\Vert x-\frac{\phi +\Phi }{2}e\right\Vert \leq \frac{1}{2}\left\vert
\Phi -\phi \right\vert .
\end{equation*}

If we choose $\phi =-1,$ $\Phi =1,$ then we have%
\begin{align*}
B_{1}\left( x,e\right) & =1-\left\langle e-x,x+e\right\rangle =1-\left(
\left\Vert e\right\Vert ^{2}-\left\Vert x\right\Vert ^{2}\right) =\left\Vert
x\right\Vert ^{2}, \\
B_{2}\left( x,e\right) & =1-\left\langle x,e\right\rangle ^{2},
\end{align*}%
provided $\left\Vert x\right\Vert \leq 1.$

Consider $H=\mathbb{R}^{2},$ $\left\langle \mathbf{x},\mathbf{y}%
\right\rangle =x_{1}y_{1}+x_{2}y_{2},\mathbf{x=}\left( x_{1},x_{2}\right) ,%
\mathbf{y=}\left( y_{1},y_{2}\right) \in \mathbb{R}^{2}$ and $\mathbf{e=}%
\left( \frac{\sqrt{2}}{2},\frac{\sqrt{2}}{2}\right) .$ Then $\left\Vert 
\mathbf{e}\right\Vert =1$ and we must compare%
\begin{equation*}
B_{1}\left( \mathbf{x}\right) =x_{1}^{2}+x_{2}^{2}
\end{equation*}%
with%
\begin{equation*}
B_{2}\left( \mathbf{x}\right) =1-\frac{\left( x_{1}+x_{2}\right) ^{2}}{2},
\end{equation*}%
provided $x_{1}^{2}+x_{2}^{2}\leq 1.$

If we choose $\mathbf{x}_{0}=\left( 1,0\right) ,$ then $\left\Vert \mathbf{x}%
_{0}\right\Vert =1$ and $B_{1}\left( \mathbf{x}_{0}\right) =1,B_{2}\left( 
\mathbf{x}_{0}\right) =\frac{1}{2}$ showing that $B_{1}>B_{2}.$ If we choose 
$\mathbf{x}_{00}=\left( -\frac{1}{2},\frac{1}{2}\right) ,$ then $B_{1}\left( 
\mathbf{x}_{00}\right) =\frac{1}{2},$ $B_{2}\left( \mathbf{x}_{00}\right)
=1, $ showing that $B_{1}<B_{2}.$
\end{remark}

We may state the following proposition.

\begin{proposition}
\label{p2.5.9}Let $\left\{ e_{i}\right\} _{i\in I}$ be a family of
orthornormal vectors in $H,$ $F$ a finite part of $I,$ $\phi _{i},$ $\Phi
_{i}\in \mathbb{K}$ $\left( i\in F\right) $. If $x\in H$ either satisfies
(i), or equivalently, (ii) of Theorem \ref{t2.2.9}, then the upper bounds%
\begin{align*}
B_{1}\left( x,e,\mathbf{\phi },\mathbf{\Phi },F\right) & :=\frac{1}{4}%
\sum_{i\in F}\left\vert \Phi _{i}-\phi _{i}\right\vert ^{2}-\func{Re}%
\left\langle \sum_{i\in F}\Phi _{i}e_{i}-x,x-\sum_{i\in F}\phi
_{i}e_{i}\right\rangle , \\
B_{2}\left( x,e,\mathbf{\phi },\mathbf{\Phi },F\right) & :=\frac{1}{4}%
\sum_{i\in F}\left\vert \Phi _{i}-\phi _{i}\right\vert ^{2}-\sum_{i\in
F}\left\vert \frac{\phi _{i}+\Phi _{i}}{2}-\left\langle x,e_{i}\right\rangle
\right\vert ^{2},
\end{align*}%
for the Bessel's difference $B_{s}\left( x,\mathbf{e},F\right) :=\left\Vert
x\right\Vert ^{2}-\sum_{i\in F}\left\vert \left\langle x,e_{i}\right\rangle
\right\vert ^{2},$ cannot be compared in general.
\end{proposition}

\subsection{A Refinement of the Gr\"{u}ss Inequality for Orthonormal Families%
}

The following result holds \cite{10NSSD}.

\begin{theorem}
\label{t3.1.9}Let $\left\{ e_{i}\right\} _{i\in I}$ be a family of
orthornormal vectors in $H,$ $F$ a finite part of $I$, $\phi _{i},\Phi _{i},$
$\gamma _{i},\Gamma _{i}\in \mathbb{K},\ i\in F$ and $x,y\in H.$ If either%
\begin{align*}
\func{Re}\left\langle \sum_{i\in F}\Phi _{i}e_{i}-x,x-\sum_{i\in F}\phi
_{i}e_{i}\right\rangle & \geq 0, \\
\func{Re}\left\langle \sum_{i\in F}\Gamma _{i}e_{i}-y,y-\sum_{i\in F}\gamma
_{i}e_{i}\right\rangle & \geq 0,
\end{align*}%
or equivalently,%
\begin{align*}
\left\Vert x-\sum_{i\in F}\frac{\Phi _{i}+\phi _{i}}{2}e_{i}\right\Vert &
\leq \frac{1}{2}\left( \sum_{i\in F}\left\vert \Phi _{i}-\phi
_{i}\right\vert ^{2}\right) ^{\frac{1}{2}}, \\
\left\Vert y-\sum_{i\in F}\frac{\Gamma _{i}+\gamma _{i}}{2}e_{i}\right\Vert
& \leq \frac{1}{2}\left( \sum_{i\in F}\left\vert \Gamma _{i}-\gamma
_{i}\right\vert ^{2}\right) ^{\frac{1}{2}},
\end{align*}%
hold, then we have the inequalities%
\begin{align}
0& \leq \left\vert \left\langle x,y\right\rangle -\sum_{i\in F}\left\langle
x,e_{i}\right\rangle \left\langle e_{i},y\right\rangle \right\vert
\label{3.3.9} \\
& \leq \frac{1}{4}\left( \sum_{i\in F}\left\vert \Phi _{i}-\phi
_{i}\right\vert ^{2}\right) ^{\frac{1}{2}}\cdot \left( \sum_{i\in
F}\left\vert \Gamma _{i}-\gamma _{i}\right\vert ^{2}\right) ^{\frac{1}{2}} 
\notag \\
& \ \ \ \ \ \ \ \ \ \ \ \ \ \ \ \ \ -\sum_{i\in F}\left\vert \frac{\Phi
_{i}+\phi _{i}}{2}-\left\langle x,e_{i}\right\rangle \right\vert \left\vert 
\frac{\Gamma _{i}+\gamma _{i}}{2}-\left\langle y,e_{i}\right\rangle
\right\vert  \notag \\
& \leq \frac{1}{4}\left( \sum_{i\in F}\left\vert \Phi _{i}-\phi
_{i}\right\vert ^{2}\right) ^{\frac{1}{2}}\cdot \left( \sum_{i\in
F}\left\vert \Gamma _{i}-\gamma _{i}\right\vert ^{2}\right) ^{\frac{1}{2}}. 
\notag
\end{align}%
The constant $\frac{1}{4}$ is best possible.
\end{theorem}

\begin{proof}
Using Schwartz's inequality in the inner product space $\left(
H,\left\langle \cdot ,\cdot \right\rangle \right) $ one has%
\begin{multline}
\left\vert \left\langle x-\sum_{i\in F}\left\langle x,e_{i}\right\rangle
e_{i},y-\sum_{i\in F}\left\langle y,e_{i}\right\rangle e_{i}\right\rangle
\right\vert ^{2}  \label{3.4.9} \\
\leq \left\Vert x-\sum_{i\in F}\left\langle x,e_{i}\right\rangle
e_{i}\right\Vert ^{2}\left\Vert y-\sum_{i\in F}\left\langle
y,e_{i}\right\rangle e_{i}\right\Vert ^{2}
\end{multline}%
and since a simple calculation shows that 
\begin{equation*}
\left\langle x-\sum_{i\in F}\left\langle x,e_{i}\right\rangle
e_{i},y-\sum_{i\in F}\left\langle y,e_{i}\right\rangle e_{i}\right\rangle
=\left\langle x,y\right\rangle -\sum_{i\in F}\left\langle
x,e_{i}\right\rangle \left\langle e_{i},y\right\rangle
\end{equation*}%
and 
\begin{equation*}
\left\Vert x-\sum_{i\in F}\left\langle x,e_{i}\right\rangle e_{i}\right\Vert
^{2}=\left\Vert x\right\Vert ^{2}-\sum_{i\in F}\left\vert \left\langle
x,e_{i}\right\rangle \right\vert ^{2}
\end{equation*}%
for any $x,y\in H,$ then by (\ref{3.4.9}) and by the reverse of Bessel's
inequality in Theorem \ref{t2.2.9}, we have%
\begin{align}
& \left\vert \left\langle x,y\right\rangle -\sum_{i\in F}\left\langle
x,e_{i}\right\rangle \left\langle e_{i},y\right\rangle \right\vert ^{2}
\label{3.5.9} \\
& \leq \left( \left\Vert x\right\Vert ^{2}-\sum_{i\in F}\left\vert
\left\langle x,e_{i}\right\rangle \right\vert ^{2}\right) \left( \left\Vert
y\right\Vert ^{2}-\sum_{i\in F}\left\vert \left\langle y,e_{i}\right\rangle
\right\vert ^{2}\right)  \notag \\
& \leq \left[ \frac{1}{4}\sum_{i\in F}\left\vert \Phi _{i}-\phi
_{i}\right\vert ^{2}-\sum_{i\in F}\left\vert \frac{\Phi _{i}+\phi _{i}}{2}%
-\left\langle x,e_{i}\right\rangle \right\vert ^{2}\right]  \notag \\
& \ \ \ \ \ \ \ \ \ \ \ \ \ \ \ \ \times \left[ \frac{1}{4}\sum_{i\in
F}\left\vert \Gamma _{i}-\gamma _{i}\right\vert ^{2}-\sum_{i\in F}\left\vert 
\frac{\Gamma _{i}+\gamma _{i}}{2}-\left\langle y,e_{i}\right\rangle
\right\vert ^{2}\right]  \notag \\
& :=K.  \notag
\end{align}%
Using Acz\'{e}l's inequality for real numbers, i.e., we recall that%
\begin{equation}
\left( a^{2}-\sum_{i\in F}a_{i}^{2}\right) \left( b^{2}-\sum_{i\in
F}b_{i}^{2}\right) \leq \left( ab-\sum_{i\in F}a_{i}b_{i}\right) ^{2},
\label{3.6.9}
\end{equation}%
provided that $a,b,a_{i},b_{i}>0,$ $i\in F,$ we may state that%
\begin{multline}
K\leq \left[ \frac{1}{4}\left( \sum_{i\in F}\left\vert \Phi _{i}-\phi
_{i}\right\vert ^{2}\right) ^{\frac{1}{2}}\cdot \left( \sum_{i\in
F}\left\vert \Gamma _{i}-\gamma _{i}\right\vert ^{2}\right) ^{\frac{1}{2}%
}\right.  \label{3.7.9} \\
\left. -\sum_{i\in F}\left\vert \frac{\Phi _{i}+\phi _{i}}{2}-\left\langle
x,e_{i}\right\rangle \right\vert \left\vert \frac{\Gamma _{i}+\gamma _{i}}{2}%
-\left\langle y,e_{i}\right\rangle \right\vert \right] ^{2}.
\end{multline}%
Using (\ref{3.5.9}) and (\ref{3.7.9}) we conclude that%
\begin{multline}
\left\vert \left\langle x,y\right\rangle -\sum_{i\in F}\left\langle
x,e_{i}\right\rangle \left\langle e_{i},y\right\rangle \right\vert ^{2}\leq 
\left[ \frac{1}{4}\left( \sum_{i\in F}\left\vert \Phi _{i}-\phi
_{i}\right\vert ^{2}\right) ^{\frac{1}{2}}\cdot \left( \sum_{i\in
F}\left\vert \Gamma _{i}-\gamma _{i}\right\vert ^{2}\right) ^{\frac{1}{2}%
}\right.  \label{3.8.9} \\
-\left. \sum_{i\in F}\left\vert \frac{\Phi _{i}+\phi _{i}}{2}-\left\langle
x,e_{i}\right\rangle \right\vert \left\vert \frac{\Gamma _{i}+\gamma _{i}}{2}%
-\left\langle y,e_{i}\right\rangle \right\vert \right] ^{2}.
\end{multline}%
Taking the square root in (\ref{3.8.9}) and taking into account that the
quantity in the last square bracket is nonnegative (this follows by (\ref%
{2.4.9}) and by the Cauchy-Bunyakovsky-Schwarz inequality), we deduce the
second inequality in (\ref{3.3.9}).

The fact that $\frac{1}{4}$ is the best possible constant follows by Theorem %
\ref{t2.2.9} and we omit the details.
\end{proof}

The following corollary may be stated \cite{10NSSD}.

\begin{corollary}
\label{c3.2.9}Let $e\in H,$ $\left\Vert e\right\Vert =1,$ $\phi ,\Phi
,\gamma ,\Gamma \in \mathbb{K}$ and $x,y\in H$ are such that either%
\begin{equation*}
\func{Re}\left\langle \Phi e-x,x-\phi e\right\rangle \geq 0\text{ \ and \ }%
\func{Re}\left\langle \Gamma e-y,y-\gamma e\right\rangle \geq 0
\end{equation*}%
or equivalently,%
\begin{equation*}
\left\Vert x-\frac{\phi +\Phi }{2}e\right\Vert \leq \frac{1}{2}\left\vert
\Phi -\phi \right\vert ,\ \ \ \left\Vert y-\frac{\gamma +\Gamma }{2}%
e\right\Vert \leq \frac{1}{2}\left\vert \Gamma -\gamma \right\vert ,
\end{equation*}%
hold. Then we have the following refinement of Gr\"{u}ss' inequality%
\begin{align*}
0& \leq \left\vert \left\langle x,y\right\rangle -\left\langle
x,e\right\rangle \left\langle e,y\right\rangle \right\vert \\
& \leq \frac{1}{4}\left\vert \Phi -\phi \right\vert \left\vert \Gamma
-\gamma \right\vert -\left\vert \frac{\phi +\Phi }{2}-\left\langle
x,e\right\rangle \right\vert \left\vert \frac{\gamma +\Gamma }{2}%
-\left\langle y,e\right\rangle \right\vert \\
& \leq \frac{1}{4}\left\vert \Phi -\phi \right\vert \left\vert \Gamma
-\gamma \right\vert .
\end{align*}%
The constant $\frac{1}{4}$ is best possible in both inequalities.
\end{corollary}

\subsection{Some Companion Inequalities}

The following companion of the Gr\"{u}ss inequality also holds \cite{10NSSD}.

\begin{theorem}
\label{t4.1.9}Let $\left\{ e_{i}\right\} _{i\in I}$ be a family of
orthornormal vectors in $H,$ $F$ a finite part of $I$ and $\phi _{i},\Phi
_{i}\in \mathbb{K},\ i\in F$, $x,y\in H$ and $\lambda \in \left( 0,1\right)
, $ such that either%
\begin{equation}
\func{Re}\left\langle \sum_{i\in F}\Phi _{i}e_{i}-\left( \lambda x+\left(
1-\lambda \right) y\right) ,\lambda x+\left( 1-\lambda \right) y-\sum_{i\in
F}\phi _{i}e_{i}\right\rangle \geq 0  \label{4.1.9}
\end{equation}%
or equivalently,%
\begin{equation*}
\left\Vert \lambda x+\left( 1-\lambda \right) y-\sum_{i\in F}\frac{\Phi
_{i}+\phi _{i}}{2}\cdot e_{i}\right\Vert \leq \frac{1}{2}\left( \sum_{i\in
F}\left\vert \Phi _{i}-\phi _{i}\right\vert ^{2}\right) ^{\frac{1}{2}},
\end{equation*}%
holds. Then we have the inequality%
\begin{align}
& \func{Re}\left[ \left\langle x,y\right\rangle -\sum_{i\in F}\left\langle
x,e_{i}\right\rangle \left\langle e_{i},y\right\rangle \right]  \label{4.3.9}
\\
& \leq \frac{1}{16}\cdot \frac{1}{\lambda \left( 1-\lambda \right) }%
\sum_{i\in F}\left\vert \Phi _{i}-\phi _{i}\right\vert ^{2}  \notag \\
& \ \ \ \ \ \ \ \ \ \ \ \ \ \ -\frac{1}{4}\frac{1}{\lambda \left( 1-\lambda
\right) }\sum_{i\in F}\left\vert \frac{\Phi _{i}+\phi _{i}}{2}-\left\langle
\lambda x+\left( 1-\lambda \right) y,e_{i}\right\rangle \right\vert ^{2} 
\notag \\
& \leq \frac{1}{16}\cdot \frac{1}{\lambda \left( 1-\lambda \right) }%
\sum_{i\in F}\left\vert \Phi _{i}-\phi _{i}\right\vert ^{2}.  \notag
\end{align}%
The constant $\frac{1}{16}$ is the best possible constant in (\ref{4.3.9})
in the sense that it cannot be replaced by a smaller constant.
\end{theorem}

\begin{proof}
We know that for any$\ z,u\in H,$ one has%
\begin{equation*}
\func{Re}\left\langle z,u\right\rangle \leq \frac{1}{4}\left\Vert
z+u\right\Vert ^{2}.
\end{equation*}%
Then for any $a,b\in H$ and $\lambda \in \left( 0,1\right) $ one has%
\begin{equation}
\func{Re}\left\langle a,b\right\rangle \leq \frac{1}{4\lambda \left(
1-\lambda \right) }\left\Vert \lambda a+\left( 1-\lambda \right)
b\right\Vert ^{2}.  \label{4.4.9}
\end{equation}%
Since%
\begin{equation*}
\left\langle x,y\right\rangle -\sum_{i\in F}\left\langle
x,e_{i}\right\rangle \left\langle e_{i},y\right\rangle =\left\langle
x-\sum_{i\in F}\left\langle x,e_{i}\right\rangle e_{i},y-\sum_{i\in
F}\left\langle y,e_{i}\right\rangle e_{i}\right\rangle ,
\end{equation*}%
for any \thinspace $x,y\in H,$ then, by (\ref{4.4.9}), we get%
\begin{align}
& \func{Re}\left[ \left\langle x,y\right\rangle -\sum_{i\in F}\left\langle
x,e_{i}\right\rangle \left\langle e_{i},y\right\rangle \right]  \label{4.5.9}
\\
& =\func{Re}\left[ \left\langle x-\sum_{i\in F}\left\langle
x,e_{i}\right\rangle e_{i},y-\sum_{i\in F}\left\langle y,e_{i}\right\rangle
e_{i}\right\rangle \right]  \notag \\
& \leq \frac{1}{4\lambda \left( 1-\lambda \right) }\left\Vert \lambda \left(
x-\sum_{i\in F}\left\langle x,e_{i}\right\rangle e_{i}\right) +\left(
1-\lambda \right) \left( y-\sum_{i\in F}\left\langle y,e_{i}\right\rangle
e_{i}\right) \right\Vert ^{2}  \notag \\
& =\frac{1}{4\lambda \left( 1-\lambda \right) }\left\Vert \lambda x+\left(
1-\lambda \right) y-\sum_{i\in F}\left\langle \lambda x+\left( 1-\lambda
\right) y,e_{i}\right\rangle e_{i}\right\Vert ^{2}  \notag \\
& =\frac{1}{4\lambda \left( 1-\lambda \right) }\left[ \left\Vert \lambda
x+\left( 1-\lambda \right) y\right\Vert ^{2}-\sum_{i\in F}\left\vert
\left\langle \lambda x+\left( 1-\lambda \right) y,e_{i}\right\rangle
\right\vert ^{2}\right] .  \notag
\end{align}%
If we apply the reverse of Bessel's inequality in Theorem \ref{t2.2.9} for $%
\lambda x+\left( 1-\lambda \right) y,$ we may state that%
\begin{align}
& \left\Vert \lambda x+\left( 1-\lambda \right) y\right\Vert ^{2}-\sum_{i\in
F}\left\vert \left\langle \lambda x+\left( 1-\lambda \right)
y,e_{i}\right\rangle \right\vert ^{2}  \label{4.6.9} \\
& \leq \frac{1}{4}\sum_{i\in F}\left\vert \Phi _{i}-\phi _{i}\right\vert
^{2}-\sum_{i\in F}\left\vert \frac{\Phi _{i}+\phi _{i}}{2}-\left\langle
\lambda x+\left( 1-\lambda \right) y,e_{i}\right\rangle \right\vert ^{2} 
\notag \\
& \leq \frac{1}{4}\sum_{i\in F}\left\vert \Phi _{i}-\phi _{i}\right\vert
^{2}.  \notag
\end{align}%
Now, by making use of (\ref{4.5.9}) and (\ref{4.6.9}), we deduce (\ref{4.3.9}%
).

The fact that $\frac{1}{16}$ is the best possible constant in (\ref{4.3.9})
follows by the fact that if in (\ref{4.1.9}) we choose $x=y,$ then it
becomes (i) of Theorem \ref{t2.2.9}, implying for $\lambda =\frac{1}{2}$ (%
\ref{2.4.9}), for which, we have shown that $\frac{1}{4}$ was the best
constant.
\end{proof}

\begin{remark}
\label{r1.9}In practical applications we may use only the inequality between
the first and the last terms in (\ref{4.3.9}).
\end{remark}

\begin{remark}
\label{r2.9}If in Theorem \ref{t4.1.9}, we choose $\lambda =\frac{1}{2},$
then we get%
\begin{align*}
& \func{Re}\left[ \left\langle x,y\right\rangle -\sum_{i\in F}\left\langle
x,e_{i}\right\rangle \left\langle e_{i},y\right\rangle \right] \\
& \leq \frac{1}{4}\sum_{i\in F}\left\vert \Phi _{i}-\phi _{i}\right\vert
^{2}-\sum_{i\in F}\left\vert \frac{\Phi _{i}+\phi _{i}}{2}-\left\langle 
\frac{x+y}{2},e_{i}\right\rangle \right\vert ^{2} \\
& \leq \frac{1}{4}\sum_{i\in F}\left\vert \Phi _{i}-\phi _{i}\right\vert
^{2},
\end{align*}%
provided%
\begin{equation*}
\func{Re}\left\langle \sum_{i\in F}\Phi _{i}e_{i}-\frac{x+y}{2},\frac{x+y}{2}%
-\sum_{i\in F}\phi _{i}e_{i}\right\rangle \geq 0
\end{equation*}%
or equivalently,%
\begin{equation*}
\left\Vert \frac{x+y}{2}-\sum_{i\in F}\frac{\Phi _{i}+\phi _{i}}{2}\cdot
e_{i}\right\Vert \leq \frac{1}{2}\left( \sum_{i\in F}\left\vert \Phi
_{i}-\phi _{i}\right\vert ^{2}\right) ^{\frac{1}{2}}.
\end{equation*}
\end{remark}

\begin{corollary}
\label{c4.2.9}With the assumptions of Theorem \ref{t4.1.9} and if%
\begin{equation*}
\func{Re}\left\langle \sum_{i\in F}\Phi _{i}e_{i}-\left( \lambda x\pm \left(
1-\lambda \right) y\right) ,\lambda x\pm \left( 1-\lambda \right)
y-\sum_{i\in F}\phi _{i}e_{i}\right\rangle \geq 0
\end{equation*}%
or, equivalently%
\begin{equation*}
\left\Vert \lambda x\pm \left( 1-\lambda \right) y-\sum_{i\in F}\frac{\Phi
_{i}+\phi _{i}}{2}\cdot e_{i}\right\Vert \leq \frac{1}{2}\left( \sum_{i\in
F}\left\vert \Phi _{i}-\phi _{i}\right\vert ^{2}\right) ^{\frac{1}{2}},
\end{equation*}%
then we have the inequality%
\begin{equation}
\left\vert \func{Re}\left[ \left\langle x,y\right\rangle -\sum_{i\in
F}\left\langle x,e_{i}\right\rangle \left\langle e_{i},y\right\rangle \right]
\right\vert \leq \frac{1}{16}\cdot \frac{1}{\lambda \left( 1-\lambda \right) 
}\sum_{i\in F}\left\vert \Phi _{i}-\phi _{i}\right\vert ^{2}.  \label{4.10.9}
\end{equation}%
The constant $\frac{1}{16}$ is best possible in (\ref{4.10.9}).
\end{corollary}

\begin{remark}
\label{r4.3.9}If $H$ is a real inner product space and $m_{i},M_{i}\in 
\mathbb{R}$ with the property%
\begin{equation*}
\left\langle \sum_{i\in F}M_{i}e_{i}-\left( \lambda x\pm \left( 1-\lambda
\right) y\right) ,\lambda x\pm \left( 1-\lambda \right) y-\sum_{i\in
F}m_{i}e_{i}\right\rangle \geq 0
\end{equation*}%
or equivalently,%
\begin{equation*}
\left\Vert \lambda x\pm \left( 1-\lambda \right) y-\sum_{i\in F}\frac{%
M_{i}+m_{i}}{2}\cdot e_{i}\right\Vert \leq \frac{1}{2}\left[ \sum_{i\in
F}\left( M_{i}-m_{i}\right) ^{2}\right] ^{\frac{1}{2}},
\end{equation*}%
then we have the Gr\"{u}ss type inequality%
\begin{equation*}
\left\vert \left\langle x,y\right\rangle -\sum_{i\in F}\left\langle
x,e_{i}\right\rangle \left\langle e_{i},y\right\rangle \right\vert \leq 
\frac{1}{16}\cdot \frac{1}{\lambda \left( 1-\lambda \right) }\sum_{i\in
F}\left( M_{i}-m_{i}\right) ^{2}.
\end{equation*}
\end{remark}

\subsection{Integral Inequalities}

The following proposition holds \cite{10NSSD}.

\begin{proposition}
\label{p5.1.9}Let $\left\{ f_{i}\right\} _{i\in I}$ be an orthornormal
family of functions in $L_{\rho }^{2}\left( \Omega ,\mathbb{K}\right) ,$ $F$
a finite subset of $I,$ $\phi _{i},\Phi _{i}\in \mathbb{K}$ $\left( i\in
F\right) $ and $f\in L_{\rho }^{2}\left( \Omega ,\mathbb{K}\right) ,$ so
that either%
\begin{equation*}
\int_{\Omega }\rho \left( s\right) \func{Re}\left[ \left( \sum_{i\in F}\Phi
_{i}f_{i}\left( s\right) -f\left( s\right) \right) \left( \overline{f}\left(
s\right) -\sum_{i\in F}\overline{\phi _{i}}\text{ }\overline{f_{i}}\left(
s\right) \right) \right] d\mu \left( s\right) \geq 0
\end{equation*}%
or equivalently,%
\begin{equation*}
\int_{\Omega }\rho \left( s\right) \left\vert f\left( s\right) -\sum_{i\in F}%
\frac{\Phi _{i}+\phi _{i}}{2}\cdot f_{i}\left( s\right) \right\vert ^{2}d\mu
\left( s\right) \leq \frac{1}{4}\sum_{i\in F}\left\vert \Phi _{i}-\phi
_{i}\right\vert ^{2}.
\end{equation*}%
Then we have the inequality%
\begin{align}
0& \leq \int_{\Omega }\rho \left( s\right) \left\vert f\left( s\right)
\right\vert ^{2}d\mu \left( s\right) -\sum_{i\in F}\left\vert \int_{\Omega
}\rho \left( s\right) f\left( s\right) \overline{f_{i}}\left( s\right) d\mu
\left( s\right) \right\vert ^{2}  \label{5.5.9} \\
& \leq \frac{1}{4}\sum_{i\in F}\left\vert \Phi _{i}-\phi _{i}\right\vert
^{2}-\sum_{i\in F}\left\vert \frac{\Phi _{i}+\phi _{i}}{2}-\int_{\Omega
}\rho \left( s\right) f\left( s\right) \overline{f_{i}}\left( s\right) d\mu
\left( s\right) \right\vert ^{2}  \notag \\
& \leq \frac{1}{4}\sum_{i\in F}\left\vert \Phi _{i}-\phi _{i}\right\vert
^{2}.  \notag
\end{align}%
The constant $\frac{1}{4}$ is best possible in both inequalities.
\end{proposition}

The proof follows by Theorem \ref{t2.2.9} applied for the Hilbert space $%
L_{\rho }^{2}\left( \Omega ,\mathbb{K}\right) $ and the orthornormal family $%
\left\{ f_{i}\right\} _{i\in I}.$

The following Gr\"{u}ss type inequality also holds \cite{10NSSD}.

\begin{proposition}
\label{p5.2.9}Let $\left\{ f_{i}\right\} _{i\in I}$ and $F$ be as in
Proposition \ref{p5.1.9}. If $\phi _{i},\Phi _{i},\gamma _{i},\Gamma _{i}\in 
\mathbb{K}$ $\left( i\in F\right) $ and $f,g\in L_{\rho }^{2}\left( \Omega ,%
\mathbb{K}\right) $ so that either%
\begin{align*}
\int_{\Omega }\rho \left( s\right) \func{Re}\left[ \left( \sum_{i\in F}\Phi
_{i}f_{i}\left( s\right) -f\left( s\right) \right) \left( \overline{f}\left(
s\right) -\sum_{i\in F}\overline{\phi _{i}}\text{ }\overline{f_{i}}\left(
s\right) \right) \right] d\mu \left( s\right) & \geq 0, \\
\int_{\Omega }\rho \left( s\right) \func{Re}\left[ \left( \sum_{i\in
F}\Gamma _{i}f_{i}\left( s\right) -g\left( s\right) \right) \left( \overline{%
g}\left( s\right) -\sum_{i\in F}\overline{\gamma _{i}}\text{ }\overline{f_{i}%
}\left( s\right) \right) \right] d\mu \left( s\right) & \geq 0,
\end{align*}%
or equivalently,%
\begin{align*}
\int_{\Omega }\rho \left( s\right) \left\vert f\left( s\right) -\sum_{i\in F}%
\frac{\Phi _{i}+\phi _{i}}{2}f_{i}\left( s\right) \right\vert ^{2}d\mu
\left( s\right) & \leq \frac{1}{4}\sum_{i\in F}\left\vert \Phi _{i}-\phi
_{i}\right\vert ^{2}, \\
\int_{\Omega }\rho \left( s\right) \left\vert g\left( s\right) -\sum_{i\in F}%
\frac{\Gamma _{i}+\gamma _{i}}{2}f_{i}\left( s\right) \right\vert ^{2}d\mu
\left( s\right) & \leq \frac{1}{4}\sum_{i\in F}\left\vert \Gamma _{i}-\gamma
_{i}\right\vert ^{2},
\end{align*}%
hold, then we have the inequalities%
\begin{align}
& \left\vert \int_{\Omega }\rho \left( s\right) f\left( s\right) \overline{%
g\left( s\right) }d\mu \left( s\right) \right.  \label{5.8.9} \\
& \qquad -\left. \sum_{i\in F}\int_{\Omega }\rho \left( s\right) f\left(
s\right) \overline{f_{i}}\left( s\right) d\mu \left( s\right) \int_{\Omega
}\rho \left( s\right) f_{i}\left( s\right) \overline{g\left( s\right) }d\mu
\left( s\right) \right\vert  \notag \\
& \leq \frac{1}{4}\left( \sum_{i\in F}\left\vert \Phi _{i}-\phi
_{i}\right\vert ^{2}\right) ^{\frac{1}{2}}\left( \sum_{i\in F}\left\vert
\Gamma _{i}-\gamma _{i}\right\vert ^{2}\right) ^{\frac{1}{2}}  \notag \\
& \qquad -\sum_{i\in F}\left\vert \frac{\Phi _{i}+\phi _{i}}{2}-\int_{\Omega
}\rho \left( s\right) f\left( s\right) \overline{f_{i}}\left( s\right) d\mu
\left( s\right) \right\vert  \notag \\
& \qquad \times \left\vert \frac{\Gamma _{i}+\gamma _{i}}{2}-\int_{\Omega
}\rho \left( s\right) g\left( s\right) \overline{f_{i}}\left( s\right) d\mu
\left( s\right) \right\vert  \notag \\
& \leq \frac{1}{4}\left( \sum_{i\in F}\left\vert \Phi _{i}-\phi
_{i}\right\vert ^{2}\right) ^{\frac{1}{2}}\left( \sum_{i\in F}\left\vert
\Gamma _{i}-\gamma _{i}\right\vert ^{2}\right) ^{\frac{1}{2}}.  \notag
\end{align}%
The constant $\frac{1}{4}$ is the best possible.
\end{proposition}

The proof follows by Theorem \ref{t3.1.9} and we omit the details.

\begin{remark}
\label{r5.3.9}Similar results may be stated if one applies the inequalities
in the above subsections. We omit the details.
\end{remark}

In the case of real spaces, the following corollaries provide much simpler
sufficient conditions for the reverse of Bessel's inequality (\ref{5.5.9})
or for the Gr\"{u}ss type inequality (\ref{5.8.9}) to hold.

\begin{corollary}
\label{c5.4.9}Let $\left\{ f_{i}\right\} _{i\in I}$ be an orthornormal
family of functions in the real Hilbert space $L_{\rho }^{2}\left( \Omega
\right) ,$ $F$ a finite part of $I,$ $M_{i},m_{i}\in \mathbb{R}$ \ $\left(
i\in F\right) $ and $f\in L_{\rho }^{2}\left( \Omega \right) $ so that%
\begin{equation*}
\sum_{i\in F}m_{i}f_{i}\left( s\right) \leq f\left( s\right) \leq \sum_{i\in
F}M_{i}f_{i}\left( s\right) \text{ \ for \ }\mu -\text{a.e. }s\in \Omega .
\end{equation*}%
Then we have the inequalities%
\begin{align*}
0& \leq \int_{\Omega }\rho \left( s\right) f^{2}\left( s\right) d\mu \left(
s\right) -\sum_{i\in F}\left[ \int_{\Omega }\rho \left( s\right) f\left(
s\right) f_{i}\left( s\right) d\mu \left( s\right) \right] ^{2} \\
& \leq \frac{1}{4}\sum_{i\in F}\left( M_{i}-m_{i}\right) ^{2}-\sum_{i\in F}%
\left[ \frac{M_{i}+m_{i}}{2}-\int_{\Omega }\rho \left( s\right) f\left(
s\right) f_{i}\left( s\right) d\mu \left( s\right) \right] ^{2} \\
& \leq \frac{1}{4}\sum_{i\in F}\left( M_{i}-m_{i}\right) ^{2}.
\end{align*}%
The constant $\frac{1}{4}$ is best possible.
\end{corollary}

\begin{corollary}
\label{c5.5.9}Let $\left\{ f_{i}\right\} _{i\in I}$ and $F$ be as in
Corollary \ref{c5.4.9}. If $M_{i},m_{i},N_{i},n_{i}\in \mathbb{R}$ $\left(
i\in F\right) $ and $f,g\in L_{\rho }^{2}\left( \Omega \right) $ are such
that%
\begin{equation*}
\sum_{i\in F}m_{i}f_{i}\left( s\right) \leq f\left( s\right) \leq \sum_{i\in
F}M_{i}f_{i}\left( s\right)
\end{equation*}%
and%
\begin{equation*}
\sum_{i\in F}n_{i}f_{i}\left( s\right) \leq g\left( s\right) \leq \sum_{i\in
F}N_{i}f_{i}\left( s\right) \text{ \ for \ }\mu -\text{a.e. }s\in \Omega ,
\end{equation*}%
hold, then we have the inequalities%
\begin{multline*}
\left\vert \int_{\Omega }\rho \left( s\right) f\left( s\right) g\left(
s\right) d\mu \left( s\right) \right. \\
-\left. \sum_{i\in F}\int_{\Omega }\rho \left( s\right) f\left( s\right)
f_{i}\left( s\right) d\mu \left( s\right) \int_{\Omega }\rho \left( s\right)
g\left( s\right) f_{i}\left( s\right) d\mu \left( s\right) \right\vert
\end{multline*}%
\begin{align*}
& \leq \frac{1}{4}\left( \sum_{i\in F}\left( M_{i}-m_{i}\right) ^{2}\right)
^{\frac{1}{2}}\left( \sum_{i\in F}\left( N_{i}-n_{i}\right) ^{2}\right) ^{%
\frac{1}{2}} \\
& \qquad -\sum_{i\in F}\left\vert \frac{M_{i}+m_{i}}{2}-\int_{\Omega }\rho
\left( s\right) f\left( s\right) f_{i}\left( s\right) d\mu \left( s\right)
\right\vert \\
& \qquad \times \left\vert \frac{N_{i}+n_{i}}{2}-\int_{\Omega }\rho \left(
s\right) g\left( s\right) f_{i}\left( s\right) d\mu \left( s\right)
\right\vert \\
& \leq \frac{1}{4}\left( \sum_{i\in F}\left( M_{i}-m_{i}\right) ^{2}\right)
^{\frac{1}{2}}\left( \sum_{i\in F}\left( N_{i}-n_{i}\right) ^{2}\right) ^{%
\frac{1}{2}}.
\end{align*}
\end{corollary}

\newpage

\section{More Reverses of Bessel's Inequality}

\subsection{A General Result}

The following reverse of Bessel's inequality holds \cite{11NSSD}.

\begin{theorem}
\label{t2.1.10}Let $\left\{ e_{i}\right\} _{i\in I}$ be a family of
orthornormal vectors in $H,$ $F$ a finite part of $I,$ and $\phi _{i},\Phi
_{i}$ $\left( i\in F\right) ,$ real or complex numbers such that $\sum_{i\in
F}\func{Re}\left( \Phi _{i}\overline{\phi _{i}}\right) >0.$ If $x\in H$ is
such that either

\begin{enumerate}
\item[(i)] $\func{Re}\left\langle \sum_{i\in F}\Phi _{i}e_{i}-x,x-\sum_{i\in
F}\phi _{i}e_{i}\right\rangle \geq 0;$\newline
or equivalently,

\item[(ii)] $\left\Vert x-\sum_{i\in F}\frac{\phi _{i}+\Phi _{i}}{2}%
e_{i}\right\Vert \leq \frac{1}{2}\left( \sum_{i\in F}\left\vert \Phi
_{i}-\phi _{i}\right\vert ^{2}\right) ^{\frac{1}{2}};$
\end{enumerate}

holds, then one has the inequality 
\begin{equation}
\left\Vert x\right\Vert ^{2}\leq \frac{1}{4}\cdot \frac{\sum_{i\in F}\left(
\left\vert \Phi _{i}\right\vert +\left\vert \phi _{i}\right\vert \right) ^{2}%
}{\sum_{i\in F}\func{Re}\left( \Phi _{i}\overline{\phi _{i}}\right) }%
\sum_{i\in F}\left\vert \left\langle x,e_{i}\right\rangle \right\vert ^{2}.
\label{2.1.10}
\end{equation}%
The constant $\frac{1}{4}$ is best possible in the sense that it cannot be
replaced by a smaller constant.
\end{theorem}

\begin{proof}
Observe that%
\begin{multline*}
\func{Re}\left\langle \sum_{i\in F}\Phi _{i}e_{i}-x,x-\sum_{i\in F}\phi
_{i}e_{i}\right\rangle \\
=\sum_{i\in F}\func{Re}\left[ \Phi _{i}\overline{\left\langle
x,e_{i}\right\rangle }+\overline{\phi _{i}}\left\langle x,e_{i}\right\rangle %
\right] -\left\Vert x\right\Vert ^{2}-\sum_{i\in F}\func{Re}\left( \Phi _{i}%
\overline{\phi _{i}}\right) ,
\end{multline*}%
giving, from (i), that%
\begin{equation}
\left\Vert x\right\Vert ^{2}+\sum_{i\in F}\func{Re}\left( \Phi _{i}\overline{%
\phi _{i}}\right) \leq \sum_{i\in F}\func{Re}\left[ \Phi _{i}\overline{%
\left\langle x,e_{i}\right\rangle }+\overline{\phi _{i}}\left\langle
x,e_{i}\right\rangle \right] .  \label{2.4.10}
\end{equation}

On the other hand, by the elementary inequality%
\begin{equation*}
\alpha p^{2}+\frac{1}{\alpha }q^{2}\geq 2pq,\ \ \alpha >0,\ p,q\geq 0;
\end{equation*}%
we deduce%
\begin{equation}
2\left\Vert x\right\Vert \leq \frac{\left\Vert x\right\Vert ^{2}}{\left[
\sum_{i\in F}\func{Re}\left( \Phi _{i}\overline{\phi _{i}}\right) \right] ^{%
\frac{1}{2}}}+\left[ \sum_{i\in F}\func{Re}\left( \Phi _{i}\overline{\phi
_{i}}\right) \right] ^{\frac{1}{2}}.  \label{2.5.10}
\end{equation}%
Dividing (\ref{2.4.10}) by $\left[ \sum_{i\in F}\func{Re}\left( \Phi _{i}%
\overline{\phi _{i}}\right) \right] ^{\frac{1}{2}}>0$ and using (\ref{2.5.10}%
), we obtain%
\begin{equation}
\left\Vert x\right\Vert \leq \frac{1}{2}\frac{\sum_{i\in F}\func{Re}\left[
\Phi _{i}\overline{\left\langle x,e_{i}\right\rangle }+\overline{\phi _{i}}%
\left\langle x,e_{i}\right\rangle \right] }{\left[ \sum_{i\in F}\func{Re}%
\left( \Phi _{i}\overline{\phi _{i}}\right) \right] ^{\frac{1}{2}}},
\label{2.6.10}
\end{equation}%
which is also an interesting inequality in itself.

Using the Cauchy-Bunyakovsky-Schwarz inequality for real numbers, we get%
\begin{align}
\sum_{i\in F}\func{Re}\left[ \Phi _{i}\overline{\left\langle
x,e_{i}\right\rangle }+\overline{\phi _{i}}\left\langle x,e_{i}\right\rangle %
\right] & \leq \sum_{i\in F}\left\vert \Phi _{i}\overline{\left\langle
x,e_{i}\right\rangle }+\overline{\phi _{i}}\left\langle x,e_{i}\right\rangle
\right\vert  \label{2.7.10} \\
& \leq \sum_{i\in F}\left( \left\vert \Phi _{i}\right\vert +\left\vert \phi
_{i}\right\vert \right) \left\vert \left\langle x,e_{i}\right\rangle
\right\vert  \notag \\
& \leq \left[ \sum_{i\in F}\left( \left\vert \Phi _{i}\right\vert
+\left\vert \phi _{i}\right\vert \right) ^{2}\right] ^{\frac{1}{2}}\left[
\sum_{i\in F}\left\vert \left\langle x,e_{i}\right\rangle \right\vert ^{2}%
\right] ^{\frac{1}{2}}.  \notag
\end{align}%
Making use of (\ref{2.6.10}) and (\ref{2.7.10}), we deduce the desired
result (\ref{2.1.10}).

To prove the sharpness of the constant $\frac{1}{4},$ let us assume that (%
\ref{2.1.10}) holds with a constant $c>0,$ i.e., 
\begin{equation}
\left\Vert x\right\Vert ^{2}\leq c\cdot \frac{\sum_{i\in F}\left( \left\vert
\Phi _{i}\right\vert +\left\vert \phi _{i}\right\vert \right) ^{2}}{%
\sum_{i\in F}\func{Re}\left( \Phi _{i}\overline{\phi _{i}}\right) }%
\sum_{i\in F}\left\vert \left\langle x,e_{i}\right\rangle \right\vert ^{2},
\label{2.8.10}
\end{equation}%
provided $x,$ $\phi _{i},\Phi _{i},i\in F$ satisfies (i).

Choose $F=\left\{ 1\right\} ,$ $e_{1}=e,$ $\left\Vert e\right\Vert =1,$ $%
\phi _{i}=m,$ $\Phi _{i}=M$ with $m,M>0,$ then, by (\ref{2.8.10}), we get%
\begin{equation}
\left\Vert x\right\Vert ^{2}\leq c\frac{\left( M+m\right) ^{2}}{mM}%
\left\vert \left\langle x,e\right\rangle \right\vert ^{2}  \label{2.9.10}
\end{equation}%
provided%
\begin{equation}
\func{Re}\left\langle Me-x,x-me\right\rangle \geq 0.  \label{2.10.10}
\end{equation}%
If $x=me,$ then obviously (\ref{2.10.10}) holds, and by (\ref{2.9.10}) we get%
\begin{equation*}
m^{2}\leq c\frac{\left( M+m\right) ^{2}}{mM}m^{2}
\end{equation*}%
giving $mM\leq c\left( M+m\right) ^{2}$ for $m,M>0.$ Now, if in this
inequality we choose $m=1-\varepsilon ,$ $M=1+\varepsilon $ $\left(
\varepsilon \in \left( 0,1\right) \right) ,$ then we get $1-\varepsilon
^{2}\leq 4c$ for $\varepsilon \in \left( 0,1\right) ,$ from where we deduce $%
c\geq \frac{1}{4}.$
\end{proof}

\begin{remark}
\label{r2.2.10}By the use of (\ref{2.6.10}), the second inequality in (\ref%
{2.7.10}) and the H\"{o}lder inequality, we may state the following reverses
of Bessel's inequality as well:%
\begin{multline*}
\left\Vert x\right\Vert ^{2}\leq \frac{1}{2}\cdot \frac{1}{\left[ \sum_{i\in
F}\func{Re}\left( \Phi _{i}\overline{\phi _{i}}\right) \right] ^{\frac{1}{2}}%
} \\
\times \left\{ 
\begin{array}{l}
\max\limits_{i\in F}\left\{ \left\vert \Phi _{i}\right\vert +\left\vert \phi
_{i}\right\vert \right\} \sum\limits_{i\in F}\left\vert \left\langle
x,e_{i}\right\rangle \right\vert ; \\ 
\\ 
\left[ \sum\limits_{i\in F}\left( \left\vert \Phi _{i}\right\vert
+\left\vert \phi _{i}\right\vert \right) ^{p}\right] ^{\frac{1}{p}}\left(
\sum\limits_{i\in F}\left\vert \left\langle x,e_{i}\right\rangle \right\vert
^{q}\right) ^{\frac{1}{q}},\text{ } \\ 
\hfill \text{for\ }p>1,\ \frac{1}{p}+\frac{1}{q}=1; \\ 
\\ 
\max_{i\in F}\left\vert \left\langle x,e_{i}\right\rangle \right\vert
\sum\limits_{i\in F}\left[ \left\vert \Phi _{i}\right\vert +\left\vert \phi
_{i}\right\vert \right] .%
\end{array}%
\right.
\end{multline*}
\end{remark}

The following corollary holds \cite{11NSSD}.

\begin{corollary}
\label{c2.3.10}With the assumption of Theorem \ref{t2.1.10} and if either
(i) or (ii) holds, then%
\begin{equation}
0\leq \left\Vert x\right\Vert ^{2}-\sum_{i\in F}\left\vert \left\langle
x,e_{i}\right\rangle \right\vert ^{2}\leq \frac{1}{4}M^{2}\left( \mathbf{%
\Phi },\mathbf{\phi },F\right) \sum_{i\in F}\left\vert \left\langle
x,e_{i}\right\rangle \right\vert ^{2},  \label{2.12.10}
\end{equation}%
where%
\begin{equation*}
M\left( \mathbf{\Phi },\mathbf{\phi },F\right) :=\left[ \frac{\sum_{i\in
F}\left\{ \left( \left\vert \Phi _{i}\right\vert -\left\vert \phi
_{i}\right\vert \right) ^{2}+4\left[ \left\vert \Phi _{i}\overline{\phi _{i}}%
\right\vert -\func{Re}\left( \Phi _{i}\overline{\phi _{i}}\right) \right]
\right\} }{\sum_{i\in F}\func{Re}\left( \Phi _{i}\overline{\phi _{i}}\right) 
}\right] ^{\frac{1}{2}}.
\end{equation*}%
The constant $\frac{1}{4}$ is best possible.
\end{corollary}

\begin{proof}
The inequality (\ref{2.12.10}) follows by (\ref{2.1.10}) on subtracting the
same quantity $\sum_{i\in F}\left\vert \left\langle x,e_{i}\right\rangle
\right\vert ^{2}$ from both sides.

To prove the sharpness of the constant $\frac{1}{4},$ assume that (\ref%
{2.12.10}) holds with $c>0,$ i.e., 
\begin{equation}
0\leq \left\Vert x\right\Vert ^{2}-\sum_{i\in F}\left\vert \left\langle
x,e_{i}\right\rangle \right\vert ^{2}\leq cM^{2}\left( \mathbf{\Phi },%
\mathbf{\phi },F\right) \sum_{i\in F}\left\vert \left\langle
x,e_{i}\right\rangle \right\vert ^{2}  \label{2.13.10}
\end{equation}%
provided the condition (i) holds.

Choose $F=\left\{ 1\right\} ,$ $e_{1}=e,$ $\left\Vert e\right\Vert =1,$ $%
\phi _{i}=\phi ,$ $\Phi _{i}=\Phi ,$ $\phi ,\Phi >0$ in (\ref{2.13.10}) to
get%
\begin{equation}
0\leq \left\Vert x\right\Vert ^{2}-\left\vert \left\langle x,e\right\rangle
\right\vert ^{2}\leq c\frac{\left( \Phi -\phi \right) ^{2}}{\phi \Phi }%
\left\vert \left\langle x,e\right\rangle \right\vert ^{2},  \label{2.14.10}
\end{equation}%
provided%
\begin{equation}
\left\langle \Phi e-x,x-\phi e\right\rangle \geq 0.  \label{2.14.1.10}
\end{equation}%
If $H=\mathbb{R}^{2},$ $x=\left( x_{1},x_{2}\right) \in \mathbb{R}^{2},$ $%
e=\left( \frac{1}{\sqrt{2}},\frac{1}{\sqrt{2}}\right) $ then we have%
\begin{align*}
\left\Vert x\right\Vert ^{2}-\left\vert \left\langle x,e\right\rangle
\right\vert ^{2}& =x_{1}^{2}+x_{2}^{2}-\frac{\left( x_{1}+x_{2}\right) ^{2}}{%
2}=\frac{1}{2}\left( x_{1}-x_{2}\right) ^{2}, \\
\left\vert \left\langle x,e\right\rangle \right\vert ^{2}& =\frac{\left(
x_{1}+x_{2}\right) ^{2}}{2}
\end{align*}%
and by (\ref{2.14.10}) we get%
\begin{equation}
\frac{\left( x_{1}-x_{2}\right) ^{2}}{2}\leq c\frac{\left( \Phi -\phi
\right) ^{2}}{\phi \Phi }\cdot \frac{\left( x_{1}+x_{2}\right) ^{2}}{2}.
\label{2.15.10}
\end{equation}%
Now, if we let $x_{1}=\frac{\phi }{\sqrt{2}},$ $x_{2}=\frac{\Phi }{\sqrt{2}}$
$\left( \phi ,\Phi >0\right) $ then obviously%
\begin{equation*}
\left\langle \Phi e-x,x-\phi e\right\rangle =\sum_{i=1}^{2}\left( \frac{\Phi 
}{\sqrt{2}}-x_{i}\right) \left( x_{i}-\frac{\phi }{\sqrt{2}}\right) =0,
\end{equation*}%
which shows that (\ref{2.14.1.10}) is fulfilled, and thus by (\ref{2.15.10})
we obtain%
\begin{equation*}
\frac{\left( \Phi -\phi \right) ^{2}}{4}\leq c\frac{\left( \Phi -\phi
\right) ^{2}}{\phi \Phi }\cdot \frac{\left( \Phi +\phi \right) ^{2}}{4}
\end{equation*}%
for any $\Phi >\phi >0.$ This implies%
\begin{equation}
c\left( \Phi +\phi \right) ^{2}\geq \phi \Phi  \label{2.16.10}
\end{equation}%
for any $\Phi >\phi >0.$

Finally, let $\phi =1-\varepsilon ,$ $\Phi =1+\varepsilon $, $\varepsilon
\in \left( 0,1\right) $. Then from (\ref{2.16.10}) we get $4c\geq
1-\varepsilon ^{2}$ for any $\varepsilon \in \left( 0,1\right) $ which
produces $c\geq \frac{1}{4}.$
\end{proof}

\begin{remark}
\label{r2.4.10}If $\left\{ e_{i}\right\} _{i\in I}$ is an orthornormal
family in the real inner product space $\left( H;\left\langle \cdot ,\cdot
\right\rangle \right) $ and $M_{i},m_{i}\in \mathbb{R}$, $i\in F$ ($F$ is a
finite part of $I$) and $x\in H$ are such that $M_{i},m_{i}\geq 0$ for $i\in
F$ with $\sum_{i\in F}M_{i}m_{i}\geq 0$ and%
\begin{equation*}
\left\langle \sum_{i\in F}M_{i}e_{i}-x,x-\sum_{i\in
F}m_{i}e_{i}\right\rangle \geq 0,
\end{equation*}%
then we have the inequality%
\begin{equation*}
0\leq \left\Vert x\right\Vert ^{2}-\sum_{i\in F}\left[ \left\langle
x,e_{i}\right\rangle \right] ^{2}\leq \frac{1}{4}\cdot \frac{\sum_{i\in
F}\left( M_{i}-m_{i}\right) ^{2}}{\sum_{i\in F}M_{i}m_{i}}\cdot \sum_{i\in F}%
\left[ \left\langle x,e_{i}\right\rangle \right] ^{2}.
\end{equation*}%
The constant $\frac{1}{4}$ is best possible.
\end{remark}

The following reverse of the Schwarz's inequality in inner product spaces
holds.

\begin{corollary}
\label{c2.5.10}Let $x,y\in H$ and $\delta ,\Delta \in \mathbb{K}$ $\left( 
\mathbb{K}=\mathbb{C},\mathbb{R}\right) $ with the property that $\func{Re}%
\left( \Delta \overline{\delta }\right) >0.$ If either%
\begin{equation*}
\func{Re}\left\langle \Delta y-x,x-\delta y\right\rangle \geq 0,
\end{equation*}%
or equivalently,%
\begin{equation*}
\left\Vert x-\frac{\delta +\Delta }{2}\cdot y\right\Vert \leq \frac{1}{2}%
\left\vert \Delta -\delta \right\vert \left\Vert y\right\Vert ,
\end{equation*}%
holds, then we have the inequalities%
\begin{align}
\left\Vert x\right\Vert \left\Vert y\right\Vert & \leq \frac{1}{2}\cdot 
\frac{\func{Re}\left[ \Delta \overline{\left\langle x,y\right\rangle }+%
\overline{\delta }\left\langle x,y\right\rangle \right] }{\sqrt{\Delta 
\overline{\delta }}}  \label{2.20.10} \\
& \leq \frac{1}{2}\cdot \frac{\left\vert \Delta \right\vert +\left\vert
\delta \right\vert }{\sqrt{\Delta \overline{\delta }}}\left\vert
\left\langle x,y\right\rangle \right\vert ,  \notag
\end{align}%
\begin{align}
0& \leq \left\Vert x\right\Vert \left\Vert y\right\Vert -\left\vert
\left\langle x,y\right\rangle \right\vert  \label{2.21.10} \\
& \leq \frac{1}{2}\cdot \frac{\left( \sqrt{\left\vert \Delta \right\vert }-%
\sqrt{\left\vert \delta \right\vert }\right) ^{2}+2\left( \sqrt{\Delta 
\overline{\delta }}-\sqrt{\func{Re}\left( \Delta \overline{\delta }\right) }%
\right) }{\sqrt{\Delta \overline{\delta }}}\left\vert \left\langle
x,y\right\rangle \right\vert ,  \notag
\end{align}%
\begin{equation}
\left\Vert x\right\Vert ^{2}\left\Vert y\right\Vert ^{2}\leq \frac{1}{4}%
\cdot \frac{\left( \left\vert \Delta \right\vert +\left\vert \delta
\right\vert \right) ^{2}}{\func{Re}\left( \Delta \overline{\delta }\right) }%
\left\vert \left\langle x,y\right\rangle \right\vert ^{2},  \label{2.22.10}
\end{equation}%
and 
\begin{align}
0& \leq \left\Vert x\right\Vert ^{2}\left\Vert y\right\Vert ^{2}-\left\vert
\left\langle x,y\right\rangle \right\vert ^{2}  \label{2.23.10} \\
& \leq \frac{1}{4}\cdot \frac{\left( \left\vert \Delta \right\vert
+\left\vert \delta \right\vert \right) ^{2}+4\left( \left\vert \Delta 
\overline{\delta }\right\vert -\func{Re}\left( \Delta \overline{\delta }%
\right) \right) }{\func{Re}\left( \Delta \overline{\delta }\right) }%
\left\vert \left\langle x,y\right\rangle \right\vert ^{2}.  \notag
\end{align}%
The constants $\frac{1}{2}$ and $\frac{1}{4}$ are best possible.
\end{corollary}

\begin{proof}
The inequality (\ref{2.20.10}) follows from (\ref{2.6.10}) on choosing $%
F=\left\{ 1\right\} ,$ $e_{1}=e=\frac{y}{\left\Vert y\right\Vert },$ $\Phi
_{1}=\Phi =\Delta \left\Vert y\right\Vert ,$ \ $\phi _{1}=\phi =\delta
\left\Vert y\right\Vert $ $\left( y\neq 0\right) .$ The inequality (\ref%
{2.21.10}) is equivalent with (\ref{2.20.10}). The inequality (\ref{2.22.10}%
) follows from (\ref{2.1.10}) for $F=\left\{ 1\right\} $ and the same
choices as above. Finally, (\ref{2.23.10}) is obviously equivalent with (\ref%
{2.22.10}).
\end{proof}

\subsection{Some Gr\"{u}ss Type Inequalities for Orthonormal Families}

The following result holds \cite{11NSSD}.

\begin{theorem}
\label{t3.1.10}Let $\left\{ e_{i}\right\} _{i\in I}$ be a family of
orthornormal vectors in $H,$ $F$ a finite part of $I$, $\phi _{i},\Phi _{i},$
$\gamma _{i},\Gamma _{i}\in \mathbb{K},\ i\in F$ and $x,y\in H.$ If either%
\begin{align*}
\func{Re}\left\langle \sum_{i\in F}\Phi _{i}e_{i}-x,x-\sum_{i\in F}\phi
_{i}e_{i}\right\rangle & \geq 0, \\
\func{Re}\left\langle \sum_{i\in F}\Gamma _{i}e_{i}-y,y-\sum_{i\in F}\gamma
_{i}e_{i}\right\rangle & \geq 0,
\end{align*}%
or equivalently,%
\begin{align*}
\left\Vert x-\sum_{i\in F}\frac{\Phi _{i}+\phi _{i}}{2}e_{i}\right\Vert &
\leq \frac{1}{2}\left( \sum_{i\in F}\left\vert \Phi _{i}-\phi
_{i}\right\vert ^{2}\right) ^{\frac{1}{2}}, \\
\left\Vert y-\sum_{i\in F}\frac{\Gamma _{i}+\gamma _{i}}{2}e_{i}\right\Vert
& \leq \frac{1}{2}\left( \sum_{i\in F}\left\vert \Gamma _{i}-\gamma
_{i}\right\vert ^{2}\right) ^{\frac{1}{2}},
\end{align*}%
hold, then we have the inequality%
\begin{align}
0& \leq \left\vert \left\langle x,y\right\rangle -\sum_{i\in F}\left\langle
x,e_{i}\right\rangle \left\langle e_{i},y\right\rangle \right\vert
\label{3.3.10} \\
& \leq \frac{1}{4}M\left( \mathbf{\Phi },\mathbf{\phi },F\right) M\left( 
\mathbf{\Gamma },\mathbf{\gamma },F\right) \left( \sum_{i\in F}\left\vert
\left\langle x,e_{i}\right\rangle \right\vert ^{2}\right) ^{\frac{1}{2}%
}\left( \sum_{i\in F}\left\vert \left\langle y,e_{i}\right\rangle
\right\vert ^{2}\right) ^{\frac{1}{2}},  \notag
\end{align}%
where $M\left( \mathbf{\Phi },\mathbf{\phi },F\right) $ is defined in
Corollary \ref{c2.3.10}.

The constant $\frac{1}{4}$ is best possible.
\end{theorem}

\begin{proof}
By the reverse of Bessel's inequality in Corollary \ref{c2.3.10}, we have%
\begin{align}
& \left\vert \left\langle x,y\right\rangle -\sum_{i\in F}\left\langle
x,e_{i}\right\rangle \left\langle e_{i},y\right\rangle \right\vert ^{2}
\label{3.5.10} \\
& \leq \left( \left\Vert x\right\Vert ^{2}-\sum_{i\in F}\left\vert
\left\langle x,e_{i}\right\rangle \right\vert ^{2}\right) \left( \left\Vert
y\right\Vert ^{2}-\sum_{i\in F}\left\vert \left\langle y,e_{i}\right\rangle
\right\vert ^{2}\right)  \notag \\
& \leq \frac{1}{4}M^{2}\left( \mathbf{\Phi },\mathbf{\phi },F\right)
\sum_{i\in F}\left\vert \left\langle x,e_{i}\right\rangle \right\vert
^{2}\cdot \frac{1}{4}M^{2}\left( \mathbf{\Gamma },\mathbf{\gamma },F\right)
\sum_{i\in F}\left\vert \left\langle y,e_{i}\right\rangle \right\vert ^{2}. 
\notag
\end{align}%
Taking the square root in (\ref{3.5.10}), we deduce (\ref{3.3.10}).

The fact that $\frac{1}{4}$ is the best possible constant follows by
Corollary \ref{c2.3.10} and we omit the details.
\end{proof}

The following corollary for real inner product spaces holds \cite{11NSSD}.

\begin{corollary}
\label{c3.2.10}Let $\left\{ e_{i}\right\} _{i\in I}$ be a family of
orthornormal vectors in $H,$ $F$ a finite part of $I$, $M_{i},m_{i},$ $%
N_{i},n_{i}\geq 0,\ i\in F$ and $x,y\in H$ such that $\sum_{i\in
F}M_{i}m_{i}>0,$ $\sum_{i\in F}N_{i}n_{i}>0$ and%
\begin{equation*}
\left\langle \sum_{i\in F}M_{i}e_{i}-x,x-\sum_{i\in
F}m_{i}e_{i}\right\rangle \geq 0,\ \ \ \ \left\langle \sum_{i\in
F}N_{i}e_{i}-y,y-\sum_{i\in F}n_{i}e_{i}\right\rangle \geq 0.
\end{equation*}%
Then we have the inequality%
\begin{align*}
0& \leq \left\vert \left\langle x,y\right\rangle -\sum_{i\in F}\left\langle
x,e_{i}\right\rangle \left\langle y,e_{i}\right\rangle \right\vert ^{2} \\
& \leq \frac{1}{16}\cdot \frac{\sum_{i\in F}\left( M_{i}-m_{i}\right)
^{2}\sum_{i\in F}\left( N_{i}-n_{i}\right) ^{2}\sum_{i\in F}\left\vert
\left\langle x,e_{i}\right\rangle \right\vert ^{2}\sum_{i\in F}\left\vert
\left\langle y,e_{i}\right\rangle \right\vert ^{2}}{\sum_{i\in
F}M_{i}m_{i}\sum_{i\in F}N_{i}n_{i}}.
\end{align*}%
The constant $\frac{1}{16}$ is best possible.
\end{corollary}

In the case where the family $\left\{ e_{i}\right\} _{i\in I}$ reduces to a
single vector, we may deduce from Theorem \ref{t3.1.10} the following
particular case first obtained in \cite{SSDc.10}.

\begin{corollary}
\label{c3.3.10}Let $e\in H,$ $\left\Vert e\right\Vert =1,$ $\phi ,\Phi
,\gamma ,\Gamma \in \mathbb{K}$ \ with $\func{Re}\left( \Phi \overline{\phi }%
\right) ,$ $\func{Re}\left( \Gamma \overline{\gamma }\right) >0$ and $x,y\in
H$ such that either 
\begin{equation*}
\func{Re}\left\langle \Phi e-x,x-\phi e\right\rangle \geq 0,\ \ \func{Re}%
\left\langle \Gamma e-y,y-\gamma e\right\rangle \geq 0,
\end{equation*}%
or equivalently,%
\begin{equation*}
\left\Vert x-\frac{\phi +\Phi }{2}e\right\Vert \leq \frac{1}{2}\left\vert
\Phi -\phi \right\vert ,\ \ \ \ \ \left\Vert y-\frac{\gamma +\Gamma }{2}%
e\right\Vert \leq \frac{1}{2}\left\vert \Gamma -\gamma \right\vert
\end{equation*}%
hold, then%
\begin{equation*}
0\leq \left\vert \left\langle x,y\right\rangle -\left\langle
x,e\right\rangle \left\langle e,y\right\rangle \right\vert \leq \frac{1}{4}%
M\left( \Phi ,\phi \right) M\left( \Gamma ,\gamma \right) \left\vert
\left\langle x,e\right\rangle \left\langle e,y\right\rangle \right\vert ,
\end{equation*}%
where%
\begin{equation*}
M\left( \Phi ,\phi \right) :=\left[ \frac{\left( \left\vert \Phi \right\vert
-\left\vert \phi \right\vert \right) ^{2}+4\left[ \left\vert \phi \Phi
\right\vert -\func{Re}\left( \Phi \overline{\phi }\right) \right] }{\func{Re}%
\left( \Phi \overline{\phi }\right) }\right] ^{\frac{1}{2}}.
\end{equation*}%
The constant $\frac{1}{4}$ is best possible.
\end{corollary}

\begin{remark}
\label{r3.4.10}If $H$ is real, $e\in H,$ $\left\Vert e\right\Vert =1$ and $%
a,b,A,B\in \mathbb{R}$ are such that $A>a>0,$ $B>b>0$ and%
\begin{equation*}
\left\Vert x-\frac{a+A}{2}e\right\Vert \leq \frac{1}{2}\left( A-a\right) ,\
\ \left\Vert y-\frac{b+B}{2}e\right\Vert \leq \frac{1}{2}\left( B-b\right) ,
\end{equation*}%
then%
\begin{equation}
\left\vert \left\langle x,y\right\rangle -\left\langle x,e\right\rangle
\left\langle e,y\right\rangle \right\vert \leq \frac{1}{4}\cdot \frac{\left(
A-a\right) \left( B-b\right) }{\sqrt{abAB}}\left\vert \left\langle
x,e\right\rangle \left\langle e,y\right\rangle \right\vert .  \label{3.12.10}
\end{equation}%
The constant $\frac{1}{4}$ is best possible.

If $\left\langle x,e\right\rangle ,$ $\left\langle y,e\right\rangle \neq 0,$
then the following equivalent form of (\ref{3.12.10}) also holds%
\begin{equation*}
\left\vert \frac{\left\langle x,y\right\rangle }{\left\langle
x,e\right\rangle \left\langle e,y\right\rangle }-1\right\vert \leq \frac{1}{4%
}\cdot \frac{\left( A-a\right) \left( B-b\right) }{\sqrt{abAB}}.
\end{equation*}
\end{remark}

\subsection{Some Companion Inequalities}

The following companion of the Gr\"{u}ss inequality also holds \cite{11NSSD}.

\begin{theorem}
\label{t4.1.10}Let $\left\{ e_{i}\right\} _{i\in I}$ be a family of
orthornormal vectors in $H,$ $F$ a finite part of $I$, $\phi _{i},\Phi
_{i}\in \mathbb{K},\ \left( i\in F\right) $, $x,y\in H$ and $\lambda \in
\left( 0,1\right) ,$ such that either%
\begin{equation}
\func{Re}\left\langle \sum_{i\in F}\Phi _{i}e_{i}-\left( \lambda x+\left(
1-\lambda \right) y\right) ,\lambda x+\left( 1-\lambda \right) y-\sum_{i\in
F}\phi _{i}e_{i}\right\rangle \geq 0,  \label{4.1.10}
\end{equation}%
or equivalently,%
\begin{equation*}
\left\Vert \lambda x+\left( 1-\lambda \right) y-\sum_{i\in F}\frac{\Phi
_{i}+\phi _{i}}{2}\cdot e_{i}\right\Vert \leq \frac{1}{2}\left( \sum_{i\in
F}\left\vert \Phi _{i}-\phi _{i}\right\vert ^{2}\right) ^{\frac{1}{2}},
\end{equation*}%
holds. Then we have the inequality%
\begin{multline}
\func{Re}\left[ \left\langle x,y\right\rangle -\sum_{i\in F}\left\langle
x,e_{i}\right\rangle \left\langle e_{i},y\right\rangle \right]
\label{4.3.10} \\
\leq \frac{1}{16}\cdot \frac{1}{\lambda \left( 1-\lambda \right) }\sum_{i\in
F}M^{2}\left( \mathbf{\Phi },\mathbf{\phi },F\right) \sum_{i\in F}\left\vert
\left\langle \lambda x+\left( 1-\lambda \right) y,e_{i}\right\rangle
\right\vert ^{2}.
\end{multline}%
The constant $\frac{1}{16}$ is the best possible constant in (\ref{4.3.10})
in the sense that it cannot be replaced by a smaller constant.
\end{theorem}

\begin{proof}
Using the known inequality%
\begin{equation*}
\func{Re}\left\langle z,u\right\rangle \leq \frac{1}{4}\left\Vert
z+u\right\Vert ^{2}
\end{equation*}%
we may state that for any $a,b\in H$ and $\lambda \in \left( 0,1\right) $ 
\begin{equation}
\func{Re}\left\langle a,b\right\rangle \leq \frac{1}{4\lambda \left(
1-\lambda \right) }\left\Vert \lambda a+\left( 1-\lambda \right)
b\right\Vert ^{2}.  \label{4.4.10}
\end{equation}%
Since%
\begin{equation*}
\left\langle x,y\right\rangle -\sum_{i\in F}\left\langle
x,e_{i}\right\rangle \left\langle e_{i},y\right\rangle =\left\langle
x-\sum_{i\in F}\left\langle x,e_{i}\right\rangle e_{i},y-\sum_{i\in
F}\left\langle y,e_{i}\right\rangle e_{i}\right\rangle ,
\end{equation*}%
for any \thinspace $x,y\in H,$ then, by (\ref{4.4.10}), we get%
\begin{align}
& \func{Re}\left[ \left\langle x,y\right\rangle -\sum_{i\in F}\left\langle
x,e_{i}\right\rangle \left\langle e_{i},y\right\rangle \right]
\label{4.5.10} \\
& =\func{Re}\left[ \left\langle x-\sum_{i\in F}\left\langle
x,e_{i}\right\rangle e_{i},y-\sum_{i\in F}\left\langle y,e_{i}\right\rangle
e_{i}\right\rangle \right]  \notag \\
& \leq \frac{1}{4\lambda \left( 1-\lambda \right) }\left\Vert \lambda \left(
x-\sum_{i\in F}\left\langle x,e_{i}\right\rangle e_{i}\right) +\left(
1-\lambda \right) \left( y-\sum_{i\in F}\left\langle y,e_{i}\right\rangle
e_{i}\right) \right\Vert ^{2}  \notag \\
& =\frac{1}{4\lambda \left( 1-\lambda \right) }\left\Vert \lambda x+\left(
1-\lambda \right) y-\sum_{i\in F}\left\langle \lambda x+\left( 1-\lambda
\right) y,e_{i}\right\rangle e_{i}\right\Vert ^{2}  \notag \\
& =\frac{1}{4\lambda \left( 1-\lambda \right) }\left[ \left\Vert \lambda
x+\left( 1-\lambda \right) y\right\Vert ^{2}-\sum_{i\in F}\left\vert
\left\langle \lambda x+\left( 1-\lambda \right) y,e_{i}\right\rangle
\right\vert ^{2}\right] .  \notag
\end{align}%
If we apply the reverse of Bessel's inequality from Corollary \ref{c2.3.10}
for $\lambda x+\left( 1-\lambda \right) y,$ we may state that%
\begin{multline}
\left\Vert \lambda x+\left( 1-\lambda \right) y\right\Vert ^{2}-\sum_{i\in
F}\left\vert \left\langle \lambda x+\left( 1-\lambda \right)
y,e_{i}\right\rangle \right\vert ^{2}  \label{4.6.10} \\
\leq \frac{1}{4}M^{2}\left( \mathbf{\Phi },\mathbf{\phi },F\right)
\sum_{i\in F}\left\vert \left\langle \lambda x+\left( 1-\lambda \right)
y,e_{i}\right\rangle \right\vert ^{2}.
\end{multline}%
Now, by making use of (\ref{4.5.10}) and (\ref{4.6.10}), we deduce (\ref%
{4.3.10}).

The fact that $\frac{1}{16}$ is the best possible constant in (\ref{4.3.10})
follows by the fact that if in (\ref{4.1.10}) we choose $x=y,$ then it
becomes (i) of Theorem \ref{t2.1.10}, implying for $\lambda =\frac{1}{2}$
the inequality (\ref{2.12.10}), for which, we have shown that $\frac{1}{4}$
is the best constant.
\end{proof}

\begin{remark}
\label{r4.2.10}If in Theorem \ref{t4.1.10}, we choose $\lambda =\frac{1}{2},$
then we get%
\begin{equation*}
\func{Re}\left[ \left\langle x,y\right\rangle -\sum_{i\in F}\left\langle
x,e_{i}\right\rangle \left\langle e_{i},y\right\rangle \right] \leq \frac{1}{%
4}M^{2}\left( \mathbf{\Phi },\mathbf{\phi },F\right) \sum_{i\in F}\left\vert
\left\langle \frac{x+y}{2},e_{i}\right\rangle \right\vert ^{2},
\end{equation*}%
provided%
\begin{equation*}
\func{Re}\left\langle \sum_{i\in F}\Phi _{i}e_{i}-\frac{x+y}{2},\frac{x+y}{2}%
-\sum_{i\in F}\phi _{i}e_{i}\right\rangle \geq 0
\end{equation*}%
or equivalently,%
\begin{equation*}
\left\Vert \frac{x+y}{2}-\sum_{i\in F}\frac{\Phi _{i}+\phi _{i}}{2}\cdot
e_{i}\right\Vert \leq \frac{1}{2}\left( \sum_{i\in F}\left\vert \Phi
_{i}-\phi _{i}\right\vert ^{2}\right) ^{\frac{1}{2}}.
\end{equation*}
\end{remark}

\subsection{Integral Inequalities}

The following proposition holds \cite{11NSSD}.

\begin{proposition}
\label{p5.1.10}Let $\left\{ f_{i}\right\} _{i\in I}$ be an orthornormal
family of functions in $L_{\rho }^{2}\left( \Omega ,\mathbb{K}\right) ,$ $F$
a finite subset of $I,$ $\phi _{i},\Phi _{i}\in \mathbb{K}$ $\left( i\in
F\right) $ such that $\sum_{i\in F}\func{Re}\left( \Phi _{i}\overline{\phi
_{i}}\right) >0$ and $f\in L_{\rho }^{2}\left( \Omega ,\mathbb{K}\right) ,$
so that either%
\begin{equation*}
\int_{\Omega }\rho \left( s\right) \func{Re}\left[ \left( \sum_{i\in F}\Phi
_{i}f_{i}\left( s\right) -f\left( s\right) \right) \left( \overline{f}\left(
s\right) -\sum_{i\in F}\overline{\phi _{i}}\text{ }\overline{f_{i}}\left(
s\right) \right) \right] d\mu \left( s\right) \geq 0,
\end{equation*}%
or equivalently,%
\begin{equation*}
\int_{\Omega }\rho \left( s\right) \left\vert f\left( s\right) -\sum_{i\in F}%
\frac{\Phi _{i}+\phi _{i}}{2}f_{i}\left( s\right) \right\vert ^{2}d\mu
\left( s\right) \leq \frac{1}{4}\sum_{i\in F}\left\vert \Phi _{i}-\phi
_{i}\right\vert ^{2}.
\end{equation*}%
Then we have the inequality%
\begin{equation*}
\left( \int_{\Omega }\rho \left( s\right) \left\vert f\left( s\right)
\right\vert ^{2}d\mu \left( s\right) \right) ^{\frac{1}{2}}\leq \frac{1}{2}%
\cdot \frac{1}{\left[ \sum_{i\in F}\func{Re}\left( \Phi _{i}\overline{\phi
_{i}}\right) \right] ^{\frac{1}{2}}}
\end{equation*}%
\begin{equation*}
\times \left\{ 
\begin{array}{l}
\max\limits_{i\in F}\left\{ \left\vert \Phi _{i}\right\vert +\left\vert \phi
_{i}\right\vert \right\} \dsum\limits_{i\in F}\left\vert \dint_{\Omega }\rho
\left( s\right) f\left( s\right) \overline{f_{i}}\left( s\right) d\mu \left(
s\right) \right\vert \\ 
\\ 
\left[ \dsum\limits_{i\in F}\left( \left\vert \Phi _{i}\right\vert
+\left\vert \phi _{i}\right\vert \right) ^{p}\right] ^{\frac{1}{p}}\left(
\dsum\limits_{i\in F}\left\vert \dint_{\Omega }\rho \left( s\right) f\left(
s\right) \overline{f_{i}}\left( s\right) d\mu \left( s\right) \right\vert
^{q}\right) ^{\frac{1}{q}},\text{ } \\ 
\hfill \text{\ for \ }p>1,\ \frac{1}{p}+\frac{1}{q}=1 \\ 
\\ 
\max\limits_{i\in F}\left\vert \dint_{\Omega }\rho \left( s\right) f\left(
s\right) \overline{f_{i}}\left( s\right) d\mu \left( s\right) \right\vert
\dsum\limits_{i\in F}\left[ \left\vert \Phi _{i}\right\vert +\left\vert \phi
_{i}\right\vert \right] .%
\end{array}%
\right.
\end{equation*}%
In particular, we have%
\begin{multline}
\int_{\Omega }\rho \left( s\right) \left\vert f\left( s\right) \right\vert
^{2}d\mu \left( s\right)  \label{5.6.10} \\
\leq \frac{1}{4}\cdot \frac{\sum_{i\in F}\left( \left\vert \Phi
_{i}\right\vert +\left\vert \phi _{i}\right\vert \right) ^{2}}{\sum_{i\in F}%
\func{Re}\left( \Phi _{i}\overline{\phi _{i}}\right) }\sum\limits_{i\in
F}\left\vert \int_{\Omega }\rho \left( s\right) f\left( s\right) \overline{%
f_{i}}\left( s\right) d\mu \left( s\right) \right\vert ^{2}.
\end{multline}%
The constant $\frac{1}{4}$ is best possible.
\end{proposition}

The proof is obvious by Theorem \ref{t2.1.10} and Remark \ref{r2.2.10}. We
omit the details.

The following proposition also holds.

\begin{proposition}
\label{p5.2.10}Assume that $f_{i},f,\phi _{i},\Phi _{i}$ and $F$ satisfy the
assumptions of Proposition \ref{p5.1.10}. Then we have the following reverse
of Bessel's inequality:%
\begin{align}
0& \leq \int_{\Omega }\rho \left( s\right) f^{2}\left( s\right) d\mu \left(
s\right) -\sum\limits_{i\in F}\left\vert \int_{\Omega }\rho \left( s\right)
f\left( s\right) \overline{f_{i}}\left( s\right) d\mu \left( s\right)
\right\vert ^{2}  \label{5.7.10} \\
& \leq \frac{1}{4}M^{2}\left( \mathbf{\Phi },\mathbf{\phi },F\right) \cdot
\sum\limits_{i\in F}\left\vert \int_{\Omega }\rho \left( s\right) f\left(
s\right) \overline{f_{i}}\left( s\right) d\mu \left( s\right) \right\vert
^{2},  \notag
\end{align}%
where, as above,%
\begin{equation}
M\left( \mathbf{\Phi },\mathbf{\phi },F\right) :=\left[ \frac{%
\sum\limits_{i\in F}\left\{ \left( \left\vert \Phi _{i}\right\vert
-\left\vert \phi _{i}\right\vert \right) ^{2}+4\left[ \left\vert \phi
_{i}\Phi _{i}\right\vert -\func{Re}\left( \Phi _{i}\overline{\phi _{i}}%
\right) \right] \right\} }{\func{Re}\left( \Phi _{i}\overline{\phi _{i}}%
\right) }\right] ^{\frac{1}{2}}.  \label{5.8.10}
\end{equation}%
The constant $\frac{1}{4}$ is the best possible.
\end{proposition}

The following Gr\"{u}ss type inequality also holds.

\begin{proposition}
\label{p5.3.10}Let $\left\{ f_{i}\right\} _{i\in I}$ and $F$ be as in
Proposition \ref{p5.1.10}. If $\phi _{i},\Phi _{i},\gamma _{i},\Gamma
_{i}\in \mathbb{K}$ $\left( i\in F\right) $ and $f,g\in L_{\rho }^{2}\left(
\Omega ,\mathbb{K}\right) $ such that either%
\begin{align*}
\int_{\Omega }\rho \left( s\right) \func{Re}\left[ \left( \sum_{i\in F}\Phi
_{i}f_{i}\left( s\right) -f\left( s\right) \right) \left( \overline{f}\left(
s\right) -\sum_{i\in F}\overline{\phi _{i}}\text{ }\overline{f_{i}}\left(
s\right) \right) \right] d\mu \left( s\right) & \geq 0, \\
\int_{\Omega }\rho \left( s\right) \func{Re}\left[ \left( \sum_{i\in
F}\Gamma _{i}f_{i}\left( s\right) -g\left( s\right) \right) \left( \overline{%
g}\left( s\right) -\sum_{i\in F}\overline{\gamma _{i}}\text{ }\overline{f_{i}%
}\left( s\right) \right) \right] d\mu \left( s\right) & \geq 0,
\end{align*}%
or equivalently,%
\begin{align*}
\int_{\Omega }\rho \left( s\right) \left\vert f\left( s\right) -\sum_{i\in F}%
\frac{\Phi _{i}+\phi _{i}}{2}\cdot f_{i}\left( s\right) \right\vert ^{2}d\mu
\left( s\right) & \leq \frac{1}{4}\sum_{i\in F}\left\vert \Phi _{i}-\phi
_{i}\right\vert ^{2}, \\
\int_{\Omega }\rho \left( s\right) \left\vert g\left( s\right) -\sum_{i\in F}%
\frac{\Gamma _{i}+\gamma _{i}}{2}\cdot f_{i}\left( s\right) \right\vert
^{2}d\mu \left( s\right) & \leq \frac{1}{4}\sum_{i\in F}\left\vert \Gamma
_{i}-\gamma _{i}\right\vert ^{2},
\end{align*}%
hold, then we have the inequality%
\begin{multline}
\left\vert \int_{\Omega }\rho \left( s\right) f\left( s\right) \overline{%
g\left( s\right) }d\mu \left( s\right) \right.  \label{5.11.10} \\
-\left. \sum_{i\in F}\int_{\Omega }\rho \left( s\right) f\left( s\right) 
\overline{f_{i}}\left( s\right) d\mu \left( s\right) \int_{\Omega }\rho
\left( s\right) f_{i}\left( s\right) \overline{g\left( s\right) }d\mu \left(
s\right) \right\vert \\
\leq \frac{1}{4}M\left( \mathbf{\Phi },\mathbf{\phi },F\right) M\left( 
\mathbf{\Gamma },\mathbf{\gamma },F\right) \left( \sum_{i\in F}\left\vert
\int_{\Omega }\rho \left( s\right) f\left( s\right) \overline{f_{i}}\left(
s\right) d\mu \left( s\right) \right\vert ^{2}\right) ^{\frac{1}{2}} \\
\times \left( \sum_{i\in F}\left\vert \rho \left( s\right) f_{i}\left(
s\right) \overline{g\left( s\right) }d\mu \left( s\right) \right\vert
^{2}\right) ^{\frac{1}{2}},
\end{multline}%
where $M\left( \mathbf{\Phi },\mathbf{\phi },F\right) $ is as defined in (%
\ref{5.8.10}).

The constant $\frac{1}{4}$ is the best possible.
\end{proposition}

The proof follows by Theorem \ref{t3.1.10} and we omit the details.

In the case of real spaces, the following corollaries provide much simpler
sufficient conditions for the reverse of Bessel's inequality (\ref{5.7.10})
or for the Gr\"{u}ss type inequality (\ref{5.11.10}) to hold.

\begin{corollary}
\label{c5.4.10}Let $\left\{ f_{i}\right\} _{i\in I}$ be an orthornormal
family of functions in the real Hilbert space $L_{\rho }^{2}\left( \Omega
\right) ,$ $F$ a finite part of $I,$ $M_{i},m_{i}\geq 0$ \ $\left( i\in
F\right) ,$ with $\sum_{i\in F}M_{i}m_{i}>0$ and $f\in L_{\rho }^{2}\left(
\Omega \right) $ so that%
\begin{equation*}
\sum_{i\in F}m_{i}f_{i}\left( s\right) \leq f\left( s\right) \leq \sum_{i\in
F}M_{i}f_{i}\left( s\right) \text{ \ for \ }\mu -\text{a.e. }s\in \Omega .
\end{equation*}%
Then we have the inequality%
\begin{align*}
0& \leq \int_{\Omega }\rho \left( s\right) f^{2}\left( s\right) d\mu \left(
s\right) -\sum_{i\in F}\left[ \int_{\Omega }\rho \left( s\right) f\left(
s\right) f_{i}\left( s\right) d\mu \left( s\right) \right] ^{2} \\
& \leq \frac{1}{4}\cdot \frac{\sum_{i\in F}\left( M_{i}-m_{i}\right) ^{2}}{%
\sum_{i\in F}M_{i}m_{i}}\cdot \sum_{i\in F}\left[ \int_{\Omega }\rho \left(
s\right) f\left( s\right) f_{i}\left( s\right) d\mu \left( s\right) \right]
^{2}.
\end{align*}%
The constant $\frac{1}{4}$ is best possible.
\end{corollary}

\begin{corollary}
\label{c5.5.10}Let $\left\{ f_{i}\right\} _{i\in I}$ and $F$ be as above. If 
$M_{i},m_{i},N_{i},n_{i}\geq 0$ $\left( i\in F\right) $ with $\sum_{i\in
F}M_{i}m_{i},\sum_{i\in F}N_{i}n_{i}>0$ and $f,g\in L_{\rho }^{2}\left(
\Omega \right) $ such that%
\begin{equation*}
\sum_{i\in F}m_{i}f_{i}\left( s\right) \leq f\left( s\right) \leq \sum_{i\in
F}M_{i}f_{i}\left( s\right)
\end{equation*}%
and%
\begin{equation*}
\sum_{i\in F}n_{i}f_{i}\left( s\right) \leq g\left( s\right) \leq \sum_{i\in
F}N_{i}f_{i}\left( s\right) \text{ \ for \ }\mu -\text{a.e. }s\in \Omega ,
\end{equation*}%
then we have the inequality%
\begin{align*}
& \left\vert \int_{\Omega }\rho \left( s\right) f\left( s\right) g\left(
s\right) d\mu \left( s\right) \right. \\
& -\left. \sum_{i\in F}\int_{\Omega }\rho \left( s\right) f\left( s\right)
f_{i}\left( s\right) d\mu \left( s\right) \int_{\Omega }\rho \left( s\right)
g\left( s\right) f_{i}\left( s\right) d\mu \left( s\right) \right\vert \\
& \leq \frac{1}{4}\left( \frac{\sum_{i\in F}\left( M_{i}-m_{i}\right) ^{2}}{%
\sum_{i\in F}M_{i}m_{i}}\right) ^{\frac{1}{2}}\left( \frac{\sum_{i\in
F}\left( N_{i}-n_{i}\right) ^{2}}{\sum_{i\in F}N_{i}n_{i}}\right) ^{\frac{1}{%
2}} \\
& \times \left[ \sum_{i\in F}\left( \int_{\Omega }\rho \left( s\right)
f\left( s\right) f_{i}\left( s\right) d\mu \left( s\right) \right)
^{2}\sum_{i\in F}\left( \int_{\Omega }\rho \left( s\right) g\left( s\right)
f_{i}\left( s\right) d\mu \left( s\right) \right) ^{2}\right] ^{\frac{1}{2}}.
\end{align*}%
\qquad
\end{corollary}

\newpage

\section{General Reverses of Bessel's Inequality in Hilbert Spaces}

\subsection{Some Reverses of Bessel's Inequality}

Let $\left( H;\left\langle \cdot ,\cdot \right\rangle \right) $ be a real or
complex infinite dimensional Hilbert space and $\left( e_{i}\right) _{i\in 
\mathbb{N}}$ an orthornormal family in $H$, i.e., we recall that $%
\left\langle e_{i},e_{j}\right\rangle =0$ if $i,j\in \mathbb{N}$, $i\neq j$
and $\left\Vert e_{i}\right\Vert =1$ for $i\in \mathbb{N}$.

It is well known that, if $x\in H,$ then the sum $\sum_{i=1}^{\infty
}\left\vert \left\langle x,e_{i}\right\rangle \right\vert ^{2}$ is
convergent and the following inequality, called \textit{Bessel's inequality}%
\begin{equation}
\sum_{i=1}^{\infty }\left\vert \left\langle x,e_{i}\right\rangle \right\vert
^{2}\leq \left\Vert x\right\Vert ^{2},  \label{5.1.11a}
\end{equation}%
holds.

If $\ell ^{2}\left( \mathbb{K}\right) :=\left\{ \mathbf{a}=\left(
a_{i}\right) _{i\in \mathbb{N}}\subset \mathbb{K}\left\vert
\sum_{i=1}^{\infty }\left\vert a_{i}\right\vert ^{2}\right. <\infty \right\}
,$ where $\mathbb{K}=\mathbb{C}$ or $\mathbb{K}=\mathbb{R}$, is the Hilbert
space of all complex or real sequences that are $2$-summable and $\mathbf{%
\lambda }=\left( \lambda _{i}\right) _{i\in \mathbb{N}}\in \ell ^{2}\left( 
\mathbb{K}\right) ,$ then the sum $\sum_{i=1}^{\infty }\lambda _{i}e_{i}$ is
convergent in $H$ and if $y:=\sum_{i=1}^{\infty }\lambda _{i}e_{i}\in H,$
then $\left\Vert y\right\Vert =\left( \sum_{i=1}^{\infty }\left\vert \lambda
_{i}\right\vert ^{2}\right) ^{\frac{1}{2}}.$

We may state the following result \cite{SSDa.11}.

\begin{theorem}
\label{t5.1.11a}Let $\left( H;\left\langle \cdot ,\cdot \right\rangle
\right) $ be an infinite dimensional Hilbert space over the real or complex
number field $\mathbb{K}$, $\left( e_{i}\right) _{i\in \mathbb{N}}$ an
orthornormal family in $H,$ $\mathbf{\lambda }=\left( \lambda _{i}\right)
_{i\in \mathbb{N}}\in \ell ^{2}\left( \mathbb{K}\right) $ and $r>0$ with the
property that%
\begin{equation*}
\sum_{i=1}^{\infty }\left\vert \lambda _{i}\right\vert ^{2}>r^{2}.
\end{equation*}%
If $x\in H$ is such that%
\begin{equation*}
\left\Vert x-\sum_{i=1}^{\infty }\lambda _{i}e_{i}\right\Vert \leq r,
\end{equation*}%
then we have the inequality%
\begin{align}
\left\Vert x\right\Vert ^{2}& \leq \frac{\left( \sum_{i=1}^{\infty }\func{Re}%
\left[ \overline{\lambda _{i}}\left\langle x,e_{i}\right\rangle \right]
\right) ^{2}}{\sum_{i=1}^{\infty }\left\vert \lambda _{i}\right\vert
^{2}-r^{2}}  \label{5.4.11a} \\
& \leq \frac{\left\vert \sum_{i=1}^{\infty }\overline{\lambda _{i}}%
\left\langle x,e_{i}\right\rangle \right\vert ^{2}}{\sum_{i=1}^{\infty
}\left\vert \lambda _{i}\right\vert ^{2}-r^{2}}  \notag \\
& \leq \frac{\sum_{i=1}^{\infty }\left\vert \lambda _{i}\right\vert ^{2}}{%
\sum_{i=1}^{\infty }\left\vert \lambda _{i}\right\vert ^{2}-r^{2}}%
\sum_{i=1}^{\infty }\left\vert \left\langle x,e_{i}\right\rangle \right\vert
^{2}  \notag
\end{align}%
and%
\begin{align}
0& \leq \left\Vert x\right\Vert ^{2}-\sum_{i=1}^{\infty }\left\vert
\left\langle x,e_{i}\right\rangle \right\vert ^{2}  \label{5.5.11a} \\
& \leq \frac{r^{2}}{\sum_{i=1}^{\infty }\left\vert \lambda _{i}\right\vert
^{2}-r^{2}}\sum_{i=1}^{\infty }\left\vert \left\langle x,e_{i}\right\rangle
\right\vert ^{2}.  \notag
\end{align}
\end{theorem}

\begin{proof}
Applying the third inequality in (\ref{2.2.3}) for $a=\sum_{i=1}^{\infty
}\lambda _{i}e_{i}\in H,$ we have%
\begin{equation}
\left\Vert x\right\Vert ^{2}\left\Vert \sum_{i=1}^{\infty }\lambda
_{i}e_{i}\right\Vert ^{2}-\left[ \func{Re}\left\langle x,\sum_{i=1}^{\infty
}\lambda _{i}e_{i}\right\rangle \right] ^{2}\leq r^{2}\left\Vert
x\right\Vert ^{2}  \label{5.6.11a}
\end{equation}%
and since%
\begin{align*}
\left\Vert \sum_{i=1}^{\infty }\lambda _{i}e_{i}\right\Vert ^{2}&
=\sum_{i=1}^{\infty }\left\vert \lambda _{i}\right\vert ^{2}, \\
\func{Re}\left\langle x,\sum_{i=1}^{\infty }\lambda _{i}e_{i}\right\rangle &
=\sum_{i=1}^{\infty }\func{Re}\left[ \overline{\lambda _{i}}\left\langle
x,e_{i}\right\rangle \right] ,
\end{align*}%
then, by (\ref{5.6.11a}), we deduce%
\begin{equation*}
\left\Vert x\right\Vert ^{2}\sum_{i=1}^{\infty }\left\vert \lambda
_{i}\right\vert ^{2}-\left[ \func{Re}\left\langle x,\sum_{i=1}^{\infty
}\lambda _{i}e_{i}\right\rangle \right] ^{2}\leq r^{2}\left\Vert
x\right\Vert ^{2},
\end{equation*}%
giving the first inequality in (\ref{5.4.11a}).

The second inequality is obvious by the modulus property.

The last inequality follows by the Cauchy-Bunyakovsky-Schwarz inequality%
\begin{equation*}
\left\vert \sum_{i=1}^{\infty }\overline{\lambda _{i}}\left\langle
x,e_{i}\right\rangle \right\vert ^{2}\leq \sum_{i=1}^{\infty }\left\vert
\lambda _{i}\right\vert ^{2}\sum_{i=1}^{\infty }\left\vert \left\langle
x,e_{i}\right\rangle \right\vert ^{2}.
\end{equation*}%
The inequality (\ref{5.5.11a}) follows by the last inequality in (\ref%
{5.4.11a}) on subtracting from both sides the quantity $\sum_{i=1}^{\infty
}\left\vert \left\langle x,e_{i}\right\rangle \right\vert ^{2}<\infty .$
\end{proof}

The following result provides a generalization for the reverse of Bessel's
inequality obtained in \cite{SSD6.11}.

\begin{theorem}
\label{t5.2.11a}Let $\left( H;\left\langle \cdot ,\cdot \right\rangle
\right) $ and $\left( e_{i}\right) _{i\in \mathbb{N}}$ be as in Theorem \ref%
{t5.1.11a}. Suppose that $\mathbf{\Gamma }=\left( \Gamma _{i}\right) _{i\in 
\mathbb{N}}\in \ell ^{2}\left( \mathbb{K}\right) ,$ $\mathbf{\gamma }=\left(
\gamma _{i}\right) _{i\in \mathbb{N}}\in \ell ^{2}\left( \mathbb{K}\right) $
are sequences of real or complex numbers such that%
\begin{equation*}
\sum_{i=1}^{\infty }\func{Re}\left( \Gamma _{i}\overline{\gamma _{i}}\right)
>0.
\end{equation*}%
If $x\in H$ is such that either%
\begin{equation}
\left\Vert x-\sum_{i=1}^{\infty }\frac{\Gamma _{i}+\gamma _{i}}{2}%
e_{i}\right\Vert \leq \frac{1}{2}\left( \sum_{i=1}^{\infty }\left\vert
\Gamma _{i}-\gamma _{i}\right\vert ^{2}\right) ^{\frac{1}{2}}
\label{5.8.11a}
\end{equation}%
or equivalently,%
\begin{equation}
\func{Re}\left\langle \sum_{i=1}^{\infty }\Gamma
_{i}e_{i}-x,x-\sum_{i=1}^{\infty }\gamma _{i}e_{i}\right\rangle \geq 0
\label{5.9.11a}
\end{equation}%
holds, then we have the inequalities%
\begin{align}
\left\Vert x\right\Vert ^{2}& \leq \frac{1}{4}\cdot \frac{\left(
\sum_{i=1}^{\infty }\func{Re}\left[ \left( \overline{\Gamma _{i}}+\overline{%
\gamma _{i}}\right) \left\langle x,e_{i}\right\rangle \right] \right) ^{2}}{%
\sum_{i=1}^{\infty }\func{Re}\left( \Gamma _{i}\overline{\gamma _{i}}\right) 
}  \label{5.10.11a} \\
& \leq \frac{1}{4}\cdot \frac{\left\vert \sum_{i=1}^{\infty }\left( 
\overline{\Gamma _{i}}+\overline{\gamma _{i}}\right) \left\langle
x,e_{i}\right\rangle \right\vert ^{2}}{\sum_{i=1}^{\infty }\func{Re}\left(
\Gamma _{i}\overline{\gamma _{i}}\right) }  \notag \\
& \leq \frac{1}{4}\cdot \frac{\sum_{i=1}^{\infty }\left\vert \Gamma
_{i}+\gamma _{i}\right\vert ^{2}}{\sum_{i=1}^{\infty }\func{Re}\left( \Gamma
_{i}\overline{\gamma _{i}}\right) }\sum_{i=1}^{\infty }\left\vert
\left\langle x,e_{i}\right\rangle \right\vert ^{2}.  \notag
\end{align}%
The constant $\frac{1}{4}$ is best possible in all inequalities in (\ref%
{5.10.11a}).

We also have the inequalities:%
\begin{equation}
0\leq \left\Vert x\right\Vert ^{2}-\sum_{i=1}^{\infty }\left\vert
\left\langle x,e_{i}\right\rangle \right\vert ^{2}\leq \frac{1}{4}\cdot 
\frac{\sum_{i=1}^{\infty }\left\vert \Gamma _{i}-\gamma _{i}\right\vert ^{2}%
}{\sum_{i=1}^{\infty }\func{Re}\left( \Gamma _{i}\overline{\gamma _{i}}%
\right) }\sum_{i=1}^{\infty }\left\vert \left\langle x,e_{i}\right\rangle
\right\vert ^{2}.  \label{5.11.11a}
\end{equation}%
Here the constant $\frac{1}{4}$ is also best possible.
\end{theorem}

\begin{proof}
Since $\mathbf{\Gamma }$, $\mathbf{\gamma }\in \ell ^{2}\left( \mathbb{K}%
\right) ,$ then also $\frac{1}{2}\left( \mathbf{\Gamma }\pm \mathbf{\gamma }%
\right) \in \ell ^{2}\left( \mathbb{K}\right) ,$ showing that the series%
\begin{equation*}
\sum_{i=1}^{\infty }\left\vert \frac{\Gamma _{i}+\gamma _{i}}{2}\right\vert
^{2},\ \sum_{i=1}^{\infty }\left\vert \frac{\Gamma _{i}-\gamma _{i}}{2}%
\right\vert ^{2}\text{ and}\ \sum_{i=1}^{\infty }\func{Re}\left( \Gamma _{i}%
\overline{\gamma _{i}}\right)
\end{equation*}%
are convergent. Also, the series 
\begin{equation*}
\sum_{i=1}^{\infty }\Gamma _{i}e_{i},\text{ }\sum_{i=1}^{\infty }\gamma
_{i}e_{i}\text{ and }\sum_{i=1}^{\infty }\frac{\gamma _{i}+\Gamma _{i}}{2}%
e_{i}
\end{equation*}%
are convergent in the Hilbert space $H.$

Now, we observe that the inequalities (\ref{5.10.11a}) and (\ref{5.11.11a})
follow from Theorem \ref{t5.1.11a} on choosing $\lambda _{i}=\frac{\gamma
_{i}+\Gamma _{i}}{2},$ $i\in \mathbb{N}$ and $r=\frac{1}{2}\left(
\sum_{i=1}^{\infty }\left\vert \Gamma _{i}-\gamma _{i}\right\vert
^{2}\right) ^{\frac{1}{2}}.$

The fact that $\frac{1}{4}$ is the best constant in both (\ref{5.10.11a})
and (\ref{5.11.11a}) follows from Theorem \ref{t2.2.3} and Corollary \ref%
{c2.3.3}, and we omit the details.
\end{proof}

For some recent results related to the Bessel inequality, see \cite{HXC.11}, 
\cite{SSD01.11}, \cite{SSDJS.11}, and \cite{GH1.11}.

\subsection{Some Gr\"{u}ss Type Inequalities for Orthonormal Families}

The following result related to the Gr\"{u}ss inequality in inner product
spaces holds \cite{SSDa.11}.

\begin{theorem}
\label{t6.1.11a}Let $\left( H;\left\langle \cdot ,\cdot \right\rangle
\right) $ be an infinite dimensional Hilbert space over the real or complex
number field $\mathbb{K}$, and $\left( e_{i}\right) _{i\in \mathbb{N}}$ an
orthornormal family in $H.$ Assume that $\mathbf{\lambda }=\left( \lambda
_{i}\right) _{i\in \mathbb{N}},\ \mathbf{\mu }=\left( \mu _{i}\right) _{i\in 
\mathbb{N}}\in \ell ^{2}\left( \mathbb{K}\right) $ and $r_{1},r_{2}>0$ with
the properties that%
\begin{equation*}
\sum_{i=1}^{\infty }\left\vert \lambda _{i}\right\vert ^{2}>r_{1}^{2},\ \ \
\sum_{i=1}^{\infty }\left\vert \mu _{i}\right\vert ^{2}>r_{2}^{2}.
\end{equation*}%
If $x,y\in H$ are such that%
\begin{equation*}
\left\Vert x-\sum_{i=1}^{\infty }\lambda _{i}e_{i}\right\Vert \leq r_{1},\ \
\ \ \ \ \left\Vert y-\sum_{i=1}^{\infty }\mu _{i}e_{i}\right\Vert \leq r_{2},
\end{equation*}%
then we have the inequalities%
\begin{align}
& \left\vert \left\langle x,y\right\rangle -\sum_{i=1}^{\infty }\left\langle
x,e_{i}\right\rangle \left\langle e_{i},y\right\rangle \right\vert
\label{6.3.11a} \\
& \leq \frac{r_{1}r_{2}}{\sqrt{\sum_{i=1}^{\infty }\left\vert \lambda
_{i}\right\vert ^{2}-r_{1}^{2}}\sqrt{\sum_{i=1}^{\infty }\left\vert \mu
_{i}\right\vert ^{2}-r_{2}^{2}}}\cdot \sqrt{\sum_{i=1}^{\infty }\left\vert
\left\langle x,e_{i}\right\rangle \right\vert ^{2}\sum_{i=1}^{\infty
}\left\vert \left\langle y,e_{i}\right\rangle \right\vert ^{2}}  \notag \\
& \leq \frac{r_{1}r_{2}\left\Vert x\right\Vert \left\Vert y\right\Vert }{%
\sqrt{\sum_{i=1}^{\infty }\left\vert \lambda _{i}\right\vert ^{2}-r_{1}^{2}}%
\sqrt{\sum_{i=1}^{\infty }\left\vert \mu _{i}\right\vert ^{2}-r_{2}^{2}}}. 
\notag
\end{align}
\end{theorem}

\begin{proof}
Applying Schwarz's inequality for the vectors $x-\sum_{i=1}^{\infty
}\left\langle x,e_{i}\right\rangle e_{i},$ $y-\sum_{i=1}^{\infty
}\left\langle y,e_{i}\right\rangle e_{i},$ we have%
\begin{multline}
\left\vert \left\langle x-\sum_{i=1}^{\infty }\left\langle
x,e_{i}\right\rangle e_{i},y-\sum_{i=1}^{\infty }\left\langle
y,e_{i}\right\rangle e_{i}\right\rangle \right\vert ^{2}  \label{6.4.11a} \\
\leq \left\Vert x-\sum_{i=1}^{\infty }\left\langle x,e_{i}\right\rangle
e_{i}\right\Vert ^{2}\left\Vert y-\sum_{i=1}^{\infty }\left\langle
y,e_{i}\right\rangle e_{i}\right\Vert ^{2}.
\end{multline}%
Since%
\begin{equation*}
\left\langle x-\sum_{i=1}^{\infty }\left\langle x,e_{i}\right\rangle
e_{i},y-\sum_{i=1}^{\infty }\left\langle y,e_{i}\right\rangle
e_{i}\right\rangle =\left\langle x,y\right\rangle -\sum_{i=1}^{\infty
}\left\langle x,e_{i}\right\rangle \left\langle e_{i},y\right\rangle
\end{equation*}%
and%
\begin{equation*}
\left\Vert x-\sum_{i=1}^{\infty }\left\langle x,e_{i}\right\rangle
e_{i}\right\Vert ^{2}=\left\Vert x\right\Vert ^{2}-\sum_{i=1}^{\infty
}\left\vert \left\langle x,e_{i}\right\rangle \right\vert ^{2},
\end{equation*}%
then by (\ref{5.5.11a}) applied for $x$ and $y,$ and from (\ref{6.4.11a}),
we deduce the first part of (\ref{6.3.11a}).

The second part follows by Bessel's inequality.
\end{proof}

The following Gr\"{u}ss type inequality may be stated as well.

\begin{theorem}
\label{t6.2.11a}Let $\left( H;\left\langle \cdot ,\cdot \right\rangle
\right) $ be an infinite dimensional Hilbert space and $\left( e_{i}\right)
_{i\in \mathbb{N}}$ an orthornormal family in $H.$ Suppose that $\left(
\Gamma _{i}\right) _{i\in \mathbb{N}},$ $\left( \gamma _{i}\right) _{i\in 
\mathbb{N}},$ $\left( \phi _{i}\right) _{i\in \mathbb{N}},$ $\left( \Phi
_{i}\right) _{i\in \mathbb{N}}\in \ell ^{2}\left( \mathbb{K}\right) $ are
sequences of real and complex numbers such that%
\begin{equation*}
\sum_{i=1}^{\infty }\func{Re}\left( \Gamma _{i}\overline{\gamma _{i}}\right)
>0,\ \ \ \sum_{i=1}^{\infty }\func{Re}\left( \Phi _{i}\overline{\phi _{i}}%
\right) >0.
\end{equation*}%
If $x,y\in H$ are such that either%
\begin{align*}
\left\Vert x-\sum_{i=1}^{\infty }\frac{\Gamma _{i}+\gamma _{i}}{2}\cdot
e_{i}\right\Vert & \leq \frac{1}{2}\left( \sum_{i=1}^{\infty }\left\vert
\Gamma _{i}-\gamma _{i}\right\vert ^{2}\right) ^{\frac{1}{2}}, \\
\left\Vert y-\sum_{i=1}^{\infty }\frac{\Phi _{i}+\phi _{i}}{2}\cdot
e_{i}\right\Vert & \leq \frac{1}{2}\left( \sum_{i=1}^{\infty }\left\vert
\Phi _{i}-\phi _{i}\right\vert ^{2}\right) ^{\frac{1}{2}},
\end{align*}%
or equivalently,%
\begin{align*}
\func{Re}\left\langle \sum_{i=1}^{\infty }\Gamma
_{i}e_{i}-x,x-\sum_{i=1}^{\infty }\gamma _{i}e_{i}\right\rangle & \geq 0, \\
\func{Re}\left\langle \sum_{i=1}^{\infty }\Phi
_{i}e_{i}-y,y-\sum_{i=1}^{\infty }\phi _{i}e_{i}\right\rangle & \geq 0,
\end{align*}%
holds, then we have the inequality%
\begin{align*}
& \left\vert \left\langle x,y\right\rangle -\sum_{i=1}^{\infty }\left\langle
x,e_{i}\right\rangle \left\langle e_{i},y\right\rangle \right\vert \\
& \leq \frac{1}{4}\cdot \frac{\left( \sum_{i=1}^{\infty }\left\vert \Gamma
_{i}-\gamma _{i}\right\vert ^{2}\right) ^{\frac{1}{2}}\left(
\sum_{i=1}^{\infty }\left\vert \Phi _{i}-\phi _{i}\right\vert ^{2}\right) ^{%
\frac{1}{2}}}{\left( \sum_{i=1}^{\infty }\func{Re}\left( \Gamma _{i}%
\overline{\gamma _{i}}\right) \right) ^{\frac{1}{2}}\left(
\sum_{i=1}^{\infty }\func{Re}\left( \Phi _{i}\overline{\phi _{i}}\right)
\right) ^{\frac{1}{2}}} \\
& \qquad \times \left( \sum_{i=1}^{\infty }\left\vert \left\langle
x,e_{i}\right\rangle \right\vert ^{2}\right) ^{\frac{1}{2}}\left(
\sum_{i=1}^{\infty }\left\vert \left\langle y,e_{i}\right\rangle \right\vert
^{2}\right) ^{\frac{1}{2}} \\
& \leq \frac{1}{4}\cdot \frac{\left( \sum_{i=1}^{\infty }\left\vert \Gamma
_{i}-\gamma _{i}\right\vert ^{2}\right) ^{\frac{1}{2}}\left(
\sum_{i=1}^{\infty }\left\vert \Phi _{i}-\phi _{i}\right\vert ^{2}\right) ^{%
\frac{1}{2}}}{\left[ \sum_{i=1}^{\infty }\func{Re}\left( \Gamma _{i}%
\overline{\gamma _{i}}\right) \right] ^{\frac{1}{2}}\left[
\sum_{i=1}^{\infty }\func{Re}\left( \Phi _{i}\overline{\phi _{i}}\right) %
\right] ^{\frac{1}{2}}}\left\Vert x\right\Vert \left\Vert y\right\Vert
\end{align*}%
The constant $\frac{1}{4}$ is best possible in the first inequality.
\end{theorem}

\begin{proof}
Follows by (\ref{5.11.11a}) and (\ref{6.4.11a}). The best constant follows
from Theorem \ref{t4.2.7a}, and we omit the details.
\end{proof}

\subsection{Other Reverses of Bessel's Inequality}

We may state the following result \cite{12NSSD}.

\begin{theorem}
\label{t5.1.11b}Let $\left( H;\left\langle \cdot ,\cdot \right\rangle
\right) $ be an infinite dimensional Hilbert space over the real or complex
number field $\mathbb{K}$, $\left( e_{i}\right) _{i\in \mathbb{N}}$ is an
orthornormal family in $H,$ $\mathbf{\lambda }=\left( \lambda _{i}\right)
_{i\in \mathbb{N}}\in \ell ^{2}\left( \mathbb{K}\right) ,$ $\mathbf{\lambda }%
\neq 0$ and $r>0.$ If $x\in H$ is such that%
\begin{equation*}
\left\Vert x-\sum_{i=1}^{\infty }\lambda _{i}e_{i}\right\Vert \leq r,
\end{equation*}%
then we have the inequality%
\begin{equation}
0\leq \left\Vert x\right\Vert -\left( \sum_{i=1}^{\infty }\left\vert
\left\langle x,e_{i}\right\rangle \right\vert ^{2}\right) ^{\frac{1}{2}}\leq 
\frac{1}{2}\cdot \frac{r^{2}}{\left( \sum_{i=1}^{\infty }\left\vert \lambda
_{i}\right\vert ^{2}\right) ^{\frac{1}{2}}}.  \label{5.3.11b}
\end{equation}%
The constant $\frac{1}{2}$ is best possible in (\ref{5.3.11b}) in the sense
that it cannot be replaced by a smaller constant.
\end{theorem}

\begin{proof}
Let $a:=\sum_{i=1}^{\infty }\lambda _{i}e_{i}\in H.$ Then by Theorem \ref%
{t2.1.4}, we have 
\begin{equation*}
\left\Vert x\right\Vert \left\Vert \sum_{i=1}^{\infty }\lambda
_{i}e_{i}\right\Vert -\left\vert \sum_{i=1}^{\infty }\bar{\lambda}%
_{i}\left\langle x,e_{i}\right\rangle \right\vert \leq \frac{1}{2}r^{2},
\end{equation*}%
giving%
\begin{equation}
\left\Vert x\right\Vert \left( \sum_{i=1}^{\infty }\left\vert \lambda
_{i}\right\vert ^{2}\right) ^{\frac{1}{2}}\leq \frac{1}{2}r^{2}+\left\vert
\sum_{i=1}^{\infty }\bar{\lambda}_{i}\left\langle x,e_{i}\right\rangle
\right\vert ,  \label{5.4.11b}
\end{equation}%
since%
\begin{equation*}
\left\Vert \sum_{i=1}^{\infty }\lambda _{i}e_{i}\right\Vert =\left(
\sum_{i=1}^{\infty }\left\vert \lambda _{i}\right\vert ^{2}\right) ^{\frac{1%
}{2}}.
\end{equation*}%
Using the Cauchy-Bunyakovsky-Schwarz inequality, we may state that%
\begin{equation}
\left\vert \sum_{i=1}^{\infty }\bar{\lambda}_{i}\left\langle
x,e_{i}\right\rangle \right\vert \leq \left( \sum_{i=1}^{\infty }\left\vert
\lambda _{i}\right\vert ^{2}\right) ^{\frac{1}{2}}\left( \sum_{i=1}^{\infty
}\left\vert \left\langle x,e_{i}\right\rangle \right\vert ^{2}\right) ^{%
\frac{1}{2}},  \label{5.5.11b}
\end{equation}%
and thus, by (\ref{5.4.11b}) and (\ref{5.5.11b}), we may state that%
\begin{equation*}
\left\Vert x\right\Vert \left( \sum_{i=1}^{\infty }\left\vert \lambda
_{i}\right\vert ^{2}\right) ^{\frac{1}{2}}\leq \frac{1}{2}r^{2}+\left(
\sum_{i=1}^{\infty }\left\vert \lambda _{i}\right\vert ^{2}\right) ^{\frac{1%
}{2}}\left( \sum_{i=1}^{\infty }\left\vert \left\langle x,e_{i}\right\rangle
\right\vert ^{2}\right) ^{\frac{1}{2}},
\end{equation*}%
from where we get the desired inequality in (\ref{5.3.11b}).

The best constant, follows by Theorem \ref{t2.1.4} on choosing $\left(
e_{i}\right) _{i\in \mathbb{N}}=\left\{ e\right\} ,$ with $\left\Vert
e\right\Vert =1$ and we omit the details.
\end{proof}

\begin{remark}
Under the assumptions of Theorem \ref{t5.1.11b}, and if we multiply by $%
\left\Vert x\right\Vert +\left( \sum_{i=1}^{\infty }\left\vert \left\langle
x,e_{i}\right\rangle \right\vert ^{2}\right) ^{\frac{1}{2}}>0,$ then we
deduce, from (\ref{5.3.11b}), that%
\begin{align}
0& \leq \left\Vert x\right\Vert ^{2}-\sum_{i=1}^{\infty }\left\vert
\left\langle x,e_{i}\right\rangle \right\vert ^{2}  \label{5.6.11b} \\
& \leq \frac{1}{2}\cdot \frac{r^{2}\left( \left\Vert x\right\Vert +\left(
\sum_{i=1}^{\infty }\left\vert \left\langle x,e_{i}\right\rangle \right\vert
^{2}\right) ^{\frac{1}{2}}\right) }{\left( \sum_{i=1}^{\infty }\left\vert
\lambda _{i}\right\vert ^{2}\right) ^{\frac{1}{2}}}  \notag \\
& \leq \frac{r^{2}\left\Vert x\right\Vert }{\left( \sum_{i=1}^{\infty
}\left\vert \lambda _{i}\right\vert ^{2}\right) ^{\frac{1}{2}}},  \notag
\end{align}%
where for the last inequality, we have used Bessel's inequality%
\begin{equation*}
\left( \sum_{i=1}^{\infty }\left\vert \left\langle x,e_{i}\right\rangle
\right\vert ^{2}\right) ^{\frac{1}{2}}\leq \left\Vert x\right\Vert ,\ \ \ \
x\in H.
\end{equation*}
\end{remark}

The following result also holds \cite{12NSSD}.

\begin{theorem}
\label{t5.2.11b}Assume that $\left( H;\left\langle \cdot ,\cdot
\right\rangle \right) $ and $\left( e_{i}\right) _{i\in \mathbb{N}}$ are as
in Theorem \ref{t5.1.11b}. If $\mathbf{\Gamma }=\left( \Gamma _{i}\right)
_{i\in \mathbb{N}},$ $\mathbf{\gamma }=\left( \gamma _{i}\right) _{i\in 
\mathbb{N}}\in \ell ^{2}\left( \mathbb{K}\right) ,$ with $\mathbf{\Gamma }%
\neq -\mathbf{\gamma }$, and $x\in H$ with the property that, either%
\begin{equation*}
\left\Vert x-\sum_{i=1}^{\infty }\frac{\Gamma _{i}+\gamma _{i}}{2}\cdot
e_{i}\right\Vert \leq \frac{1}{2}\left( \sum_{i=1}^{\infty }\left\vert
\Gamma _{i}-\gamma _{i}\right\vert ^{2}\right) ^{\frac{1}{2}},
\end{equation*}%
or equivalently,%
\begin{equation*}
\func{Re}\left\langle \sum_{i=1}^{\infty }\Gamma
_{i}e_{i}-x,x-\sum_{i=1}^{\infty }\gamma _{i}e_{i}\right\rangle \geq 0,
\end{equation*}%
holds, then we have the inequality%
\begin{equation}
0\leq \left\Vert x\right\Vert -\left( \sum_{i=1}^{\infty }\left\vert
\left\langle x,e_{i}\right\rangle \right\vert ^{2}\right) ^{\frac{1}{2}}\leq 
\frac{1}{4}\cdot \frac{\sum_{i=1}^{\infty }\left\vert \Gamma _{i}-\gamma
_{i}\right\vert ^{2}}{\left( \sum_{i=1}^{\infty }\left\vert \Gamma
_{i}+\gamma _{i}\right\vert ^{2}\right) ^{\frac{1}{2}}}.  \label{5.9.11b}
\end{equation}%
The constant $\frac{1}{4}$ is best possible in the sense that it cannot be
replaced by a smaller constant.
\end{theorem}

\begin{proof}
Since $\mathbf{\Gamma },$ $\mathbf{\gamma }\in \ell ^{2}\left( \mathbb{K}%
\right) ,$ then we have that $\frac{1}{2}\left( \mathbf{\Gamma }\pm \mathbf{%
\gamma }\right) \in \ell ^{2}\left( \mathbb{K}\right) ,$ showing that the
series%
\begin{equation*}
\sum_{i=1}^{\infty }\left\vert \frac{\Gamma _{i}+\gamma _{i}}{2}\right\vert
^{2},\ \ \ \sum_{i=1}^{\infty }\left\vert \frac{\Gamma _{i}-\gamma _{i}}{2}%
\right\vert ^{2}
\end{equation*}%
are convergent. In addition, the series $\sum_{i=1}^{\infty }\Gamma
_{i}e_{i} $, $\sum_{i=1}^{\infty }\gamma _{i}e_{i}$ and $\sum_{i=1}^{\infty }%
\frac{\Gamma _{i}+\gamma _{i}}{2}e_{i}$ are also convergent in the Hilbert
space $H.$

Now, we observe that the inequality (\ref{5.9.11b}) follows from Theorem \ref%
{t5.1.11b} on choosing $\lambda _{i}=\frac{\Gamma _{i}+\gamma _{i}}{2},$ $%
i\in \mathbb{N}$ and $r=\frac{1}{2}\left( \sum_{i=1}^{\infty }\left\vert
\Gamma _{i}-\gamma _{i}\right\vert ^{2}\right) ^{\frac{1}{2}}.$

The fact that $\frac{1}{4}$ is the best possible constant in (\ref{5.9.11b})
follows from Theorem \ref{t2.2.4}, and we omit the details.
\end{proof}

\begin{remark}
With the assumptions of Theorem \ref{t5.2.11b}, we have 
\begin{align}
0& \leq \left\Vert x\right\Vert ^{2}-\sum_{i=1}^{\infty }\left\vert
\left\langle x,e_{i}\right\rangle \right\vert ^{2}  \label{5.10.11b} \\
& \leq \frac{1}{4}\cdot \frac{\sum_{i=1}^{\infty }\left\vert \Gamma
_{i}-\gamma _{i}\right\vert ^{2}}{\left( \sum_{i=1}^{\infty }\left\vert
\Gamma _{i}+\gamma _{i}\right\vert ^{2}\right) ^{\frac{1}{2}}}\left[
\left\Vert x\right\Vert +\left( \sum_{i=1}^{\infty }\left\vert \left\langle
x,e_{i}\right\rangle \right\vert ^{2}\right) ^{\frac{1}{2}}\right]  \notag \\
& \leq \frac{1}{2}\cdot \frac{\sum_{i=1}^{\infty }\left\vert \Gamma
_{i}-\gamma _{i}\right\vert ^{2}}{\left( \sum_{i=1}^{\infty }\left\vert
\Gamma _{i}+\gamma _{i}\right\vert ^{2}\right) ^{\frac{1}{2}}}\left\Vert
x\right\Vert .  \notag
\end{align}
\end{remark}

\subsection{More Gr\"{u}ss Type Inequalities for Orthonormal Families}

The following result holds \cite{12NSSD}.

\begin{theorem}
\label{t6.1.11b}Let $\left( H;\left\langle \cdot ,\cdot \right\rangle
\right) $ be an infinite dimensional Hilbert space over the real or complex
number field $\mathbb{K}$ and $\left( e_{i}\right) _{i\in \mathbb{N}}$ an
orthornormal family in $H.$ If $\mathbf{\lambda }=\left( \lambda _{i}\right)
_{i\in \mathbb{N}},\ \mathbf{\mu }=\left( \mu _{i}\right) _{i\in \mathbb{N}%
}\in \ell ^{2}\left( \mathbb{K}\right) ,$ $\mathbf{\lambda }$, $\mathbf{\mu }%
\neq 0,$ $r_{1},r_{2}>0$ and $x,y\in H$ are such that%
\begin{equation*}
\left\Vert x-\sum_{i=1}^{\infty }\lambda _{i}e_{i}\right\Vert \leq r_{1},\ \
\ \left\Vert y-\sum_{i=1}^{\infty }\mu _{i}e_{i}\right\Vert \leq r_{2},
\end{equation*}%
then we have the inequality%
\begin{align*}
& \left\vert \left\langle x,y\right\rangle -\sum_{i=1}^{\infty }\left\langle
x,e_{i}\right\rangle \left\langle e_{i},y\right\rangle \right\vert \\
& \leq \frac{1}{2}r_{1}r_{2}\frac{\left[ \left\Vert x\right\Vert +\left(
\sum_{i=1}^{\infty }\left\vert \left\langle x,e_{i}\right\rangle \right\vert
^{2}\right) ^{\frac{1}{2}}\right] ^{\frac{1}{2}}\left[ \left\Vert
y\right\Vert +\left( \sum_{i=1}^{\infty }\left\vert \left\langle
y,e_{i}\right\rangle \right\vert ^{2}\right) ^{\frac{1}{2}}\right] ^{\frac{1%
}{2}}}{\left( \sum_{i=1}^{\infty }\left\vert \lambda _{i}\right\vert
^{2}\right) ^{\frac{1}{4}}\left( \sum_{i=1}^{\infty }\left\vert \mu
_{i}\right\vert ^{2}\right) ^{\frac{1}{4}}} \\
& \leq r_{1}r_{2}\frac{\left\Vert x\right\Vert ^{\frac{1}{2}}\left\Vert
y\right\Vert ^{\frac{1}{2}}}{\left( \sum_{i=1}^{\infty }\left\vert \lambda
_{i}\right\vert ^{2}\right) ^{\frac{1}{4}}\left( \sum_{i=1}^{\infty
}\left\vert \mu _{i}\right\vert ^{2}\right) ^{\frac{1}{4}}}.
\end{align*}
\end{theorem}

\begin{proof}
It follows by (\ref{5.6.11b}) applied for $x$ and $y.$ We omit the details.
\end{proof}

Finally we may state the following theorem \cite{12NSSD}.

\begin{theorem}
\label{t6.2.11b}Assume that $\left( H;\left\langle \cdot ,\cdot
\right\rangle \right) $ and $\left( e_{i}\right) _{i\in \mathbb{N}}$ are as
in Theorem \ref{t6.1.11b}. If $\mathbf{\Gamma }=\left( \Gamma _{i}\right)
_{i\in \mathbb{N}},$ $\mathbf{\Gamma }=\left( \Gamma _{i}\right) _{i\in 
\mathbb{N}},$ $\phi =\left( \mathbf{\phi }_{i}\right) _{i\in \mathbb{N}},$ $%
\mathbf{\Phi }=\left( \Phi _{i}\right) _{i\in \mathbb{N}}\in \ell ^{2}\left( 
\mathbb{K}\right) ,$ with $\mathbf{\Gamma }\neq -\mathbf{\gamma }$, $\mathbf{%
\Phi }\neq -\mathbf{\phi }$, and $x,y\in H$ are such that, either%
\begin{align*}
\left\Vert x-\sum_{i=1}^{\infty }\frac{\Gamma _{i}+\gamma _{i}}{2}\cdot
e_{i}\right\Vert & \leq \frac{1}{2}\left( \sum_{i=1}^{\infty }\left\vert
\Gamma _{i}-\gamma _{i}\right\vert ^{2}\right) ^{\frac{1}{2}}, \\
\left\Vert y-\sum_{i=1}^{\infty }\frac{\Phi _{i}+\phi _{i}}{2}\cdot
e_{i}\right\Vert & \leq \frac{1}{2}\left( \sum_{i=1}^{\infty }\left\vert
\Phi _{i}-\phi _{i}\right\vert ^{2}\right) ^{\frac{1}{2}},
\end{align*}%
or equivalently,%
\begin{align*}
\func{Re}\left\langle \sum_{i=1}^{\infty }\Gamma
_{i}e_{i}-x,x-\sum_{i=1}^{\infty }\gamma _{i}e_{i}\right\rangle & \geq 0, \\
\func{Re}\left\langle \sum_{i=1}^{\infty }\Phi
_{i}e_{i}-y,y-\sum_{i=1}^{\infty }\phi _{i}e_{i}\right\rangle & \geq 0,
\end{align*}%
holds, then we have the inequality%
\begin{align*}
& \left\vert \left\langle x,y\right\rangle -\sum_{i=1}^{\infty }\left\langle
x,e_{i}\right\rangle \left\langle e_{i},y\right\rangle \right\vert \\
& \leq \frac{1}{4}\cdot \left( \sum_{i=1}^{\infty }\left\vert \Phi _{i}-\phi
_{i}\right\vert ^{2}\right) ^{\frac{1}{2}}\left( \sum_{i=1}^{\infty
}\left\vert \Gamma _{i}-\gamma _{i}\right\vert ^{2}\right) ^{\frac{1}{2}} \\
& \qquad \qquad \times \frac{\left[ \left\Vert x\right\Vert +\left(
\sum_{i=1}^{\infty }\left\vert \left\langle x,e_{i}\right\rangle \right\vert
^{2}\right) ^{\frac{1}{2}}\right] ^{\frac{1}{2}}\left[ \left\Vert
y\right\Vert +\left( \sum_{i=1}^{\infty }\left\vert \left\langle
y,e_{i}\right\rangle \right\vert ^{2}\right) ^{\frac{1}{2}}\right] ^{\frac{1%
}{2}}}{\left( \sum_{i=1}^{\infty }\left\vert \Phi _{i}+\phi _{i}\right\vert
^{2}\right) ^{\frac{1}{4}}\left( \sum_{i=1}^{\infty }\left\vert \Gamma
_{i}+\gamma _{i}\right\vert ^{2}\right) ^{\frac{1}{4}}} \\
& \leq \frac{1}{2}\cdot \frac{\left( \sum_{i=1}^{\infty }\left\vert \Phi
_{i}-\phi _{i}\right\vert ^{2}\right) ^{\frac{1}{2}}\left(
\sum_{i=1}^{\infty }\left\vert \Gamma _{i}-\gamma _{i}\right\vert
^{2}\right) ^{\frac{1}{2}}}{\left( \sum_{i=1}^{\infty }\left\vert \Phi
_{i}+\phi _{i}\right\vert ^{2}\right) ^{\frac{1}{4}}\left(
\sum_{i=1}^{\infty }\left\vert \Gamma _{i}+\gamma _{i}\right\vert
^{2}\right) ^{\frac{1}{4}}}\left\Vert x\right\Vert ^{\frac{1}{2}}\left\Vert
y\right\Vert ^{\frac{1}{2}}.
\end{align*}
\end{theorem}

The proof follow by (\ref{5.10.11b}) applied for $x$ and $y.$ We omit the
details.

\newpage

\chapter{Generalisations of Bessel's Inequality\label{chap4}}

\section{Boas-Bellman Type Inequalities}

\subsection{Introduction}

Let $\left( H;\left( \cdot ,\cdot \right) \right) $ be an inner product
space over the real or complex number field $\mathbb{K}$. If $\left(
e_{i}\right) _{1\leq i\leq n}$ are orthonormal vectors in the inner product
space $H,$ i.e., $\left( e_{i},e_{j}\right) =\delta _{ij}$ for all $i,j\in
\left\{ 1,\dots ,n\right\} ,$ where $\delta _{ij}$ is the Kronecker delta,
then we have the following inequality is well known in the literature as
Bessel's inequality (see for example \cite[p. 391]{6b.12}): 
\begin{equation*}
\sum_{i=1}^{n}\left\vert \left( x,e_{i}\right) \right\vert ^{2}\leq
\left\Vert x\right\Vert ^{2}\text{ \ for any \ }x\in H.
\end{equation*}

For other results related to Bessel's inequality, see \cite{3b.12} -- \cite%
{5b.12} and Chapter XV in the book \cite{6b.12}.

In 1941, R.P. Boas \cite{2b.12} and in 1944, independently, R. Bellman \cite%
{1b.12} proved the following generalisation of Bessel's inequality (see also 
\cite[p. 392]{6b.12}).

\begin{theorem}
\label{t1.1.12}If $x,y_{1},\dots ,y_{n}$ are elements of an inner product
space $\left( H;\left( \cdot ,\cdot \right) \right) ,$ then the following
inequality: 
\begin{equation}
\sum_{i=1}^{n}\left\vert \left( x,y_{i}\right) \right\vert ^{2}\leq
\left\Vert x\right\Vert ^{2}\left[ \max_{1\leq i\leq n}\left\Vert
y_{i}\right\Vert ^{2}+\left( \sum_{1\leq i\neq j\leq n}\left\vert \left(
y_{i},y_{j}\right) \right\vert ^{2}\right) ^{\frac{1}{2}}\right] ,
\label{1.2.12}
\end{equation}%
holds.
\end{theorem}

A recent generalisation of the Boas-Bellman result was given in Mitrinovi%
\'{c}-Pe\v{c}ari\'{c}-Fink \cite[p. 392]{6b.12} where they proved the
following:

\begin{theorem}
\label{t1.2.12}If $x,y_{1},\dots ,y_{n}$ are as in Theorem \ref{t1.1.12} and 
$c_{1},\dots ,c_{n}\in \mathbb{K}$, then one has the inequality: 
\begin{multline}
\left\vert \sum_{i=1}^{n}c_{i}\left( x,y_{i}\right) \right\vert ^{2}
\label{1.3.12} \\
\leq \left\Vert x\right\Vert ^{2}\sum_{i=1}^{n}\left\vert c_{i}\right\vert
^{2}\left[ \max_{1\leq i\leq n}\left\Vert y_{i}\right\Vert ^{2}+\left(
\sum_{1\leq i\neq j\leq n}\left\vert \left( y_{i},y_{j}\right) \right\vert
^{2}\right) ^{\frac{1}{2}}\right] .
\end{multline}
\end{theorem}

They also noted that if in (\ref{1.3.12}) one chooses $c_{i}=\overline{%
\left( x,y_{i}\right) }$, then this inequality becomes (\ref{1.2.12}).

For other results related to the Boas-Bellman inequality, see \cite{4b.12}.

In this section, by following \cite{13NSSD}, we point out some new results
that may be related to both the Mitrinovi\'{c}-Pe\v{c}ari\'{c}-Fink and
Boas-Bellman inequalities.

\subsection{Preliminary Results}

We start with the following lemma which is also interesting in itself \cite%
{13NSSD}.

\begin{lemma}
\label{l2.1.12}Let $z_{1},\dots ,z_{n}\in H$ and $\alpha _{1},\dots ,\alpha
_{n}\in \mathbb{K}$. Then one has the inequality: 
\begin{multline}
\left\Vert \sum_{i=1}^{n}\alpha _{i}z_{i}\right\Vert ^{2}\leq \left\{ 
\begin{array}{l}
\max\limits_{1\leq i\leq n}\left\vert \alpha _{i}\right\vert
^{2}\sum\limits_{i=1}^{n}\left\Vert z_{i}\right\Vert ^{2}; \\ 
\\ 
\left( \sum\limits_{i=1}^{n}\left\vert \alpha _{i}\right\vert ^{2\alpha
}\right) ^{\frac{1}{\alpha }}\left( \sum\limits_{i=1}^{n}\left\Vert
z_{i}\right\Vert ^{2\beta }\right) ^{\frac{1}{\beta }},\ \ \ \text{where \ }%
\alpha >1, \\ 
\hfill \frac{1}{\alpha }+\frac{1}{\beta }=1; \\ 
\sum\limits_{i=1}^{n}\left\vert \alpha _{i}\right\vert
^{2}\max\limits_{1\leq i\leq n}\left\Vert z_{i}\right\Vert ^{2},%
\end{array}%
\right.  \label{2.1.12} \\
+\left\{ 
\begin{array}{l}
\max\limits_{1\leq i\neq j\leq n}\left\{ \left\vert \alpha _{i}\alpha
_{j}\right\vert \right\} \sum\limits_{1\leq i\neq j\leq n}\left\vert \left(
z_{i},z_{j}\right) \right\vert ; \\ 
\\ 
\left[ \left( \sum\limits_{i=1}^{n}\left\vert \alpha _{i}\right\vert
^{\gamma }\right) ^{2}-\left( \sum\limits_{i=1}^{n}\left\vert \alpha
_{i}\right\vert ^{2\gamma }\right) \right] ^{\frac{1}{\gamma }}\left(
\sum\limits_{1\leq i\neq j\leq n}\left\vert \left( z_{i},z_{j}\right)
\right\vert ^{\delta }\right) ^{\frac{1}{\delta }}, \\ 
\hfill \ \ \ \text{where \ }\gamma >1,\ \ \frac{1}{\gamma }+\frac{1}{\delta }%
=1; \\ 
\\ 
\left[ \left( \sum\limits_{i=1}^{n}\left\vert \alpha _{i}\right\vert \right)
^{2}-\sum\limits_{i=1}^{n}\left\vert \alpha _{i}\right\vert ^{2}\right]
\max\limits_{1\leq i\neq j\leq n}\left\vert \left( z_{i},z_{j}\right)
\right\vert .%
\end{array}%
\right.
\end{multline}
\end{lemma}

\begin{proof}
We observe that 
\begin{align}
\left\Vert \sum_{i=1}^{n}\alpha _{i}z_{i}\right\Vert ^{2}& =\left(
\sum_{i=1}^{n}\alpha _{i}z_{i},\sum_{j=1}^{n}\alpha _{j}z_{j}\right)
\label{2.2.12} \\
& =\sum_{i=1}^{n}\sum_{j=1}^{n}\alpha _{i}\overline{\alpha _{j}}\left(
z_{i},z_{j}\right) =\left\vert \sum_{i=1}^{n}\sum_{j=1}^{n}\alpha _{i}%
\overline{\alpha _{j}}\left( z_{i},z_{j}\right) \right\vert  \notag \\
& \leq \sum_{i=1}^{n}\sum_{j=1}^{n}\left\vert \alpha _{i}\right\vert
\left\vert \overline{\alpha _{j}}\right\vert \left\vert \left(
z_{i},z_{j}\right) \right\vert  \notag \\
& =\sum_{i=1}^{n}\left\vert \alpha _{i}\right\vert ^{2}\left\Vert
z_{i}\right\Vert ^{2}+\sum\limits_{1\leq i\neq j\leq n}\left\vert \alpha
_{i}\right\vert \left\vert \alpha _{j}\right\vert \left\vert \left(
z_{i},z_{j}\right) \right\vert .  \notag
\end{align}%
Using H\"{o}lder's inequality, we may write that 
\begin{multline}
\sum_{i=1}^{n}\left\vert \alpha _{i}\right\vert ^{2}\left\Vert
z_{i}\right\Vert ^{2}  \label{2.3.12} \\
\leq \left\{ 
\begin{array}{l}
\max\limits_{1\leq i\leq n}\left\vert \alpha _{i}\right\vert
^{2}\sum\limits_{i=1}^{n}\left\Vert z_{i}\right\Vert ^{2}; \\ 
\\ 
\left( \sum\limits_{i=1}^{n}\left\vert \alpha _{i}\right\vert ^{2\alpha
}\right) ^{\frac{1}{\alpha }}\left( \sum\limits_{i=1}^{n}\left\Vert
z_{i}\right\Vert ^{2\beta }\right) ^{\frac{1}{\beta }},\ \ \ \text{where \ }%
\alpha >1,\frac{1}{\alpha }+\frac{1}{\beta }=1; \\ 
\\ 
\sum\limits_{i=1}^{n}\left\vert \alpha _{i}\right\vert
^{2}\max\limits_{1\leq i\leq n}\left\Vert z_{i}\right\Vert ^{2}.%
\end{array}%
\right.
\end{multline}%
By H\"{o}lder's inequality for double sums we also have 
\begin{multline}
\sum\limits_{1\leq i\neq j\leq n}\left\vert \alpha _{i}\right\vert
\left\vert \alpha _{j}\right\vert \left\vert \left( z_{i},z_{j}\right)
\right\vert  \label{2.4.12} \\
\leq \left\{ 
\begin{array}{l}
\max\limits_{1\leq i\neq j\leq n}\left\vert \alpha _{i}\alpha
_{j}\right\vert \sum\limits_{1\leq i\neq j\leq n}\left\vert \left(
z_{i},z_{j}\right) \right\vert ; \\ 
\\ 
\left( \sum\limits_{1\leq i\neq j\leq n}\left\vert \alpha _{i}\right\vert
^{\gamma }\left\vert \alpha _{j}\right\vert ^{\gamma }\right) ^{\frac{1}{%
\gamma }}\left( \sum\limits_{1\leq i\neq j\leq n}\left\vert \left(
z_{i},z_{j}\right) \right\vert ^{\delta }\right) ^{\frac{1}{\delta }}, \\ 
\hfill \ \ \ \text{where \ }\gamma >1,\ \ \frac{1}{\gamma }+\frac{1}{\delta }%
=1; \\ 
\\ 
\sum\limits_{1\leq i\neq j\leq n}\left\vert \alpha _{i}\right\vert
\left\vert \alpha _{j}\right\vert \max\limits_{1\leq i\neq j\leq
n}\left\vert \left( z_{i},z_{j}\right) \right\vert ,%
\end{array}%
\right.
\end{multline}%
\begin{equation*}
=\left\{ 
\begin{array}{l}
\max\limits_{1\leq i\neq j\leq n}\left\{ \left\vert \alpha _{i}\alpha
_{j}\right\vert \right\} \sum\limits_{1\leq i\neq j\leq n}\left\vert \left(
z_{i},z_{j}\right) \right\vert ; \\ 
\\ 
\left[ \left( \sum\limits_{i=1}^{n}\left\vert \alpha _{i}\right\vert
^{\gamma }\right) ^{2}-\left( \sum\limits_{i=1}^{n}\left\vert \alpha
_{i}\right\vert ^{2\gamma }\right) \right] ^{\frac{1}{\gamma }}\left(
\sum\limits_{1\leq i\neq j\leq n}\left\vert \left( z_{i},z_{j}\right)
\right\vert ^{\delta }\right) ^{\frac{1}{\delta }}, \\ 
\hfill \ \ \ \text{where \ }\gamma >1,\ \ \frac{1}{\gamma }+\frac{1}{\delta }%
=1; \\ 
\\ 
\left[ \left( \sum\limits_{i=1}^{n}\left\vert \alpha _{i}\right\vert \right)
^{2}-\sum\limits_{i=1}^{n}\left\vert \alpha _{i}\right\vert ^{2}\right]
\max\limits_{1\leq i\neq j\leq n}\left\vert \left( z_{i},z_{j}\right)
\right\vert .%
\end{array}%
\right.
\end{equation*}%
Utilising (\ref{2.3.12}) and (\ref{2.4.12}) in (\ref{2.2.12}), we may deduce
the desired result (\ref{2.1.12}).
\end{proof}

\begin{remark}
\label{r2.2.12}Inequality (\ref{2.1.12}) contains in fact 9 different
inequalities which may be obtained by combining the first 3 ones with the
last 3 ones.
\end{remark}

A particular case that may be related to the Boas-Bellman result is embodied
in the following inequality \cite{13NSSD}.

\begin{corollary}
\label{c2.3.12}With the assumptions in Lemma \ref{l2.1.12}, we have 
\begin{equation}
\left\Vert \sum_{i=1}^{n}\alpha _{i}z_{i}\right\Vert ^{2}  \label{2.5.12}
\end{equation}
\end{corollary}

\begin{align}
& \leq \sum_{i=1}^{n}\left| \alpha _{i}\right| ^{2}\left\{
\max\limits_{1\leq i\leq n}\left\| z_{i}\right\| ^{2}+\frac{\left[ \left(
\sum_{i=1}^{n}\left| \alpha _{i}\right| ^{2}\right)
^{2}-\sum_{i=1}^{n}\left| \alpha _{i}\right| ^{4}\right] ^{\frac{1}{2}}}{%
\sum_{i=1}^{n}\left| \alpha _{i}\right| ^{2}}\left( \sum\limits_{1\leq i\neq
j\leq n}\left| \left( z_{i},z_{j}\right) \right| ^{2}\right) ^{\frac{1}{2}%
}\right\}  \notag \\
& \leq \sum_{i=1}^{n}\left| \alpha _{i}\right| ^{2}\left\{
\max\limits_{1\leq i\leq n}\left\| z_{i}\right\| ^{2}+\left(
\sum\limits_{1\leq i\neq j\leq n}\left| \left( z_{i},z_{j}\right) \right|
^{2}\right) ^{\frac{1}{2}}\right\} .  \notag
\end{align}
The first inequality follows by taking the third branch in the first curly
bracket with the second branch in the second curly bracket for $\gamma
=\delta =2.$

The second inequality in (\ref{2.5.12}) follows by the fact that 
\begin{equation*}
\left[ \left( \sum_{i=1}^{n}\left\vert \alpha _{i}\right\vert ^{2}\right)
^{2}-\sum_{i=1}^{n}\left\vert \alpha _{i}\right\vert ^{4}\right] ^{\frac{1}{2%
}}\leq \sum_{i=1}^{n}\left\vert \alpha _{i}\right\vert ^{2}.
\end{equation*}%
Applying the following Cauchy-Bunyakovsky-Schwarz type inequality 
\begin{equation*}
\left( \sum_{i=1}^{n}a_{i}\right) ^{2}\leq n\sum_{i=1}^{n}a_{i}^{2},\ \ \
a_{i}\in \mathbb{R}_{+},\ \ 1\leq i\leq n,
\end{equation*}%
we may write that 
\begin{equation}
\left( \sum\limits_{i=1}^{n}\left\vert \alpha _{i}\right\vert ^{\gamma
}\right) ^{2}-\sum\limits_{i=1}^{n}\left\vert \alpha _{i}\right\vert
^{2\gamma }\leq \left( n-1\right) \sum\limits_{i=1}^{n}\left\vert \alpha
_{i}\right\vert ^{2\gamma }\ \ \ \ \ \ \left( n\geq 1\right)  \label{2.7.12}
\end{equation}%
and 
\begin{equation}
\left( \sum\limits_{i=1}^{n}\left\vert \alpha _{i}\right\vert \right)
^{2}-\sum\limits_{i=1}^{n}\left\vert \alpha _{i}\right\vert ^{2}\leq \left(
n-1\right) \sum\limits_{i=1}^{n}\left\vert \alpha _{i}\right\vert ^{2}\ \ \
\ \ \ \left( n\geq 1\right) .  \label{2.8.12}
\end{equation}%
Also, it is obvious that: 
\begin{equation}
\max\limits_{1\leq i\neq j\leq n}\left\{ \left\vert \alpha _{i}\alpha
_{j}\right\vert \right\} \leq \max\limits_{1\leq i\leq n}\left\vert \alpha
_{i}\right\vert ^{2}.  \label{2.9.12}
\end{equation}%
Consequently, we may state the following coarser upper bounds for $%
\left\Vert \sum_{i=1}^{n}\alpha _{i}z_{i}\right\Vert ^{2}$ that may be
useful in applications \cite{13NSSD}.

\begin{corollary}
\label{c2.4.12}With the assumptions in Lemma \ref{l2.1.12}, we have the
inequalities: 
\begin{multline}
\left\Vert \sum_{i=1}^{n}\alpha _{i}z_{i}\right\Vert ^{2}  \label{2.10.12} \\
\leq \left\{ 
\begin{array}{l}
\max\limits_{1\leq i\leq n}\left\vert \alpha _{i}\right\vert
^{2}\sum\limits_{i=1}^{n}\left\Vert z_{i}\right\Vert ^{2}; \\ 
\\ 
\left( \sum\limits_{i=1}^{n}\left\vert \alpha _{i}\right\vert ^{2\alpha
}\right) ^{\frac{1}{\alpha }}\left( \sum\limits_{i=1}^{n}\left\Vert
z_{i}\right\Vert ^{2\beta }\right) ^{\frac{1}{\beta }},\ \ \ \text{where \ }%
\alpha >1,\frac{1}{\alpha }+\frac{1}{\beta }=1; \\ 
\\ 
\sum\limits_{i=1}^{n}\left\vert \alpha _{i}\right\vert
^{2}\max\limits_{1\leq i\leq n}\left\Vert z_{i}\right\Vert ^{2},%
\end{array}%
\right. \\
+\left\{ 
\begin{array}{l}
\max\limits_{1\leq i\leq n}\left\vert \alpha _{i}\right\vert
^{2}\sum\limits_{1\leq i\neq j\leq n}\left\vert \left( z_{i},z_{j}\right)
\right\vert ; \\ 
\\ 
\left( n-1\right) ^{\frac{1}{\gamma }}\left( \sum\limits_{i=1}^{n}\left\vert
\alpha _{i}\right\vert ^{2\gamma }\right) ^{\frac{1}{\gamma }}\left(
\sum\limits_{1\leq i\neq j\leq n}\left\vert \left( z_{i},z_{j}\right)
\right\vert ^{\delta }\right) ^{\frac{1}{\delta }}, \\ 
\hfill \ \ \ \text{where \ }\gamma >1,\ \ \frac{1}{\gamma }+\frac{1}{\delta }%
=1; \\ 
\\ 
\left( n-1\right) \sum\limits_{i=1}^{n}\left\vert \alpha _{i}\right\vert
^{2}\max\limits_{1\leq i\neq j\leq n}\left\vert \left( z_{i},z_{j}\right)
\right\vert .%
\end{array}%
\right.
\end{multline}
\end{corollary}

The proof is obvious by Lemma \ref{l2.1.12} in applying the inequalities (%
\ref{2.7.12}) -- (\ref{2.9.12}).

\begin{remark}
\label{r2.5.12}The following inequalities which are incorporated in (\ref%
{2.10.12}) are of special interest: 
\begin{equation}
\left\Vert \sum_{i=1}^{n}\alpha _{i}z_{i}\right\Vert ^{2}\leq
\max\limits_{1\leq i\leq n}\left\vert \alpha _{i}\right\vert ^{2}\left[
\sum\limits_{i=1}^{n}\left\Vert z_{i}\right\Vert ^{2}+\sum\limits_{1\leq
i\neq j\leq n}\left\vert \left( z_{i},z_{j}\right) \right\vert \right] ;
\label{2.11.12}
\end{equation}%
\begin{multline}
\left\Vert \sum_{i=1}^{n}\alpha _{i}z_{i}\right\Vert ^{2}\leq \left(
\sum\limits_{i=1}^{n}\left\vert \alpha _{i}\right\vert ^{2p}\right) ^{\frac{1%
}{p}}\left[ \left( \sum\limits_{i=1}^{n}\left\Vert z_{i}\right\Vert
^{2q}\right) ^{\frac{1}{q}}\right.  \label{2.12.12} \\
+\left( n-1\right) \left. ^{\frac{1}{p}}\left( \sum\limits_{1\leq i\neq
j\leq n}\left\vert \left( z_{i},z_{j}\right) \right\vert ^{q}\right) ^{\frac{%
1}{q}}\right] ,
\end{multline}%
where $p>1,$ $\frac{1}{p}+\frac{1}{q}=1;$ and 
\begin{equation}
\left\Vert \sum_{i=1}^{n}\alpha _{i}z_{i}\right\Vert ^{2}\leq
\sum\limits_{i=1}^{n}\left\vert \alpha _{i}\right\vert ^{2}\left[
\max\limits_{1\leq i\leq n}\left\Vert z_{i}\right\Vert ^{2}+\left(
n-1\right) \max\limits_{1\leq i\neq j\leq n}\left\vert \left(
z_{i},z_{j}\right) \right\vert \right] .  \label{2.13.12}
\end{equation}
\end{remark}

\subsection{Mitrinovi\'{c}-Pe\v{c}ari\'{c}-Fink Type Inequalities}

We are now able to present the following result obtained in \cite{13NSSD},
which complements the inequality (\ref{1.3.12}) due to Mitrinovi\'{c}, Pe%
\v{c}ari\'{c} and Fink \cite[p. 392]{6b.12}.

\begin{theorem}
\label{t3.1.12}Let $x,y_{1},\dots ,y_{n}$ be vectors of an inner product
space $\left( H;\left( \cdot ,\cdot \right) \right) $ and $c_{1},\dots
,c_{n}\in \mathbb{K}$ $\left( \mathbb{K}=\mathbb{C},\mathbb{R}\right) .$
Then one has the inequalities: 
\begin{multline*}
\left\vert \sum_{i=1}^{n}c_{i}\left( x,y_{i}\right) \right\vert ^{2} \\
\leq \left\Vert x\right\Vert ^{2}\left\{ 
\begin{array}{l}
\max\limits_{1\leq i\leq n}\left\vert c_{i}\right\vert
^{2}\sum_{i=1}^{n}\left\Vert y_{i}\right\Vert ^{2}; \\ 
\\ 
\left( \sum_{i=1}^{n}\left\vert c_{i}\right\vert ^{2\alpha }\right) ^{\frac{1%
}{\alpha }}\left( \sum_{i=1}^{n}\left\Vert y_{i}\right\Vert ^{2\beta
}\right) ^{\frac{1}{\beta }},\ \ \ \text{where \ }\alpha >1, \\ 
\hfill \frac{1}{\alpha }+\frac{1}{\beta }=1; \\ 
\sum_{i=1}^{n}\left\vert c_{i}\right\vert ^{2}\max\limits_{1\leq i\leq
n}\left\Vert y_{i}\right\Vert ^{2},%
\end{array}%
\right.
\end{multline*}%
\begin{equation}
+\left\Vert x\right\Vert ^{2}\left\{ 
\begin{array}{l}
\max\limits_{1\leq i\neq j\leq n}\left\{ \left\vert c_{i}c_{j}\right\vert
\right\} \sum_{1\leq i\neq j\leq n}\left\vert \left( y_{i},y_{j}\right)
\right\vert ; \\ 
\\ 
\left[ \left( \sum_{i=1}^{n}\left\vert c_{i}\right\vert ^{\gamma }\right)
^{2}-\left( \sum_{i=1}^{n}\left\vert c_{i}\right\vert ^{2\gamma }\right) %
\right] ^{\frac{1}{\gamma }} \\ 
\qquad \times \left( \sum_{1\leq i\neq j\leq n}\left\vert \left(
y_{i},y_{j}\right) \right\vert ^{\delta }\right) ^{\frac{1}{\delta }}, \\ 
\hfill \ \ \text{where \ }\gamma >1,\ \ \frac{1}{\gamma }+\frac{1}{\delta }%
=1; \\ 
\\ 
\left[ \left( \sum_{i=1}^{n}\left\vert c_{i}\right\vert \right)
^{2}-\sum_{i=1}^{n}\left\vert c_{i}\right\vert ^{2}\right]
\max\limits_{1\leq i\neq j\leq n}\left\vert \left( y_{i},y_{j}\right)
\right\vert .%
\end{array}%
\right.  \label{3.1.12}
\end{equation}
\end{theorem}

\begin{proof}
We note that 
\begin{equation*}
\sum_{i=1}^{n}c_{i}\left( x,y_{i}\right) =\left( x,\sum_{i=1}^{n}\overline{%
c_{i}}y_{i}\right) .
\end{equation*}%
Using Schwarz's inequality in inner product spaces, we have: 
\begin{equation*}
\left\vert \sum_{i=1}^{n}c_{i}\left( x,y_{i}\right) \right\vert ^{2}\leq
\left\Vert x\right\Vert ^{2}\left\Vert \sum_{i=1}^{n}\overline{c_{i}}%
y_{i}\right\Vert ^{2}.
\end{equation*}%
Now using Lemma \ref{l2.1.12} with $\alpha _{i}=\overline{c_{i}},$ $%
z_{i}=y_{i}$ $\left( i=1,\dots ,n\right) ,$ we deduce the desired inequality
(\ref{3.1.12}).
\end{proof}

The following particular inequalities that may be obtained by the
Corollaries \ref{c2.3.12} and \ref{c2.4.12} and Remark \ref{r2.5.12} hold 
\cite{13NSSD}.

\begin{corollary}
\label{c3.2.12}With the assumptions in Theorem \ref{t3.1.12}, one has the
inequalities: 
\begin{multline}
\left\vert \sum_{i=1}^{n}c_{i}\left( x,y_{i}\right) \right\vert ^{2}
\label{3.2.12} \\
\leq \left\Vert x\right\Vert ^{2}\left\{ 
\begin{array}{l}
\sum_{i=1}^{n}\left\vert c_{i}\right\vert ^{2}\left\{ \max\limits_{1\leq
i\leq n}\left\Vert y_{i}\right\Vert ^{2}+\left( \sum_{1\leq i\neq j\leq
n}\left\vert \left( y_{i},y_{j}\right) \right\vert ^{2}\right) ^{\frac{1}{2}%
}\right\} ; \\ 
\\ 
\max\limits_{1\leq i\leq n}\left\vert c_{i}\right\vert ^{2}\left\{
\sum_{i=1}^{n}\left\Vert y_{i}\right\Vert ^{2}+\sum_{1\leq i\neq j\leq
n}\left\vert \left( y_{i},y_{j}\right) \right\vert \right\} \\ 
\\ 
\left( \sum_{i=1}^{n}\left\vert c_{i}\right\vert ^{2p}\right) ^{\frac{1}{p}%
}\left\{ \left( \sum_{i=1}^{n}\left\Vert y_{i}\right\Vert ^{2q}\right) ^{%
\frac{1}{q}}\right. \\ 
\qquad +\left. \left( n-1\right) ^{\frac{1}{p}}\left( \sum_{1\leq i\neq
j\leq n}\left\vert \left( y_{i},y_{j}\right) \right\vert ^{q}\right) ^{\frac{%
1}{q}}\right\} , \\ 
\hfill \ \ \ \text{where \ }p>1,\frac{1}{p}+\frac{1}{q}=1; \\ 
\sum_{i=1}^{n}\left\vert c_{i}\right\vert ^{2}\left\{ \max\limits_{1\leq
i\leq n}\left\Vert y_{i}\right\Vert ^{2}+\left( n-1\right)
\max\limits_{1\leq i\neq j\leq n}\left\vert \left( y_{i},y_{j}\right)
\right\vert \right\} .%
\end{array}%
\right.
\end{multline}
\end{corollary}

\begin{remark}
\label{r3.3.12}Note that the first inequality in (\ref{3.2.12}) is the
result obtained by Mitrinovi\'{c}-Pe\v{c}ari\'{c}-Fink in \cite{6b.12}. The
other 3 provide similar bounds in terms of the $p-$norms of the vector $%
\left( \left\vert c_{1}\right\vert ^{2},\dots ,\left\vert c_{n}\right\vert
^{2}\right) .$
\end{remark}

\subsection{Boas-Bellman Type Inequalities}

If one chooses $c_{i}=\overline{\left( x,y_{i}\right) }$ $\left( i=1,\dots
,n\right) $ in (\ref{3.1.12}), then it is possible to obtain 9 different
inequalities between the Fourier coefficients $\left( x,y_{i}\right) $ and
the norms and inner products of the vectors $y_{i}$ $\left( i=1,\dots
,n\right) .$ We restrict ourselves only to those inequalities that may be
obtained from (\ref{3.2.12}).

As Mitrinovi\'{c}, Pe\v{c}ari\'{c} and Fink noted in \cite[p. 392]{6b.12},
the first inequality in (\ref{3.2.12}) for the above selection of $c_{i}$
will produce the Boas-Bellman inequality (\ref{1.2.12}).

From the second inequality in (\ref{3.2.12}) for $c_{i}=\overline{\left(
x,y_{i}\right) }$ we get 
\begin{equation*}
\left( \sum_{i=1}^{n}\left\vert \left( x,y_{i}\right) \right\vert
^{2}\right) ^{2}\leq \left\Vert x\right\Vert ^{2}\max_{1\leq i\leq
n}\left\vert \left( x,y_{i}\right) \right\vert ^{2}\left\{
\sum_{i=1}^{n}\left\Vert y_{i}\right\Vert ^{2}+\sum_{1\leq i\neq j\leq
n}\left\vert \left( y_{i},y_{j}\right) \right\vert \right\} .
\end{equation*}%
Taking the square root in this inequality we obtain: 
\begin{multline}
\sum_{i=1}^{n}\left\vert \left( x,y_{i}\right) \right\vert ^{2}
\label{4.1.12} \\
\leq \left\Vert x\right\Vert \max_{1\leq i\leq n}\left\vert \left(
x,y_{i}\right) \right\vert \left\{ \sum_{i=1}^{n}\left\Vert y_{i}\right\Vert
^{2}+\sum_{1\leq i\neq j\leq n}\left\vert \left( y_{i},y_{j}\right)
\right\vert \right\} ^{\frac{1}{2}},
\end{multline}%
for any $x,y_{1},\dots ,y_{n}$ vectors in the inner product space $\left(
H;\left( \cdot ,\cdot \right) \right) .$

If we assume that $\left( e_{i}\right) _{1\leq i\leq n}$ is an orthonormal
family in $H,$ then by (\ref{4.1.12}) we have 
\begin{equation*}
\sum_{i=1}^{n}\left\vert \left( x,e_{i}\right) \right\vert ^{2}\leq \sqrt{n}%
\left\Vert x\right\Vert \max_{1\leq i\leq n}\left\vert \left( x,e_{i}\right)
\right\vert ,\ \ \ x\in H.
\end{equation*}%
From the third inequality in (\ref{3.2.12}) for $c_{i}=\overline{\left(
x,y_{i}\right) }$ we deduce 
\begin{multline*}
\left( \sum_{i=1}^{n}\left\vert \left( x,y_{i}\right) \right\vert
^{2}\right) ^{2}\leq \left\Vert x\right\Vert ^{2}\left(
\sum_{i=1}^{n}\left\vert \left( x,y_{i}\right) \right\vert ^{2p}\right) ^{%
\frac{1}{p}} \\
\times \left\{ \left( \sum\limits_{i=1}^{n}\left\Vert y_{i}\right\Vert
^{2q}\right) ^{\frac{1}{q}}+\left( n-1\right) ^{\frac{1}{p}}\left(
\sum\limits_{1\leq i\neq j\leq n}\left\vert \left( y_{i},y_{j}\right)
\right\vert ^{q}\right) ^{\frac{1}{q}}\right\} ,
\end{multline*}%
for $p>1,$ $\frac{1}{p}+\frac{1}{q}=1.$

Taking the square root in this inequality we get 
\begin{multline}
\sum_{i=1}^{n}\left\vert \left( x,y_{i}\right) \right\vert ^{2}\leq
\left\Vert x\right\Vert \left( \sum_{i=1}^{n}\left\vert \left(
x,y_{i}\right) \right\vert ^{2p}\right) ^{\frac{1}{2p}}  \label{4.3.12} \\
\times \left\{ \left( \sum\limits_{i=1}^{n}\left\Vert y_{i}\right\Vert
^{2q}\right) ^{\frac{1}{q}}+\left( n-1\right) ^{\frac{1}{p}}\left(
\sum\limits_{1\leq i\neq j\leq n}\left\vert \left( y_{i},y_{j}\right)
\right\vert ^{q}\right) ^{\frac{1}{q}}\right\} ^{\frac{1}{2}},
\end{multline}%
for any $x,y_{1},\dots ,y_{n}\in H,$ $p>1,$ $\frac{1}{p}+\frac{1}{q}=1.$

The above inequality (\ref{4.3.12}) becomes, for an orthornormal family $%
\left( e_{i}\right) _{1\leq i\leq n},$%
\begin{equation*}
\sum_{i=1}^{n}\left\vert \left( x,e_{i}\right) \right\vert ^{2}\leq n^{\frac{%
1}{q}}\left\Vert x\right\Vert \left( \sum_{i=1}^{n}\left\vert \left(
x,e_{i}\right) \right\vert ^{2p}\right) ^{\frac{1}{2p}},\ \ \ x\in H.
\end{equation*}%
Finally, the choice $c_{i}=\overline{\left( x,y_{i}\right) }$ $\left(
i=1,\dots ,n\right) $ will produce in the last inequality in (\ref{3.2.12}) 
\begin{equation*}
\left( \sum_{i=1}^{n}\left\vert \left( x,y_{i}\right) \right\vert
^{2}\right) ^{2}\leq \left\Vert x\right\Vert ^{2}\sum_{i=1}^{n}\left\vert
\left( x,y_{i}\right) \right\vert ^{2}\left\{ \max\limits_{1\leq i\leq
n}\left\Vert y_{i}\right\Vert ^{2}+\left( n-1\right) \max\limits_{1\leq
i\neq j\leq n}\left\vert \left( y_{i},y_{j}\right) \right\vert \right\}
\end{equation*}%
giving the following Boas-Bellman type inequality 
\begin{equation}
\sum_{i=1}^{n}\left\vert \left( x,y_{i}\right) \right\vert ^{2}\leq
\left\Vert x\right\Vert ^{2}\left\{ \max\limits_{1\leq i\leq n}\left\Vert
y_{i}\right\Vert ^{2}+\left( n-1\right) \max\limits_{1\leq i\neq j\leq
n}\left\vert \left( y_{i},y_{j}\right) \right\vert \right\} ,  \label{4.5.12}
\end{equation}%
for any $x,y_{1},\dots ,y_{n}\in H.$

It is obvious that (\ref{4.5.12}) will give for orthonormal families the
well known Bessel inequality.

\begin{remark}
In order the compare the Boas-Bellman result with our result (\ref{4.5.12}),
it is enough to compare the quantities 
\begin{equation*}
A:=\left( \sum\limits_{1\leq i\neq j\leq n}\left\vert \left(
y_{i},y_{j}\right) \right\vert ^{2}\right) ^{\frac{1}{2}}
\end{equation*}%
and 
\begin{equation*}
B:=\left( n-1\right) \max\limits_{1\leq i\neq j\leq n}\left\vert \left(
y_{i},y_{j}\right) \right\vert .
\end{equation*}%
Consider the inner product space $H=\mathbb{R}$ with $\left( x,y\right) =xy,$
and choose $n=3,$ $y_{1}=a>0$, $y_{2}=b>0,$ $y_{3}=c>0.$ Then 
\begin{equation*}
A=\sqrt{2}\left( a^{2}b^{2}+b^{2}c^{2}+c^{2}a^{2}\right) ^{\frac{1}{2}},\ \
\ \ \ \ \ B=2\max \left( ab,ac,bc\right) .
\end{equation*}%
Denote $ab=p,$ $bc=q,$ $ca=r.$ Then 
\begin{equation*}
A=\sqrt{2}\left( p^{2}+q^{2}+r^{2}\right) ^{\frac{1}{2}},\ \ \ \ \ \ \
B=2\max \left( p,q,r\right) .
\end{equation*}%
Firstly, if we assume that $p=q=r,$ then $A=\sqrt{6}p,$ $B=2p$ which shows
that $A>B.$

Now choose $r=1$ and $p,q=\frac{1}{2}.$ Then $A=\sqrt{3}$ and $B=2$ showing
that $B>A.$

Consequently, in general, the Boas-Bellman inequality and our inequality (%
\ref{4.5.12}) cannot be compared.
\end{remark}

\newpage

\section{Bombieri Type Inequalities}

\subsection{Introduction}

In 1971, E. Bombieri \cite{2ab.13} (see also \cite[p. 394]{6b.13}) gave the
following generalisation of Bessel's inequality.

\begin{theorem}
\label{t1.1.13}If $x,y_{1},\dots ,y_{n}$ are vectors in the inner product
space $\left( H;\left( \cdot ,\cdot \right) \right) ,$ then the following
inequality: 
\begin{equation}
\sum_{i=1}^{n}\left\vert \left( x,y_{i}\right) \right\vert ^{2}\leq
\left\Vert x\right\Vert ^{2}\max_{1\leq i\leq n}\left\{
\sum_{j=1}^{n}\left\vert \left( y_{i},y_{j}\right) \right\vert \right\} ,
\label{1.2.13}
\end{equation}%
holds.
\end{theorem}

It is obvious that if $\left( y_{i}\right) _{1\leq i\leq n}$ are
orthonormal, then from (\ref{1.2.13}) one can deduce Bessel's inequality.

Another generalisation of Bessel's inequality was obtained by A. Selberg
(see for example \cite[p. 394]{6b.13}):

\begin{theorem}
\label{t1.2.13}Let $x,y_{1},\dots ,y_{n}$ be vectors in $H$ with $y_{i}\neq
0 $ $\left( i=1,\dots ,n\right) .$ Then one has the inequality: 
\begin{equation}
\sum_{i=1}^{n}\frac{\left\vert \left( x,y_{i}\right) \right\vert ^{2}}{%
\sum_{j=1}^{n}\left\vert \left( y_{i},y_{j}\right) \right\vert }\leq
\left\Vert x\right\Vert ^{2}.  \label{1.3.13}
\end{equation}
\end{theorem}

In this case, also, if $\left( y_{i}\right) _{1\leq i\leq n}$ are
orthonormal, then from (\ref{1.3.13}) one may deduce Bessel's inequality.

Another type of inequality related to Bessel's result, was discovered in
1958 by H. Heilbronn \cite{5ab.13} (see also \cite[p. 395]{6b.13}).

\begin{theorem}
\label{t1.3.13}With the assumptions in Theorem \ref{t1.1.13}, one has 
\begin{equation}
\sum_{i=1}^{n}\left\vert \left( x,y_{i}\right) \right\vert \leq \left\Vert
x\right\Vert \left( \sum_{i,j=1}^{n}\left\vert \left( y_{i},y_{j}\right)
\right\vert \right) ^{\frac{1}{2}}.  \label{1.4.13}
\end{equation}
\end{theorem}

If in (\ref{1.4.13}) one chooses $y_{i}=e_{i}$ $\left( i=1,\dots ,n\right) ,$
where $\left( e_{i}\right) _{1\leq i\leq n}$ are orthonormal vectors in $H,$
then 
\begin{equation*}
\sum_{i=1}^{n}\left\vert \left( x,e_{i}\right) \right\vert \leq \sqrt{n}%
\left\Vert x\right\Vert ,\text{ \ for any \ }x\in H.
\end{equation*}

In 1992, J.E. Pe\v{c}ari\'{c} \cite{7b.13} (see also \cite[p. 394]{6b.13})
proved the following general inequality in inner product spaces.

\begin{theorem}
\label{t1.4.13}Let $x,y_{1},\dots ,y_{n}\in H$ and $c_{1},\dots ,c_{n}\in 
\mathbb{K}$. Then 
\begin{align}
\left\vert \sum_{i=1}^{n}c_{i}\left( x,y_{i}\right) \right\vert ^{2}& \leq
\left\Vert x\right\Vert ^{2}\sum_{i=1}^{n}\left\vert c_{i}\right\vert
^{2}\left( \sum_{j=1}^{n}\left\vert \left( y_{i},y_{j}\right) \right\vert
\right)  \label{1.6.13} \\
& \leq \left\Vert x\right\Vert ^{2}\sum_{i=1}^{n}\left\vert c_{i}\right\vert
^{2}\max_{1\leq i\leq n}\left\{ \sum_{j=1}^{n}\left\vert \left(
y_{i},y_{j}\right) \right\vert \right\} .  \notag
\end{align}
\end{theorem}

He showed that the Bombieri inequality (\ref{1.2.13}) may be obtained from (%
\ref{1.6.13}) for the choice $c_{i}=\overline{\left( x,y_{i}\right) }$
(using the second inequality), the Selberg inequality (\ref{1.3.13}) may be
obtained from the first part of (\ref{1.6.13}) for the choice 
\begin{equation*}
c_{i}=\frac{\overline{\left( x,y_{i}\right) }}{\sum_{j=1}^{n}\left\vert
\left( y_{i},y_{j}\right) \right\vert },\ \ \ i\in \left\{ 1,\dots
,n\right\} ;
\end{equation*}%
while the Heilbronn inequality (\ref{1.4.13}) may be obtained from the first
part of (\ref{1.6.13}) if one chooses $c_{i}=\frac{\overline{\left(
x,y_{i}\right) }}{\left\vert \left( x,y_{i}\right) \right\vert },$ for any $%
i\in \left\{ 1,\dots ,n\right\} .$

For other results connected with the above ones, see \cite{4b.13} and \cite%
{5b.13}.

\subsection{Some Norm Inequalities}

We start with the following lemma which is also interesting in itself \cite%
{14NSSD}.

\begin{lemma}
\label{l2.1.13}Let $z_{1},\dots ,z_{n}\in H$ and $\alpha _{1},\dots ,\alpha
_{n}\in \mathbb{K}.$ Then one has the inequality: 
\begin{equation}
\left\Vert \sum_{i=1}^{n}\alpha _{i}z_{i}\right\Vert ^{2}\leq \left\{ 
\begin{array}{c}
A \\ 
B \\ 
C%
\end{array}%
\right. ,  \label{2.1.13}
\end{equation}%
where 
\begin{equation*}
A:=\left\{ 
\begin{array}{l}
\max\limits_{1\leq k\leq n}\left\vert \alpha _{k}\right\vert
^{2}\sum\limits_{i,j=1}^{n}\left\vert \left( z_{i},z_{j}\right) \right\vert ;
\\ 
\\ 
\max\limits_{1\leq k\leq n}\left\vert \alpha _{k}\right\vert \left(
\sum\limits_{i=1}^{n}\left\vert \alpha _{i}\right\vert ^{r}\right) ^{\frac{1%
}{r}}\left( \sum\limits_{i=1}^{n}\left( \sum\limits_{j=1}^{n}\left\vert
\left( z_{i},z_{j}\right) \right\vert \right) ^{s}\right) ^{\frac{1}{s}},\ \
\ r>1,\ \frac{1}{r}+\frac{1}{s}=1; \\ 
\\ 
\max\limits_{1\leq k\leq n}\left\vert \alpha _{k}\right\vert
\sum\limits_{k=1}^{n}\left\vert \alpha _{k}\right\vert \max\limits_{1\leq
i\leq n}\left( \sum\limits_{j=1}^{n}\left\vert \left( z_{i},z_{j}\right)
\right\vert \right) ;%
\end{array}%
\right.
\end{equation*}%
\begin{equation*}
B:=\left\{ 
\begin{array}{l}
\left( \sum\limits_{k=1}^{n}\left\vert \alpha _{k}\right\vert ^{p}\right) ^{%
\frac{1}{p}}\max\limits_{1\leq i\leq n}\left\vert \alpha _{i}\right\vert
\left( \sum\limits_{i=1}^{n}\left( \sum\limits_{j=1}^{n}\left\vert \left(
z_{i},z_{j}\right) \right\vert \right) ^{q}\right) ^{\frac{1}{q}},\ \ \
p>1,\ \frac{1}{p}+\frac{1}{q}=1; \\ 
\\ 
\left( \sum\limits_{k=1}^{n}\left\vert \alpha _{k}\right\vert ^{p}\right) ^{%
\frac{1}{p}}\left( \sum\limits_{i=1}^{n}\left\vert \alpha _{i}\right\vert
^{t}\right) ^{\frac{1}{t}}\left[ \sum\limits_{i=1}^{n}\left(
\sum\limits_{j=1}^{n}\left\vert \left( z_{i},z_{j}\right) \right\vert
^{q}\right) ^{\frac{u}{q}}\right] ^{\frac{1}{u}},\ \ \ p>1,\ \frac{1}{p}+%
\frac{1}{q}=1; \\ 
\hfill \hfill t>1,\ \frac{1}{t}+\frac{1}{u}=1; \\ 
\left( \sum\limits_{k=1}^{n}\left\vert \alpha _{k}\right\vert ^{p}\right) ^{%
\frac{1}{p}}\sum\limits_{i=1}^{n}\left\vert \alpha _{i}\right\vert
\max\limits_{1\leq i\leq n}\left\{ \left( \sum\limits_{j=1}^{n}\left\vert
\left( z_{i},z_{j}\right) \right\vert ^{q}\right) ^{\frac{1}{q}}\right\} ,\
\ \ p>1,\ \frac{1}{p}+\frac{1}{q}=1;%
\end{array}%
\right.
\end{equation*}%
and%
\begin{equation*}
C:=\left\{ 
\begin{array}{l}
\sum\limits_{k=1}^{n}\left\vert \alpha _{k}\right\vert \ \max\limits_{1\leq
i\leq n}\left\vert \alpha _{i}\right\vert \ \sum\limits_{i=1}^{n}\left[
\max\limits_{1\leq j\leq n}\left\vert \left( z_{i},z_{j}\right) \right\vert %
\right] ; \\ 
\\ 
\sum\limits_{k=1}^{n}\left\vert \alpha _{k}\right\vert \ \left(
\sum\limits_{i=1}^{n}\left\vert \alpha _{i}\right\vert ^{m}\right) ^{\frac{1%
}{m}}\left( \sum\limits_{i=1}^{n}\left[ \max\limits_{1\leq j\leq
n}\left\vert \left( z_{i},z_{j}\right) \right\vert \right] ^{l}\right) ^{%
\frac{1}{l}},\ \ \ m>1,\ \frac{1}{m}+\frac{1}{l}=1; \\ 
\\ 
\left( \sum\limits_{k=1}^{n}\left\vert \alpha _{k}\right\vert \right)
^{2}\max\limits_{i,1\leq j\leq n}\left\vert \left( z_{i},z_{j}\right)
\right\vert .%
\end{array}%
\right.
\end{equation*}
\end{lemma}

\begin{proof}
We observe that 
\begin{align*}
\left\Vert \sum_{i=1}^{n}\alpha _{i}z_{i}\right\Vert ^{2}& =\left(
\sum_{i=1}^{n}\alpha _{i}z_{i},\sum_{j=1}^{n}\alpha _{j}z_{j}\right) \\
& =\sum_{i=1}^{n}\sum_{j=1}^{n}\alpha _{i}\overline{\alpha _{j}}\left(
z_{i},z_{j}\right) =\left\vert \sum_{i=1}^{n}\sum_{j=1}^{n}\alpha _{i}%
\overline{\alpha _{j}}\left( z_{i},z_{j}\right) \right\vert \\
& \leq \sum_{i=1}^{n}\sum_{j=1}^{n}\left\vert \alpha _{i}\right\vert
\left\vert \alpha _{j}\right\vert \left\vert \left( z_{i},z_{j}\right)
\right\vert =\sum_{i=1}^{n}\left\vert \alpha _{i}\right\vert \left(
\sum_{j=1}^{n}\left\vert \alpha _{j}\right\vert \left\vert \left(
z_{i},z_{j}\right) \right\vert \right) \\
& :=M.
\end{align*}%
Using H\"{o}lder's inequality, we may write that 
\begin{equation*}
\sum_{j=1}^{n}\left\vert \alpha _{j}\right\vert \left\vert \left(
z_{i},z_{j}\right) \right\vert \leq \left\{ 
\begin{array}{l}
\max\limits_{1\leq k\leq n}\left\vert \alpha _{k}\right\vert
\sum\limits_{j=1}^{n}\left\vert \left( z_{i},z_{j}\right) \right\vert \\ 
\\ 
\left( \sum\limits_{k=1}^{n}\left\vert \alpha _{k}\right\vert ^{p}\right) ^{%
\frac{1}{p}}\left( \sum\limits_{j=1}^{n}\left\vert \left( z_{i},z_{j}\right)
\right\vert ^{q}\right) ^{\frac{1}{q}},\ \ p>1,\ \frac{1}{p}+\frac{1}{q}=1;
\\ 
\\ 
\sum\limits_{k=1}^{n}\left\vert \alpha _{k}\right\vert \ \max\limits_{1\leq
j\leq n}\left\vert \left( z_{i},z_{j}\right) \right\vert%
\end{array}%
\right.
\end{equation*}%
for any $i\in \left\{ 1,\dots ,n\right\} ,$ giving 
\begin{equation*}
M\leq \left\{ 
\begin{array}{l}
\max\limits_{1\leq k\leq n}\left\vert \alpha _{k}\right\vert
\sum\limits_{i=1}^{n}\left\vert \alpha _{i}\right\vert
\sum\limits_{j=1}^{n}\left\vert \left( z_{i},z_{j}\right) \right\vert
=:M_{1}; \\ 
\\ 
\left( \sum\limits_{k=1}^{n}\left\vert \alpha _{k}\right\vert ^{p}\right) ^{%
\frac{1}{p}}\sum\limits_{i=1}^{n}\left\vert \alpha _{i}\right\vert \left(
\sum\limits_{j=1}^{n}\left\vert \left( z_{i},z_{j}\right) \right\vert
^{q}\right) ^{\frac{1}{q}}:=M_{p},\  \\ 
\hfill \hfill p>1,\ \frac{1}{p}+\frac{1}{q}=1; \\ 
\sum\limits_{k=1}^{n}\left\vert \alpha _{k}\right\vert
\sum\limits_{i=1}^{n}\left\vert \alpha _{i}\right\vert \ \max\limits_{1\leq
j\leq n}\left\vert \left( z_{i},z_{j}\right) \right\vert =:M_{\infty }.%
\end{array}%
\right.
\end{equation*}%
By H\"{o}lder's inequality we also have: 
\begin{multline*}
\sum\limits_{i=1}^{n}\left\vert \alpha _{i}\right\vert \left(
\sum\limits_{j=1}^{n}\left\vert \left( z_{i},z_{j}\right) \right\vert \right)
\\
\leq \left\{ 
\begin{array}{l}
\max\limits_{1\leq i\leq n}\left\vert \alpha _{i}\right\vert
\sum\limits_{i,j=1}^{n}\left\vert \left( z_{i},z_{j}\right) \right\vert ; \\ 
\\ 
\left( \sum\limits_{i=1}^{n}\left\vert \alpha _{i}\right\vert ^{r}\right) ^{%
\frac{1}{r}}\left( \sum\limits_{i=1}^{n}\left(
\sum\limits_{j=1}^{n}\left\vert \left( z_{i},z_{j}\right) \right\vert
\right) ^{s}\right) ^{\frac{1}{s}},\ \ \ r>1,\ \frac{1}{r}+\frac{1}{s}=1; \\ 
\\ 
\sum\limits_{i=1}^{n}\left\vert \alpha _{i}\right\vert \ \max\limits_{1\leq
i\leq n}\left( \sum\limits_{j=1}^{n}\left\vert \left( z_{i},z_{j}\right)
\right\vert \right) ;%
\end{array}%
\right.
\end{multline*}%
and thus 
\begin{equation*}
M_{1}\leq \left\{ 
\begin{array}{l}
\max\limits_{1\leq k\leq n}\left\vert \alpha _{k}\right\vert
^{2}\sum\limits_{i,j=1}^{n}\left\vert \left( z_{i},z_{j}\right) \right\vert ;
\\ 
\\ 
\max\limits_{1\leq k\leq n}\left\vert \alpha _{k}\right\vert \left(
\sum\limits_{i=1}^{n}\left\vert \alpha _{i}\right\vert ^{r}\right) ^{\frac{1%
}{r}}\left( \sum\limits_{i=1}^{n}\left( \sum\limits_{j=1}^{n}\left\vert
\left( z_{i},z_{j}\right) \right\vert \right) ^{s}\right) ^{\frac{1}{s}},\ \
\ r>1,\ \frac{1}{r}+\frac{1}{s}=1; \\ 
\\ 
\max\limits_{1\leq k\leq n}\left\vert \alpha _{k}\right\vert
\sum\limits_{i=1}^{n}\left\vert \alpha _{i}\right\vert \ \max\limits_{1\leq
i\leq n}\left( \sum\limits_{j=1}^{n}\left\vert \left( z_{i},z_{j}\right)
\right\vert \right) ;%
\end{array}%
\right.
\end{equation*}%
and the first 3 inequalities in (\ref{2.1.13}) are obtained.

By H\"{o}lder's inequality we also have: 
\begin{multline*}
M_{p}\leq \left( \sum\limits_{k=1}^{n}\left\vert \alpha _{k}\right\vert
^{p}\right) ^{\frac{1}{p}} \\
\times \left\{ 
\begin{array}{l}
\max\limits_{1\leq i\leq n}\left\vert \alpha _{i}\right\vert
\sum\limits_{i=1}^{n}\left( \sum\limits_{j=1}^{n}\left\vert \left(
z_{i},z_{j}\right) \right\vert ^{q}\right) ^{\frac{1}{q}}; \\ 
\\ 
\left( \sum\limits_{i=1}^{n}\left\vert \alpha _{i}\right\vert ^{t}\right) ^{%
\frac{1}{t}}\left( \sum\limits_{i=1}^{n}\left(
\sum\limits_{j=1}^{n}\left\vert \left( z_{i},z_{j}\right) \right\vert
^{q}\right) ^{\frac{u}{q}}\right) ^{\frac{1}{u}},\ \ \ \hfill t>1,\ \frac{1}{%
t}+\frac{1}{u}=1; \\ 
\\ 
\sum\limits_{i=1}^{n}\left\vert \alpha _{i}\right\vert \ \max\limits_{1\leq
i\leq n}\left\{ \left( \sum\limits_{j=1}^{n}\left\vert \left(
z_{i},z_{j}\right) \right\vert ^{q}\right) ^{\frac{1}{q}}\right\} ;%
\end{array}%
\right.
\end{multline*}%
and the next 3 inequalities in (\ref{2.1.13}) are proved.

Finally, by the same H\"{o}lder inequality we may state that: 
\begin{equation*}
M_{\infty }\leq \sum\limits_{k=1}^{n}\left\vert \alpha _{k}\right\vert
\times \left\{ 
\begin{array}{l}
\max\limits_{1\leq i\leq n}\left\vert \alpha _{i}\right\vert
\sum\limits_{i=1}^{n}\left( \max\limits_{1\leq j\leq n}\left\vert \left(
z_{i},z_{j}\right) \right\vert \right) ; \\ 
\\ 
\left( \sum\limits_{i=1}^{n}\left\vert \alpha _{i}\right\vert ^{m}\right) ^{%
\frac{1}{m}}\left( \sum\limits_{i=1}^{n}\left( \max\limits_{1\leq j\leq
n}\left\vert \left( z_{i},z_{j}\right) \right\vert \right) ^{l}\right) ^{%
\frac{1}{l}},\ \ \ \hfill m>1,\ \frac{1}{m}+\frac{1}{l}=1; \\ 
\\ 
\sum\limits_{i=1}^{n}\left\vert \alpha _{i}\right\vert \ \max\limits_{1\leq
i,j\leq n}\left\vert \left( z_{i},z_{j}\right) \right\vert ;%
\end{array}%
\right.
\end{equation*}%
and the last 3 inequalities in (\ref{2.1.13}) are proved.
\end{proof}

If we would like to have some bounds for $\left\| \sum_{i=1}^{n}\alpha
_{i}z_{i}\right\| ^{2}$ in terms of $\sum_{i=1}^{n}\left| \alpha _{i}\right|
^{2},$ then the following corollaries may be used.

\begin{corollary}
\label{c2.1.a.13}Let $z_{1},\dots ,z_{n}$ and $\alpha _{1},\dots ,\alpha
_{n} $ be as in Lemma \ref{l2.1.13}. If $1<p\leq 2$, $1<t\leq 2,$ then one
has the inequality 
\begin{equation}
\left\Vert \sum_{i=1}^{n}\alpha _{i}z_{i}\right\Vert ^{2}\leq n^{\frac{1}{p}+%
\frac{1}{t}-1}\sum\limits_{k=1}^{n}\left\vert \alpha _{k}\right\vert ^{2}%
\left[ \sum\limits_{i=1}^{n}\left( \sum\limits_{j=1}^{n}\left\vert \left(
z_{i},z_{j}\right) \right\vert ^{q}\right) ^{\frac{u}{q}}\right] ^{\frac{1}{u%
}},  \label{2.5.a.13}
\end{equation}%
where $\frac{1}{p}+\frac{1}{q}=1,$ $\frac{1}{t}+\frac{1}{u}=1.$
\end{corollary}

\begin{proof}
Observe, by the monotonicity of power means, we may write that 
\begin{align*}
\left( \frac{\sum_{k=1}^{n}\left\vert \alpha _{k}\right\vert ^{p}}{n}\right)
^{\frac{1}{p}}& \leq \left( \frac{\sum_{k=1}^{n}\left\vert \alpha
_{k}\right\vert ^{2}}{n}\right) ^{\frac{1}{2}};\ \ 1<p\leq 2, \\
\left( \frac{\sum_{k=1}^{n}\left\vert \alpha _{k}\right\vert ^{t}}{n}\right)
^{\frac{1}{t}}& \leq \left( \frac{\sum_{k=1}^{n}\left\vert \alpha
_{k}\right\vert ^{2}}{n}\right) ^{\frac{1}{2}};\ \ 1<t\leq 2,
\end{align*}%
from where we get 
\begin{align*}
\left( \sum_{k=1}^{n}\left\vert \alpha _{k}\right\vert ^{p}\right) ^{\frac{1%
}{p}}& \leq n^{\frac{1}{p}-\frac{1}{2}}\left( \sum_{k=1}^{n}\left\vert
\alpha _{k}\right\vert ^{2}\right) ^{\frac{1}{2}}, \\
\left( \sum_{k=1}^{n}\left\vert \alpha _{k}\right\vert ^{t}\right) ^{\frac{1%
}{t}}& \leq n^{\frac{1}{t}-\frac{1}{2}}\left( \sum_{k=1}^{n}\left\vert
\alpha _{k}\right\vert ^{2}\right) ^{\frac{1}{2}}.
\end{align*}%
Using the fifth inequality in (\ref{2.1.13}), we then deduce (\ref{2.5.a.13}%
).
\end{proof}

\begin{remark}
\label{r2.1.b.13}An interesting particular case is the one for $p=q=t=u=2,$
giving 
\begin{equation*}
\left\Vert \sum_{i=1}^{n}\alpha _{i}z_{i}\right\Vert ^{2}\leq
\sum_{k=1}^{n}\left\vert \alpha _{k}\right\vert ^{2}\left(
\sum\limits_{i,j=1}^{n}\left\vert \left( z_{i},z_{j}\right) \right\vert
^{2}\right) ^{\frac{1}{2}}.
\end{equation*}
\end{remark}

\begin{corollary}
\label{c2.1.c.13}With the assumptions of Lemma \ref{l2.1.13} and if $1<p\leq
2,$ then 
\begin{equation}
\left\Vert \sum_{i=1}^{n}\alpha _{i}z_{i}\right\Vert ^{2}\leq n^{\frac{1}{p}%
}\sum_{k=1}^{n}\left\vert \alpha _{k}\right\vert ^{2}\max\limits_{1\leq
i\leq n}\left[ \left( \sum\limits_{j=1}^{n}\left\vert \left(
z_{i},z_{j}\right) \right\vert ^{q}\right) ^{\frac{1}{q}}\right] ,
\label{2.5.c.13}
\end{equation}%
where $\frac{1}{p}+\frac{1}{q}=1.$
\end{corollary}

\begin{proof}
Since 
\begin{equation*}
\left( \sum_{k=1}^{n}\left\vert \alpha _{k}\right\vert ^{p}\right) ^{\frac{1%
}{p}}\leq n^{\frac{1}{p}-\frac{1}{2}}\left( \sum_{k=1}^{n}\left\vert \alpha
_{k}\right\vert ^{2}\right) ^{\frac{1}{2}},
\end{equation*}%
and 
\begin{equation*}
\sum_{k=1}^{n}\left\vert \alpha _{k}\right\vert \leq n^{\frac{1}{2}}\left(
\sum_{k=1}^{n}\left\vert \alpha _{k}\right\vert ^{2}\right) ^{\frac{1}{2}},
\end{equation*}%
then by the sixth inequality in (\ref{2.1.13}) we deduce (\ref{2.5.c.13}).
\end{proof}

In a similar fashion, one may prove the following two corollaries.

\begin{corollary}
\label{c2.1.1.13}With the assumptions of Lemma \ref{l2.1.13} and if $1<m\leq
2,$ then 
\begin{equation*}
\left\Vert \sum_{i=1}^{n}\alpha _{i}z_{i}\right\Vert ^{2}\leq n^{\frac{1}{m}%
}\sum_{k=1}^{n}\left\vert \alpha _{k}\right\vert ^{2}\left(
\sum\limits_{i=1}^{n}\left[ \max\limits_{1\leq j\leq n}\left\vert \left(
z_{i},z_{j}\right) \right\vert \right] ^{l}\right) ^{\frac{1}{l}},
\end{equation*}%
where $\frac{1}{m}+\frac{1}{l}=1.$
\end{corollary}

\begin{corollary}
\label{c2.1.e.13}With the assumptions of Lemma \ref{l2.1.13}, we have: 
\begin{equation*}
\left\Vert \sum_{i=1}^{n}\alpha _{i}z_{i}\right\Vert ^{2}\leq
n\sum_{k=1}^{n}\left\vert \alpha _{k}\right\vert ^{2}\ \max\limits_{1\leq
i,j\leq n}\left\vert \left( z_{i},z_{j}\right) \right\vert .
\end{equation*}
\end{corollary}

The following lemma may be of interest as well \cite{14NSSD}.

\begin{lemma}
\label{l2.2.13}With the assumptions of Lemma \ref{l2.1.13}, one has the
inequalities 
\begin{align}
\left\Vert \sum_{i=1}^{n}\alpha _{i}z_{i}\right\Vert ^{2}& \leq
\sum_{i=1}^{n}\left\vert \alpha _{i}\right\vert
^{2}\sum\limits_{j=1}^{n}\left\vert \left( z_{i},z_{j}\right) \right\vert
\label{2.6.13} \\
& \leq \left\{ 
\begin{array}{l}
\sum_{i=1}^{n}\left\vert \alpha _{i}\right\vert ^{2}\ \max\limits_{1\leq
i\leq n}\left[ \sum_{j=1}^{n}\left\vert \left( z_{i},z_{j}\right)
\right\vert \right] ; \\ 
\\ 
\left( \sum_{i=1}^{n}\left\vert \alpha _{i}\right\vert ^{2p}\right) ^{\frac{1%
}{p}}\left( \left( \sum_{j=1}^{n}\left\vert \left( z_{i},z_{j}\right)
\right\vert \right) ^{q}\right) ^{\frac{1}{q}},\  \\ 
\hfill p>1,\ \frac{1}{p}+\frac{1}{q}=1; \\ 
\\ 
\max\limits_{1\leq i\leq n}\left\vert \alpha _{i}\right\vert ^{2}\
\sum_{i,j=1}^{n}\left\vert \left( z_{i},z_{j}\right) \right\vert .%
\end{array}%
\right.  \notag
\end{align}
\end{lemma}

\begin{proof}
As in Lemma \ref{l2.1.13}, we know that 
\begin{equation*}
\left\Vert \sum_{i=1}^{n}\alpha _{i}z_{i}\right\Vert ^{2}\leq
\sum_{i=1}^{n}\sum_{j=1}^{n}\left\vert \alpha _{i}\right\vert \left\vert
\alpha _{j}\right\vert \left\vert \left( z_{i},z_{j}\right) \right\vert .
\end{equation*}%
Using the simple observation that (see also \cite[p. 394]{6b.13}) 
\begin{equation*}
\left\vert \alpha _{i}\right\vert \left\vert \alpha _{j}\right\vert \leq 
\frac{1}{2}\left( \left\vert \alpha _{i}\right\vert ^{2}+\left\vert \alpha
_{j}\right\vert ^{2}\right) ,\ \ \ i,j\in \left\{ 1,\dots ,n\right\} ,
\end{equation*}%
we have 
\begin{align*}
\sum_{i=1}^{n}\sum_{j=1}^{n}\left\vert \alpha _{i}\right\vert \left\vert
\alpha _{j}\right\vert \left\vert \left( z_{i},z_{j}\right) \right\vert &
\leq \frac{1}{2}\sum\limits_{i,j=1}^{n}\left( \left\vert \alpha
_{i}\right\vert ^{2}+\left\vert \alpha _{j}\right\vert ^{2}\right)
\left\vert \left( z_{i},z_{j}\right) \right\vert \\
& =\frac{1}{2}\left[ \sum\limits_{i,j=1}^{n}\left\vert \alpha
_{i}\right\vert ^{2}\left\vert \left( z_{i},z_{j}\right) \right\vert
+\sum\limits_{i,j=1}^{n}\left\vert \alpha _{j}\right\vert ^{2}\left\vert
\left( z_{i},z_{j}\right) \right\vert \right] \\
& =\sum\limits_{i,j=1}^{n}\left\vert \alpha _{i}\right\vert ^{2}\left\vert
\left( z_{i},z_{j}\right) \right\vert ,
\end{align*}%
which proves the first inequality in (\ref{2.6.13}).

The second part follows by H\"{o}lder's inequality and we omit the details.
\end{proof}

\begin{remark}
The first part in (\ref{2.6.13}) is the inequality obtained by Pe\v{c}ari%
\'{c} in \cite{7b.13}.
\end{remark}

\subsection{Pe\v{c}ari\'{c} Type Inequalities}

We are now able to present the following result obtainend in \cite{14NSSD},
which complements the inequality (\ref{1.6.13}) due to J.E. Pe\v{c}ari\'{c} 
\cite{7b.13} (see also \cite[p. 394]{6b.13}).

\begin{theorem}
\label{t3.1.13}Let $x,y_{1},\dots ,y_{n}$ be vectors of an inner product
space $\left( H;\left( \cdot ,\cdot \right) \right) $ and $c_{1},\dots
,c_{n}\in \mathbb{K}$. Then one has the inequalities: 
\begin{equation}
\left\vert \sum\limits_{i=1}^{n}c_{i}\left( x,y_{i}\right) \right\vert
^{2}\leq \left\Vert x\right\Vert ^{2}\times \left\{ 
\begin{array}{c}
D \\ 
E \\ 
F%
\end{array}%
,\right.  \label{3.1.13}
\end{equation}%
where%
\begin{equation*}
D:=\left\{ 
\begin{array}{l}
\max\limits_{1\leq k\leq n}\left\vert c_{k}\right\vert
^{2}\sum\limits_{i,j=1}^{n}\left\vert \left( y_{i},y_{j}\right) \right\vert ;
\\ 
\\ 
\max\limits_{1\leq k\leq n}\left\vert c_{k}\right\vert \left(
\sum\limits_{i=1}^{n}\left\vert c_{i}\right\vert ^{r}\right) ^{\frac{1}{r}}%
\left[ \sum\limits_{i=1}^{n}\left( \sum\limits_{j=1}^{n}\left\vert \left(
y_{i},y_{j}\right) \right\vert \right) ^{s}\right] ^{\frac{1}{s}},\ \ \
r>1,\ \frac{1}{r}+\frac{1}{s}=1; \\ 
\\ 
\max\limits_{1\leq k\leq n}\left\vert c_{k}\right\vert
\sum\limits_{k=1}^{n}\left\vert c_{k}\right\vert \max\limits_{1\leq i\leq
n}\left( \sum\limits_{j=1}^{n}\left\vert \left( y_{i},y_{j}\right)
\right\vert \right) ;%
\end{array}%
\right.
\end{equation*}%
\begin{equation*}
E:=\left\{ 
\begin{array}{ll}
\left( \sum\limits_{k=1}^{n}\left\vert c_{k}\right\vert ^{p}\right) ^{\frac{1%
}{p}}\max\limits_{1\leq i\leq n}\left\vert c_{i}\right\vert \left(
\sum\limits_{i=1}^{n}\left( \sum\limits_{j=1}^{n}\left\vert \left(
y_{i},y_{j}\right) \right\vert \right) ^{q}\right) ^{\frac{1}{q}}, & p>1,\ 
\frac{1}{p}+\frac{1}{q}=1; \\ 
&  \\ 
\left( \sum\limits_{k=1}^{n}\left\vert c_{k}\right\vert ^{p}\right) ^{\frac{1%
}{p}}\left( \sum\limits_{i=1}^{n}\left\vert c_{i}\right\vert ^{t}\right) ^{%
\frac{1}{t}}\left[ \sum\limits_{i=1}^{n}\left(
\sum\limits_{j=1}^{n}\left\vert \left( y_{i},y_{j}\right) \right\vert
^{q}\right) ^{\frac{u}{q}}\right] ^{\frac{1}{u}}, & p>1,\ \frac{1}{p}+\frac{1%
}{q}=1; \\ 
& t>1,\ \frac{1}{t}+\frac{1}{u}=1; \\ 
&  \\ 
\left( \sum\limits_{k=1}^{n}\left\vert c_{k}\right\vert ^{p}\right) ^{\frac{1%
}{p}}\sum\limits_{i=1}^{n}\left\vert c_{i}\right\vert \max\limits_{1\leq
i\leq n}\left\{ \left( \sum\limits_{j=1}^{n}\left\vert \left(
y_{i},y_{j}\right) \right\vert ^{q}\right) ^{\frac{1}{q}}\right\} , & p>1,\ 
\frac{1}{p}+\frac{1}{q}=1;%
\end{array}%
\right.
\end{equation*}%
and%
\begin{equation*}
F:=\left\{ 
\begin{array}{l}
\sum_{k=1}^{n}\left\vert c_{k}\right\vert \ \max\limits_{1\leq i\leq
n}\left\vert c_{i}\right\vert \ \sum_{i=1}^{n}\left[ \max\limits_{1\leq
j\leq n}\left\vert \left( y_{i},y_{j}\right) \right\vert \right] ; \\ 
\\ 
\sum_{k=1}^{n}\left\vert c_{k}\right\vert \ \left( \sum_{i=1}^{n}\left\vert
c_{i}\right\vert ^{m}\right) ^{\frac{1}{m}}\left( \sum_{i=1}^{n}\left[
\max\limits_{1\leq j\leq n}\left\vert \left( y_{i},y_{j}\right) \right\vert %
\right] ^{l}\right) ^{\frac{1}{l}}, \\ 
\hfill \ \ \ m>1,\ \frac{1}{m}+\frac{1}{l}=1; \\ 
\\ 
\left( \sum_{k=1}^{n}\left\vert c_{k}\right\vert \right)
^{2}\max\limits_{i,1\leq j\leq n}\left\vert \left( y_{i},y_{j}\right)
\right\vert .%
\end{array}%
\right.
\end{equation*}
\end{theorem}

\begin{proof}
We note that 
\begin{equation*}
\sum\limits_{i=1}^{n}c_{i}\left( x,y_{i}\right) =\left(
x,\sum\limits_{i=1}^{n}\overline{c_{i}}y_{i}\right) .
\end{equation*}%
Using Schwarz's inequality in inner product spaces, we have 
\begin{equation*}
\left\vert \sum\limits_{i=1}^{n}c_{i}\left( x,y_{i}\right) \right\vert
^{2}\leq \left\Vert x\right\Vert ^{2}\left\Vert \sum\limits_{i=1}^{n}%
\overline{c_{i}}y_{i}\right\Vert ^{2}.
\end{equation*}%
Finally, using Lemma \ref{l2.1.13} with $\alpha _{i}=\overline{c_{i}},$ $%
z_{i}=y_{i}$ $\left( i=1,\dots ,n\right) ,$ we deduce the desired inequality
(\ref{3.1.13}). We omit the details.
\end{proof}

The following corollaries may be useful if one needs bounds in terms of $%
\sum_{i=1}^{n}\left| c_{i}\right| ^{2}.$

\begin{corollary}
\label{c3.2.13}With the assumptions in Theorem \ref{t3.1.13} and if $1<p\leq
2,$ $1<t\leq 2,$ $\frac{1}{p}+\frac{1}{q}=1,$ $\frac{1}{t}+\frac{1}{u}=1,$
one has the inequality: 
\begin{multline}
\left\vert \sum\limits_{i=1}^{n}c_{i}\left( x,y_{i}\right) \right\vert ^{2}
\label{3.3.13} \\
\leq \left\Vert x\right\Vert ^{2}n^{\frac{1}{p}+\frac{1}{t}%
-1}\sum\limits_{i=1}^{n}\left\vert c_{i}\right\vert ^{2}\left[
\sum\limits_{i=1}^{n}\left( \sum\limits_{j=1}^{n}\left\vert \left(
y_{i},y_{j}\right) \right\vert ^{q}\right) ^{\frac{u}{q}}\right] ^{\frac{1}{u%
}},
\end{multline}%
and, in particular, for $p=q=t=u=2,$%
\begin{equation*}
\left\vert \sum\limits_{i=1}^{n}c_{i}\left( x,y_{i}\right) \right\vert
^{2}\leq \left\Vert x\right\Vert ^{2}\sum\limits_{i=1}^{n}\left\vert
c_{i}\right\vert ^{2}\left( \sum\limits_{i,j=1}^{n}\left\vert \left(
y_{i},y_{j}\right) \right\vert ^{2}\right) ^{\frac{1}{2}}.
\end{equation*}
\end{corollary}

The proof is similar to the one in Corollary \ref{c2.1.a.13}.

\begin{corollary}
\label{c3.3.13}With the assumptions in Theorem \ref{t3.1.13} and if $1<p\leq
2,$ then 
\begin{equation*}
\left\vert \sum\limits_{i=1}^{n}c_{i}\left( x,y_{i}\right) \right\vert
^{2}\leq \left\Vert x\right\Vert ^{2}n^{\frac{1}{p}}\sum\limits_{k=1}^{n}%
\left\vert c_{k}\right\vert ^{2}\ \max\limits_{1\leq i\leq n}\left[
\sum\limits_{j=1}^{n}\left\vert \left( y_{i},y_{j}\right) \right\vert ^{q}%
\right] ^{\frac{1}{q}},
\end{equation*}%
where $\frac{1}{p}+\frac{1}{q}=1.$
\end{corollary}

The proof is similar to the one in Corollary \ref{c2.1.c.13}.

The following two inequalities also hold.

\begin{corollary}
\label{c3.4.13}With the above assumptions for $x,y_{i},c_{i}$ and if $%
1<m\leq 2,$ then 
\begin{equation}
\left\vert \sum\limits_{i=1}^{n}c_{i}\left( x,y_{i}\right) \right\vert
^{2}\leq \left\Vert x\right\Vert ^{2}n^{\frac{1}{m}}\sum\limits_{k=1}^{n}%
\left\vert c_{k}\right\vert ^{2}\left( \sum\limits_{i=1}^{n}\left[
\max\limits_{1\leq j\leq n}\left\vert \left( y_{i},y_{j}\right) \right\vert %
\right] ^{l}\right) ^{\frac{1}{l}},  \label{3.6.13}
\end{equation}%
where $\frac{1}{m}+\frac{1}{l}=1.$
\end{corollary}

\begin{corollary}
\label{c3.5.13}With the above assumptions for $x,y_{i},c_{i},$ one has 
\begin{equation}
\left\vert \sum\limits_{i=1}^{n}c_{i}\left( x,y_{i}\right) \right\vert
^{2}\leq \left\Vert x\right\Vert ^{2}n\sum\limits_{k=1}^{n}\left\vert
c_{k}\right\vert ^{2}\max\limits_{1\leq j\leq n}\left\vert \left(
y_{i},y_{j}\right) \right\vert .  \label{3.7.13}
\end{equation}
\end{corollary}

Using Lemma \ref{l2.2.13}, we may state the following result as well.

\begin{remark}
\label{t3.6.13}With the assumptions of Theorem \ref{t3.1.13}, one has the
inequalities: 
\begin{align*}
\left\vert \sum\limits_{i=1}^{n}c_{i}\left( x,y_{i}\right) \right\vert ^{2}&
\leq \left\Vert x\right\Vert ^{2}\sum\limits_{i=1}^{n}\left\vert
c_{i}\right\vert ^{2}\sum\limits_{j=1}^{n}\left\vert \left(
y_{i},y_{j}\right) \right\vert \\
& \leq \left\Vert x\right\Vert ^{2}\left\{ 
\begin{array}{l}
\sum_{i=1}^{n}\left\vert c_{i}\right\vert ^{2}\ \max\limits_{1\leq i\leq n} 
\left[ \sum_{j=1}^{n}\left\vert \left( y_{i},y_{j}\right) \right\vert \right]
; \\ 
\\ 
\left( \sum_{i=1}^{n}\left\vert c_{i}\right\vert ^{2p}\right) ^{\frac{1}{p}%
}\left( \sum_{i=1}^{n}\left( \sum_{j=1}^{n}\left\vert \left(
y_{i},y_{j}\right) \right\vert \right) ^{q}\right) ^{\frac{1}{q}}, \\ 
\hfill \ p>1,\ \frac{1}{p}+\frac{1}{q}=1; \\ 
\\ 
\max\limits_{1\leq i\leq n}\left\vert c_{i}\right\vert ^{2}\
\sum_{i,j=1}^{n}\left\vert \left( y_{i},y_{j}\right) \right\vert ;%
\end{array}%
\right.
\end{align*}%
that provide some alternatives to Pe\v{c}ari\'{c}'s result (\ref{1.6.13}).
\end{remark}

\subsection{Inequalities of Bombieri Type}

In this section we point out some inequalities of Bombieri type that may be
obtained from (\ref{3.1.13}) on choosing $c_{i}=\overline{\left(
x,y_{i}\right) }$ $\left( i=1,\dots ,n\right) .$

If the above choice was made in the first inequality in (\ref{3.1.13}), then
one can obtain: 
\begin{equation*}
\left( \sum\limits_{i=1}^{n}\left\vert \left( x,y_{i}\right) \right\vert
^{2}\right) ^{2}\leq \left\Vert x\right\Vert ^{2}\max\limits_{1\leq i\leq
n}\left\vert \left( x,y_{i}\right) \right\vert
^{2}\sum\limits_{i,j=1}^{n}\left\vert \left( y_{i},y_{j}\right) \right\vert
\end{equation*}%
giving, by taking the square root, 
\begin{equation}
\sum\limits_{i=1}^{n}\left\vert \left( x,y_{i}\right) \right\vert ^{2}\leq
\left\Vert x\right\Vert \max\limits_{1\leq i\leq n}\left\vert \left(
x,y_{i}\right) \right\vert \left( \sum\limits_{i,j=1}^{n}\left\vert \left(
y_{i},y_{j}\right) \right\vert \right) ^{\frac{1}{2}},\ \ x\in H.
\label{4.1.13}
\end{equation}%
If the same choice for $c_{i}$ is made in the second inequality in (\ref%
{3.1.13}), then one can get 
\begin{multline*}
\left( \sum\limits_{i=1}^{n}\left\vert \left( x,y_{i}\right) \right\vert
^{2}\right) ^{2}\leq \left\Vert x\right\Vert ^{2}\max\limits_{1\leq i\leq
n}\left\vert \left( x,y_{i}\right) \right\vert \left(
\sum\limits_{i=1}^{n}\left\vert \left( x,y_{i}\right) \right\vert
^{r}\right) ^{\frac{1}{r}} \\
\times \left[ \sum\limits_{i=1}^{n}\left( \sum\limits_{j=1}^{n}\left\vert
\left( y_{i},y_{j}\right) \right\vert \right) ^{s}\right] ^{\frac{1}{s}},
\end{multline*}%
implying 
\begin{multline}
\sum\limits_{i=1}^{n}\left\vert \left( x,y_{i}\right) \right\vert ^{2}\leq
\left\Vert x\right\Vert \max\limits_{1\leq i\leq n}\left\vert \left(
x,y_{i}\right) \right\vert ^{\frac{1}{2}}\left(
\sum\limits_{i=1}^{n}\left\vert \left( x,y_{i}\right) \right\vert
^{r}\right) ^{\frac{1}{2r}}  \label{4.2.13} \\
\times \left[ \sum\limits_{i=1}^{n}\left( \sum\limits_{j=1}^{n}\left\vert
\left( y_{i},y_{j}\right) \right\vert \right) ^{s}\right] ^{\frac{1}{2s}},
\end{multline}%
where $\frac{1}{r}+\frac{1}{s}=1,$ $s>1.$

The other inequalities in (\ref{3.1.13}) will produce the following results,
respectively 
\begin{multline}
\sum\limits_{i=1}^{n}\left\vert \left( x,y_{i}\right) \right\vert ^{2}\leq
\left\Vert x\right\Vert \max\limits_{1\leq i\leq n}\left\vert \left(
x,y_{i}\right) \right\vert ^{\frac{1}{2}}\left(
\sum\limits_{i=1}^{n}\left\vert \left( x,y_{i}\right) \right\vert \right) ^{%
\frac{1}{2}}  \label{4.3.13} \\
\times \left[ \max\limits_{1\leq i\leq n}\left(
\sum\limits_{j=1}^{n}\left\vert \left( y_{i},y_{j}\right) \right\vert
\right) \right] ;
\end{multline}%
\begin{multline}
\sum\limits_{i=1}^{n}\left\vert \left( x,y_{i}\right) \right\vert ^{2}\leq
\left\Vert x\right\Vert \max\limits_{1\leq i\leq n}\left\vert \left(
x,y_{i}\right) \right\vert ^{\frac{1}{2}}\left(
\sum\limits_{i=1}^{n}\left\vert \left( x,y_{i}\right) \right\vert
^{p}\right) ^{\frac{1}{2p}}  \label{4.4.13} \\
\times \left[ \sum\limits_{i=1}^{n}\left( \sum\limits_{j=1}^{n}\left\vert
\left( y_{i},y_{j}\right) \right\vert ^{q}\right) ^{\frac{1}{q}}\right] ^{%
\frac{1}{2}},
\end{multline}%
where $p>1,$ $\frac{1}{p}+\frac{1}{q}=1;$%
\begin{multline}
\sum\limits_{i=1}^{n}\left\vert \left( x,y_{i}\right) \right\vert ^{2}\leq
\left\Vert x\right\Vert \left( \sum\limits_{i=1}^{n}\left\vert \left(
x,y_{i}\right) \right\vert ^{p}\right) ^{\frac{1}{2p}}\left(
\sum\limits_{i=1}^{n}\left\vert \left( x,y_{i}\right) \right\vert
^{t}\right) ^{\frac{1}{2t}}  \label{4.5.13} \\
\times \left[ \sum\limits_{i=1}^{n}\left( \sum\limits_{j=1}^{n}\left\vert
\left( y_{i},y_{j}\right) \right\vert ^{q}\right) ^{\frac{u}{q}}\right] ^{%
\frac{1}{2u}},
\end{multline}%
where $p>1,$ $\frac{1}{p}+\frac{1}{q}=1,$ $t>1,$ $\frac{1}{t}+\frac{1}{u}=1;$%
\begin{multline}
\sum\limits_{i=1}^{n}\left\vert \left( x,y_{i}\right) \right\vert ^{2}\leq
\left\Vert x\right\Vert \left( \sum\limits_{i=1}^{n}\left\vert \left(
x,y_{i}\right) \right\vert ^{p}\right) ^{\frac{1}{2p}}\left(
\sum\limits_{i=1}^{n}\left\vert \left( x,y_{i}\right) \right\vert \right) ^{%
\frac{1}{2}}  \label{4.6.13} \\
\times \max\limits_{1\leq i\leq n}\left\{ \left(
\sum\limits_{j=1}^{n}\left\vert \left( y_{i},y_{j}\right) \right\vert
^{q}\right) ^{\frac{1}{2q}}\right\} ,
\end{multline}%
where $p>1,$ $\frac{1}{p}+\frac{1}{q}=1;$%
\begin{multline}
\sum\limits_{i=1}^{n}\left\vert \left( x,y_{i}\right) \right\vert ^{2}\leq
\left\Vert x\right\Vert \left[ \sum\limits_{i=1}^{n}\left\vert \left(
x,y_{i}\right) \right\vert \right] ^{\frac{1}{2}}\max\limits_{1\leq i\leq
n}\left\vert \left( x,y_{i}\right) \right\vert ^{\frac{1}{2}}  \label{4.7.13}
\\
\times \left( \sum\limits_{i=1}^{n}\left[ \max\limits_{1\leq j\leq
n}\left\vert \left( y_{i},y_{j}\right) \right\vert \right] \right) ^{\frac{1%
}{2}};
\end{multline}%
\begin{equation}
\sum\limits_{i=1}^{n}\left\vert \left( x,y_{i}\right) \right\vert ^{2}\leq
\left\Vert x\right\Vert \left[ \sum\limits_{i=1}^{n}\left\vert \left(
x,y_{i}\right) \right\vert ^{m}\right] ^{\frac{1}{2m}}\left[
\sum\limits_{i=1}^{n}\left[ \max\limits_{1\leq j\leq n}\left\vert \left(
y_{i},y_{j}\right) \right\vert ^{l}\right] \right] ^{\frac{1}{2l}},
\label{4.8.13}
\end{equation}%
where $m>1,$ $\frac{1}{m}+\frac{1}{l}=1;$ and 
\begin{equation}
\sum\limits_{i=1}^{n}\left\vert \left( x,y_{i}\right) \right\vert ^{2}\leq
\left\Vert x\right\Vert \sum\limits_{i=1}^{n}\left\vert \left(
x,y_{i}\right) \right\vert \max\limits_{1\leq j\leq n}\left\vert \left(
y_{i},y_{j}\right) \right\vert ^{\frac{1}{2}}.  \label{4.9.13}
\end{equation}

If in the above inequalities we assume that $\left( y_{i}\right) _{1\leq
i\leq n}=\left( e_{i}\right) _{1\leq i\leq n},$ where $\left( e_{i}\right)
_{1\leq i\leq n}$ are orthonormal vectors in the inner product space $\left(
H,\left( \cdot ,\cdot \right) \right) ,$ then from (\ref{4.1.13}) -- (\ref%
{4.9.13}) we may deduce the following inequalities similar in a sense to
Bessel's inequality: 
\begin{equation*}
\sum\limits_{i=1}^{n}\left\vert \left( x,e_{i}\right) \right\vert ^{2}\leq 
\sqrt{n}\left\Vert x\right\Vert \max\limits_{1\leq i\leq n}\left\{
\left\vert \left( x,e_{i}\right) \right\vert \right\} ;
\end{equation*}%
\begin{equation*}
\sum\limits_{i=1}^{n}\left\vert \left( x,e_{i}\right) \right\vert ^{2}\leq
n^{\frac{1}{2s}}\left\Vert x\right\Vert \max\limits_{1\leq i\leq n}\left\{
\left\vert \left( x,e_{i}\right) \right\vert ^{\frac{1}{2}}\right\} \left(
\sum\limits_{i=1}^{n}\left\vert \left( x,e_{i}\right) \right\vert
^{r}\right) ^{\frac{1}{2r}},
\end{equation*}%
where $r>1,$ $\frac{1}{r}+\frac{1}{s}=1;$%
\begin{equation*}
\sum\limits_{i=1}^{n}\left\vert \left( x,e_{i}\right) \right\vert ^{2}\leq
\left\Vert x\right\Vert \max\limits_{1\leq i\leq n}\left\{ \left\vert \left(
x,e_{i}\right) \right\vert ^{\frac{1}{2}}\right\} \left(
\sum\limits_{i=1}^{n}\left\vert \left( x,e_{i}\right) \right\vert \right) ^{%
\frac{1}{2}},
\end{equation*}%
\begin{equation*}
\sum\limits_{i=1}^{n}\left\vert \left( x,e_{i}\right) \right\vert ^{2}\leq 
\sqrt{n}\left\Vert x\right\Vert \max\limits_{1\leq i\leq n}\left\{
\left\vert \left( x,e_{i}\right) \right\vert ^{\frac{1}{2}}\right\} \left(
\sum\limits_{i=1}^{n}\left\vert \left( x,e_{i}\right) \right\vert
^{p}\right) ^{\frac{1}{2p}},
\end{equation*}%
where $p>1;$%
\begin{equation*}
\sum\limits_{i=1}^{n}\left\vert \left( x,e_{i}\right) \right\vert ^{2}\leq
n^{\frac{1}{2u}}\left\Vert x\right\Vert \left(
\sum\limits_{i=1}^{n}\left\vert \left( x,e_{i}\right) \right\vert
^{p}\right) ^{\frac{1}{2p}}\left( \sum\limits_{i=1}^{n}\left\vert \left(
x,e_{i}\right) \right\vert ^{t}\right) ^{\frac{1}{2t}},
\end{equation*}%
where $p>1,$ $t>1,$ $\frac{1}{t}+\frac{1}{u}=1;$%
\begin{equation*}
\sum\limits_{i=1}^{n}\left\vert \left( x,e_{i}\right) \right\vert ^{2}\leq
\left\Vert x\right\Vert \left( \sum\limits_{i=1}^{n}\left\vert \left(
x,e_{i}\right) \right\vert ^{p}\right) ^{\frac{1}{2p}}\left(
\sum\limits_{i=1}^{n}\left\vert \left( x,e_{i}\right) \right\vert \right) ^{%
\frac{1}{2}},\ \ p>1;
\end{equation*}%
\begin{equation*}
\sum\limits_{i=1}^{n}\left\vert \left( x,e_{i}\right) \right\vert ^{2}\leq 
\sqrt{n}\left\Vert x\right\Vert \left( \sum\limits_{i=1}^{n}\left\vert
\left( x,e_{i}\right) \right\vert \right) ^{\frac{1}{2}}\max\limits_{1\leq
i\leq n}\left\{ \left\vert \left( x,e_{i}\right) \right\vert ^{\frac{1}{2}%
}\right\} ;
\end{equation*}%
\begin{equation*}
\sum\limits_{i=1}^{n}\left\vert \left( x,e_{i}\right) \right\vert ^{2}\leq
n^{\frac{1}{2l}}\left\Vert x\right\Vert \left[ \sum\limits_{i=1}^{n}\left%
\vert \left( x,e_{i}\right) \right\vert ^{m}\right] ^{\frac{1}{m}},\ \ \
m>1,\ \frac{1}{m}+\frac{1}{l}=1;
\end{equation*}%
\begin{equation*}
\sum\limits_{i=1}^{n}\left\vert \left( x,e_{i}\right) \right\vert ^{2}\leq
\left\Vert x\right\Vert \sum\limits_{i=1}^{n}\left\vert \left(
x,e_{i}\right) \right\vert .
\end{equation*}

Corollaries \ref{c3.2.13} -- \ref{c3.5.13} will produce the following
results which do not contain the Fourier coefficients in the right side of
the inequality.

Indeed, if one chooses $c_{i}=\overline{\left( x,y_{i}\right) }$ in (\ref%
{3.3.13}), then 
\begin{equation*}
\left( \sum\limits_{i=1}^{n}\left\vert \left( x,y_{i}\right) \right\vert
^{2}\right) ^{2}\leq \left\Vert x\right\Vert ^{2}n^{\frac{1}{p}+\frac{1}{t}%
-1}\sum\limits_{i=1}^{n}\left\vert \left( x,y_{i}\right) \right\vert ^{2}%
\left[ \sum\limits_{i=1}^{n}\left( \sum\limits_{j=1}^{n}\left\vert \left(
y_{i},y_{j}\right) \right\vert ^{q}\right) ^{\frac{u}{q}}\right] ^{\frac{1}{u%
}},
\end{equation*}%
giving the following Bombieri type inequality: 
\begin{equation*}
\sum\limits_{i=1}^{n}\left\vert \left( x,y_{i}\right) \right\vert ^{2}\leq
n^{\frac{1}{p}+\frac{1}{t}-1}\left\Vert x\right\Vert ^{2}\left[
\sum\limits_{i=1}^{n}\left( \sum\limits_{j=1}^{n}\left\vert \left(
y_{i},y_{j}\right) \right\vert ^{q}\right) ^{\frac{u}{q}}\right] ^{\frac{1}{u%
}},
\end{equation*}%
where $1<p\leq 2,$ $1<t\leq 2,$ $\frac{1}{p}+\frac{1}{q}=1,$ $\frac{1}{t}+%
\frac{1}{u}=1.$

If in this inequality we consider $p=q=t=u=2,$ then 
\begin{equation*}
\sum\limits_{i=1}^{n}\left\vert \left( x,y_{i}\right) \right\vert ^{2}\leq
\left\Vert x\right\Vert ^{2}\left( \sum\limits_{i,j=1}^{n}\left\vert \left(
y_{i},y_{j}\right) \right\vert ^{2}\right) ^{\frac{1}{2}}.
\end{equation*}%
For a different proof of this result see also \cite{5b.13}.

In a similar way, if $c_{i}=\overline{\left( x,y_{i}\right) }$ in (\ref%
{3.6.13}), then 
\begin{equation*}
\sum\limits_{i=1}^{n}\left\vert \left( x,y_{i}\right) \right\vert ^{2}\leq
n^{\frac{1}{m}}\left\Vert x\right\Vert ^{2}\left( \sum\limits_{i=1}^{n}\left[
\max\limits_{1\leq j\leq n}\left\vert \left( y_{i},y_{j}\right) \right\vert %
\right] ^{l}\right) ^{\frac{1}{l}},
\end{equation*}%
where $m>1,$ $\frac{1}{m}+\frac{1}{l}=1.$

Finally, if $c_{i}=\overline{\left( x,y_{i}\right) }$ $\left( i=1,\dots
,n\right) ,$ is taken in (\ref{3.7.13}), then 
\begin{equation*}
\sum\limits_{i=1}^{n}\left\vert \left( x,y_{i}\right) \right\vert ^{2}\leq
n\left\Vert x\right\Vert ^{2}\max\limits_{1\leq i,j\leq n}\left\vert \left(
y_{i},y_{j}\right) \right\vert .
\end{equation*}

\begin{remark}
\label{r4.1.13}Let us compare Bombieri's result 
\begin{equation}
\sum\limits_{i=1}^{n}\left\vert \left( x,y_{i}\right) \right\vert ^{2}\leq
\left\Vert x\right\Vert ^{2}\max\limits_{1\leq i\leq n}\left\{
\sum\limits_{j=1}^{n}\left\vert \left( y_{i},y_{j}\right) \right\vert
\right\}  \label{4.23.13}
\end{equation}%
with our result 
\begin{equation}
\sum\limits_{i=1}^{n}\left\vert \left( x,y_{i}\right) \right\vert ^{2}\leq
\left\Vert x\right\Vert ^{2}\left\{ \sum\limits_{i,j=1}^{n}\left\vert \left(
y_{i},y_{j}\right) \right\vert ^{2}\right\} ^{\frac{1}{2}}.  \label{4.24.13}
\end{equation}%
Denote 
\begin{equation*}
M_{1}:=\max_{1\leq i\leq n}\left\{ \sum_{j=1}^{n}\left\vert \left(
y_{i},y_{j}\right) \right\vert \right\}
\end{equation*}%
and 
\begin{equation*}
M_{2}:=\left[ \sum_{i,j=1}^{n}\left\vert \left( y_{i},y_{j}\right)
\right\vert ^{2}\right] ^{\frac{1}{2}}.
\end{equation*}

If we choose the inner product space $H=\mathbb{R}$, $\left( x,y\right) :=xy$
and $n=2,$ then for $y_{1}=a,$ $y_{2}=b,$ $a,b>0,$ we have 
\begin{align*}
M_{1}& =\max \left\{ a^{2}+ab,ab+b^{2}\right\} =\left( a+b\right) \max
\left( a,b\right) , \\
M_{2}& =\left( a^{4}+a^{2}b^{2}+a^{2}b^{2}+b^{4}\right) ^{\frac{1}{2}%
}=a^{2}+b^{2}.
\end{align*}

Assume that $a\geq b.$ Then $M_{1}=a^{2}+ab\geq a^{2}+b^{2}=M_{2},$ showing
that, in this case, the bound provided by (\ref{4.24.13}) is better than the
bound provided by (\ref{4.23.13}). If $\left( y_{i}\right) _{1\leq i\leq n}$
are orthonormal vectors, then $M_{1}=1,$ $M_{2}=\sqrt{n}$, showing that in
this case the Bombieri inequality (which becomes Bessel's inequality)
provides a better bound than (\ref{4.24.13}).
\end{remark}

\newpage

\section{Pe\v{c}ari\'{c} Type Inequalities}

\subsection{Introduction}

In 1992, J.E. Pe\v{c}ari\'{c} \cite{3b.14} proved the following inequality
for vectors in complex inner product spaces $\left( H;\left( \cdot ,\cdot
\right) \right) $.

\begin{theorem}
\label{t1.1.14}Suppose that $x,y_{1},\dots ,y_{n}$ are vectors in $H$ and $%
c_{1},\dots ,c_{n}$ are complex numbers. Then the following inequalities 
\begin{align}
\left\vert \sum\limits_{i=1}^{n}c_{i}\left( x,y_{i}\right) \right\vert ^{2}&
\leq \left\Vert x\right\Vert ^{2}\sum\limits_{i=1}^{n}\left\vert
c_{i}\right\vert ^{2}\left( \sum\limits_{j=1}^{n}\left\vert \left(
y_{i},y_{j}\right) \right\vert \right)  \label{1.1.14} \\
& \leq \left\Vert x\right\Vert ^{2}\sum\limits_{i=1}^{n}\left\vert
c_{i}\right\vert ^{2}\max_{1\leq i\leq n}\left(
\sum\limits_{j=1}^{n}\left\vert \left( y_{i},y_{j}\right) \right\vert
\right) ,  \notag
\end{align}%
hold.
\end{theorem}

He also showed that for $c_{i}=\overline{\left( x,y_{i}\right) },$ $i\in
\left\{ 1,\dots ,n\right\} ,$ one gets 
\begin{align*}
\left( \sum\limits_{i=1}^{n}\left\vert \left( x,y_{i}\right) \right\vert
^{2}\right) ^{2}& \leq \left\Vert x\right\Vert
^{2}\sum\limits_{i=1}^{n}\left\vert \left( x,y_{i}\right) \right\vert
^{2}\left( \sum\limits_{j=1}^{n}\left\vert \left( y_{i},y_{j}\right)
\right\vert \right) \\
& \leq \left\Vert x\right\Vert ^{2}\sum\limits_{i=1}^{n}\left\vert \left(
x,y_{i}\right) \right\vert ^{2}\max_{1\leq i\leq n}\left(
\sum\limits_{j=1}^{n}\left\vert \left( y_{i},y_{j}\right) \right\vert
\right) ,
\end{align*}%
which improves Bombieri's result \cite{1b.14} (see also \cite[p. 394]{2b.14}%
) 
\begin{equation}
\sum\limits_{i=1}^{n}\left\vert \left( x,y_{i}\right) \right\vert ^{2}\leq
\left\Vert x\right\Vert ^{2}\max_{1\leq i\leq n}\left(
\sum\limits_{j=1}^{n}\left\vert \left( y_{i},y_{j}\right) \right\vert
\right) .  \label{1.3.14}
\end{equation}%
Note that (\ref{1.3.14}) is in its turn a natural generalisation of \textit{%
Bessel's inequality} 
\begin{equation}
\sum\limits_{i=1}^{n}\left\vert \left( x,e_{i}\right) \right\vert ^{2}\leq
\left\Vert x\right\Vert ^{2},\ \ x\in H,  \label{1.4.14}
\end{equation}%
which holds for the orthornormal vectors $\left( e_{i}\right) _{1\leq i\leq
n}.$

In this section, by following \cite{15NSSD}, we point out other related
results to Pe\v{c}ari\'{c}'s inequality (\ref{1.1.14}) than the ones stated
in the previous sections. Some results of Bombieri type are also mentioned.

\subsection{Some Norm Inequalities}

We start with the following lemma that is interesting in its own right \cite%
{15NSSD}.

\begin{lemma}
\label{l2.1.14}Let $z_{1},\dots ,z_{n}\in H$ and $\alpha _{1},\dots ,\alpha
_{n}\in \mathbb{K}$. Then one has the inequalities: 
\begin{multline}
\left\Vert \sum_{i=1}^{n}\alpha _{i}z_{i}\right\Vert ^{2}\leq \left(
\sum\limits_{i=1}^{n}\left\vert \alpha _{i}\right\vert ^{p}\left(
\sum\limits_{j=1}^{n}\left\vert \left( z_{i},z_{j}\right) \right\vert
\right) \right) ^{\frac{1}{p}}  \label{2.1.14} \\
\times \left( \sum\limits_{i=1}^{n}\left\vert \alpha _{i}\right\vert
^{q}\left( \sum\limits_{j=1}^{n}\left\vert \left( z_{i},z_{j}\right)
\right\vert \right) \right) ^{\frac{1}{q}}\leq \left\{ 
\begin{array}{c}
A \\ 
B \\ 
C%
\end{array}%
\right. ,
\end{multline}%
where%
\begin{equation*}
A:=\left\{ 
\begin{array}{l}
\max\limits_{1\leq i\leq n}\left\vert \alpha _{i}\right\vert
^{2}\sum_{i,j=1}^{n}\left\vert \left( z_{i},z_{j}\right) \right\vert ; \\ 
\vspace{-0.05in} \\ 
\max\limits_{1\leq i\leq n}\left\vert \alpha _{i}\right\vert \left(
\sum_{i=1}^{n}\left\vert \alpha _{i}\right\vert ^{\gamma q}\right) ^{\frac{1%
}{\gamma q}}\left( \sum_{i,j=1}^{n}\left\vert \left( z_{i},z_{j}\right)
\right\vert \right) ^{\frac{1}{p}}\left( \sum_{i=1}^{n}\left(
\sum_{j=1}^{n}\left\vert \left( z_{i},z_{j}\right) \right\vert \right)
^{\delta }\right) ^{\frac{1}{\delta q}}, \\ 
\hfill \hfill \ \ \text{if}\ \gamma >1,\ \frac{1}{\gamma }+\frac{1}{\delta }%
=1; \\ 
\vspace{-0.05in} \\ 
\max\limits_{1\leq i\leq n}\left\vert \alpha _{i}\right\vert \left(
\sum_{i=1}^{n}\left\vert \alpha _{i}\right\vert ^{q}\right) ^{\frac{1}{q}%
}\left( \sum_{i,j=1}^{n}\left\vert \left( z_{i},z_{j}\right) \right\vert
\right) ^{\frac{1}{p}}\max\limits_{1\leq i\leq n}\left(
\sum_{j=1}^{n}\left\vert \left( z_{i},z_{j}\right) \right\vert \right) ^{%
\frac{1}{q}};%
\end{array}%
\right.
\end{equation*}%
\vspace{-0.1in}%
\begin{equation*}
B:=\left\{ 
\begin{array}{l}
\max\limits_{1\leq i\leq n}\left\vert \alpha _{i}\right\vert \left(
\sum_{i=1}^{n}\left\vert \alpha _{i}\right\vert ^{\alpha p}\right) ^{\frac{1%
}{\alpha p}}\left( \sum_{i,j=1}^{n}\left\vert \left( z_{i},z_{j}\right)
\right\vert \right) ^{\frac{1}{q}}\left( \sum_{i=1}^{n}\left(
\sum_{j=1}^{n}\left\vert \left( z_{i},z_{j}\right) \right\vert \right)
^{\beta }\right) ^{\frac{1}{\beta q}}, \\ 
\hfill \hfill \ \ \text{if}\ \alpha >1,\ \frac{1}{\alpha }+\frac{1}{\beta }%
=1; \\ 
\vspace{-0.05in} \\ 
\left( \sum_{i=1}^{n}\left\vert \alpha _{i}\right\vert ^{\alpha p}\right) ^{%
\frac{1}{\alpha p}}\left( \sum_{i=1}^{n}\left\vert \alpha _{i}\right\vert
^{\gamma q}\right) ^{\frac{1}{\gamma q}} \\ 
\ \ \ \ \ \times \left( \sum_{i=1}^{n}\left( \sum_{j=1}^{n}\left\vert \left(
z_{i},z_{j}\right) \right\vert \right) ^{\beta }\right) ^{\frac{1}{p\beta }%
}\left( \sum_{i=1}^{n}\left( \sum_{j=1}^{n}\left\vert \left(
z_{i},z_{j}\right) \right\vert \right) ^{\delta }\right) ^{\frac{1}{\delta q}%
}\hfill \\ 
\hfill \text{if}\ \alpha >1,\ \frac{1}{\alpha }+\frac{1}{\beta }=1\text{ and 
}\ \gamma >1,\ \frac{1}{\gamma }+\frac{1}{\delta }=1; \\ 
\vspace{-0.05in} \\ 
\left( \sum_{i=1}^{n}\left\vert \alpha _{i}\right\vert ^{q}\right) ^{\frac{1%
}{q}}\left( \sum_{i=1}^{n}\left\vert \alpha _{i}\right\vert ^{\alpha
p}\right) ^{\frac{1}{\alpha p}}\max\limits_{1\leq i\leq n}\left(
\sum_{j=1}^{n}\left\vert \left( z_{i},z_{j}\right) \right\vert \right) ^{%
\frac{1}{q}} \\ 
\qquad \times \left( \sum_{i=1}^{n}\left( \sum_{j=1}^{n}\left\vert \left(
z_{i},z_{j}\right) \right\vert \right) ^{\beta }\right) ^{\frac{1}{p\beta }%
},\hfill \hfill \ \ \text{if}\ \alpha >1,\ \frac{1}{\alpha }+\frac{1}{\beta }%
=1;%
\end{array}%
\right.
\end{equation*}%
and%
\begin{equation*}
C:=\left\{ 
\begin{array}{l}
\max\limits_{1\leq i\leq n}\left\vert \alpha _{i}\right\vert \left(
\sum_{i=1}^{n}\left\vert \alpha _{i}\right\vert ^{p}\right) ^{\frac{1}{p}%
}\max\limits_{1\leq i\leq n}\left( \sum_{j=1}^{n}\left\vert \left(
z_{i},z_{j}\right) \right\vert \right) ^{\frac{1}{p}}\left(
\sum_{i,j=1}^{n}\left\vert \left( z_{i},z_{j}\right) \right\vert \right) ^{%
\frac{1}{q}}; \\ 
\vspace{-0.05in} \\ 
\left( \sum_{i=1}^{n}\left\vert \alpha _{i}\right\vert ^{p}\right) ^{\frac{1%
}{p}}\left( \sum_{i=1}^{n}\left\vert \alpha _{i}\right\vert ^{\gamma
q}\right) ^{\frac{1}{\gamma q}}\max\limits_{1\leq i\leq n}\left(
\sum_{j=1}^{n}\left\vert \left( z_{i},z_{j}\right) \right\vert \right) ^{%
\frac{1}{p}} \\ 
\qquad \times \left( \sum_{i=1}^{n}\left( \sum_{j=1}^{n}\left\vert \left(
z_{i},z_{j}\right) \right\vert \right) ^{\delta }\right) ^{\frac{1}{\delta q}%
},\hfill \ \ \text{if}\ \gamma >1,\ \frac{1}{\gamma }+\frac{1}{\delta }=1;
\\ 
\vspace{-0.05in} \\ 
\left( \sum_{i=1}^{n}\left\vert \alpha _{i}\right\vert ^{p}\right) ^{\frac{1%
}{p}}\left( \sum_{i=1}^{n}\left\vert \alpha _{i}\right\vert ^{q}\right) ^{%
\frac{1}{q}}\max\limits_{1\leq i\leq n}\left( \sum_{j=1}^{n}\left\vert
\left( z_{i},z_{j}\right) \right\vert \right) ,%
\end{array}%
\right.
\end{equation*}%
where $p>1,$ $\frac{1}{p}+\frac{1}{q}=1.$
\end{lemma}

\begin{proof}
We observe that 
\begin{align}
\left\Vert \sum_{i=1}^{n}\alpha _{i}z_{i}\right\Vert ^{2}& =\left(
\sum_{i=1}^{n}\alpha _{i}z_{i},\sum_{j=1}^{n}\alpha _{j}z_{j}\right)
\label{2.2.14} \\
& =\sum_{i=1}^{n}\sum_{j=1}^{n}\alpha _{i}\overline{\alpha _{j}}\left(
z_{i},z_{j}\right) =\left\vert \sum_{i=1}^{n}\sum_{j=1}^{n}\alpha _{i}%
\overline{\alpha _{j}}\left( z_{i},z_{j}\right) \right\vert  \notag \\
& \leq \sum_{i=1}^{n}\sum_{j=1}^{n}\left\vert \alpha _{i}\right\vert
\left\vert \alpha _{j}\right\vert \left\vert \left( z_{i},z_{j}\right)
\right\vert =:M.  \notag
\end{align}%
If one uses the H\"{o}lder inequality for double sums, i.e., we recall it 
\begin{equation*}
\sum\limits_{i,j=1}^{n}m_{ij}a_{ij}b_{ij}\leq \left(
\sum\limits_{i,j=1}^{n}m_{ij}a_{ij}^{p}\right) ^{\frac{1}{p}}\left(
\sum\limits_{i,j=1}^{n}m_{ij}b_{ij}^{q}\right) ^{\frac{1}{q}},
\end{equation*}%
where $m_{ij},a_{ij},b_{ij}\geq 0,$ $\frac{1}{p}+\frac{1}{q}=1,$ $p>1;$ then 
\begin{align}
M& \leq \left( \sum\limits_{i,j=1}^{n}\left\vert \left( z_{i},z_{j}\right)
\right\vert \left\vert \alpha _{i}\right\vert ^{p}\right) ^{\frac{1}{p}%
}\left( \sum\limits_{i,j=1}^{n}\left\vert \left( z_{i},z_{j}\right)
\right\vert \left\vert \alpha _{i}\right\vert ^{q}\right) ^{\frac{1}{q}}
\label{2.4.14} \\
& =\left( \sum_{i=1}^{n}\left\vert \alpha _{i}\right\vert ^{p}\left(
\sum_{j=1}^{n}\left\vert \left( z_{i},z_{j}\right) \right\vert \right)
\right) ^{\frac{1}{p}}\left( \sum_{i=1}^{n}\left\vert \alpha _{i}\right\vert
^{q}\left( \sum_{j=1}^{n}\left\vert \left( z_{i},z_{j}\right) \right\vert
\right) \right) ^{\frac{1}{q}},  \notag
\end{align}%
and the first inequality in (\ref{2.1.14}) is proved.

Observe that 
\begin{equation*}
\sum_{i=1}^{n}\left\vert \alpha _{i}\right\vert ^{p}\left(
\sum_{j=1}^{n}\left\vert \left( z_{i},z_{j}\right) \right\vert \right) \leq
\left\{ 
\begin{array}{l}
\max\limits_{1\leq i\leq n}\left\vert \alpha _{i}\right\vert
^{p}\sum_{i,j=1}^{n}\left\vert \left( z_{i},z_{j}\right) \right\vert ; \\ 
\\ 
\left( \sum_{i=1}^{n}\left\vert \alpha _{i}\right\vert ^{\alpha p}\right) ^{%
\frac{1}{\alpha }}\left( \sum_{i=1}^{n}\left( \sum_{j=1}^{n}\left\vert
\left( z_{i},z_{j}\right) \right\vert \right) ^{\beta }\right) ^{\frac{1}{%
\beta }} \\ 
\hfill \ \text{if}\ \alpha >1,\ \frac{1}{\alpha }+\frac{1}{\beta }=1; \\ 
\\ 
\sum_{i=1}^{n}\left\vert \alpha _{i}\right\vert ^{p}\max\limits_{1\leq i\leq
n}\left( \sum_{j=1}^{n}\left\vert \left( z_{i},z_{j}\right) \right\vert
\right) ;%
\end{array}%
\right.
\end{equation*}%
giving 
\begin{multline}
\left( \sum_{i=1}^{n}\left\vert \alpha _{i}\right\vert ^{p}\left(
\sum_{j=1}^{n}\left\vert \left( z_{i},z_{j}\right) \right\vert \right)
\right) ^{\frac{1}{p}}  \label{2.5.14} \\
\leq \left\{ 
\begin{array}{l}
\max\limits_{1\leq i\leq n}\left\vert \alpha _{i}\right\vert \left(
\sum_{i,j=1}^{n}\left\vert \left( z_{i},z_{j}\right) \right\vert \right) ^{%
\frac{1}{p}};\hfill \\ 
\\ 
\left( \sum_{i=1}^{n}\left\vert \alpha _{i}\right\vert ^{\alpha p}\right) ^{%
\frac{1}{\alpha p}}\left( \sum_{i=1}^{n}\left( \sum_{j=1}^{n}\left\vert
\left( z_{i},z_{j}\right) \right\vert \right) ^{\beta }\right) ^{\frac{1}{%
\beta p}}\hfill \  \\ 
\hfill \text{if}\ \alpha >1,\ \frac{1}{\alpha }+\frac{1}{\beta }=1; \\ 
\\ 
\left( \sum_{i=1}^{n}\left\vert \alpha _{i}\right\vert ^{p}\right) ^{\frac{1%
}{p}}\max\limits_{1\leq i\leq n}\left( \sum_{j=1}^{n}\left\vert \left(
z_{i},z_{j}\right) \right\vert \right) ^{\frac{1}{p}}.%
\end{array}%
\right.
\end{multline}%
Similarly, we have 
\begin{multline}
\left( \sum_{i=1}^{n}\left\vert \alpha _{i}\right\vert ^{q}\left(
\sum_{j=1}^{n}\left\vert \left( z_{i},z_{j}\right) \right\vert \right)
\right) ^{\frac{1}{q}}  \label{2.6.14} \\
\leq \left\{ 
\begin{array}{l}
\max\limits_{1\leq i\leq n}\left\vert \alpha _{i}\right\vert \left(
\sum_{i,j=1}^{n}\left\vert \left( z_{i},z_{j}\right) \right\vert \right) ^{%
\frac{1}{q}}\hfill \\ 
\\ 
\left( \sum_{i=1}^{n}\left\vert \alpha _{i}\right\vert ^{\gamma q}\right) ^{%
\frac{1}{\gamma q}}\left( \sum_{i=1}^{n}\left( \sum_{j=1}^{n}\left\vert
\left( z_{i},z_{j}\right) \right\vert \right) ^{\delta }\right) ^{\frac{1}{%
\delta q}}\hfill \  \\ 
\hfill \text{if}\ \gamma >1,\ \frac{1}{\gamma }+\frac{1}{\delta }=1; \\ 
\\ 
\left( \sum_{i=1}^{n}\left\vert \alpha _{i}\right\vert ^{q}\right) ^{\frac{1%
}{q}}\max\limits_{1\leq i\leq n}\left( \sum_{j=1}^{n}\left\vert \left(
z_{i},z_{j}\right) \right\vert \right) ^{\frac{1}{q}}.%
\end{array}%
\right.
\end{multline}%
Using (\ref{2.2.14}) and (\ref{2.5.14}) -- (\ref{2.6.14}), we deduce the 9
inequalities in the second part of (\ref{2.1.14}).
\end{proof}

If we choose $p=q=2,$ then the following result holds \cite{15NSSD}.

\begin{corollary}
\label{c2.2.14}If $z_{1},\dots ,z_{n}\in H$ and $\alpha _{1},\dots ,\alpha
_{n}\in \mathbb{K}$, then one has 
\begin{align}
\left\Vert \sum_{i=1}^{n}\alpha _{i}z_{i}\right\Vert ^{2}& \leq
\sum\limits_{i=1}^{n}\left\vert \alpha _{i}\right\vert ^{2}\left(
\sum\limits_{j=1}^{n}\left\vert \left( z_{i},z_{j}\right) \right\vert \right)
\label{2.7.14} \\
& \leq \left\{ 
\begin{array}{c}
D \\ 
E \\ 
F%
\end{array}%
\right. ,  \notag
\end{align}%
where%
\begin{equation*}
D:=\left\{ 
\begin{array}{l}
\max\limits_{1\leq i\leq n}\left\vert \alpha _{i}\right\vert
^{2}\sum\limits_{i,j=1}^{n}\left\vert \left( z_{i},z_{j}\right) \right\vert ;
\\ 
\\ 
\max\limits_{1\leq i\leq n}\left\vert \alpha _{i}\right\vert \left(
\sum\limits_{i=1}^{n}\left\vert \alpha _{i}\right\vert ^{2\gamma }\right) ^{%
\frac{1}{2\gamma }}\left( \sum\limits_{i,j=1}^{n}\left\vert \left(
z_{i},z_{j}\right) \right\vert \right) ^{\frac{1}{2}}\left(
\sum\limits_{i=1}^{n}\left( \sum\limits_{j=1}^{n}\left\vert \left(
z_{i},z_{j}\right) \right\vert \right) ^{\delta }\right) ^{\frac{1}{2\delta }%
}, \\ 
\hfill \ \ \text{if}\ \gamma >1,\ \frac{1}{\gamma }+\frac{1}{\delta }=1; \\ 
\\ 
\max\limits_{1\leq i\leq n}\left\vert \alpha _{i}\right\vert \left(
\sum\limits_{i=1}^{n}\left\vert \alpha _{i}\right\vert ^{2}\right) ^{\frac{1%
}{2}}\left( \sum\limits_{i,j=1}^{n}\left\vert \left( z_{i},z_{j}\right)
\right\vert \right) ^{\frac{1}{2}}\max\limits_{1\leq i\leq n}\left(
\sum\limits_{j=1}^{n}\left\vert \left( z_{i},z_{j}\right) \right\vert
\right) ^{\frac{1}{2}};%
\end{array}%
\right.
\end{equation*}%
\begin{equation*}
E:=\left\{ 
\begin{array}{l}
\max\limits_{1\leq i\leq n}\left\vert \alpha _{i}\right\vert \left(
\sum\limits_{i=1}^{n}\left\vert \alpha _{i}\right\vert ^{2\alpha }\right) ^{%
\frac{1}{2\alpha }}\left( \sum\limits_{i,j=1}^{n}\left\vert \left(
z_{i},z_{j}\right) \right\vert \right) ^{\frac{1}{2}}\left(
\sum\limits_{i=1}^{n}\left( \sum\limits_{j=1}^{n}\left\vert \left(
z_{i},z_{j}\right) \right\vert \right) ^{\beta }\right) ^{\frac{1}{2\beta }},
\\ 
\hfill \ \ \text{if}\ \alpha >1,\ \frac{1}{\alpha }+\frac{1}{\beta }=1; \\ 
\\ 
\left( \sum\limits_{i=1}^{n}\left\vert \alpha _{i}\right\vert ^{2\alpha
}\right) ^{\frac{1}{2\alpha }}\left( \sum\limits_{i=1}^{n}\left\vert \alpha
_{i}\right\vert ^{2\gamma }\right) ^{\frac{1}{2\gamma }}\left(
\sum\limits_{i=1}^{n}\left( \sum\limits_{j=1}^{n}\left\vert \left(
z_{i},z_{j}\right) \right\vert \right) ^{\beta }\right) ^{\frac{1}{2\beta }}
\\ 
\ \ \ \ \ \times \left( \sum\limits_{i=1}^{n}\left(
\sum\limits_{j=1}^{n}\left\vert \left( z_{i},z_{j}\right) \right\vert
\right) ^{\delta }\right) ^{\frac{1}{2\delta }}\hfill \ \text{if}\ \alpha
>1,\ \frac{1}{\alpha }+\frac{1}{\beta }=1\text{ and }\ \gamma >1,\ \frac{1}{%
\gamma }+\frac{1}{\delta }=1; \\ 
\\ 
\left( \sum\limits_{i=1}^{n}\left\vert \alpha _{i}\right\vert ^{2}\right) ^{%
\frac{1}{2}}\left( \sum\limits_{i=1}^{n}\left\vert \alpha _{i}\right\vert
^{2\alpha }\right) ^{\frac{1}{2\alpha }}\max\limits_{1\leq i\leq n}\left(
\sum\limits_{j=1}^{n}\left\vert \left( z_{i},z_{j}\right) \right\vert
\right) ^{\frac{1}{2}}\left( \sum\limits_{i=1}^{n}\left(
\sum\limits_{j=1}^{n}\left\vert \left( z_{i},z_{j}\right) \right\vert
\right) ^{\beta }\right) ^{\frac{1}{2\beta }}, \\ 
\hfill \ \ \text{if}\ \alpha >1,\ \frac{1}{\alpha }+\frac{1}{\beta }=1;%
\end{array}%
\right.
\end{equation*}%
and%
\begin{equation*}
F:=\left\{ 
\begin{array}{l}
\max\limits_{1\leq i\leq n}\left\vert \alpha _{i}\right\vert \left(
\sum\limits_{i=1}^{n}\left\vert \alpha _{i}\right\vert ^{2}\right) ^{\frac{1%
}{2}}\max\limits_{1\leq i\leq n}\left( \sum\limits_{j=1}^{n}\left\vert
\left( z_{i},z_{j}\right) \right\vert \right) ^{\frac{1}{2}}\left(
\sum\limits_{i,j=1}^{n}\left\vert \left( z_{i},z_{j}\right) \right\vert
\right) ^{\frac{1}{2}}; \\ 
\\ 
\left( \sum\limits_{i=1}^{n}\left\vert \alpha _{i}\right\vert ^{2}\right) ^{%
\frac{1}{2}}\left( \sum\limits_{i=1}^{n}\left\vert \alpha _{i}\right\vert
^{2\gamma }\right) ^{\frac{1}{2\gamma }}\max\limits_{1\leq i\leq n}\left(
\sum\limits_{j=1}^{n}\left\vert \left( z_{i},z_{j}\right) \right\vert
\right) ^{\frac{1}{2}}\left( \sum\limits_{i=1}^{n}\left(
\sum\limits_{j=1}^{n}\left\vert \left( z_{i},z_{j}\right) \right\vert
\right) ^{\delta }\right) ^{\frac{1}{2\delta }}, \\ 
\hfill \ \ \text{if}\ \gamma >1,\ \frac{1}{\gamma }+\frac{1}{\delta }=1; \\ 
\\ 
\sum\limits_{i=1}^{n}\left\vert \alpha _{i}\right\vert
^{2}\max\limits_{1\leq i\leq n}\left( \sum\limits_{j=1}^{n}\left\vert \left(
z_{i},z_{j}\right) \right\vert \right) .%
\end{array}%
\right.
\end{equation*}
\end{corollary}

\subsection{Some Pe\v{c}ari\'{c} Type Inequalities}

We are now able to point out the following result obtained in \cite{15NSSD},
which complements and generalises the inequality (\ref{1.1.14}) due to J. Pe%
\v{c}ari\'{c}.

\begin{theorem}
\label{t3.1.14}Let $x,y_{1},\dots ,y_{n}$ be vectors of an inner product
space $\left( H;\left( \cdot ,\cdot \right) \right) $ and $c_{1},\dots
,c_{n}\in \mathbb{K}$. Then one has the inequalities: 
\begin{multline}
\left\vert \sum\limits_{i=1}^{n}c_{i}\left( x,y_{i}\right) \right\vert ^{2}
\label{3.1.14} \\
\leq \left\Vert x\right\Vert ^{2}\left( \sum\limits_{i=1}^{n}\left\vert
c_{i}\right\vert ^{p}\left( \sum\limits_{j=1}^{n}\left\vert \left(
y_{i},y_{j}\right) \right\vert \right) \right) ^{\frac{1}{p}}\left(
\sum\limits_{i=1}^{n}\left\vert c_{i}\right\vert ^{q}\left(
\sum\limits_{j=1}^{n}\left\vert \left( y_{i},y_{j}\right) \right\vert
\right) \right) ^{\frac{1}{q}} \\
\leq \left\Vert x\right\Vert ^{2}\times \left\{ 
\begin{array}{c}
G \\ 
H \\ 
I%
\end{array}%
\right. ,
\end{multline}%
where%
\begin{equation*}
G:=\left\{ 
\begin{array}{l}
\max\limits_{1\leq i\leq n}\left\vert c_{i}\right\vert
^{2}\sum_{i,j=1}^{n}\left\vert \left( y_{i},y_{j}\right) \right\vert ; \\ 
\\ 
\max\limits_{1\leq i\leq n}\left\vert c_{i}\right\vert \left(
\sum_{i=1}^{n}\left\vert c_{i}\right\vert ^{\gamma q}\right) ^{\frac{1}{%
\gamma q}}\left( \sum_{i,j=1}^{n}\left\vert \left( y_{i},y_{j}\right)
\right\vert \right) ^{\frac{1}{p}} \\ 
\qquad \times \left( \sum_{i=1}^{n}\left( \sum_{j=1}^{n}\left\vert \left(
y_{i},y_{j}\right) \right\vert \right) ^{\delta }\right) ^{\frac{1}{\delta q}%
},\hfill \ \ \text{if}\ \gamma >1,\ \frac{1}{\gamma }+\frac{1}{\delta }=1;
\\ 
\\ 
\max\limits_{1\leq i\leq n}\left\vert c_{i}\right\vert \left(
\sum_{i=1}^{n}\left\vert c_{i}\right\vert ^{q}\right) ^{\frac{1}{q}}\left(
\sum_{i,j=1}^{n}\left\vert \left( y_{i},y_{j}\right) \right\vert \right) ^{%
\frac{1}{p}}\max\limits_{1\leq i\leq n}\left( \sum_{j=1}^{n}\left\vert
\left( y_{i},y_{j}\right) \right\vert \right) ^{\frac{1}{q}};%
\end{array}%
\right.
\end{equation*}%
\begin{equation*}
H:=\left\{ 
\begin{array}{l}
\max\limits_{1\leq i\leq n}\left\vert c_{i}\right\vert \left(
\sum_{i=1}^{n}\left\vert c_{i}\right\vert ^{\alpha p}\right) ^{\frac{1}{%
\alpha p}}\left( \sum_{i,j=1}^{n}\left\vert \left( y_{i},y_{j}\right)
\right\vert \right) ^{\frac{1}{q}} \\ 
\qquad \times \left( \sum_{i=1}^{n}\left( \sum_{j=1}^{n}\left\vert \left(
y_{i},y_{j}\right) \right\vert \right) ^{\beta }\right) ^{\frac{1}{p\beta }%
},\hfill \ \ \text{if}\ \alpha >1,\ \frac{1}{\alpha }+\frac{1}{\beta }=1; \\ 
\\ 
\left( \sum_{i=1}^{n}\left\vert c_{i}\right\vert ^{\alpha p}\right) ^{\frac{1%
}{\alpha p}}\left( \sum_{i=1}^{n}\left\vert c_{i}\right\vert ^{\gamma
q}\right) ^{\frac{1}{\gamma q}}\left( \sum_{i=1}^{n}\left(
\sum_{j=1}^{n}\left\vert \left( y_{i},y_{j}\right) \right\vert \right)
^{\beta }\right) ^{\frac{1}{p\beta }} \\ 
\qquad \times \left( \sum_{i=1}^{n}\left( \sum_{j=1}^{n}\left\vert \left(
y_{i},y_{j}\right) \right\vert \right) ^{\delta }\right) ^{\frac{1}{\delta q}%
}\  \\ 
\hfill \text{if}\ \alpha >1,\ \frac{1}{\alpha }+\frac{1}{\beta }=1\text{ and 
}\gamma >1,\ \frac{1}{\gamma }+\frac{1}{\delta }=1; \\ 
\\ 
\left( \sum_{i=1}^{n}\left\vert c_{i}\right\vert ^{q}\right) ^{\frac{1}{q}%
}\left( \sum_{i=1}^{n}\left\vert c_{i}\right\vert ^{\alpha p}\right) ^{\frac{%
1}{\alpha p}}\max\limits_{1\leq i\leq n}\left( \sum_{j=1}^{n}\left\vert
\left( y_{i},y_{j}\right) \right\vert \right) ^{\frac{1}{q}} \\ 
\ \ \ \ \ \times \left( \sum_{i=1}^{n}\left( \sum_{j=1}^{n}\left\vert \left(
y_{i},y_{j}\right) \right\vert \right) ^{\beta }\right) ^{\frac{1}{p\beta }%
},\hfill \ \ \text{if}\ \alpha >1,\ \frac{1}{\alpha }+\frac{1}{\beta }=1;%
\end{array}%
\right.
\end{equation*}%
and%
\begin{equation*}
I:=\left\{ 
\begin{array}{l}
\max\limits_{1\leq i\leq n}\left\vert c_{i}\right\vert \left(
\sum_{i=1}^{n}\left\vert c_{i}\right\vert ^{p}\right) ^{\frac{1}{p}%
}\max\limits_{1\leq i\leq n}\left( \sum_{j=1}^{n}\left\vert \left(
y_{i},y_{j}\right) \right\vert \right) ^{\frac{1}{p}}\left(
\sum_{i,j=1}^{n}\left\vert \left( y_{i},y_{j}\right) \right\vert \right) ^{%
\frac{1}{q}}; \\ 
\\ 
\left( \sum_{i=1}^{n}\left\vert c_{i}\right\vert ^{p}\right) ^{\frac{1}{p}%
}\left( \sum_{i=1}^{n}\left\vert c_{i}\right\vert ^{\gamma q}\right) ^{\frac{%
1}{\gamma q}}\max\limits_{1\leq i\leq n}\left( \sum_{j=1}^{n}\left\vert
\left( y_{i},y_{j}\right) \right\vert \right) ^{\frac{1}{p}} \\ 
\ \ \ \ \ \ \times \left( \sum_{i=1}^{n}\left( \sum_{j=1}^{n}\left\vert
\left( y_{i},y_{j}\right) \right\vert \right) ^{\delta }\right) ^{\frac{1}{%
\delta q}},\ \hfill \ \ \text{if}\ \gamma >1,\ \frac{1}{\gamma }+\frac{1}{%
\delta }=1; \\ 
\\ 
\left( \sum_{i=1}^{n}\left\vert c_{i}\right\vert ^{p}\right) ^{\frac{1}{p}%
}\left( \sum_{i=1}^{n}\left\vert c_{i}\right\vert ^{q}\right) ^{\frac{1}{q}%
}\max\limits_{1\leq i\leq n}\left( \sum_{j=1}^{n}\left\vert \left(
y_{i},y_{j}\right) \right\vert \right) ;%
\end{array}%
\right.
\end{equation*}%
for $p>1,$ $\frac{1}{p}+\frac{1}{q}=1.$
\end{theorem}

\begin{proof}
We note that 
\begin{equation*}
\sum\limits_{i=1}^{n}c_{i}\left( x,y_{i}\right) =\left(
x,\sum\limits_{i=1}^{n}\overline{c_{i}}y_{i}\right) .
\end{equation*}%
Using Schwarz's inequality in inner product spaces, we have 
\begin{equation*}
\left\vert \sum\limits_{i=1}^{n}c_{i}\left( x,y_{i}\right) \right\vert
^{2}\leq \left\Vert x\right\Vert ^{2}\left\Vert \sum\limits_{i=1}^{n}%
\overline{c_{i}}y_{i}\right\Vert ^{2}.
\end{equation*}%
Finally, using Lemma \ref{l2.1.14} with $\alpha _{i}=\overline{c_{i}},$ $%
z_{i}=y_{i}$ $\left( i=1,\dots ,n\right) ,$ we deduce the desired inequality
(\ref{3.1.14}).
\end{proof}

\begin{remark}
\label{r3.2.14}If in (\ref{3.1.14}) we choose $p=q=2,$ we obtain amongst
others, the result $\left( \ref{1.1.14}\right) $ due to J. Pe\v{c}ari\'{c}.
\end{remark}

\subsection{More Results of Bombieri Type}

The following results of Bombieri type hold \cite{15NSSD}.

\begin{theorem}
\label{t4.1.14}Let $x,y_{1},\dots ,y_{n}\in H.$ Then one has the inequality: 
\begin{align}
& \sum\limits_{i=1}^{n}\left\vert \left( x,y_{i}\right) \right\vert ^{2}
\label{4.1.14} \\
& \leq \left\Vert x\right\Vert \left[ \sum_{i=1}^{n}\left\vert \left(
x,y_{i}\right) \right\vert ^{p}\left( \sum\limits_{j=1}^{n}\left\vert \left(
y_{i},y_{j}\right) \right\vert \right) \right] ^{\frac{1}{2p}}\left[
\sum\limits_{i=1}^{n}\left\vert \left( x,y_{i}\right) \right\vert ^{q}\left(
\sum\limits_{j=1}^{n}\left\vert \left( y_{i},y_{j}\right) \right\vert
\right) \right] ^{\frac{1}{2q}}  \notag \\
& \leq \left\Vert x\right\Vert \times \left\{ 
\begin{array}{c}
J \\ 
K \\ 
L%
\end{array}%
\right. ,  \notag
\end{align}%
where%
\begin{equation*}
J:=\left\{ 
\begin{array}{l}
\max\limits_{1\leq i\leq n}\left\vert \left( x,y_{i}\right) \right\vert
\left( \sum_{i,j=1}^{n}\left\vert \left( y_{i},y_{j}\right) \right\vert
\right) ^{\frac{1}{2}}; \\ 
\\ 
\max\limits_{1\leq i\leq n}\left\vert \left( x,y_{i}\right) \right\vert ^{%
\frac{1}{2}}\left( \sum_{i=1}^{n}\left\vert \left( x,y_{i}\right)
\right\vert ^{\gamma q}\right) ^{\frac{1}{2\gamma q}}\left(
\sum_{i,j=1}^{n}\left\vert \left( y_{i},y_{j}\right) \right\vert \right) ^{%
\frac{1}{2p}} \\ 
\ \ \ \ \ \ \times \left( \sum_{i=1}^{n}\left( \sum_{j=1}^{n}\left\vert
\left( y_{i},y_{j}\right) \right\vert \right) ^{\delta }\right) ^{\frac{1}{%
2\delta q}},\hfill \ \ \text{if}\ \gamma >1,\ \frac{1}{\gamma }+\frac{1}{%
\delta }=1; \\ 
\\ 
\max\limits_{1\leq i\leq n}\left\vert \left( x,y_{i}\right) \right\vert ^{%
\frac{1}{2}}\left( \sum_{i=1}^{n}\left\vert \left( x,y_{i}\right)
\right\vert ^{q}\right) ^{\frac{1}{2q}}\left( \sum_{i,j=1}^{n}\left\vert
\left( y_{i},y_{j}\right) \right\vert \right) ^{\frac{1}{2p}%
}\max\limits_{1\leq i\leq n}\left( \sum_{j=1}^{n}\left\vert \left(
y_{i},y_{j}\right) \right\vert \right) ^{\frac{1}{2q}};%
\end{array}%
\right.
\end{equation*}%
\begin{equation*}
K:=\left\{ 
\begin{array}{l}
\max\limits_{1\leq i\leq n}\left\vert \left( x,y_{i}\right) \right\vert ^{%
\frac{1}{2}}\left( \sum_{i=1}^{n}\left\vert \left( x,y_{i}\right)
\right\vert ^{\alpha p}\right) ^{\frac{1}{2\alpha \beta }}\left(
\sum_{i,j=1}^{n}\left\vert \left( y_{i},y_{j}\right) \right\vert \right) ^{%
\frac{1}{2q}} \\ 
\ \ \ \ \ \ \times \left( \sum_{i=1}^{n}\left( \sum_{j=1}^{n}\left\vert
\left( y_{i},y_{j}\right) \right\vert \right) ^{\beta }\right) ^{\frac{1}{%
p\beta }},\hfill \ \ \text{if}\ \alpha >1,\ \frac{1}{\alpha }+\frac{1}{\beta 
}=1; \\ 
\\ 
\left( \sum_{i=1}^{n}\left\vert \left( x,y_{i}\right) \right\vert ^{\alpha
p}\right) ^{\frac{1}{2\alpha p}}\left( \sum\limits_{i=1}^{n}\left\vert
\left( x,y_{i}\right) \right\vert ^{\gamma q}\right) ^{\frac{1}{2\gamma q}%
}\left( \sum_{i=1}^{n}\left( \sum_{j=1}^{n}\left\vert \left(
y_{i},y_{j}\right) \right\vert \right) ^{\beta }\right) ^{\frac{1}{2p\beta }}
\\ 
\times \left( \sum_{i=1}^{n}\left( \sum_{j=1}^{n}\left\vert \left(
y_{i},y_{j}\right) \right\vert \right) ^{\delta }\right) ^{\frac{1}{2\delta q%
}}\hfill \ \text{if}\ \alpha >1,\ \frac{1}{\alpha }+\frac{1}{\beta }=1\text{ 
} \\ 
\hfill \text{and }\ \gamma >1,\ \frac{1}{\gamma }+\frac{1}{\delta }=1; \\ 
\\ 
\left( \sum_{i=1}^{n}\left\vert \left( x,y_{i}\right) \right\vert
^{q}\right) ^{\frac{1}{2q}}\left( \sum\limits_{i=1}^{n}\left\vert \left(
x,y_{i}\right) \right\vert ^{\alpha p}\right) ^{\frac{1}{2\alpha p}%
}\max\limits_{1\leq i\leq n}\left( \sum_{j=1}^{n}\left\vert \left(
y_{i},y_{j}\right) \right\vert \right) ^{\frac{1}{2p}} \\ 
\ \ \ \ \ \ \ \times \left( \sum_{i=1}^{n}\left( \sum_{j=1}^{n}\left\vert
\left( y_{i},y_{j}\right) \right\vert \right) ^{\beta }\right) ^{\frac{1}{%
2p\beta }},\hfill \ \ \text{if}\ \alpha >1,\ \frac{1}{\alpha }+\frac{1}{%
\beta }=1;%
\end{array}%
\right.
\end{equation*}%
and%
\begin{equation*}
L:=\left\{ 
\begin{array}{l}
\max\limits_{1\leq i\leq n}\left\vert \left( x,y_{i}\right) \right\vert ^{%
\frac{1}{2}}\left( \sum_{i=1}^{n}\left\vert \left( x,y_{i}\right)
\right\vert ^{p}\right) ^{\frac{1}{2p}}\max\limits_{1\leq i\leq n}\left(
\sum_{j=1}^{n}\left\vert \left( y_{i},y_{j}\right) \right\vert \right) ^{%
\frac{1}{2p}} \\ 
\qquad \times \left( \sum_{i,j=1}^{n}\left\vert \left( y_{i},y_{j}\right)
\right\vert \right) ^{\frac{1}{2q}}; \\ 
\\ 
\left( \sum_{i=1}^{n}\left\vert \left( x,y_{i}\right) \right\vert
^{p}\right) ^{\frac{1}{2p}}\left( \sum_{i=1}^{n}\left\vert \left(
x,y_{i}\right) \right\vert ^{\gamma q}\right) ^{\frac{1}{2\gamma q}%
}\max\limits_{1\leq i\leq n}\left( \sum_{j=1}^{n}\left\vert \left(
y_{i},y_{j}\right) \right\vert \right) ^{\frac{1}{2p}} \\ 
\ \ \ \ \ \ \ \times \left( \sum_{i=1}^{n}\left( \sum_{j=1}^{n}\left\vert
\left( y_{i},y_{j}\right) \right\vert \right) ^{\delta }\right) ^{\frac{1}{%
2\delta q}},\hfill \ \ \text{if}\ \gamma >1,\ \frac{1}{\gamma }+\frac{1}{%
\delta }=1; \\ 
\\ 
\left( \sum_{i=1}^{n}\left\vert \left( x,y_{i}\right) \right\vert
^{p}\right) ^{\frac{1}{2p}}\left( \sum_{i=1}^{n}\left\vert \left(
x,y_{i}\right) \right\vert ^{q}\right) ^{\frac{1}{2q}}\max\limits_{1\leq
i\leq n}\left( \sum_{j=1}^{n}\left\vert \left( y_{i},y_{j}\right)
\right\vert \right) ^{\frac{1}{2}},%
\end{array}%
\right.
\end{equation*}%
for $p>1,\frac{1}{p}+\frac{1}{q}=1.$
\end{theorem}

\begin{proof}
The proof follows by Theorem \ref{t3.1.14} on choosing $c_{i}=\overline{%
\left( x,y_{i}\right) },i\in \left\{ 1,\dots ,n\right\} $ and taking the
square root in both sides of the inequalities involved. We omit the details.
\end{proof}

\begin{remark}
We observe, by the last inequality in $\left( \ref{4.1.14}\right) $, that 
\begin{equation*}
\frac{\left( \sum\limits_{i=1}^{n}\left\vert \left( x,y_{i}\right)
\right\vert ^{2}\right) ^{2}}{\left( \sum\limits_{i=1}^{n}\left\vert \left(
x,y_{i}\right) \right\vert ^{p}\right) ^{\frac{1}{p}}\left(
\sum\limits_{i=1}^{n}\left\vert \left( x,y_{i}\right) \right\vert
^{q}\right) ^{\frac{1}{q}}}\leq \left\Vert x\right\Vert
^{2}\max\limits_{1\leq i\leq n}\left( \sum\limits_{j=1}^{n}\left\vert \left(
y_{i},y_{j}\right) \right\vert \right) ,
\end{equation*}%
where $p>1,\frac{1}{p}+\frac{1}{q}=1.$

If in this inequality we choose $p=q=2,$ then we recapture Bombieri's result 
$\left( \ref{1.3.14}\right) .$
\end{remark}

\newpage

\chapter{Some Gr\"{u}ss' Type Inequalities for $n$-Tuples of Vectors}

\section{Introduction}

We start by recalling some of the most important Gr\"{u}ss type discrete
inequalities for real numbers that are available in the literature.

\begin{enumerate}
\item (1950) \textit{Biernacki, Pidek, Ryll-Nardzewski} \cite{1b.15.0}.

If $\mathbf{\bar{a}}=\left( a_{1},\dots ,a_{n}\right) $ and $\mathbf{\bar{b}}%
=\left( b_{1},\dots ,b_{n}\right) $ are $n$-tuples of real numbers such that
there exists the real numbers $a,A,b,B$ \hspace{0.05in}with 
\begin{equation}
a\leq a_{i}\leq A,\;\;b\leq b_{i}\leq B,\;\;i\in \left\{ 1,\dots ,n\right\} ,
\label{b.1.15.0}
\end{equation}%
then 
\begin{eqnarray*}
\left\vert C_{n}\left( \mathbf{\bar{a}},\mathbf{\bar{b}}\right) \right\vert
&\leq &\frac{1}{n}\left[ \frac{n}{2}\right] \left( 1-\frac{1}{n}\left[ \frac{%
n}{2}\right] \right) \left( A-a\right) \left( B-b\right) \\
&=&\frac{1}{n^{2}}\left[ \frac{n^{4}}{4}\right] \left( A-a\right) \left(
B-b\right) \\
&\leq &\frac{1}{4}\left( A-a\right) \left( B-b\right) .
\end{eqnarray*}%
\bigskip

\item[2.] (1988) \textit{Andrica-Badea} \cite{2b.15.0}.

Let $\mathbf{\bar{a}}$, $\mathbf{\bar{b}}$ satisfy (\ref{b.1.15.0}) and $%
\mathbf{\bar{p}}=\left( p_{1},\dots ,p_{n}\right) $ be an $n-$tuple of
nonnegative numbers with $P_{n}>0.$ If $S$ is a subset of $\left\{ 1,\dots
,n\right\} $ that minimises the expression 
\begin{equation*}
\left\vert \sum_{i\in S}q_{i}-\frac{1}{2}Q_{n}\right\vert ,
\end{equation*}%
then 
\begin{align*}
C_{n}\left( \mathbf{\bar{p}},\mathbf{\bar{a}},\mathbf{\bar{b}}\right) & \leq 
\frac{Q_{S}}{Q_{n}}\left( 1-\frac{Q_{S}}{Q_{n}}\right) \left( A-a\right)
\left( B-b\right) \\
& \leq \frac{1}{4}\left( A-a\right) \left( B-b\right) ,
\end{align*}%
where $Q_{S}:=\sum_{i\in S}Q_{i}.$

\item[3.] (2000) \textit{Dragomir-Booth} \cite{3b.15.0}.

If $\mathbf{\bar{a}}$, $\mathbf{\bar{b}}$ are real $n-$tuples and $\mathbf{%
\bar{p}}$ is nonnegative with $P_{n}>0,$ then 
\begin{equation*}
\left\vert C_{n}\left( \mathbf{\bar{p}},\mathbf{\bar{a}},\mathbf{\bar{b}}%
\right) \right\vert \leq \max\limits_{1\leq j\leq n-1}\left\vert \Delta
a_{j}\right\vert \max\limits_{1\leq j\leq n-1}\left\vert \Delta
b_{j}\right\vert C_{n}\left( \mathbf{\bar{p}},\mathbf{\bar{e}},\mathbf{\bar{e%
}}\right) ,
\end{equation*}%
where $\mathbf{\bar{e}}=\left( 1,2,\dots ,n\right) $ and $\Delta
a_{j}:=a_{j+1}-a_{j}$ is the forward difference, $j=1,\dots ,n-1.$ Note that 
\begin{equation*}
C_{n}\left( \mathbf{\bar{p}},\mathbf{\bar{e}},\mathbf{\bar{e}}\right) =\frac{%
1}{P_{n}^{2}}\sum_{i=1}^{n}i^{2}p_{i}-\left( \frac{1}{P_{n}}%
\sum_{i=1}^{n}ip_{i}\right) ^{2}.
\end{equation*}%
In particular, we have 
\begin{equation*}
\left\vert C_{n}\left( \mathbf{\bar{a}},\mathbf{\bar{b}}\right) \right\vert
\leq \frac{1}{12}\left( n^{2}-1\right) \max\limits_{1\leq j\leq
n-1}\left\vert \Delta a_{j}\right\vert \max\limits_{1\leq j\leq
n-1}\left\vert \Delta b_{j}\right\vert .
\end{equation*}%
The constant $\frac{1}{12}$ is best possible.\bigskip

\item[4.] (2002) \textit{Dragomir }\cite{4b.15.0}.

With the assumptions in 3, the following inequality holds 
\begin{equation*}
\left\vert C_{n}\left( \mathbf{\bar{p}},\mathbf{\bar{a}},\mathbf{\bar{b}}%
\right) \right\vert \leq \frac{1}{P_{n}^{2}}\sum_{1\leq j<i\leq n}\left(
i-j\right) \left( \sum_{k=1}^{n-1}\left\vert \Delta a_{k}\right\vert
^{p}\right) ^{\frac{1}{p}}\left( \sum_{k=1}^{n-1}\left\vert \Delta
b_{k}\right\vert ^{q}\right) ^{\frac{1}{q}},
\end{equation*}%
where $p>1,$ $\frac{1}{p}+\frac{1}{q}=1.$

In particular, we have 
\begin{equation*}
\left\vert C_{n}\left( \mathbf{\bar{a}},\mathbf{\bar{b}}\right) \right\vert
\leq \frac{1}{6}\cdot \frac{n^{2}-1}{n}\left( \sum_{k=1}^{n-1}\left\vert
\Delta a_{k}\right\vert ^{p}\right) ^{\frac{1}{p}}\left(
\sum_{k=1}^{n-1}\left\vert \Delta b_{k}\right\vert ^{q}\right) ^{\frac{1}{q}%
}.
\end{equation*}%
The constant $\frac{1}{6}$ is best possible.\bigskip

\item[5.] (2002) \textit{Dragomir} \cite{5b.15.0}.

The following inequality holds, where $\mathbf{\bar{p}},\mathbf{\bar{a}},%
\mathbf{\bar{b}}$ and $P_{n}$ are as in assumption 3, 
\begin{equation*}
\left\vert C_{n}\left( \mathbf{\bar{p}},\mathbf{\bar{a}},\mathbf{\bar{b}}%
\right) \right\vert \leq \frac{1}{2}\cdot \frac{1}{P_{n}^{2}}%
\sum_{i=1}^{n}p_{i}\left( P_{n}-p_{i}\right) \sum_{k=1}^{n-1}\left\vert
\Delta a_{k}\right\vert \sum_{k=1}^{n-1}\left\vert \Delta b_{k}\right\vert .
\end{equation*}%
In particular, we have 
\begin{equation*}
\left\vert C_{n}\left( \mathbf{\bar{a}},\mathbf{\bar{b}}\right) \right\vert
\leq \frac{1}{2}\left( 1-\frac{1}{n}\right) \sum_{k=1}^{n-1}\left\vert
\Delta a_{k}\right\vert \sum_{k=1}^{n-1}\left\vert \Delta b_{k}\right\vert .
\end{equation*}%
The constant $\frac{1}{2}$ is sharp.\bigskip

\item[6.] (2002) \textit{Cerone-Dragomir} \cite{6b.15.0}.

If $\mathbf{\bar{a}}$, $\mathbf{\bar{b}}$ are real $n-$tuples and $\mathbf{%
\bar{p}}$ is a positive $n-$tuple and there exists $m,M\in \mathbb{R}$ such
that 
\begin{equation*}
m\leq a_{i}\leq M,
\end{equation*}%
then one has the inequality 
\begin{equation*}
\left\vert C_{n}\left( \mathbf{\bar{p}},\mathbf{\bar{a}},\mathbf{\bar{b}}%
\right) \right\vert \leq \frac{1}{2}\left( M-m\right) \frac{1}{P_{n}}%
\sum_{i=1}^{n}p_{i}\left\vert b_{i}-\frac{1}{P_{n}}\sum_{j=1}^{n}p_{j}b_{j}%
\right\vert .
\end{equation*}%
The constant $\frac{1}{2}$\hspace{0.05in}is best possible.

In particular, we have 
\begin{equation*}
\left\vert C_{n}\left( \mathbf{\bar{a}},\mathbf{\bar{b}}\right) \right\vert
\leq \frac{1}{2}\left( M-m\right) \cdot \frac{1}{n}\sum_{i=1}^{n}\left\vert
b_{i}-\frac{1}{n}\sum_{j=1}^{n}b_{j}\right\vert .
\end{equation*}%
The constant $\frac{1}{2}$ is best possible.
\end{enumerate}

The main aim of this chapter is to present some extensions of the above
results holding in the general setting of $n$-tuples of vectors in an inner
product space.

\newpage

\section{The Version for Norms}

\subsection{Preliminary Results}

The following lemma is of interest in itself \cite{16NSSD}.

\begin{lemma}
\label{l2.1.15}Let $\left( H;\left\langle \cdot ,\cdot \right\rangle \right) 
$ be an inner product space over the real or complex number field $\mathbb{K}%
,$ $x_{i}\in H$ and $p_{i}\geq 0$ $\left( i=1,\dots ,n\right) $ such that $%
\sum_{i=1}^{n}p_{i}=1$ $\left( n\geq 2\right) .$\ If $x,X\in H$ are such
that 
\begin{equation}
\func{Re}\left\langle X-x_{i},x_{i}-x\right\rangle \geq 0\text{ for all }%
i\in \left\{ 1,\dots ,n\right\} ,  \label{2.1.15}
\end{equation}%
or equivalently,%
\begin{equation*}
\left\Vert x_{i}-\frac{x+X}{2}\right\Vert \leq \frac{1}{2}\left\Vert
X-x\right\Vert \text{ for all }i\in \left\{ 1,\dots ,n\right\} ,
\end{equation*}%
then we have the inequality 
\begin{equation}
0\leq \sum_{i=1}^{n}p_{i}\left\Vert x_{i}\right\Vert ^{2}-\left\Vert
\sum_{i=1}^{n}p_{i}x_{i}\right\Vert ^{2}\leq \frac{1}{4}\left\Vert
X-x\right\Vert ^{2}.  \label{2.2.15}
\end{equation}%
The constant $\frac{1}{4}$ is sharp.
\end{lemma}

\begin{proof}
Define 
\begin{equation*}
I_{1}:=\left\langle
X-\sum_{i=1}^{n}p_{i}x_{i},\sum_{i=1}^{n}p_{i}x_{i}-x\right\rangle
\end{equation*}%
and 
\begin{equation*}
I_{2}:=\sum_{i=1}^{n}p_{i}\left\langle X-x_{i},x_{i}-x\right\rangle .
\end{equation*}%
Then 
\begin{equation*}
I_{1}=\sum_{i=1}^{n}p_{i}\left\langle X,x_{i}\right\rangle -\left\langle
X,x\right\rangle -\left\Vert \sum_{i=1}^{n}p_{i}x_{i}\right\Vert
^{2}+\sum_{i=1}^{n}p_{i}\left\langle x_{i},x\right\rangle
\end{equation*}%
and 
\begin{equation*}
I_{2}=\sum_{i=1}^{n}p_{i}\left\langle X,x_{i}\right\rangle -\left\langle
X,x\right\rangle -\sum_{i=1}^{n}p_{i}\left\Vert x_{i}\right\Vert
^{2}+\sum_{i=1}^{n}p_{i}\left\langle x_{i},x\right\rangle .
\end{equation*}%
Consequently, 
\begin{equation}
I_{1}-I_{2}=\sum_{i=1}^{n}p_{i}\left\Vert x_{i}\right\Vert ^{2}-\left\Vert
\sum_{i=1}^{n}p_{i}x_{i}\right\Vert ^{2}.  \label{2.3.15}
\end{equation}%
Taking the real value in $\left( \ref{2.3.15}\right) ,$ we can state that 
\begin{multline}
\sum_{i=1}^{n}p_{i}\left\Vert x_{i}\right\Vert ^{2}-\left\Vert
\sum_{i=1}^{n}p_{i}x_{i}\right\Vert ^{2}  \label{2.4.15} \\
=\func{Re}\left\langle
X-\sum_{i=1}^{n}p_{i}x_{i},\sum_{i=1}^{n}p_{i}x_{i}-x\right\rangle
-\sum_{i=1}^{n}p_{i}\func{Re}\left\langle X-x_{i},x_{i}-x\right\rangle ,
\end{multline}%
which is also an identity of interest in itself.\newline
Using the assumption $\left( \ref{2.1.15}\right) ,$ we can conclude, by $%
\left( \ref{2.4.15}\right) ,$ that 
\begin{equation}
\sum_{i=1}^{n}p_{i}\left\Vert x_{i}\right\Vert ^{2}-\left\Vert
\sum_{i=1}^{n}p_{i}x_{i}\right\Vert ^{2}\leq \func{Re}\left\langle
X-\sum_{i=1}^{n}p_{i}x_{i},\sum_{i=1}^{n}p_{i}x_{i}-x\right\rangle .
\label{2.5.15}
\end{equation}%
It is known that if $y,z\in H,$ then 
\begin{equation}
4\func{Re}\left\langle z,y\right\rangle \leq \left\Vert z+y\right\Vert ^{2},
\label{2.6.15}
\end{equation}%
with equality iff $z=y.$

Now, by $\left( \ref{2.6.15}\right) $ we can state that 
\begin{align*}
\func{Re}\left\langle
X-\sum_{i=1}^{n}p_{i}x_{i},\sum_{i=1}^{n}p_{i}x_{i}-x\right\rangle & \leq 
\frac{1}{4}\left\Vert
X-\sum_{i=1}^{n}p_{i}x_{i}+\sum_{i=1}^{n}p_{i}x_{i}-x\right\Vert ^{2} \\
& =\frac{1}{4}\left\Vert X-x\right\Vert ^{2}.
\end{align*}%
Using $\left( \ref{2.5.15}\right) ,$ we obtain $\left( \ref{2.2.15}\right) .$

To prove the sharpness of the constant $\frac{1}{4},$ let us assume that the
inequality $\left( \ref{2.2.15}\right) $ holds with a constant $c>0,$ i.e., 
\begin{equation}
0\leq \sum_{i=1}^{n}p_{i}\left\Vert x_{i}\right\Vert ^{2}-\left\Vert
\sum_{i=1}^{n}p_{i}x_{i}\right\Vert ^{2}\leq c\left\Vert X-x\right\Vert ^{2}
\label{2.7.15}
\end{equation}%
for all $p_{i},x_{i}$ and $n$ as in the hypothesis of Lemma \ref{l2.1.15}.

Assume that $n=2,$ $p_{1}=p_{2}=\frac{1}{2},$ $x_{1}=x$ and $x_{2}=X$ with $%
x,X\in H$ and $x\neq X.$ Then, obviously, 
\begin{equation*}
\left\langle X-x_{1},x_{1}-x\right\rangle =\left\langle
X-x_{2},x_{2}-x\right\rangle =0,
\end{equation*}%
which shows that the condition $\left( \ref{2.1.15}\right) $ holds.

If we replace $n,p_{1},p_{2},x_{1},x_{2}$ in $\left( \ref{2.7.15}\right) ,$
we obtain 
\begin{align*}
\sum_{i=1}^{2}p_{i}\left\Vert x_{i}\right\Vert ^{2}-\left\Vert
\sum_{i=1}^{2}p_{i}x_{i}\right\Vert ^{2}& =\frac{1}{2}\left( \left\Vert
x\right\Vert ^{2}+\left\Vert X\right\Vert ^{2}\right) -\left\Vert \frac{x+X}{%
2}\right\Vert ^{2} \\
& =\frac{1}{4}\left\Vert x-X\right\Vert ^{2} \\
& \leq c\left\Vert x-X\right\Vert ^{2},
\end{align*}%
from where we deduce that $c\geq \frac{1}{4},$ which proves the sharpness of
the constant $\frac{1}{4}.$
\end{proof}

\begin{remark}
\label{r2.2.15}The assumption $\left( \ref{2.1.15}\right) $ can be replaced
by the more general condition 
\begin{equation}
\sum_{i=1}^{n}p_{i}\func{Re}\left\langle X-x_{i},x_{i}-x\right\rangle \geq 0,
\label{2.8.15}
\end{equation}%
and the conclusion $\left( \ref{2.2.15}\right) $ will still remain valid.
\end{remark}

The following corollary is natural.

\begin{corollary}
\label{c2.3.15}Let $a_{i}\in \mathbb{K},$ $p_{i}\geq 0,$ $\left( i=1,\dots
,n\right) $ $\left( n\geq 2\right) $ with $\sum_{i=1}^{n}p_{i}=1.$ If $%
a,A\in \mathbb{K}$ are such that 
\begin{equation}
\func{Re}\left[ \left( A-a_{i}\right) \left( \bar{a}_{i}-\bar{a}\right) %
\right] \geq 0\text{ for all }i\in \left\{ 1,\dots ,n\right\} ,
\label{2.9.15}
\end{equation}%
then we have the inequality 
\begin{equation}
0\leq \sum_{i=1}^{n}p_{i}\left\vert a_{i}\right\vert ^{2}-\left\vert
\sum_{i=1}^{n}p_{i}a_{i}\right\vert ^{2}\leq \frac{1}{4}\left\vert
A-a\right\vert ^{2}.  \label{2.10.15}
\end{equation}%
The constant $\frac{1}{4}$ is sharp.
\end{corollary}

The proof follows by the above lemma by choosing $H=\mathbb{K},$ $%
\left\langle x,y\right\rangle :=x\bar{y},$ $x_{i}=a_{i},$ $x=a$ and $X=A.$
We omit the details.

\begin{remark}
\label{r2.4.15}The condition $\left( \ref{2.9.15}\right) $ can be replaced
by the more general assumption 
\begin{equation*}
\sum_{i=1}^{n}p_{i}\func{Re}\left[ \left( A-a_{i}\right) \left( \bar{a}_{i}-%
\bar{a}\right) \right] \geq 0.
\end{equation*}
\end{remark}

\subsection{A Discrete Inequality of Gr\"{u}ss' Type}

The following Gr\"{u}ss type inequality holds \cite{16NSSD}.

\begin{theorem}
\label{t3.1.15}Let $\left( H;\left\langle \cdot ,\cdot \right\rangle \right) 
$ be an inner product space over $\mathbb{K},$ $\mathbb{K}=\mathbb{R},%
\mathbb{C},\;x_{i}\in H,\;a_{i}\in \mathbb{K},\;p_{i}\geq 0$ $\left(
i=1,\dots ,n\right) $ $\left( n\geq 2\right) $ with $\sum_{i=1}^{n}p_{i}=1.$
If $a,A\in \mathbb{K}$ and $x,X\in H$ are such that 
\begin{equation}
\func{Re}\left[ \left( A-a_{i}\right) \left( \bar{a}_{i}-\bar{a}\right) %
\right] \geq 0,\text{ }\func{Re}\left\langle X-x_{i},x_{i}-x\right\rangle
\geq 0\text{ }  \label{3.1.15a}
\end{equation}%
for all $i\in \left\{ 1,\dots ,n\right\} ;$ then we have the inequality 
\begin{equation}
0\leq \left\Vert \sum_{i=1}^{n}p_{i}a_{i}x_{i}-\sum_{i=1}^{n}p_{i}a_{i}\cdot
\sum_{i=1}^{n}p_{i}x_{i}\right\Vert \leq \frac{1}{4}\left\vert
A-a\right\vert \left\Vert X-x\right\Vert .  \label{3.2.15}
\end{equation}%
The constant $\frac{1}{4}$ is sharp.
\end{theorem}

\begin{proof}
A simple computation shows that 
\begin{equation}
\sum_{i=1}^{n}p_{i}a_{i}x_{i}-\sum_{i=1}^{n}p_{i}a_{i}%
\sum_{i=1}^{n}p_{i}x_{i}=\frac{1}{2}\sum_{i,j=1}^{n}p_{i}p_{j}\left(
a_{i}-a_{j}\right) \left( x_{i}-x_{j}\right) .  \label{3.3.15}
\end{equation}%
Taking the norm in both parts of $\left( \ref{3.3.15}\right) $ and using the
generalized triangle inequality, we obtain 
\begin{equation}
\left\Vert
\sum_{i=1}^{n}p_{i}a_{i}x_{i}-\sum_{i=1}^{n}p_{i}a_{i}%
\sum_{i=1}^{n}p_{i}x_{i}\right\Vert \leq \frac{1}{2}%
\sum_{i,j=1}^{n}p_{i}p_{j}\left\vert a_{i}-a_{j}\right\vert \left\Vert
x_{i}-x_{j}\right\Vert .  \label{3.4.15}
\end{equation}%
By the Cauchy-Bunyakovsky-Schwartz discrete inequality for double sums, we
obtain 
\begin{multline}
\left( \frac{1}{2}\sum_{i,j=1}^{n}p_{i}p_{j}\left\vert
a_{i}-a_{j}\right\vert \left\Vert x_{i}-x_{j}\right\Vert \right) ^{2}
\label{3.5.15} \\
\leq \left( \frac{1}{2}\sum_{i,j=1}^{n}p_{i}p_{j}\left\vert
a_{i}-a_{j}\right\vert ^{2}\right) \left( \frac{1}{2}%
\sum_{i,j=1}^{n}p_{i}p_{j}\left\Vert x_{i}-x_{j}\right\Vert ^{2}\right) .
\end{multline}%
As a simple calculation reveals that 
\begin{equation*}
\frac{1}{2}\sum_{i,j=1}^{n}p_{i}p_{j}\left\vert a_{i}-a_{j}\right\vert
^{2}=\sum_{i=1}^{n}p_{i}\left\vert a_{i}\right\vert ^{2}-\left\vert
\sum_{i=1}^{n}p_{i}a_{i}\right\vert ^{2}
\end{equation*}%
and 
\begin{equation*}
\frac{1}{2}\sum_{i,j=1}^{n}p_{i}p_{j}\left\Vert x_{i}-x_{j}\right\Vert
^{2}=\sum_{i=1}^{n}p_{i}\left\Vert x_{i}\right\Vert ^{2}-\left\Vert
\sum_{i=1}^{n}p_{i}x_{i}\right\Vert ^{2},
\end{equation*}%
then, by $\left( \ref{3.4.15}\right) $ and $\left( \ref{3.5.15}\right) ,$ we
conclude that 
\begin{multline}
\left\Vert
\sum_{i=1}^{n}p_{i}a_{i}x_{i}-\sum_{i=1}^{n}p_{i}a_{i}%
\sum_{i=1}^{n}p_{i}x_{i}\right\Vert  \label{3.6.15} \\
\leq \left( \sum_{i=1}^{n}p_{i}\left\vert a_{i}\right\vert ^{2}-\left\vert
\sum_{i=1}^{n}p_{i}a_{i}\right\vert ^{2}\right) ^{\frac{1}{2}}\left(
\sum_{i=1}^{n}p_{i}\left\Vert x_{i}\right\Vert ^{2}-\left\Vert
\sum_{i=1}^{n}p_{i}x_{i}\right\Vert ^{2}\right) ^{\frac{1}{2}}.
\end{multline}%
However, from Lemma \ref{l2.1.15} and Corollary \ref{c2.3.15}, we know that 
\begin{equation}
\left( \sum_{i=1}^{n}p_{i}\left\Vert x_{i}\right\Vert ^{2}-\left\Vert
\sum_{i=1}^{n}p_{i}x_{i}\right\Vert ^{2}\right) ^{\frac{1}{2}}\leq \frac{1}{2%
}\left\Vert X-x\right\Vert  \label{3.7.15}
\end{equation}%
and 
\begin{equation}
\left( \sum_{i=1}^{n}p_{i}\left\vert a_{i}\right\vert ^{2}-\left\vert
\sum_{i=1}^{n}p_{i}a_{i}\right\vert ^{2}\right) ^{\frac{1}{2}}\leq \frac{1}{2%
}\left\vert A-a\right\vert .  \label{3.8.15}
\end{equation}%
Consequently, by using $\left( \ref{3.6.15}\right) -\left( \ref{3.8.15}%
\right) ,$ we deduce the desired estimate $\left( \ref{3.2.15}\right) .$

To prove the sharpness of the constant $\frac{1}{4},$ assume that $\left( %
\ref{3.2.15}\right) $ holds with a constant $c>0,$ i.e., 
\begin{equation}
\left\Vert
\sum_{i=1}^{n}p_{i}a_{i}x_{i}-\sum_{i=1}^{n}p_{i}a_{i}%
\sum_{i=1}^{n}p_{i}x_{i}\right\Vert \leq c\left\vert A-a\right\vert
\left\Vert X-x\right\Vert  \label{3.9.15}
\end{equation}%
for all $p_{i},a_{i},x_{i},a,A,x,X$ and $n$ as in the hypothesis of Theorem %
\ref{t3.1.15}.

If we choose $n=2,$ $a_{1}=a,$ $a_{2}=A,$ $x_{1}=x,$ $x_{2}=X$ $\left( a\neq
A,x\neq X\right) $ and $p_{1}=p_{2}=\frac{1}{2},$ then 
\begin{align*}
\sum_{i=1}^{2}p_{i}a_{i}x_{i}-\sum_{i=1}^{2}p_{i}a_{i}%
\sum_{i=1}^{2}p_{i}x_{i}& =\frac{1}{2}\sum_{i,j=1}^{2}p_{i}p_{j}\left(
a_{i}-a_{j}\right) \left( x_{i}-x_{j}\right) \\
& =\frac{1}{4}\left( a-A\right) \left( x-X\right) .
\end{align*}%
Consequently, from $\left( \ref{3.9.15}\right) ,$ we deduce 
\begin{equation*}
\frac{1}{4}\left\vert a-A\right\vert \left\Vert X-x\right\Vert \leq
c\left\vert A-a\right\vert \left\Vert X-x\right\Vert ,
\end{equation*}%
which implies that $c\geq \frac{1}{4},$ and the theorem is completely proved.
\end{proof}

\begin{remark}
\label{r3.2.15}The condition $\left( \ref{3.1.15a}\right) $ can be replaced
by the more general assumption 
\begin{equation*}
\sum_{i=1}^{n}p_{i}\func{Re}\left[ \left( A-a_{i}\right) \left( \bar{a}_{i}-%
\bar{a}\right) \right] \geq 0\text{, }\sum_{i=1}^{n}p_{i}\func{Re}%
\left\langle X-x_{i},x_{i}-x\right\rangle \geq 0
\end{equation*}%
and the conclusion $\left( \ref{3.2.15}\right) $ will still be valid.
\end{remark}

The following corollary for real or complex numbers holds.

\begin{corollary}
\label{c3.3.15}Let $a_{i},b_{i}\in \mathbb{K}$ $\left( \mathbb{K}=\mathbb{C},%
\mathbb{R}\right) ,$ $p_{i}\geq 0$ $\left( i=1,\dots ,n\right) $ with $%
\sum_{i=1}^{n}p_{i}=1.$ If $a,A,b,B\in \mathbb{K}$ are such that 
\begin{equation}
\func{Re}\left[ \left( A-a_{i}\right) \left( \bar{a}_{i}-\bar{a}\right) %
\right] \geq 0,\text{ \ }\func{Re}\left[ \left( B-b_{i}\right) \left( \bar{b}%
_{i}-\bar{b}\right) \right] \geq 0,  \label{3.11.15}
\end{equation}%
then we have the inequality 
\begin{equation}
0\leq \left\vert
\sum_{i=1}^{n}p_{i}a_{i}b_{i}-\sum_{i=1}^{n}p_{i}a_{i}%
\sum_{i=1}^{n}p_{i}b_{i}\right\vert \leq \frac{1}{4}\left\vert
A-a\right\vert \left\vert B-b\right\vert ,  \label{3.12.15}
\end{equation}%
where the constant $\frac{1}{4}$ is sharp.
\end{corollary}

\begin{remark}
\label{r3.5.15}If we assume that $a_{i},b_{i},a,A,b,B$ are real numbers,
then $\left( \ref{3.11.15}\right) $ is equivalent to 
\begin{equation*}
a\leq a_{i}\leq A,b\leq b_{i}\leq B\text{ for all }i\in \left\{ 1,\dots
,n\right\} ,
\end{equation*}%
and $\left( \ref{3.12.15}\right) $ becomes 
\begin{equation*}
0\leq \left\vert
\sum_{i=1}^{n}p_{i}a_{i}b_{i}-\sum_{i=1}^{n}p_{i}a_{i}%
\sum_{i=1}^{n}p_{i}b_{i}\right\vert \leq \frac{1}{4}\left( A-a\right) \left(
B-b\right) ,
\end{equation*}%
which is the classical Gr\"{u}ss inequality for sequences of real numbers.
\end{remark}

\subsection{Applications for Discrete Fourier Transforms}

Let $\left( H;\left\langle \cdot ,\cdot \right\rangle \right) $ be an inner
product space over $\mathbb{K}$ and $\mathbf{\bar{x}}=\left( x_{1},\dots
,x_{n}\right) $ be a sequence of vectors in $H.$

For a given $w\in \mathbb{R},$ define the \textit{discrete Fourier transform}
as 
\begin{equation}
\mathcal{F}_{w}\left( \mathbf{\bar{x}}\right) \left( m\right)
:=\sum_{k=1}^{n}\exp \left( 2wimk\right) \times x_{k},\;m=1,\dots ,n.
\label{4.1.15}
\end{equation}%
The following approximation result for the Fourier transform $\left( \ref%
{4.1.15}\right) $ holds \cite{16NSSD}.

\begin{theorem}
\label{t4.1.15}Let $\left( H;\left\langle \cdot ,\cdot \right\rangle \right) 
$ and $\mathbf{\bar{x}}\in H^{n}$ be as above. If there exists the vectors $%
x,X\in H$ such that 
\begin{equation}
\func{Re}\left\langle X-x_{k},x_{k}-x\right\rangle \geq 0\text{ for all }%
k\in \left\{ 1,\dots ,n\right\} ,  \label{4.2.15}
\end{equation}%
then we have the inequality 
\begin{multline}
\left\Vert \mathcal{F}_{w}\left( \mathbf{\bar{x}}\right) \left( m\right) -%
\frac{\sin \left( wmn\right) }{\sin \left( wm\right) }\exp \left[ w\left(
n+1\right) im\right] \times \frac{1}{n}\sum_{k=1}^{n}x_{k}\right\Vert
\label{4.3.15} \\
\leq \frac{1}{2}\left\Vert X-x\right\Vert \left[ n^{2}-\frac{\sin ^{2}\left(
wmn\right) }{\sin ^{2}\left( wm\right) }\right] ^{\frac{1}{2}},
\end{multline}%
for all $m\in \left\{ 1,\dots ,n\right\} $ and $w\in \mathbb{R},\;w\neq 
\frac{l}{m}\pi ,\;l\in \mathbb{Z}.$
\end{theorem}

\begin{proof}
From the inequality $\left( \ref{3.6.15}\right) $ in Theorem \ref{t3.1.15},
we can state that 
\begin{multline}
\left\Vert \frac{1}{n}\sum_{k=1}^{n}a_{k}x_{k}-\frac{1}{n}%
\sum_{k=1}^{n}a_{k}\cdot \frac{1}{n}\sum_{k=1}^{n}x_{k}\right\Vert
\label{4.4.15} \\
\leq \left( \frac{1}{n}\sum_{k=1}^{n}\left\vert a_{k}\right\vert
^{2}-\left\vert \frac{1}{n}\sum_{k=1}^{n}a_{k}\right\vert ^{2}\right) ^{%
\frac{1}{2}}\left( \frac{1}{n}\sum_{k=1}^{n}\left\Vert x_{k}\right\Vert
^{2}-\left\Vert \frac{1}{n}\sum_{k=1}^{n}x_{k}\right\Vert ^{2}\right) ^{%
\frac{1}{2}}
\end{multline}%
for all $a_{k}\in \mathbb{K},\;x_{k}\in H\;\left( k=1,\dots ,n\right) .$

However, the $x_{k}\;\left( k=1,\dots ,n\right) $ satisfy $\left( \ref%
{4.2.15}\right) ,$ and therefore, by Lemma \ref{l2.1.15}, we have 
\begin{equation}
0\leq \frac{1}{n}\sum_{k=1}^{n}\left\Vert x_{k}\right\Vert ^{2}-\left\Vert 
\frac{1}{n}\sum_{k=1}^{n}x_{k}\right\Vert ^{2}\leq \frac{1}{4}\left\Vert
X-x\right\Vert ^{2}.  \label{4.5.15}
\end{equation}%
Consequently, by $\left( \ref{4.4.15}\right) $ and $\left( \ref{4.5.15}%
\right) ,$ we conclude that 
\begin{multline}
\left\Vert \sum_{k=1}^{n}a_{k}x_{k}-\sum_{k=1}^{n}a_{k}\cdot \frac{1}{n}%
\sum_{k=1}^{n}x_{k}\right\Vert  \label{4.6.15} \\
\leq \frac{1}{2}\left\Vert X-x\right\Vert \left( n\sum_{k=1}^{n}\left\vert
a_{k}\right\vert ^{2}-\left\vert \sum_{k=1}^{n}a_{k}\right\vert ^{2}\right)
^{\frac{1}{2}}
\end{multline}%
for all $a_{k}\in \mathbb{K\;}\left( k=1,\dots ,n\right) .$

We now choose in $\left( \ref{4.6.15}\right) ,$ $a_{k}=\exp \left(
2wimk\right) $ to obtain 
\begin{multline}
\left\Vert \mathcal{F}_{w}\left( \mathbf{\bar{x}}\right) \left( m\right)
-\sum_{k=1}^{n}\exp \left( 2wimk\right) \times \frac{1}{n}%
\sum_{k=1}^{n}x_{k}\right\Vert  \label{4.7.15} \\
\leq \frac{1}{2}\left\Vert X-x\right\Vert \left( n\sum_{k=1}^{n}\left\vert
\exp \left( 2wimk\right) \right\vert ^{2}-\left\vert \sum_{k=1}^{n}\exp
\left( 2wimk\right) \right\vert ^{2}\right) ^{\frac{1}{2}}
\end{multline}%
for all $m\in \left\{ 1,\dots ,n\right\} .$

As a simple calculation reveals that 
\begin{align*}
\sum_{k=1}^{n}\exp \left( 2wimk\right) & =\exp \left( 2wim\right) \times 
\left[ \frac{\exp \left( 2wimn\right) -1}{\exp \left( 2wim\right) -1}\right]
\\
& =\exp \left( 2wim\right) \times \left[ \frac{\cos \left( 2wmn\right)
+i\sin \left( 2wmn\right) -1}{\cos \left( 2wm\right) +i\sin \left(
2wm\right) -1}\right] \\
& =\exp \left( 2wim\right) \times \frac{\sin \left( wmn\right) }{\sin \left(
wm\right) }\left[ \frac{\cos \left( wmn\right) +i\sin \left( wmn\right) }{%
\cos \left( wm\right) +i\sin \left( wm\right) }\right] \\
& =\frac{\sin \left( wmn\right) }{\sin \left( wm\right) }\times \exp \left(
2wim\right) \left[ \frac{\exp \left( iwmn\right) }{\exp \left( iwm\right) }%
\right] \\
& =\frac{\sin \left( wmn\right) }{\sin \left( wm\right) }\times \exp \left[
w\left( n+1\right) im\right] ,
\end{align*}%
\begin{equation*}
\sum_{k=1}^{n}\left\vert \exp \left( 2wimk\right) \right\vert ^{2}=n
\end{equation*}%
and 
\begin{equation*}
\left\vert \sum_{k=1}^{n}\exp \left( 2wimk\right) \right\vert ^{2}=\frac{%
\sin ^{2}\left( wmn\right) }{\sin ^{2}\left( wm\right) },\text{ for }w\neq 
\frac{l}{m}\pi ,\;l\in \mathbb{Z},
\end{equation*}%
thus, from $\left( \ref{4.7.15}\right) ,$ we deduce the desired inequality $%
\left( \ref{4.3.15}\right) .$
\end{proof}

\begin{remark}
\label{r4.2.15}The assumption $\left( \ref{4.2.15}\right) $ can be replaced
by the more general condition 
\begin{equation*}
\sum_{i=1}^{n}\func{Re}\left\langle X-x_{i},x_{i}-x\right\rangle \geq 0,
\end{equation*}%
and the conclusion $\left( \ref{4.3.15}\right) $ will still remain valid.
\end{remark}

The following corollary is an obvious consequence of $\left( \ref{4.3.15}%
\right) .$

\begin{corollary}
\label{c4.3.15}Let $a_{i}\in \mathbb{K}$ $\left( i=1,\dots ,n\right) .$ If $%
a,A\in \mathbb{K}$ are such that 
\begin{equation}
\func{Re}\left[ \left( A-a_{i}\right) \left( \bar{a}_{i}-\bar{a}\right) %
\right] \geq 0\text{ for all }i\in \left\{ 1,\dots ,n\right\} ,
\label{4.9.15}
\end{equation}%
then we have an approximation of the Fourier transform for the vector $%
\overline{\mathbf{a}}=\left( a_{1},\dots ,a_{n}\right) \in \mathbb{K}^{n}:$%
\begin{multline}
\left\Vert \mathcal{F}_{w}\left( \overline{\mathbf{a}}\right) \left(
m\right) -\frac{\sin \left( wmn\right) }{\sin \left( wm\right) }\exp \left[
w\left( n+1\right) im\right] \times \frac{1}{n}\sum_{k=1}^{n}a_{k}\right\Vert
\label{4.10.15} \\
\leq \frac{1}{2}\left\vert A-a\right\vert \left[ n^{2}-\frac{\sin ^{2}\left(
wmn\right) }{\sin ^{2}\left( wm\right) }\right] ^{\frac{1}{2}},
\end{multline}%
for all $m\in \left\{ 1,\dots ,n\right\} $ and $w\in \mathbb{R}$ so that $%
w\neq \frac{l}{m}\pi ,\;l\in \mathbb{Z}.$
\end{corollary}

\begin{remark}
\label{r4.5.15}If we assume that $\mathbb{K=R},$ then $\left( \ref{4.9.15}%
\right) $ is equivalent to 
\begin{equation}
a\leq a_{i}\leq A\text{ for all }i\in \left\{ 1,\dots ,n\right\} .
\label{4.12.15}
\end{equation}%
Consequently, with the assumption $\left( \ref{4.12.15}\right) ,$ we obtain
the following approximation of the Fourier transform 
\begin{multline*}
\left\Vert \mathcal{F}_{w}\left( \overline{\mathbf{a}}\right) \left(
m\right) -\frac{\sin \left( wmn\right) }{\sin \left( wm\right) }\exp \left[
w\left( n+1\right) im\right] \times \frac{1}{n}\sum_{k=1}^{n}a_{k}\right\Vert
\\
\leq \frac{1}{2}\left( A-a\right) \left[ n^{2}-\frac{\sin ^{2}\left(
wmn\right) }{\sin ^{2}\left( wm\right) }\right] ^{\frac{1}{2}},
\end{multline*}%
for all $m\in \left\{ 1,\dots ,n\right\} $ and $w\neq \frac{l}{m}\pi ,\;l\in 
\mathbb{Z}.$
\end{remark}

\subsection{Applications for the Discrete Mellin Transform}

Let $\left( H;\left\langle \cdot ,\cdot \right\rangle \right) $ be an inner
product over $\mathbb{R}$ and $\mathbf{\bar{x}}=\left( x_{1},\dots
,x_{n}\right) $ be a sequence of vectors in $H.$

Define the \textit{Mellin transform}: 
\begin{equation*}
\mathcal{M}\left( \mathbf{\bar{x}}\right) \left( m\right)
:=\sum_{k=1}^{n}k^{m-1}x_{k},\;m=1,\dots ,n;
\end{equation*}%
of the sequence $\mathbf{\bar{x}}\in H^{n}.$

The following approximation result holds \cite{16NSSD}.

\begin{theorem}
\label{t5.1.15}Let $H$ and $\mathbf{\bar{x}}$ be as above. If there exist
the vectors $x,X\in H$ such that 
\begin{equation}
\func{Re}\left\langle X-x_{k},x_{k}-x\right\rangle \geq 0\text{ for all }%
k=1,\dots ,n;  \label{5.2.15}
\end{equation}%
then we have the inequality 
\begin{multline}
\left\Vert \mathcal{M}\left( \mathbf{\bar{x}}\right) \left( m\right)
-S_{m-1}\left( n\right) \cdot \frac{1}{n}\sum_{k=1}^{n}x_{k}\right\Vert
\label{5.3.15} \\
\leq \frac{1}{2}\left\Vert X-x\right\Vert \left[ nS_{2m-2}\left( n\right)
-S_{m-1}^{2}\left( n\right) \right] ^{\frac{1}{2}},\;m\in \left\{ 1,\dots
,n\right\} ,
\end{multline}%
where $S_{p}\left( n\right) ,\;p\in \mathbb{R},\;n\in \mathbb{N}$ is the $p-$%
powered sum of the first $n$ natural numbers, i.e., 
\begin{equation*}
S_{p}\left( n\right) :=\sum_{k=1}^{n}k^{p}.
\end{equation*}
\end{theorem}

\begin{proof}
We apply the inequality $\left( \ref{4.6.15}\right) $ to obtain 
\begin{align*}
& \left\Vert \sum_{k=1}^{n}k^{m-1}x_{k}-\sum_{k=1}^{n}k^{m-1}\cdot \frac{1}{n%
}\sum_{k=1}^{n}x_{k}\right\Vert \\
& \leq \frac{1}{2}\left\Vert X-x\right\Vert \left[ n\sum_{k=1}^{n}k^{2\left(
m-1\right) }-\left( \sum_{k=1}^{n}k^{m-1}\right) ^{2}\right] ^{\frac{1}{2}}
\\
& =\frac{1}{2}\left\Vert X-x\right\Vert \left[ nS_{2m-2}\left( n\right)
-S_{m-1}^{2}\left( n\right) \right] ^{\frac{1}{2}},
\end{align*}%
and the inequality $\left( \ref{5.3.15}\right) $ is proved.
\end{proof}

Consider the following particular values of Mellin Transform 
\begin{equation*}
\mu _{1}\left( \mathbf{\bar{x}}\right) :=\sum_{k=1}^{n}kx_{k}
\end{equation*}%
and 
\begin{equation*}
\mu _{2}\left( \mathbf{\bar{x}}\right) :=\sum_{k=1}^{n}k^{2}x_{k}.
\end{equation*}%
The following corollary holds.

\begin{corollary}
\label{c5.2.15}Let $H$ and $\mathbf{\bar{x}}$ be as in Theorem \ref{t5.1.15}%
. Then we have the inequalities: 
\begin{equation}
\left\Vert \mu _{1}\left( \mathbf{\bar{x}}\right) -\frac{n+1}{2}\cdot
\sum_{k=1}^{n}x_{k}\right\Vert \leq \frac{1}{2}\left\Vert X-x\right\Vert n%
\left[ \frac{n\left( n+1\right) }{2}\right] ^{\frac{1}{2}}  \label{5.5.15}
\end{equation}%
and 
\begin{multline}
\left\Vert \mu _{2}\left( \bar{x}\right) -\frac{\left( n+1\right) \left(
2n+1\right) }{6}\cdot \sum_{k=1}^{n}x_{k}\right\Vert  \label{5.6.15} \\
\leq \frac{1}{12\sqrt{5}}\left\Vert X-x\right\Vert n\sqrt{\left( n-1\right)
\left( n+1\right) \left( 2n+1\right) \left( 8n+1\right) }.
\end{multline}
\end{corollary}

\begin{remark}
\label{r5.3.15}If we assume that $\overline{\mathbf{p}}=\left( p_{1},\dots
,p_{n}\right) $ is a probability distribution, i.e., $p_{k}\geq 0$ $\left(
k=1,\dots ,n\right) $ and $\sum_{k=1}^{n}p_{k}=1$ and $p\leq p_{k}\leq P$ $%
\left( k=1,\dots ,n\right) ,$ then by $\left( \ref{5.5.15}\right) $ and $%
\left( \ref{5.6.15}\right) ,$ we get the inequalities 
\begin{equation*}
\left\vert \sum_{k=1}^{n}kp_{k}-\frac{n+1}{2}\right\vert \leq \frac{1}{2}%
\left( P-p\right) n\left[ \frac{n\left( n+1\right) }{2}\right] ^{\frac{1}{2}}
\end{equation*}%
and 
\begin{multline*}
\left\vert \sum_{k=1}^{n}k^{2}p_{k}-\frac{\left( n+1\right) \left(
2n+1\right) }{6}\right\vert \\
\leq \frac{1}{12\sqrt{5}}\left( P-p\right) n\sqrt{\left( n-1\right) \left(
n+1\right) \left( 2n+1\right) \left( 8n+1\right) }.
\end{multline*}
\end{remark}

\subsection{Applications for Polynomials}

Let $\left( H;\left\langle \cdot ,\cdot \right\rangle \right) $ be an inner
product space over $\mathbb{K}$ and $\mathbf{\bar{c}}=\left( c_{0},\dots
,c_{n}\right) $ be a sequence of vectors in $H.$

Define the polynomial $P:\mathbb{C\rightarrow }H$ with the coefficients $%
\mathbf{\bar{c}}=\left( c_{0},\dots ,c_{n}\right) $ by 
\begin{equation*}
P\left( z\right) =c_{0}+zc_{1}+\cdots +z^{n}c_{n},\;z\in \mathbb{C}%
,\;c_{n}\neq 0.
\end{equation*}%
The following approximation result for the polynomial $P$ holds \cite{16NSSD}%
.

\begin{theorem}
\label{t6.1.15}Let $H,\mathbf{\bar{c}}$ and $P$ be as above. If there exist
the vectors $c,C\in H$ such that 
\begin{equation}
\func{Re}\left\langle C-c_{k},c_{k}-c\right\rangle \geq 0\text{ for all }%
k\in \left\{ 0,\dots ,n\right\} ,  \label{6.1.15}
\end{equation}%
then we have the inequality 
\begin{multline}
\left\Vert P\left( z\right) -\frac{z^{n+1}-1}{z-1}\times \frac{%
c_{0}+c_{1}+\cdots +c_{n}}{n+1}\right\Vert  \label{6.2.15} \\
\leq \frac{1}{2}\left\Vert C-c\right\Vert \left[ \left( n+1\right) \frac{%
\left\vert z\right\vert ^{2n+2}-1}{\left\vert z\right\vert ^{2}-1}-\frac{%
\left\vert z\right\vert ^{2n+2}-2\func{Re}\left( z^{n+1}\right) +1}{%
\left\vert z\right\vert ^{2}-2\func{Re}\left( z\right) +1}\right] ^{\frac{1}{%
2}}
\end{multline}%
for all $z\in \mathbb{C},$ $\left\vert z\right\vert \neq 1.$
\end{theorem}

\begin{proof}
Using the inequality $\left( \ref{4.6.15}\right) ,$ we can state that 
\begin{align}
& \left\Vert \sum_{k=0}^{n}z^{k}c_{k}-\sum_{k=0}^{n}z^{k}\cdot \frac{1}{n+1}%
\sum_{k=0}^{n}c_{k}\right\Vert  \label{6.3.15} \\
& \leq \frac{1}{2}\left\Vert C-c\right\Vert \left( \left( n+1\right)
\sum_{k=0}^{n}\left\vert z\right\vert ^{2k}-\left\vert
\sum_{k=0}^{n}z^{k}\right\vert ^{2}\right) ^{\frac{1}{2}}  \notag \\
& =\frac{1}{2}\left\Vert C-c\right\Vert \left[ \left( n+1\right) \frac{%
\left\vert z\right\vert ^{2n+2}-1}{\left\vert z\right\vert ^{2}-1}%
-\left\vert \frac{z^{n+1}-1}{z-1}\right\vert ^{2}\right] ^{\frac{1}{2}} 
\notag \\
& =\frac{1}{2}\left\Vert C-c\right\Vert \left[ \left( n+1\right) \frac{%
\left\vert z\right\vert ^{2n+2}-1}{\left\vert z\right\vert ^{2}-1}-\frac{%
\left\vert z\right\vert ^{2n+2}-2\func{Re}\left( z^{n+1}\right) +1}{%
\left\vert z\right\vert ^{2}-2\func{Re}\left( z\right) +1}\right] ^{\frac{1}{%
2}}  \notag
\end{align}%
and the inequality $\left( \ref{6.2.15}\right) $ is proved.
\end{proof}

The following result for the complex roots of the unity also holds \cite%
{16NSSD}.

\begin{theorem}
\label{t6.2.15}Let $z_{k}:=\cos \left( \frac{k\pi }{n+1}\right) +i\sin
\left( \frac{k\pi }{n+1}\right) ,\;k\in \left\{ 0,\dots ,n\right\} $ be the
complex $\left( n+1\right) -$roots of the unity. Then we have the inequality 
\begin{equation}
\left\Vert P\left( z_{k}\right) \right\Vert \leq \frac{1}{2}\left(
n+1\right) \left\Vert C-c\right\Vert ,\;k\in \left\{ 1,\dots ,n\right\} ;
\label{6.4.15}
\end{equation}%
where the coefficients $\mathbf{\bar{c}}=\left( c_{0},\dots ,c_{n}\right)
\in H^{n+1}$ satisfy the assumption $\left( \ref{6.1.15}\right) .$
\end{theorem}

\begin{proof}
From the inequality $\left( \ref{6.3.15}\right) ,$ we can state that 
\begin{multline}
\left\Vert P\left( z_{k}\right) -\frac{z^{n+1}-1}{z-1}\cdot \frac{1}{n+1}%
\sum_{k=0}^{n}c_{k}\right\Vert  \notag \\
\leq \frac{1}{2}\left\Vert C-c\right\Vert \left[ \left( n+1\right)
\sum_{k=0}^{n}\left\vert z\right\vert ^{2k}-\left\vert \frac{z^{n+1}-1}{z-1}%
\right\vert ^{2}\right] ^{\frac{1}{2}}
\end{multline}%
for all $z\in \mathbb{C},\;z\neq 1.$

Putting $z=z_{k},\;k\in \left\{ 1,\dots ,n\right\} $ and taking into account
that $z_{k}^{n+1}=1,\;\left\vert z_{k}\right\vert =1,$ we get the desired
result $\left( \ref{6.4.15}\right) .$
\end{proof}

The following corollary is a natural consequence of Theorem \ref{t6.2.15}.

\begin{corollary}
\label{c6.3.15}Let $P\left( z\right) :=\sum_{k=0}^{n}a_{k}z^{k}$ be a
polynomial with real coefficients and $z_{k}$ the $\left( n+1\right) $-roots
of the unity as defined above. If $a\leq a_{k}\leq A,\;k=0,\dots ,n,$ then
we have the inequality: 
\begin{equation*}
\left\vert P\left( z_{k}\right) \right\vert \leq \frac{1}{2}\left(
n+1\right) \left( A-a\right) .
\end{equation*}
\end{corollary}

\subsection{Applications for Lipschitzian Mappings}

Let $\left( H;\left\langle \cdot ,\cdot \right\rangle \right) $ be as above
and $F:H\rightarrow B$ a mapping defined on the inner product space $H$ with
values in the normed linear space $B$ which satisfy the \textit{Lipschitzian
condition:} 
\begin{equation}
\left\vert F\left( x\right) -F\left( y\right) \right\vert \leq L\left\Vert
x-y\right\Vert ,\;\text{for\ all\ }x,y\in H,  \label{7.1.15}
\end{equation}%
where $\left\vert \,\mathbf{\cdot }\,\right\vert $ denotes the norm on $B$
and $\left\Vert \,\mathbf{\cdot }\,\right\Vert $ is the Euclidean norm on $H$%
.

The following theorem holds \cite{16NSSD}.

\begin{theorem}
\label{t7.1.15}Let $F:H\rightarrow B$ be as above and $x_{i}\in
H,\;p_{i}\geq 0$ $\left( i=1,\dots ,n\right) $ with $P_{n}:=%
\sum_{i=1}^{n}p_{i}>0$. If there exists two vectors $x,X\in H$ such that 
\begin{equation}
\func{Re}\left\langle X-x_{i},x_{i}-x\right\rangle \geq 0\text{ for all }%
i\in \left\{ 1,\dots ,n\right\} ,  \label{7.2.15}
\end{equation}%
then we have the inequality 
\begin{equation}
\left\vert \frac{1}{P_{n}}\sum_{i=1}^{n}p_{i}F\left( x_{i}\right) -F\left( 
\frac{1}{P_{n}}\sum_{i=1}^{n}p_{i}x_{i}\right) \right\vert \leq \frac{1}{2}%
\cdot L\left\Vert X-x\right\Vert .  \label{7.3.15}
\end{equation}
\end{theorem}

\begin{proof}
As $F$ is Lipschitzian, we have $\left( \ref{7.1.15}\right) $ for all $%
x,y\in H.$ Choose $x=\frac{1}{P_{n}}\sum_{i=1}^{n}p_{i}x_{i}$ and $%
y=x_{j}\;\left( j=1,\dots ,n\right) ,$ to get 
\begin{equation}
\left\vert F\left( \frac{1}{P_{n}}\sum_{i=1}^{n}p_{i}x_{i}\right) -F\left(
x_{j}\right) \right\vert \leq L\left\Vert \frac{1}{P_{n}}%
\sum_{i=1}^{n}p_{i}x_{i}-x_{j}\right\Vert ,  \label{7.4.15}
\end{equation}%
for all $j\in \left\{ 1,\dots ,n\right\} $.

If we multiply $\left( \ref{7.4.15}\right) $ by $p_{j}\geq 0$ and sum over $%
j $ from $1$ to $n$, we obtain 
\begin{equation}
\sum_{j=1}^{n}p_{j}\left\vert F\left( \frac{1}{P_{n}}%
\sum_{i=1}^{n}p_{i}x_{i}\right) -F\left( x_{j}\right) \right\vert \leq
L\sum_{j=1}^{n}p_{j}\left\Vert \frac{1}{P_{n}}\sum_{i=1}^{n}p_{i}x_{i}-x_{j}%
\right\Vert .  \label{7.5.15}
\end{equation}%
Using the generalized triangle inequality, we have 
\begin{multline}
\sum_{j=1}^{n}p_{j}\left\vert F\left( \frac{1}{P_{n}}%
\sum_{i=1}^{n}p_{i}x_{i}\right) -F\left( x_{j}\right) \right\vert
\label{7.6.15} \\
\geq \left\vert P_{n}F\left( \frac{1}{P_{n}}\sum_{i=1}^{n}p_{i}x_{i}\right)
-\sum_{j=1}^{n}p_{j}F\left( x_{j}\right) \right\vert .
\end{multline}

By the Cauchy-Buniakowsky-Schwarz inequality, we also have 
\begin{align}
& \sum_{j=1}^{n}p_{j}\left\Vert \frac{1}{P_{n}}%
\sum_{i=1}^{n}p_{i}x_{i}-x_{j}\right\Vert  \label{7.7.15} \\
& \leq \left[ \sum_{j=1}^{n}p_{j}\left\Vert \frac{1}{P_{n}}%
\sum_{i=1}^{n}p_{i}x_{i}-x_{j}\right\Vert ^{2}\right] ^{\frac{1}{2}}P_{n}^{%
\frac{1}{2}}  \notag \\
& =P_{n}^{\frac{1}{2}}\left[ \sum_{j=1}^{n}p_{j}\left[ \left\Vert \frac{1}{%
P_{n}}\sum_{i=1}^{n}p_{i}x_{i}\right\Vert ^{2}-2\func{Re}\left\langle \frac{1%
}{P_{n}}\sum_{i=1}^{n}p_{i}x_{i},x_{j}\right\rangle +\left\Vert
x_{j}\right\Vert ^{2}\right] \right] ^{\frac{1}{2}}  \notag \\
& =P_{n}^{\frac{1}{2}}\left[ P_{n}\left\Vert \frac{1}{P_{n}}%
\sum_{i=1}^{n}p_{i}x_{i}\right\Vert ^{2}-2\func{Re}\left\langle \frac{1}{%
P_{n}}\sum_{i=1}^{n}p_{i}x_{i},\sum_{j=1}^{n}p_{j}x_{j}\right\rangle
+\sum_{j=1}^{n}p_{j}\left\Vert x_{j}\right\Vert ^{2}\right] ^{\frac{1}{2}} 
\notag \\
& =P_{n}\left[ \frac{1}{P_{n}}\sum_{i=1}^{n}p_{i}\left\Vert x_{i}\right\Vert
^{2}-\left\Vert \frac{1}{P_{n}}\sum_{i=1}^{n}p_{i}x_{i}\right\Vert ^{2}%
\right] ^{\frac{1}{2}}.  \notag
\end{align}%
Combining the above inequalities $\left( \ref{7.5.15}\right) -\left( \ref%
{7.7.15}\right) $ we deduce, by dividing with $P_{n}>0$, that 
\begin{multline}
\left\vert F\left( \frac{1}{P_{n}}\sum_{i=1}^{n}p_{i}x_{i}\right) -\frac{1}{%
P_{n}}\sum_{i=1}^{n}p_{i}F\left( x_{i}\right) \right\vert  \label{7.8.15} \\
\leq L\cdot \left[ \frac{1}{P_{n}}\sum_{i=1}^{n}p_{i}\left\Vert
x_{i}\right\Vert ^{2}-\left\Vert \frac{1}{P_{n}}\sum_{i=1}^{n}p_{i}x_{i}%
\right\Vert ^{2}\right] ^{\frac{1}{2}}.
\end{multline}%
Finally, using Lemma \ref{l2.1.15}, we obtain the desired result.
\end{proof}

\begin{remark}
\label{r7.9.15}The condition $\left( \ref{7.2.15}\right) $ can be
substituted by the more general condition 
\begin{equation*}
\sum_{i=1}^{n}p_{i}\func{Re}\left\langle X-x_{i},x_{i}-x\right\rangle \geq 0,
\end{equation*}%
and the conclusion $\left( \ref{7.3.15}\right) $ will still remain valid.
\end{remark}

The following corollary is a natural consequence of the above findings.

\begin{corollary}
\label{c7.10.15}Let $x_{i}\in H$ $\left( i=1,\dots ,n\right) $ and $x,X\in H$
be such that the condition $\left( \ref{7.2.15}\right) $ holds. Then we have
the inequality 
\begin{equation*}
0\leq \frac{1}{P_{n}}\sum_{i=1}^{n}p_{i}\left\Vert x_{i}\right\Vert
-\left\Vert \frac{1}{P_{n}}\sum_{i=1}^{n}p_{i}x_{i}\right\Vert \leq \frac{1}{%
2}\left\Vert X-x\right\Vert .
\end{equation*}
\end{corollary}

The proof follows by Theorem \ref{t7.1.15} by choosing $F:H\rightarrow 
\mathbb{R}$, $F\left( x\right) =\left\Vert x\right\Vert $ which is
Lipschitzian with the constant $L=1$, as 
\begin{equation*}
\left\vert F\left( x\right) -F\left( y\right) \right\vert =\left\vert
\left\Vert x\right\Vert -\left\Vert y\right\Vert \right\vert \leq \left\Vert
x-y\right\Vert ,
\end{equation*}
for all $x,y\in H$.

\newpage

\section{The Version for Inner-Products}

\subsection{A Discrete Inequality of Gr\"{u}ss Type}

The following Gr\"{u}ss type inequality holds \cite{17NSSD}.

\begin{theorem}
\label{gtisvipsAppt3.1}Let $\left( H;\left\langle \cdot ,\cdot \right\rangle
\right) $ be an inner product space over $\mathbb{K};$ $\mathbb{K=C},\mathbb{%
R},$ $x_{i},y_{i}\in H$, $p_{i}\geq 0$ $\left( i=0,\dots ,n\right) $ \ $%
\left( n\geq 2\right) $ with $\sum_{i=1}^{n}p_{i}=1$. If $x,X,y,Y\in H$ are
such that 
\begin{equation}
\func{Re}\left\langle X-x_{i},x_{i}-x\right\rangle \geq 0\text{ and\ }\func{%
Re}\left\langle Y-y_{i},y_{i}-y\right\rangle \geq 0\text{ }
\label{gtisvipsApp3.1}
\end{equation}%
\ for all $i\in \left\{ 1,\dots ,n\right\} ,$ then we have the inequality 
\begin{equation}
\left\vert \sum_{i=1}^{n}p_{i}\left\langle x_{i},y_{i}\right\rangle
-\left\langle \sum_{i=1}^{n}p_{i}x_{i},\sum_{i=1}^{n}p_{i}y_{i}\right\rangle
\right\vert \leq \frac{1}{4}\left\Vert X-x\right\Vert \left\Vert
Y-y\right\Vert .  \label{gtisvipsApp3.2}
\end{equation}%
The constant $\frac{1}{4}$ is sharp.
\end{theorem}

\begin{proof}
A simple calculation shows that 
\begin{multline}
\sum_{i=1}^{n}p_{i}\left\langle x_{i},y_{i}\right\rangle -\left\langle
\sum_{i=1}^{n}p_{i}x_{i},\sum_{i=1}^{n}p_{i}y_{i}\right\rangle
\label{gtisvipsApp3.3} \\
=\frac{1}{2}\sum_{i,j=1}^{n}p_{i}p_{j}\left\langle
x_{i}-x_{j},y_{i}-y_{j}\right\rangle .
\end{multline}%
Taking the modulus in both parts of (\ref{gtisvipsApp3.3}) and using the
generalized triangle inequality, we obtain 
\begin{equation*}
\left\vert \sum_{i=1}^{n}p_{i}\left\langle x_{i},y_{i}\right\rangle
-\left\langle \sum_{i=1}^{n}p_{i}x_{i},\sum_{i=1}^{n}p_{i}y_{i}\right\rangle
\right\vert \leq \frac{1}{2}\sum_{i,j=1}^{n}p_{i}p_{j}\left\vert
\left\langle x_{i}-x_{j},y_{i}-y_{j}\right\rangle \right\vert .
\end{equation*}%
By Schwarz's inequality in inner product spaces we have 
\begin{equation*}
\left\vert \left\langle x_{i}-x_{j},y_{i}-y_{j}\right\rangle \right\vert
\leq \left\Vert x_{i}-x_{j}\right\Vert \left\Vert y_{i}-y_{j}\right\Vert ,
\end{equation*}%
for all $i,j\in \left\{ 1,\dots ,n\right\} ,$ and therefore 
\begin{equation*}
\left\vert \sum_{i=1}^{n}p_{i}\left\langle x_{i},y_{i}\right\rangle
-\left\langle \sum_{i=1}^{n}p_{i}x_{i},\sum_{i=1}^{n}p_{i}y_{i}\right\rangle
\right\vert \leq \frac{1}{2}\sum_{i,j=1}^{n}p_{i}p_{j}\left\Vert
x_{i}-x_{j}\right\Vert \left\Vert y_{i}-y_{j}\right\Vert .
\end{equation*}%
Using the Cauchy-Bunyakovsky-Schwarz inequality for double sums, we can
state that 
\begin{multline*}
\frac{1}{2}\sum_{i,j=1}^{n}p_{i}p_{j}\left\Vert x_{i}-x_{j}\right\Vert
\left\Vert y_{i}-y_{j}\right\Vert \\
\leq \left( \frac{1}{2}\sum_{i,j=1}^{n}p_{i}p_{j}\left\Vert
x_{i}-x_{j}\right\Vert ^{2}\right) ^{\frac{1}{2}}\left( \frac{1}{2}%
\sum_{i,j=1}^{n}p_{i}p_{j}\left\Vert y_{i}-y_{j}\right\Vert ^{2}\right) ^{%
\frac{1}{2}}
\end{multline*}%
and, as a simple calculation shows that, 
\begin{equation*}
\frac{1}{2}\sum_{i,j=1}^{n}p_{i}p_{j}\left\Vert x_{i}-x_{j}\right\Vert
^{2}=\sum_{i=1}^{n}p_{i}\left\Vert x_{i}\right\Vert ^{2}-\left\Vert
\sum_{i=1}^{n}p_{i}x_{i}\right\Vert ^{2}
\end{equation*}%
and 
\begin{equation*}
\frac{1}{2}\sum_{i,j=1}^{n}p_{i}p_{j}\left\Vert y_{i}-y_{j}\right\Vert
^{2}=\sum_{i=1}^{n}p_{i}\left\Vert y_{i}\right\Vert ^{2}-\left\Vert
\sum_{i=1}^{n}p_{i}y_{i}\right\Vert ^{2},
\end{equation*}%
we obtain 
\begin{multline}
\left\vert \sum_{i=1}^{n}p_{i}\left\langle x_{i},y_{i}\right\rangle
-\left\langle \sum_{i=1}^{n}p_{i}x_{i},\sum_{i=1}^{n}p_{i}y_{i}\right\rangle
\right\vert  \label{gtisvipsApp3.8} \\
\leq \left( \sum_{i=1}^{n}p_{i}\left\Vert x_{i}\right\Vert ^{2}-\left\Vert
\sum_{i=1}^{n}p_{i}x_{i}\right\Vert ^{2}\right) ^{\frac{1}{2}}\left(
\sum_{i=1}^{n}p_{i}\left\Vert y_{i}\right\Vert ^{2}-\left\Vert
\sum_{i=1}^{n}p_{i}y_{i}\right\Vert ^{2}\right) ^{\frac{1}{2}}.
\end{multline}%
Using Lemma \ref{l2.1.15}, we know that 
\begin{equation*}
\left( \sum_{i=1}^{n}p_{i}\left\Vert x_{i}\right\Vert ^{2}-\left\Vert
\sum_{i=1}^{n}p_{i}x_{i}\right\Vert ^{2}\right) ^{\frac{1}{2}}\leq \frac{1}{2%
}\left\Vert X-x\right\Vert
\end{equation*}%
and 
\begin{equation*}
\left( \sum_{i=1}^{n}p_{i}\left\Vert y_{i}\right\Vert ^{2}-\left\Vert
\sum_{i=1}^{n}p_{i}y_{i}\right\Vert ^{2}\right) ^{\frac{1}{2}}\leq \frac{1}{2%
}\left\Vert Y-y\right\Vert ,
\end{equation*}%
and then, by (\ref{gtisvipsApp3.8}), we deduce the desired inequality (\ref%
{gtisvipsApp3.3}).

To prove the sharpness of the constant $\frac{1}{4}$, let us assume that (%
\ref{gtisvipsApp3.2}) holds with a constant $c>0$, i.e., 
\begin{equation}
\left| \sum_{i=1}^{n}p_{i}\left\langle x_{i},y_{i}\right\rangle
-\left\langle \sum_{i=1}^{n}p_{i}x_{i},\sum_{i=1}^{n}p_{i}y_{i}\right\rangle
\right| \leq c\left\| X-x\right\| \left\| Y-y\right\|  \label{gtisvipsApp3.9}
\end{equation}
under the above assumptions for $p_{i},\;x_{i},\;y_{i},\;x,\;X,\;y,\;Y$ and $%
n\geq 2$.

If we choose $n=2$, $x_{1}=x,\;x_{2}=X,$ $y_{1}=y,\;y_{2}=Y$ \ $\left( x\neq
X,\;y\neq Y\right) $ and $p_{1}=p_{2}=\frac{1}{2}$, then 
\begin{align*}
\sum_{i=1}^{2}p_{i}\left\langle x_{i},y_{i}\right\rangle -\left\langle
\sum_{i=1}^{2}p_{i}x_{i},\sum_{i=1}^{2}p_{i}y_{i}\right\rangle & =\frac{1}{2}%
\sum_{i,j=1}^{2}p_{i}p_{j}\left\langle x_{i}-x_{j},y_{i}-y_{j}\right\rangle
\\
& =\sum_{1\leq i<j\leq 2}p_{i}p_{j}\left\langle
x_{i}-x_{j},y_{i}-y_{j}\right\rangle \\
& =\frac{1}{4}\left\langle x-X,y-Y\right\rangle
\end{align*}%
and then 
\begin{equation*}
\left\vert \sum_{i=1}^{2}p_{i}\left\langle x_{i},y_{i}\right\rangle
-\left\langle \sum_{i=1}^{2}p_{i}x_{i},\sum_{i=1}^{2}p_{i}y_{i}\right\rangle
\right\vert =\frac{1}{4}\left\vert \left\langle x-X,y-Y\right\rangle
\right\vert .
\end{equation*}%
Choose $X-x=z,$ $Y-y=z,$ $z\neq 0$. Then using (\ref{gtisvipsApp3.9}), we
derive 
\begin{equation*}
\frac{1}{4}\left\Vert z\right\Vert ^{2}\leq c\left\Vert z\right\Vert
^{2},\;\;z\neq 0
\end{equation*}%
which implies that $c\geq \frac{1}{4}$, and the theorem is proved.
\end{proof}

\begin{remark}
\label{gtisvipsAppr3.2}The condition (\ref{gtisvipsApp3.1}) can be replaced
by the more general assumption 
\begin{equation*}
\sum_{i=1}^{n}p_{i}\func{Re}\left\langle X-x_{i},x_{i}-x\right\rangle \geq
0,\;\;\sum_{i=1}^{n}p_{i}\func{Re}\left\langle Y-y_{i},y_{i}-y\right\rangle
\geq 0
\end{equation*}%
and the conclusion (\ref{gtisvipsApp3.2}) still remains valid.
\end{remark}

The following corollary for real or complex numbers holds.

\begin{corollary}
\label{gtisvipsAppc3.3}Let $a_{i},$ $b_{i}\in \mathbb{K}$ $\left( \mathbb{K=C%
},\mathbb{R}\right) ,$ $p_{i}\geq 0$ \ $\left( i=1,\dots ,n\right) $ with $%
\sum_{i=1}^{n}p_{i}=1.$ If $a,A,b,B\in \mathbb{K}$ are such that 
\begin{equation}
\func{Re}\left[ \left( A-a_{i}\right) \left( \bar{a}_{i}-\bar{a}\right) %
\right] \geq 0,\;\;\func{Re}\left[ \left( B-b_{i}\right) \left( \bar{b}_{i}-%
\bar{b}\right) \right] \geq 0,  \label{gtisvipsApp3.11}
\end{equation}%
then we have the inequality 
\begin{equation}
\left\vert \sum_{i=1}^{n}p_{i}a_{i}\bar{b}_{i}-\sum_{i=1}^{n}p_{i}a_{i}%
\sum_{i=1}^{n}p_{i}\bar{b}_{i}\right\vert \leq \frac{1}{4}\left\vert
A-a\right\vert \left\vert B-b\right\vert  \label{gtisvipsApp3.12}
\end{equation}%
and the constant $\frac{1}{4}$ is sharp.
\end{corollary}

The proof is obvious by Theorem \ref{gtisvipsAppt3.1} applied for the inner
product space $\left( \mathbb{C},\left\langle \cdot ,\cdot \right\rangle
\right) ,$ where $\left\langle x,y\right\rangle =x\cdot \bar{y}$. We omit
the details.

\begin{remark}
\label{gtisvipsAppr3.4}The condition (\ref{gtisvipsApp3.11}) can be replaced
by the more general condition 
\begin{equation*}
\sum_{i=1}^{n}p_{i}\func{Re}\left[ \left( A-a_{i}\right) \left( \bar{a}_{i}-%
\bar{a}\right) \right] \geq 0,\;\;\sum_{i=1}^{n}p_{i}\func{Re}\left[ \left(
B-b_{i}\right) \left( \bar{b}_{i}-\bar{b}\right) \right] \geq 0
\end{equation*}%
and the conclusion of the above corollary will still remain valid.
\end{remark}

\begin{remark}
\label{gtisvipsAppr3.5}If we assume that $a_{i},$ $b_{i}$, $a$, $b$, $A$, $B$
are real numbers, then (\ref{gtisvipsApp3.11}) is equivalent to 
\begin{equation*}
a\leq a_{i}\leq A,\;\;b\leq b_{i}\leq B\text{ \ for all }i\in \left\{
1,\dots ,n\right\}
\end{equation*}%
and (\ref{gtisvipsApp3.12}) becomes 
\begin{equation*}
0\leq \left\vert
\sum_{i=1}^{n}p_{i}a_{i}b_{i}-\sum_{i=1}^{n}p_{i}a_{i}%
\sum_{i=1}^{n}p_{i}b_{i}\right\vert \leq \frac{1}{4}\left( A-a\right) \left(
B-b\right) ,
\end{equation*}%
which is the classical Gr\"{u}ss inequality for sequences of real numbers.
\end{remark}

\subsection{Applications for Convex Functions}

Let $\left( H;\left\langle \cdot ,\cdot \right\rangle \right) $ be a real
inner product space and $F:H\rightarrow \mathbb{R}$ a Fr\'{e}chet
differentiable convex mapping on $H$. Then we have the \textquotedblleft 
\textit{gradient inequality}\textquotedblright\ 
\begin{equation}
F\left( x\right) -F\left( y\right) \geq \left\langle \triangledown F\left(
y\right) ,x-y\right\rangle  \label{gtisvipsApp4.1}
\end{equation}%
for all $x,y\in H$, where $\nabla F:H\rightarrow H$ is the gradient operator
associated to the differentiable convex function $F$.

The following theorem holds \cite{17NSSD}.

\begin{theorem}
\label{gtisvipsAppt4.1}Let $F:H\rightarrow \mathbb{R}$ be as above and $%
x_{i}\in H$ $\left( i=1,\dots ,n\right) $. Suppose that there exists the
vectors $x,X\in H$ such that $\left\langle x_{i}-x,X-x_{i}\right\rangle \geq
0$ for all $i\in \left\{ 1,\dots ,m\right\} $ and $y,Y\in H$ such that $%
\left\langle \triangledown F\left( x_{i}\right) -y,\;Y-\triangledown F\left(
x_{i}\right) \right\rangle \geq 0$ for all $i\in \left\{ 1,\dots ,m\right\} $%
. Then for all $p_{i}\geq 0$ $\left( i=1,\dots ,m\right) $ with $%
P_{m}:=\sum_{i=1}^{m}p_{i}>0$, we have the inequality 
\begin{equation}
0\leq \frac{1}{P_{m}}\sum_{i=1}^{m}p_{i}F\left( x_{i}\right) -F\left( \frac{1%
}{P_{m}}\sum_{i=1}^{m}p_{i}x_{i}\right) \leq \frac{1}{4}\left\Vert
X-x\right\Vert \left\Vert Y-y\right\Vert .  \label{gtisvipsApp4.2}
\end{equation}
\end{theorem}

\begin{proof}
Choose in (\ref{gtisvipsApp4.1}), $x=\frac{1}{P_{M}}\sum_{i=1}^{m}p_{i}x_{i}$
and $y=x_{j}$ to obtain 
\begin{equation}
F\left( \frac{1}{P_{m}}\sum_{i=1}^{m}p_{i}x_{i}\right) -F\left( x_{j}\right)
\geq \left\langle \triangledown F\left( x_{j}\right) ,\frac{1}{P_{m}}%
\sum_{i=1}^{m}p_{i}x_{i}-x_{j}\right\rangle  \label{gtisvipsApp4.3}
\end{equation}%
for all $j\in \left\{ 1,\dots ,n\right\} $.

If we multiply (\ref{gtisvipsApp4.3}) by $p_{j}\geq 0$ and sum over $j$ from 
$1$ to $m$, we have 
\begin{multline*}
P_{m}F\left( \frac{1}{P_{m}}\sum_{i=1}^{m}p_{i}x_{i}\right)
-\sum_{j=1}^{m}p_{j}F\left( x_{j}\right) \\
\geq \frac{1}{P_{m}}\left\langle \sum_{j=1}^{m}\triangledown F\left(
x_{j}\right) ,\sum_{i=1}^{m}p_{i}x_{i}\right\rangle
-\sum_{i=1}^{m}\left\langle \triangledown F\left( x_{j}\right)
,x_{j}\right\rangle .
\end{multline*}%
Dividing by $P_{m}>0,$ we obtain the inequality 
\begin{align}
0& \leq \frac{1}{P_{m}}\sum_{i=1}^{m}p_{i}F\left( x_{i}\right) -F\left( 
\frac{1}{P_{m}}\sum_{i=1}^{m}p_{i}x_{i}\right)  \label{gtisvipsApp4.4} \\
& \leq \frac{1}{P_{m}}\sum_{i=1}^{m}p_{i}\left\langle \triangledown F\left(
x_{i}\right) ,x_{i}\right\rangle -\left\langle \frac{1}{P_{m}}%
\sum_{i=1}^{m}p_{i}\triangledown F\left( x_{i}\right) ,\frac{1}{P_{m}}%
\sum_{i=1}^{m}p_{i}x_{i}\right\rangle ,  \notag
\end{align}%
which is a generalisation for the case of inner product spaces of the result
by Dragomir-Goh established in 1996 for the case of differentiable mappings
defined on $\mathbb{R}^{n}$ \cite{gtisvipsApp8a}.

Applying Theorem \ref{gtisvipsAppt3.1} for real inner product spaces, and $%
y_{i}=\nabla F\left( x_{i}\right) $, we easily deduce 
\begin{multline}
\frac{1}{P_{m}}\sum_{i=1}^{m}p_{i}\left\langle x_{i},\triangledown F\left(
x_{i}\right) \right\rangle -\left\langle \frac{1}{P_{m}}%
\sum_{i=1}^{m}p_{i}x_{i},\frac{1}{P_{m}}\sum_{i=1}^{m}p_{i}\triangledown
F\left( x_{i}\right) \right\rangle  \label{gtisvipsApp4.5} \\
\leq \frac{1}{4}\left\Vert X-x\right\Vert \left\Vert Y-y\right\Vert
\end{multline}%
and then, by (\ref{gtisvipsApp4.4}) and (\ref{gtisvipsApp4.5}) we can
conclude that the desired inequality (\ref{gtisvipsApp4.2}) holds.
\end{proof}

\subsection{Applications for Some Discrete Transforms}

Let $\left( H;\left\langle \cdot ,\cdot \right\rangle \right) $ be an inner
product space over $\mathbb{K}$, $\mathbb{K=C}$, $\mathbb{R}$ and $\mathbf{%
\bar{x}}=\left( x_{1},\dots ,x_{n}\right) $ be a sequence of vectors in $H$.

For a given $m\in \mathbb{K}$, define the \textit{discrete Fourier Transform}
\begin{equation*}
\mathcal{F}_{w}\left( \mathbf{\bar{x}}\right) \left( m\right)
=\sum_{k=1}^{n}\exp \left( 2wimk\right) \times x_{k},\;\;m=1,\dots ,n.
\end{equation*}%
The complex number $\sum_{k=1}^{n}\exp \left( 2wimk\right) \left\langle
x_{k},y_{k}\right\rangle $ is actually the usual Fourier transform of the
vector $\left( \left\langle x_{1},y_{1}\right\rangle ,\dots ,\left\langle
x_{n},y_{n}\right\rangle \right) \in \mathbb{K}^{n}$ and will be denoted by 
\begin{equation*}
\mathcal{F}_{w}\left( \mathbf{\bar{x}}\cdot \mathbf{\bar{y}}\right) \left(
m\right) =\sum_{k=1}^{n}\exp \left( 2wimk\right) \left\langle
x_{k},y_{k}\right\rangle ,\;\;m=1,\dots ,n.
\end{equation*}

The following result holds \cite{17NSSD}.

\begin{theorem}
\label{gtisvipsAppt5.1}Let $\mathbf{\bar{x}}$, $\mathbf{\bar{y}}\in H^{n}$
be sequences of vectors such that there exists the vectors $c,C,y,Y\in H$
with the properties 
\begin{equation*}
\func{Re}\left\langle C-\exp \left( 2wimk\right) x_{k},\exp \left(
2wimk\right) x_{k}-c\right\rangle \geq 0,\;\;k,m=1,\dots ,n
\end{equation*}%
and 
\begin{equation}
\func{Re}\left\langle Y-y_{k},y_{k}-y\right\rangle \geq 0,\;\;k=1,\dots ,n.
\label{gtisvipsApp5.4}
\end{equation}%
Then we have the inequality 
\begin{equation*}
\left\vert \mathcal{F}_{w}\left( \mathbf{\bar{x}}\cdot \mathbf{\bar{y}}%
\right) \left( m\right) -\left\langle \mathcal{F}_{w}\left( \mathbf{\bar{x}}%
\right) \left( m\right) ,\frac{1}{n}\sum_{k=1}^{n}y_{k}\right\rangle
\right\vert \leq \frac{n}{4}\left\Vert C-c\right\Vert \left\Vert
Y-y\right\Vert ,
\end{equation*}%
for all $m\in \left\{ 1,\dots ,n\right\} $.
\end{theorem}

The proof follows by Theorem \ref{gtisvipsAppt3.1} applied for $p_{k}=\frac{1%
}{n}$ and for the sequences $x_{k}\rightarrow c_{k}=\exp \left( 2wimk\right)
x_{k}$ and $y_{k}$ \ $\left( k=1,\dots ,n\right) $. We omit the details.

We can also consider the \textit{Mellin transform} 
\begin{equation*}
\mathcal{M}\left( \mathbf{\bar{x}}\right) \left( m\right)
:=\sum_{k=1}^{n}k^{m-1}x_{k},\;\;m=1,\dots ,n,
\end{equation*}%
of the sequence $\mathbf{\bar{x}}=\left( x_{1},\dots ,x_{n}\right) \in H^{n}$%
.

We remark that the complex number $\sum_{k=1}^{n}k^{m-1}\left\langle
x_{k},y_{k}\right\rangle $ is actually the Mellin transform of the vector $%
\left( \left\langle x_{1},y_{1}\right\rangle ,\dots ,\left\langle
x_{n},y_{n}\right\rangle \right) \in \mathbb{K}^{n}$ and will be denoted by 
\begin{equation*}
\mathcal{M}\left( \mathbf{\bar{x}}\cdot \mathbf{\bar{y}}\right) \left(
m\right) :=\sum_{k=1}^{n}k^{m-1}\left\langle x_{k},y_{k}\right\rangle .
\end{equation*}

The following theorem holds \cite{17NSSD}.

\begin{theorem}
\label{gtisvipsAppt5.2}Let $\mathbf{\bar{x}}$, $\mathbf{\bar{y}}\in H^{n}$
be sequences of vectors such that there exist the vectors $d,D,y,Y\in H$
with the properties 
\begin{equation*}
\func{Re}\left\langle D-k^{m-1}x_{k},k^{m-1}x_{k}-d\right\rangle \geq 0
\end{equation*}%
for all $k,m\in \left\{ 1,\dots ,n\right\} ,$ and (\ref{gtisvipsApp5.4}) is
fulfilled.

Then we have the inequality 
\begin{equation*}
\left\vert \mathcal{M}\left( \mathbf{\bar{x}}\cdot \mathbf{\bar{y}}\right)
\left( m\right) -\left\langle \mathcal{M}\left( \bar{x}\right) \left(
m\right) ,\frac{1}{n}\sum_{k=1}^{n}y_{k}\right\rangle \right\vert \leq \frac{%
n}{4}\left\Vert D-d\right\Vert \left\Vert Y-y\right\Vert
\end{equation*}%
for all $m\in \left\{ 1,\dots ,n\right\} $.
\end{theorem}

The proof follows by Theorem \ref{gtisvipsAppt3.1} applied for $p_{k}=\frac{1%
}{n}$ and for the sequences $x_{k}\rightarrow d_{k}=kx_{k}$ and $y_{k}$ $%
\left( k=1,\dots ,n\right) $. We omit the details.

Another result which connects the Fourier transforms for different
parameters $w$ also holds \cite{17NSSD}.

\begin{theorem}
\label{gtisvipsAppt5.3}Let $\mathbf{\bar{x}}$, $\mathbf{\bar{y}}\in H^{n}$
and $w,z\in \mathbb{K}$. If there exist the vectors $e,E,f,F\in H$ such that 
\begin{equation*}
\func{Re}\left\langle E-\exp \left( 2wimk\right) x_{k},\exp \left(
2wimk\right) x_{k}-e\right\rangle \geq 0,\;\;k,m=1,\dots ,n
\end{equation*}%
and 
\begin{equation*}
\func{Re}\left\langle F-\exp \left( 2zimk\right) y_{k},\exp \left(
2zimk\right) y_{k}-f\right\rangle \geq 0,\;\;k,m=1,\dots ,n
\end{equation*}%
then we have the inequality: 
\begin{equation*}
\left\vert \frac{1}{n}\mathcal{F}_{w+z}\left( \mathbf{\bar{x}}\cdot \mathbf{%
\bar{y}}\right) \left( m\right) -\left\langle \frac{1}{n}\mathcal{F}%
_{w}\left( \mathbf{\bar{x}}\right) \left( m\right) ,\frac{1}{n}\mathcal{F}%
_{z}\left( \mathbf{\bar{y}}\right) \left( m\right) \right\rangle \right\vert
\leq \frac{1}{4}\left\Vert E-e\right\Vert \left\Vert F-f\right\Vert ,
\end{equation*}%
for all $m\in \left\{ 1,\dots ,n\right\} $.
\end{theorem}

The proof follows by Theorem \ref{gtisvipsAppt3.1} for the sequences $\exp
\left( 2wimk\right) x_{k}$, $\exp \left( 2zimk\right) y_{k}$ \ $\left(
k=1,\dots ,n\right) $. We omit the details.

\newpage

\section{More Gr\"{u}ss' Type Inequalities}

\subsection{Introduction}

In the recent paper \cite{D3.16}, the author has obtained the following Gr%
\"{u}ss type inequality for forward difference.

\begin{theorem}
\label{t1.3.16}Let $\overline{\mathbf{x}}=\left( x_{1},\dots ,x_{n}\right) ,$
$\overline{\mathbf{y}}=\left( y_{1},\dots ,y_{n}\right) \in H^{n}$ and $%
\overline{\mathbf{p}}\in \mathbb{R}_{+}^{n}$ be a probability sequence. Then
one has the inequalities%
\begin{multline}
\left\vert \sum_{i=1}^{n}p_{i}\left\langle x_{i},y_{i}\right\rangle
-\left\langle \sum_{i=1}^{n}p_{i}x_{i},\sum_{i=1}^{n}p_{i}y_{i}\right\rangle
\right\vert  \label{1.6.16} \\
\leq \left\{ 
\begin{array}{l}
\left[ \sum\limits_{i=1}^{n}i^{2}p_{i}-\left(
\sum\limits_{i=1}^{n}ip_{i}\right) ^{2}\right] \max\limits_{k=1,\dots
,n-1}\left\Vert \Delta x_{k}\right\Vert \max\limits_{k=1,\dots
,n-1}\left\Vert \Delta y_{k}\right\Vert ; \\ 
\\ 
\sum\limits_{1\leq j<i\leq n}p_{i}p_{j}\left( i-j\right) \left(
\sum\limits_{k=1}^{n-1}\left\Vert \Delta x_{k}\right\Vert ^{p}\right) ^{%
\frac{1}{p}}\left( \sum\limits_{k=1}^{n-1}\left\Vert \Delta y_{k}\right\Vert
^{q}\right) ^{\frac{1}{q}} \\ 
\hfill \text{if }p>1,\ \frac{1}{p}+\frac{1}{q}=1; \\ 
\\ 
\dfrac{1}{2}\left[ \sum\limits_{i=1}^{n}p_{i}\left( 1-p_{i}\right) \right]
\sum\limits_{k=1}^{n-1}\left\Vert \Delta x_{k}\right\Vert
\sum\limits_{k=1}^{n-1}\left\Vert \Delta y_{k}\right\Vert .%
\end{array}%
\right.
\end{multline}%
The constants $1,$ $1$ and $\frac{1}{2}$ in the right hand side of the
inequality (\ref{1.6.16}) are best in the sense that they cannot be replaced
by smaller constants.
\end{theorem}

If one chooses $p_{i}=\frac{1}{n}$ $\left( i=1,\dots ,n\right) $ in (\ref%
{1.6.16}), then the following unweighted inequalities hold:%
\begin{multline}
\left\vert \frac{1}{n}\sum_{i=1}^{n}\left\langle x_{i},y_{i}\right\rangle
-\left\langle \frac{1}{n}\sum_{i=1}^{n}x_{i},\frac{1}{n}\sum_{i=1}^{n}y_{i}%
\right\rangle \right\vert  \label{1.7.16} \\
\leq \left\{ 
\begin{array}{l}
\dfrac{n^{2}-1}{12}\max\limits_{k=1,\dots ,n-1}\left\Vert \Delta
x_{k}\right\Vert \max\limits_{k=1,\dots ,n-1}\left\Vert \Delta
y_{k}\right\Vert ; \\ 
\\ 
\dfrac{n^{2}-1}{6n}\left( \sum\limits_{k=1}^{n-1}\left\Vert \Delta
x_{k}\right\Vert ^{p}\right) ^{\frac{1}{p}}\left(
\sum\limits_{k=1}^{n-1}\left\Vert \Delta y_{k}\right\Vert ^{q}\right) ^{%
\frac{1}{q}} \\ 
\hfill \text{if }p>1,\ \frac{1}{p}+\frac{1}{q}=1; \\ 
\\ 
\dfrac{n-1}{2n}\sum\limits_{k=1}^{n-1}\left\Vert \Delta x_{k}\right\Vert
\sum\limits_{k=1}^{n-1}\left\Vert \Delta y_{k}\right\Vert .%
\end{array}%
\right.
\end{multline}%
Here, the constants $\frac{1}{12},$ $\frac{1}{6}$ and $\frac{1}{2}$ are also
best possible in the above sense.

The following reverse of the Cauchy-Bunyakovsky-Schwarz inequality for
sequences of vectors in inner product spaces holds.

\begin{corollary}
\label{c1.4.16}With the assumptions in Theorem \ref{t1.3.16} for $\overline{%
\mathbf{x}}$ and $\overline{\mathbf{p}}$ one has the inequalities%
\begin{align}
0& \leq \sum_{i=1}^{n}p_{i}\left\Vert x_{i}\right\Vert ^{2}-\left\Vert
\sum_{i=1}^{n}p_{i}x_{i}\right\Vert ^{2}  \label{1.8.16} \\
& \leq \left\{ 
\begin{array}{l}
\left[ \sum\limits_{i=1}^{n}i^{2}p_{i}-\left(
\sum\limits_{i=1}^{n}ip_{i}\right) ^{2}\right] \max\limits_{k=\overline{1,n-1%
}}\left\Vert \Delta x_{k}\right\Vert ^{2}; \\ 
\\ 
\sum\limits_{1\leq j<i\leq n}p_{i}p_{j}\left( i-j\right) \left(
\sum\limits_{k=1}^{n-1}\left\Vert \Delta x_{k}\right\Vert ^{p}\right) ^{%
\frac{1}{p}}\left( \sum\limits_{k=1}^{n-1}\left\Vert \Delta x_{k}\right\Vert
^{q}\right) ^{\frac{1}{q}} \\ 
\hfill \text{if }p>1,\ \frac{1}{p}+\frac{1}{q}=1; \\ 
\\ 
\dfrac{1}{2}\left[ \sum\limits_{i=1}^{n}p_{i}\left( 1-p_{i}\right) \right]
\left( \sum\limits_{k=1}^{n-1}\left\Vert \Delta x_{k}\right\Vert \right)
^{2}.%
\end{array}%
\right.  \notag
\end{align}%
The constants $1,$ $1$ and $\frac{1}{2}$ are best possible in the above
sense.
\end{corollary}

The following particular inequalities that may be deduced from (\ref{1.8.16}%
) on choosing the equal weights $p_{i}=\frac{1}{n},$ $i=1,\dots ,n$ are also
of interest%
\begin{align}
0& \leq \frac{1}{n}\sum_{i=1}^{n}\left\Vert x_{i}\right\Vert ^{2}-\left\Vert 
\frac{1}{n}\sum_{i=1}^{n}x_{i}\right\Vert ^{2}  \label{1.9.16} \\
& \leq \left\{ 
\begin{array}{l}
\dfrac{n^{2}-1}{12}\max\limits_{k=\overline{1,n-1}}\left\Vert \Delta
x_{k}\right\Vert ^{2}; \\ 
\\ 
\dfrac{n^{2}-1}{6n}\left( \sum\limits_{k=1}^{n-1}\left\Vert \Delta
x_{k}\right\Vert ^{p}\right) ^{\frac{1}{p}}\left(
\sum\limits_{k=1}^{n-1}\left\Vert \Delta x_{k}\right\Vert ^{q}\right) ^{%
\frac{1}{q}} \\ 
\hfill \text{if }p>1,\ \frac{1}{p}+\frac{1}{q}=1; \\ 
\\ 
\dfrac{n-1}{2n}\left( \sum\limits_{k=1}^{n-1}\left\Vert \Delta
x_{k}\right\Vert \right) ^{2}.%
\end{array}%
\right.  \notag
\end{align}%
Here the constants $\frac{1}{12},$ $\frac{1}{6}$ and $\frac{1}{2}$ are also
best possible.

The main aim of this section is to present, by following \cite{18NSSD}, a
different class of Gr\"{u}ss type inequalities for sequences of vectors in
inner product spaces and to apply them for obtaining a reverse of Jensen's
inequality for convex functions defined on such spaces.

\subsection{More Gr\"{u}ss Type Inequalities}

The following lemma holds (see also \cite{D4.16}).

\begin{lemma}
\label{l2.1.16}Let $a,x,A$ be vectors in the inner product space $\left(
H;\left\langle \cdot ,\cdot \right\rangle \right) $ over the real or complex
number field $\mathbb{K}$ $\left( \mathbb{K=R},\mathbb{C}\right) $ with $%
a\neq A.$ The following statements are equivalent:

\begin{enumerate}
\item[(i)] $\func{Re}\left\langle A-x,x-a\right\rangle \geq 0;$

\item[(ii)] $\left\Vert x-\frac{a+A}{2}\right\Vert \leq \frac{1}{2}%
\left\Vert A-a\right\Vert .$
\end{enumerate}
\end{lemma}

The following inequality of Gr\"{u}ss type for sequences of vectors in inner
product spaces holds \cite{18NSSD}.

\begin{theorem}
\label{t2.3.16}Let $\left( H;\left\langle \cdot ,\cdot \right\rangle \right) 
$ be an inner product over $\mathbb{K}$ $\left( \mathbb{K=C},\mathbb{R}%
\right) $, and $\overline{\mathbf{x}}=\left( x_{1},\dots ,x_{n}\right) ,$ $%
\overline{\mathbf{y}}=\left( y_{1},\dots ,y_{n}\right) \in H^{n},$ $%
\overline{\mathbf{p}}\in \mathbb{R}_{+}^{n}$ with $\sum_{i=1}^{n}p_{i}=1.$
If $x,X\in H$ are such that%
\begin{equation}
\func{Re}\left\langle X-x_{i},x_{i}-x\right\rangle \geq 0\ \text{\ for each
\ }i\in \left\{ 1,\dots ,n\right\} ,  \label{2.1.16}
\end{equation}%
or equivalently,%
\begin{equation}
\left\Vert x_{i}-\frac{x+X}{2}\right\Vert \leq \frac{1}{2}\left\Vert
X-x\right\Vert \ \text{\ for each \ }i\in \left\{ 1,\dots ,n\right\} ,
\label{2.2.16}
\end{equation}%
then one has the inequality%
\begin{align}
& \left\vert \sum_{i=1}^{n}p_{i}\left\langle x_{i},y_{i}\right\rangle
-\left\langle \sum_{i=1}^{n}p_{i}x_{i},\sum_{i=1}^{n}p_{i}y_{i}\right\rangle
\right\vert  \label{2.3.16} \\
& \leq \frac{1}{2}\left\Vert X-x\right\Vert \sum_{i=1}^{n}p_{i}\left\Vert
y_{i}-\sum_{j=1}^{n}p_{j}y_{j}\right\Vert  \notag \\
& \leq \frac{1}{2}\left\Vert X-x\right\Vert \left[ \sum_{i=1}^{n}p_{i}\left%
\Vert y_{i}\right\Vert ^{2}-\left\Vert \sum_{i=1}^{n}p_{i}y_{i}\right\Vert
^{2}\right] ^{\frac{1}{2}}.  \notag
\end{align}%
The constant $\frac{1}{2}$ is best possible in the first and second
inequality in the sense that it cannot be replaced by a smaller constant.
\end{theorem}

\begin{proof}
It is easy to see that the following identity holds true%
\begin{multline}
\sum_{i=1}^{n}p_{i}\left\langle x_{i},y_{i}\right\rangle -\left\langle
\sum_{i=1}^{n}p_{i}x_{i},\sum_{i=1}^{n}p_{i}y_{i}\right\rangle
\label{2.4.16} \\
=\sum_{i=1}^{n}p_{i}\left\langle x_{i}-\frac{x+X}{2},y_{i}-%
\sum_{j=1}^{n}p_{j}y_{j}\right\rangle .
\end{multline}%
Taking the modulus in (\ref{2.4.16}) and using the Schwarz inequality in the
inner product space $\left( H;\left\langle \cdot ,\cdot \right\rangle
\right) ,$ we have%
\begin{align*}
& \left\vert \sum_{i=1}^{n}p_{i}\left\langle x_{i},y_{i}\right\rangle
-\left\langle \sum_{i=1}^{n}p_{i}x_{i},\sum_{i=1}^{n}p_{i}y_{i}\right\rangle
\right\vert \\
& \leq \sum_{i=1}^{n}p_{i}\left\vert \left\langle x_{i}-\frac{x+X}{2}%
,y_{i}-\sum_{j=1}^{n}p_{j}y_{j}\right\rangle \right\vert \\
& \leq \sum_{i=1}^{n}p_{i}\left\Vert x_{i}-\frac{x+X}{2}\right\Vert
\left\Vert y_{i}-\sum_{j=1}^{n}p_{j}y_{j}\right\Vert \\
& \leq \frac{1}{2}\left\Vert X-x\right\Vert \sum_{i=1}^{n}p_{i}\left\Vert
y_{i}-\sum_{j=1}^{n}p_{j}y_{j}\right\Vert ,
\end{align*}%
and the first inequality in (\ref{2.3.16}) is proved.

Using the Cauchy-Bunyakovsky-Schwarz inequality for positive sequences and
the calculation rules in inner product spaces, we have%
\begin{equation*}
\sum_{i=1}^{n}p_{i}\left\Vert y_{i}-\sum_{j=1}^{n}p_{j}y_{j}\right\Vert \leq 
\left[ \sum_{i=1}^{n}p_{i}\left\Vert
y_{i}-\sum_{j=1}^{n}p_{j}y_{j}\right\Vert ^{2}\right] ^{\frac{1}{2}}
\end{equation*}%
and%
\begin{equation*}
\sum_{i=1}^{n}p_{i}\left\Vert y_{i}-\sum_{j=1}^{n}p_{j}y_{j}\right\Vert
^{2}=\sum_{i=1}^{n}p_{i}\left\Vert y_{i}\right\Vert ^{2}-\left\Vert
\sum_{i=1}^{n}p_{i}y_{i}\right\Vert ^{2}
\end{equation*}%
giving the second part of (\ref{2.3.16}).

To prove the sharpness of the constant $\frac{1}{2}$ in the first inequality
in (\ref{2.3.16}), let us assume that, under the assumptions of the theorem,
the inequality holds with a constant $C>0,$ i.e., 
\begin{multline}
\left\vert \sum_{i=1}^{n}p_{i}\left\langle x_{i},y_{i}\right\rangle
-\left\langle \sum_{i=1}^{n}p_{i}x_{i},\sum_{i=1}^{n}p_{i}y_{i}\right\rangle
\right\vert  \label{2.4a.16} \\
\leq C\left\Vert X-x\right\Vert \sum_{i=1}^{n}p_{i}\left\Vert
y_{i}-\sum_{j=1}^{n}p_{j}y_{j}\right\Vert .
\end{multline}%
Consider $n=2$ and observe that%
\begin{align*}
\sum_{i=1}^{2}p_{i}\left\langle x_{i},y_{i}\right\rangle -\left\langle
\sum_{i=1}^{2}p_{i}x_{i},\sum_{i=1}^{2}p_{i}y_{i}\right\rangle &
=p_{2}p_{1}\left\langle x_{2}-x_{1},y_{2}-y_{1}\right\rangle , \\
\sum_{i=1}^{2}p_{i}\left\Vert y_{i}-\sum_{j=1}^{2}p_{j}y_{j}\right\Vert &
=2p_{2}p_{1}\left\Vert y_{2}-y_{1}\right\Vert
\end{align*}%
and then, by (\ref{2.4a.16}), we deduce%
\begin{equation}
p_{2}p_{1}\left\vert \left\langle x_{2}-x_{1},y_{2}-y_{1}\right\rangle
\right\vert \leq 2C\left\Vert X-x\right\Vert p_{2}p_{1}\left\Vert
y_{2}-y_{1}\right\Vert .  \label{2.5.16}
\end{equation}%
If we choose $p_{1},p_{2}>0,$ $y_{2}=x_{2},$ $y_{1}=x_{1}$ and $x_{2}=X,$ $%
x_{1}=x$ with $x\neq X,$ then (\ref{2.2.16}) holds and from (\ref{2.5.16})
we deduce $C\geq \frac{1}{2}.$

The fact that $\frac{1}{2}$ is best possible in the second inequality may be
proven in a similar manner and we omit the details.
\end{proof}

\begin{remark}
\label{r2.4.16}If $\overline{\mathbf{x}}$ and $\overline{\mathbf{y}}$
satisfy the assumptions of Theorem \ref{t2.3.16}, or equivalently%
\begin{equation}
\left\Vert x_{i}-\frac{x+X}{2}\right\Vert \leq \frac{1}{2}\left\Vert
X-x\right\Vert ,\ \ \ \ \left\Vert y_{i}-\frac{y+Y}{2}\right\Vert \leq \frac{%
1}{2}\left\Vert Y-y\right\Vert ,  \label{2.6.16}
\end{equation}%
for each $i\in \left\{ 1,\dots ,n\right\} ,$ then by Theorem \ref{t2.3.16}
we may state the following sequence of inequalities improving the Gr\"{u}ss
inequality (\ref{2.4.16})%
\begin{align}
0& \leq \left\vert \sum_{i=1}^{n}p_{i}\left\langle x_{i},y_{i}\right\rangle
-\left\langle \sum_{i=1}^{n}p_{i}x_{i},\sum_{i=1}^{n}p_{i}y_{i}\right\rangle
\right\vert  \label{2.7.16} \\
& \leq \frac{1}{2}\left\Vert X-x\right\Vert \sum_{i=1}^{n}p_{i}\left\Vert
y_{i}-\sum_{j=1}^{n}p_{j}y_{j}\right\Vert  \notag \\
& \leq \frac{1}{2}\left\Vert X-x\right\Vert \left(
\sum_{i=1}^{n}p_{i}\left\Vert y_{i}\right\Vert ^{2}-\left\Vert
\sum_{i=1}^{n}p_{i}y_{i}\right\Vert ^{2}\right) ^{\frac{1}{2}}  \notag \\
& \leq \frac{1}{4}\left\Vert X-x\right\Vert \left\Vert Y-y\right\Vert . 
\notag
\end{align}%
In particular, for $x_{i}=y_{i}$ $\left( i=1,\dots ,n\right) ,$ one has%
\begin{equation}
0\leq \sum_{i=1}^{n}p_{i}\left\Vert x_{i}\right\Vert ^{2}-\left\Vert
\sum_{i=1}^{n}p_{i}x_{i}\right\Vert ^{2}\leq \frac{1}{2}\left\Vert
X-x\right\Vert \sum_{i=1}^{n}p_{i}\left\Vert
x_{i}-\sum_{j=1}^{n}p_{j}x_{j}\right\Vert  \label{2.8.16}
\end{equation}%
and the constant $\frac{1}{2}$ is best possible.
\end{remark}

The following result also holds \cite{18NSSD}.

\begin{theorem}
\label{t2.5.16}Let $\left( H;\left\langle \cdot ,\cdot \right\rangle \right) 
$ and $\mathbb{K}$ be as above and $\overline{\mathbf{x}}=\left( x_{1},\dots
,x_{n}\right) \in H^{n},$ $\overline{\mathbf{\alpha }}=\left( \alpha
_{1},\dots ,\alpha _{n}\right) \in \mathbb{K}^{n}$ and $\overline{\mathbf{p}}
$ a probability vector. If $x,X\in H$ are such that (\ref{2.1.16}) or
equivalently, (\ref{2.2.16}) holds, then we have the inequality%
\begin{align}
0& \leq \left\Vert \sum_{i=1}^{n}p_{i}\alpha
_{i}x_{i}-\sum_{i=1}^{n}p_{i}\alpha _{i}\cdot
\sum_{i=1}^{n}p_{i}x_{i}\right\Vert  \label{2.9.16} \\
& \leq \frac{1}{2}\left\Vert X-x\right\Vert \sum_{i=1}^{n}p_{i}\left\vert
\alpha _{i}-\sum_{j=1}^{n}p_{j}\alpha _{j}\right\vert  \notag \\
& \leq \frac{1}{2}\left\Vert X-x\right\Vert \left[ \sum_{i=1}^{n}p_{i}\left%
\vert \alpha _{i}\right\vert ^{2}-\left\vert \sum_{i=1}^{n}p_{i}\alpha
_{i}\right\vert ^{2}\right] ^{\frac{1}{2}}.  \notag
\end{align}%
The constant $\frac{1}{2}$ in the first and second inequalities is best
possible in the sense that it cannot be replaced by a smaller constant.
\end{theorem}

\begin{proof}
We start with the following equality that may be easily verified by direct
calculation%
\begin{multline}
\sum_{i=1}^{n}p_{i}\alpha _{i}x_{i}-\sum_{i=1}^{n}p_{i}\alpha _{i}\cdot
\sum_{i=1}^{n}p_{i}x_{i}  \label{2.10.16} \\
=\sum_{i=1}^{n}p_{i}\left( \alpha _{i}-\sum_{j=1}^{n}p_{j}\alpha _{j}\right)
\left( x_{i}-\frac{x+X}{2}\right) .
\end{multline}%
If we take the norm in (\ref{2.10.16}), we deduce%
\begin{align*}
\left\Vert \sum_{i=1}^{n}p_{i}\alpha _{i}x_{i}-\sum_{i=1}^{n}p_{i}\alpha
_{i}\cdot \sum_{i=1}^{n}p_{i}x_{i}\right\Vert & \leq
\sum_{i=1}^{n}p_{i}\left\vert \alpha _{i}-\sum_{j=1}^{n}p_{j}\alpha
_{j}\right\vert \left\Vert x_{i}-\frac{x+X}{2}\right\Vert \\
& \leq \frac{1}{2}\left\Vert X-x\right\Vert \sum_{i=1}^{n}p_{i}\left\vert
\alpha _{i}-\sum_{j=1}^{n}p_{j}\alpha _{j}\right\vert \\
& \leq \frac{1}{2}\left\Vert X-x\right\Vert \left( \sum_{i=1}^{n}p_{i}\left(
\alpha _{i}-\sum_{j=1}^{n}p_{j}\alpha _{j}\right) ^{2}\right) ^{\frac{1}{2}}
\\
& =\frac{1}{2}\left\Vert X-x\right\Vert \left( \sum_{i=1}^{n}p_{i}\left\vert
\alpha _{i}\right\vert ^{2}-\left\vert \sum_{i=1}^{n}p_{i}\alpha
_{i}\right\vert ^{2}\right) ^{\frac{1}{2}},
\end{align*}%
proving the inequality (\ref{2.9.16}).

The fact that the constant $\frac{1}{2}$ is sharp may be proven in a similar
manner to the one embodied in the proof of Theorem \ref{t2.3.16}. We omit
the details.
\end{proof}

\begin{remark}
\label{r2.6.16}If $\overline{\mathbf{x}}$ and $\overline{\mathbf{\alpha }}$
satisfy the assumption 
\begin{equation*}
\left\Vert \alpha _{i}-\frac{a+A}{2}\right\Vert \leq \frac{1}{2}\left\vert
A-a\right\vert ,\ \ \ \ \left\Vert x_{i}-\frac{x+X}{2}\right\Vert \leq \frac{%
1}{2}\left\Vert X-x\right\Vert ,
\end{equation*}%
for each $i\in \left\{ 1,\dots ,n\right\} ,$ then by Theorem \ref{t2.3.16}
we may state the following sequence of inequalities improving the Gr\"{u}ss
inequality 
\begin{align*}
0& \leq \left\Vert \sum_{i=1}^{n}p_{i}\alpha
_{i}x_{i}-\sum_{i=1}^{n}p_{i}\alpha _{i}\cdot
\sum_{i=1}^{n}p_{i}x_{i}\right\Vert \\
& \leq \frac{1}{2}\left\Vert X-x\right\Vert \sum_{i=1}^{n}p_{i}\left\vert
\alpha _{i}-\sum_{j=1}^{n}p_{j}\alpha _{j}\right\vert \\
& \leq \frac{1}{2}\left\Vert X-x\right\Vert \left(
\sum_{i=1}^{n}p_{i}\left\vert \alpha _{i}\right\vert ^{2}-\left\vert
\sum_{i=1}^{n}p_{i}\alpha _{i}\right\vert ^{2}\right) ^{\frac{1}{2}} \\
& \leq \frac{1}{4}\left\vert A-a\right\vert \left\Vert X-x\right\Vert .
\end{align*}
\end{remark}

\begin{remark}
\label{r2.7.16}If in (\ref{2.9.16}) we choose $x_{i}=\alpha _{i}\in \mathbb{C%
}$ and assume that $\left\vert \alpha _{i}-\frac{a+A}{2}\right\vert \leq 
\frac{1}{2}\left\vert A-a\right\vert ,$ where $a,A\in \mathbb{C}$, then we
get the following interesting inequality for complex numbers%
\begin{align*}
0& \leq \left\vert \sum_{i=1}^{n}p_{i}\alpha _{i}^{2}-\left(
\sum_{i=1}^{n}p_{i}\alpha _{i}\right) ^{2}\right\vert \\
& \leq \frac{1}{2}\left\vert A-a\right\vert \sum_{i=1}^{n}p_{i}\left\vert
\alpha _{i}-\sum_{j=1}^{n}p_{j}\alpha _{j}\right\vert \\
& \leq \frac{1}{2}\left\vert A-a\right\vert \left[ \sum_{i=1}^{n}p_{i}\left%
\vert \alpha _{i}\right\vert ^{2}-\left\vert \sum_{i=1}^{n}p_{i}\alpha
_{i}\right\vert ^{2}\right] ^{\frac{1}{2}}.
\end{align*}
\end{remark}

\subsection{Applications for Convex Functions}

Let $\left( H;\left\langle \cdot ,\cdot \right\rangle \right) $ be a real
inner product space and $F:H\rightarrow \mathbb{R}$ a Fr\'{e}chet
differentiable convex function on $H.$ If $\triangledown F:H\rightarrow H$
denotes the gradient operator associated to $F,$ then we have the inequality%
\begin{equation}
F\left( x\right) -F\left( y\right) \geq \left\langle \triangledown F\left(
y\right) ,x-y\right\rangle  \label{3.1.16}
\end{equation}%
for each $x,y\in H.$

The following result holds \cite{18NSSD}.

\begin{theorem}
\label{t3.1.16}Let $F:H\rightarrow \mathbb{R}$ be as above and $z_{i}\in H,$ 
$i\in \left\{ 1,\dots ,n\right\} .$ Suppose that there exist the vectors $%
m,M\in H$ such that either%
\begin{equation*}
\left\langle \triangledown F\left( z_{i}\right) -m,M-\triangledown F\left(
z_{i}\right) \right\rangle \geq 0\text{ \ for each \ }i\in \left\{ 1,\dots
,n\right\} ;
\end{equation*}%
or equivalently,%
\begin{equation*}
\left\Vert \triangledown F\left( z_{i}\right) -\frac{m+M}{2}\right\Vert \leq 
\frac{1}{2}\left\Vert M-m\right\Vert \text{ \ for each \ }i\in \left\{
1,\dots ,n\right\} .
\end{equation*}%
If $q_{i}\geq 0$ \ $\left( i\in \left\{ 1,\dots ,n\right\} \right) $ with $%
Q_{n}:=\sum_{i=1}^{n}q_{i}>0,$ then we have the following converse of
Jensen's inequality%
\begin{align}
0& \leq \frac{1}{Q_{n}}\sum_{i=1}^{n}q_{i}F\left( z_{i}\right) -F\left( 
\frac{1}{Q_{n}}\sum_{i=1}^{n}q_{i}z_{i}\right)  \label{3.4.16} \\
& \leq \frac{1}{2}\left\Vert M-m\right\Vert \frac{1}{Q_{n}}%
\sum_{i=1}^{n}q_{i}\left\Vert z_{i}-\frac{1}{Q_{n}}\sum_{j=1}^{n}q_{j}z_{j}%
\right\Vert  \notag \\
& \leq \frac{1}{2}\left\Vert M-m\right\Vert \left[ \frac{1}{Q_{n}}%
\sum_{i=1}^{n}q_{i}\left\Vert z_{i}\right\Vert ^{2}-\left\Vert \frac{1}{Q_{n}%
}\sum_{i=1}^{n}q_{i}z_{i}\right\Vert ^{2}\right] ^{\frac{1}{2}}.  \notag
\end{align}
\end{theorem}

\begin{proof}
We know, see for example \cite[Eq. (4.4)]{D2.16}, that the following reverse
of Jensen's inequality for Fr\'{e}chet differentiable convex functions%
\begin{align}
0& \leq \frac{1}{Q_{n}}\sum_{i=1}^{n}q_{i}F\left( z_{i}\right) -F\left( 
\frac{1}{Q_{n}}\sum_{i=1}^{n}q_{i}z_{i}\right)  \label{3.5.16} \\
& \leq \frac{1}{Q_{n}}\sum_{i=1}^{n}q_{i}\left\langle \triangledown F\left(
z_{i}\right) ,z_{i}\right\rangle -\left\langle \frac{1}{Q_{n}}%
\sum_{i=1}^{n}q_{i}\triangledown F\left( z_{i}\right) ,\frac{1}{Q_{n}}%
\sum_{i=1}^{n}q_{i}z_{i}\right\rangle  \notag
\end{align}%
holds.

Now, if we use Theorem \ref{t2.3.16} for the choices $x_{i}=\nabla F\left(
z_{i}\right) ,$ $y_{i}=z_{i}$ and $p_{i}=\frac{1}{Q_{n}}q_{i},$ $i\in
\left\{ 1,\dots ,n\right\} ,$ we can state the inequality%
\begin{align}
& \frac{1}{Q_{n}}\sum_{i=1}^{n}q_{i}\left\langle \triangledown F\left(
z_{i}\right) ,z_{i}\right\rangle -\left\langle \frac{1}{Q_{n}}%
\sum_{i=1}^{n}q_{i}\triangledown F\left( z_{i}\right) ,\frac{1}{Q_{n}}%
\sum_{i=1}^{n}q_{i}z_{i}\right\rangle  \label{3.6.16} \\
& \leq \frac{1}{2}\left\Vert M-m\right\Vert \frac{1}{Q_{n}}%
\sum_{i=1}^{n}q_{i}\left\Vert z_{i}-\frac{1}{Q_{n}}\sum_{j=1}^{n}q_{j}z_{j}%
\right\Vert  \notag \\
& \leq \frac{1}{2}\left\Vert M-m\right\Vert \left[ \frac{1}{Q_{n}}%
\sum_{i=1}^{n}q_{i}\left\Vert z_{i}\right\Vert ^{2}-\left\Vert \frac{1}{Q_{n}%
}\sum_{i=1}^{n}q_{i}z_{i}\right\Vert ^{2}\right] ^{\frac{1}{2}}.  \notag
\end{align}%
Utilizing (\ref{3.5.16}) and (\ref{3.6.16}), we deduce the desired result (%
\ref{3.4.16}).
\end{proof}

If more information is available about the vector sequence $\overline{%
\mathbf{z}}=\left( z_{1},\dots ,z_{n}\right) \in H^{n},$ then we may state
the following corollary.

\begin{corollary}
\label{c3.2.16}With the assumptions in Theorem \ref{t3.1.16} and if there
exist the vectors $z,Z\in H$ such that either%
\begin{equation}
\left\langle z_{i}-z,Z-z_{i}\right\rangle \geq 0\text{ \ for each \ }i\in
\left\{ 1,\dots ,n\right\} ;  \label{3.7.16}
\end{equation}%
or, equivalently%
\begin{equation}
\left\Vert z_{i}-\frac{z+Z}{2}\right\Vert \leq \frac{1}{2}\left\Vert
Z-z\right\Vert \text{ \ for each \ }i\in \left\{ 1,\dots ,n\right\} ,
\label{3.8.16}
\end{equation}%
then we have the inequality%
\begin{align}
0& \leq \frac{1}{Q_{n}}\sum_{i=1}^{n}q_{i}F\left( z_{i}\right) -F\left( 
\frac{1}{Q_{n}}\sum_{i=1}^{n}q_{i}z_{i}\right)  \label{3.9.16} \\
& \leq \frac{1}{2}\left\Vert M-m\right\Vert \frac{1}{Q_{n}}%
\sum_{i=1}^{n}q_{i}\left\Vert z_{i}-\frac{1}{Q_{n}}\sum_{j=1}^{n}q_{j}z_{j}%
\right\Vert  \notag \\
& \leq \frac{1}{2}\left\Vert M-m\right\Vert \left[ \frac{1}{Q_{n}}%
\sum_{i=1}^{n}q_{i}\left\Vert z_{i}\right\Vert ^{2}-\left\Vert \frac{1}{Q_{n}%
}\sum_{i=1}^{n}q_{i}z_{i}\right\Vert ^{2}\right] ^{\frac{1}{2}}  \notag \\
& \leq \frac{1}{4}\left\Vert M-m\right\Vert \left\Vert Z-z\right\Vert . 
\notag
\end{align}
\end{corollary}

\begin{remark}
Note that the inequality between the first term and the last term in (\ref%
{3.9.16}) was first proved in \cite[Theorem 4.1]{D2.16}. Consequently, the
above corollary provides an improvement of the reverse of Jensen's
inequality established in \cite{D2.16}.
\end{remark}

\newpage

\section{Some Inequalities for Forward Difference}

\subsection{Introduction}

In \cite{ab.16.0}, we have proved the following generalisation of the Gr\"{u}%
ss inequality.

\begin{theorem}
\label{ta.16.0} Let $\left( H,\left\langle \cdot ,\cdot \right\rangle
\right) $ be an inner product space over $\mathbb{K},$ $\mathbb{K}=\mathbb{C}%
,\mathbb{R}$ and $e\in H,$ $\left\Vert e\right\Vert =1.$ If $\phi ,\Phi
,\gamma ,\Gamma \in \mathbb{K}$ and $x,y\in H$ are such that%
\begin{equation*}
\func{Re}\left\langle \Phi e-x,x-\phi e\right\rangle \geq 0\text{ \ and \ \ }%
\func{Re}\left\langle \Gamma e-y,y-\gamma e\right\rangle \geq 0
\end{equation*}%
hold, then we have the inequality%
\begin{equation*}
\left\vert \left\langle x,y\right\rangle -\left\langle x,e\right\rangle
\left\langle e,y\right\rangle \right\vert \leq \frac{1}{4}\left\vert \Phi
-\phi \right\vert \left\vert \Gamma -\gamma \right\vert .
\end{equation*}%
The constant $\frac{1}{4}$ is the best possible.
\end{theorem}

A Gr\"{u}ss type inequality for sequences of vectors in inner product spaces
was pointed out in \cite{bb.16.0}.

\begin{theorem}
\label{tb.16.0}Let $H$ and $\mathbb{K}$ be as in Theorem \ref{ta.16.0} and $%
x_{i}\in H,$ $a_{i}\in \mathbb{K},$ $p_{i}\geq 0$ \ $\left( i=1,\dots
,n\right) $ $\left( n\geq 2\right) $ with $\sum_{i=1}^{n}p_{i}=1$. If $%
a,A\in \mathbb{K}$ and $x,X\in H$ are such that:%
\begin{equation*}
\func{Re}\left[ \left( A-a_{i}\right) \left( \overline{a_{i}}-\overline{a}%
\right) \right] \geq 0,\ \ \ \func{Re}\left\langle
X-x_{i},x_{i}-x\right\rangle \geq 0
\end{equation*}%
for any $i\in \left\{ 1,\dots ,n\right\} ,$ then we have the inequality%
\begin{align*}
0& \leq \left\Vert
\sum_{i=1}^{n}p_{i}a_{i}x_{i}-\sum_{i=1}^{n}p_{i}a_{i}\cdot
\sum_{i=1}^{n}p_{i}x_{i}\right\Vert \\
& \leq \frac{1}{4}\left\vert A-a\right\vert \left\Vert X-x\right\Vert .
\end{align*}%
The constant $\frac{1}{4}$ is best possible.
\end{theorem}

A complementary result for two sequences of vectors in inner product spaces
is the following result that has been obtained in \cite{cb.16.0}.

\begin{theorem}
\label{tc.16.0}Let $H$ and $\mathbb{K}$ be as above, $x_{i},y_{i}\in H,$ $%
p_{i}\geq 0$ \ $\left( i=1,\dots ,n\right) $ $\left( n\geq 2\right) $ with $%
\sum_{i=1}^{n}p_{i}=1$. If $x,X,y,Y\in H$ are such that:%
\begin{equation*}
\func{Re}\left\langle X-x_{i},x_{i}-x\right\rangle \geq 0\text{ \ \ and \ \ }%
\func{Re}\left\langle Y-y_{i},y_{i}-y\right\rangle \geq 0\text{ \ }
\end{equation*}%
for all \ $i\in \left\{ 1,\dots ,n\right\} ,$ then we have the inequality%
\begin{equation*}
0\leq \left\vert \sum_{i=1}^{n}p_{i}\left\langle x_{i},y_{i}\right\rangle
-\left\langle \sum_{i=1}^{n}p_{i}x_{i},\sum_{i=1}^{n}p_{i}y_{i}\right\rangle
\right\vert \leq \frac{1}{4}\left\Vert X-x\right\Vert \left\Vert
Y-y\right\Vert .
\end{equation*}%
The constant $\frac{1}{4}$ is best possible.
\end{theorem}

In the general case of normed linear spaces, the following Gr\"{u}ss type
inequality in terms of the forward difference is known, see \cite{db.16.0}.

\begin{theorem}
\label{td.16.0}Let $\left( E,\left\Vert \cdot \right\Vert \right) $ be a
normed linear space over $\mathbb{K}=\mathbb{C},\mathbb{R},$ $x_{i}\in E,$ $%
\alpha _{i}\in \mathbb{K}$ and $p_{i}\geq 0$ \ $\left( i=1,\dots ,n\right) $
such that $\sum_{i=1}^{n}p_{i}=1$. Then we have the inequality%
\begin{eqnarray}
0 &\leq &\left\Vert \sum_{i=1}^{n}p_{i}\alpha
_{i}x_{i}-\sum_{i=1}^{n}p_{i}\alpha _{i}\cdot
\sum_{i=1}^{n}p_{i}x_{i}\right\Vert  \label{1.7.16.0} \\
&\leq &\max_{1\leq j\leq n-1}\left\vert \Delta \alpha _{j}\right\vert
\max_{1\leq j\leq n-1}\left\Vert \Delta x_{j}\right\Vert \left[
\sum_{i=1}^{n}i^{2}p_{i}-\left( \sum_{i=1}^{n}ip_{i}\right) ^{2}\right] , 
\notag
\end{eqnarray}%
where $\Delta \alpha _{j}=\alpha _{j+1}-\alpha _{j}$ and $\Delta
x_{j}=x_{j+1}-x_{j}\left( j=1,\dots ,n-1\right) $ are the forward
differences of the vectors having the components $\alpha _{j}$ and $%
x_{j}\left( j=1,\dots ,n-1\right) ,$ respectively.

The inequality (\ref{1.7.16.0}) is sharp in the sense that the
multiplicative constant $C=1$ in the right hand side cannot be replaced by a
smaller one.
\end{theorem}

An important particular case is the one where all the weights are equal,
giving the following corollary \cite{db.16.0}.

\begin{corollary}
\label{ce.16.0}Under the above assumptions for $\alpha _{i},x_{i}$ $\left(
i=1,\dots ,n\right) $ we have the inequality%
\begin{align}
0& \leq \left\Vert \frac{1}{n}\sum_{i=1}^{n}\alpha _{i}x_{i}-\frac{1}{n}%
\sum_{i=1}^{n}\alpha _{i}\cdot \frac{1}{n}\sum_{i=1}^{n}x_{i}\right\Vert
\label{1.8.16.0} \\
& \leq \frac{n^{2}-1}{12}\max_{1\leq j\leq n-1}\left\vert \Delta \alpha
_{j}\right\vert \max_{1\leq j\leq n-1}\left\Vert \Delta x_{j}\right\Vert . 
\notag
\end{align}%
The constant $\frac{1}{12}$ is best possible.
\end{corollary}

Another result of this type was proved in \cite{e1b.16.0}.

\begin{theorem}
\label{tf.16.0}With the assumptions of Theorem \ref{td.16.0}, one has the
inequality%
\begin{align}
0& \leq \left\Vert \sum_{i=1}^{n}p_{i}\alpha
_{i}x_{i}-\sum_{i=1}^{n}p_{i}\alpha _{i}\cdot
\sum_{i=1}^{n}p_{i}x_{i}\right\Vert  \label{1.9.16.0} \\
& \leq \frac{1}{2}\sum_{j=1}^{n-1}\left\vert \Delta \alpha _{j}\right\vert
\sum_{j=1}^{n-1}\left\Vert \Delta x_{j}\right\Vert \sum_{i=1}^{n}p_{i}\left(
1-p_{i}\right) .  \notag
\end{align}%
The constant $\frac{1}{2}$ is best possible.
\end{theorem}

As a useful particular case, we have the following corollary \cite{e1b.16.0}.

\begin{corollary}
\label{ch.16.0}If $\alpha _{i},x_{i}$ $\left( i=1,\dots ,n\right) $ are as
in Theorem \ref{td.16.0}, then%
\begin{align*}
0& \leq \left\Vert \frac{1}{n}\sum_{i=1}^{n}\alpha _{i}x_{i}-\frac{1}{n}%
\sum_{i=1}^{n}\alpha _{i}\cdot \frac{1}{n}\sum_{i=1}^{n}x_{i}\right\Vert \\
& \leq \frac{1}{2}\left( 1-\frac{1}{n}\right) \sum_{i=1}^{n-1}\left\vert
\Delta \alpha _{i}\right\vert \sum_{i=1}^{n-1}\left\Vert \Delta
x_{i}\right\Vert .
\end{align*}%
The constant $\frac{1}{2}$ is the best possible.
\end{corollary}

Finally, the following result is also known \cite{eb.16.0}.

\begin{theorem}
\label{tg.16.0}With the assumptions in Theorem \ref{td.16.0}, we have the
inequality:%
\begin{align}
0& \leq \left\Vert \sum_{i=1}^{n}p_{i}\alpha
_{i}x_{i}-\sum_{i=1}^{n}p_{i}\alpha _{i}\cdot
\sum_{i=1}^{n}p_{i}x_{i}\right\Vert  \label{1.11.16.0} \\
& \leq \left( \sum_{j=1}^{n-1}\left\vert \Delta \alpha _{j}\right\vert
^{p}\right) ^{\frac{1}{p}}\left( \sum_{j=1}^{n-1}\left\Vert \Delta
x_{j}\right\Vert ^{q}\right) ^{\frac{1}{q}}\sum_{1\leq i<j\leq n}\left(
j-i\right) p_{i}p_{j},  \notag
\end{align}%
where $p>1,$ $\frac{1}{p}+\frac{1}{q}=1.$

The constant $c=1$ in the right hand side of (\ref{1.11.16.0}) is sharp.
\end{theorem}

The case of equal weights is embodied in the following\ corollary \cite%
{eb.16.0}.

\begin{corollary}
\label{cl.16.0}With the above assumptions for $\alpha _{i},x_{i}$ $\left(
i=1,\dots ,n\right) $ one has%
\begin{align*}
0& \leq \left\Vert \frac{1}{n}\sum_{i=1}^{n}\alpha _{i}x_{i}-\frac{1}{n}%
\sum_{i=1}^{n}\alpha _{i}\cdot \frac{1}{n}\sum_{i=1}^{n}x_{i}\right\Vert \\
& \leq \frac{n^{2}-1}{6n}\left( \sum_{j=1}^{n-1}\left\vert \Delta \alpha
_{j}\right\vert ^{p}\right) ^{\frac{1}{p}}\left( \sum_{j=1}^{n-1}\left\Vert
\Delta x_{j}\right\Vert ^{q}\right) ^{\frac{1}{q}},
\end{align*}%
where $p>1,$ $\frac{1}{p}+\frac{1}{q}=1.$

The constant $\frac{1}{6}$ is the best possible.
\end{corollary}

The main aim of this section is to establish some similar bounds for the
absolute value of the difference%
\begin{equation*}
\sum_{i=1}^{n}p_{i}\left\langle x_{i},y_{i}\right\rangle -\left\langle
\sum_{i=1}^{n}p_{i}x_{i},\sum_{i=1}^{n}p_{i}y_{i}\right\rangle
\end{equation*}%
provided that $x_{i},y_{i}$ $\left( i=1,\dots ,n\right) $ are vectors in an
inner product space $H,$ and $p_{i}\geq 0$ $\left( i=1,\dots ,n\right) $
with $\sum_{i=1}^{n}p_{i}=1.$

\subsection{The Main Results}

We assume that $\left( H,\left\langle \cdot ,\cdot \right\rangle \right) $
is an inner product space over $\mathbb{K},$ $\mathbb{K}=\mathbb{C}$ or $%
\mathbb{K}=\mathbb{R}.$ The following discrete inequality of Gr\"{u}ss' type
holds.

\begin{theorem}
\label{t2.1.16.0}If $x_{i},y_{i}\in H,$ $p_{i}\geq 0$ $\left( i=1,\dots
,n\right) $ with $\sum_{i=1}^{n}p_{i}=1$, then one has the inequalities:%
\begin{align}
& \left\vert \sum_{i=1}^{n}p_{i}\left\langle x_{i},y_{i}\right\rangle
-\left\langle \sum_{i=1}^{n}p_{i}x_{i},\sum_{i=1}^{n}p_{i}y_{i}\right\rangle
\right\vert  \label{2.1.16.0} \\
& \leq \left\{ 
\begin{array}{l}
\left[ \sum\limits_{i=1}^{n}i^{2}p_{i}-\left(
\sum\limits_{i=1}^{n}ip_{i}\right) ^{2}\right] \max\limits_{k=1,\dots
,n-1}\left\Vert \Delta x_{k}\right\Vert \max\limits_{k=1,\dots
,n-1}\left\Vert \Delta y_{k}\right\Vert ; \\ 
\\ 
\left[ \sum\limits_{1\leq j<i\leq n}p_{i}p_{j}\left( i-j\right) \right]
\left( \sum\limits_{k=1}^{n-1}\left\Vert \Delta x_{k}\right\Vert ^{p}\right)
^{\frac{1}{p}}\left( \sum\limits_{k=1}^{n-1}\left\Vert \Delta
y_{k}\right\Vert ^{q}\right) ^{\frac{1}{q}} \\ 
\hfill \text{if }p>1,\ \frac{1}{p}+\frac{1}{q}=1; \\ 
\\ 
\dfrac{1}{2}\left[ \sum\limits_{i=1}^{n}p_{i}\left( 1-p_{i}\right) \right]
\sum\limits_{k=1}^{n-1}\left\Vert \Delta x_{k}\right\Vert
\sum\limits_{k=1}^{n-1}\left\Vert \Delta y_{k}\right\Vert .%
\end{array}%
\right.  \notag
\end{align}%
All the inequalities in (\ref{2.1.16.0}) are sharp.
\end{theorem}

The following particular case for equal vectors holds.

\begin{corollary}
\label{c2.1.a.16.0}With the assumptions of Theorem \ref{t2.1.16.0}, one has
the inequalities%
\begin{align*}
0& \leq \sum_{i=1}^{n}p_{i}\left\Vert x_{i}\right\Vert ^{2}-\left\Vert
\sum_{i=1}^{n}p_{i}x_{i}\right\Vert ^{2} \\
& \leq \left\{ 
\begin{array}{l}
\left[ \sum\limits_{i=1}^{n}i^{2}p_{i}-\left(
\sum\limits_{i=1}^{n}ip_{i}\right) ^{2}\right] \max\limits_{k=1,\dots
,n-1}\left\Vert \Delta x_{k}\right\Vert ^{2}; \\ 
\\ 
\sum\limits_{1\leq j<i\leq n}p_{i}p_{j}\left( i-j\right) \left(
\sum\limits_{k=1}^{n-1}\left\Vert \Delta x_{k}\right\Vert ^{p}\right) ^{%
\frac{1}{p}}\left( \sum\limits_{k=1}^{n-1}\left\Vert \Delta x_{k}\right\Vert
^{q}\right) ^{\frac{1}{q}} \\ 
\hfill \text{if }p>1,\ \frac{1}{p}+\frac{1}{q}=1; \\ 
\\ 
\dfrac{1}{2}\sum\limits_{i=1}^{n}p_{i}\left( 1-p_{i}\right) \left(
\sum\limits_{k=1}^{n-1}\left\Vert \Delta x_{k}\right\Vert \right) ^{2}.%
\end{array}%
\right.
\end{align*}
\end{corollary}

The following particular case for equal weights may be useful in practice.

\begin{corollary}
\label{c2.2.16.0}If $x_{i},y_{i}\in H$ $\left( i=1,\dots ,n\right) $, then
one has the inequalities:%
\begin{equation*}
\left\vert \frac{1}{n}\sum_{i=1}^{n}\left\langle x_{i},y_{i}\right\rangle
-\left\langle \frac{1}{n}\sum_{i=1}^{n}x_{i},\frac{1}{n}\sum_{i=1}^{n}y_{i}%
\right\rangle \right\vert
\end{equation*}%
\begin{equation*}
\leq \left\{ 
\begin{array}{l}
\dfrac{n^{2}-1}{12}\max\limits_{k=1,\dots ,n-1}\left\Vert \Delta
x_{k}\right\Vert \max\limits_{k=1,\dots ,n-1}\left\Vert \Delta
y_{k}\right\Vert ; \\ 
\\ 
\dfrac{n^{2}-1}{6n}\left( \sum\limits_{k=1}^{n-1}\left\Vert \Delta
x_{k}\right\Vert ^{p}\right) ^{\frac{1}{p}}\left(
\sum\limits_{k=1}^{n-1}\left\Vert \Delta y_{k}\right\Vert ^{q}\right) ^{%
\frac{1}{q}} \\ 
\hfill \text{if }p>1,\ \frac{1}{p}+\frac{1}{q}=1; \\ 
\\ 
\dfrac{n-1}{2n}\sum\limits_{k=1}^{n-1}\left\Vert \Delta x_{k}\right\Vert
\sum\limits_{k=1}^{n-1}\left\Vert \Delta y_{k}\right\Vert .%
\end{array}%
\right.
\end{equation*}%
The constants $\frac{1}{12}$, $\frac{1}{6}$ and $\frac{1}{2}$ are best
possible.
\end{corollary}

In particular, the following corollary holds.

\begin{corollary}
\label{c2.3.16.0}If $x_{i}\in H$ $\left( i=1,\dots ,n\right) ,$ then one has
the inequality%
\begin{align*}
0& \leq \frac{1}{n}\sum_{i=1}^{n}\left\Vert x_{i}\right\Vert ^{2}-\left\Vert 
\frac{1}{n}\sum_{i=1}^{n}x_{i}\right\Vert ^{2} \\
& \leq \left\{ 
\begin{array}{l}
\dfrac{n^{2}-1}{12}\max\limits_{k=\overline{1,n}}\left\Vert \Delta
x_{k}\right\Vert ^{2}; \\ 
\\ 
\dfrac{n^{2}-1}{6n}\left( \sum\limits_{k=1}^{n-1}\left\Vert \Delta
x_{k}\right\Vert ^{p}\right) ^{\frac{1}{p}}\left(
\sum\limits_{k=1}^{n-1}\left\Vert \Delta x_{k}\right\Vert ^{q}\right) ^{%
\frac{1}{q}} \\ 
\hfill \text{if }p>1,\ \frac{1}{p}+\frac{1}{q}=1; \\ 
\\ 
\dfrac{n-1}{2n}\left( \sum\limits_{k=1}^{n-1}\left\Vert \Delta
x_{k}\right\Vert \right) ^{2}.%
\end{array}%
\right.
\end{align*}%
The constants $\frac{1}{12}$, $\frac{1}{6}$ and $\frac{1}{2}$ are best
possible.
\end{corollary}

\subsection{Proof of the Main Result}

It is well known that, the following identity holds in inner product spaces:%
\begin{align}
& \sum_{i=1}^{n}p_{i}\left\langle x_{i},y_{i}\right\rangle -\left\langle
\sum_{i=1}^{n}p_{i}x_{i},\sum_{i=1}^{n}p_{i}y_{i}\right\rangle
\label{3.1.16.0} \\
& =\frac{1}{2}\sum_{i,j=1}^{n}p_{i}p_{j}\left\langle
x_{i}-x_{j},y_{i}-y_{j}\right\rangle  \notag \\
& =\sum\limits_{1\leq j<i\leq n}p_{i}p_{j}\left\langle
x_{i}-x_{j},y_{i}-y_{j}\right\rangle .  \notag
\end{align}%
We observe, for $i>j,$ we can write that 
\begin{equation}
x_{i}-x_{j}=\sum\limits_{k=j}^{i-1}\Delta x_{k},\ \ \ \ \ \ \
y_{i}-y_{j}=\sum\limits_{k=j}^{i-1}\Delta y_{k}.  \label{3.2.16.0}
\end{equation}%
Taking the modulus in (\ref{3.1.16.0}) and by the use of (\ref{3.2.16.0})
and Schwarz's inequality in inner product spaces, i.e., we recall that $%
\left\vert \left\langle z,u\right\rangle \right\vert \leq \left\Vert
z\right\Vert \left\Vert u\right\Vert ,$ $z,u\in H,$ we have:%
\begin{align*}
& \left\vert \sum_{i=1}^{n}p_{i}\left\langle x_{i},y_{i}\right\rangle
-\left\langle \sum_{i=1}^{n}p_{i}x_{i},\sum_{i=1}^{n}p_{i}y_{i}\right\rangle
\right\vert \\
& \leq \sum\limits_{1\leq j<i\leq n}p_{i}p_{j}\left\vert \left\langle
x_{i}-x_{j},y_{i}-y_{j}\right\rangle \right\vert \\
& \leq \sum\limits_{1\leq j<i\leq n}p_{i}p_{j}\left\Vert
x_{i}-x_{j}\right\Vert \left\Vert y_{i}-y_{j}\right\Vert \\
& =\sum\limits_{1\leq j<i\leq n}p_{i}p_{j}\left\Vert
\sum\limits_{k=j}^{i-1}\Delta x_{k}\right\Vert \left\Vert
\sum\limits_{l=j}^{i-1}\Delta y_{l}\right\Vert \\
& \leq \sum\limits_{1\leq j<i\leq
n}p_{i}p_{j}\sum\limits_{k=j}^{i-1}\left\Vert \Delta x_{k}\right\Vert
\sum\limits_{l=j}^{i-1}\left\Vert \Delta y_{l}\right\Vert \\
& :=M.
\end{align*}%
It is obvious that%
\begin{equation*}
\sum\limits_{k=j}^{i-1}\left\Vert \Delta x_{k}\right\Vert \leq \left(
i-j\right) \max\limits_{k=j,\dots ,i-1}\left\Vert \Delta x_{k}\right\Vert
\leq \left( i-j\right) \max\limits_{k=1,\dots ,n}\left\Vert \Delta
x_{k}\right\Vert
\end{equation*}%
and%
\begin{equation*}
\sum\limits_{k=j}^{i-1}\left\Vert \Delta y_{k}\right\Vert \leq \left(
i-j\right) \max\limits_{k=j,\dots ,i-1}\left\Vert \Delta y_{k}\right\Vert
\leq \left( i-j\right) \max\limits_{k=1,\dots ,n}\left\Vert \Delta
y_{k}\right\Vert ,
\end{equation*}%
giving that%
\begin{equation*}
M\leq \sum\limits_{1\leq j<i\leq n}p_{i}p_{j}\left( i-j\right) ^{2}\cdot
\max\limits_{k=1,\dots ,n}\left\Vert \Delta x_{k}\right\Vert
\max\limits_{k=1,\dots ,n}\left\Vert \Delta y_{k}\right\Vert ,
\end{equation*}%
and since%
\begin{equation*}
\sum\limits_{1\leq j<i\leq n}p_{i}p_{j}\left( i-j\right) ^{2}=\frac{1}{2}%
\sum_{i,j=1}^{n}p_{i}p_{j}\left( i-j\right)
^{2}=\sum\limits_{i=1}^{n}p_{i}i^{2}-\left(
\sum\limits_{i=1}^{n}ip_{i}\right) ^{2},
\end{equation*}%
the first inequality in (\ref{2.1.16.0}) is proved.

Using H\"{o}lder's discrete inequality, we can state that%
\begin{equation*}
\sum\limits_{k=j}^{i-1}\left\Vert \Delta x_{k}\right\Vert \leq \left(
i-j\right) ^{\frac{1}{q}}\left( \sum\limits_{k=j}^{i-1}\left\Vert \Delta
x_{k}\right\Vert ^{p}\right) ^{\frac{1}{p}}\leq \left( i-j\right) ^{\frac{1}{%
q}}\left( \sum\limits_{k=1}^{n-1}\left\Vert \Delta x_{k}\right\Vert
^{p}\right) ^{\frac{1}{p}}
\end{equation*}%
and%
\begin{equation*}
\sum\limits_{k=j}^{i-1}\left\Vert \Delta y_{k}\right\Vert \leq \left(
i-j\right) ^{\frac{1}{p}}\left( \sum\limits_{k=j}^{i-1}\left\Vert \Delta
y_{k}\right\Vert ^{q}\right) ^{\frac{1}{q}}\leq \left( i-j\right) ^{\frac{1}{%
p}}\left( \sum\limits_{k=1}^{n-1}\left\Vert \Delta y_{k}\right\Vert
^{q}\right) ^{\frac{1}{q}},
\end{equation*}%
for $p>1,$ $\frac{1}{p}+\frac{1}{q}=1,$ giving that:%
\begin{equation*}
M\leq \sum\limits_{1\leq j<i\leq n}p_{i}p_{j}\left( i-j\right) \cdot \left(
\sum\limits_{k=1}^{n-1}\left\Vert \Delta x_{k}\right\Vert ^{p}\right) ^{%
\frac{1}{p}}\left( \sum\limits_{k=1}^{n-1}\left\Vert \Delta y_{k}\right\Vert
^{q}\right) ^{\frac{1}{q}}
\end{equation*}%
and the second inequality in (\ref{2.1.16.0}) is proved.

Also, observe that%
\begin{equation*}
\sum\limits_{k=j}^{i-1}\left\Vert \Delta x_{k}\right\Vert \leq
\sum\limits_{k=1}^{n-1}\left\Vert \Delta x_{k}\right\Vert \text{ \ and \ }%
\sum\limits_{k=j}^{i-1}\left\Vert \Delta y_{k}\right\Vert \leq
\sum\limits_{k=1}^{n-1}\left\Vert \Delta y_{k}\right\Vert
\end{equation*}%
and thus%
\begin{equation*}
M\leq \sum\limits_{1\leq j<i\leq n}p_{i}p_{j}\left( i-j\right)
\sum\limits_{k=1}^{n-1}\left\Vert \Delta x_{k}\right\Vert
\sum\limits_{k=1}^{n-1}\left\Vert \Delta y_{k}\right\Vert .
\end{equation*}%
Since%
\begin{align*}
\sum\limits_{1\leq j<i\leq n}p_{i}p_{j}& =\frac{1}{2}\left[
\sum_{i,j=1}^{n}p_{i}p_{j}-\sum\limits_{k=1}^{n}p_{k}^{2}\right] \\
& =\frac{1}{2}\left( 1-\sum\limits_{k=1}^{n}p_{k}^{2}\right) \\
& =\frac{1}{2}\sum\limits_{i=1}^{n}p_{i}\left( 1-p_{i}\right) ,
\end{align*}%
the last part of (\ref{2.1.16.0}) is also proved.

Now, assume that the first inequality in (\ref{2.1.16.0}) holds with a
constant $c>0,$ i.e.,%
\begin{align*}
& \sum_{i=1}^{n}p_{i}\left\langle x_{i},y_{i}\right\rangle -\left\langle
\sum_{i=1}^{n}p_{i}x_{i},\sum_{i=1}^{n}p_{i}y_{i}\right\rangle \\
& \leq c\left[ \sum\limits_{i=1}^{n}i^{2}p_{i}-\left(
\sum\limits_{i=1}^{n}ip_{i}\right) ^{2}\right] \max\limits_{k=1,\dots
,n-1}\left\Vert \Delta x_{k}\right\Vert \max\limits_{k=1,\dots
,n-1}\left\Vert \Delta y_{k}\right\Vert
\end{align*}%
and choose $n=2$ to get%
\begin{equation}
p_{1}p_{2}\left\vert \left\langle x_{2}-x_{1},y_{2}-y_{1}\right\rangle
\right\vert \leq cp_{1}p_{2}\left\Vert x_{2}-x_{1}\right\Vert \left\Vert
y_{2}-y_{1}\right\Vert  \label{3.5.16.0}
\end{equation}%
for any $p_{1},p_{2}>0$ and $x_{1},x_{2},y_{1},y_{2}\in H.$

If in (\ref{3.5.16.0}) we choose $y_{2}=x_{2},$ $y_{1}=x_{1}$ and $x_{2}\neq
x_{1},$ then we deduce $c\geq 1,$ which proves the sharpness of the constant
in the first inequality in (\ref{2.1.16.0}).

In a similar way one may show that the other two inequalities are sharp, and
the theorem is completely proved.

\subsection{A Reverse for Jensen's Inequality}

Let $\left( H;\left\langle \cdot ,\cdot \right\rangle \right) $ be a real
inner product space and $F:H\rightarrow \mathbb{R}$ a Fr\'{e}chet
differentiable convex function on $H.$ If $\triangledown F:H\rightarrow H$
denotes the gradient operator associated to $F,$ then we have the inequality%
\begin{equation*}
F\left( x\right) -F\left( y\right) \geq \left\langle \triangledown F\left(
y\right) ,x-y\right\rangle
\end{equation*}%
for each $x,y\in H.$

The following result holds.

\begin{theorem}
\label{t4.1.16.0}Let $F:H\rightarrow \mathbb{R}$ be as above and $z_{i}\in
H, $ $i\in \left\{ 1,\dots ,n\right\} .$ If $q_{i}\geq 0$ \ $\left( i\in
\left\{ 1,\dots ,n\right\} \right) $ with $\sum_{i=1}^{n}q_{i}=1,$ then we
have the following reverse of Jensen's inequality%
\begin{align}
0& \leq \sum_{i=1}^{n}q_{i}F\left( z_{i}\right) -F\left(
\sum_{i=1}^{n}q_{i}z_{i}\right)  \label{4.2.16.0} \\
& \leq \left\{ 
\begin{array}{l}
\left[ \sum\limits_{i=1}^{n}i^{2}q_{i}-\left(
\sum\limits_{i=1}^{n}iq_{i}\right) ^{2}\right] \max\limits_{k=1,\dots
,n-1}\left\Vert \Delta \left( \triangledown F\left( z_{i}\right) \right)
\right\Vert \max\limits_{k=1,\dots ,n-1}\left\Vert \Delta z_{i}\right\Vert ;
\\ 
\\ 
\left[ \sum\limits_{1\leq j<i\leq n}q_{i}q_{j}\left( i-j\right) \right]
\left( \sum\limits_{i=1}^{n-1}\left\Vert \Delta \left( \triangledown F\left(
z_{i}\right) \right) \right\Vert ^{p}\right) ^{\frac{1}{p}}\left(
\sum\limits_{i=1}^{n-1}\left\Vert \Delta z_{i}\right\Vert ^{q}\right) ^{%
\frac{1}{q}} \\ 
\hfill \text{if }p>1,\ \frac{1}{p}+\frac{1}{q}=1; \\ 
\\ 
\dfrac{1}{2}\left[ \sum\limits_{i=1}^{n}q_{i}\left( 1-q_{i}\right) \right]
\sum\limits_{i=1}^{n-1}\left\Vert \Delta \left( \triangledown F\left(
z_{i}\right) \right) \right\Vert \sum\limits_{i=1}^{n-1}\left\Vert \Delta
z_{i}\right\Vert .%
\end{array}%
\right.  \notag \\
&  \notag
\end{align}
\end{theorem}

\begin{proof}
We know, see for example \cite[Eq. (4.4)]{cb.16.0}, that the following
reverse of Jensen's inequality for Fr\'{e}chet differentiable convex
functions%
\begin{align}
0& \leq \sum_{i=1}^{n}q_{i}F\left( z_{i}\right) -F\left(
\sum_{i=1}^{n}q_{i}z_{i}\right)  \label{4.3.16.0} \\
& \leq \sum_{i=1}^{n}q_{i}\left\langle \triangledown F\left( z_{i}\right)
,z_{i}\right\rangle -\left\langle \sum_{i=1}^{n}q_{i}\triangledown F\left(
z_{i}\right) ,\sum_{i=1}^{n}q_{i}z_{i}\right\rangle  \notag
\end{align}%
holds.

Now, if we apply Theorem \ref{t2.1.16.0} for the choices $%
x_{i}=\triangledown F\left( z_{i}\right) ,y_{i}=z_{i}$ and $%
p_{i}=q_{i}\left( i=1,...,n\right) ,$ then we may state%
\begin{align}
& \left\vert \sum_{i=1}^{n}q_{i}\left\langle \triangledown F\left(
z_{i}\right) ,z_{i}\right\rangle -\left\langle
\sum_{i=1}^{n}q_{i}\triangledown F\left( z_{i}\right)
,\sum_{i=1}^{n}q_{i}z_{i}\right\rangle \right\vert  \label{4.4.16.0} \\
& \leq \left\{ 
\begin{array}{l}
\left[ \sum\limits_{i=1}^{n}i^{2}q_{i}-\left(
\sum\limits_{i=1}^{n}iq_{i}\right) ^{2}\right] \max\limits_{k=1,\dots
,n-1}\left\Vert \Delta \left( \triangledown F\left( z_{k}\right) \right)
\right\Vert \max\limits_{k=1,\dots ,n-1}\left\Vert \Delta z_{k}\right\Vert ;
\\ 
\\ 
\left[ \sum\limits_{1\leq j<i\leq n}q_{i}q_{j}\left( i-j\right) \right]
\left( \sum\limits_{k=1}^{n-1}\left\Vert \Delta \left( \triangledown F\left(
z_{k}\right) \right) \right\Vert ^{p}\right) ^{\frac{1}{p}}\left(
\sum\limits_{k=1}^{n-1}\left\Vert \Delta z_{k}\right\Vert ^{q}\right) ^{%
\frac{1}{q}} \\ 
\hfill \text{if }p>1,\ \frac{1}{p}+\frac{1}{q}=1; \\ 
\\ 
\dfrac{1}{2}\left[ \sum\limits_{i=1}^{n}q_{i}\left( 1-p_{i}\right) \right]
\sum\limits_{k=1}^{n-1}\left\Vert \Delta \left( \triangledown F\left(
z_{k}\right) \right) \right\Vert \sum\limits_{k=1}^{n-1}\left\Vert \Delta
z_{k}\right\Vert .%
\end{array}%
\right.  \notag
\end{align}%
Finally, on making use of the inequalities (\ref{4.3.16.0}) and (\ref%
{4.4.16.0}), we deduce the desired result (\ref{4.2.16.0}).
\end{proof}

The unweighted case may useful in application and is incorporated in the
following corollary.

\begin{corollary}
Let $F:H\rightarrow \mathbb{R}$ be as above and $z_{i}\in H,$ $i\in \left\{
1,\dots ,n\right\} .$ Then we have the inequalities%
\begin{align*}
0& \leq \frac{1}{n}\sum_{i=1}^{n}F\left( z_{i}\right) -F\left( \frac{1}{n}%
\sum_{i=1}^{n}z_{i}\right) \\
& \leq \left\{ 
\begin{array}{l}
\dfrac{n^{2}-1}{12}\max\limits_{k=1,\dots ,n-1}\left\Vert \Delta \left(
\triangledown F\left( z_{k}\right) \right) \right\Vert
\max\limits_{k=1,\dots ,n-1}\left\Vert \Delta z_{k}\right\Vert ; \\ 
\\ 
\dfrac{n^{2}-1}{6n}\left( \sum\limits_{k=1}^{n-1}\left\Vert \Delta \left(
\triangledown F\left( z_{k}\right) \right) \right\Vert ^{p}\right) ^{\frac{1%
}{p}}\left( \sum\limits_{k=1}^{n-1}\left\Vert \Delta z_{k}\right\Vert
^{q}\right) ^{\frac{1}{q}} \\ 
\hfill \text{if }p>1,\ \frac{1}{p}+\frac{1}{q}=1; \\ 
\\ 
\dfrac{n-1}{2n}\sum\limits_{k=1}^{n-1}\left\Vert \Delta \left( \triangledown
F\left( z_{k}\right) \right) \right\Vert \sum\limits_{k=1}^{n-1}\left\Vert
\Delta z_{k}\right\Vert .%
\end{array}%
\right. \\
&
\end{align*}
\end{corollary}

\newpage

\section{Bounds for a Pair of n-Tuples of Vectors}

\subsection{Introduction}

Let $\left( H;\left\langle \cdot ,\cdot \right\rangle \right) $ be an inner
product over the real or complex number field $\mathbb{K}$. For $\mathbf{p}%
=\left( p_{1},\dots ,p_{n}\right) \in \mathbb{R}^{n}$ and $\mathbf{x}=\left(
x_{1},\dots ,x_{n}\right) ,$ $\mathbf{y}=\left( y_{1},\dots ,y_{n}\right)
\in H^{n},$ define the \textit{\v{C}eby\v{s}ev functional}%
\begin{equation}
T_{n}\left( \mathbf{\bar{p}};\mathbf{\bar{x}},\mathbf{\bar{y}}\right)
:=P_{n}\sum_{i=1}^{n}p_{i}\left\langle x_{i},y_{i}\right\rangle
-\left\langle \sum_{i=1}^{n}p_{i}x_{i},\sum_{i=1}^{n}p_{i}y_{i}\right\rangle
,  \label{1.1.16.1}
\end{equation}%
where $P_{n}:=\sum_{i=1}^{n}p_{i}.$

The following Gr\"{u}ss type inequality has been obtained in \cite{SSD1.16.1}%
.

\begin{theorem}
\label{t1.1.16.1}Let $H,$ $\mathbf{x},\mathbf{y}$ be as above and $p_{i}\geq
0$ \ $\left( i\in \left\{ 1,\dots ,n\right\} \right) $ with $%
\sum_{i=1}^{n}p_{i}=1,$ i.e., $\mathbf{p}$ is a probability sequence. If $%
x,X,y,Y\in H$ are such that%
\begin{equation}
\func{Re}\left\langle X-x_{i},x_{i}-x\right\rangle \geq 0,\ \ \ \func{Re}%
\left\langle Y-y_{i},y_{i}-y\right\rangle \geq 0  \label{1.2.16.1}
\end{equation}%
for each $i\in \left\{ 1,\dots ,n\right\} ,$ or, equivalently, (see \cite%
{SSD2.16.1})%
\begin{equation}
\left\Vert x_{i}-\frac{x+X}{2}\right\Vert \leq \frac{1}{2}\left\Vert
X-x\right\Vert ,\ \ \ \ \ \ \left\Vert y_{i}-\frac{y+Y}{2}\right\Vert \leq 
\frac{1}{2}\left\Vert Y-y\right\Vert  \label{1.2a.16.1}
\end{equation}%
for each $i\in \left\{ 1,\dots ,n\right\} ,$ then we have the inequality%
\begin{equation}
\left\vert T_{n}\left( \mathbf{\bar{p}};\mathbf{\bar{x}},\mathbf{\bar{y}}%
\right) \right\vert \leq \frac{1}{4}\left\Vert X-x\right\Vert \left\Vert
Y-y\right\Vert .  \label{1.3.16.1}
\end{equation}%
The constant $\frac{1}{4}$ is best possible in the sense that it cannot be
replaced by a smaller constant.
\end{theorem}

In \cite{SSD3.16.1}, the following Gr\"{u}ss type inequality for the forward
difference of vectors was established.

\begin{theorem}
\label{t1.2.16.1}Let $\mathbf{x}=\left( x_{1},\dots ,x_{n}\right) ,$ $%
\mathbf{y}=\left( y_{1},\dots ,y_{n}\right) \in H^{n}$ and $\mathbf{p}\in 
\mathbb{R}_{+}^{n}$ be a probability sequence. Then one has the inequality:%
\begin{align}
& \left\vert T_{n}\left( \mathbf{\bar{p}};\mathbf{\bar{x}},\mathbf{\bar{y}}%
\right) \right\vert  \label{1.4.16.1} \\
& \leq \left\{ 
\begin{array}{l}
\left[ \sum\limits_{i=1}^{n}i^{2}p_{i}-\left(
\sum\limits_{i=1}^{n}ip_{i}\right) ^{2}\right] \max\limits_{1\leq k\leq
n-1}\left\Vert \Delta x_{k}\right\Vert \max\limits_{1\leq k\leq
n-1}\left\Vert \Delta y_{k}\right\Vert ; \\ 
\\ 
\sum\limits_{1\leq j<i\leq n}p_{i}p_{j}\left( i-j\right) \left(
\sum\limits_{k=1}^{n-1}\left\Vert \Delta x_{k}\right\Vert ^{p}\right) ^{%
\frac{1}{p}}\left( \sum\limits_{k=1}^{n-1}\left\Vert \Delta y_{k}\right\Vert
^{q}\right) ^{\frac{1}{q}} \\ 
\hfill \text{ \ if \ }p>1,\ \frac{1}{p}+\frac{1}{q}+1 \\ 
\\ 
\dfrac{1}{2}\left[ \sum\limits_{i=1}^{n}p_{i}\left( 1-p_{i}\right) \right]
\sum\limits_{k=1}^{n-1}\left\Vert \Delta x_{k}\right\Vert
\sum\limits_{k=1}^{n-1}\left\Vert \Delta y_{k}\right\Vert .%
\end{array}%
\right.  \notag
\end{align}%
The constants $1,1$ and $\frac{1}{2}$ in the right hand side of inequality (%
\ref{1.4.16.1}) are best in the sense that they cannot be replaced by
smaller constants.
\end{theorem}

Another result is incorporated in the following theorem (see \cite{SSD2.16.1}%
).

\begin{theorem}
\label{t1.3.16.1}Let $\mathbf{x},\mathbf{y}$ and $\mathbf{p}$ be as in
Theorem \ref{t1.2.16.1}. If there exist $x,X\in H$ such that%
\begin{equation}
\func{Re}\left\langle X-x_{i},x_{i}-x\right\rangle \geq 0\ \ \text{ for each
\ }i\in \left\{ 1,\dots ,n\right\} ,  \label{1.5.16.1}
\end{equation}%
or, equivalently,%
\begin{equation}
\left\Vert x_{i}-\frac{x+X}{2}\right\Vert \leq \frac{1}{2}\left\Vert
X-x\right\Vert \ \ \text{ for each \ }i\in \left\{ 1,\dots ,n\right\} ,
\label{1.6.16.1}
\end{equation}%
then one has the inequality%
\begin{align}
\left\vert T_{n}\left( \mathbf{\bar{p}};\mathbf{\bar{x}},\mathbf{\bar{y}}%
\right) \right\vert & \leq \frac{1}{2}\left\Vert X-x\right\Vert
\sum\limits_{i=1}^{n}p_{i}\left\Vert
y_{i}-\sum_{j=1}^{n}p_{j}y_{j}\right\Vert  \label{1.7.16.1} \\
& \leq \frac{1}{2}\left\Vert X-x\right\Vert \left[ \sum_{i=1}^{n}p_{i}\left%
\Vert y_{i}\right\Vert ^{2}-\left\Vert
\sum\limits_{i=1}^{n}p_{i}y_{i}\right\Vert ^{2}\right] ^{\frac{1}{2}}. 
\notag
\end{align}%
The constant $\frac{1}{2}$ is best possible in the first and second
inequalities in the sense that it cannot be replaced by a smaller constant.
\end{theorem}

\begin{remark}
If $\mathbf{x}$ and $\mathbf{y}$ satisfy the assumptions of Theorem \ref%
{t1.1.16.1}, then we have the following sequence of inequalities improving
the Gr\"{u}ss inequality (\ref{1.3.16.1}):%
\begin{align}
\left\vert T_{n}\left( \mathbf{\bar{p}};\mathbf{\bar{x}},\mathbf{\bar{y}}%
\right) \right\vert & \leq \frac{1}{2}\left\Vert X-x\right\Vert
\sum_{i=1}^{n}p_{i}\left\Vert y_{i}-\sum_{j=1}^{n}p_{j}y_{j}\right\Vert
\label{1.8.16.1} \\
& \leq \frac{1}{2}\left\Vert X-x\right\Vert \left(
\sum_{i=1}^{n}p_{i}\left\Vert y_{i}\right\Vert ^{2}-\left\Vert
\sum\limits_{i=1}^{n}p_{i}y_{i}\right\Vert ^{2}\right) ^{\frac{1}{2}}  \notag
\\
& \leq \frac{1}{4}\left\Vert X-x\right\Vert \left\Vert Y-y\right\Vert . 
\notag
\end{align}
\end{remark}

Now, if we consider the \v{C}eby\v{s}ev functional for the uniform
probability distribution $u=\left( \frac{1}{n},\dots ,\frac{1}{n}\right) ,$%
\begin{equation*}
T_{n}\left( \mathbf{\bar{x}},\mathbf{\bar{y}}\right) :=\frac{1}{n}%
\sum_{i=1}^{n}\left\langle x_{i},y_{i}\right\rangle -\left\langle \frac{1}{n}%
\sum_{i=1}^{n}x_{i},\frac{1}{n}\sum_{i=1}^{n}y_{i}\right\rangle ,
\end{equation*}%
then, with the assumptions of Theorem \ref{t1.1.16.1}, we have%
\begin{equation}
\left\vert T_{n}\left( \mathbf{\bar{x}},\mathbf{\bar{y}}\right) \right\vert
\leq \frac{1}{4}\left\Vert X-x\right\Vert \left\Vert Y-y\right\Vert .
\label{1.9.16.1}
\end{equation}

Theorem \ref{t1.2.16.1} will provide the following inequalities%
\begin{multline}
\left\vert T_{n}\left( \mathbf{\bar{x}},\mathbf{\bar{y}}\right) \right\vert
\label{1.10.16.1} \\
\leq \left\{ 
\begin{array}{l}
\dfrac{1}{12}\left( n^{2}-1\right) \max\limits_{1\leq k\leq n-1}\left\Vert
\Delta x_{k}\right\Vert \max\limits_{1\leq k\leq n-1}\left\Vert \Delta
y_{k}\right\Vert ; \\ 
\\ 
\dfrac{1}{6}\left( n-\dfrac{1}{n}\right) \left(
\sum\limits_{k=1}^{n-1}\left\Vert \Delta x_{k}\right\Vert ^{p}\right) ^{%
\frac{1}{p}}\left( \sum\limits_{k=1}^{n-1}\left\Vert \Delta y_{k}\right\Vert
^{q}\right) ^{\frac{1}{q}}\hfill \text{ \ if \ }p>1,\ \frac{1}{p}+\frac{1}{q}%
+1; \\ 
\\ 
\dfrac{1}{2}\left( 1-\dfrac{1}{n}\right) \sum\limits_{k=1}^{n-1}\left\Vert
\Delta x_{k}\right\Vert \sum\limits_{k=1}^{n-1}\left\Vert \Delta
y_{k}\right\Vert .%
\end{array}%
\right.
\end{multline}%
Here the constants $\frac{1}{12},$ $\frac{1}{6}$ and $\frac{1}{2}$ are best
possible in the above sense.

Finally, from (\ref{1.8.16.1}), we have%
\begin{align}
\left\vert T_{n}\left( \mathbf{\bar{x}},\mathbf{\bar{y}}\right) \right\vert
& \leq \frac{1}{2n}\left\Vert X-x\right\Vert \sum_{i=1}^{n}\left\Vert y_{i}-%
\dfrac{1}{n}\sum_{j=1}^{n}y_{j}\right\Vert  \label{1.11.16.1} \\
& \leq \frac{1}{2}\left\Vert X-x\right\Vert \left( \frac{1}{n}%
\sum_{i=1}^{n}\left\Vert y_{i}\right\Vert ^{2}-\left\Vert \frac{1}{n}%
\sum_{i=1}^{n}y_{i}\right\Vert ^{2}\right) ^{\frac{1}{2}}  \notag \\
& \leq \frac{1}{4}\left\Vert X-x\right\Vert \left\Vert Y-y\right\Vert . 
\notag
\end{align}

It is the main aim of this section to point out other bounds for the \v{C}eby%
\v{s}ev functionals $T_{n}\left( \mathbf{p},\mathbf{x},\mathbf{y}\right) $
and $T_{n}\left( \mathbf{x},\mathbf{y}\right) .$

\subsection{Identities for Inner Products}

For $\mathbf{p}=\left( p_{1},\dots ,p_{n}\right) \in \mathbb{R}^{n}$ and $%
\mathbf{a}=\left( a_{1},\dots ,a_{n}\right) \in H^{n}$ we define%
\begin{equation*}
P_{i}:=\sum_{k=1}^{i}p_{k},\ \ \ \ \ \ \bar{P}_{i}=P_{n}-P_{i},\ \ \ \ \
i\in \left\{ 1,\dots ,n-1\right\}
\end{equation*}%
and the vectors%
\begin{equation*}
A_{i}\left( \mathbf{p}\right) =\sum_{k=1}^{i}p_{k}a_{k},\ \ \ \ \ \ \bar{A}%
_{i}\left( \mathbf{p}\right) =A_{n}\left( \mathbf{p}\right) -A_{i}\left( 
\mathbf{p}\right)
\end{equation*}%
for $i\in \left\{ 1,\dots ,n-1\right\} .$

The following result holds.

\begin{theorem}
\label{t2.1.16.1}Let $\left( H;\left\langle \cdot ,\cdot \right\rangle
\right) $ be an inner product space over $\mathbb{K}$, $\mathbf{p}=\left(
p_{1},\dots ,p_{n}\right) \in \mathbb{R}^{n}$ and $\mathbf{a}=\left(
a_{1},\dots ,a_{n}\right) ,\mathbf{b}=\left( b_{1},\dots ,b_{n}\right) \in
H^{n}.$ Then we have the identities%
\begin{align}
T_{n}\left( \mathbf{p};\mathbf{a},\mathbf{b}\right) &
=\sum_{i=1}^{n-1}\left\langle P_{i}A_{n}\left( \mathbf{p}\right)
-P_{n}A_{i}\left( \mathbf{p}\right) ,\Delta b_{i}\right\rangle
\label{2.1.16.1} \\
& =P_{n}\sum_{i=1}^{n-1}P_{i}\left\langle \frac{1}{P_{n}}A_{n}\left( \mathbf{%
p}\right) -\frac{1}{P_{i}}A_{i}\left( \mathbf{p}\right) ,\Delta
b_{i}\right\rangle  \notag \\
\text{(if }P_{i}& \neq 0,\ i\in \left\{ 1,\dots ,n\right\} \text{)}  \notag
\\
& =\sum_{i=1}^{n-1}P_{i}\bar{P}_{i}\left\langle \frac{1}{\bar{P}_{i}}\bar{A}%
_{i}\left( \mathbf{p}\right) -\frac{1}{P_{i}}A_{i}\left( \mathbf{p}\right)
,\Delta b_{i}\right\rangle  \notag \\
\text{(if }P_{i},\bar{P}_{i}& \neq 0,\ i\in \left\{ 1,\dots ,n-1\right\} 
\text{),}  \notag
\end{align}%
where $\Delta x_{i}=x_{i+1}-x_{i}$ $\left( i\in \left\{ 1,\dots ,n-1\right\}
\right) $ is the forward difference.
\end{theorem}

\begin{proof}
We use the following summation by parts formula for vectors in inner product
spaces%
\begin{equation}
\sum_{l=p}^{q-1}\left\langle d_{l},\Delta v_{l}\right\rangle =\left\langle
d_{l},v_{l}\right\rangle \big|_{p}^{q}-\sum_{l=p}^{q-1}\left\langle
v_{l+1},\Delta d_{l}\right\rangle  \label{2.2.16.1}
\end{equation}%
where $d_{l},$ $v_{l}$ are vectors in $H,$ $l=p,\dots ,q$ ($q>p;$ $p,q$ are
natural numbers).

If we choose in (\ref{2.2.16.1}), $p=1,$ $q=n,$ $d_{i}=P_{i}A_{n}\left( 
\mathbf{p}\right) -P_{n}A_{i}\left( \mathbf{p}\right) $ and $v_{i}=b_{i}$ $%
\left( i\in \left\{ 1,\dots ,n-1\right\} \right) ,$ then we get%
\begin{eqnarray*}
&&\sum_{i=1}^{n-1}\left\langle P_{i}A_{n}\left( \mathbf{p}\right)
-P_{n}A_{i}\left( \mathbf{p}\right) ,\Delta b_{i}\right\rangle \\
&=&\left\langle P_{i}A_{n}\left( \mathbf{p}\right) -P_{n}A_{i}\left( \mathbf{%
p}\right) ,b_{i}\right\rangle \big|_{1}^{n}-\sum_{i=1}^{n-1}\left\langle
\Delta \left( P_{i}A_{n}\left( \mathbf{p}\right) -P_{n}A_{i}\left( \mathbf{p}%
\right) \right) ,b_{i+1}\right\rangle \\
&=&\left\langle P_{n}A_{n}\left( \mathbf{p}\right) -P_{n}A_{n}\left( \mathbf{%
p}\right) ,b_{n}\right\rangle -\left\langle P_{1}A_{n}\left( \mathbf{p}%
\right) -P_{n}A_{1}\left( \mathbf{p}\right) ,b_{1}\right\rangle \\
&&-\sum_{i=1}^{n-1}\left\langle P_{i+1}A_{n}\left( \mathbf{p}\right)
-P_{n}A_{i+1}\left( \mathbf{p}\right) -P_{i}A_{n}\left( \mathbf{p}\right)
+P_{n}A_{i}\left( \mathbf{p}\right) ,b_{i+1}\right\rangle \\
&=&P_{n}p_{1}\left\langle a_{1},x_{1}\right\rangle -p_{1}\left\langle
A_{n}\left( \mathbf{p}\right) ,b_{1}\right\rangle -\left\langle A_{n}\left( 
\mathbf{p}\right) ,\sum_{i=1}^{n-1}p_{i+1}b_{i+1}\right\rangle \\
&&+P_{n}\sum_{i=1}^{n-1}p_{i+1}\left\langle a_{i+1},b_{i+1}\right\rangle \\
&=&P_{n}\sum_{i=1}^{n}p_{i}\left\langle a_{i},b_{i}\right\rangle
-\left\langle \sum_{i=1}^{n}p_{i}a_{i},\sum_{i=1}^{n}p_{i}b_{i}\right\rangle
\\
&=&T_{n}\left( \mathbf{p};\mathbf{a},\mathbf{b}\right) ,
\end{eqnarray*}%
proving the first identity in (\ref{2.1.16.1}).

The second and third identities are obvious and we omit the details.
\end{proof}

The following lemma is of interest in itself.

\begin{lemma}
\label{l2.2.16.1}Let $\mathbf{p}=\left( p_{1},\dots ,p_{n}\right) \in 
\mathbb{R}^{n}$ and $\mathbf{a}=\left( a_{1},\dots ,a_{n}\right) \in H.$
Then we have the equality%
\begin{equation}
P_{i}A_{n}\left( \mathbf{p}\right) -P_{n}A_{i}\left( \mathbf{p}\right)
=\sum_{j=1}^{n-1}P_{\min \left\{ i,j\right\} }\bar{P}_{\max \left\{
i,j\right\} }\Delta a_{j}  \label{2.3.16.1}
\end{equation}%
for each $i\in \left\{ 1,\dots ,n-1\right\} .$
\end{lemma}

\begin{proof}
Define, for $i\in \left\{ 1,\dots ,n-1\right\} ,$ the vector%
\begin{equation*}
K\left( i\right) :=\sum_{j=1}^{n-1}P_{\min \left\{ i,j\right\} }\bar{P}%
_{\max \left\{ i,j\right\} }\cdot \Delta a_{j}.
\end{equation*}%
We have%
\begin{align}
K\left( i\right) & =\sum_{j=1}^{i}P_{\min \left\{ i,j\right\} }\bar{P}_{\max
\left\{ i,j\right\} }\cdot \Delta a_{j}+\sum_{j=i+1}^{n-1}P_{\min \left\{
i,j\right\} }\bar{P}_{\max \left\{ i,j\right\} }\cdot \Delta a_{j}
\label{2.4.16.1} \\
& =\sum_{j=1}^{i}P_{j}\bar{P}_{i}\cdot \Delta a_{j}+\sum_{j=i+1}^{n-1}P_{i}%
\bar{P}_{j}\cdot \Delta a_{j}  \notag \\
& =\bar{P}_{i}\sum_{j=1}^{i}P_{j}\cdot \Delta a_{j}+P_{i}\sum_{j=i+1}^{n-1}%
\bar{P}_{j}\cdot \Delta a_{j}.  \notag
\end{align}%
Using the summation by parts formula, we have%
\begin{align}
\sum_{j=1}^{i}P_{j}\cdot \Delta a_{j}& =P_{j}a_{j}\big|_{1}^{i+1}-%
\sum_{j=1}^{i}\left( P_{j+1}-P_{j}\right) a_{j+1}  \label{2.5.16.1} \\
& =P_{i+1}a_{i+1}-p_{1}a_{1}-\sum_{j=1}^{i}p_{j+1}a_{j+1}  \notag \\
& =P_{i+1}a_{i+1}-\sum_{j=1}^{i+1}p_{j}a_{j}  \notag
\end{align}%
and%
\begin{align}
\sum_{j=i+1}^{n-1}\bar{P}_{j}\cdot \Delta a_{j}& =\bar{P}_{j}a_{j}\big|%
_{i+1}^{n}-\sum_{j=i+1}^{n-1}\left( \bar{P}_{j+1}-\bar{P}_{j}\right) a_{j+1}
\label{2.6.16.1} \\
& =\bar{P}_{n}a_{n}-\bar{P}_{i+1}a_{i+1}-\sum_{j=i+1}^{n-1}\left(
P_{n}-P_{j+1}-P_{n}+P_{j}\right) a_{j+1}  \notag \\
& =-\bar{P}_{i+1}a_{i+1}+\sum_{j=i+1}^{n-1}p_{j+1}a_{j+1}.  \notag
\end{align}%
Using (\ref{2.5.16.1}) and (\ref{2.6.16.1}), we have%
\begin{align*}
K\left( i\right) & =\bar{P}_{i}\left(
P_{i+1}a_{i+1}-\sum_{j=1}^{i+1}p_{j}a_{j}\right) +P_{i}\left(
\sum_{j=i+1}^{n-1}p_{j+1}a_{j+1}-\bar{P}_{i+1}a_{i+1}\right) \\
& =\bar{P}_{i}P_{i+1}a_{i+1}-\bar{P}_{i}\bar{P}_{i+1}a_{i+1}-\bar{P}%
_{i}\sum_{j=1}^{i+1}p_{j}a_{j}+P_{i}\sum_{j=i+1}^{n-1}p_{j+1}a_{j+1} \\
& =\left[ \left( P_{n}-P_{i}\right) P_{i+1}-P_{i}\left( P_{n}-P_{i+1}\right) %
\right] a_{i+1}+P_{i}\sum_{j=i+1}^{n-1}p_{j+1}a_{j+1}-\bar{P}%
_{i}\sum_{j=1}^{i+1}p_{j}a_{j} \\
& =P_{n}p_{i+1}a_{i+1}+P_{i}\sum_{j=i+1}^{n-1}p_{j+1}a_{j+1}-\bar{P}%
_{i}\sum_{j=1}^{i+1}p_{j}a_{j} \\
& =\left( P_{i}+\bar{P}_{i}\right)
p_{i+1}a_{i+1}+P_{i}\sum_{j=i+1}^{n-1}p_{j+1}a_{j+1}-\bar{P}%
_{i}\sum_{j=1}^{i+1}p_{j}a_{j} \\
& =P_{i}\sum_{j=i+1}^{n-1}p_{j}a_{j}-\bar{P}_{i}\sum_{j=1}^{i}p_{j}a_{j} \\
& =P_{i}\bar{A}_{i}\left( \mathbf{p}\right) -\bar{P}_{i}A_{i}\left( \mathbf{p%
}\right) \\
& =P_{i}A_{n}\left( \mathbf{p}\right) -P_{n}A_{i}\left( \mathbf{p}\right) ,
\end{align*}%
and the identity is proved.
\end{proof}

We are able now to state and prove the second identity for the \v{C}eby\v{s}%
ev functional.

\begin{theorem}
\label{t2.3.16.1}With the assumptions of Theorem \ref{t2.1.16.1}, we have
the identity%
\begin{equation}
T_{n}\left( \mathbf{p};\mathbf{a},\mathbf{b}\right)
=\sum_{i=1}^{n-1}\sum_{j=1}^{n-1}P_{\min \left\{ i,j\right\} }\bar{P}_{\max
\left\{ i,j\right\} }\cdot \left\langle \Delta a_{j},\Delta
b_{i}\right\rangle .  \label{2.7.16.1}
\end{equation}
\end{theorem}

\begin{proof}
Follows by Theorem \ref{t2.1.16.1} and Lemma \ref{l2.2.16.1} and we omit the
details.
\end{proof}

\subsection{New Inequalities}

The following result holds.

\begin{theorem}
\label{t3.1.16.1}Let $\left( H;\left\langle \cdot ,\cdot \right\rangle
\right) $ be an inner product space over the real or complex number field $%
\mathbb{K}$; $\mathbf{p}=\left( p_{1},\dots ,p_{n}\right) \in \mathbb{R}^{n}$
and $\mathbf{a}=\left( a_{1},\dots ,a_{n}\right) ,\mathbf{b}=\left(
b_{1},\dots ,b_{n}\right) \in H^{n}.$ Then we have the inequalities%
\begin{equation}
\left\vert T_{n}\left( \mathbf{p};\mathbf{a},\mathbf{b}\right) \right\vert
\leq \left\{ 
\begin{array}{l}
\max\limits_{1\leq i\leq n-1}\left\Vert P_{i}A_{n}\left( \mathbf{p}\right)
-P_{n}A_{i}\left( \mathbf{p}\right) \right\Vert
\sum\limits_{j=1}^{n-1}\left\Vert \Delta b_{j}\right\Vert ; \\ 
\\ 
\left( \sum\limits_{i=1}^{n-1}\left\Vert P_{i}A_{n}\left( \mathbf{p}\right)
-P_{n}A_{i}\left( \mathbf{p}\right) \right\Vert ^{q}\right) ^{\frac{1}{q}%
}\left( \sum\limits_{j=1}^{n-1}\left\Vert \Delta b_{j}\right\Vert
^{p}\right) ^{\frac{1}{p}} \\ 
\hfill \text{for \ }p>1,\ \frac{1}{p}+\frac{1}{q}=1; \\ 
\\ 
\sum\limits_{i=1}^{n-1}\left\Vert P_{i}A_{n}\left( \mathbf{p}\right)
-P_{n}A_{i}\left( \mathbf{p}\right) \right\Vert \cdot \max\limits_{1\leq
j\leq n-1}\left\Vert \Delta b_{j}\right\Vert .%
\end{array}%
\right.  \label{3.1.16.1}
\end{equation}%
All the inequalities in (\ref{3.1.16.1}) are sharp in the sense that the
constants $1$ cannot be replaced by smaller constants.
\end{theorem}

\begin{proof}
Using the first identity in (\ref{2.1.16.1}) and Schwarz's inequality in $H,$
i.e., $\left\vert \left\langle u,v\right\rangle \right\vert \leq \left\Vert
u\right\Vert \left\Vert v\right\Vert ,$ $u,v\in H,$ we have successively:%
\begin{align*}
\left\vert T_{n}\left( \mathbf{p};\mathbf{a},\mathbf{b}\right) \right\vert &
\leq \sum\limits_{i=1}^{n-1}\left\vert \left\langle P_{i}A_{n}\left( \mathbf{%
p}\right) -P_{n}A_{i}\left( \mathbf{p}\right) ,\Delta b_{i}\right\rangle
\right\vert \\
& \leq \sum\limits_{i=1}^{n-1}\left\Vert P_{i}A_{n}\left( \mathbf{p}\right)
-P_{n}A_{i}\left( \mathbf{p}\right) \right\Vert \left\Vert \Delta
b_{i}\right\Vert .
\end{align*}%
Using H\"{o}lder's inequality, we deduce the desired result (\ref{3.1.16.1}).

Let us prove, for instance, that the constant 1 in the second inequality is
best possible.

Assume, for $c>0,$ we have that%
\begin{equation}
\left\vert T_{n}\left( \mathbf{p};\mathbf{a},\mathbf{b}\right) \right\vert
\leq c\left( \sum\limits_{i=1}^{n-1}\left\Vert P_{i}A_{n}\left( \mathbf{p}%
\right) -P_{n}A_{i}\left( \mathbf{p}\right) \right\Vert ^{q}\right) ^{\frac{1%
}{q}}\left( \sum\limits_{j=1}^{n-1}\left\Vert \Delta b_{j}\right\Vert
^{p}\right) ^{\frac{1}{p}}  \label{3.2.16.1}
\end{equation}%
for $p>1,$ $\frac{1}{p}+\frac{1}{q}=1,$ $n\geq 2.$

If we choose $n=2,$ then we get%
\begin{equation*}
\left\vert T_{2}\left( \mathbf{p};\mathbf{a},\mathbf{b}\right) \right\vert
\leq p_{1}p_{2}\left\langle a_{2}-a_{1},b_{2}-b_{1}\right\rangle .
\end{equation*}%
Also, for $n=2,$%
\begin{equation*}
\left( \sum\limits_{i=1}^{n-1}\left\Vert P_{i}A_{n}\left( \mathbf{p}\right)
-P_{n}A_{i}\left( \mathbf{p}\right) \right\Vert ^{q}\right) ^{\frac{1}{q}%
}=\left\vert p_{1}p_{2}\right\vert \left\Vert a_{2}-a_{1}\right\Vert
\end{equation*}%
and 
\begin{equation*}
\left( \sum\limits_{j=1}^{n-1}\left\Vert \Delta b_{j}\right\Vert ^{p}\right)
^{\frac{1}{p}}=\left\Vert b_{2}-b_{1}\right\Vert ,
\end{equation*}%
and then, from (\ref{3.2.16.1}), for $n=2,$ we deduce%
\begin{equation}
\left\vert p_{1}p_{2}\right\vert \left\vert \left\langle
a_{2}-a_{1},b_{2}-b_{1}\right\rangle \right\vert \leq c\left\vert
p_{1}p_{2}\right\vert \left\Vert a_{2}-a_{1}\right\Vert \left\Vert
b_{2}-b_{1}\right\Vert .  \label{3.3.16.1}
\end{equation}%
If in (\ref{3.3.16.1}) we choose $a_{2}=b_{2},$ $a_{2}=b_{1}$ and $b_{2}\neq
b_{1},$ $p_{1},p_{2}\neq 0,$ we deduce $c\geq 1,$ proving that 1 is the best
possible constant in that inequality.
\end{proof}

The following corollary for the uniform distribution of the probability $%
\mathbf{p}$ holds.

\begin{corollary}
\label{c3.2.16.1}With the assumptions of Theorem \ref{t3.1.16.1} for $%
\mathbf{a}$ and $\mathbf{b},$ we have the inequalities%
\begin{equation}
0\leq \left\vert T_{n}\left( \mathbf{a},\mathbf{b}\right) \right\vert \leq 
\frac{1}{n^{2}}\left\{ 
\begin{array}{l}
\max\limits_{1\leq i\leq n-1}\left\Vert
i\sum\limits_{k=1}^{n}a_{k}-n\sum\limits_{k=1}^{i}a_{k}\right\Vert
\sum\limits_{j=1}^{n-1}\left\Vert \Delta b_{j}\right\Vert ; \\ 
\\ 
\left( \sum\limits_{i=1}^{n-1}\left\Vert
i\sum\limits_{k=1}^{n}a_{k}-n\sum\limits_{k=1}^{i}a_{k}\right\Vert
^{q}\right) ^{\frac{1}{q}}\left( \sum\limits_{j=1}^{n-1}\left\Vert \Delta
b_{j}\right\Vert ^{p}\right) ^{\frac{1}{p}} \\ 
\hfill \text{for \ }p>1,\ \frac{1}{p}+\frac{1}{q}=1; \\ 
\\ 
\sum\limits_{i=1}^{n-1}\left\Vert
i\sum\limits_{k=1}^{n}a_{k}-n\sum\limits_{k=1}^{i}a_{k}\right\Vert \cdot
\max\limits_{1\leq j\leq n-1}\left\Vert \Delta b_{j}\right\Vert .%
\end{array}%
\right.  \label{3.4.16.1}
\end{equation}
\end{corollary}

The following result may be stated as well.

\begin{theorem}
\label{t3.3.16.1}With the assumptions of Theorem \ref{t3.1.16.1} and if $%
P_{i}\neq 0$ $\left( i=1,\dots ,n\right) ,$ then we have the inequalities%
\begin{multline}
\left\vert T_{n}\left( \mathbf{p};\mathbf{a},\mathbf{b}\right) \right\vert
\label{3.5.16.1} \\
\leq \left\vert P_{n}\right\vert \times \left\{ 
\begin{array}{l}
\max\limits_{1\leq i\leq n-1}\left\Vert \dfrac{1}{P_{n}}A_{n}\left( \mathbf{p%
}\right) -\dfrac{1}{P_{i}}A_{i}\left( \mathbf{p}\right) \right\Vert
\sum\limits_{i=1}^{n-1}\left\vert P_{i}\right\vert \left\Vert \Delta
b_{i}\right\Vert ; \\ 
\\ 
\left( \sum\limits_{i=1}^{n-1}\left\vert P_{i}\right\vert \left\Vert \dfrac{1%
}{P_{n}}A_{n}\left( \mathbf{p}\right) -\dfrac{1}{P_{i}}A_{i}\left( \mathbf{p}%
\right) \right\Vert ^{q}\right) ^{\frac{1}{q}}\left(
\sum\limits_{i=1}^{n}\left\vert P_{i}\right\vert \left\Vert \Delta
b_{i}\right\Vert ^{p}\right) ^{\frac{1}{p}} \\ 
\hfill \text{for \ }p>1,\ \frac{1}{p}+\frac{1}{q}=1; \\ 
\\ 
\sum\limits_{i=1}^{n-1}\left\vert P_{i}\right\vert \left\Vert \dfrac{1}{P_{n}%
}A_{n}\left( \mathbf{p}\right) -\dfrac{1}{P_{i}}A_{i}\left( \mathbf{p}%
\right) \right\Vert \cdot \max\limits_{1\leq i\leq n-1}\left\Vert \Delta
b_{i}\right\Vert .%
\end{array}%
\right.
\end{multline}%
All the inequalities in (\ref{3.5.16.1}) are sharp in the sense that the
constant 1 cannot be replaced by a smaller constant.
\end{theorem}

\begin{proof}
Using the second equality in (\ref{2.1.16.1}) and Schwarz's inequality, we
have%
\begin{align*}
\left\vert T_{n}\left( \mathbf{p};\mathbf{a},\mathbf{b}\right) \right\vert &
\leq \left\vert P_{n}\right\vert \sum\limits_{i=1}^{n-1}\left\vert
P_{i}\right\vert \left\vert \left\langle \dfrac{1}{P_{n}}A_{n}\left( \mathbf{%
p}\right) -\dfrac{1}{P_{i}}A_{i}\left( \mathbf{p}\right) ,\Delta
b_{i}\right\rangle \right\vert \\
& \leq \left\vert P_{n}\right\vert \sum\limits_{i=1}^{n-1}\left\vert
P_{i}\right\vert \left\Vert \dfrac{1}{P_{n}}A_{n}\left( \mathbf{p}\right) -%
\dfrac{1}{P_{i}}A_{i}\left( \mathbf{p}\right) \right\Vert \left\Vert \Delta
b_{i}\right\Vert .
\end{align*}%
Using H\"{o}lder's weighted inequality, we deduce (\ref{3.5.16.1}).

The sharpness of the constant may be proven in a similar manner to the one
in Theorem \ref{t3.1.16.1}. We omit the details.
\end{proof}

The following corollary containing the unweighted inequalities holds.

\begin{corollary}
\label{c3.4.16.1}With the above assumptions for $\mathbf{a}$ and $\mathbf{b}%
, $ one has%
\begin{equation}
\left\vert T_{n}\left( \mathbf{a},\mathbf{b}\right) \right\vert \leq \frac{1%
}{n}\left\{ 
\begin{array}{l}
\max\limits_{1\leq i\leq n-1}\left\Vert \dfrac{1}{n}\sum%
\limits_{k=1}^{n}a_{k}-\dfrac{1}{i}\sum\limits_{k=1}^{i}a_{k}\right\Vert
\sum\limits_{i=1}^{n-1}i\left\Vert \Delta b_{i}\right\Vert ; \\ 
\\ 
\left( \sum\limits_{i=1}^{n-1}i\left\Vert \dfrac{1}{n}\sum%
\limits_{k=1}^{n}a_{k}-\dfrac{1}{i}\sum\limits_{k=1}^{i}a_{k}\right\Vert
^{q}\right) ^{\frac{1}{q}}\left( \sum\limits_{i=1}^{n-1}i\left\Vert \Delta
b_{i}\right\Vert ^{p}\right) ^{\frac{1}{p}} \\ 
\hfill \text{for \ }p>1,\ \frac{1}{p}+\frac{1}{q}=1; \\ 
\\ 
\sum\limits_{i=1}^{n-1}i\left\Vert \dfrac{1}{n}\sum\limits_{k=1}^{n}a_{k}-%
\dfrac{1}{i}\sum\limits_{k=1}^{i}a_{k}\right\Vert \cdot \max\limits_{1\leq
i\leq n-1}\left\Vert \Delta b_{i}\right\Vert .%
\end{array}%
\right.  \label{3.6.16.1}
\end{equation}%
The inequalities (\ref{3.6.16.1}) are sharp in the sense mentioned above.
\end{corollary}

Another type of inequality may be stated if ones used the third identity in (%
\ref{2.1.16.1}) and H\"{o}lder's weighted inequality with the weights: $%
\left\vert P_{i}\right\vert \left\vert \bar{P}_{i}\right\vert ,$ $i\in
\left\{ 1,\dots ,n-1\right\} .$

\begin{theorem}
\label{t3.5.16.1}With the assumptions in Theorem \ref{t3.1.16.1} and if $%
P_{i},$ $\bar{P}_{i}\neq 0,$ $i\in \left\{ 1,\dots ,n-1\right\} ,$ then we
have the inequalities%
\begin{multline}
\left\vert T_{n}\left( \mathbf{p};\mathbf{a},\mathbf{b}\right) \right\vert
\label{3.7.16.1} \\
\leq \left\vert P_{n}\right\vert \times \left\{ 
\begin{array}{l}
\max\limits_{1\leq i\leq n-1}\left\Vert \dfrac{1}{\bar{P}_{i}}\bar{A}%
_{i}\left( \mathbf{p}\right) -\dfrac{1}{P_{i}}A_{i}\left( \mathbf{p}\right)
\right\Vert \sum\limits_{i=1}^{n-1}\left\vert P_{i}\right\vert \left\vert 
\bar{P}_{i}\right\vert \left\Vert \Delta b_{i}\right\Vert ; \\ 
\\ 
\left( \sum\limits_{i=1}^{n-1}\left\vert P_{i}\right\vert \left\vert \bar{P}%
_{i}\right\vert \left\Vert \dfrac{1}{\bar{P}_{i}}\bar{A}_{i}\left( \mathbf{p}%
\right) -\dfrac{1}{P_{i}}A_{i}\left( \mathbf{p}\right) \right\Vert
^{q}\right) ^{\frac{1}{q}}\left( \sum\limits_{i=1}^{n-1}\left\vert
P_{i}\right\vert \left\vert \bar{P}_{i}\right\vert \left\Vert \Delta
b_{i}\right\Vert ^{p}\right) ^{\frac{1}{p}} \\ 
\hfill \text{for \ }p>1,\ \frac{1}{p}+\frac{1}{q}=1; \\ 
\\ 
\sum\limits_{i=1}^{n-1}\left\vert P_{i}\right\vert \left\vert \bar{P}%
_{i}\right\vert \left\Vert \dfrac{1}{\bar{P}_{i}}\bar{A}_{i}\left( \mathbf{p}%
\right) -\dfrac{1}{P_{i}}A_{i}\left( \mathbf{p}\right) \right\Vert \cdot
\max\limits_{1\leq i\leq n-1}\left\Vert \Delta b_{i}\right\Vert .%
\end{array}%
\right.
\end{multline}%
In particular, if $p_{i}=\frac{1}{n},$ $i\in \left\{ 1,\dots ,n\right\} ,$
then we have%
\begin{multline}
\left\vert T_{n}\left( \mathbf{a},\mathbf{b}\right) \right\vert
\label{3.8.16.1} \\
\leq \frac{1}{n^{2}}\left\{ 
\begin{array}{l}
\max\limits_{1\leq i\leq n-1}\left\Vert \dfrac{1}{n-i}\sum%
\limits_{k=i+1}^{n}a_{k}-\dfrac{1}{i}\sum\limits_{k=1}^{i}a_{k}\right\Vert
\sum\limits_{i=1}^{n-1}i\left( n-i\right) \left\Vert \Delta b_{i}\right\Vert
; \\ 
\\ 
\left( \sum\limits_{i=1}^{n-1}i\left( n-i\right) \left\Vert \dfrac{1}{n-i}%
\sum\limits_{k=i+1}^{n}a_{k}-\dfrac{1}{i}\sum\limits_{k=1}^{i}a_{k}\right%
\Vert ^{q}\right) ^{\frac{1}{q}}\left( \sum\limits_{i=1}^{n-1}i\left(
n-i\right) \left\Vert \Delta b_{i}\right\Vert ^{p}\right) ^{\frac{1}{p}} \\ 
\hfill \text{for \ }p>1,\ \frac{1}{p}+\frac{1}{q}=1; \\ 
\\ 
\sum\limits_{i=1}^{n-1}i\left( n-i\right) \left\Vert \dfrac{1}{n-i}%
\sum\limits_{k=i+1}^{n}a_{k}-\dfrac{1}{i}\sum\limits_{k=1}^{i}a_{k}\right%
\Vert \cdot \max\limits_{1\leq i\leq n-1}\left\Vert \Delta b_{i}\right\Vert .%
\end{array}%
\right.
\end{multline}%
The inequalities in (\ref{3.7.16.1}) and (\ref{3.8.16.1}) are sharp in the
above mentioned sense.
\end{theorem}

A different approach may be considered if one uses the representation in
terms of double sums for the \v{C}eby\v{s}ev functional provided by Theorem %
\ref{t2.3.16.1}.

The following result holds.

\begin{theorem}
\label{t3.6.16.1}With the above assumptions of Theorem \ref{t3.1.16.1}, we
have the inequalities%
\begin{multline}
\left\vert T_{n}\left( \mathbf{p};\mathbf{a},\mathbf{b}\right) \right\vert
\label{3.9.16.1} \\
\leq \left\vert P_{n}\right\vert \times \left\{ 
\begin{array}{l}
\max\limits_{1\leq i,j\leq n-1}\left\{ \left\vert P_{\min \left\{
i,j\right\} }\right\vert ,\left\vert \bar{P}_{\max \left\{ i,j\right\}
}\right\vert \right\} \sum\limits_{i=1}^{n-1}\left\Vert \Delta
a_{i}\right\Vert \sum\limits_{i=1}^{n-1}\left\Vert \Delta b_{i}\right\Vert ;
\\ 
\\ 
\left( \sum\limits_{i=1}^{n-1}\sum\limits_{j=1}^{n-1}\left\vert P_{\min
\left\{ i,j\right\} }\right\vert ^{q}\left\vert \bar{P}_{\max \left\{
i,j\right\} }\right\vert ^{q}\right) ^{\frac{1}{q}}\left(
\sum\limits_{i=1}^{n-1}\left\Vert \Delta a_{i}\right\Vert ^{p}\right) ^{%
\frac{1}{p}}\left( \sum\limits_{i=1}^{n-1}\left\Vert \Delta b_{i}\right\Vert
^{p}\right) ^{\frac{1}{p}} \\ 
\hfill \text{for \ }p>1,\ \frac{1}{p}+\frac{1}{q}=1; \\ 
\\ 
\sum\limits_{i=1}^{n-1}\sum\limits_{j=1}^{n-1}\left\vert P_{\min \left\{
i,j\right\} }\right\vert \left\vert \bar{P}_{\max \left\{ i,j\right\}
}\right\vert \max\limits_{1\leq i\leq n-1}\left\Vert \Delta a_{i}\right\Vert
\max\limits_{1\leq i\leq n-1}\left\Vert \Delta b_{i}\right\Vert .%
\end{array}%
\right.
\end{multline}%
The inequalities are sharp in the sense mentioned above.
\end{theorem}

The proof follows by the identity (\ref{2.7.16.1}) on using H\"{o}lder's
inequality for double sums and we omit the details.

Now, define%
\begin{equation*}
k_{\infty }:=\max_{1\leq i,j\leq n-1}\left\{ \frac{\min \left\{ i,j\right\} 
}{n}\left( 1-\frac{\max \left\{ i,j\right\} }{n}\right) \right\} ,\ \ n\geq
2.
\end{equation*}%
Using the elementary inequality%
\begin{equation*}
ab\leq \frac{1}{4}\left( a+b\right) ^{2},\ \ \ \ a,b\in \mathbb{R};
\end{equation*}%
we deduce%
\begin{equation*}
\min \left\{ i,j\right\} \left( n-\max \left\{ i,j\right\} \right) \leq 
\frac{1}{4}\left( n-\left\vert i-j\right\vert \right) ^{2}
\end{equation*}%
for $1\leq i,j\leq n-2.$ Consequently, we deduce%
\begin{equation*}
k_{\infty }\leq \frac{1}{4n^{2}}\max_{1\leq i,j\leq n-1}\left\{ \left(
n-\left\vert i-j\right\vert \right) ^{2}\right\} =\frac{1}{4}.
\end{equation*}

We may now state the following corollary of Theorem \ref{t3.6.16.1}.

\begin{corollary}
\label{c3.7.16.1}With the assumptions of Theorem \ref{t3.1.16.1} for $%
\mathbf{a}$ and $\mathbf{b},$ we have the inequality%
\begin{align}
\left\vert T_{n}\left( \mathbf{a},\mathbf{b}\right) \right\vert & \leq
k_{\infty }\sum\limits_{i=1}^{n-1}\left\Vert \Delta a_{i}\right\Vert
\sum\limits_{i=1}^{n-1}\left\Vert \Delta b_{i}\right\Vert  \label{3.10.16.1}
\\
& \leq \frac{1}{4}\sum\limits_{i=1}^{n-1}\left\Vert \Delta a_{i}\right\Vert
\sum\limits_{i=1}^{n-1}\left\Vert \Delta b_{i}\right\Vert .  \notag
\end{align}%
The constant $\frac{1}{4}$ cannot be replaced in general by a smaller
constant.
\end{corollary}

\begin{remark}
\label{r3.8.16.1}The inequality (\ref{3.10.16.1}) is better than the third
inequality in (\ref{1.10.16.1}).
\end{remark}

Consider now, for $q>1,$ the number%
\begin{equation*}
k_{q}:=\frac{1}{n^{2}}\left( \sum\limits_{i=1}^{n-1}\sum\limits_{j=1}^{n-1}%
\left[ \min \left\{ i,j\right\} \left( n-\max \left\{ i,j\right\} \right) %
\right] ^{q}\right) ^{\frac{1}{q}}.
\end{equation*}%
We observe, by symmetry of the terms under the summation symbol, we have that%
\begin{equation*}
k_{q}=\frac{1}{n^{2}}\left( 2\sum_{1\leq i<j\leq n-1}i^{q}\left( n-j\right)
^{q}+\sum_{i=1}^{n-1}i^{q}\left( n-i\right) ^{q}\right) ^{\frac{1}{q}},
\end{equation*}%
that may be computed exactly if $q=2$ or another natural number.

Since, as above,%
\begin{equation*}
\left[ \min \left\{ i,j\right\} \left( n-\max \left\{ i,j\right\} \right) %
\right] ^{q}\leq \frac{1}{4^{q}}\left( n-\left\vert i-j\right\vert \right)
^{2q},
\end{equation*}%
we deduce%
\begin{align*}
k_{q}& \leq \frac{1}{4n^{2}}\left(
\sum\limits_{i=1}^{n-1}\sum\limits_{j=1}^{n-1}\left( n-\left\vert
i-j\right\vert \right) ^{2q}\right) ^{\frac{1}{q}} \\
& \leq \frac{1}{4n^{2}}\left[ \left( n-1\right) ^{2}n^{2q}\right] ^{\frac{1}{%
q}} \\
& =\frac{1}{4}\left( n-1\right) ^{\frac{2}{q}}.
\end{align*}%
Consequently, we may state the following corollary as well.

\begin{corollary}
\label{c3.9.16.1}With the assumptions of Theorem \ref{t3.1.16.1} for $%
\mathbf{a}$ and $\mathbf{b},$ we have the inequalities%
\begin{align}
\left\vert T_{n}\left( \mathbf{a},\mathbf{b}\right) \right\vert & \leq
k_{q}\left( \sum\limits_{i=1}^{n-1}\left\Vert \Delta a_{i}\right\Vert
^{p}\right) ^{\frac{1}{p}}\left( \sum\limits_{i=1}^{n-1}\left\Vert \Delta
b_{i}\right\Vert ^{p}\right) ^{\frac{1}{p}}  \label{3.11.16.1} \\
& \leq \frac{1}{4}\left( n-1\right) ^{\frac{2}{q}}\left(
\sum\limits_{i=1}^{n-1}\left\Vert \Delta a_{i}\right\Vert ^{p}\right) ^{%
\frac{1}{p}}\left( \sum\limits_{i=1}^{n-1}\left\Vert \Delta b_{i}\right\Vert
^{p}\right) ^{\frac{1}{p}},  \notag
\end{align}%
provided $p>1$ , $\frac{1}{p}+\frac{1}{q}=1.$ The constant $\frac{1}{4}$
cannot be replaced in general by a smaller constant.
\end{corollary}

Finally, if we denote%
\begin{equation*}
k_{1}:=\frac{1}{n^{2}}\sum\limits_{i=1}^{n-1}\sum\limits_{j=1}^{n-1}\left[
\min \left\{ i,j\right\} \left( n-\max \left\{ i,j\right\} \right) \right] ,
\end{equation*}%
then we observe, for $\mathbf{u}=\left( \frac{1}{n},\dots ,\frac{1}{n}%
\right) ,$ $\mathbf{e}=\left( 1,2,\dots ,n\right) ,$ that%
\begin{equation*}
k_{1}=\left\vert T_{n}\left( \mathbf{u};\mathbf{e},\mathbf{e}\right)
\right\vert =\frac{1}{n}\sum\limits_{i=1}^{n}i^{2}-\left( \frac{1}{n}%
\sum\limits_{i=1}^{n}i\right) ^{2}=\frac{1}{12}\left( n^{2}-1\right) ,
\end{equation*}%
and by Theorem \ref{t3.6.16.1}, we deduce the inequality%
\begin{equation*}
\left\vert T_{n}\left( \mathbf{a},\mathbf{b}\right) \right\vert \leq \frac{1%
}{12}\left( n^{2}-1\right) \max_{1\leq j\leq n-1}\left\Vert \Delta
a_{j}\right\Vert \max_{1\leq j\leq n-1}\left\Vert \Delta b_{j}\right\Vert .
\end{equation*}%
Note that, the above inequality has been discovered using a different method
in \cite{SSD3.16.1}. The constant $\frac{1}{12}$ is best possible.

\newpage

\chapter{Other Inequalities in Inner Product Spaces}

\section{The Ostrowski Inequality}

\subsection{Introduction}

In 1951, A.M. Ostrowski \cite[p. 289]{OS.17} proved the following result
(see also \cite[p. 92]{MPF.17}):

\begin{theorem}
\label{t.1.17} Suppose that $\mathbf{a},\mathbf{b}$ and $\mathbf{x}$ are
real $n$-tuples such that $\mathbf{a}\neq 0$ and%
\begin{equation*}
\sum_{i=1}^{n}a_{i}x_{i}=0\text{ and }\sum_{i=1}^{n}b_{i}x_{i}=1.
\end{equation*}%
Then%
\begin{equation*}
\sum_{i=1}^{n}x_{i}^{2}\geq \frac{\sum_{i=1}^{n}a_{i}^{2}}{%
\sum_{i=1}^{n}a_{i}^{2}\sum_{i=1}^{n}b_{i}^{2}-\left(
\sum_{i=1}^{n}a_{i}b_{i}\right) ^{2}}
\end{equation*}%
with equality if and only if 
\begin{equation*}
x_{k}=\frac{b_{k}\sum_{i=1}^{n}a_{i}^{2}-a_{k}\sum_{i=1}^{n}a_{i}b_{i}}{%
\sum_{i=1}^{n}a_{i}^{2}\sum_{i=1}^{n}b_{i}^{2}-\left(
\sum_{i=1}^{n}a_{i}b_{i}\right) ^{2}},
\end{equation*}%
for $k\in \left\{ 1,\dots ,n\right\} .$
\end{theorem}

An integral version of this inequality was obtained by Pearce, Pe\v{c}ari%
\'{c} and Varo\v{s}anec in 1998, \cite{PPV.17}.

H. \v{S}iki\'{c} and T. \v{S}iki\'{c} in 2001, \cite{SS.17}, by the use of
an argument based on orthogonal projections in inner product spaces have
observed that Ostrowski's inequality may be naturally stated in an abstract
setting as follows:

\begin{theorem}
\label{t.2.17} Let $\left( H;\left\langle \cdot ,\cdot \right\rangle \right) 
$ be a real or complex inner product space and $a,b\in H$ two linearly
independent vectors. If $x\in H$ is such that%
\begin{equation*}
\left\langle x,a\right\rangle =0\text{ and }\left\langle x,b\right\rangle =1,
\end{equation*}%
then one has the inequality%
\begin{equation}
\left\Vert x\right\Vert ^{2}\geq \frac{\left\Vert a\right\Vert ^{2}}{%
\left\Vert a\right\Vert ^{2}\left\Vert b\right\Vert ^{2}-\left\vert
\left\langle a,b\right\rangle \right\vert ^{2}},  \label{1.5.17}
\end{equation}%
with equality if and only if%
\begin{equation*}
x=\frac{\left\Vert a\right\Vert ^{2}b-\overline{\left\langle
a,b\right\rangle }\cdot a}{\left\Vert a\right\Vert ^{2}\left\Vert
b\right\Vert ^{2}-\left\vert \left\langle a,b\right\rangle \right\vert ^{2}}.
\end{equation*}
\end{theorem}

In the present section, by the use of elementary arguments and Schwarz's
inequality in inner product spaces, we show that Ostrowski's inequality $%
\left( \ref{1.5.17}\right) $ holds true for a larger class of elements $x\in
H.$ The case of equality is analyzed. Applications for complex sequences and
integrals are also provided.

\subsection{The General Inequality}

The following theorem holds \cite{19NSSD}.

\begin{theorem}
\label{t.3.17} Let $\left( H;\left\langle \cdot ,\cdot \right\rangle \right) 
$ be a real or complex inner product space and $a,b\in H$ two linearly
independent vectors. If $x\in H$ is such that 
\begin{equation}
\left\langle x,a\right\rangle =0,\text{ and }\left\vert \left\langle
x,b\right\rangle \right\vert =1;  \label{2.0.17}
\end{equation}%
then one has the inequality%
\begin{equation}
\left\Vert x\right\Vert ^{2}\geq \frac{\left\Vert a\right\Vert ^{2}}{%
\left\Vert a\right\Vert ^{2}\left\Vert b\right\Vert ^{2}-\left\vert
\left\langle a,b\right\rangle \right\vert ^{2}}.  \label{2.1.17}
\end{equation}%
The equality holds in $\left( \ref{2.1.17}\right) $ if and only if%
\begin{equation*}
x=\mu \left( b-\frac{\overline{\left\langle a,b\right\rangle }}{\left\Vert
a\right\Vert ^{2}}\cdot a\right)
\end{equation*}%
where $\mu \in \mathbb{K}\left( \mathbb{K=R}\text{,}\mathbb{C}\right) $ is
such that 
\begin{equation}
\left\vert \mu \right\vert =\frac{\left\Vert a\right\Vert ^{2}}{\left\Vert
a\right\Vert ^{2}\left\Vert b\right\Vert ^{2}-\left\vert \left\langle
a,b\right\rangle \right\vert ^{2}}.  \label{2.3.17}
\end{equation}
\end{theorem}

\begin{proof}
We use Schwarz's inequality in the inner product space $\left(
H;\left\langle \cdot ,\cdot \right\rangle \right) ,$ i.e.,%
\begin{equation}
\left\Vert u\right\Vert ^{2}\left\Vert v\right\Vert ^{2}\geq \left\vert
\left\langle u,v\right\rangle \right\vert ^{2};\ \ \ \ \ u,v\in H
\label{2.4.17}
\end{equation}%
with equality iff there exists a scalar $\alpha \in \mathbb{K}$ such that $%
u=\alpha v.$

If we apply $\left( \ref{2.4.17}\right) $ for 
\begin{equation*}
u=z-\frac{\left\langle z,c\right\rangle }{\left\Vert c\right\Vert ^{2}}\cdot
c,\ \ \ \ \ v=d-\frac{\left\langle d,c\right\rangle }{\left\Vert
c\right\Vert ^{2}}\cdot c,
\end{equation*}%
where $c\neq 0$ and $c,d,z\in H,$ we have%
\begin{equation}
\left\Vert z-\frac{\left\langle z,c\right\rangle }{\left\Vert c\right\Vert
^{2}}\cdot c\right\Vert ^{2}\left\Vert d-\frac{\left\langle d,c\right\rangle 
}{\left\Vert c\right\Vert ^{2}}\cdot c\right\Vert ^{2}\geq \left\vert
\left\langle z-\frac{\left\langle z,c\right\rangle }{\left\Vert c\right\Vert
^{2}}\cdot c,d-\frac{\left\langle d,c\right\rangle }{\left\Vert c\right\Vert
^{2}}\cdot c\right\rangle \right\vert ^{2}  \label{2.5.17}
\end{equation}%
with equality iff there is a scalar $\beta \in \mathbb{K}$ such that%
\begin{equation}
z=\frac{\left\langle z,c\right\rangle }{\left\Vert c\right\Vert ^{2}}\cdot
c+\beta \left( d-\frac{\left\langle d,c\right\rangle }{\left\Vert
c\right\Vert ^{2}}\cdot c\right) .  \label{2.6.17}
\end{equation}%
Since simple calculations show that%
\begin{align*}
\left\Vert z-\frac{\left\langle z,c\right\rangle }{\left\Vert c\right\Vert
^{2}}\cdot c\right\Vert ^{2}& =\frac{\left\Vert z\right\Vert ^{2}\left\Vert
c\right\Vert ^{2}-\left\vert \left\langle z,c\right\rangle \right\vert ^{2}}{%
\left\Vert c\right\Vert ^{2}}, \\
\left\Vert d-\frac{\left\langle d,c\right\rangle }{\left\Vert c\right\Vert
^{2}}\cdot c\right\Vert ^{2}& =\frac{\left\Vert d\right\Vert ^{2}\left\Vert
c\right\Vert ^{2}-\left\vert \left\langle d,c\right\rangle \right\vert ^{2}}{%
\left\Vert c\right\Vert ^{2}},
\end{align*}%
and%
\begin{equation*}
\left\langle z-\frac{\left\langle z,c\right\rangle }{\left\Vert c\right\Vert
^{2}}\cdot c,d-\frac{\left\langle d,c\right\rangle }{\left\Vert c\right\Vert
^{2}}\cdot c\right\rangle =\frac{\left\langle z,d\right\rangle \left\Vert
c\right\Vert ^{2}-\left\langle z,c\right\rangle \left\langle
c,d\right\rangle }{\left\Vert c\right\Vert ^{2}},
\end{equation*}%
then, by $\left( \ref{2.5.17}\right) $, we deduce%
\begin{multline}
\left[ \left\Vert z\right\Vert ^{2}\left\Vert c\right\Vert ^{2}-\left\vert
\left\langle z,c\right\rangle \right\vert ^{2}\right] \left[ \left\Vert
d\right\Vert ^{2}\left\Vert c\right\Vert ^{2}-\left\vert \left\langle
d,c\right\rangle \right\vert ^{2}\right]  \label{2.7.17} \\
\geq \left\vert \left\langle z,d\right\rangle \left\Vert c\right\Vert
^{2}-\left\langle z,c\right\rangle \left\langle c,d\right\rangle \right\vert
^{2},  \notag
\end{multline}%
with equality if and only if there is a $\beta \in \mathbb{K}$ such that $%
\left( \ref{2.6.17}\right) $ holds.

If $a,x,b$ satisfy $\left( \ref{2.0.17}\right) $ then by $\left( \ref{2.7.17}%
\right) $ and $\left( \ref{2.6.17}\right) $ for the choices $z=x,c=a$ and $%
d=b$ we deduce the inequality $\left( \ref{2.1.17}\right) $ with equality
iff there exists a $\mu \in \mathbb{K}$ such that 
\begin{equation*}
x=\mu \left( b-\frac{\overline{\left\langle a,b\right\rangle }}{\left\Vert
a\right\Vert ^{2}}\cdot a\right)
\end{equation*}%
and, by the second condition in $\left( \ref{2.0.17}\right) ,$%
\begin{equation}
\left\vert \mu \left\langle b-\frac{\overline{\left\langle a,b\right\rangle }%
}{\left\Vert a\right\Vert ^{2}}\cdot a,b\right\rangle \right\vert =1.
\label{2.8.17}
\end{equation}%
Since $\left( \ref{2.8.17}\right) $ is clearly equivalent to $\left( \ref%
{2.3.17}\right) $, the theorem is completely proved.
\end{proof}

\subsection{Applications for Sequences and Integrals}

The following particular cases hold.

\begin{enumerate}
\item[\textbf{1.}] If $\mathbf{a},\mathbf{b},\mathbf{x}\in \ell ^{2}\left( 
\mathbb{K}\right) ,$ where $\ell ^{2}\left( \mathbb{K}\right) :=\left\{ 
\mathbf{x=}\left( x_{i}\right) _{i\in \mathbb{N}},\sum_{i=1}^{\infty
}\left\vert x_{i}\right\vert ^{2}<\infty \right\} ,$ with $\mathbf{a},%
\mathbf{b}$ linearly independent and%
\begin{equation*}
\sum_{i=1}^{\infty }x_{i}\overline{a_{i}}=0\text{ and }\left\vert
\sum_{i=1}^{\infty }x_{i}\overline{b_{i}}\right\vert =1,
\end{equation*}%
then one has the inequality%
\begin{equation*}
\sum_{i=1}^{\infty }\left\vert x_{i}\right\vert ^{2}\geq \frac{%
\sum_{i=1}^{\infty }\left\vert a_{i}\right\vert ^{2}}{\sum_{i=1}^{\infty
}\left\vert a_{i}\right\vert ^{2}\sum_{i=1}^{\infty }\left\vert
b_{i}\right\vert ^{2}-\left\vert \sum_{i=1}^{\infty }a_{i}\overline{b_{i}}%
\right\vert ^{2}},
\end{equation*}%
with equality iff 
\begin{equation*}
x_{i}=\mu \left[ b_{i}-\frac{\sum_{k=1}^{\infty }\overline{a_{k}}b_{k}}{%
\sum_{k=1}^{\infty }\left\vert a_{k}\right\vert ^{2}}\cdot a_{i}\right] ,\ \
\ i\in \mathbb{N}
\end{equation*}%
and $\mu \in \mathbb{K}$ with the property%
\begin{equation*}
\left\vert \mu \right\vert =\frac{\sum_{i=1}^{\infty }\left\vert
a_{i}\right\vert ^{2}\sum_{i=1}^{\infty }\left\vert b_{i}\right\vert
^{2}-\left\vert \sum_{i=1}^{\infty }a_{i}\overline{b_{i}}\right\vert ^{2}}{%
\sum_{i=1}^{\infty }\left\vert a_{i}\right\vert ^{2}}.
\end{equation*}

\item[\textbf{2.}] If $f,g,h\in L^{2}\left( \Omega ,m\right) $, where $%
\Omega $ is a measurable space and $L^{2}\left( \Omega ,m\right) :=\left\{
f:\Omega \rightarrow \mathbb{K}\text{, }\int_{\Omega }\left\vert f\left(
x\right) \right\vert ^{2}dm\left( x\right) <\infty \right\} ,$ with $f,g$
linearly independent and%
\begin{equation*}
\int_{\Omega }h\left( x\right) \overline{f\left( x\right) }dm\left( x\right)
=0,\text{ }\left\vert \int_{\Omega }h\left( x\right) \overline{g\left(
x\right) }dm\left( x\right) \right\vert =1,
\end{equation*}%
then one has the inequality%
\begin{multline*}
\int_{\Omega }\left\vert h\left( x\right) \right\vert ^{2}dm\left( x\right)
\\
\geq \frac{\int_{\Omega }\left\vert f\left( x\right) \right\vert
^{2}dm\left( x\right) }{\int_{\Omega }\left\vert f\left( x\right)
\right\vert ^{2}dm\left( x\right) \int_{\Omega }\left\vert g\left( x\right)
\right\vert ^{2}dm\left( x\right) -\left\vert \int_{\Omega }f\left( x\right) 
\overline{g\left( x\right) }dm\left( x\right) \right\vert ^{2}}
\end{multline*}%
with equality iff%
\begin{equation*}
h\left( x\right) =\nu \left[ g\left( x\right) -\frac{\int_{\Omega }f%
\overline{\left( x\right) }g\left( x\right) dm\left( x\right) }{\int_{\Omega
}\left\vert f\left( x\right) \right\vert ^{2}dm\left( x\right) }\cdot
f\left( x\right) \right]
\end{equation*}%
for $m-$ a.e. $x\in \Omega ,$ and $\nu \in \mathbb{K}$ with%
\begin{equation*}
\left\vert \nu \right\vert =\frac{\int_{\Omega }\left\vert f\left( x\right)
\right\vert ^{2}dm\left( x\right) }{\int_{\Omega }\left\vert f\left(
x\right) \right\vert ^{2}dm\left( x\right) \int_{\Omega }\left\vert g\left(
x\right) \right\vert ^{2}dm\left( x\right) -\left\vert \int_{\Omega }f\left(
x\right) \overline{g\left( x\right) }dm\left( x\right) \right\vert }.
\end{equation*}
\end{enumerate}

\newpage

\section{Another Ostrowski Type Inequality}

\subsection{Introduction}

Another result due to Ostrowski which is far less known than the one
incorporated in Theorem \ref{t.1.17}\ and obtained in the same work \cite[p.
130]{2b.18} (see also \cite[p. 94]{1b.18}), is the following one.

\begin{theorem}
\label{t1.2.18}Let $\mathbf{a}$, $\mathbf{b}$ and $\mathbf{x}$ be $n-$tuples
of real numbers with $\mathbf{a}\neq 0$ and%
\begin{equation*}
\sum_{i=1}^{n}a_{i}x_{i}=0\text{ \ and \ }\sum_{i=1}^{n}x_{i}^{2}=1.
\end{equation*}%
Then%
\begin{equation}
\frac{\sum_{i=1}^{n}a_{i}^{2}\sum_{i=1}^{n}b_{i}^{2}-\left(
\sum_{i=1}^{n}a_{i}b_{i}\right) ^{2}}{\sum_{i=1}^{n}a_{i}^{2}}\geq \left(
\sum_{i=1}^{n}b_{i}x_{i}\right) ^{2}.  \label{1.5.18}
\end{equation}%
If $\mathbf{a}$ and $\mathbf{b}$ are not proportional, then the equality
holds in (\ref{1.5.18}) iff%
\begin{equation*}
x_{k}=q\cdot \frac{b_{k}\sum_{i=1}^{n}a_{i}^{2}-a_{k}\sum_{i=1}^{n}a_{i}b_{i}%
}{\left( \sum_{k=1}^{n}a_{k}^{2}\right) ^{\frac{1}{2}}\left[
\sum_{i=1}^{n}a_{i}^{2}\sum_{i=1}^{n}b_{i}^{2}-\left(
\sum_{i=1}^{n}a_{i}b_{i}\right) ^{2}\right] ^{\frac{1}{2}}},\ \ k\in \left\{
1,\dots ,n\right\} ,
\end{equation*}%
with $q\in \left\{ -1,1,\right\} .$
\end{theorem}

The case of equality which was neither mentioned in \cite{1b.18} nor in \cite%
{2b.18} is considered in Remark \ref{r1.18}.

In the present section, by the use of an elementary argument based on\
Schwarz's inequality, a natural generalisation in inner-product spaces of (%
\ref{1.5.18}) is given. The case of equality is analyzed. Applications for
sequences and integrals are also provided.

\subsection{The General Result}

The following theorem holds \cite{20NSSD}.

\begin{theorem}
\label{t2.1.18}Let $\left( H,\left\langle \cdot ,\cdot \right\rangle \right) 
$ be a real or complex inner product space and $a,b\in H$ two linearly
independent vectors. If $x\in H$ is such that%
\begin{equation}
\left\langle x,a\right\rangle =0\ \text{and\ }\left\Vert x\right\Vert =1, 
\tag{i}  \label{i}
\end{equation}%
then%
\begin{equation}
\frac{\left\Vert a\right\Vert ^{2}\left\Vert b\right\Vert ^{2}-\left\vert
\left\langle a,b\right\rangle \right\vert ^{2}}{\left\Vert a\right\Vert ^{2}}%
\geq \left\vert \left\langle x,b\right\rangle \right\vert ^{2}.
\label{2.1.18}
\end{equation}%
The equality holds in (\ref{2.1.18}) iff%
\begin{equation*}
x=\nu \left( b-\frac{\overline{\left\langle a,b\right\rangle }}{\left\Vert
a\right\Vert ^{2}}\cdot a\right) ,
\end{equation*}%
where $\nu \in \mathbb{K}$ $\left( \mathbb{C},\mathbb{R}\right) $ is such
that%
\begin{equation*}
\left\vert \nu \right\vert =\frac{\left\Vert a\right\Vert }{\left[
\left\Vert a\right\Vert ^{2}\left\Vert b\right\Vert ^{2}-\left\vert
\left\langle a,b\right\rangle \right\vert ^{2}\right] ^{\frac{1}{2}}}.
\end{equation*}
\end{theorem}

\begin{proof}
We use Schwarz's inequality in the inner product space $H,$ i.e.,%
\begin{equation}
\left\Vert u\right\Vert ^{2}\left\Vert v\right\Vert ^{2}\geq \left\vert
\left\langle u,v\right\rangle \right\vert ^{2},\ \ u,v\in H  \label{2.4.18}
\end{equation}%
with equality iff there is a scalar $\alpha \in \mathbb{K}$ such that%
\begin{equation*}
u=\alpha v.
\end{equation*}%
If we apply (\ref{2.4.18}) for $u=z-\frac{\left\langle z,c\right\rangle }{%
\left\Vert c\right\Vert ^{2}}\cdot c,$ $v=d-\frac{\left\langle
d,c\right\rangle }{\left\Vert c\right\Vert ^{2}}\cdot c$, where $c\neq 0$
and $c,d,z\in H,$ then we deduce the inequality%
\begin{multline}
\left[ \left\Vert z\right\Vert ^{2}\left\Vert c\right\Vert ^{2}-\left\vert
\left\langle z,c\right\rangle \right\vert ^{2}\right] \left[ \left\Vert
d\right\Vert ^{2}\left\Vert c\right\Vert ^{2}-\left\vert \left\langle
d,c\right\rangle \right\vert ^{2}\right]  \label{2.6.18} \\
\geq \left\vert \left\langle z,d\right\rangle \left\Vert c\right\Vert
^{2}-\left\langle z,c\right\rangle \left\langle c,d\right\rangle \right\vert
^{2}
\end{multline}%
with equality iff there is a $\beta \in \mathbb{K}$ such that%
\begin{equation*}
z=\frac{\left\langle z,c\right\rangle }{\left\Vert c\right\Vert ^{2}}\cdot
c+\beta \left( d-\frac{\left\langle d,c\right\rangle }{\left\Vert
c\right\Vert ^{2}}\cdot c\right) .
\end{equation*}%
If in (\ref{2.6.18}) we choose $z=x,$ $c=a$ and $d=b,$ where $a$ and $x$
statisfy (i), then we deduce%
\begin{equation*}
\left\Vert a\right\Vert ^{2}\left[ \left\Vert a\right\Vert ^{2}\left\Vert
b\right\Vert ^{2}-\left\vert \left\langle a,b\right\rangle \right\vert ^{2}%
\right] \geq \left[ \left\langle x,b\right\rangle \left\Vert a\right\Vert
^{2}\right] ^{2}
\end{equation*}%
which is clearly equivalent to (\ref{2.1.18}).

The equality holds in (\ref{2.1.18}) iff%
\begin{equation*}
x=\nu \left( b-\frac{\overline{\left\langle a,b\right\rangle }}{\left\Vert
a\right\Vert ^{2}}\cdot a\right) ,
\end{equation*}%
where $\nu \in \mathbb{K}$ satisfies the condition%
\begin{equation*}
1=\left\Vert x\right\Vert =\left\vert \nu \right\vert \left\Vert b-\frac{%
\overline{\left\langle a,b\right\rangle }}{\left\Vert a\right\Vert ^{2}}%
\cdot a\right\Vert =\left\vert \nu \right\vert \left[ \frac{\left\Vert
a\right\Vert ^{2}\left\Vert b\right\Vert ^{2}-\left\vert \left\langle
a,b\right\rangle \right\vert ^{2}}{\left\Vert a\right\Vert ^{2}}\right] ^{%
\frac{1}{2}},
\end{equation*}%
and the theorem is thus proved.
\end{proof}

\subsection{Applications for Sequences and Integrals}

The following particular cases hold.

\begin{enumerate}
\item[\textbf{1.}] If $\mathbf{a}$, $\mathbf{b}$, $\mathbf{x}\in \ell
^{2}\left( \mathbb{K}\right) ,$ $\mathbb{K}=\mathbb{C},\mathbb{R},$ where 
\begin{equation*}
\ell ^{2}\left( \mathbb{K}\right) :=\left\{ \mathbf{x}=\left( x_{i}\right)
_{i\in \mathbb{N}},\ \sum_{i=1}^{\infty }\left\vert x_{i}\right\vert
^{2}<\infty \right\}
\end{equation*}%
with $\mathbf{a}$, $\mathbf{b}$ linearly independent and%
\begin{equation*}
\sum_{i=1}^{\infty }x_{i}\overline{a_{i}}=0,\ \ \ \sum_{i=1}^{\infty
}\left\vert x_{i}\right\vert ^{2}=1,
\end{equation*}%
then%
\begin{equation}
\frac{\sum_{i=1}^{\infty }\left\vert a_{i}\right\vert ^{2}\sum_{i=1}^{\infty
}\left\vert b_{i}\right\vert ^{2}-\left\vert \sum_{i=1}^{\infty }a_{i}%
\overline{b_{i}}\right\vert ^{2}}{\sum_{i=1}^{\infty }\left\vert
a_{i}\right\vert ^{2}}\geq \left\vert \sum_{i=1}^{\infty }x_{i}\overline{%
b_{i}}\right\vert ^{2}.  \label{2.9.18}
\end{equation}%
The equality holds in (\ref{2.9.18}) iff%
\begin{equation*}
x_{i}=\nu \left[ b_{i}-\frac{\sum_{k=1}^{\infty }a_{k}\overline{b_{k}}}{%
\sum_{k=1}^{\infty }\left\vert a_{k}\right\vert ^{2}}\cdot a_{i}\right] ,\ \
\ i\in \left\{ 1,2,\dots \right\}
\end{equation*}%
with $\nu \in \mathbb{K}$ such that%
\begin{equation*}
\left\vert \nu \right\vert =\frac{\left( \sum_{k=1}^{\infty }\left\vert
a_{k}\right\vert ^{2}\right) ^{\frac{1}{2}}}{\left[ \sum_{k=1}^{\infty
}\left\vert a_{k}\right\vert ^{2}\sum_{k=1}^{\infty }\left\vert
b_{k}\right\vert ^{2}-\left\vert \sum_{k=1}^{\infty }a_{k}\overline{b_{k}}%
\right\vert ^{2}\right] ^{\frac{1}{2}}}.
\end{equation*}
\end{enumerate}

\begin{remark}
\label{r1.18}The case of equality in (\ref{1.5.18}) is obviously a
particular case of the above. We omit the details.
\end{remark}

\begin{enumerate}
\item[\textbf{2.}] If $f,g,h\in L^{2}\left( \Omega ,m\right) ,$ where $%
\Omega $ is an $m-$measurable space and 
\begin{equation*}
L^{2}\left( \Omega ,m\right) :=\left\{ f:\Omega \rightarrow \mathbb{K},\
\int_{\Omega }\left\vert f\left( x\right) \right\vert ^{2}dm\left( x\right)
<\infty \right\} ,
\end{equation*}%
with $f,g$ being linearly independent and%
\begin{equation*}
\int_{\Omega }h\left( x\right) \overline{f\left( x\right) }dm\left( x\right)
=0,\ \ \ \int_{\Omega }\left\vert h\left( x\right) \right\vert ^{2}dm\left(
x\right) =1,
\end{equation*}%
then%
\begin{multline}
\frac{\int_{\Omega }\left\vert f\left( x\right) \right\vert ^{2}dm\left(
x\right) \int_{\Omega }\left\vert g\left( x\right) \right\vert ^{2}dm\left(
x\right) -\left\vert \int_{\Omega }f\left( x\right) \overline{g\left(
x\right) }dm\left( x\right) \right\vert ^{2}}{\int_{\Omega }\left\vert
f\left( x\right) \right\vert ^{2}dm\left( x\right) }  \label{2.13.18} \\
\geq \left\vert \int_{\Omega }h\left( x\right) \overline{g\left( x\right) }%
dm\left( x\right) \right\vert ^{2}.
\end{multline}%
The equality holds in (\ref{2.13.18}) iff%
\begin{equation*}
h\left( x\right) =\nu \left[ g\left( x\right) -\frac{\int_{\Omega }g\left(
x\right) \overline{f\left( x\right) }dm\left( x\right) }{\int_{\Omega
}\left\vert f\left( x\right) \right\vert ^{2}dm\left( x\right) }f\left(
x\right) \right] \text{ \ for \ a.e. \ }x\in \Omega
\end{equation*}%
and\ $\nu \in \mathbb{K}$ with%
\begin{equation*}
\left\vert \nu \right\vert =\frac{\left( \int_{\Omega }\left\vert f\left(
x\right) \right\vert ^{2}dm\left( x\right) \right) ^{\frac{1}{2}}}{\left[
\int_{\Omega }\left\vert f\left( x\right) \right\vert ^{2}dm\left( x\right)
\int_{\Omega }\left\vert g\left( x\right) \right\vert ^{2}dm\left( x\right)
-\left\vert \int_{\Omega }f\left( x\right) \overline{g\left( x\right) }%
dm\left( x\right) \right\vert ^{2}\right] ^{\frac{1}{2}}}.
\end{equation*}
\end{enumerate}

\newpage

\section{The Wagner Inequality in Inner Product Spaces}

\subsection{Introduction}

In 1965, S.S. Wagner \cite{3b.19} (see also \cite{1b.19} or \cite[p. 85]%
{2b.19}) pointed out the following generalisation of
Cauchy-Bunyakovsky-Schwarz's inequality for real numbers.

\begin{theorem}
\label{1.1.19} Let $\mathbf{a}=\left( a_{1},\dots ,a_{n}\right) $ and $%
\mathbf{b}=\left( b_{1},\dots ,b_{n}\right) $ be two $n$-tuples of real
numbers. Then for any $x\in \left[ 0,1\right] ,$ one has the inequality 
\begin{multline*}
\left( \sum_{k=1}^{n}a_{k}b_{k}+x\cdot \sum_{1\leq i\neq j\leq
n}a_{i}b_{j}\right) ^{2} \\
\leq \left( \sum_{k=1}^{n}a_{k}^{2}+2x\cdot \sum_{1\leq i<j\leq
n}a_{i}a_{j}\right) \cdot \left( \sum_{k=1}^{n}b_{k}^{2}+2x\cdot \sum_{1\leq
i<j\leq n}b_{i}b_{j}\right) .
\end{multline*}
\end{theorem}

For $x=0,$ we recapture the Cauchy-Bunyakovsky-Schwarz's inequality, i.e.,
(see for example \cite[p. 84]{2b.19}) 
\begin{equation*}
\left( \sum_{k=1}^{n}a_{k}b_{k}\right) ^{2}\leq
\sum_{k=1}^{n}a_{k}^{2}\sum_{k=1}^{n}b_{k}^{2},
\end{equation*}%
with equality if and only if there exists a real number $r$ such that $%
a_{k}=rb_{k}$ for each $k\in \{1,\dots ,n\}.$

In this section we extend the above result for sequences of vectors in real
or complex inner product spaces.

\subsection{The Results}

Let $\left( H;\left\langle \cdot ,\cdot \right\rangle \right) $ be an inner
product space over $\mathbb{K}$, where $\mathbb{K}=\mathbb{R}$ or $\mathbb{K}%
=\mathbb{C}$. The following result holds \cite{21NSSD}.

\begin{theorem}
\label{2.1.19} Let $\mathbf{x}=\left( x_{1},\dots ,x_{n}\right) $ and $%
\mathbf{y}=\left( y_{1},\dots ,y_{n}\right) $ be two $n$-tuples of vectors
in $H.$ Then for any $\alpha \in \left[ 0,1\right] $ one has the inequality 
\begin{equation}
\left[ \sum_{k=1}^{n}\func{Re}\left\langle x_{k},y_{k}\right\rangle +\alpha
\cdot \sum_{1\leq i\neq j\leq n}\func{Re}\left\langle
x_{i},y_{j}\right\rangle \right] ^{2}  \label{a.2.19}
\end{equation}%
\begin{equation*}
\leq \left[ \sum_{k=1}^{n}\left\Vert x_{k}\right\Vert ^{2}+2\alpha \cdot
\sum_{1\leq i<j\leq n}\func{Re}\left\langle x_{i},x_{j}\right\rangle \right] %
\left[ \sum_{k=1}^{n}\left\Vert y_{k}\right\Vert ^{2}+2\alpha \cdot
\sum_{1\leq i<j\leq n}\func{Re}\left\langle y_{i},y_{j}\right\rangle \right]
.
\end{equation*}
\end{theorem}

\begin{proof}
Following the proof by P. Flor \cite{1b.19}, we may consider the function $f:%
\mathbb{R}\rightarrow \mathbb{R}$, given by 
\begin{equation}
f\left( t\right) =\left( 1-\alpha \right) \cdot \sum_{k=1}^{n}\left\Vert
tx_{k}-y_{k}\right\Vert ^{2}+\alpha \cdot \left\Vert \sum_{k=1}^{n}\left(
tx_{k}-y_{k}\right) \right\Vert ^{2}.  \label{a.3.19}
\end{equation}%
Then 
\begin{multline}
f\left( t\right) =\left[ \left( 1-\alpha \right) \cdot
\sum_{k=1}^{n}\left\Vert x_{k}\right\Vert ^{2}+\alpha \cdot \left\Vert
\sum_{k=1}^{n}x_{k}\right\Vert ^{2}\right] t^{2}  \label{a.4.19} \\
+2\left[ \left( 1-\alpha \right) \cdot \sum_{k=1}^{n}\func{Re}\left\langle
x_{k},y_{k}\right\rangle +\alpha \cdot \func{Re}\left\langle
\sum_{k=1}^{n}x_{k},\sum_{k=1}^{n}y_{k}\right\rangle \right] t \\
+\left[ \left( 1-\alpha \right) \cdot \sum_{k=1}^{n}\left\Vert
y_{k}\right\Vert ^{2}+\alpha \cdot \left\Vert \sum_{k=1}^{n}y_{k}\right\Vert
^{2}\right] .
\end{multline}%
Observe that 
\begin{equation}
\left\Vert \sum_{k=1}^{n}x_{k}\right\Vert ^{2}=\sum_{k=1}^{n}\left\Vert
x_{k}\right\Vert ^{2}+2\cdot \sum_{1\leq i<j\leq n}\func{Re}\left\langle
x_{i},x_{j}\right\rangle  \label{a.5.19}
\end{equation}%
and 
\begin{equation}
\left\Vert \sum_{k=1}^{n}y_{k}\right\Vert ^{2}=\sum_{k=1}^{n}\left\Vert
y_{k}\right\Vert ^{2}+2\cdot \sum_{1\leq i<j\leq n}\func{Re}\left\langle
y_{i},y_{j}\right\rangle .  \label{a.6.19}
\end{equation}%
Also 
\begin{equation}
\func{Re}\left\langle \sum_{k=1}^{n}x_{k},\sum_{k=1}^{n}y_{k}\right\rangle
=\sum_{k=1}^{n}\func{Re}\left\langle x_{k},y_{k}\right\rangle +\sum_{1\leq
i\neq j\leq n}\func{Re}\left\langle x_{i},y_{j}\right\rangle .
\label{a.7.19}
\end{equation}%
Using $\left( \ref{a.4.19}\right) -\left( \ref{a.7.19}\right) ,$ we deduce
that 
\begin{multline}
f\left( t\right) =\left[ \sum_{k=1}^{n}\left\Vert x_{k}\right\Vert
^{2}+2\alpha \cdot \sum_{1\leq i<j\leq n}\func{Re}\left\langle
x_{i},x_{j}\right\rangle \right] t^{2}  \label{a.8.19} \\
+2\left[ \sum_{k=1}^{n}\func{Re}\left\langle x_{k},y_{k}\right\rangle
+\alpha \cdot \sum_{1\leq i\neq j\leq n}\func{Re}\left\langle
x_{i},y_{j}\right\rangle \right] t \\
+\left[ \sum_{k=1}^{n}\left\Vert y_{k}\right\Vert ^{2}+2\alpha \cdot
\sum_{1\leq i<j\leq n}\func{Re}\left\langle y_{i},y_{j}\right\rangle \right]
.
\end{multline}%
Since, by $\left( \ref{a.3.19}\right) ,$ $f\left( t\right) \geq 0$ for any $%
t\in \mathbb{R}\mathbf{,}$ it follows that the discriminant of the quadratic
function given by $\left( \ref{a.8.19}\right) $ is negative, which is
clearly equivalent with the desired inequality $\left( \ref{a.2.19}\right) .$
\end{proof}

One may obtain an interesting inequality if $\mathbf{x}$ and $\mathbf{y}$
are assumed to incorporate orthogonal vectors.

\begin{corollary}
\label{2.2.19} Assume that $\{x_{i}\}_{i=1,\dots ,n}$ are orthogonal, i.e., $%
x_{i}\perp x_{j}$ for any $i,j\in \{1,\dots ,n\},$ $i\neq j;$ and $%
\{y_{i}\}_{i=1,\dots ,n}$ are also orthogonal in the real inner product
space $\left( H;\left\langle \cdot ,\cdot \right\rangle \right) .$ Then 
\begin{equation*}
\sup_{\alpha \in \left[ 0,1\right] }\left[ \sum_{k=1}^{n}\left\langle
x_{k},y_{k}\right\rangle +\alpha \cdot \sum_{1\leq i\neq j\leq
n}\left\langle x_{i},y_{j}\right\rangle \right] ^{2}\leq
\sum_{k=1}^{n}\left\Vert x_{k}\right\Vert ^{2}\sum_{k=1}^{n}\left\Vert
y_{k}\right\Vert ^{2}.
\end{equation*}
\end{corollary}

\subsection{Applications}

\begin{enumerate}
\item[\textbf{1}.] If we assume that $H=\mathbb{C}\mathbf{,}$ with the inner
product $\left\langle x,y\right\rangle =x\cdot \bar{y},$ then by $\left( \ref%
{a.2.19}\right) $ we may deduce the following Wagner type inequality for
complex numbers 
\begin{multline*}
\left[ \sum_{k=1}^{n}\func{Re}\left( a_{k}\bar{b}_{k}\right) +\alpha \cdot
\sum_{1\leq i\neq j\leq n}\func{Re}\left( a_{i}\bar{b}_{j}\right) \right]
^{2} \\
\leq \left[ \sum_{k=1}^{n}\left\vert a_{k}\right\vert ^{2}+2\alpha \cdot
\sum_{1\leq i<j\leq n}\func{Re}\left( a_{i}\bar{a}_{j}\right) \right] \\
\times \left[ \sum_{k=1}^{n}\left\vert b_{k}\right\vert ^{2}+2\alpha \cdot
\sum_{1\leq i<j\leq n}\func{Re}\left( b_{i}\bar{b}_{j}\right) \right] ,
\end{multline*}%
where $\alpha \in \left[ 0,1\right] $ and $\mathbf{a}=\left( a_{1},\dots
,a_{n}\right) $, $\mathbf{b}=\left( b_{1},\dots ,b_{n}\right) \in \mathbb{C}%
^{n}$.

\item[\textbf{2}.] Consider the Hilbert space $L_{2}\left( \Omega ,\mu
\right) :=\left\{ f:\Omega \rightarrow \mathbb{C}\mathbf{,\ }\int_{\Omega
}\left\vert f\left( x\right) \right\vert ^{2}d\mu \left( x\right) <\infty
\right\} ,$ where $\Omega $ is a $\mu $-measurable space and $\mu :\Omega
\rightarrow \left[ 0,\infty \right] $ is a positive measure on $\Omega .$
Then for $H=L_{2}\left( \Omega ,\mu \right) $ and since the inner product
generating the norm is given by 
\begin{equation*}
\left\langle f,g\right\rangle =\int_{\Omega }f\left( x\right) \bar{g}\left(
x\right) d\mu \left( x\right) ,
\end{equation*}%
we get the inequality 
\begin{multline*}
\left[ \sum_{k=1}^{n}\int_{\Omega }\func{Re}\left( f_{k}\left( x\right) \bar{%
g}_{k}\left( x\right) \right) d\mu \left( x\right) +\alpha \cdot \sum_{1\leq
i\neq j\leq n}\int_{\Omega }\func{Re}\left( f_{i}\left( x\right) \bar{g}%
_{j}\left( x\right) \right) d\mu \left( x\right) \right] ^{2} \\
\leq \left[ \sum_{k=1}^{n}\int_{\Omega }\left\vert f_{k}\left( x\right)
\right\vert ^{2}d\mu \left( x\right) +2\alpha \cdot \sum_{1\leq i<j\leq
n}\int_{\Omega }\func{Re}\left( f_{i}\left( x\right) \bar{f}_{j}\left(
x\right) \right) d\mu \left( x\right) \right] \\
\times \left[ \sum_{k=1}^{n}\int_{\Omega }\left\vert g_{k}\left( x\right)
\right\vert ^{2}d\mu \left( x\right) +2\alpha \cdot \sum_{1\leq i<j\leq
n}\int_{\Omega }\func{Re}\left( g_{i}\left( x\right) \bar{g}_{j}\left(
x\right) \right) d\mu \left( x\right) \right] ,
\end{multline*}%
where $f_{i},g_{i}\in L_{2}\left( \Omega ,\mu \right) ,i\in \{1,\dots ,n\}$
and $\alpha \in \left[ 0,1\right] .$
\end{enumerate}

\newpage

\section{A Monotoniciy Property of Bessel's Inequality}

Let $X$ be a linear space over the real or complex number field $\mathbb{K}$%
. A mapping $\left( \cdot ,\cdot \right) :X\times X\rightarrow \mathbb{K}$
is said to be of \textit{positive hermitian form }if the following
conditions are satisfied:

\begin{enumerate}
\item[(i)] $\left( \alpha x+\beta y,z\right) =\alpha \left( x,z\right)
+\beta \left( y,z\right) $ for all $x,y,z\in X$ and $\alpha ,\beta \in 
\mathbb{K}$;

\item[(ii)] $\left( y,x\right) =\overline{\left( x,y\right) }$ for all $%
x,y\in X$;

\item[(iii)] $\left( x,x\right) \geq 0$ for all $x\in X.$
\end{enumerate}

If $\left\Vert x\right\Vert :=\left( x,x\right) ^{\frac{1}{2}}$ denotes the
semi-norm associated to this form and $\left( e_{i}\right) _{i\in I}$ is an
orthonormal family of vectors in $X$, that is, $\left( e_{i},e_{j}\right)
=\delta _{ij}\;\left( i,j\in I\right) $, then one has \cite{Y.20}: 
\begin{equation}
\left\Vert x\right\Vert ^{2}\geq \sum_{i\in I}\left\vert \left(
x,e_{i}\right) \right\vert ^{2}\;\;\;\text{for all }x\in X,  \label{e1.20}
\end{equation}%
which is well known in the literature as Bessel's inequality.

The main aim of the section is to point out an improvement for this result
as follows \cite{22NSSD}.

\begin{theorem}
\label{tnew}Let $X$ be a linear space and $\left( \cdot ,\cdot \right)
_{2},\left( \cdot ,\cdot \right) _{1}$ two hermitian forms on $X$ such that $%
\left\Vert \cdot \right\Vert _{2}$ is greater than or equal to $\left\Vert
\cdot \right\Vert _{1}$, that is, $\left\Vert x\right\Vert _{2}\geq
\left\Vert x\right\Vert _{1}$ for all $x\in X$. Assume that $\left(
e_{i}\right) _{i\in I}$ is an orthonormal family in $\left( X;\left( \cdot
,\cdot \right) _{2}\right) $ and $\left( f_{j}\right) _{j\in J}$ is an
orthonormal family in $\left( X;\left( \cdot ,\cdot \right) _{1}\right) $
such that for any $i\in I$ there exists a finite $K\subset J$ such that 
\begin{equation}
e_{i}=\sum_{j\in K}\alpha _{j}f_{j},\;\;\alpha _{j}\in \mathbb{K},\ \ \ \
\left( j\in K\right) ,  \label{e2.20}
\end{equation}%
then one has the inequality: 
\begin{equation}
\left\Vert x\right\Vert _{2}^{2}-\sum_{i\in I}\left\vert \left(
x,e_{i}\right) _{2}\right\vert ^{2}\geq \left\Vert x\right\Vert
_{1}^{2}-\sum_{j\in J}\left\vert \left( x,f_{j}\right) _{1}\right\vert
^{2}\geq 0,\ \text{for\ all}\ x\in X.  \label{e3.20}
\end{equation}
\end{theorem}

In order to prove this, we require the following lemma.

\begin{lemma}
\label{lemnew}Let $X$ be a linear space endowed with a positive hermitian
form $\left( \cdot ,\cdot \right) $ and $\left( g_{k}\right) ,$ $k\in
\{1,\dots ,n\}$ be an orthonormal family in $\left( X;\left( \cdot ,\cdot
\right) \right) .$ Then 
\begin{equation*}
\left\Vert x-\sum_{k=1}^{n}\lambda _{k}g_{k}\right\Vert ^{2}\geq \left\Vert
x\right\Vert ^{2}-\sum_{k=1}^{n}\left\vert \left( x,g_{k}\right) \right\vert
^{2}\geq 0,
\end{equation*}%
for\ all\thinspace \thinspace $\lambda _{k}\in \mathbb{K},$ $k\in \{1\ldots
,n\}$\thinspace \thinspace and\thinspace \thinspace \thinspace $x\in X.$
\end{lemma}

The proof follows by mathematical induction.

\begin{proof}[Proof of Theorem \protect\ref{tnew}]
Let $H$ be a finite subset of $I$. Since $\left\Vert \cdot \right\Vert _{2}$
is greater than $\left\Vert \cdot \right\Vert _{1}$, we have: 
\begin{equation*}
\left\Vert x\right\Vert _{2}^{2}-\sum_{i\in H}\left\vert \left(
x,e_{i}\right) _{2}\right\vert ^{2}=\left\Vert x-\sum_{i\in H}\left(
x,e_{i}\right) _{2}e_{i}\right\Vert _{2}^{2}\geq \left\Vert x-\sum_{i\in
H}\left( x,e_{i}\right) _{2}e_{i}\right\Vert _{1}^{2},\;\;x\in X.
\end{equation*}%
Since, by (\ref{e2.20}), we may state that for any $i\in H$ there exists a
finite $K\subset J$ with 
\begin{equation*}
e_{i}=\sum_{j\in K}\left( e_{i},f_{j}\right) _{1}f_{j},
\end{equation*}%
we have, for all $x\in X$ 
\begin{align*}
\left\Vert x-\sum_{i\in H}\left( x,e_{i}\right) _{2}e_{i}\right\Vert
_{1}^{2}& =\left\Vert x-\sum_{i\in H}\left( x,e_{i}\right) _{2}\sum_{j\in
K}\left( e_{i},f_{j}\right) _{1}f_{j}\right\Vert _{1}^{2} \\
& =\left\Vert x-\sum_{j\in K}\left( \sum_{i\in H}\left( x,e_{i}\right)
_{2}e_{i},f_{j}\right) _{1}f_{j}\right\Vert _{1}^{2}.
\end{align*}%
Applying Lemma \ref{lemnew}, we can conclude that 
\begin{equation*}
\left\Vert x-\sum_{j\in K}\lambda _{j}f_{j}\right\Vert _{1}^{2}\geq
\left\Vert x\right\Vert _{1}^{2}-\sum_{j\in K}\left\vert \left(
x,f_{j}\right) _{1}\right\vert ^{2},\;\;x\in X,
\end{equation*}%
where 
\begin{equation*}
\lambda _{j}=\left( \sum_{i\in H}\left( x,e_{i}\right)
_{2}e_{i},f_{j}\right) _{1}\in \mathbb{K},\left( j\in K\right) .
\end{equation*}%
Consequently, we have 
\begin{equation*}
\left\Vert x\right\Vert _{2}^{2}-\sum_{i\in H}\left\vert \left(
x,e_{i}\right) _{2}\right\vert ^{2}\geq \left\Vert x\right\Vert
_{1}^{2}-\sum_{j\in K}\left\vert \left( x,f_{j}\right) _{1}\right\vert
^{2}\geq \left\Vert x\right\Vert _{1}^{2}-\sum_{j\in J}\left\vert \left(
x,f_{j}\right) _{1}\right\vert ^{2}
\end{equation*}%
for all $x\in X$ and $H$ a finite subset of $I$, from which (\ref{e3.20})
results.
\end{proof}

\begin{corollary}
\label{cnew}Let $\left\Vert \cdot \right\Vert _{1},\left\Vert \cdot
\right\Vert _{2}:X\rightarrow \mathbb{R}_{+}$ be as above. Then for all $%
x,y\in X$, we have the inequality: 
\begin{equation}
\left\Vert x\right\Vert _{2}^{2}\left\Vert y\right\Vert _{2}^{2}-\left\vert
\left( x,y\right) _{2}\right\vert ^{2}\geq \left\Vert x\right\Vert
_{1}^{2}\left\Vert y\right\Vert _{1}^{2}-\left\vert \left( x,y\right)
_{1}\right\vert ^{2}\geq 0,  \label{e4.20}
\end{equation}%
which is an improvement of the well known Cauchy-Schwarz inequality.
\end{corollary}

\begin{remark}
\label{rnew}For a different proof of (\ref{e4.20}), see also \cite{DM.20} or 
\cite{DMP.20}.
\end{remark}

Now, we will give some natural applications of the above theorem.

\begin{enumerate}
\item Let $\left( X;\left( \cdot ,\cdot \right) \right) $ be an inner
product space and $\left( e_{i}\right) _{i\in I}$ an orthonormal family in $%
X $. Assume that $A:X\rightarrow X$ is a linear operator such that $%
\left\Vert Ax\right\Vert \leq \left\Vert x\right\Vert $ for all $x\in X$ and 
$\left( Ae_{i},Ae_{j}\right) =\delta _{ij}$ for all $i,j\in I$. Then one has
the inequality 
\begin{equation*}
\left\Vert x\right\Vert ^{2}-\sum_{i\in I}\left\vert \left( x,e_{i}\right)
\right\vert ^{2}\geq \left\Vert Ax\right\Vert ^{2}-\sum_{i\in I}\left\vert
\left( Ax,Ae_{i}\right) \right\vert ^{2}\geq 0,\ \text{for\ all}\ x\in X.
\end{equation*}

\item If $A:X\rightarrow X$ is such that $\left\Vert Ax\right\Vert \geq
\left\Vert x\right\Vert $ for all $x\in X$, then, with the previous
assumptions, we also have 
\begin{equation*}
0\leq \left\Vert x\right\Vert ^{2}-\sum_{i\in I}\left\vert \left(
x,e_{i}\right) \right\vert ^{2}\leq \left\Vert Ax\right\Vert ^{2}-\sum_{i\in
I}\left\vert \left( Ax,Ae_{i}\right) \right\vert ^{2},\ \text{for\ all}\
x\in X.
\end{equation*}

\item Suppose that $A:X\rightarrow X$ is a symmetric positive definite
operator with $\left( Ax,x\right) \geq \left\Vert x\right\Vert ^{2}$ for all 
$x\in X$. If $\left( e_{i}\right) _{i\in I}$ is an orthonormal family in $X$
such that $\left( Ae_{i},Ae_{j}\right) =\delta _{ij}$ for all $i,j\in I$,
then one has the inequality 
\begin{equation*}
0\leq \left\Vert x\right\Vert ^{2}-\sum_{i\in I}\left\vert \left(
x,e_{i}\right) \right\vert ^{2}\leq \left( Ax,x\right) -\sum_{i\in
I}\left\vert \left( Ax,e_{i}\right) \right\vert ^{2},
\end{equation*}%
for any $x\in X.$
\end{enumerate}

\newpage

\section{Other Bombieri Type Inequalities}

\subsection{Introduction}

In 1971, E. Bombieri \cite{1b.23} gave the following generalisation of
Bessel's inequality: 
\begin{equation}
\sum_{i=1}^{n}\left\vert \left( x,y_{i}\right) \right\vert ^{2}\leq
\left\Vert x\right\Vert ^{2}\max_{1\leq i\leq n}\left\{
\sum_{j=1}^{n}\left\vert \left( y_{i},y_{j}\right) \right\vert \right\} ,
\label{1.1.23}
\end{equation}%
where $x,y_{1},\dots ,y_{n}$ are vectors in the inner product space $\left(
H;\left( \cdot ,\cdot \right) \right) .$

It is obvious that if $\left( y_{i}\right) _{1\leq i\leq n}=\left(
e_{i}\right) _{1\leq i\leq n},$ where $\left( e_{i}\right) _{1\leq i\leq n}$
are orthornormal vectors in $H,$ i.e., $\left( e_{i},e_{j}\right) =\delta
_{ij}$ $\left( i,j=1,\dots ,n\right) ,$ where $\delta _{ij}$ is the
Kronecker delta, then (\ref{1.1.23}) provides Bessel's inequality 
\begin{equation*}
\sum_{i=1}^{n}\left\vert \left( x,e_{i}\right) \right\vert ^{2}\leq
\left\Vert x\right\Vert ^{2},\ \ \ x\in H.
\end{equation*}

In this section we point out some Bombieri type inequalities that complement
the results obtained in Chapter \ref{chap4}.

\subsection{The Results}

The following lemma, which is of interest in itself, holds \cite{23NSSD}.

\begin{lemma}
\label{l2.1.23}Let $z_{1},\dots ,z_{n}\in H$ and $\alpha _{1},\dots ,\alpha
_{n}\in \mathbb{K}$. Then one has the inequalities: 
\begin{align*}
& \left\Vert \sum_{i=1}^{n}\alpha _{i}z_{i}\right\Vert ^{2} \\
& \leq \left\{ 
\begin{array}{l}
\max\limits_{1\leq i\leq n}\left\vert \alpha _{i}\right\vert
^{2}\sum_{i,j=1}^{n}\left\vert \left( z_{i},z_{j}\right) \right\vert ; \\ 
\\ 
\left( \sum_{i=1}^{n}\left\vert \alpha _{i}\right\vert ^{p}\right) ^{\frac{2%
}{p}}\left( \sum_{i,j=1}^{n}\left\vert \left( z_{i},z_{j}\right) \right\vert
^{q}\right) ^{\frac{1}{q}},\ \ \text{where \ }p>1,\ \frac{1}{p}+\frac{1}{q}%
=1; \\ 
\\ 
\left( \sum_{i=1}^{n}\left\vert \alpha _{i}\right\vert \right) ^{2}\
\max\limits_{1\leq i,j\leq n}\left\vert \left( z_{i},z_{j}\right)
\right\vert ;%
\end{array}%
\right. \\
& \leq \left\{ 
\begin{array}{l}
\max\limits_{1\leq i\leq n}\left\vert \alpha _{i}\right\vert ^{2}\left(
\sum_{i=1}^{n}\left\Vert z_{i}\right\Vert \right) ^{2}; \\ 
\\ 
\left( \sum_{i=1}^{n}\left\vert \alpha _{i}\right\vert ^{p}\right) ^{\frac{2%
}{p}}\left( \sum_{i=1}^{n}\left\Vert z_{i}\right\Vert ^{q}\right) ^{\frac{2}{%
q}},\ \ \text{where \ }p>1,\ \frac{1}{p}+\frac{1}{q}=1; \\ 
\\ 
\left( \sum_{i=1}^{n}\left\vert \alpha _{i}\right\vert \right) ^{2}\
\max\limits_{1\leq i\leq n}\left\Vert z_{i}\right\Vert ^{2}.%
\end{array}%
\right.
\end{align*}
\end{lemma}

\begin{proof}
We observe that 
\begin{align*}
\left\Vert \sum_{i=1}^{n}\alpha _{i}z_{i}\right\Vert ^{2}& =\left(
\sum_{i=1}^{n}\alpha _{i}z_{i},\sum_{j=1}^{n}\alpha _{j}z_{j}\right)
=\sum_{i=1}^{n}\sum_{j=1}^{n}\alpha _{i}\overline{\alpha _{j}}\left(
z_{i},z_{j}\right) \\
& =\left\vert \sum_{i=1}^{n}\sum_{j=1}^{n}\alpha _{i}\overline{\alpha _{j}}%
\left\vert \left( z_{i},z_{j}\right) \right\vert \right\vert \leq
\sum_{i=1}^{n}\sum_{j=1}^{n}\left\vert \alpha _{i}\right\vert \left\vert
\alpha _{j}\right\vert \left\vert \left( z_{i},z_{j}\right) \right\vert =:M.
\end{align*}%
Firstly, we have 
\begin{align*}
M& \leq \max\limits_{1\leq i,j\leq n}\left\{ \left\vert \alpha
_{i}\right\vert \left\vert \alpha _{j}\right\vert \right\}
\sum\limits_{i,j=1}^{n}\left\vert \left( z_{i},z_{j}\right) \right\vert \\
& =\max\limits_{1\leq i\leq n}\left\vert \alpha _{i}\right\vert
^{2}\sum\limits_{i,j=1}^{n}\left\vert \left( z_{i},z_{j}\right) \right\vert .
\end{align*}%
Secondly, by the H\"{o}lder inequality for double sums, we obtain 
\begin{align*}
M& \leq \left[ \sum\limits_{i,j=1}^{n}\left( \left\vert \alpha
_{i}\right\vert \left\vert \alpha _{j}\right\vert \right) ^{p}\right] ^{%
\frac{1}{p}}\left( \sum\limits_{i,j=1}^{n}\left\vert \left(
z_{i},z_{j}\right) \right\vert ^{q}\right) ^{\frac{1}{q}} \\
& =\left( \sum_{i=1}^{n}\left\vert \alpha _{i}\right\vert
^{p}\sum_{j=1}^{n}\left\vert \alpha _{j}\right\vert ^{p}\right) ^{\frac{1}{p}%
}\left( \sum\limits_{i,j=1}^{n}\left\vert \left( z_{i},z_{j}\right)
\right\vert ^{q}\right) ^{\frac{1}{q}} \\
& =\left( \sum\limits_{i=1}^{n}\left\vert \alpha _{i}\right\vert ^{p}\right)
^{\frac{2}{p}}\left( \sum\limits_{i,j=1}^{n}\left\vert \left(
z_{i},z_{j}\right) \right\vert ^{q}\right) ^{\frac{1}{q}},
\end{align*}%
where $p>1,$ $\frac{1}{p}+\frac{1}{q}=1.$

Finally, we have 
\begin{equation*}
M\leq \max\limits_{1\leq i,j\leq n}\left| \left( z_{i},z_{j}\right) \right|
\sum\limits_{i,j=1}^{n}\left| \alpha _{i}\right| \left| \alpha _{j}\right|
=\left( \sum\limits_{i=1}^{n}\left| \alpha _{i}\right| \right)
^{2}\max\limits_{1\leq i,j\leq n}\left| \left( z_{i},z_{j}\right) \right|
\end{equation*}
and the first part of the lemma is proved.

The second part is obvious on taking into account, by Schwarz's inequality
in $H$, that we have 
\begin{equation*}
\left\vert \left( z_{i},z_{j}\right) \right\vert \leq \left\Vert
z_{i}\right\Vert \left\Vert z_{j}\right\Vert ,
\end{equation*}%
for any $i,j\in \left\{ 1,\dots ,n\right\} .$ We omit the details.
\end{proof}

\begin{corollary}
\label{c2.2.23}With the assumptions in Lemma \ref{l2.1.23}, one has 
\begin{align}
\left\Vert \sum_{i=1}^{n}\alpha _{i}z_{i}\right\Vert ^{2}& \leq
\sum\limits_{i=1}^{n}\left\vert \alpha _{i}\right\vert ^{2}\left(
\sum\limits_{i,j=1}^{n}\left\vert \left( z_{i},z_{j}\right) \right\vert
^{2}\right) ^{\frac{1}{2}}  \label{2.3.23} \\
& \leq \sum\limits_{i=1}^{n}\left\vert \alpha _{i}\right\vert
^{2}\sum_{i=1}^{n}\left\Vert z_{i}\right\Vert ^{2}.  \notag
\end{align}
\end{corollary}

The proof follows by Lemma \ref{l2.1.23} on choosing $p=q=2.$

Note also that $\left( \ref{2.3.23}\right) $ provides a refinement of the
well known Cauchy-Bunyakovsky-Schwarz inequality for sequences of vectors in
inner product spaces, namely 
\begin{equation*}
\left\Vert \sum_{i=1}^{n}\alpha _{i}z_{i}\right\Vert ^{2}\leq
\sum\limits_{i=1}^{n}\left\vert \alpha _{i}\right\vert
^{2}\sum_{i=1}^{n}\left\Vert z_{i}\right\Vert ^{2}.
\end{equation*}

The following lemma also holds \cite{23NSSD}.

\begin{lemma}
\label{l2.4.23}Let $x,y_{1},\dots ,y_{n}\in H$ and $c_{1},\dots ,c_{n}\in 
\mathbb{K}$. Then one has the inequalities: 
\begin{align}
& \left\vert \sum\limits_{i=1}^{n}c_{i}\left( x,y_{i}\right) \right\vert ^{2}
\label{2.4a.23} \\
& \leq \left\Vert x\right\Vert ^{2}\times \left\{ 
\begin{array}{l}
\ \max\limits_{1\leq i\leq n}\left\vert c_{i}\right\vert
^{2}\sum_{i,j=1}^{n}\left\vert \left( y_{i},y_{j}\right) \right\vert ; \\ 
\\ 
\left( \sum_{i=1}^{n}\left\vert c_{i}\right\vert ^{p}\right) ^{\frac{2}{p}%
}\left( \sum_{i,j=1}^{n}\left\vert \left( y_{i},y_{j}\right) \right\vert
^{q}\right) ^{\frac{1}{q}},\ \ \  \\ 
\hfill \ \text{where}\ \ p>1,\ \frac{1}{p}+\frac{1}{q}=1; \\ 
\\ 
\left( \sum_{i=1}^{n}\left\vert c_{i}\right\vert \right)
^{2}\max\limits_{1\leq i,j\leq n}\left\vert \left( y_{i},y_{j}\right)
\right\vert ;%
\end{array}%
\right.  \notag \\
& \leq \left\Vert x\right\Vert ^{2}\times \left\{ 
\begin{array}{l}
\ \max\limits_{1\leq i\leq n}\left\vert c_{i}\right\vert ^{2}\left(
\sum_{i=1}^{n}\left\Vert y_{i}\right\Vert \right) ^{2}; \\ 
\\ 
\left( \sum_{i=1}^{n}\left\vert c_{i}\right\vert ^{p}\right) ^{\frac{2}{p}%
}\left( \sum_{i=1}^{n}\left\Vert y_{i}\right\Vert ^{q}\right) ^{\frac{2}{q}%
},\ \ \ \  \\ 
\hfill \text{where}\ \ p>1,\ \frac{1}{p}+\frac{1}{q}=1; \\ 
\\ 
\left( \sum_{i=1}^{n}\left\vert c_{i}\right\vert \right)
^{2}\max\limits_{1\leq i\leq n}\left\Vert y_{i}\right\Vert ^{2}.%
\end{array}%
\right.  \notag
\end{align}
\end{lemma}

\begin{proof}
We have, by Schwarz's inequality in the inner product $\left( H;\left( \cdot
,\cdot \right) \right) ,$ that 
\begin{equation*}
\left\vert \sum\limits_{i=1}^{n}c_{i}\left( x,y_{i}\right) \right\vert
^{2}=\left\vert \left( x,\sum\limits_{i=1}^{n}\overline{c_{i}}y_{i}\right)
\right\vert ^{2}\leq \left\Vert x\right\Vert ^{2}\left\Vert
\sum\limits_{i=1}^{n}\overline{c_{i}}y_{i}\right\Vert ^{2}.
\end{equation*}%
Now, applying Lemma \ref{l2.1.23} for $\alpha _{i}=\overline{c_{i}}$, $%
z_{i}=y_{i}$ $\left( i=1,\dots ,n\right) ,$ the inequality (\ref{2.4a.23})
is proved.
\end{proof}

\begin{corollary}
\label{c2.5.23}With the assumptions in Lemma \ref{l2.4.23}, one has 
\begin{align}
\left\vert \sum\limits_{i=1}^{n}c_{i}\left( x,y_{i}\right) \right\vert ^{2}&
\leq \left\Vert x\right\Vert ^{2}\sum\limits_{i=1}^{n}\left\vert
c_{i}\right\vert ^{2}\left( \sum\limits_{i,j=1}^{n}\left\vert \left(
y_{i},y_{j}\right) \right\vert ^{2}\right) ^{\frac{1}{2}}  \label{2.5.23} \\
& \leq \left\Vert x\right\Vert ^{2}\sum\limits_{i=1}^{n}\left\vert
c_{i}\right\vert ^{2}\sum\limits_{i=1}^{n}\left\Vert y_{i}\right\Vert ^{2}. 
\notag
\end{align}
\end{corollary}

The proof follows by Lemma \ref{l2.4.23}, on choosing $p=q=2.$

\begin{remark}
\label{r2.6.23}The inequality (\ref{2.5.23}) was firstly obtained in \cite[%
Inequality (7)]{2b.23}.
\end{remark}

The following theorem incorporating three Bombieri type inequalities holds 
\cite{23NSSD}.

\begin{theorem}
\label{t2.7.23}Let $x,y_{1},\dots ,y_{n}\in H.$ Then one has the
inequalities: 
\begin{multline}
\sum\limits_{i=1}^{n}\left\vert \left( x,y_{i}\right) \right\vert ^{2}
\label{2.6.23} \\
\leq \left\Vert x\right\Vert \times \left\{ 
\begin{array}{l}
\ \max\limits_{1\leq i\leq n}\left\vert \left( x,y_{i}\right) \right\vert
\left( \sum_{i,j=1}^{n}\left\vert \left( y_{i},y_{j}\right) \right\vert
\right) ^{\frac{1}{2}}; \\ 
\\ 
\left( \sum_{i=1}^{n}\left\vert \left( x,y_{i}\right) \right\vert
^{p}\right) ^{\frac{1}{p}}\left( \sum_{i,j=1}^{n}\left\vert \left(
y_{i},y_{j}\right) \right\vert ^{q}\right) ^{\frac{1}{2q}},\ \ \ \  \\ 
\hfill \text{where}\ \ p>1,\ \frac{1}{p}+\frac{1}{q}=1; \\ 
\\ 
\sum_{i=1}^{n}\left\vert \left( x,y_{i}\right) \right\vert \
\max\limits_{1\leq i,j\leq n}\left\vert \left( y_{i},y_{j}\right)
\right\vert ^{\frac{1}{2}}.%
\end{array}%
\right.
\end{multline}
\end{theorem}

\begin{proof}
Choosing $c_{i}=\overline{\left( x,y_{i}\right) }$ $\left( i=1,\dots
,n\right) $ in (\ref{2.4a.23}) we deduce 
\begin{multline}
\left( \sum\limits_{i=1}^{n}\left\vert \left( x,y_{i}\right) \right\vert
^{2}\right) ^{2}  \label{2.7.23} \\
\leq \left\Vert x\right\Vert ^{2}\times \left\{ 
\begin{array}{l}
\max\limits_{1\leq i\leq n}\left\vert \left( x,y_{i}\right) \right\vert
^{2}\left( \sum_{i,j=1}^{n}\left\vert \left( y_{i},y_{j}\right) \right\vert
\right) ; \\ 
\\ 
\left( \sum_{i=1}^{n}\left\vert \left( x,y_{i}\right) \right\vert
^{p}\right) ^{\frac{2}{p}}\left( \sum_{i,j=1}^{n}\left\vert \left(
y_{i},y_{j}\right) \right\vert ^{q}\right) ^{\frac{1}{q}}, \\ 
\hfill \ \ \ \ \text{where}\ \ p>1,\ \frac{1}{p}+\frac{1}{q}=1; \\ 
\\ 
\left( \sum_{i=1}^{n}\left\vert \left( x,y_{i}\right) \right\vert \right)
^{2}\ \max\limits_{1\leq i,j\leq n}\left\vert \left( y_{i},y_{j}\right)
\right\vert ;%
\end{array}%
\right.
\end{multline}%
which, by taking the square root, is clearly equivalent to (\ref{2.6.23}).
\end{proof}

\begin{remark}
If $\left( y_{i}\right) _{1\leq i\leq n}=\left( e_{i}\right) _{1\leq i\leq
n},$ where $\left( e_{i}\right) _{1\leq i\leq n}$ are orthornormal vectors
in $H,$ then by (\ref{2.6.23}) we deduce 
\begin{equation}
\sum\limits_{i=1}^{n}\left\vert \left( x,e_{i}\right) \right\vert ^{2}\leq
\left\Vert x\right\Vert \left\{ 
\begin{array}{l}
\sqrt{n}\ \max\limits_{1\leq i\leq n}\left\vert \left( x,e_{i}\right)
\right\vert ; \\ 
\\ 
n^{\frac{1}{2q}}\left( \sum_{i=1}^{n}\left\vert \left( x,e_{i}\right)
\right\vert ^{p}\right) ^{\frac{1}{p}},\  \\ 
\hfill \text{where}\ \ p>1,\ \frac{1}{p}+\frac{1}{q}=1; \\ 
\\ 
\sum_{i=1}^{n}\left\vert \left( x,e_{i}\right) \right\vert .%
\end{array}%
\right.  \label{2.7.a.23}
\end{equation}
\end{remark}

If in (\ref{2.7.23}) we take $p=q=2,$ then we obtain the following
inequality which was formulated in \cite[p. 81]{2b.23}.

\begin{corollary}
\label{c2.8.23}With the assumptions in Theorem \ref{t2.7.23}, we have: 
\begin{equation}
\sum\limits_{i=1}^{n}\left\vert \left( x,y_{i}\right) \right\vert ^{2}\leq
\left\Vert x\right\Vert ^{2}\left( \sum\limits_{i,j=1}^{n}\left\vert \left(
y_{i},y_{j}\right) \right\vert ^{2}\right) ^{\frac{1}{2}}.  \label{2.8.23}
\end{equation}
\end{corollary}

\begin{remark}
\label{r2.9.23}Observe, that by the monotonicity of power means, we may
write 
\begin{equation*}
\left( \frac{\sum_{i=1}^{n}\left\vert \left( x,y_{i}\right) \right\vert ^{p}%
}{n}\right) ^{\frac{1}{p}}\leq \left( \frac{\sum_{i=1}^{n}\left\vert \left(
x,y_{i}\right) \right\vert ^{2}}{n}\right) ^{\frac{1}{2}},\ \ 1<p\leq 2.
\end{equation*}%
Taking the square in both sides, one has 
\begin{equation*}
\left( \frac{\sum_{i=1}^{n}\left\vert \left( x,y_{i}\right) \right\vert ^{p}%
}{n}\right) ^{\frac{2}{p}}\leq \frac{\sum_{i=1}^{n}\left\vert \left(
x,y_{i}\right) \right\vert ^{2}}{n},
\end{equation*}%
giving 
\begin{equation}
\left( \sum_{i=1}^{n}\left\vert \left( x,y_{i}\right) \right\vert
^{p}\right) ^{\frac{2}{p}}\leq n^{\frac{2}{p}-1}\sum_{i=1}^{n}\left\vert
\left( x,y_{i}\right) \right\vert ^{2}.  \label{2.10.23}
\end{equation}%
Using (\ref{2.10.23}) and the second inequality in (\ref{2.7.23}) we may
deduce the following result 
\begin{equation}
\sum_{i=1}^{n}\left\vert \left( x,y_{i}\right) \right\vert ^{2}\leq n^{\frac{%
2}{p}-1}\left\Vert x\right\Vert ^{2}\left( \sum\limits_{i,j=1}^{n}\left\vert
\left( y_{i},y_{j}\right) \right\vert ^{q}\right) ^{\frac{1}{q}},
\label{2.11.23}
\end{equation}%
for $1<p\leq 2,$ $\frac{1}{p}+\frac{1}{q}=1.$

Note that for $p=2$ $\left( q=2\right) $ we recapture (\ref{2.8.23}).
\end{remark}

\begin{remark}
Let us compare Bombieri's result 
\begin{equation}
\sum\limits_{i=1}^{n}\left\vert \left( x,y_{i}\right) \right\vert ^{2}\leq
\left\Vert x\right\Vert ^{2}\max\limits_{1\leq i\leq n}\left\{
\sum\limits_{j=1}^{n}\left\vert \left( y_{i},y_{j}\right) \right\vert
\right\}  \label{2.12.23}
\end{equation}%
with our general result (\ref{2.11.23}).

To do that, denote 
\begin{equation*}
M_{1}:=\max\limits_{1\leq i\leq n}\left\{ \sum\limits_{j=1}^{n}\left\vert
\left( y_{i},y_{j}\right) \right\vert \right\}
\end{equation*}%
and 
\begin{equation*}
M_{2}:=n^{\frac{2}{p}-1}\left( \sum\limits_{i,j=1}^{n}\left\vert \left(
y_{i},y_{j}\right) \right\vert ^{q}\right) ^{\frac{1}{q}},\ \ \ 1<p\leq 2,\ 
\frac{1}{p}+\frac{1}{q}=1.
\end{equation*}%
Consider the inner product space $H=\mathbb{R}$, $\left( x,y\right) =x\cdot
y $, $n=2$ and $y_{1}=a>0,$ $y_{2}=b>0.$ Then 
\begin{align*}
M_{1}& =\max \left\{ a^{2}+ab,ab+b^{2}\right\} =\left( a+b\right) \max
\left\{ a,b\right\} , \\
M_{2}& =2^{\frac{2}{p}-1}\left( a^{q}+b^{q}\right) ^{\frac{2}{q}}=2^{\frac{2%
}{p}-1}\left( a^{\frac{p}{p-1}}+b^{\frac{p}{p-1}}\right) ^{\frac{2\left(
p-1\right) }{p}},\ 1<p\leq 2.
\end{align*}%
Assume that $a=1,$ $b\in \left[ 0,1\right] ,$ $p\in (1,2].$ Utilizing Maple
6, one may easily see by plotting the function 
\begin{equation*}
f\left( b,p\right) :=M_{2}-M_{1}=2^{\frac{2}{p}-1}\left( 1+b^{\frac{p}{p-1}%
}\right) ^{\frac{2\left( p-1\right) }{p}}-1-b
\end{equation*}%
that it has positive and negative values in the box $\left[ 0,1\right]
\times \left[ 1,2\right] $, showing that the inequalities (\ref{2.11.23})
and (\ref{2.12.23}) cannot be compared. This means that one is not always
better than the other.
\end{remark}

\newpage

\section{Some Pre-Gr\"{u}ss Inequalities}

\subsection{Introduction}

Let $f,\,g$ be two functions defined and integrable on $\left[ a{,}b\right] $%
. Assume that%
\begin{equation*}
\varphi \leq f\left( x\right) \leq \Phi \ \ \text{and\ \ }\gamma \leq
g\left( x\right) \leq \Gamma
\end{equation*}%
for each $x\in \left[ {a,\,b}\right] $, where $\varphi ,\,\Phi ,\,\gamma
,\,\Gamma $ are given real constants. Then we have the following inequality
which is well known in the literature as the Gr\"{u}ss inequality (\cite[pp.
296]{1b.24})%
\begin{multline*}
\left\vert {{\frac{1}{{b-a}}}\int_{a}^{b}{f\left( x\right) g\left( x\right)
dx}}-\frac{1}{b-a}{\int_{a}^{b}}f\left( x\right) dx{\cdot }\frac{1}{b-a}{%
\int_{a}^{b}{g\left( x\right) dx}}\right\vert \\
\leq {\frac{1}{4}}\left\vert {\Phi -\varphi }\right\vert \cdot \left\vert {%
\Gamma -\gamma }\right\vert .
\end{multline*}%
In this inequality, G. Gr\"{u}ss has proven that the constant ${\frac{1}{4}}$
is the best possible in the sense that it cannot be replaced by a smaller
one, and this is achieved when%
\begin{equation*}
f\left( x\right) =g\left( x\right) =\func{sgn}\left( {x-{\frac{{a+b}}{2}}}%
\right) .
\end{equation*}

Recently, S.S. Dragomir proved the following Gr\"{u}ss' type inequality in
real or complex inner product spaces \cite{2b.24}.

\begin{theorem}
\label{ta.24}Let $\left( {H,\left\langle {\cdot ,\cdot }\right\rangle }%
\right) $ be an inner product space over $\mathbb{K}\ \left( \mathbb{K}=%
\mathbb{R},\mathbb{C}\right) $ and $e\in H,\,\left\Vert e\right\Vert =1$. If 
$\varphi ,\gamma ,\Phi ,\Gamma $ are real or complex numbers and $x,y$ are
vectors in $H$ such that the conditions%
\begin{equation*}
\func{Re}\left\langle {\Phi e-x,x-\varphi e}\right\rangle \geq 0\ \ \text{%
and\ \ }\func{Re}\left\langle {\Gamma e-x,x-\gamma e}\right\rangle \geq 0
\end{equation*}%
hold, then we have the inequality%
\begin{equation}
\left\vert {\left\langle {x,y}\right\rangle -\left\langle {x,e}\right\rangle
\,\left\langle {e,y}\right\rangle }\right\vert \leq {\frac{1}{4}}\left\vert {%
\Phi -\varphi }\right\vert \cdot \left\vert {\Gamma -\gamma }\right\vert .
\label{4.24}
\end{equation}%
The constant ${\frac{1}{4}}$ is best possible in sense that it cannot be
replaced by a smaller constant.
\end{theorem}

In \cite{3b.24}, by using the following lemmas

\begin{lemma}
\label{la.24}Let $x,e\in H$ with $\,\left\Vert e\right\Vert =1$ and $\delta
,\,\Delta \in \mathbb{K}$ with $\delta \neq \Delta $. Then%
\begin{equation*}
\func{Re}\left\langle {\Delta e-x,x-\delta e}\right\rangle \geq 0
\end{equation*}%
if and only if%
\begin{equation*}
\left\Vert {x-{\frac{{\delta +\Delta }}{2}}e}\right\Vert \leq {\frac{1}{2}}%
\left\vert {\Delta -\delta }\right\vert .
\end{equation*}
\end{lemma}

and

\begin{lemma}
\label{lb.24}Let $x,e\in H$ with $\,\left\Vert e\right\Vert =1$. Then one
has the following representation%
\begin{equation*}
0\leq \left\Vert x\right\Vert ^{2}-\left\vert {\left\langle {x,\,e}%
\right\rangle }\right\vert ^{2}=\inf_{\lambda \in K}\left\Vert {x-\lambda e}%
\right\Vert ^{2}.
\end{equation*}
\end{lemma}

the author gave an alternative proof for (\ref{4.24}) and also obtained the
following refinement of it, namely

\begin{theorem}
\label{t3.1.5.24}Let $\left( H,\left\langle \cdot ,\cdot \right\rangle
\right) $ be an inner product space over $\mathbb{K\ }\left( \mathbb{K}=%
\mathbb{R}\text{,}\mathbb{C}\right) $ and $e\in H,\left\Vert e\right\Vert
=1. $ If $\varphi ,\gamma ,\Phi ,\Gamma $ are real or complex numbers and $%
x,y$ are vectors in $H$ such that either the conditions%
\begin{equation*}
\func{Re}\left\langle \Phi e-x,x-\varphi e\right\rangle \geq 0\text{, }\func{%
Re}\left\langle \Gamma e-y,y-\gamma e\right\rangle \geq 0,
\end{equation*}%
or equivalently, 
\begin{equation*}
\left\Vert x-\frac{\varphi +\Phi }{2}\cdot e\right\Vert \leq \frac{1}{2}%
\left\vert \Phi -\varphi \right\vert ,\text{ \ }\left\Vert y-\frac{\gamma
+\Gamma }{2}\cdot e\right\Vert \leq \frac{1}{2}\left\vert \Gamma -\gamma
\right\vert ,
\end{equation*}%
hold, then we have the inequality 
\begin{align*}
& \left\vert \left\langle x,y\right\rangle -\left\langle x,e\right\rangle
\left\langle e,y\right\rangle \right\vert \\
& \leq \frac{1}{4}\left\vert \Phi -\varphi \right\vert \cdot \left\vert
\Gamma -\gamma \right\vert -\left[ \func{Re}\left\langle \Phi e-x,x-\varphi
e\right\rangle \right] ^{\frac{1}{2}}\left[ \func{Re}\left\langle \Gamma
e-y,y-\gamma e\right\rangle \right] ^{\frac{1}{2}} \\
& \leq \left( \frac{1}{4}\left\vert \Phi -\varphi \right\vert \cdot
\left\vert \Gamma -\gamma \right\vert \right) .
\end{align*}%
The constant $\frac{1}{4}$ is best possible.
\end{theorem}

Further, as a generalization for orthonormal families of vectors in inner
product spaces, S.S. Dragomir proved, in \cite{4b.24}, the following reverse
of Bessel's inequality:

\begin{theorem}
\label{tb.24}Let $\left\{ {e_{i}}\right\} ,\,i\in I$ be a family of
orthonormal vectors in $H$, $F$ a finite part of $I$, $\varphi _{i},\,\Phi
_{i}\in \mathbb{K}$, $i\in F$ and $x$ a vector in $H$ such that either the
condition%
\begin{equation*}
\func{Re}\left\langle {\sum\limits_{i\in F}{\Phi _{i}e_{i}}%
-x,\,x-\sum\limits_{i\in F}{\varphi _{i}e_{i}}}\right\rangle \geq 0,
\end{equation*}%
or equivalently,%
\begin{equation*}
\left\Vert {x-\sum\limits_{i\in F}{{\frac{{\Phi _{i}+\varphi _{i}}}{2}}}e_{i}%
}\right\Vert \leq {\frac{1}{2}}\left( {\sum\limits_{i\in F}{\left\vert {\Phi
_{i}-\varphi _{i}}\right\vert ^{2}}}\right) ^{{\frac{1}{2}}},
\end{equation*}%
holds, then we have the following reverse of Bessel's inequality%
\begin{multline}
\quad \left\Vert x\right\Vert ^{2}-\sum\limits_{i\in F}{\left\vert {%
\left\langle {x,\,e_{i}}\right\rangle }\right\vert }^{2}  \label{10.24} \\
\leq {\frac{1}{4}}\sum\limits_{i\in F}{\left\vert {\Phi _{i}-\varphi _{i}}%
\right\vert ^{2}}-\sum\limits_{i\in F}{\left\vert {{\frac{{\varphi _{i}+\Phi
_{i}}}{2}}-\left\langle {x,\,e_{i}}\right\rangle }\right\vert ^{2}.}\quad
\end{multline}%
The constant ${\frac{1}{4}}$ is best possible.
\end{theorem}

The corresponding Gr\"{u}ss type inequality is embodied in the following
theorem:

\begin{theorem}
\label{tc.24}Let $\left\{ e_{i}\right\} _{i\in I}$ be a family of
orthornormal vectors in $H,$ $F$ a finite part of $I,$ $\phi _{i},\gamma
_{i},\Phi _{i},\Gamma _{i}\in \mathbb{R}$ $\left( i\in F\right) $, and $%
x,y\in H.$ If either%
\begin{align*}
\func{Re}\left\langle \sum_{i=1}^{n}\Phi _{i}e_{i}-x,x-\sum_{i=1}^{n}\phi
_{i}e_{i}\right\rangle & \geq 0,\  \\
\func{Re}\left\langle \sum_{i=1}^{n}\Gamma
_{i}e_{i}-y,y-\sum_{i=1}^{n}\gamma _{i}e_{i}\right\rangle & \geq 0,
\end{align*}%
or equivalently,%
\begin{align*}
\left\Vert x-\sum_{i\in F}\frac{\Phi _{i}+\phi _{i}}{2}e_{i}\right\Vert &
\leq \frac{1}{2}\left( \sum_{i\in F}\left\vert \Phi _{i}-\phi
_{i}\right\vert ^{2}\right) ^{\frac{1}{2}}, \\
\left\Vert y-\sum_{i\in F}\frac{\Gamma _{i}+\gamma _{i}}{2}e_{i}\right\Vert
& \leq \frac{1}{2}\left( \sum_{i\in F}\left\vert \Gamma _{i}-\gamma
_{i}\right\vert ^{2}\right) ^{\frac{1}{2}},
\end{align*}%
hold true, then%
\begin{align*}
0& \leq \left\vert \left\langle x,y\right\rangle -\sum_{i=1}^{n}\left\langle
x,e_{i}\right\rangle \left\langle e_{i},y\right\rangle \right\vert \\
& \leq \frac{1}{4}\left( \sum_{i=1}^{n}\left\vert \Phi _{i}-\phi
_{i}\right\vert ^{2}\right) ^{\frac{1}{2}}\cdot \left(
\sum_{i=1}^{n}\left\vert \Gamma _{i}-\gamma _{i}\right\vert ^{2}\right) ^{%
\frac{1}{2}} \\
& \ \ \ \ \ \ \ \ \ \ \ \ \ \ \ \ \ -\sum_{i\in F}\left\vert \frac{\Phi
_{i}+\phi _{i}}{2}-\left\langle x,e_{i}\right\rangle \right\vert \left\vert 
\frac{\Gamma _{i}+\gamma _{i}}{2}-\left\langle y,e_{i}\right\rangle
\right\vert \\
& \left( \leq \frac{1}{4}\left( \sum_{i=1}^{n}\left\vert \Phi _{i}-\phi
_{i}\right\vert ^{2}\right) ^{\frac{1}{2}}\cdot \left(
\sum_{i=1}^{n}\left\vert \Gamma _{i}-\gamma _{i}\right\vert ^{2}\right) ^{%
\frac{1}{2}}\right) .
\end{align*}%
The constant $\frac{1}{4}$ is best possible in the sense that it cannot be
replaced by a smaller constant.
\end{theorem}

The main aim here is to provide some similar inequalities which, giving
refinements of the usual Gr\"{u}ss' inequality, are known in the literature
as pre-Gr\"{u}ss type inequalities. Applications for Lebesgue integrals in
general measure spaces are also given.

\subsection{Pre-Gr\"{u}ss Inequalities in Inner Product Spaces}

We start with the following result \cite{24NSSD}:

\begin{theorem}
\label{t1.24}Let $\left( {H,\left\langle {\cdot ,\cdot }\right\rangle }%
\right) $ be an inner product space over $\mathbb{K},\ \left( \mathbb{K}=%
\mathbb{R},\mathbb{C}\right) $ and $e\in H,\,\left\Vert e\right\Vert =1$. If 
$\varphi ,{\Phi }$ are real or complex numbers and $x,y$ are vectors in $H$
such that either the condition%
\begin{equation*}
\func{Re}\left\langle {\Phi e-x,x-\varphi e}\right\rangle \geq 0,\text{\ }
\end{equation*}%
or equivalently,%
\begin{equation}
\left\Vert {x-{\frac{{\varphi +\Phi }}{2}}e}\right\Vert \leq {\frac{1}{2}}%
\left\vert {\Phi -\varphi }\right\vert ,  \label{11.24}
\end{equation}%
holds true, then we have the inequalities%
\begin{equation}
\left\vert {\left\langle {x,y}\right\rangle -\left\langle {x,e}\right\rangle
\,\left\langle {e,y}\right\rangle }\right\vert \leq {\frac{1}{2}}\left\vert {%
\Phi -\varphi }\right\vert \cdot \sqrt{\left( {\left\Vert y\right\Vert
^{2}-\left\vert {\left\langle {y,\,e}\right\rangle }\right\vert ^{2}}\right) 
}  \label{12.24}
\end{equation}%
and%
\begin{multline}
\left\vert {\left\langle {x,y}\right\rangle -\left\langle {x,e}\right\rangle
\,\left\langle {e,y}\right\rangle }\right\vert  \label{13.24} \\
\leq {\frac{1}{2}}\left\vert {\Phi -\varphi }\right\vert \cdot \left\Vert
y\right\Vert -\left( {\func{Re}\left\langle {\Phi e-x,x-\varphi e}%
\right\rangle }\right) ^{{\frac{1}{2}}}\cdot \left\vert {\left\langle {y,\,e}%
\right\rangle }\right\vert .
\end{multline}
\end{theorem}

\begin{proof}
It is obvious that:%
\begin{equation*}
\left\langle {x,y}\right\rangle -\left\langle {x,e}\right\rangle
\left\langle {e,y}\right\rangle =\left\langle {x-\left\langle {x,e}%
\right\rangle e,y-\left\langle {y,e}\right\rangle e}\right\rangle .
\end{equation*}%
Using Schwarz's inequality in inner product spaces $\left\vert {\left\langle 
{u,\,v}\right\rangle }\right\vert \leq \left\Vert u\right\Vert \cdot
\left\Vert v\right\Vert $ for the vectors ${x-\left\langle {x,e}%
\right\rangle e}$ and ${y-\left\langle {y,e}\right\rangle e,}$ we deduce:%
\begin{equation}
\left\vert \left\langle {x,y}\right\rangle -\left\langle {x,e}\right\rangle
\left\langle {e,y}\right\rangle \right\vert ^{2}\leq \left( {\left\Vert
x\right\Vert ^{2}-\left\vert {\left\langle {x,\,e}\right\rangle }\right\vert
^{2}}\right) \cdot \left( {\left\Vert y\right\Vert ^{2}-\left\vert {%
\left\langle {y,\,e}\right\rangle }\right\vert ^{2}}\right) {.}
\label{15a.24}
\end{equation}%
Now, the inequality (\ref{12.24}) is a simple consequence of (\ref{4.24})
for $x=y$, or of Lemma \ref{lb.24} and (\ref{11.24}).

Since (see for instance \cite{2b.24}), 
\begin{multline}
\quad \left\Vert x\right\Vert ^{2}-\left\vert {\left\langle {x,\,e}%
\right\rangle }\right\vert ^{2}  \label{16.24} \\
=\func{Re}\left( {\left( {\Phi -\left\langle {x,\,e}\right\rangle }\right)
\cdot \left( {\left\langle {e,\,x}\right\rangle -\bar{\varphi}}\right) }%
\right) -\func{Re}\left\langle {\Phi e-x,\,x-\varphi e}\right\rangle ,\quad
\end{multline}%
then making use of the elementary inequality $4\func{Re}\left( {a\bar{b}}%
\right) \leq \left\vert {a+b}\right\vert ^{2}$ with $a,b\in \mathbb{K\ }%
\,\left( \mathbb{K}=\mathbb{R}{,}\mathbb{C}\right) ,$ we can state that%
\begin{equation}
\func{Re}\left( {\left( {\Phi -\left\langle {x,\,e}\right\rangle }\right)
\cdot \left( {\left\langle {e,\,x}\right\rangle -\bar{\varphi}}\right) }%
\right) \leq {\frac{1}{4}}\left\vert {\Phi -\varphi }\right\vert ^{2}.
\label{17.24}
\end{equation}%
Using (\ref{16.24}) and (\ref{17.24}) we have%
\begin{equation}
\left\Vert x\right\Vert ^{2}-\left\vert {\left\langle {x,\,e}\right\rangle }%
\right\vert ^{2}\leq \left( {{\frac{1}{2}}\left\vert {\Phi -\varphi }%
\right\vert }\right) ^{2}-\left( {\left( {\func{Re}\left\langle {\Phi
e-x,\,x-\varphi e}\right\rangle }\right) ^{{\frac{1}{2}}}}\right) ^{2}.
\label{18.24}
\end{equation}%
Taking into account the inequalities (\ref{15a.24}) and (\ref{18.24}), we
get that%
\begin{multline*}
\left\vert {\left\langle {x-\left\langle {x,\,e}\right\rangle
e,\,y-\left\langle {y,\,e}\right\rangle e}\right\rangle }\right\vert ^{2} \\
\leq \left( {\left( {{\frac{1}{2}}\left\vert {\Phi -\varphi }\right\vert }%
\right) ^{2}-\left( {\left( {\func{Re}\left\langle {\Phi e-x,x-\varphi e}%
\right\rangle }\right) ^{{\frac{1}{2}}}}\right) ^{2}}\right) \cdot \left( {%
\left\Vert y\right\Vert ^{2}-\left\vert {\left\langle {y,e}\right\rangle }%
\right\vert ^{2}}\right) .
\end{multline*}%
Finally, using the elementary inequality for positive real numbers: 
\begin{equation}
\left( {m^{2}-n^{2}}\right) \cdot \left( {p^{2}-q^{2}}\right) \leq \left( {%
mp-nq}\right) ^{2},  \label{20.24}
\end{equation}%
we have: 
\begin{multline*}
\left( {\left( {{\frac{1}{2}}\left\vert {\Phi -\varphi }\right\vert }\right)
^{2}-\left( {\left( {\func{Re}\left\langle {\Phi e-x,\,x-\varphi e}%
\right\rangle }\right) ^{{\frac{1}{2}}}}\right) ^{2}}\right) \cdot \left( {%
\left\Vert y\right\Vert ^{2}-\left\vert {\left\langle {y,e}\right\rangle }%
\right\vert ^{2}}\right) \\
\leq \left( {{\frac{1}{2}}\left\vert {\Phi -\varphi }\right\vert \cdot
\left\Vert y\right\Vert -\left( {\func{Re}\left\langle {\Phi e-x,\,x-\varphi
e}\right\rangle }\right) ^{{\frac{1}{2}}}\cdot \left\vert {\left\langle {%
y,\,e}\right\rangle }\right\vert }\right) ^{2},
\end{multline*}%
giving the desired inequality (\ref{13.24}).
\end{proof}

A similar version for Bessel's inequality is incorporated in the following
theorem \cite{24NSSD}:

\begin{theorem}
\label{t2.24}Let $\left\{ {e_{i}}\right\} _{i\in I},$ be a family of
orthonormal vectors in $H$, $F$ a finite part of $I$, $\varphi _{i},\,\Phi
_{i}\in \mathbb{K}$, $i\in F$ and $x,y$ vectors in $H$ such that either the
condition%
\begin{equation*}
\func{Re}\left\langle {\sum\limits_{i\in F}{\Phi _{i}e_{i}}%
-x,\,x-\sum\limits_{i\in F}{\varphi _{i}e_{i}}}\right\rangle \geq 0,
\end{equation*}%
or equivalently,%
\begin{equation*}
\left\Vert {x-\sum\limits_{i\in F}{{\frac{{\Phi _{i}+\varphi _{i}}}{2}}}e_{i}%
}\right\Vert \leq {\frac{1}{2}}\left( {\sum\limits_{i\in F}{\left\vert {\Phi
_{i}-\varphi _{i}}\right\vert ^{2}}}\right) ^{{\frac{1}{2}}}
\end{equation*}%
holds. Then we have inequalities%
\begin{multline}
\left\vert {\left\langle {x,y}\right\rangle -\sum\limits_{i\in F}{%
\left\langle {x,e_{i}}\right\rangle \,\left\langle {e_{i},y}\right\rangle }}%
\right\vert  \label{23.24} \\
\leq {\frac{1}{2}}\left( {\sum\limits_{i\in F}{\left\vert {\Phi _{i}-\varphi
_{i}}\right\vert ^{2}}}\right) ^{{\frac{1}{2}}}\sqrt{\left( {\left\Vert
y\right\Vert ^{2}-\sum\limits_{i\in F}{\left\vert {\left\langle {y,\,e_{i}}%
\right\rangle }\right\vert ^{2}}}\right) }
\end{multline}%
and%
\begin{multline}
\left\vert {\left\langle {x,y}\right\rangle -\sum\limits_{i\in F}{%
\left\langle {x,e_{i}}\right\rangle \,\left\langle {e_{i},y}\right\rangle }}%
\right\vert \leq {\frac{1}{2}}\left( {\sum\limits_{i\in F}{\left\vert {\Phi
_{i}-\varphi _{i}}\right\vert ^{2}}}\right) ^{{\frac{1}{2}}}\cdot \left\Vert
y\right\Vert  \label{24.24} \\
-\left( {\sum\limits_{i\in F}{\left\vert {{\frac{{\Phi _{i}+\varphi _{i}}}{2}%
}-\left\langle {x,\,e_{i}}\right\rangle }\right\vert ^{2}}}\right) ^{{\frac{1%
}{2}}}\cdot \left( {\sum\limits_{i\in F}{\left\vert {\left\langle {y,\,e_{i}}%
\right\rangle }\right\vert ^{2}}}\right) ^{{\frac{1}{2}}}.
\end{multline}
\end{theorem}

\begin{proof}
It is obvious (see for example \cite{4b.24}) that:%
\begin{equation*}
\left\langle {x,\,y}\right\rangle -\sum\limits_{i\in F}{\left\langle {%
x,\,e_{i}}\right\rangle \left\langle {e_{i},\,y}\right\rangle }=\left\langle 
{x-\sum\limits_{i\in F}{\left\langle {x,\,e_{i}}\right\rangle e_{i}}%
,\,y-\sum\limits_{i\in F}{\left\langle {y,\,e_{i}}\right\rangle e_{i}}}%
\right\rangle .
\end{equation*}%
Using Schwarz's inequality in inner product spaces, we have:%
\begin{align}
& \left\vert {\left\langle {x-\sum\limits_{i\in F}{\left\langle {x,\,e_{i}}%
\right\rangle e_{i}},\,y-\sum\limits_{i\in F}{\left\langle {y,\,e_{i}}%
\right\rangle e_{i}}}\right\rangle }\right\vert ^{2}  \label{26.24} \\
& \leq \left\Vert {x-\sum\limits_{i\in F}{\left\langle {x,\,e_{i}}%
\right\rangle e_{i}}}\right\Vert ^{2}\cdot \left\Vert {x-\sum\limits_{i\in F}%
{\left\langle {y,\,e_{i}}\right\rangle e_{i}}}\right\Vert ^{2}  \notag \\
& =\left( {\left\Vert x\right\Vert ^{2}-\sum\limits_{i\in F}{\left\vert {%
\left\langle {x,\,e_{i}}\right\rangle }\right\vert ^{2}}}\right) \cdot
\left( {\left\Vert y\right\Vert ^{2}-\sum\limits_{i\in F}{\left\vert {%
\left\langle {y,\,e_{i}}\right\rangle }\right\vert ^{2}}}\right) .  \notag
\end{align}%
In a similar manner to the one in the proof of Theorem \ref{t1.24} we may
conclude that (\ref{23.24}) holds true.

Now, using (\ref{10.24}) and (\ref{26.24}) we also have:%
\begin{multline*}
\left\vert {\left\langle {x-\sum\limits_{i\in F}{\left\langle {x,\,e_{i}}%
\right\rangle e_{i}},\,y-\sum\limits_{i\in F}{\left\langle {y,\,e_{i}}%
\right\rangle e_{i}}}\right\rangle }\right\vert ^{2} \\
\leq \left( {{\frac{1}{2}}\left( {\left( {\sum\limits_{i\in F}{\left\vert {%
\Phi _{i}-\varphi _{i}}\right\vert ^{2}}}\right) ^{{\frac{1}{2}}}}\right)
^{2}-\left( {\left( {\sum\limits_{i\in F}{\left\vert {{\frac{{\varphi
_{i}+\Phi _{i}}}{2}}-\left\langle {x,\,e_{i}}\right\rangle }\right\vert ^{2}}%
}\right) ^{{\frac{1}{2}}}}\right) ^{2}}\right) \\
\times \left( {\left\Vert y\right\Vert ^{2}-\sum\limits_{i\in F}{\left\vert {%
\left\langle {y,\,e_{i}}\right\rangle }\right\vert ^{2}}}\right) .
\end{multline*}%
Finally, utilizing the elementary inequality (\ref{20.24}), we have%
\begin{multline}
\left( {{\frac{1}{2}}\left( {\left( {\sum\limits_{i\in F}{\left\vert {\Phi
_{i}-\varphi _{i}}\right\vert ^{2}}}\right) ^{{\frac{1}{2}}}}\right)
^{2}-\left( {\left( {\sum\limits_{i\in F}{{}{{\frac{{\varphi _{i}+\Phi _{i}}%
}{2}}-\left\langle {x,\,e_{i}}\right\rangle }^{2}}}\right) ^{{\frac{1}{2}}}}%
\right) ^{2}}\right)  \label{28.24} \\
\times \left( {\left\Vert y\right\Vert ^{2}-\sum\limits_{i\in F}{\left\vert {%
\left\langle {y,\,e_{i}}\right\rangle }\right\vert ^{2}}}\right) \leq \left( 
{{\frac{1}{2}}\left( {\sum\limits_{i\in F}{\left\vert {\Phi _{i}-\varphi _{i}%
}\right\vert ^{2}}}\right) ^{{\frac{1}{2}}}\cdot \left\Vert y\right\Vert ^{2}%
}\right. \\
\left. {-\left( {\sum\limits_{i\in F}{\left\vert {{\frac{{\varphi _{i}+\Phi
_{i}}}{2}}-\left\langle {x,\,e_{i}}\right\rangle }\right\vert ^{2}}}\right)
^{{\frac{1}{2}}}\cdot \sum\limits_{i\in F}{\left\vert {\left\langle {%
y,\,e_{i}}\right\rangle }\right\vert ^{2}}}\right) ^{2},
\end{multline}%
which gives the desired result (\ref{24.24}).
\end{proof}

\subsection{Applications for Integrals}

Let $\left( {\Omega ,\,\Sigma ,\,\mu }\right) $ be a measure space
consisting of a set $\Omega $, $\Sigma $ a $\sigma -$algebra of parts and $%
\mu $ a countably additive and positive measure on $\Sigma $ with values in $%
\mathbb{R}\cup \left\{ \infty \right\} $. Denote by $L^{2}\left( {\Omega ,}%
\mathbb{K}\right) $ the Hilbert space of all real or complex valued
functions $f$ defined on $\Omega $ and 2--integrable on $\Omega $, i. e. 
\begin{equation*}
\int_{\Omega }{\left\vert {f\left( s\right) }\right\vert ^{2}d\mu \left(
s\right) <\infty .}
\end{equation*}%
The following proposition holds \cite{24NSSD}.

\begin{proposition}
\label{c1.24}If $f,g,h\in L^{2}\left( {\Omega ,}\mathbb{K}\right) $ and $%
\varphi ,\,\Phi \in \mathbb{K}$ are such that $\int_{\Omega }{\left\vert {%
h\left( s\right) }\right\vert ^{2}d\mu \left( s\right) =1}$ and, either%
\begin{equation}
\int\limits_{\Omega }{\func{Re}\left( {\left( {\Phi h\left( s\right)
-f\left( s\right) }\right) \,\left( {\bar{f}\left( s\right) -}\overline{{%
\varphi }}{\bar{h}\left( s\right) }\right) }\right) d\mu \left( s\right) }%
\geq 0,  \label{29.24}
\end{equation}%
or equivalently,%
\begin{equation*}
\left( {\int_{\Omega }{\left\vert {f\left( s\right) -{\frac{{\Phi +\varphi }%
}{2}}h\left( s\right) }\right\vert ^{2}d\mu \left( s\right) }}\right) ^{{%
\frac{1}{2}}}\leq {\frac{1}{2}}\left\vert {\Phi -\varphi }\right\vert
\end{equation*}%
holds, then we have the inequalities%
\begin{multline*}
\left\vert {\int_{\Omega }{f\left( s\right) \bar{g}\left( s\right) d\mu
\left( s\right) }-\int_{\Omega }{f\left( s\right) \bar{h}\left( s\right)
d\mu \left( s\right) }\int_{\Omega }{h\left( s\right) \bar{g}\left( s\right)
d\mu \left( s\right) }}\right\vert \\
\leq {\frac{1}{2}}\left\vert {\Phi -\varphi }\right\vert \cdot \sqrt{\left( {%
\int_{\Omega }}\left\vert {g\left( s\right) }\right\vert ^{2}{{d\mu \left(
s\right) }-\left\vert {\int_{\Omega }{h\left( s\right) \bar{g}\left(
s\right) d\mu \left( s\right) }}\right\vert ^{2}}\right) }
\end{multline*}%
and%
\begin{multline*}
\left\vert {\int_{\Omega }{f\left( s\right) \bar{g}\left( s\right) d\mu
\left( s\right) }-\int_{\Omega }{f\left( s\right) \bar{h}\left( s\right)
d\mu \left( s\right) }\int_{\Omega }{h\left( s\right) \bar{g}\left( s\right)
d\mu \left( s\right) }}\right\vert \\
\leq {\frac{1}{2}}\left\vert {\Phi -\varphi }\right\vert \cdot \left( {%
\int_{\Omega }}\left\vert {g\left( s\right) }\right\vert ^{2}{d\mu \left(
s\right) }\right) ^{\frac{1}{2}} \\
-\left( {\int_{\Omega }{\func{Re}\left( {\left( {\Phi h\left( s\right)
-f\left( s\right) }\right) \left( {h\left( s\right) \bar{f}\left( s\right)
-\varphi h\left( s\right) }\right) }\right) d\mu \left( s\right) }}\right) ^{%
{\frac{1}{2}}}\left\vert {\int_{\Omega }{h\left( s\right) \bar{g}\left(
s\right) d\mu \left( s\right) }}\right\vert .
\end{multline*}
\end{proposition}

\begin{proof}
The proof follows by Theorem \ref{t1.24} on choosing $H=L^{2}\left( {\Omega
,\,K}\right) $ with the inner product 
\begin{equation*}
\left\langle {f,g}\right\rangle =\int\limits_{{\Omega }}f\left( s\right) 
\bar{g}\left( s\right) d\mu \left( s\right) .
\end{equation*}
\end{proof}

\begin{remark}
We observe that a sufficient condition for (\ref{29.24}) to hold is:%
\begin{equation}
\func{Re}{\left( {\Phi h\left( s\right) -f\left( s\right) }\right) \,\left( {%
\bar{f}\left( s\right) -}\overline{{\varphi }}{\bar{h}\left( s\right) }%
\right) \geq 0,}  \label{33.24}
\end{equation}%
for $\mu -$a.e. $s\in {\Omega .}$

If the functions are real-valued, then, for ${\Phi }$ and $\varphi $ real
numbers, a sufficient condition for (\ref{33.24}) to hold is%
\begin{equation*}
{\Phi h\left( s\right) \geq f\left( s\right) \geq \varphi h\left( s\right) }
\end{equation*}%
for $\mu -$a.e. $s\in {\Omega .}$

In this way we can see the close connection that exists between the
classical Gr\"{u}ss inequality and the results obtained above.
\end{remark}

Now, consider the family $\left\{ f_{i}\right\} _{i\in I}$ of functions in $%
L^{2}\left( \Omega ,\mathbb{K}\right) $ with the properties that%
\begin{equation*}
\int_{\Omega }f_{i}\left( s\right) \overline{f_{j}}\left( s\right) d\mu
\left( s\right) =\delta _{ij},\ \ \ i,j\in I,
\end{equation*}%
where $\delta _{ij}$ is $0$ if $i\neq j$ and $\delta _{ij}=1$ if $i=j.$ $%
\left\{ f_{i}\right\} _{i\in I}$ is an orthornormal family in $L^{2}\left(
\Omega ,\mathbb{K}\right) .$

The following proposition holds \cite{24NSSD}.

\begin{proposition}
\label{p5.1.24}Let $\left\{ f_{i}\right\} _{i\in I}$ be an orthornormal
family of functions in $L^{2}\left( \Omega ,\mathbb{K}\right) ,$ $F$ a
finite subset of $I,$ $\phi _{i},\Phi _{i}\in \mathbb{K}$ $\left( i\in
F\right) $ and $f\in L^{2}\left( \Omega ,\mathbb{K}\right) ,$ such that
either%
\begin{equation}
\int_{\Omega }\func{Re}\left[ \left( \sum_{i\in F}\Phi _{i}f_{i}\left(
s\right) -f\left( s\right) \right) \left( \overline{f}\left( s\right)
-\sum_{i\in F}\overline{\phi _{i}}\text{ }\overline{f_{i}}\left( s\right)
\right) \right] d\mu \left( s\right) \geq 0  \label{34.24}
\end{equation}%
or equivalently,%
\begin{equation*}
\int_{\Omega }\left\vert f\left( s\right) -\sum_{i\in F}\frac{\Phi _{i}+\phi
_{i}}{2}f_{i}\left( s\right) \right\vert ^{2}d\mu \left( s\right) \leq \frac{%
1}{4}\sum_{i\in F}\left\vert \Phi _{i}-\phi _{i}\right\vert ^{2}.
\end{equation*}%
holds. Then we have the inequalities%
\begin{multline*}
\left\vert \int_{\Omega }f\left( s\right) \overline{g\left( s\right) }d\mu
\left( s\right) -\sum_{i\in F}\int_{\Omega }f\left( s\right) \overline{f_{i}}%
\left( s\right) d\mu \left( s\right) \int_{\Omega }f_{i}\left( s\right) 
\overline{g\left( s\right) }d\mu \left( s\right) \right\vert \\
\leq \frac{1}{2}\left( \sum_{i\in F}\left\vert \Phi _{i}-\phi
_{i}\right\vert ^{2}\right) ^{\frac{1}{2}}\left( \int_{\Omega }\left\vert
g\left( s\right) \right\vert ^{2}d\mu \left( s\right) -\sum_{i\in
F}\left\vert \int_{\Omega }g\left( s\right) \overline{f_{i}\left( s\right) }%
d\mu \left( s\right) \right\vert ^{2}\right) ^{\frac{1}{2}}
\end{multline*}%
and%
\begin{multline*}
\left\vert \int_{\Omega }f\left( s\right) \overline{g\left( s\right) }d\mu
\left( s\right) -\sum_{i\in F}\int_{\Omega }f\left( s\right) \overline{f_{i}}%
\left( s\right) d\mu \left( s\right) \int_{\Omega }f_{i}\left( s\right) 
\overline{g\left( s\right) }d\mu \left( s\right) \right\vert \\
\leq \frac{1}{2}\left( \sum_{i\in F}\left\vert \Phi _{i}-\phi
_{i}\right\vert ^{2}\right) ^{\frac{1}{2}}\left( \int_{\Omega }\left\vert
g\left( s\right) \right\vert ^{2}d\mu \left( s\right) \right) ^{\frac{1}{2}}
\\
-\left( \sum_{i\in F}\left\vert \frac{\Phi _{i}+\phi _{i}}{2}-\int_{\Omega
}f\left( s\right) \overline{f_{i}}\left( s\right) d\mu \left( s\right)
\right\vert ^{2}\right) ^{\frac{1}{2}}\left( \sum_{i\in F}\left\vert
\int_{\Omega }f\left( s\right) \overline{f_{i}}\left( s\right) d\mu \left(
s\right) \right\vert ^{2}\right) ^{\frac{1}{2}}.
\end{multline*}
\end{proposition}

The proof is obvious by Theorem \ref{t1.24} and we omit the details.

\begin{remark}
In the real case, we observe that a sufficient condition for (\ref{34.24})
to hold is%
\begin{equation*}
\sum_{i\in F}\Phi _{i}f_{i}\left( s\right) \geq f\left( s\right) \geq
\sum_{i\in F}\varphi _{i}f_{i}\left( s\right)
\end{equation*}%
for $\mu -$a.e. $s\in {\Omega .}$
\end{remark}

\end{document}